\documentclass[11pt,a4paper, dvips,twoside]{book}
\usepackage{amsmath,amsthm,amssymb,amsfonts,upref,cite,epsf,color}
\usepackage{pst-all,pstricks,pst-node,pst-text,pst-3d}
\usepackage{hyperref}
\usepackage{pstricks-add}
\usepackage{multido}
\usepackage{graphicx,epsfig,rotating} 
\usepackage{psfrag,texdraw,bm}
\usepackage{nomencl}
\usepackage{amsmath}
\usepackage{subfigure}
\usepackage[titletoc]{appendix}
\usepackage{amssymb}
\usepackage{xspace}
\usepackage{dsfont}
\usepackage{trsym}
\usepackage{pifont}
\usepackage{listings,color}   
\usepackage[T1]{fontenc}
\usepackage{textcomp}

\newtheorem{theorem}{Theorem}[section]
\newtheorem{definition}{Definition}[section]
\newtheorem{lemma}[theorem]{Lemma}
\newtheorem{example}[theorem]{Example}

\newtheorem{postulate}[theorem]{Postulate}
\newtheorem{corollary}[theorem]{Corollary}

\DeclareMathOperator{\supp}{supp}
\DeclareMathOperator*{\argmax}{argmax}
\DeclareMathOperator*{\spark}{spark}
\DeclareMathOperator*{\rank}{rank}
\DeclareMathOperator*{\vectorizesymb}{vec}

\DeclareMathOperator*{\arginf}{arginf}
\DeclareMathOperator*{\trace}{Tr}
\DeclareMathOperator*{\deto}{det}

\DeclareMathOperator{\linspan}{span}
\DeclareMathOperator{\closure}{cl}
\def\ML_est{\hat{\mathbf{x}}_{\text{ML}}}
\newcommand{\CRBfull}{{C}ram\'{e}r--{R}ao bound\xspace}
\newcommand{\minachievevar}{L_{\mathcal{M}}}
\newcommand{\minvarproblem}{\big( \mathcal{E}, \mathbf{c}(\cdot),\mathbf{x}_{0} \big)}
\newcommand{\minvarproblemscalar}{\big( \mathcal{E},c(\cdot),\mathbf{x}_{0} \big)}
\newcommand{\be}{\begin{equation}}
\newcommand{\ee}{\end{equation}}
\newcommand{\ist}{\hspace*{.2mm}}
\newcommand{\rmv}{\hspace*{-.2mm}}

\newcommand{\detm}[1]{\det\{#1\}}
\newcommand{\tracem}[1]{\trace\{#1\}}
\newcommand{\vectorize}[1]{\vectorizesymb\{#1\}}
\newcommand{\detmb}[1]{\det \big \{#1 \big\}}

\newcommand{\equref}[1]{(\ref{#1})}
\newcommand{\vecestproblem}{\left(  \mathcal{X}, f(\mathbf{y};\mathbf{x}), \mathbf{g}(\cdot) \right)}
\newcommand{\scalarestproblem}{\left(  \mathcal{X}, f(\mathbf{y};\mathbf{x}),g(\cdot) \right)}
\definecolor{hl}{rgb}{0.7,0,0}
\definecolor{fr}{rgb}{0.9,0.5,0}
\definecolor{ret}{rgb}{0,.5,0}

\makenomenclature
\allowdisplaybreaks
\begin{document}
\begin{titlepage}
\begin{center}
{\bf \Large Dissertation}\\[25mm]
{\bf \huge An RKHS Approach to Estimation} \\ [5mm]
{\bf \huge with Sparsity Constraints} \\ [30mm]
ausgef\"uhrt zum Zwecke der Erlangung des akademischen Grades eines \\[1mm]
Doktors der technischen Wissenschaften\\[15mm]
unter der Leitung von\\[1mm]
Ao. Univ.-Prof. Dipl.-Ing. Dr. Franz Hlawatsch\\[1mm]
Insitute of Telecommunications\\[15mm]
eingereicht an der Technischen Universit\"at Wien \\[1mm]
Fakult\"at f\"ur Elektrotechnik\\[15mm]
von\\[0mm]
Dipl.-Ing.\ Alexander Jung\\
Vorderer {\"O}dhof 1/1\\
A-3062 Kirchstetten\\[15mm]
\flushleft{Kirchstetten, Mai 2011%
\hspace{40mm} \underline{\hspace{50mm}}}
\end{center}
\end{titlepage}
\pagenumbering{roman}%


\setcounter{page}{2} 

\chapter*{Acknowledgements}

First of all, I would like to thank Prof.\ Franz Hlawatsch, whose spirit of work will always be a motivating but unreachable ideal for me.
\\[0.5cm]
I am furthermore grateful to Prof.\ Helmut B{\"o}lcskei, who very kindly agreed to act as a referee.
\\[0.5cm]
I would like to express my thanks to my parents. Their support made my studies of electrical engineering possible.
\\[0.5cm]
Special thanks to Dipl.-Ing.\ Sebastian Schmutzhard for sharing his impressive expertise on the mathematical theory of reproducing kernel Hilbert spaces. 
\\[0.5cm]
Finally and most importantly I would like to express my deep gratefulness for Julia, Helene and Emily for being with me. While
continuously supporting and motivating me throughout the writing of this thesis, they regularly reminded me about what is really important in life.


\begin{titlepage}
\thispagestyle{empty}
\vspace*{50mm}
\emph{To Julia, Helene and Emily}
\end{titlepage}


\thispagestyle{plain} 

\chapter*{Abstract}

The investigation of the effects of sparsity or sparsity constraints in signal processing problems has received considerable attention recently. 
Sparsity constraints refer to the a priori information that the object or signal of interest can be represented by using only few elements 
of a predefined dictionary. Within this thesis, sparsity refers to the fact that a vector to be estimated has only few nonzero entries.  

One specific field concerned with sparsity constraints has become popular under the name ``Compressed Sensing'' (CS). Within CS, the sparsity is exploited in 
order to perform (nearly) lossless compression. Moreover, this compression is carried out jointly or simultaneously with the process of sensing a 
physical quantity. 
 
In contrast to CS, one can alternatively use sparsity to enhance signal processing methods. 
Obviously, sparsity constraints can only improve the obtainable estimation performance since 
the constraints can be interpreted as an additional prior information about the unknown parameter vector which is to be estimated. 
Our main focus will be on this aspect of sparsity, i.e., we analyze 
how much we can gain in estimation performance due to the sparsity constraints.

We will study in detail two specific estimation problems that are already well investigated in the absence of sparsity constraints. 
First, we consider an estimation problem coined the ``sparse linear model'' (SLM). Here, the unknown parameter determines the mean 
of an observed Gaussian random vector. Second, we consider an estimation problem that we refer to as the ``sparse parametric covariance model'' (SPCM). 
In this case, the unknown parameter determines the covariance matrix of an observed Gaussian random vector. 

The main results of this thesis will be lower bounds on the variance of estimators for these two sparse estimation problems. 
We will demonstrate that the mathematical framework of reproducing kernel Hilbert spaces (RKHS) allows a simple derivation of lower bounds on the estimator variance 
and, moreover, a natural comparison of estimation problems with and without sparsity constraints.


\chapter*{Zusammenfassung}
Die Untersuchung der Auswirkungen von Sp{\"a}rlichkeitsbedingungen (SB) in verschiedenen Problemen der Signalverarbeitung ist ein wichtiges Ziel der aktuellen 
Forschung. Unter SB versteht man das Vorwissen, dass das Nutzsignal als Superposition weniger ``Elementarsignale'' dargestellt werden kann.
In dieser Arbeit beziehen sich die SB auf das Vorwissen, dass ein zu sch{\"a}tzender Parametervektor wenige (im Vergleich zu seiner L{\"a}nge) Eintr{\"a}ge 
aufweist, die von Null verschieden sind. Diese Vektoren werde dann als ``sp{\"a}rlich'' (sparse) bezeichnet. 

Solche SB spielen eine zentrale Rolle in der Theorie von ``Compressed Sensing'' (CS). Im Rahmen von CS werden die SB verwendet, um eine (nahezu) verlustfreie 
Datenkompression durchzuf{\"u}ren. Diese Datenkompression wird dabei mit dem Prozess des Messens einer physikalischen Gr{\"o}\ss e kombiniert. 

Im Gegensatz zu CS kann man die SB auch zur Versbesserung bestehender Signalverarbeitungsmethoden benutzen. 
In dieser Arbeit steht diese zweite M{\"o}glichkeit im Vordergrund. Ein wesentliches Ziel ist die Quantifizierung der durch SB bewirkten potenziellen Leistungsverbesserung. 

Es werden zwei verschiedene Sch{\"a}tzprobleme mit SB im Detail analysiert. Im ersten Problem, dem ``sparse linear model'' (SLM), bestimmt der unbekannte sp{\"a}rliche 
Parametervektor in linearer Weise den Mittelwert eines beobachteten Gau{\ss}schen Zufallsvektors. 
Im zweiten Problem, hier als ``sparse parametric covariance model'' (SPCM) bezeichnet, bestimmt der unbekannte sp{\"a}rliche 
Parametervektor die Kovarianzmatrix eines beobachteten Gau{\ss}schen Zufallsvektors. 

Die Hauptergebnisse dieser Arbeit sind untere Schranken f{\"u}r die Varianz von Sch{\"a}tzern f{\"u}r diese beiden Sch{\"a}tzprobleme. Es wird gezeigt, dass die Theorie der 
``reproducing kernel Hilbert spaces'' (RKHS) eine elegante Ableitung von unteren Varianzschranken erlaubt. 
Die gefundenen Schranken gestatten einen nat{\"u}rlichen Vergleich von Sch{\"a}tzproblemen mit und ohne SB.




\tableofcontents%
\newpage

\pagenumbering{arabic}%
\begin{titlepage}
$ $ 
\end{titlepage}
\setcounter{page}{1}




\chapter{Introduction}

\emph{Sparsity} or \emph{sparsity constraints} refer to the a-priori information that the parameter or signal vector of interest can be represented by using only few elements 
of a predefined dictionary. In particular, within this thesis, sparsity refers to the fact that a vector has only few non-zero entries. 
This thesis is concerned with the theory of classical (i.e., non-Bayesian) estimation problems where it is known that the unknown parameter vector is sparse. 

Given a general classical estimation problem,\footnote{We will make the notion of a classical estimation problem precise in Chapter \ref{chap_classic_est_fund}.} 
in particular one with sparsity constraints, two main questions are: 
\begin{itemize}
\item How to design accurate and efficient estimators? 
\item Given a specific estimator, how far is it away from being optimal? 
\end{itemize}
This work will mainly address the second question, i.e., the assessment of the quality of a given estimator. To that end we have to define an optimality criterion. In this thesis we  
will use the estimator variance (cf. Chapter \ref{chap_classic_est_fund}) as the primary performance criterion. If we have a lower bound on the 
set of all realizable values of the performance measure, we can assess the optimality of a given estimator by comparing its performance measure, i.e., its variance, with this lower bound. A central part of this 
work is concerned with the derivation of such lower bounds on the minimum achievable variance of estimators with a prescribed bias function. We will make these ideas precise in Chapter \ref{chap_classic_est_fund}. 

This thesis considers the application of the well-known theory of \emph{reproducing kernel Hilbert spaces} (RKHS) to specific problems within classical estimation theory. 
A RKHS is a special type of Hilbert space which has specific properties that will be discussed in detail in Chapter \ref{chap_RKHS_Fund}. The application of RKHS theory to 
statistical signal processing proved to be very successful, in particular in the context of machine learning \cite{Cucker02onthe,sun_jfaa}. 
While the application of RKHS to general classical estimation theory dates back to the late $1950$s \cite{Parzen59}, it seems that so far the RKHS framework has not been used for the 
study of estimation problems with sparsity constraints. It will be demonstrated in this thesis that RKHS theory is a well-suited tool for the analysis 
of the effect of sparsity constraints within classical estimation problems. 

Using the theory of RKHS, we will study in detail two specific estimation problems that have already been well investigated in the absence of sparsity constraints. First, 
we consider an estimation problem coined the ``sparse linear model'' (SLM),\footnote{Again we justify the unusual convention to name estimation problems after models by the fact, that an estimation problem is 
defined to a large extend by an observation model which connects the unknown parameter to be estimation with the observation.} in which the unknown parameter determines the mean 
of an observed Gaussian random vector. Second, we consider the estimation problem coined the ``sparse parametric covariance model'' (SPCM), 
in which the unknown parameter determines the covariance matrix of an observed Gaussian random vector. 
We will show how the SLM and the SPCM can be obtained as special cases of a larger class of sparse estimation problems where the statistic of the observation 
is given by an \emph{exponential family} (cf. Section \ref{sec_exp_family}). 
The main results of this thesis will be lower bounds on the variance of estimators for the sparse estimation problems that have a prescribed bias or mean. 

\section{Sparsity in Statistical Signal Processing}

Consider a vector $\mathbf{x} \in \mathbb{R}^{N}$ in the $N$-dimensional real Euclidean space \cite{HalmosFiniteDimVecSpace,RudinBookPrinciplesMatheAnalysis, RudinBook}. 
We then call the vector $\mathbf{x}$ \emph{sparse} if it has only few nonzero entries, i.e., $\| \mathbf{x} \|_{0} \ll N$ where $\| \mathbf{x} \|_{0}$ denotes the number of nonzero entries of $\mathbf{x}$. 
More precisely, the set of all vectors $\mathbf{x} \inÊ\mathbb{R}^{N}$ with \emph{at most} $S$ nonzero entries will be called the set of \emph{strictly} $S$\emph{-sparse vectors} and denoted by 
\begin{equation}
\mathcal{X}_{S} \triangleq \{ \mathbf{x} \in \mathbb{R}^{N} \big| \| \mathbf{x} \|_{0} \leq S \},
\end{equation}
where $S \in [N] \triangleq \{1,\ldots,N\}$, is called the \emph{sparsity degree}. 
Beside the notion of strict sparsity, there are also weaker notions of sparsity or sparse signals. E.g., the class of \emph{compressible} signals \cite{Cevher09ITW}
consists of vectors whose sorted magnitudes of its entires exhibit a power-law decay so that the signal vector can be well approximated by a strictly $S$-sparse vector with a suitable sparsity degree $S$. 
Alternatively, we call a vector $\mathbf{x} \in \mathbb{R}^{N}$ approximately $S$-sparse if it consists of $S$ ``large'' entries and $N-S$ ``small'' entries. However, for 
a rigorous analysis, one has to make the meaning of a ``large'' and a ``small'' entry precise. 
A specific approach to characterize the set of approximately $S$-sparse vectors is to measure the norm of the ``tail,'' i.e., the 
norm of the vector that is obtained from $\mathbf{x}$ by setting the $S$ largest (in magnitude) entries to zero. This quantification 
of approximate sparsity is used in \cite{Can06a} to characterize the performance of a signal processing algorithm. 
Another approach to model approximately sparse vectors will be discussed in detail in Section \ref{sec_strict_sparstiy_SLM}. 

If for a given application, e.g. radar \cite{CSRadar} or channel estimation \cite{GT_icassp08,Carbonelli06}, it is known that the object of interest can be represented by a sparse vector, then we should exploit this a-priori information. This can 
be done in two alternative ways. On the one hand, we can use sparsity to perform a nearly lossless compression. On the other hand, we could use the sparsity information in order
to improve the performance of existing signal processing schemes. 

The former approach is taken within \emph{Compressed Sensing} (CS) \cite{Don06,CandesCSTutorial,BaraniukCSTutorial}, i.e., CS refers to signal processing applications where the sparsity is exploited in order to perform a nearly 
lossless compression simultaneously with the sensing process. The output of the sensing process is then referred to as compressive measurements \cite{DavSigProcCompMeas}. 
From a more formal viewpoint, the theory of CS is concerned with two fundamental problems: 

\begin{itemize}
\item
The first problem considered in CS is that of finding good methods for the joint compression and sensing process. 
In particular, if the sensing process is linear and modeled by a matrix-vector multiplication between a CS measurement matrix $\mathbf{M} \in \mathbb{R}^{M \times N}$ and the sparse signal vector $\mathbf{x} \in \mathbb{R}^{N}$, 
we would like to construct measurement matrices which yield both a good compression performance (i.e., a small number of compressive measurements $z_{m}$, with $\mathbf{z} = \mathbf{M} \mathbf{x}$, which is equivalent to the requirement $M \ll N$) 
and an accurate reconstruction of the signal vector $\mathbf{x}$ from the compressive measurements $\mathbf{z} = \mathbf{M}\mathbf{x}$ at a later point. These two desiderata are contradicting goals and thus 
there is a tradeoff between compression performance and achievable reconstruction accuracy.  
So far, the best known performance guarantees are obtained by random constructions of the CS measurement matrix \cite{CandesCSTutorial,Don06}.
\item
The second main problem faced in CS is how to exploit the sparsity of the sparse signal vector $\mathbf{x} \in \mathbb{R}^{N}$ in order to estimate or recover it in an efficient and accurate way from given compressive measurements $\mathbf{z} = \mathbf{M} \mathbf{x}$.
Existing results show that the sparsity allows for computationally efficient and accurate estimation or recovery procedures \cite{JustRelax,GreedisGood,CoSAMP}.  
\end{itemize}

By contrast to CS, one may choose not to perform any compression but instead use the full available amount of raw data and exploit the sparsity constraints in order to obtain a better performance 
of signal processing algorithms (compared to the performance obtained in the absence of any sparsity constraints). Our thesis is concerned with this rationale, i.e., we analyze how much one can 
potentially gain in performance using the a priori information represented by the sparsity constraints. 

\section{Outline of the Thesis} 

In Chapter \ref{chap_classic_est_fund}, we present a brief review of some important concepts used within classical (non-Bayesian) estimation theory. 
We will discuss two major rationales within classical estimation theory, i.e., the minimax and the minimum variance rationales. 
The mathematical tools used in later chapters are presented in Chapter \ref{chap_RKHS_Fund}. In Chapter \ref{chap_RKHS_MVE}, we 
introduce the RKHS approach to minimum variance estimation and give detailed interpretations of some well-known lower bounds on the estimator variance  
from the (geometric) viewpoint of RKHS. 
The two core chapters of this thesis are Chapter \ref{chap_SLM} and Chapter \ref{chap_SCM}, where the RKHS approach to minimum variance estimation is specialized to 
two specific estimation problems, i.e., the SLM and the SPCM. In Chapter \ref{chap_SLM}, we give a detailed characterization of the RKHS associated to the SLM and use this characterization for the 
derivation of novel lower bounds on the estimator variance. We will also show in Chapter \ref{chap_SLM} that the (linear) CS recovery problem can be interpreted as a specific instance of the SLM. 
In Chapter \ref{chap_SCM}, we will introduce an RKHS approach to the SPCM, which is quite different from the RKHS approach to the SLM. 
This RKHS approach will then be used to derive novel lower bounds on the estimator variance for the SPCM. 

\section{Contribution and Related Work} 

This thesis develops the application of the theory of RKHS to classical estimation problems with sparsity constraints. 
The main part of the thesis is concerned with the analysis of two specific estimation problems. 

The first problem concerns the SLM. The goal is the estimation of a sparse parameter vector of which a noisy and linearly distorted version is observed. 
The SLM is obtained from the well-known linear Gaussian model (LGM) \cite{poorsspbook,kay,scharf91} by adding sparsity constraints on the unknown parameter vector. 
A special case of the SLM has been considered in \cite{DonohoJohnstone94, Donoho94idealspatial, Donoho98minimaxestimation}, where the authors investigate asymptotic minimax estimation. 
By contrast, we consider a finite-dimensional setting. Furthermore, we do not consider minimax estimation but a different estimation rationale, i.e., minimum variance estimation, where the bias or mean of an estimator is prescribed. 
We give results concerning the minimum variance that can be obtained by estimators with a prescribed bias. In this spirit, 
the authors of \cite{ZvikaCRB, ZvikaSSP} derive lower bounds on the variance of estimators for the SLM 
with a locally prescribed bias. By locally we mean that the bias has to be defined only at a given, fixed parameter value and in a small neighborhood around this point. 
Their bound is an adaption of the well-known \CRBfull (CRB) \cite{cramer45} to the sparse setting and depends only on a first-order characterization (related to a truncated 
Taylor expansion of an analytic function) of the estimation problem and bias. In a similar setting, the authors of \cite{BabAympAchievCRBSLM} prove the asymptotic achievability of the 
CRB for the CS recovery problem. Whereas \cite{BabAympAchievCRBSLM} considers the SLM with a
random system matrix, we consider a fixed deterministic system matrix as in \cite{ZvikaCRB}. 

By contrast to \cite{ZvikaCRB, ZvikaSSP}, we derive lower bounds on the estimator variance using a higher-order characterization of the estimation problem and estimator bias. 
This implies that our lower bounds will be in general tighter, i.e., higher,  than those found in \cite{ZvikaCRB, ZvikaSSP}. Another major difference is that the results of \cite{ZvikaCRB,ZvikaSSP} 
and, in a more general setting, those of \cite{GormanHero} are stated in terms of generalized (matrix) inequalities \cite{BoydConvexBook} involving the covariance matrix of an estimator. 
However, in this work, we present bounds on the variance of the individual entries of an estimator. In principle both type of bounds are equivalent, i.e., we can in general reformulate our 
bounds in terms of generalized inequalities and vice versa. In particular, the bounds in \cite{ZvikaCRB,ZvikaSSP,GormanHero} can be reformulated such that they can be compared with our bounds. 

The second estimation problem considered in this thesis concerns the SPCM. The goal is to estimate a sparse parameter vector that determines the covariance matrix $\mathbf{C}$ of a Gaussian (random) signal vector of which a noisy version is observed. While covariance estimation is an established topic of research (see, e.g., the recent work \cite{CovEstDecompGraphModel, RueBue09, ShrinkageMMSECovEst}), sparsity constraints 
are mostly placed directly on the inverse covariance matrix, i.e., the \emph{precision matrix} $\mathbf{C}^{-1}$, by assuming that most entries of $\mathbf{C}^{-1}$ are zero. This is mainly because the precision matrix of a Gaussian random vector is directly related to a graphical model \cite{GraphModExpFamVarInfWainJor} describing the statistical structure of the random vector. A sparse precision matrix $\mathbf{C}^{-1}$ thus corresponds to a sparse graph, i.e., a graph with relatively few edges compared to the number of nodes. 

From a methodological perspective, this thesis specializes the theory of RKHS for the application to minimum variance estimation with sparsity constraints. By sparsity constraint, we mean that the unknown parameter vector $\mathbf{x} \in \mathbb{R}^{N}$ 
to be estimated is known to be $S$-sparse, i.e., $\mathbf{x} \in \mathcal{X}_{S}$ with a known sparsity degree $S \in \mathbb{N}$. The application of Hilbert space, in particular RKHS, techniques to minimum variance 
estimation dates back to the seminal work of Parzen \cite{Parzen59}, which introduced the theory of RKHS to general estimation problems. The approach of \cite{Parzen59} is further developed in \cite{Duttweiler73b} for specific 
estimation problems. We see our thesis as another specialization of the work \cite{Parzen59} to estimation problems with sparsity constraints. 

Differently from this ``RKHS line'' of work, a Hilbert space approach has been proposed for general minimum variance estimation problems in the recent paper \cite{TodrosLowerBounds2010}. The Hilbert space used in \cite{TodrosLowerBounds2010} for the derivation of lower bounds on the estimator variance is different from the RKHS used in this thesis and based on integral operators associated to a minimum variance estimation problem. Furthermore, the authors in \cite{TodrosLowerBounds2010} do not consider estimation with sparsity constraints. 

The main contributions of this thesis (listed according to their order of appearance) are the following.
\begin{itemize} 

\item Chapter \ref{chap_RKHS_MVE}: We show in Section \ref{sec_classes_min_var_problems_con_var_kernel} that under mild conditions on the underlying estimation problem and prescribed, 
the minimum achievable variance $L_{\mathcal{M}}$ of the minimum variance problem $\mathcal{M} = \minvarproblem$, 
viewed as a function of the parameter vector $\mathbf{x}_{0}$, is lower semi-continuous. 

\item Chapter \ref{chap_RKHS_MVE}: We prove in Section \ref{sec_suff_stat_RKHS}, that the RKHS associated to a minimum variance problem remains unchanged if the original observation is replaced by any sufficient statistic. This result 
could be regarded as an RKHS analogue of the famous ``Rao-Blackwell-Lehmann-Scheff{\'e}'' (RBLS) theorem \cite{kay}.\footnote{The RBLS theorem basically states that given an arbitrary estimator one can always 
find an estimator that (i) depends on the observation only via a sufficient statistic, (ii) has the same bias as the original estimator and (iii) its variance does not exceed those of the original estimator. Thus, using a sufficient statistic as the new observation means 
no loss of information in terms of minimum variance estimation.}

\item Chapter \ref{chap_RKHS_MVE}: We provide a detailed derivation of some well-known variance bounds using exclusively RKHS theory. This yields an intuitive geometric interpretation of these bounds. 

\item Section \ref{sec_RKHS_SLM}: We present a derivation and complete characterization of the RKHS associated with the SLM. This includes as a special case a novel characterization of the RKHS associated with the LGM. Our characterization will 
be quite different from those given \cite{Duttweiler73b} which already considered the LGM and associated RKHS. 

\item Section \ref{sec_RKHS_basic_facts_SLM}: We discuss minimum variance estimation for the SLM. We characterize the class of valid bias functions, i.e., 
those bias functions for which there exists at least one finite-variance estimator having this bias function. In particular, we show that there does not exist a finite-variance unbiased estimator of the support set of the unknown sparse parameter vector. 

\item Section \ref{sec_RKHS_Lower_Bounds_SLM}: We interpret existing lower bounds on the estimator variance for the SLM from the RKHS viewpoint and then derive novel lower bounds on the estimator variance based on the RKHS derived in Section \ref{sec_RKHS_SLM}. 

\item Section \ref{sec_SSNM}: We focus on a special case of the SLM termed the \emph{sparse signal in noise model} (SSNM). By exploiting the specific structure of the SSNM, it is possible to derive stronger 
results than for the general SLM. In particular, we obtain closed-form expressions for the minimum achievable variance (i.e., the Barankin bound) 
and the corresponding estimator which achieves this minimum variance (i.e., the locally minimum variance (LMV) estimator, which will be defined in Section \ref{sec_min_var_est}) for a wide class of bias functions, which 
includes the bias function of the hard-thresholding (HT) estimator \cite{MallatBook}.  

\item  Section \ref{sec_strict_sparstiy_SLM}: We show that the strict sparsity constraints of the SLM are necessary in order for the minimum achievable 
variance of the estimation problem to be strictly lower than that of the minimum achievable variance for the LGM.

\item Section \ref{sec_slm_viewpoint_CS}: Based on the RKHS approach, we present an SLM-based analysis of the CS recovery problem. 

\item Section \ref{sec_numerical_SLM}: We compare the actual variance behaviors of some popular estimation schemes for the SLM and SSNM with the corresponding lower bounds. 

\item Section \ref{sec_intro_SPCM} and \ref{sec_D_restricted_SPCM}: We present a specific RKHS approach to the SPCM and discuss the differences between the RKHS approach to the SPCM and SLM. 

\item Section \ref{sec_lower_bounds_SPCM}: We derive novel lower bounds on the estimator variance for the SPCM based on the RKHS defined in Section \ref{sec_D_restricted_SPCM}. 

\item Section \ref{sec_SDPCM}: We focus on a special case of the SPCM termed the \emph{sparse diagonalizable parametric covariance model} (SDPCM). By exploiting the specific structure of the SDPCM, it is possible to derive stronger 
results than for the general SPCM. In particular, we will show in Section \ref{sec_unbiasecd_SDPCM} that for unbiased estimation for the SDPCM, 
the strict sparsity constraints are necessary in order for the minimum achievable variance to be strictly lower 
than for the case where no sparsity constraints are present. 

\item Section \ref{sec_comp_bounds_actual_var_SPCM}: We compare our lower bounds on the variance with the actual variance behavior of two specific estimation schemes for the SDPCM. 


\end{itemize}

\section{Notation} 

In the following we will define the notation used in this thesis. The notation related to probability theory will be introduced separately in Chapter \ref{chap_classic_est_fund}. 

The set of (nonnegative) real numbers is denoted by ($\mathbb{R}_{+}$) $\mathbb{R}$, the set of natural numbers is denoted by $\mathbb{N} = \{1,2,3, \ldotsÊ\}$ 
and the set of nonnegative integers by $\mathbb{Z}_{+} = \{0,1,2,\ldots\}$. 
Given a nonnegative integer $L \in \mathbb{Z}_{+}$, we denote by $[L]$ the set of all numbers up to $L$, i.e., $[L] \triangleq \{1,\ldots,L\}$ if $L>0$ and for the special case $L=0$ we set $[0] \triangleq \emptyset$ where $\emptyset$ denotes the empty set. 
Given an arbitrary set $\mathcal{A}$, we define the Kronecker delta $\delta_{l,l'} \in \{0,1\}$ for $l,l' \in \mathcal{A}$ by setting $\delta_{l,l'} = 0$ if and only if $l \neq l'$. Given the two sets $\mathcal{A}$ and $\mathcal{B}$, we denote 
by $\mathcal{A} \setminus \mathcal{B}$ the set which consists of those elements of $\mathcal{A}$ that are not an element of $\mathcal{B}$ and by $\mathcal{A} \cup \mathcal{B}$ and $\mathcal{A} \cap \mathcal{B}$ the union and intersection, 
respectively of $\mathcal{A}$ and $\mathcal{B}$ \cite{RudinBook}. 
Furthermore we denote by $|\mathcal{A}|$ the cardinality, i.e., the number of elements of a finite set $\mathcal{A}$. 
We denote by $\log x$ the natural logarithm of the positive real number $x$. 

Uppercase boldface letters denote matrices and lowercase boldface letters denote vectors.\footnote{Since a matrix can be viewed as a special case of 
a linear operator, we will use uppercase boldface letters also to denote general linear operators.} 
The superscript $^{T}$ stands for transposition. Given a matrix $\mathbf{H} \in \mathbb{R}^{M \times N}$ we denote 
by $\vectorize{\mathbf{H}} \in \mathbb{R}^{MN}$ the certain vector that is obtained by vertically stacking the columns of $\mathbf{H}$. 
For a square matrix $\mathbf{H} \in \mathbb{R}^{N \times N}$, we denote by 
$\tracem{ \mathbf{H}}$, $\detm{\mathbf{H}}$ and $\mathbf{H}^{-1}$ its trace, determinant and inverse (if it exists), respectively.
The kernel or null space $\mathcal{N}(\mathbf{H})$ of a matrix $\mathbf{H} \in \mathbb{R}^{M \times N}$ is defined as 
$\mathcal{N}(\mathbf{H}) \triangleq \{ \mathbf{x} \in \mathbb{R}^{N} | \mathbf{H} \mathbf{x} = \mathbf{0} \}$. We denote by $\linspan(\mathbf{H})$ the column span of a matrix $\mathbf{H} \in \mathbb{R}^{M \times N}$, i.e. 
$\linspan(\mathbf{H}) \triangleq \{ \mathbf{y} \in \mathbb{R}^{M} | \exists \mathbf{x} \in \mathbb{R}^{N}: \mathbf{y} = \mathbf{H} \mathbf{x} \}$. 
 
The $k$th entry of a vector $\mathbf{x} \in \mathbb{R}^{N}$ and the entry in the $k$th row and $l$th column of a matrix $\mathbf{H} \in \mathbb{R}^{M \times N}$ are denoted 
by $x_{k}$ and $(\mathbf{H})_{k,l}$, respectively.
We designate the identity matrix of dimension $N \times N$ as $\mathbf{I}_{N}$ or just as $\mathbf{I}$ if the dimension is clear from the context. The $k$th column of the identity matrix will 
be denoted by $\mathbf{e}_{k}$. 
By $\mathbf{0}$ we denote a vector or matrix consisting only of zeros whereby the dimension of the vector or matrix should by clear from context. 
Given two vectors $\mathbf{a}, \mathbf{b} \in \mathbb{R}^{N}$ and two matrices $\mathbf{A}, \mathbf{B} \in \mathbb{R}^{M \times N}$ we denote by $\mathbf{a} \neq \mathbf{b}$ and $\mathbf{A} \neq \mathbf{B}$ the fact 
that the two vectors and matrices, respectively differ in at least one entry. 
The support, i.e. the set of indices that correspond to non-zero entries, and number of non-zero entries of a vector $\mathbf{x} \in \mathbb{R}^{N}$ are denoted by $\supp(\mathbf{x})$ and $\| \mathbf{x} \|_{0}$, respectively.
Given a vector $\mathbf{x} \in \mathbb{R}^{N}$, we define its $p$ norm as $\| \mathbf{x} \|_{p} \triangleq \left( \sum_{k \in [N]} x_{k}^{p}\right)^{1/p}$. 
We denote by  $\| \mathbf{a} \|_{\infty} \triangleq \max_{k \in [N]} |a_{k}|$ the largest magnitude of the entries of the vector $\mathbf{a} \in \mathbb{R}^{N}$.
Given a square matrix $\mathbf{H} \in \mathbb{R}^{M \times M}$ we denote by $\| \mathbf{H} \|_{p}$ its matrix $p$-norm defined as $\| \mathbf{H} \|_{p} \triangleq \sup_{\| \mathbf{x} \|_{p}=1} \| \mathbf{H} \mathbf{x} \|_{p}$, 
where $\mathbf{x} \in \mathbb{R}^{M}$. 

Given an index set $\mathcal{K} \subseteq [N]$, a vector $\mathbf{x} \in \mathbb{R}^{N}$, and a matrix $\mathbf{H} \in \mathbb{R}^{M \times N}$, we denote by $\mathbf{x}^{\mathcal{K}} \in \mathbb{R}^{N}$ and 
$\mathbf{H}_{\mathcal{K}} \in \mathbb{R}^{M \times |\mathcal{K}|}$ the vector which is obtained from $\mathbf{x}$ by zeroing all entries except those indexed by $\mathcal{K}$ and the matrix that is obtained by selecting only the columns of $\mathbf{H}$ indexed by $\mathcal{K}$, respectively. 

The fact that a square matrix $\mathbf{R} \in \mathbb{R}^{N \times N}$ is positive semi-definite (psd), i.e., 
\begin{equation} 
\label{equ_def_psd_matrix}
\mathbf{x}^{T} \mathbf{R} \mathbf{x} \geq 0 \quad \forall \mathbf{x} \in \mathbb{R}^{N},
\end{equation} 
is denoted by $\mathbf{R} \geq \mathbf{0}$. Given two psd matrices $\mathbf{A}, \mathbf{B} \in \mathbb{R}^{N \times N}$, we write $\mathbf{A} \geq \mathbf{B}$ if the matrix $\mathbf{A} - \mathbf{B}$, given 
elementwise as the difference between the elements of $\mathbf{A}$ and $\mathbf{B}$, is psd, i.e., $\mathbf{A} - \mathbf{B} \geq \mathbf{0}$ \cite{BoydConvexBook}. 
If the inequality in \eqref{equ_def_psd_matrix} is replaced by a strict inequality, we call the matrix $\mathbf{R}$ positive definite, in symbols $\mathbf{R} > \mathbf{0}$. 

Given an arbitrary matrix $\mathbf{H} \in \mathbb{R}^{M \times N}$, we call any decomposition of the form 
\begin{equation} 
\label{equ_thin_SVD}
\mathbf{H} = \mathbf{U} \mathbf{\Sigma} \mathbf{V}^{T} 
\end{equation} 
with $\mathbf{U} \in \mathbb{R}^{M \times r}, \mathbf{V} \in \mathbb{R}^{N \times r}, \mathbf{\Sigma} \in \mathbb{R}^{r \times r}$ the \emph{thin singular value decomposition} (SVD) \cite{golub96} of 
$\mathbf{H}$ if $\mathbf{U}^{T}\mathbf{U}=\mathbf{I}$, $\mathbf{V}^{T} \mathbf{V} = \mathbf{I}$ and the matrix $\mathbf{\Sigma}$ is a diagonal matrix 
with strictly positive real-valued main diagonal elements, i.e., $(\mathbf{\Sigma})_{k,k} > 0$. These main diagonal entries are called the \emph{singular values} of $\mathbf{H}$. 
Note that it can be shown that every nonzero matrix $\mathbf{H} \in \mathbb{R}^{M \times N}$ has a thin SVD \cite{golub96}. The rank or column rank of $\mathbf{H}$, denoted by $\rank(\mathbf{H})$, is defined as the 
unique dimension $r$ of the matrix $\mathbf{\Sigma}$ appearing in \eqref{equ_thin_SVD}. We also set $\rank(\mathbf{0}) = 0$. 

Given a psd matrix $\mathbf{H} \in \mathbb{R}^{N \times N}$, i.e., $\mathbf{H} \geq \mathbf{0}$, we call any decomposition of the form 
\begin{equation} 
\label{equ_EVD}
\mathbf{H} = \mathbf{U} \mathbf{\Sigma} \mathbf{U}^{T} 
\end{equation} 
with $\mathbf{U} \in \mathbb{R}^{N \times \rank(\mathbf{H})}, \mathbf{\Sigma} \in \mathbb{R}^{\rank(\mathbf{H}) \times \rank(\mathbf{H})}$ the \emph{thin eigenvalue decomposition} 
(EVD) (or \emph{thin symmetric Schur decomposition}) \cite{golub96} of $\mathbf{H}$ 
if $\mathbf{U}^{T}\mathbf{U}=\mathbf{I}$ and the matrix $\mathbf{\Sigma}$ is a diagonal matrix with positive real-valued main diagonal elements, i.e., $(\mathbf{\Sigma})_{k,k} > 0$. 
These main diagonal elements are called the \emph{eigenvalues} of $\mathbf{H}$. Note that it can be shown that every psd matrix $\mathbf{H} \in \mathbb{R}^{N \times N}$, 
where $\mathbf{H} \neq \mathbf{0}$, has a thin EVD \cite{golub96}.

Given a nonzero matrix $\mathbf{H} \in \mathbb{R}^{M \times N}$ and its thin SVD \eqref{equ_thin_SVD}, we denote by $\mathbf{H}^{\dagger} \in \mathbb{R}^{N
\times M}$ its Moore Penrose pseudo inverse \cite{golub96} defined by $\mathbf{H}^{\dagger} \triangleq \mathbf{V} \mathbf{\Sigma}^{-1} \mathbf{U}^{T}$. Note that 
due to the properties of the thin SVD, the inverse $\mathbf{\Sigma}^{-1}$ is guaranteed to exist. For $\mathbf{H} = \mathbf{0}$, the Moore Penrose inverse is defined as $\mathbf{H}^{\dagger} = \mathbf{0}$. 

A sequence of elements $f_{l} \in \mathcal{M}$, where $l \in \mathbb{N}$, that belong to a certain set $\mathcal{M}$ is denoted by $\{ f_{l} \in \mathcal{M} \}_{l \rightarrow \infty}$ or, if the underlying set 
$\mathcal{M}$ is clear from the context, by $\{ f_{l} \}_{l \rightarrow \infty}$. 

We call a map or function $f(\cdot): \mathcal{A} \rightarrow \mathcal{B}$ invertible if there exists a map $f^{-1}(\cdot)$, called the inverse map or inverse function, such that $f^{-1}(f(\mathbf{x})) = \mathbf{x}$ for every $\mathbf{x} \in \mathcal{A}$. 
A map $f(\cdot): \mathcal{A} \rightarrow \mathcal{B}$ is called bijective if $f(\mathcal{A}) = \mathcal{B}$ and if it is invertible, i.e., there exists an inverse map $f^{-1}(\cdot)$ \cite{RudinBookPrinciplesMatheAnalysis}. 
If a real-valued function $f(\cdot): \mathcal{D} \rightarrow \mathbb{R}$ is identically zero on its domain, i.e., $f(\mathcal{D}) = \{ 0 \}$, we will denote this by $f(\cdot) \equiv 0$. 
Given a function $f(\cdot): \mathcal{D} \rightarrow \mathbb{R}$ defined on some domain $\mathcal{D}$ and a subdomain $\mathcal{D}_{1} \subseteq \mathcal{D}$, we denote by 
$f\big|_{\mathcal{D}_{1}}(\cdot)$ or $f(\cdot)\big|_{\mathcal{D}_{1}}$ its restriction to $\mathcal{D}_{1}$, i.e., the specific function defined on $\mathcal{D}_{1}$ that coincides with $f(\cdot)$ on $\mathcal{D}_{1}$. More 
explicitly, the function $f\big|_{\mathcal{D}_{1}}(\cdot): \mathcal{D}_{1} \rightarrow \mathbb{R}$ satisfies $f\big|_{\mathcal{D}_{1}}(\mathbf{x}) = f(\mathbf{x})$ for every $\mathbf{x} \in \mathcal{D}_{1}$. 
For the special case when the subdomain consists of a single element, i.e., $\mathcal{D}_{1}=\{ \mathbf{x} \}$, we denote 
by $f(\cdot)\big|_{\mathbf{x}}$ the function value $f(\mathbf{x})$ at $\mathbf{x}$. 

The closed ball centered at $\mathbf{x}_{0} \in \mathbb{R}^{N}$ with radius $r >0$ is denoted by $\mathcal{B}(\mathbf{x}_{0},r) \triangleq \{ \mathbf{x} \in \mathbb{R}^{N} \big| \| \mathbf{x} - \mathbf{x}_{0} \|_{2} \leq r \}$.

Given a real-valued function $f(\cdot): \mathcal{D}\subseteq \mathbb{R}^{N} \rightarrow \mathbb{R}$, we say that $f(\cdot)$ is continuous at the point $\mathbf{x}_{0} \in \mathcal{D}$ if $\lim\limits_{\mathbf{x} \rightarrow \mathbf{x}_{0}} f(\mathbf{x}) = f(\mathbf{x}_{0})$. A function is said to be continuous if it is continuous at every point of its domain. 
Furthermore, we say that a function $f(\cdot): \mathcal{D}\subseteq \mathbb{R}^{N} \rightarrow \mathbb{R}$ is lower (upper) semi-continuous at the point $\mathbf{x}_{0}\in \mathcal{D}$ if for every $\varepsilon >0$ there exists a radius $r$ such that 
$f(\mathbf{x}) \geq f(\mathbf{x}_{0}) - \varepsilon$ ($f(\mathbf{x}) \leq f(\mathbf{x}_{0}) + \varepsilon$) for every $\mathbf{x} \in \mathcal{B}(\mathbf{x}_{0}, r)$ and denote this fact by 
$\liminf\limits_{\mathbf{x} \rightarrow \mathbf{x}_{0}} f(\mathbf{x}) \geq f(\mathbf{x}_{0})$ ($\limsup\limits_{\mathbf{x} \rightarrow \mathbf{x}_{0}} f(\mathbf{x}) \leq f(\mathbf{x}_{0})$) \cite{CounterExamplesAnalysis,RudinBookPrinciplesMatheAnalysis}, where 
\begin{equation} 
\liminf\limits_{\mathbf{x} \rightarrow \mathbf{x}_{0}} f(\mathbf{x}) \triangleq \sup_{r >0} \bigg( \inf_{\mathbf{x} \in \mathcal{D} \cap (\mathcal{B}(\mathbf{x}_{0},r) \setminus \{ \mathbf{x}_{0}\})} f(\mathbf{x}) \bigg)\mbox{, }
\quad 
\limsup\limits_{\mathbf{x} \rightarrow \mathbf{x}_{0}} f(\mathbf{x}) \triangleq \inf_{r >0} \bigg( \sup_{\mathbf{x} \in \mathcal{D} \cap (\mathcal{B}(\mathbf{x}_{0},r) \setminus \{ \mathbf{x}_{0}\})} f(\mathbf{x}) \bigg). 
\end{equation} 

Given a real-valued function $R(\cdot,\cdot): \mathcal{D} \times \mathcal{D} \rightarrow \mathbb{R}$ and 
an element $\mathbf{x} \in \mathcal{D}$, we denote by $R(\cdot, \mathbf{x})$ the specific function $f(\cdot): \mathcal{D} \rightarrow \mathbb{R}$ given by $f(\mathbf{x}') = R(\mathbf{x}', \mathbf{x})$ for every $\mathbf{x}' \in \mathcal{D}$. 
Given two functions $f(\cdot): \mathcal{A} \rightarrow \mathcal{B}$, $g(\cdot): \mathcal{C} \rightarrow \mathcal{D}$, we denote by $f(\cdot)=g(\cdot)$ the fact that $\mathcal{A} = \mathcal{C}$, $\mathcal{B} = \mathcal{D}$ and 
for every $\mathbf{x} \inÊ\mathcal{A}$ we have $f(\mathbf{x}) = g(\mathbf{x})$. 

Given a real-valued function $f(\cdot): \mathcal{D} \rightarrow \mathbb{R}$ where $\mathcal{D} \subseteq \mathbb{R}^{N}$, we denote by $\frac{\partial^{\mathbf{e}_{k}} f(\mathbf{x})}{\partial \mathbf{x}^{\mathbf{e}_{k}}}$ or 
$\frac{\partial f(\mathbf{x})}{\partial x_{k} }$ the partial derivative (if it exists) of $f(\cdot)$ with respect to the $k$th entry, i.e., $\frac{\partial^{\mathbf{e}_{k}} f(\mathbf{x})}{\partial \mathbf{x}^{\mathbf{e}_{k}}}$ denotes the function 
$g(\mathbf{x}) = \lim_{h \rightarrow 0} \frac{f(\mathbf{x}+h \mathbf{e}_{k}) - f(\mathbf{x})}{h}$ \cite{RudinBookPrinciplesMatheAnalysis}. 
This definition is then recursively extended to an arbitrary multi-index $\mathbf{p} \in \mathbb{Z}_{+}^{N}$ \cite{KranzPrimerAnalytic} as follows: 
Denoting the partial derivative of order $\mathbf{p}$ by $h(\mathbf{x})$, i.e., $ h(\mathbf{x}) = \frac{\partial^{\mathbf{p}} f(\mathbf{x})}{\partial \mathbf{x}^{\mathbf{p}}}$, the 
partial derivative of order $\mathbf{p}' = \mathbf{p}+ \mathbf{e}_{k}$ denoted by $\frac{\partial^{\mathbf{p}'} f(\mathbf{x})}{\partial \mathbf{x}^{\mathbf{p}'}}$ is defined as $\frac{\partial^{\mathbf{e}_{k}} h(\mathbf{x})}{\partial \mathbf{x}^{\mathbf{e}_{k}}}$. 
Strictly speaking, this definition is not consistent when used for an arbitrary function $f(\cdot)$, since the partial derivative of order $\mathbf{p}$ might depend on the order of adding the $\mathbf{e}_{k}$. However, within this thesis we exclusively 
deal with functions $f(\cdot)$ for which this definition of the partial derivative of order $\mathbf{p} \in \mathbb{Z}_{+}^{N}$ is unambiguous. 
Given a real-valued function $g(\cdot): \mathcal{D} \rightarrow \mathbb{R}$ where $\mathcal{D} \subseteq \mathbb{R}^{N}$, we denote by $\frac{\partial g(\mathbf{x})} {\partial \mathbf{x}}\big|_{\mathbf{x}_{0}}$ 
the vector $\mathbf{a} \in \mathbb{R}^{N}$ given elementwise by $a_{k} = \frac{\partial g(\mathbf{x})} {\partial x_{k}}\big|_{\mathbf{x}_{0}}$.
Given a real-valued function $R(\cdot,\cdot): \mathcal{D} \times \mathcal{D} \rightarrow \mathbb{R}$ where $\mathcal{D} \subseteq \mathbb{R}^{N}$, we denote by 
$\frac{\partial^{\mathbf{e}_{k}} R(\mathbf{x}_{1},\mathbf{x}_{2})}{\partial \mathbf{x}_{2}^{\mathbf{e}_{k}}}$ or $\frac{\partial R(\mathbf{x}_{1},\mathbf{x}_{2})}{\partial x_{2,k}}$ 
its partial derivative of order $\mathbf{e}_{k}$ w.r.t.\ $\mathbf{x}_{2}$, i.e., the function $g(\cdot,\cdot) : \mathcal{D} \times \mathcal{D} \rightarrow \mathbb{R}$ 
defined as $g(\mathbf{x}_{1},\mathbf{x}_{2}) =   \lim_{h \rightarrow 0} \frac{R(\mathbf{x}_{1},\mathbf{x}_{2}+h \mathbf{e}_{k}) - R(\mathbf{x}_{1}, \mathbf{x}_{2})}{h}$. 
The generalization to partial derivatives of the form $\frac{\partial^{\mathbf{p}_{1}} \partial^{\mathbf{p}_{2}} R(\mathbf{x}_{1},\mathbf{x}_{2})}{\partial \mathbf{x}_{1}^{\mathbf{p}_{1}}\mathbf{x}_{2}^{\mathbf{p}_{2}}}$
with $\mathbf{p}_{1}, \mathbf{p}_{2} \in \mathbb{Z}_{+}^{N}$ is a straightforward extension of the definition of $\frac{\partial^{\mathbf{e}_{k}} R(\mathbf{x}_{1},\mathbf{x}_{2})}{\partial \mathbf{x}_{2}^{\mathbf{e}_{k}}}$ obtained 
by considering the pair $(\mathbf{x}_{1},\mathbf{x}_{2})$ as a ``super-vector'' of length $2N$.



For a finite set $\mathcal{T}$, the sum $\sum_{l \in \mathcal{T}} f[l]$ and product $\prod_{\l \in \mathcal{T}} f[l]$ of a real valued function $f[\cdot]: \mathcal{T} \rightarrow \mathbb{R}$ are defined 
in the usual sense. For the special case given by the empty set, i.e., $\mathcal{T} = \emptyset$ we set $\sum_{l \in \mathcal{T}} f[l] = 0$ and $\prod_{\l \in \mathcal{T}} f[l]=1$ respectively. 
Given an arbitrary set $\mathcal{T}$, we define the sum $\sum_{l \in \mathcal{T}} f[l]$ as the supremum 
of all finite sums $\sum_{l \in \mathcal{T}' \subseteq \mathcal{T}} f[l]$ where $\mathcal{T}'$ is an arbitrary finite subset of $\mathcal{T}$ \cite[p. 83]{RudinBook}. 
As noted in \cite{RudinBook}, the so defined sum is nothing but the Lesbegue integral with respect to the counting measure on $\mathcal{T}$. 

We denote by $\ell^{2}(\mathcal{T})$ the set of all real-valued functions $f(\cdot): \mathcal{T} \rightarrow \mathbb{R}$ defined on $\mathcal{T}$ for which 
the sum $\sum_{l \in \mathcal{T}}| f[l]|^{2}$ is finite. It can be shown \cite{RudinBook} that the set $\ell^{2}(\mathcal{T})$ is a normed vector space with norm $\| f[\cdot]\|_{\mathcal{T}} \triangleq  \sqrt{\sum_{l \in \mathcal{T}}| f[l]|^{2}}$.

Given a nonnegative integer $p \in \mathbb{Z}_{+}$ we denote by $p!$ the factorial, defined recursively by $0!  \triangleq 1$ and $p! \triangleq  p(p-1)!$. 
Given a multi-index $\mathbf{p} \in \mathbb{Z}_{+}^{N}$ \cite{KranzPrimerAnalytic}, i.e., an $N$-tuple of nonnegative integers, we denote by $\mathbf{p}!$ and $| \mathbf{p}|$ the product $\mathbf{p}! \triangleq \prod_{l \in [N]} p_{l}!$ 
and the sum $| \mathbf{p}|\triangleq \sum_{l \in [N]} p_{l}$, respectively. For a vector $\mathbf{x} \in \mathbb{R}^{N}$ 
and a multi-index $\mathbf{p} \in \mathbb{Z}_{+}^{N}$, we denote by $\mathbf{x}^{\mathbf{p}}$ the product $\mathbf{x}^{\mathbf{p}} \triangleq \prod_{l \in [N]} (x_{l})^{p_l}$. 
Given two multi-indices $\mathbf{p}_{1}, \mathbf{p}_{2} \in \mathbb{Z}_{+}^{N}$, we mean by $\mathbf{p}_{1} \leq \mathbf{p}_{2}$ that this inequality holds separately for every single entry, i.e., $p_{1,l} \leq p_{2,l}$ for every $l \in [N]$. 


\chapter{Elements of Classical Estimation Theory}
\label{chap_classic_est_fund}


Classical or non-Bayesian estimation theory \cite{kay,scharf91,LC,poorsspbook} is concerned with the problem of inferring the value of an unknown but deterministic parameter $\mathbf{x} \in \mathbb{R}^{N}$ based on the observation 
of a random vector $\mathbf{y} \in \mathbb{R}^{M}$.\footnote{Of course $\mathbf{x}$ and $\mathbf{y}$ do not necessarily have to be real-valued vectors, but we will assume this without loss of generality for our purposes.} By contrast 
to Bayesian estimation theory \cite{kay,scharf91,LC,poorsspbook}, we do not model the parameter $\mathbf{x}$ as a random vector that has a certain statistical characterization. 
Instead, it is assumed that the parameter vector has a fixed value that is however unknown. 

\section{Basic Concepts}
\label{sec_basic_concepts} 

The elements of a classical estimation problem are the following.
\begin{itemize}

\item {\bf Parameter set.} We assume that the parameter $\mathbf{x}$ is an element of the parameter set $\mathcal{X}$, i.e., 
\begin{equation} 
\label{equ_parameter_set}
\mathbf{x} \in \mathcal{X}.
\end{equation} 
Note that the knowledge of \eqref{equ_parameter_set} also expresses some prior information. Intuitively speaking, the smaller the parameter set $\mathcal{X}$, the more prior information is available. In the extreme case (which is practically irrelevant) 
of a singleton $\mathcal{X} = \{ \mathbf{x}_{0} \}$, i.e., the parameter set consists of a single element, the estimation problem becomes obviously trivial. A main part of this thesis is concerned with the question of how much certain reductions of the parameter set 
$\mathcal{X}$ help in terms of estimation quality.

\item {\bf Statistical model.} In order to infer the value of the parameter $\mathbf{x}$ from an observed random quantity $\mathbf{y}$, there has to be some relationship between 
$\mathbf{x}$ and $\mathbf{y}$. Within classical estimation theory, this relationship is modeled by a family (a set) of probability density functions (pdf) denoted by $f(\mathbf{y}; \mathbf{x})$ and called the statistical model. 
If the true parameter vector is $\mathbf{x}_{0}$, then the probability of the event $\{ \mathbf{y} \in \mathcal{A} \}$, where $\mathcal{A} \subseteq \mathbb{R}^{M}$ 
is an arbitrary measurable\footnote{We consider the Lebesgue measure for the Euclidean space $\mathbb{R}^{N}$.} set, is given by 
\begin{equation} 
\label{equ_int_prob_event}
\mathsf{P} \{ \mathbf{y} \in \mathcal{A} \} = \int_{\mathbf{y} \in \mathcal{A}} f(\mathbf{y}; \mathbf{x}_{0}) d \mathbf{y}. 
\end{equation} 
We will consider only statistical models for which one can define a pdf $f(\mathbf{y}; \mathbf{x})$ w.r.t.\ to the 
Lebesgue measure in $\mathbb{R}^{M}$ such that the relation \eqref{equ_int_prob_event} can be interpreted in the 
sense of integration w.r.t.\ the Lebesgue measure \cite{RudinBook}. However, the main concepts developed in this thesis could be in principle also applied to a more general case where 
the statistics of the observation $\mathbf{y}$ are described by a \emph{probability measure} \cite{AshProbMeasure,BillingsleyProbMeasure}.\footnote{A pdf is nothing but a convenient representation of probability measures that satisfy some technical conditions \cite{AshProbMeasure,BillingsleyProbMeasure}.}
We will use the symbol $\mathsf{E}_{\mathbf{x}_{0}} \{ \mathbf{T}(\mathbf{y}) \}$, where $\mathbf{T}(\mathbf{y}) \in \mathbb{R}^{K}$ is a (possibly random) function of the observation, for 
the expectation (cf. \cite{AshProbMeasure,BillingsleyProbMeasure}) of $\mathbf{T}(\mathbf{y})$ w.r.t.\ the pdf $f(\mathbf{y}; \mathbf{x}_{0})$, i.e.,
\begin{equation}
\label{equ_def_expec_random_func_observation} 
\mathsf{E}_{\mathbf{x}_{0}} \{ \mathbf{T}(\mathbf{y}) \} \triangleq  \int_{\mathbf{z},\mathbf{y}} \mathbf{z} f(\mathbf{z}=\mathbf{T}(\mathbf{y}) \big| \mathbf{y}) f(\mathbf{y}; \mathbf{x}_{0}) d \mathbf{z} \, d \mathbf{y}. 
\end{equation} 
Note that the conditional pdf $f(\mathbf{z}=\mathbf{T}(\mathbf{y}) \big| \mathbf{y})$ of the function value $\mathbf{z}=\mathbf{T}(\mathbf{y})$, 
given the observation $\mathbf{y}$, does not depend on the parameter vector $\mathbf{x}_{0}$. This corresponds to the fact 
that the random function $\mathbf{T}(\cdot)$ is statistically independent of the observation $\mathbf{y}$ for every parameter vector $\mathbf{x} \in \mathcal{X}$.  
For the special case of a deterministic mapping $\mathbf{T}(\mathbf{y})$, the definition \eqref{equ_def_expec_random_func_observation} reduces to 
\begin{equation}
\label{equ_def_expec_deterministic_func_observation} 
\mathsf{E}_{\mathbf{x}_{0}} \{ \mathbf{T}(\mathbf{y}) \} \triangleq  \int_{\mathbf{y}}Ê\mathbf{T}(\mathbf{y}) f(\mathbf{y}; \mathbf{x}_{0}) d \mathbf{y}. 
\end{equation} 

\item {\bf Estimation problem.} In some situations, we are not interested in the parameter $\mathbf{x}$ directly but rather the value $\mathbf{g}(\mathbf{x})$ of a vector-valued function 
\begin{equation} 
\label{equ_def_par_func}
\mathbf{g}(\cdot): \mathcal{X} \rightarrow \mathbb{R}^{P}
\end{equation}
of the true parameter. 
Consider, e.g., the situation where the vector $\mathbf{x}$ represents a transmitted signal and one is not interested in the signal itself but only if a signal has been transmitted, i.e., 
if the signal energy $g(\mathbf{x}) = \| \mathbf{x} \|^{2}_{2}$ is above a threshold.  

Let us formalize the notion of a \emph{classical estimation problem} by
\begin{definition} 
A (classical) estimation problem $\mathcal{E}$ is a triplet $\mathcal{E} = \vecestproblem$ consisting of the parameter set $\mathcal{X} \subseteq \mathbb{R}^{N}$, 
a statistical model $f(\mathbf{y};\mathbf{x})$, and the function $\mathbf{g}(\cdot): \mathcal{X} \rightarrow \mathbb{R}^{P}$ of which we would like to estimate the value $\mathbf{g}(\mathbf{x})$. 
\end{definition}

\item {\bf Estimator.} The estimation or inference of the value $\mathbf{g}(\mathbf{x})$ based on the observation $\mathbf{y}$ is performed by applying an estimator function $\hat{\mathbf{g}}(\cdot): \mathbb{R}^{M} \rightarrow \mathbb{R}^{P}$ to the 
observation $\mathbf{y}$. For a specific realization of the observation $\mathbf{y}$, the function value $\hat{\mathbf{g}}(\mathbf{y}) \in \mathbb{R}^{P}$ is an estimate of the value $\mathbf{g}(\mathbf{x})$. 
Note that the estimator function value $\hat{\mathbf{g}}(\mathbf{y})$ is a random quantity since it is obtained by a mapping applied to the random vector $\mathbf{y}$. 
In this thesis, we deal mainly with estimators defined via a deterministic mapping $\hat{\mathbf{g}}(\cdot)$. However, it is also possible to consider random or randomized estimators which are defined 
by a conditional pdf $f(\hat{\mathbf{g}} \big| \mathbf{y})$, i.e., even if we make the same observation $\mathbf{y}$ twice, the corresponding outputs of the randomized estimator may be different. 
Such randomized estimators may be useful especially if the parameter set $\mathcal{X}$ is finite. 

Implicit in the notation $f(\hat{\mathbf{g}} \big| \mathbf{y})$ for the conditional pdf of the estimator outcome $\hat{\mathbf{g}}$ given 
the observation $\mathbf{y}$, is the requirement that it \emph{cannot} depend on the parameter vector $\mathbf{x}$, which determines the statistics of $\mathbf{y}$ via the statistical model $f(\mathbf{y}; \mathbf{x})$ (see \cite[p.\ 33]{LC}). Strictly speaking, the correct notation 
for the conditional pdf of the estimator would be $f(\hat{\mathbf{g}}\big| \mathbf{y}; \mathbf{x})$, but since its values cannot depend on $\mathbf{x}$ we are allowed to write $f(\hat{\mathbf{g}}\big| \mathbf{y}; \mathbf{x})=f(\hat{\mathbf{g}} \big| \mathbf{y})$ with 
a suitable conditional pdf $f(\hat{\mathbf{g}} \big| \mathbf{y})$. Equivalently, a randomized estimator $\hat{\mathbf{g}}(\cdot)$ is a random function (a realization of $\hat{\mathbf{g}}(\cdot)$ is then a deterministic mapping $\mathbb{R}^{M} \rightarrow \mathbb{R}^{P}$), 
which is statistically independent of the observation $\mathbf{y}$ for every $\mathbf{x} \in \mathcal{X}$. 

We emphasize that all concepts and facts developed in this thesis that relate to the notion of an estimator, hold for deterministic as well as random(ized) estimators. 

\item {\bf Mean squared error.} We also need a measure to compare the quality of different estimators $\hat{\mathbf{g}}(\cdot)$. A popular choice for the performance measure that will also be used 
in this work is the \emph{mean squared error} (MSE) 
\begin{equation}
\varepsilon(\hat{\mathbf{g}}(\cdot);\mathbf{x}) \triangleq \mathsf{E}_{\mathbf{x}} \{ \| \hat{\mathbf{g}}(\mathbf{y}) - \mathbf{g}(\mathbf{x}) \|^{2}_{2} \}.
\end{equation} 

\item{\bf Bias.} 
Given an estimation problem $\mathcal{E}=\vecestproblem$, we define the \emph{bias} $\mathbf{b}(\hat{\mathbf{g}}(\cdot);\mathbf{x})$ of an estimator $\hat{\mathbf{g}}(\cdot)$ as 
\begin{equation}
\label{equ_def_bias} 
\mathbf{b}(\hat{\mathbf{g}}(\cdot);\mathbf{x}) \triangleq \mathsf{E}_{\mathbf{x}}  \{ \hat{\mathbf{g}}(\mathbf{y})  \} - \mathbf{g}(\mathbf{x}). 
\end{equation} 
An estimator $\hat{\mathbf{g}}(\cdot)$ whose bias vanishes identically, i.e., 
\begin{equation} 
\mathbf{b}(\hat{\mathbf{g}}(\cdot);\mathbf{x})  = \mathbf{0} \quad \quad \forall \mathbf{x} \in \mathcal{X}
\end{equation}
is called \emph{unbiased}.
The bias of an estimator represents the ``systematic'' or deterministic error which is incurred by the estimator. 

\item{\bf Variance and stochastic power.} 
Given an estimation problem $\mathcal{E}=\vecestproblem$, we define the \emph{variance} $v(\hat{\mathbf{g}}(\cdot); \mathbf{x})$ of an estimator $\hat{\mathbf{g}}(\cdot)$ as
\begin{equation} 
\label{equ_def_var} 
v(\hat{\mathbf{g}}(\cdot); \mathbf{x}) \triangleq \mathsf{E}_{\mathbf{x}} \big \{ \big \| \hat{\mathbf{g}}(\mathbf{y}) - \mathsf{E}_{\mathbf{x}}  \{ \hat{\mathbf{g}}(\mathbf{y})  \}  \big \|^{2}_{2} \big\}. 
\end{equation}
Similarly, we define the \emph{stochastic power} $P(\hat{\mathbf{g}}(\cdot); \mathbf{x})$ as 
\begin{equation} 
P(\hat{\mathbf{g}}(\cdot); \mathbf{x}) \triangleq \mathsf{E}_{\mathbf{x}} \big \{ \big \| \hat{\mathbf{g}}(\mathbf{y}) \}  \big \|^{2}_{2} \big\}. 
\end{equation}
The following relation between the estimator's mean, variance, and stochastic power can be verified easily: 
\begin{align}
\label{equ_rel_power_var}
P(\hat{\mathbf{g}}(\cdot); \mathbf{x})  = v(\hat{\mathbf{g}}(\cdot); \mathbf{x}) + \big \| \mathsf{E}_{\mathbf{x}} \big\{ \hat{\mathbf{g}}(\mathbf{y}) \}\|^{2}_{2}. 
\end{align} 
Note that two estimators that differ only on a set $\mathcal{A} \subseteq \mathbb{R}^{M}$ of measure zero yield the same bias, variance, MSE, and stochastic power at every $\mathbf{x} \in \mathcal{X}$. 
Therefore, when we say that two estimators are identical, we mean actually that they may differ only on a set of measure zero. 

\item {\bf Bias - variance tradeoff.} The MSE $\varepsilon(\hat{\mathbf{g}}(\cdot);\mathbf{x})$ of an estimator $\hat{\mathbf{g}}(\cdot)$ can be decomposed as the sum of  two nonnegative terms: 
\begin{equation} 
\label{equ_bias_var_decomp}
\varepsilon(\hat{\mathbf{g}}(\cdot);\mathbf{x}) = \| \mathbf{b} (\hat{\mathbf{g}}(\cdot);\mathbf{x}) \|^{2}_{2} + v(\hat{\mathbf{g}}(\cdot); \mathbf{x}), 
\end{equation}
where we used the estimator bias as defined in \eqref{equ_def_bias} and the estimator variance as defined in \eqref{equ_def_var}. 
For a small MSE $\varepsilon(\hat{\mathbf{g}}(\cdot);\mathbf{x})$, according to \eqref{equ_bias_var_decomp}, it is desirable to construct estimators that have both a small bias and a small variance. However, in general 
these are conflicting desiderata, i.e., an estimator with a small variance tends to have a large bias and vice versa. Therefore, to obtain a small MSE, one has to carefully balance between these two terms. It 
might be, e.g., that if one allows a small bias the variance can be reduced significantly, resulting in a smaller MSE as compared to unbiased estimators \cite{RethinkingBiasedEldar}. 
However, there are at least two situations where unbiased estimators are preferable. First, if the statistical model $f(\mathbf{y};\mathbf{x})$, 
viewed as a function of $\mathbf{y}$, is highly concentrated around its mean, then the bias effectively determines the MSE and thus a small MSE requires the estimator to be effectively unbiased. 
The second situation where unbiasedness is advantageous is when one can observe a large number $N$ of independent identically distributed (i.i.d.) realizations $\{ \mathbf{y}_{l} \}_{l \in [N]}$ of the observation $\mathbf{y}$. It 
can be shown that, without loss of generality, the estimation can then be based on the sample mean $\bar{\mathbf{y}} = \frac{1}{N} \sum\limits_{l \in [N]} \mathbf{y}_{l}$.\footnote{This is due to the fact that the 
sample mean $\bar{\mathbf{y}}$ is a sufficient statistic for the problem of estimating $\mathbf{x}$ from the observation $\{ \mathbf{y}_{l} \}_{l \in [N]}$ 
(cf.\ Section \ref{sec_suff_stat_RKHS}).}Under mild conditions on the statistical model $f(\mathbf{y}; \mathbf{x})$, it follows from the central limit theorem \cite{papoulis} that the pdf of $\bar{\mathbf{y}}$ 
becomes more and more concentrated around its mean as $N$ increases and therefore we have the same situation as before, i.e., 
the MSE is dominated by the bias and thus a small MSE requires an estimator to be effectively unbiased. 

\item {\bf Optimality of estimators.} 
Given a specific estimation problem $\mathcal{E}=\vecestproblem$, it is natural to look for an \emph{optimum} estimator $\hat{\mathbf{g}}(\cdot)$. Since our primary objective is a 
small MSE $\varepsilon(\hat{\mathbf{g}}(\cdot);\mathbf{x})$, we could be tempted to search for the specific estimator that minimizes the MSE, i.e., to find the specific estimator (if it exists) that solves the problem
\begin{equation} 
\label{equ_min_uniformly_MSE}
\min_{\hat{\mathbf{g}}(\cdot)} \,\, \varepsilon(\hat{\mathbf{g}}(\cdot);\mathbf{x}) 
\end{equation} 
for all $\mathbf{x} \in \mathcal{X}$ \emph{simultaneously}. 

However, it turns out that in general such an optimality criterion is not meaningful since there does not exist a single estimator 
that minimizes the MSE at all parameter vectors $\mathbf{x} \in \mathcal{X}$ simultaneously \cite{LC,RethinkingBiasedEldar}, i.e., the minimization problem in \eqref{equ_min_uniformly_MSE} has in general no solution.
This can be verified easily since the minimum MSE obtainable at any fixed parameter vector $\mathbf{x}_{0}$ is zero and 
achieved by the dumb estimator $\hat{\mathbf{g}}(\mathbf{y}) = \mathbf{g}(\mathbf{x}_{0})$ which always yields a constant value since it completely ignores the observation $\mathbf{y}$.  
Therefore, if there existed an optimum estimator in the sense of minimum MSE, it would have to yield zero MSE for all parameter vectors $\mathbf{x}$ and this is in general impossible. 
We thus have to consider alternative optimality criteria, two of which will be discussed in the following. 
\end{itemize}


\section{Minimax Estimation} 
\label{sec_minimax_est}

A widely used notion of estimator optimality is based on the concept of \emph{robustness}. By a robust estimator we mean an estimator that performs well under any operating condition \cite{VerduPoorMinimax,BlindMinimaxZvika}. 
For our scope, this means that we would like to keep control of the worst-case MSE of an estimator $\hat{\mathbf{g}}(\cdot)$ for a given estimation problem $\mathcal{E}=\vecestproblem$. This worst-case MSE is defined as 
\begin{equation}
R_{\mathcal{E}}(\hat{\mathbf{g}}(\cdot)) \triangleq \sup_{\mathbf{x} \in \mathcal{X}} \varepsilon(\hat{\mathbf{g}}(\cdot); \mathbf{x}). 
\end{equation}
The optimal estimator in terms of robustness, called the \emph{minimax estimator}, is defined as the solution to the problem \cite{LC}
\begin{equation}
\label{equ_opt_minimax}
\arginf_{\hat{\mathbf{g}}(\cdot)} R_{\mathcal{E}} (\hat{\mathbf{g}}(\cdot)) =  \arginf_{\hat{\mathbf{g}}(\cdot)} \sup_{\mathbf{x} \in \mathcal{X}} \varepsilon(\hat{\mathbf{g}}(\cdot); \mathbf{x}).
\end{equation} 
The infimum value $R_{\mathcal{E}} \triangleq \inf_{\hat{\mathbf{g}}(\cdot)} R_{\mathcal{E}} (\hat{\mathbf{g}}(\cdot))$ is called the \emph{minimax risk} of the estimation problem $\mathcal{E}$. 
Finding minimax estimators or even characterizing the minimax risk (e.g., by lower bounds) is often difficult. 
However, for the estimation problem corresponding to the \emph{sparse linear model}, which will be discussed in detail in Chapter \ref{chap_SLM}, some recent work (see e.g. \cite{RaskuttiMinmaxSLM,VerzelenMinmaxSLM,DantzigCandes,Meinshausen09lasso-typerecovery}) 
derives sharp lower bounds on the minimax risk and presents practical estimation schemes that are close to the minimax estimator. 
In this regard, also the authors of \cite{UnifiedMEstMinimax,ZvikaCoherenceTSP} present performance evaluations of existing estimation schemes revealing that they are close to the minimax estimator. 


\section{Minimum Variance Estimation} 
\label{sec_min_var_est}

An optimality criterion different from the minimax criterion is based on the bias-variance decomposition of the MSE in \eqref{equ_bias_var_decomp}.
In particular, it is common to fix the estimator bias, i.e., to require that
\begin{equation} 
\mathbf{b}(\hat{\mathbf{g}}(\cdot);\mathbf{x}) \stackrel{!}{=} \mathbf{c}(\mathbf{x}) \quad \quad \forall \mathbf{x} \in \mathcal{X},
\end{equation} 
with a given function $\mathbf{c}(\cdot): \mathcal{X} \rightarrow \mathbb{R}^{P}$ (called the prescribed bias function), and then look for estimators that minimize the MSE.

If, for a given estimation problem $\mathcal{E}= \vecestproblem$, one considers only estimators with the same prescribed bias function $\mathbf{c}(\mathbf{x})$, 
minimizing the MSE $\varepsilon(\hat{\mathbf{g}}(\cdot);\mathbf{x})$ is completely equivalent to minimizing the variance $v(\hat{\mathbf{g}}(\cdot);\mathbf{x})$. 
It is often impossible to perform a uniform minimization of the variance, and therefore we will pursue a local approach: 
We consider a specific parameter vector $\mathbf{x}_{0} \in \mathcal{X}$ and try to minimize the variance $v(\hat{\mathbf{g}}(\cdot);\mathbf{x}_{0})$ at $\mathbf{x}_{0}$ of estimators $\hat{\mathbf{g}}(\cdot)$ whose 
bias is equal to the prescribed function $\mathbf{c}(\cdot): \mathcal{X} \rightarrow \mathbb{R}^{P}$, i.e., $\mathbf{b}(\hat{\mathbf{g}}(\cdot);\mathbf{x}_{0}) = \mathbf{c}(\mathbf{x})$ for all $\mathbf{x} \in \mathcal{X}$. 

We make some important definitions: 
\begin{definition}
A ``minimum variance estimation problem'' or ``minimum variance problem'' $\mathcal{M}$ is the triplet
\begin{equation}
\mathcal{M} \triangleq \big( \mathcal{E}, \mathbf{c}(\mathbf{x}), \mathbf{x}_{0} \big) 
\end{equation} 
consisting of an estimation problem $\mathcal{E} = \vecestproblem$, a prescribed bias function $\mathbf{c}(\cdot): \mathcal{X} \rightarrow \mathbb{R}^{P}$, and a fixed
parameter vector $\mathbf{x}_{0} \in \mathcal{X}$ at which we will try to minimize the variance of estimators whose bias is equal to $\mathbf{c}(\mathbf{x})$ for all $\mathbf{x} \in \mathcal{X}$. 
\end{definition}
\begin{definition} 
\label{def_est_finite_var_prescr_bias}
Given a minimum variance problem $\mathcal{M}= \minvarproblem$, we denote by $\mathcal{F}(\mathcal{M})$ the set of all (possibly randomized) estimators with the prescribed bias function $\mathbf{c}(\mathbf{x})$ and 
with finite variance at $\mathbf{x}_{0}$, i.e.,  
\begin{equation}
\label{equ_est_finite_var_prescr_bias}
\mathcal{F}(\mathcal{M}) \triangleq \{ \hat{\mathbf{g}}(\cdot) \big| v(\hat{\mathbf{g}}(\cdot);\mathbf{x}_{0}) < \infty \mbox{, } \mathbf{b}(\hat{\mathbf{g}}(\cdot);\mathbf{x}) = \mathbf{c}(\mathbf{x}) \,\, \forall \mathbf{x} \in \mathcal{X}  \}. 
\end{equation}
We will refer to the set $\mathcal{F}(\mathcal{M})$ as the set of allowed estimators for the minimum variance problem $\mathcal{M}$. 
\end{definition} 

The theoretically achievable performance of any estimator for a minimum variance problem is characterized by 
\begin{definition}
\label{def_min_ach_var}
Given a minimum variance problem $\mathcal{M} = \minvarproblem$, 
the minimum achievable variance at the parameter vector $\mathbf{x}_{0}$ is defined as 
\begin{equation}
\label{equ_def_min_ach_var}
\minachievevar \triangleq \inf_{\hat{\mathbf{g}}(\cdot) \in \mathcal{F}(\mathcal{M})}   v(\hat{\mathbf{g}}(\cdot);\mathbf{x}_{0}).
\end{equation} 
If there is exists no allowed estimator for $\mathcal{M}$, i.e., the set $\mathcal{F}(\mathcal{M})$ is empty, we define $\minachievevar \triangleq \infty$. 
\end{definition} 
In the literature, the minimum achievable variance is sometimes referred to as the \emph{Barankin bound} (in particular if unbiased estimation is considered). 

The most important notion of estimator optimality in the context of minimum variance estimation is stated in 
\begin{definition}[Locally Minimum Variance Estimation]
\label{def_LMV}
Given a minimum variance problem $\mathcal{M}=\minvarproblem$, we call any (possibly randomized) estimator $\hat{\mathbf{g}}^{(\mathbf{x}_{0})}(\cdot) \in  \mathcal{F}(\mathcal{M})$  
whose variance at $\mathbf{x}_{0}$ attains $L_{\mathcal{M}}$, i.e., 
\begin{equation} 
v(\hat{\mathbf{g}}^{(\mathbf{x}_{0})}(\cdot); \mathbf{x}_{0}) = L_{\mathcal{M}},
\end{equation} 
a locally minimum variance (LMV) estimator for $\mathcal{M}$. For the special case where the prescribed bias function is identically zero, i.e., $\mathbf{c}(\cdot) \equiv \mathbf{0}$, 
such an estimator is called a locally minimum variance unbiased (LMVU) estimator.
\end{definition} 
While an LMV estimator in general performs well (has a small variance) only locally around a specific parameter vector $\mathbf{x}_{0}$, 
there may be estimators which have the prescribed bias $\mathbf{c}(\cdot)$ and moreover perform well (have a small variance) everywhere, i.e., for any parameter value $\mathbf{x} \in \mathcal{X}$: 
\begin{definition}[Uniformly Minimum Variance Estimation] 
Given an estimation problem $\mathcal{E}=\vecestproblem$, a prescribed bias function $\mathbf{c}(\cdot): \mathcal{X} \rightarrow \mathbb{R}^{P}$ and a parameter vector $\mathbf{x}_{0} \in \mathcal{X}$, 
we denote by $\mathcal{M}(\mathbf{x}_{0})$ the specific minimum variance problem that is given by $\mathcal{M}(\mathbf{x}_{0})=\left(\mathcal{E},\mathbf{c}(\cdot), \mathbf{x}_{0} \right)$.
Then, we call any estimator  $\hat{\mathbf{g}}(\cdot): \mathbb{R}^{M} \rightarrow \mathbb{R}^{P}$ which is an LMV estimator for $\mathcal{M}(\mathbf{x}_{0})$ simultaneously  for every $\mathbf{x}_{0} \in \mathcal{X}$, 
a uniformly minimum variance (UMV) estimator for the given estimation problem $\mathcal{E}$ and  bias function $\mathbf{c}(\cdot)$. 
For the special case where the bias function is identically zero, i.e., $\mathbf{c}(\cdot) \equiv \mathbf{0}$, 
such an estimator is called a uniformly minimum variance unbiased (UMVU) estimator for the estimation problem $\mathcal{E}$.
\end{definition} 

\subsection{Vector-Valued vs. Scalar-Valued Parameter Function}
\label{sec_sep_vector_scalar_param_function} 
So far, we considered the estimation of the value $\mathbf{g}(\mathbf{x})$ of a vector-valued function $\mathbf{g}(\cdot): \mathcal{X} \rightarrow \mathbb{R}^{P}$ of the parameter vector $\mathbf{x}$. 
However, the minimum variance estimation of a $P$-dimensional vector-valued function $\mathbf{g}(\mathbf{x}) = \big( g_{1}(\mathbf{x}), \ldots, g_{P}(\mathbf{x}) \big)^{T}$ 
can be reduced to separately estimating the $P$ scalar functions $\{ g_{k}(\mathbf{x})\}_{k \in [P]}$. This is possible because on the one hand we have that 
the variance $v(\hat{\mathbf{g}}(\cdot);\mathbf{x})$ can be decomposed as $v(\hat{\mathbf{g}}(\cdot);\mathbf{x}) = \sum_{k \in [P]} v(\hat{g}_{k}(\cdot);\mathbf{x})$, 
and on the other hand the bias constraint decomposes elementwise, i.e., 
$\mathbf{b}(\hat{\mathbf{g}}(\cdot);\mathbf{x}) = \mathbf{c}(\mathbf{x})$ if and only if $b(\hat{g}_{k}(\cdot);\mathbf{x}) =  c_{k}(\mathbf{x})$ for every $k \in [P]$. 
Thus, given a minimum variance problem $\mathcal{M}=\minvarproblem$ with prescribed bias $\mathbf{c}(\cdot): \mathcal{X} \rightarrow \mathbb{R}^{P}$, 
we have that the minimum achievable variance at $\mathbf{x}_{0}$ is given as 
\begin{equation}
\label{equ_sum_min_var_scalar_min_var}
L_{\mathcal{M}} = \sum_{k \in [P]} L_{\mathcal{M}_{k}},
\end{equation} 
where $L_{\mathcal{M}_{k}}$ denotes the minimum achievable variance for the ``scalar'' minimum variance problem 
$\mathcal{M}_{k} = \left( \mathcal{E}_{k}, c_{k}(\cdot),  \mathbf{x}_{0} \right)$, 
with $\mathcal{E}_{k}$ denoting the associated scalar estimation problems, i.e, 
$\mathcal{E}_{k} \triangleq  \left( \mathcal{X}, f(\mathbf{y}; \mathbf{x}), g_{k} (\cdot) \right)$.
Here, $g_{k}(\cdot)$ and $c_{k}(\cdot)$ denote the $k$th component of the parameter function $\mathbf{g}(\cdot)$ and the prescribed bias function $\mathbf{c}(\cdot)$, respectively.
Furthermore, if the estimators $\{ \hat{g}_{k}(\cdot) \}_{k \in [P]}$ are LMV estimators for the minimum variance problems $\mathcal{M}_{k}$, then 
the vector-valued estimator given by $\hat{\mathbf{g}}(\cdot) =  \big( \hat{g}_{1}(\cdot), \ldots, \hat{g}_{P}(\cdot) \big)^{T}$ is the LMV estimator for $\mathcal{M}$, and vice versa. 

Therefore, in the context of minimum variance estimation, we can restrict ourselves without loss of generality to estimation problems $\mathcal{E} = \scalarestproblem$ 
with a scalar-valued parameter function $g(\cdot): \mathcal{X} \rightarrow \mathbb{R}$, i.e., we can assume $P=1$ in \eqref{equ_def_par_func}. 
Note that the restriction to scalar-valued parameter functions is in general not possible in the context of minimax estimation 
since the optimization problem \eqref{equ_opt_minimax} does not decompose in general. 

\subsection{Lower Bounds via Generalized Inequalities}

Some results on minimum variance estimation (e.g., \cite{GormanHero,ZvikaCRB}) are stated in terms of generalized (matrix) inequalities \cite{BoydConvexBook} involving the estimator covariance matrix 
$\mbox{cov}_{\mathbf{x}_{0}} \{ \hat{\mathbf{g}}(\mathbf{y}) \}  \in  \mathbb{R}^{P \times P}$ defined as
\begin{equation} 
\mbox{cov}_{\mathbf{x}_{0}} \{ \hat{\mathbf{g}}(\mathbf{y}) \} 
\triangleq \mathsf{E}_{\mathbf{x}_{0}} \big\{ \big[ \hat{\mathbf{g}}(\mathbf{y}) - \mathsf{E}_{\mathbf{x}_{0}}Ê\{ \hat{\mathbf{g}}(\mathbf{y}) \} \big] \big[ \hat{\mathbf{g}}(\mathbf{y}) - \mathsf{E}_{\mathbf{x}_{0}}Ê\{ \hat{\mathbf{g}}(\mathbf{y}) \}\big]^{T} \big \}.
\end{equation}  
Note that $\mbox{cov}_{\mathbf{x}_{0}} \{ \hat{\mathbf{g}}(\mathbf{y}) \}$ is a psd matrix by its very definition.
A typical result is then a lower bound of the form 
\begin{equation} 
\label{equ_bound_gen_inequ_cov}
\mbox{cov}_{\mathbf{x}_{0}} \{ \hat{\mathbf{g}}(\mathbf{y}) \} \geq \mathbf{L},
\end{equation} 
where the lower bound $\mathbf{L} \in \mathbb{R}^{P \times P}$ is a psd matrix that depends of course on the underlying estimation problem.  
By contrast, we derive exclusively lower bounds on the variance $v(\hat{\mathbf{g}}(\cdot);\mathbf{x}_{0})$, which is always a scalar quantity. 

However, since $v(\hat{\mathbf{g}}(\cdot);\mathbf{x}_{0}) = \tracem {\mbox{cov}_{\mathbf{x}_{0}} \{ \hat{\mathbf{g}}(\mathbf{y}) \} }$, we have that 
any bound of the form \eqref{equ_bound_gen_inequ_cov} induces also a lower bound 
on $v(\hat{\mathbf{g}}(\cdot);\mathbf{x}_{0})$ via 
\begin{equation}
v(\hat{\mathbf{g}}(\cdot);\mathbf{x}_{0}) = \tracem{ \mbox{cov}_{\mathbf{x}_{0}} \{ \hat{\mathbf{g}}(\mathbf{y}) \} } \geq \tracem{\mathbf{L}},
\end{equation} 
where the last inequality follows from the fact that for two psd matrices $\mathbf{A}$, $\mathbf{B} \in \mathbb{R}^{P \times P}$, 
we have $\mathbf{A} \geq \mathbf{B} \Rightarrow \tracem { \mathbf{A} } \geq \tracem { \mathbf{B} }$ \cite{BoydConvexBook}. 

\subsection{Exchangeability of Prescribed Bias and Parameter Function} 
\label{sec_equ_bias_param_function}

A minimum variance problem $\mathcal{M}$ depends on the prescribed bias $c(\cdot): \mathcal{X} \rightarrow \mathbb{R}$ 
and the parameter function $g(\cdot): \mathcal{X} \rightarrow \mathbb{R}$. 
However, the minimum achievable variance $L_{\mathcal{M}}$ as well 
as the corresponding LMV estimator do not depend \emph{independently} on $c(\cdot)$ and $g(\cdot)$, as stated in 
\begin{theorem} 
\label{thm_equ_bias_param_function}
Consider two estimation problems $\mathcal{E}=\left( \mathcal{X},f(\mathbf{y};\mathbf{x}),g(\cdot)\right)$ and 
$\mathcal{E}'=\left( \mathcal{X},f(\mathbf{y};\mathbf{x}),g'(\cdot) \right)$ which share a common parameter set $\mathcal{X}$ and statistical model $f(\mathbf{y}; \mathbf{x})$. 
We have that any estimator $\hat{g}(\cdot)$ for $\mathcal{E}$ with bias $c(\cdot)$ is also an estimator for $\mathcal{E}'$ with bias $c'(\cdot)=c(\cdot) + g(\cdot) - g'(\cdot)$. 
For any two minimum variance problems $\mathcal{M}= \minvarproblemscalar$, $\mathcal{M}' = \left( \mathcal{E}', c'(\cdot), \mathbf{x}_{0}\right)$ 
associated with $\mathcal{E}$ and $\mathcal{E}'$, respectively, with a common fixed parameter vector $\mathbf{x}_{0} \in \mathcal{X}$ but 
different prescribed bias functions $c(\cdot)$ and $c'(\cdot)$, we have that if 
\begin{equation}
\label{equ_same_mean_two_min_var_problems}
g(\cdot) + c(\cdot) = g'(\cdot) + c'(\cdot),
\end{equation} 
then the minimum achievable variance is the same for $\mathcal{M}$ and $\mathcal{M}'$, i.e., 
\begin{equation}
L_{\mathcal{M}} = L_{\mathcal{M}'}.  
\end{equation} 
Furthermore, if an estimator $\hat{g}(\cdot)$ is an LMV estimator for $\mathcal{M}$, it is also an LMV estimator for $\mathcal{M}'$. 
\end{theorem} 
\begin{proof} 
Consider an estimator $\hat{g}(\cdot)$ for $\mathcal{E}$ with bias $c(\cdot)$, i.e., $b(\hat{g}(\cdot); \mathbf{x}) = c(\mathbf{x})$ for any $\mathbf{x} \in \mathcal{X}$. 
This implies by \eqref{equ_def_bias} that $\mathsf{E}_{\mathbf{x}_{0}} \{ \hat{g}( \mathbf{y}) \} = g(\mathbf{x}) + c(\mathbf{x})$ for any $\mathbf{x} \in \mathcal{X}$. 
However, if we use the same estimator for $\mathcal{E}'$ we have that its bias is $b(\hat{g}(\cdot);\mathbf{x}) = \mathsf{E}_{\mathbf{x}_{0}} \{ \hat{g}( \mathbf{y}) \} - g'(\mathbf{x}) =  g(\mathbf{x}) + c(\mathbf{x}) - g'(\mathbf{x})$. 
Furthermore, since both estimation problems $\mathcal{E}$ and $\mathcal{E}'$ share the same statistical model $f(\mathbf{y}; \mathbf{x})$, 
the variance $v(\hat{g}(\cdot);\mathbf{x})$ (cf.\ \eqref{equ_def_var}) is the same for both estimation problems. 
It follows that if the two minimum variance problems $\mathcal{M}$ and $\mathcal{M}'$ satisfy \eqref{equ_same_mean_two_min_var_problems}, 
the associated sets of allowed estimators coincide, i.e., $\mathcal{F}(\mathcal{M}) = \mathcal{F}(\mathcal{M}')$, and moreover for each allowed estimator, 
the corresponding value of the objective in \eqref{equ_def_min_ach_var} is the same for $\mathcal{M}$ and $\mathcal{M}'$, i.e., $L_{\mathcal{M}} = L_{\mathcal{M}'}$. 
\end{proof} 

According to the definition of the minimum achievable variance in \eqref{equ_def_min_ach_var}, it is possible that for a given minimum variance problem $\mathcal{M}= \minvarproblem$ 
there is no LMV estimator, i.e., no estimator whose bias is equal to $\mathbf{c}(\mathbf{x})$ for all $\mathbf{x} \in \mathcal{X}$ and whose variance at $\mathbf{x}_{0}$ is finite and equal to $L_{\mathcal{M}}$. 
This motivates
\begin{definition} 
\label{def_valid_bias_func_classic_est}
Given a minimum variance problem $\mathcal{M} = \minvarproblem$, we call the prescribed bias function $\mathbf{c}(\cdot): \mathcal{X} \rightarrow \mathbb{R}^{P}$ valid for $\mathcal{M}$ if 
$L_{\mathcal{M}}$ is finite and there exists at least one estimator with bias $\mathbf{c}(\mathbf{x})$ whose variance at $\mathbf{x}_{0}$ equals $L_{\mathcal{M}}$.
\end{definition} 
Similarly, we have \cite{Duttweiler73b,Parzen59}
\begin{definition} 
\label{def_estimable_par_function}
Given a minimum variance problem with zero bias, $\mathcal{M}=\left( \mathcal{E}, \mathbf{c}(\cdot) \equiv 0,\mathbf{x}_{0} \right)$, with associated estimation problem $\mathcal{E} = \vecestproblem$, 
we call a parameter function $\mathbf{g}(\mathbf{x})$ estimable for $\mathcal{M}$ if 
$L_{\mathcal{M}}$ is finite and there exists at least one unbiased estimator for $g(\mathbf{x})$ whose variance at $\mathbf{x}_{0}$ equals $L_{\mathcal{M}}$.
\end{definition} 
Due to Theorem \ref{thm_equ_bias_param_function}, we have 
\begin{corollary}
\label{cor_estimable_equiv_valid}
Given the minimum variance problem $\mathcal{M}= \minvarproblemscalar$ associated with the estimation problem $\mathcal{E} = \scalarestproblem$, 
the prescribed bias function $c(\cdot): \mathcal{X} \rightarrow \mathbb{R}$ 
is valid for $\mathcal{M}$ if and only if the parameter function $g'(\cdot) = g(\cdot) + c(\cdot)$ of the estimation problem $\mathcal{E}'= \left(\mathcal{X}, f(\mathbf{y}; \mathbf{x}), g(\cdot) \right)$ 
is estimable for the minimum variance problem $\mathcal{M}' = \left( \mathcal{E}', c'(\cdot) \equiv 0, \mathbf{x}_{0} \right)$. 
\end{corollary}

\subsection{Existence and Uniqueness of Minimum Variance Estimators}
\label{sec_esist_uniqu_mve}

As already mentioned above, it is not guaranteed that there exists an LMV estimator for a given minimum variance problem $\mathcal{M}=\minvarproblem$ associated to an estimation problem 
$\mathcal{E} = \vecestproblem$. Discussions on the conditions for the existence of an LMV estimator can be found in the seminal works \cite{stein50,Barankin49}. 
It is interesting to note that there need not exist an LMV estimator even if the set $\mathcal{F}(\mathcal{M})$ is nonempty, i.e., even if there is at least one estimator $\hat{\mathbf{g}}(\cdot)$ with bias $\mathbf{c}(\mathbf{x})$ 
and whose variance at $\mathbf{x}_{0}$ is finite. An example of a minimum variance problem where this is the case can be found in \cite[p. 407]{stein50}. 
The difficulty in establishing the existence of an LMV estimator, i.e., an allowed estimator whose variance at $\mathbf{x}_{0}$ attains the infimum in \eqref{equ_def_min_ach_var}, 
results from the nontrivial topological characterization of the set $\mathcal{F}(\mathcal{M})$ of allowed estimators. 
While one can easily verify that this set of estimators is an affine subset \cite[p. 21]{BoydConvexBook} of a specific Hilbert space, it is in general not true that this set is closed, i.e., there might exist a Cauchy sequence of estimators 
belonging to $\mathcal{F}(\mathcal{M})$ which has no limit that belongs to $\mathcal{F}(\mathcal{M})$.\footnote{We will introduce the concepts 
of closedness, limits, and Cauchy sequences in Chapter \ref{chap_RKHS_Fund}.} We will discuss this issue in more detail in Section \ref{sec_hilbert_space_est}. 
In any case, if a LMV estimator exists it is unique, as stated in 
\begin{theorem} 
\label{thm_uniqueness_LMV}
Consider a minimum variance problem $\mathcal{M}=\minvarproblem$.
If an LMV estimator exists for $\mathcal{M}$, it is unique. 
\end{theorem} 
\begin{proof} 
The uniqueness of an LMVU estimator, i.e., the uniqueness of an LMV estimator for the special case of a minimum variance problem with zero bias, $\mathbf{c}(\mathbf{x}) = \mathbf{0}$, has been proven 
in the seminal works \cite{Barankin49,stein50}. The statement follows then for an arbitrary prescribed bias function $\mathbf{c}(\cdot): \mathcal{X} \rightarrow \mathbb{R}^{P}$ by Theorem \ref{thm_equ_bias_param_function}. 
\end{proof} 
A direct consequence of Theorem \ref{thm_uniqueness_LMV}, which is also proved in \cite{papoulis} is stated 
in 
\begin{theorem}
Consider an estimation problem $\mathcal{E}=\vecestproblem$ and a prescribed bias function $\mathbf{c}(\cdot): \mathcal{X} \rightarrow \mathbb{R}^{P}$. 
If an UMV estimator exists for $\mathcal{E}$ and $\mathbf{c}(\cdot)$, it is unique. 
\end{theorem} 
\begin{proof}
This result follows straightforwardly from Theorem \ref{thm_uniqueness_LMV} since the UMV is necessarily also the unique 
LMV estimator for any minimum variance problem $\mathcal{M}=\minvarproblem$, which is obtained from $\mathcal{E}$ and $\mathbf{c}(\cdot)$ by choosing an arbitrary value for $\mathbf{x}_{0} \in \mathcal{X}$. 
If there would be two different UMV estimators, they would be at the same time two different LMV estimators for $\mathcal{M}$ which contradicts Theorem \ref{thm_uniqueness_LMV}. 
\end{proof}

\subsection{Transformations of the Parameter Function} 

Consider an estimation problem $\mathcal{E}_{1}= \left(\mathcal{X},f(\mathbf{y}; \mathbf{x}),g_{1}(\cdot) \right)$ and a modified 
estimation problem $\mathcal{E}_{2}= \left( \mathcal{X}, f(\mathbf{y}; \mathbf{x}),g_{2}(\cdot) Ê\right)$ that is obtained from $\mathcal{E}_{1}$ 
by using the parameter function $g_{2}(\mathbf{x})$ instead of $g_{1}(\mathbf{x})$. If the two parameter functions are related by $g_{2}(\mathbf{x}) =g_{1}(\mathbf{x}) + a$, where $a \in \mathbb{R}$ is a fixed constant, 
we have obviously that if an estimator $\hat{g}_{1}(\cdot)$ is an unbiased estimator for the estimation problem $\mathcal{E}_{1}$, then the estimator $\hat{g}_{2}(\cdot) \triangleq \hat{g}_{1}(\cdot) +a$ is an unbiased estimator for the estimation problem $\mathcal{E}_{2}$ (cf.\ Theorem \ref{thm_equ_bias_param_function}).
Moreover, this estimator has the same variance as the original estimator, i.e., 
$v_{\mathcal{E}_{2}}(\hat{g}_{2}(\cdot);\mathbf{x}) = v_{\mathcal{E}_{1}}(\hat{g}_{1}(\cdot); \mathbf{x})$
for every $\mathbf{x} \in \mathcal{X}$, where $v_{\mathcal{E}}(\hat{g}(\cdot); \mathbf{x})$ denotes the variance $v(\hat{g}(\cdot); \mathbf{x})$ at $\mathbf{x}$ of an estimator $\hat{g}(\cdot)$ that uses an observation $\mathbf{y}$ whose statistics 
are determined by the statistical model $f(\mathbf{y}; \mathbf{x})$ of the estimation problem $\mathcal{E}$.
By a slight generalization of this simple example, we arrive at  
\begin{theorem}
\label{thm_aff_trafo_parameter_function}
Consider two estimation problems $\mathcal{E}_{1}= \left(\mathcal{X},f(\mathbf{y}; \mathbf{x}),\mathbf{g}_{1}(\cdot) \right)$, 
$\mathcal{E}_{2}= \left(\mathcal{X},f(\mathbf{y}; \mathbf{x}),\mathbf{g}_{2}(\cdot) \right)$ and the associated minimum variance problems 
$\mathcal{M}_{1} = \left( \mathcal{E}_{1}, \mathbf{c}(\cdot), \mathbf{x}_{0} \right)$, $\mathcal{M}_{2} = \left( \mathcal{E}_{2}, \mathbf{c}(\cdot), \mathbf{x}_{0} \right)$ that are identical except 
for the parameter functions $\mathbf{g}_{1}(\mathbf{x})$ and $\mathbf{g}_{2}(\mathbf{x})$, respectively. 
If the parameter functions are related by 
\begin{equation}
\mathbf{g}_{2}(\cdot) = \mathbf{U} \mathbf{g}_{1}(\cdot) + \mathbf{a} 
\end{equation}
where $\mathbf{a} \in \mathbb{R}^{P}$ is a constant vector and $\mathbf{U} \in \mathbb{R}^{P \times P}$ is an orthonormal matrix, i.e., $\mathbf{U}^{T} \mathbf{U} = \mathbf{I}$, then we have 
that 
\begin{equation}
L_{\mathcal{M}_{1}} = L_{\mathcal{M}_{2}}.
\end{equation}
Moreover, if an estimator $\hat{\mathbf{g}}_{1}(\cdot)$ is the LMV estimator for $\mathcal{M}_{1}$, then the estimator $\hat{\mathbf{g}}_{2} (\cdot) \triangleq \mathbf{U}\hat{\mathbf{g}}_{1} (\cdot)+ \mathbf{a}$ is 
the LMV estimator for $\mathcal{M}_{2}$. 
\end{theorem} 
\begin{proof} 
The statement is proven by showing that there is a one-to-one correspondence between the sets $\mathcal{F}(\mathcal{M}_{1})$ and $\mathcal{F}(\mathcal{M}_{2})$. Indeed, 
to any estimator $\hat{\mathbf{g}}_{1}(\cdot) \in \mathcal{F}(\mathcal{M}_{1})$, we can associate uniquely the estimator $\hat{\mathbf{g}}_{2}(\cdot) \triangleq \mathbf{U} \hat{\mathbf{g}}_{1} (\cdot)+ \mathbf{a} \in \mathcal{F}(\mathcal{M}_{2})$. 
Moreover, the value of the objective in \eqref{equ_def_min_ach_var}, i.e., $v(\hat{\mathbf{g}}(\cdot);\mathbf{x}_{0})$, is the same for $\hat{\mathbf{g}}_{1}(\cdot)$ and $\hat{\mathbf{g}}_{2}(\cdot)$.
\end{proof}

\subsection{Relaxing The Bias Equality Constraints}

Instead of requiring the estimator bias to equal a given function $\mathbf{c}(\mathbf{x})$ for every $\mathbf{x} \in \mathcal{X}$, one may pose inequality constraints on the bias. 
This is pursued by the authors of \cite{HeroUniformCRB} and \cite{EldarUniformCRB}, who consider minimum variance estimation with inequality constraints placed on the bias gradient. 
More specifically, they require the norm of the bias gradient matrix to be upper bounded by a prescribed maximum value.


\section{Transformation of the Observation}  

In some situations it might be desirable not to work with the observation $\mathbf{y}$ of a given estimation problem directly, but to perform some 
``preprocessing'' yielding a modified or transformed observation $\mathbf{z}$ which is easier to process for some reason. As an example for such a preprocessing consider e.g. a compression or dimension reduction 
as is done within CS by multiplying the observation $\mathbf{y} \in \mathbb{R}^{M}$ with a CS measurement matrix $\mathbf{M} \in \mathbb{R}^{L \times M}$ (typically $L \ll M$) 
to yield a small number of compressive measurements $z_{k}$, with $\mathbf{z} = \mathbf{M} \mathbf{y}$. 
However, in general we cannot expect that the achievable estimation accuracy remains unchanged if we use 
the transformed observation. The next two subsections discuss the influence of a transformation of the observation, in the context of minimum variance estimation, for the two cases 
where the transformation is invertible and where it is not. 

\subsection[Invertible Transformation]{Invertible Transformation -- Invariance of Classical Estimation Problems} 
\label{sec_invariance_mve}
Consider the problem of estimating the value $g(\mathbf{x})$ of the  parameter function $g(\cdot): \mathcal{X} \rightarrow \mathbb{R}$ using the observation $\mathbf{y} \in \mathbb{R}^{M}$. Then, intuitively, 
the estimation problem does not become any harder if we use not directly $\mathbf{y}$ but a transformed observation $\mathbf{z} =Ê\mathbf{T}(\mathbf{y})$, where 
$\mathbf{T}(\cdot): \mathbb{R}^{M} \rightarrow \mathbb{R}^{K}$ denotes an invertible deterministic function (this requires that $K \geq M$), as the new observation. Indeed, we have 
\begin{theorem}
\label{thm_invariance_classical_est}
Consider an estimation problem $\mathcal{E} = \vecestproblem$ with a statistical model $f(\mathbf{y}; \mathbf{x})$ and an invertible deterministic mapping $\mathbf{T}(\cdot): \mathbb{R}^{M} \rightarrow \mathbb{R}^{K}$. 
We denote by $\mathcal{E}'$ the estimation problem that is obtained from $\mathcal{E}$ by using the vector $\mathbf{z} =Ê\mathbf{T}(\mathbf{y})$ as the observation, i.e., 
$\mathcal{E}' = \left( \mathcal{X}, f(\mathbf{z} ; \mathbf{x}), \mathbf{g}(\cdot) \right)$ where $f(\mathbf{z};\mathbf{x})$ denotes the pdf of the new observation $\mathbf{z} = \mathbf{T}(\mathbf{y})$.
We then have that for any two minimum variance problems $\mathcal{M} = \left( \mathcal{E}, \mathbf{c}(\cdot), \mathbf{x}_{0}) \right)$, $\mathcal{M}' = \left( \mathcal{E}', \mathbf{c}(\cdot), \mathbf{x}_{0}) \right)$ which 
share a common prescribed bias function $\mathbf{c}(\cdot): \mathcal{X} \rightarrow \mathbb{R}^{P}$ and a fixed parameter vector $\mathbf{x}_{0} \in \mathcal{X}$,
\begin{equation}
L_{\mathcal{M}} = L_{\mathcal{M}'}. 
\end{equation}  
Moreover, the minimax risks of $\mathcal{E}$ and $\mathcal{E}'$ satisfy 
\begin{equation} 
R_{\mathcal{E}} = R_{\mathcal{E}'}.
\end{equation} 
\end{theorem}
\begin{proof}
This statement if proven by showing that the optimization problems \eqref{equ_def_min_ach_var} and \eqref{equ_opt_minimax} which 
define the quantities $L_{\mathcal{M}}$, $L_{\mathcal{M}'}$ and $R_{\mathcal{E}}$, $R_{\mathcal{E}'}$ respectively are equivalent for $\mathcal{E}$, $\mathcal{E}'$ 
and $\mathcal{M}$, $\mathcal{M}'$. Indeed, for any estimator $\hat{\mathbf{g}}(\mathbf{y})$ for $\mathcal{E}$ and $\mathcal{M}$ 
that uses the observation $\mathbf{y}$ we have that the specific estimator $\hat{\mathbf{g}}'(\mathbf{z}) \triangleq \hat{\mathbf{g}}(\mathbf{T}^{-1}(\mathbf{z}))$ that uses the observation 
$\mathbf{z}=\mathbf{T}(\mathbf{y})$ has the same bias and variance at any $\mathbf{x} \in \mathcal{X}$. 
The same goes the other way, i.e., given any estimator $\hat{\mathbf{g}}'(\mathbf{z})$ for $\mathcal{E}'$ and $\mathcal{M}'$ 
that uses the observation $\mathbf{z}$ we have that the specific estimator $\hat{\mathbf{g}}(\mathbf{y}) \triangleq \hat{\mathbf{g}}'(\mathbf{T}(\mathbf{y}))$ that uses the observation $\mathbf{y}$ 
has the same bias and variance at any $\mathbf{x} \in \mathcal{X}$.
Thus we have that $\hat{\mathbf{g}}'(\mathbf{z})  \in \mathcal{F}(\mathcal{M}') \Leftrightarrow \hat{\mathbf{g}}(\mathbf{y}) \inÊ\mathcal{F}(\mathcal{M})$, i.e., 
$b(\hat{\mathbf{g}}(\mathbf{y}); \mathbf{x}) = b(\hat{\mathbf{g}}'(\mathbf{z}); \mathbf{x})$ for every 
$\mathbf{x} \in \mathcal{X}$, and $v(\hat{\mathbf{g}}(\mathbf{y});\mathbf{x}_{0}) = v(\hat{\mathbf{g}}'(\mathbf{z}); \mathbf{x}_{0})$. 
\end{proof}


\subsection[Non-invertible Transformation]{Non-invertible Transformation -- Data Processing Inequality for Classical Estimation} 
\label{sec_data_proc_inequ_classical_est}

As discussed in the previous section, replacing the observation $\mathbf{y}$ by $\mathbf{z} = \mathbf{T}(\mathbf{y})$ where $\mathbf{T}(\cdot): \mathbb{R}^{M} \rightarrow \mathbb{R}^{K}$ 
is an arbitrary (but known) invertible deterministic function is irrelevant w.r.t.\ minimum variance and minimax estimation. Loosely speaking, an invertible transformation is just a ``relabeling'' of the observed data, i.e., the realization of the observation, which does not 
incur any gain or loss of information. 

However, if the transformation $\mathbf{T}(\cdot)$ applied to the observation is not invertible or even random, the situation changes. While in general the achievable performance for minimum variance and minimax estimators 
is different for the transformed estimation problem $\mathcal{E}'$ whose statistical model is given by $f( \mathbf{z} ; \mathbf{x})$, we still have a relation between the modified and the original estimation problem. 
This relation is due to the fact that any estimator $\hat{\mathbf{g}}(\cdot)$ that uses the new observation $\mathbf{z}$ is 
also an (possibly random) estimator for the original estimation problem with the observation $\mathbf{y}$ since $\hat{\mathbf{g}}(\mathbf{z}) = \hat{\mathbf{g}}( \mathbf{T}(\mathbf{y}))$. 
More precisely, we have 
\begin{theorem}
\label{thm_data_proc_inequ_classical_est}
Consider an estimation problem $\mathcal{E} = \vecestproblem$ with a statistical model $f(\mathbf{y}; \mathbf{x})$ and a (possibly random) function $\mathbf{T}(\cdot): \mathbb{R}^{M} \rightarrow \mathbb{R}^{K}$. 
If the function $\mathbf{T}(\cdot)$ is random we assume that it is statistically independent of $\mathbf{y}$ for every $\mathbf{x} \in \mathcal{X}$. 
We denote by $\mathcal{E}'$ the estimation problem that is obtained from $\mathcal{E}$ by using the vector $\mathbf{z} =Ê\mathbf{T}(\mathbf{y})$ as the observation, i.e., 
$\mathcal{E}' = \left( \mathcal{X}, f(\mathbf{z} ; \mathbf{x}), \mathbf{g}(\cdot) \right)$ where $f(\mathbf{z};\mathbf{x})$ denotes the pdf of the new observation $\mathbf{z} = \mathbf{T}(\mathbf{y})$.
We then have that for any two minimum variance problems given by $\mathcal{M} = \minvarproblem$, $\mathcal{M}' = \left( \mathcal{E}', \mathbf{c}(\cdot), \mathbf{x}_{0}) \right)$ which 
share a common prescribed bias function $\mathbf{c}(\cdot): \mathcal{X} \rightarrow \mathbb{R}^{P}$ and a fixed parameter vector $\mathbf{x}_{0} \in \mathcal{X}$, 
\begin{equation}
\label{equ_data_proc_inequ_classic_est_min_var}
L_{\mathcal{M}} \leq L_{\mathcal{M}'}. 
\end{equation}  
Moreover, the minimax risks of $\mathcal{E}$ and $\mathcal{E}'$ satisfy 
\begin{equation} 
R_{\mathcal{E}} \leq R_{\mathcal{E}'}.
\end{equation} 
\end{theorem}
\begin{proof}
The statement follows from the fact that for any estimator $\hat{\mathbf{g}}'(\mathbf{z})$ for $\mathcal{E}'$ and $\mathcal{M}'$ 
that uses the observation $\mathbf{z}=\mathbf{T}(\mathbf{y})$, we have that the specific (possibly random) estimator $\hat{\mathbf{g}}(\mathbf{y}) \triangleq \hat{\mathbf{g}}'(\mathbf{T}(\mathbf{y}))$ for 
$\mathcal{E}$ and $\mathcal{M}$ that uses the observation $\mathbf{y}$ has the same bias and variance at any $\mathbf{x} \in \mathcal{X}$. 

In particular, assume that $L_{\mathcal{M}'} < L_{\mathcal{M}}$ which would imply, via Definition \ref{def_min_ach_var}, that there must be an allowed estimator $\hat{\mathbf{g}}'_{0}(\mathbf{z})$ for $\mathcal{M}'$ whose 
variance at $\mathbf{x}_{0}$ is strictly smaller than $L_{\mathcal{M}}$, i.e., $v(\hat{\mathbf{g}}'(\mathbf{z});\mathbf{x}_{0}) < L_{\mathcal{M}}$. We can then construct an allowed estimator for $\mathcal{M}$ as 
$\hat{\mathbf{g}}_{0}(\mathbf{y}) \triangleq \hat{\mathbf{g}}'_{0}(\mathbf{T}(\mathbf{y}))$ which has the same variance $v(\hat{\mathbf{g}}'(\mathbf{z});\mathbf{x}_{0})< L_{\mathcal{M}}$ at $\mathbf{x}_{0}$. This, however, is 
a contradiction due to the definition of $L_{\mathcal{M}}$. A similar argument holds for the minimax risks of the estimation problems $\mathcal{E}$ and $\mathcal{E}'$.  

Thus the infimum in \eqref{equ_def_min_ach_var} and \eqref{equ_opt_minimax} cannot decrease by moving from $\mathcal{E}$ to $\mathcal{E}'$ or from $\mathcal{M}$ to $\mathcal{M}'$, respectively.
\end{proof}

\section{Reducing the Parameter Set}
\label{sec_effect_red_par_set}

The following definition will prove handy: 
\begin{definition}
\label{def_min_var_problem_reduc_par_set}
Consider a minimum variance problem $\mathcal{M}=\minvarproblem$ associated with the estimation problem $\mathcal{E} = \vecestproblem$ 
and an arbitrary subset $\mathcal{X}' \subseteq \mathcal{X}$ of the parameter set of $\mathcal{M}$. 
We then denote by 
\begin{equation}
\mathcal{M}\big|_{\mathcal{X}'} \triangleq \bigg( \mathcal{E}', \mathbf{c}(\cdot) \big|_{\mathcal{X}'}, \mathbf{x}_{0} \bigg)
\end{equation}
the specific minimum variance problem associated with the estimation problem 
\begin{equation}
\mathcal{E}' \triangleq \bigg( \mathcal{X}', f(\mathbf{y}; \mathbf{x}) , \mathbf{g}(\cdot)\big|_{\mathcal{X}'} \bigg)
\end{equation}
which is identical to $\mathcal{M}$ except for the parameter set. 
\end{definition}

Given an estimation problem $\mathcal{E}=\vecestproblem$, it may be known beforehand that the true parameter vector $\mathbf{x}$ 
belongs to a subset $\mathcal{X}' \subseteq \mathcal{X}$ of the ``nominal'' parameter set $\mathcal{X}$. 
Intuitively, this a priori knowledge should be exploitable in order to improve the performance of estimators for $\mathbf{g}(\mathbf{x})$. 
A precise statement in this regard for minimax estimation is  
\begin{theorem} 
Given two estimation problems $\mathcal{E}=\vecestproblem$, $\mathcal{E}'=\left( \mathcal{X}', f(\mathbf{y};\mathbf{x}),\mathbf{g}(\cdot) \right)$ which differ only in their parameter sets $\mathcal{X}$ and $\mathcal{X}'$, we have 
\begin{equation} 
\label{equ_reduc_para_set_minimax_risk_gets_smaller}
\mathcal{X}' \subseteq \mathcal{X} \Rightarrow R_{\mathcal{E}'} \leq R_{\mathcal{E}}, 
\end{equation} 
i.e., reducing the parameter set can never increase the minimax risk of an estimation problem. 
\end{theorem}
\begin{proof}
For any estimator $\hat{\mathbf{g}}(\cdot)$, we have due to $\mathcal{X}' \subseteq \mathcal{X}$ the following relation for the worst MSE's associated with $\mathcal{E}$ and $\mathcal{E}'$: 
\begin{equation}
R_{\mathcal{E}'}(\hat{\mathbf{g}}(\cdot)) =  \sup_{\mathbf{x} \in \mathcal{X}'} \varepsilon(\hat{\mathbf{g}}(\cdot), \mathbf{x}) \leq  \sup_{\mathbf{x} \in \mathcal{X}} \varepsilon(\hat{\mathbf{g}}(\cdot), \mathbf{x})=R_{\mathcal{E}} (\hat{\mathbf{g}}(\cdot)).
\end{equation}  
Since this inequality holds for every estimator and since $R_{\mathcal{E}}(\hat{\mathbf{g}}(\cdot)) \triangleq \sup_{\mathbf{x} \in \mathcal{X}} \varepsilon(\hat{\mathbf{g}}(\cdot); \mathbf{x})$, we obtain \eqref{equ_reduc_para_set_minimax_risk_gets_smaller}.
\end{proof} 

Compared to minimax estimation, the effect of reducing the parameter set of an estimation problem in the context of minimum variance estimation is more subtle \cite{ZvikaCRB,ZvikaSSP,StoicaNgCCRB,ZvikaCCRB,MooreCCRB,GormanHero}. 
As one may intuitively expect, the achievable performance as quantified by the minimum achievable variance cannot become worse. However, this improved performance bound may not be achieved by incorporating 
the a-priori information of a reduced parameter set directly into the design of an estimator. E.g., if we are given a minimum variance problem $\mathcal{M} = \minvarproblem$ 
associated with the estimation problem  $\mathcal{E}=\vecestproblem$ and it is known beforehand that the true parameter vector $\mathbf{x}$ 
is contained within $\mathcal{X}' \subseteq \mathcal{X}$, then it might be reasonable to constrain an estimator $\hat{\mathbf{g}}(\mathbf{y})$ for $\mathcal{M}$ 
to take on values only within the set $\mathbf{g}(\mathcal{X}') \subseteq \mathbb{R}^{P}$. However, for minimum variance estimation it may not be optimal to place such an additional constraint on an estimator. 
The real effect of a reduced parameter set for minimum variance estimation is that the set of allowed estimators $\mathcal{F}(\mathcal{M})$ becomes larger 
if the parameter set $\mathcal{X}$ is reduced to $\mathcal{X}'$ since the bias constraint defining $\mathcal{F}(\mathcal{M})$ (see \eqref{equ_est_finite_var_prescr_bias}) 
has to be satisfied only on the smaller set $\mathcal{X}'$ instead of the original parameter set $\mathcal{X}$. 
Since the set of allowed estimators can only become larger, the minimum achievable variance (cf. \eqref{equ_def_min_ach_var}) can never increase, as is stated in 
\begin{theorem} 
\label{thm_par_set_reduction_classic_est_mve}
Consider a minimum variance problem $\mathcal{M}=\minvarproblem$ and an arbitrary subset $\mathcal{X}' \subseteq \mathcal{X}$ of the parameter set of $\mathcal{M}$. 
For any modified minimum variance problem $\mathcal{M}' = \mathcal{M}\big|_{\mathcal{X}'}$, associated with the smaller parameter set $\mathcal{X}' \subseteq \mathcal{X}$, we have
\begin{equation} 
L_{\mathcal{M}'} \leq L_{\mathcal{M}}, 
\end{equation}
i.e., a reduction of the parameter set can never result in an increase of the minimum achievable variance. 

Furthermore, if we consider a minimum variance problem $\mathcal{M}$ with zero bias, $\mathbf{c}(\cdot) \equiv 0$, and an arbitrary parameter function $\mathbf{g}(\cdot): \mathcal{X} \rightarrow \mathbb{R}^{P}$, we 
have that if the parameter function $\mathbf{g}(\cdot)\big|_{\mathcal{X}'}$ is not estimable for $\mathcal{M}'$, then the parameter function $\mathbf{g}(\cdot)$ is also not estimable for $\mathcal{M}$. 
\end{theorem} 
\begin{proof} 
This result follows from \eqref{equ_def_min_ach_var} by the obvious fact that under the condition 
$\mathcal{X}' \subseteq \mathcal{X}$ we have that $\mathcal{F}(\mathcal{M}) \subseteq \mathcal{F}(\mathcal{M}')$, i.e., any allowed estimator for $\mathcal{M}$ is necessarily also 
an allowed estimator for $\mathcal{M}'$. 
\end{proof}

\section{Statistical Models Belonging to an Exponential Family} 
\label{sec_exp_family} 
An important class of estimation problems is obtained by assuming that the statistical model $f(\mathbf{y}; \mathbf{x})$ belongs to an exponential family  \cite{LC,FundmentExpFamBrown,GraphModExpFamVarInfWainJor}. 
An exponential family is a class of pdfs that is associated with two fixed vector-valued functions ${\bm \Phi}(\cdot): \mathbb{R}^{M} \rightarrow \mathbb{R}^{K}$, $\mathbf{u}(\cdot): \mathbb{R}^{N} \rightarrow \mathbb{R}^{K}$ of the observation $\mathbf{y}$ 
and the parameter vector $\mathbf{x}$, respectively: 
\begin{equation} 
\label{equ_def_pdf_exp_fam}
f^{( {\bm \Phi},\mathbf{u})}(\mathbf{y}; \mathbf{x}) = \exp \left( \big[ {\bm \Phi}(\mathbf{y}) \big]^{T} \mathbf{u}(\mathbf{x}) - A^{( {\bm \Phi})}(\mathbf{x}) \right) h(\mathbf{y})
\end{equation} 
where $h(\cdot): \mathbb{R}^{M} \rightarrow \mathbb{R}$ is a fixed nonnegative weight function. 
Associated with an exponential family is the \emph{natural parameter space} \cite{FundmentExpFamBrown} denoted $N^{( {\bm \Phi})}$ and defined as 
\begin{equation}
 N^{( {\bm \Phi})} \triangleq \left\{ \mathbf{x} \in \mathbb{R}^{N} \bigg| \int_{\mathbf{y}}  \exp \left( \big[ {\bm \Phi}(\mathbf{y}) \big]^{T}  \mathbf{u}(\mathbf{x}) \right) h(\mathbf{y})d\mathbf{y} < \infty \right\}. 
\end{equation}
The function $A^{( {\bm \Phi})}(\mathbf{x}):  N^{( {\bm \Phi})}  \rightarrow \mathbb{R}$ is defined by the requirement that $\int_{\mathbf{y}}f^{( {\bm \Phi},\mathbf{u})}(\mathbf{y}; \mathbf{x}) d\mathbf{y} = 1$ for every $\mathbf{x} \in N^{( {\bm \Phi})}$, i.e., 
\begin{equation} 
A^{( {\bm \Phi})}(\mathbf{x}) \triangleq \log \int_{\mathbf{y}}\exp \left( \big[ {\bm \Phi}(\mathbf{y}) \big]^{T}  \mathbf{u}(\mathbf{x})  \right) h(\mathbf{y}) d\mathbf{y}.
\end{equation}  
In the literature, the function ${\bm \Phi}(\cdot)$ is known as the \emph{potential function} or \emph{sufficient statistic} and the function $A^{( {\bm \Phi})}(\cdot)$ is known as the \emph{log partition
function} or \emph{cumulant (generating) function} \cite{GraphModExpFamVarInfWainJor,FundmentExpFamBrown}.  

It can be shown \cite{FundmentExpFamBrown} that under some weak technical requirements on the parameter function $\mathbf{u}(\cdot)$, every exponential family can be reduced to 
a standard form, with parameter function $\mathbf{u}'(\mathbf{x}) = \mathbf{x}$. In particular, any estimation problem with parameter set $\mathcal{X}$ and statistical model 
$f^{( {\bm \Phi},\mathbf{u})}(\mathbf{y}; \mathbf{x})$ is equivalent to an estimation problem with parameter set $\mathcal{X}' = \mathbf{u}(\mathcal{X})$ and statistical model 
$f^{( {\bm \Phi},\mathbf{u}'(\mathbf{x}) = \mathbf{x})}(\mathbf{y}; \mathbf{x})$.

By definition, every exponential family in standard form, i.e., $\mathbf{u}(\mathbf{x}) = \mathbf{x})$, is completely characterized by the log partition function $A^{( {\bm \Phi})}(\cdot):  N^{( {\bm \Phi})}  \rightarrow \mathbb{R}$ which 
has some interesting properties. In particular, for an exponential family in standard form, according to \cite{FundmentExpFamBrown,GraphModExpFamVarInfWainJor} for every $\mathbf{x}_{0}$ in the interior of $N^{( {\bm \Phi})}$, 
i.e., there exists a radius $r >0$ such that $\mathcal{B}(\mathbf{x}_{0},r) \subseteq N^{( {\bm \Phi})}$, the partial derivatives 
$\frac{ \partial^{\mathbf{p}} A^{( {\bm \Phi})}(\mathbf{x})}{\partial \mathbf{x}^{\mathbf{p}}}\bigg|_{\mathbf{x} = \mathbf{x}_{0}}$ up to any order $\mathbf{p} \in \mathbb{Z}_{+}^{N}$ exist and 
are given by the cumulants of the pdf $f^{( {\bm \Phi})}(\mathbf{y}; \mathbf{x}_{0})$. In particular, we obtain for the first order derivatives
\begin{equation}
 \frac{ \partial^{\mathbf{e}_{l}} A^{( {\bm \Phi})}(\mathbf{x})}{\partial \mathbf{x}^{\mathbf{e}_{l}}}\bigg|_{\mathbf{x} = \mathbf{x}_{0}} = \mathsf{E}_{\mathbf{x}_{0}} \{ \Phi_{l}(\mathbf{y})  \}. 
\end{equation}

Another useful fact is that the partial derivatives $\frac{\partial^{\mathbf{p}} f(\mathbf{y}; \mathbf{x})}{\partial \mathbf{x}^{\mathbf{p}}}\big|_{\mathbf{x} = \mathbf{x}_{0}}$ exist for any order $\mathbf{p} \in \mathbb{Z}_{+}^{N}$ at every 
$\mathbf{x}_{0}$ in the interior of $N^{( {\bm \Phi})}$. Moreover, it holds that  \cite{FundmentExpFamBrown,GraphModExpFamVarInfWainJor}
\begin{equation}
\mathsf{E}_{\mathbf{x}_{0}} \bigg\{ \left( \frac{1}{f(\mathbf{y};\mathbf{x})} \frac{\partial^{\mathbf{p}} f(\mathbf{y}; \mathbf{x})}{\partial \mathbf{x}^{\mathbf{p}}}Ê\right)^{2} \bigg\}  < \infty
\end{equation}
for any multi-index $\mathbf{p} \in \mathbb{Z}_{+}^{N}$. 

Maybe the most important instance of a pdf belonging to an exponential family is the multivariate normal distribution $f(\mathbf{y})$ of a Gaussian random vector $\mathbf{y} \in \mathbb{R}^{M}$ 
with mean ${\bm \mu} \in \mathbb{R}^{M}$ and the positive definite covariance matrix $\mathbf{C}\in \mathbb{R}^{M \times M}$, i.e., 
$f(\mathbf{y}) = \frac{1}{\sqrt{ (2 \pi)^{M} \detm{Ê\mathbf{C}}}} \exp \left( - \frac{1}{2} (\mathbf{y} - {\bm \mu})^{T} \mathbf{C}^{-1} (\mathbf{y} - {\bm \mu}) \right)$. 
This distribution is obtained for the choice $\mathbf{x} = {\bm \mu}$ from the exponential family in standard form, i.e., $\mathbf{u}(\mathbf{x}) = \mathbf{x}$, with sufficient statistic ${\bm \Phi}(\mathbf{y}) = \mathbf{C}^{-1} \mathbf{y}$, weight function
$h(\mathbf{y}) =  \frac{1}{\sqrt{ (2 \pi)^{M} \detm{Ê\mathbf{C}}}} \exp\big( - \frac{1}{2} \mathbf{y}^{T} \mathbf{C}^{-1}\mathbf{y}\big)$ and cumulant function $A^{( {\bm \Phi})}(\mathbf{x})= \frac{1}{2} \mathbf{x}^{T} \mathbf{C}^{-1}\mathbf{x}$. 
If a random vector $\mathbf{y} \in \mathbb{R}^{M}$ is distributed normally with mean ${\bm \mu}$ and covariance matrix $\mathbf{C} \in \mathbb{R}^{M \times M}$ we denote this by $\mathbf{y} \sim \mathcal{N}({\bm \mu}, \mathbf{C})$. 

Finally, we note that the concept of exponential families provides a natural link between classical estimation theory and the theory of graphical models. Indeed, the structure of the sufficient statistics translates directly 
to a specific connectivity structure of the associated graphical model \cite{GraphModExpFamVarInfWainJor,CovEstDecompGraphModel}.

\chapter{Mathematical Tools}
\label{chap_RKHS_Fund}
This chapter is a self-contained presentation of the mathematical tools that will be required later on. 
All the mathematical concepts discussed in this chapter are based on sets of real-valued functions $f(\cdot): \mathcal{D} \rightarrow \mathbb{R}$ which share a common domain $\mathcal{D}$. 
We define the multiplication of a function with a real number and the sum of two functions over the same domain in the usual pointwise sense. The following presentation is based 
mainly on \cite{RudinBook, RudinBookPrinciplesMatheAnalysis} for general Hilbert space theory and \cite{Parzen59,aronszajn1950} for the theory of RKHS. 
 
\section{Hilbert Spaces}

\subsection{Some Basic Definitions} 

\begin{definition}
An inner-product function space $\mathcal{I}$ over a domain $\mathcal{D}$ is a linear space of real-valued functions $f(\cdot): \mathcal{D} \rightarrow \mathbb{R}$. By a linear space we mean that 
\begin{equation} 
f(\cdot) , g(\cdot) \in \mathcal{I}\,\mbox{ and } \, a, b \in \mathbb{R} \quad  \Rightarrow  \quad a f(\cdot) + b g(\cdot) \in \mathcal{I},
\end{equation}
i.e., any linear combination of two elements of $\mathcal{I}$ is also an element of $\mathcal{I}$. 
Moreover, an inner product is defined on $\mathcal{I}$, i.e., a mapping that associates to two functions $f(\cdot),g(\cdot) \in \mathcal{I}$ the real number $\langle f(\cdot), g(\cdot) \rangle_{\mathcal{I}} \in \mathbb{R}$. 
This inner product has the properties 
\begin{align}
(a) \,\,\,\, &  \langle f(\cdot), g(\cdot) \rangle_{\mathcal{I}}= \langle g(\cdot), f(\cdot) \rangle_{\mathcal{I}} \nonumber\\ 
(b) \,\,\,\, &  \langle a f(\cdot)+ bh(\cdot), g(\cdot) \rangle_{\mathcal{I}} = a\langle f(\cdot), g(\cdot) \rangle_{\mathcal{I}} + b\langle h(\cdot), g(\cdot) \rangle_{\mathcal{I}}  \quad\quad \forall a,b \in \mathbb{R} \mbox{ and }  \forall f(\cdot),h(\cdot),g(\cdot)Ê\in \mathcal{I} \nonumber \\ 
(c) \,\,\,\, &  \langle f(\cdot), f(\cdot) \rangle_{\mathcal{I}} \geq 0 \quad \forall f(\cdot) \in \mathcal{I} \nonumber \\
(d) \,\,\,\, &  \langle f(\cdot), f(\cdot) \rangle_{\mathcal{I}} = 0 \quad \mbox{ only if } \quad f(\cdot) \equiv 0.   \label{equ_cond_inner_product}
\end{align}
\end{definition} 

In what follows we will need the famous \emph{Cauchy-Schwarz inequality} \cite{RudinBook} as stated in 
 \begin{theorem}
 \label{thm_cauchy_schwarz}
 Consider an function space $\mathcal{I}'$ together with a linear mapping that assigns to any two functions $f(\cdot) ,g(\cdot) \in \mathcal{I}'$ the real number $\langle f(\cdot), g(\cdot) \rangle_{\mathcal{I}'}$.
 If the mapping is such that $\langle f(\cdot), g(\cdot) \rangle_{\mathcal{I}'}$ has the properties $(a)$-$(c)$ of \eqref{equ_cond_inner_product}, then we have that 
 \begin{equation}
 \label{equ_cauchy_schwarz}
 | \langle f(\cdot), g(\cdot) \rangle_{\mathcal{I}'}| \leq \sqrt{\langle f(\cdot), f(\cdot) \rangle_{\mathcal{I}'}  \langle g(\cdot),  g(\cdot) \rangle_{\mathcal{I}'}},
 \end{equation} 
 where the equality sign holds if and only if there is a number $c \in \mathbb{R}$ such that $\langle f(\cdot)-cg(\cdot), f(\cdot) - cg(\cdot) \rangle_{\mathcal{I}'} =0$. 
 \end{theorem}
 \begin{proof} 
 \cite[p.\ 77]{RudinBook}
 \end{proof} 
 
For every inner product space $\mathcal{I}$, one can naturally define an \emph{induced norm} $\| f(\cdot) \|_{\mathcal{I}}$ via 
\begin{equation}
 \| f(\cdot) \|_{\mathcal{I}} \triangleq \sqrt{ \langle f(\cdot), f(\cdot) \rangle_{\mathcal{I}}}.
\end{equation}
As can be verified easily, the induced norm satisfies the axioms of a norm, i.e., 
\begin{itemize}
 \item Positive definiteness:
 \begin{equation} 
 \label{equ_pos_definite_norm}
 \forall f(\cdot) \in \mathcal{H}: \quad \| f(\cdot) \|_{\mathcal{I}} \geq 0, \quad \quad  \| f(\cdot) \|_{\mathcal{I}} = 0 \,\,\,\, \Leftrightarrow \,\,\,\, f(\cdot) \equiv 0.
 \end{equation}
 \item Homogeneity:
 \begin{equation} 
 a \in \mathbb{R}, f(\cdot) \in \mathcal{I} \Rightarrow \| a f(\cdot) \|_{\mathcal{I}} = |a|  \| f(\cdot) \|_{\mathcal{I}}.
 \end{equation}  
 \item Triangle inequality: 
 \begin{equation} 
 \label{equ_triangle_inequality}
 \|f(\cdot) + g(\cdot)\|_{\mathcal{I}} \leq \|f(\cdot)\|_{\mathcal{I}} + \|g(\cdot)\|_{\mathcal{I}}.
 \end{equation}
\end{itemize}
It may happen that for a given linear function space $\mathcal{I}'$ one can define naturally a real-valued function $\langle f(\cdot), g(\cdot) \rangle_{\mathcal{I}'}$ of pairs of elements of the function space, 
that satisfies all requirements for an inner product in \eqref{equ_cond_inner_product} except condition $(d)$.\footnote{Maybe the best known example of this situation is given by the set of square-integrable real-valued functions.} 
It is then customary (cf. \cite{RudinBook}) to consider two functions $f(\cdot),g(\cdot) \in \mathcal{I}'$ as identical 
if  $\langle f(\cdot) - g(\cdot), f(\cdot) - g(\cdot) \rangle_{\mathcal{I}'} = 0$. Technically speaking, one constructs an inner-product space whose elements are the equivalence classes 
\begin{equation} 
\big[f(\cdot)\big] \triangleq \bigg\{ f'(\cdot) \in \mathcal{I}' \big|  \big\langle f(\cdot) - f'(\cdot), f(\cdot) - f'(\cdot) \big\rangle_{\mathcal{I}'} = 0 \bigg\} 
\end{equation} 
where $f(\cdot) \in \mathcal{I}'$. 
Based on Theorem \ref{thm_cauchy_schwarz}, one can then verify that this inner-product space is well defined, i.e., if $f'(\cdot)\in [f(\cdot)]$ and $g'(\cdot) \in [g(\cdot)]$ then 
$\langle f(\cdot) , g(\cdot) \rangle_{\mathcal{I}'}=\langle f'(\cdot), g'(\cdot) \rangle_{\mathcal{I}'}$. 
However, we will follow the usual convention of considering the space $\mathcal{I}'$ itself as an inner-product space \cite{RudinBook} and denote by a function $f(\cdot)$ not only the function itself but 
the whole equivalence class $[f(\cdot)]$. 

Therefore, if we write $f(\cdot) \equiv 0$, we mean not necessarily that $f(\mathbf{x}) = 0$ for every $\mathbf{x} \in \mathcal{D}$ but 
that the function $f(\cdot)$ is in the same class as the function which is identically zero. However, within this 
thesis we will work primarily with a special kind of Hilbert spaces, i.e., a RKHS, and in that case 
we have that every class $\big[ f(\cdot) \big]$ consists of one and only one function, i.e., $\big[ f(\cdot) \big] = \{ f(\cdot) \}$. 

The induced norm $\| \cdot \|_{\mathcal{I}}$ also defines the notion of convergence in an inner-product function space $\mathcal{I}$. 
\begin{definition}
A sequence $\{ f_{l} (\cdot) \}_{l \rightarrow \infty}$ of functions $f_{l}(\cdot)$ belonging to an inner-product function space $\mathcal{I}$, i.e., $f_{l}(\cdot) \in \mathcal{I}$, is called 
convergent to a function $f(\cdot) \in \mathcal{I}$ if to any $\varepsilon > 0$ there exists a $n_{0} \in \mathbb{N}$ such that 
\begin{equation}
n \geq n_{0} \,\, \Rightarrow \,\, \| f_{n}(\cdot) - f(\cdot) \|_{\mathcal{I}} \leq \varepsilon.
\end{equation}
The function $f(\cdot)$ is then called the limit of the sequence $\{ f_{l} (\cdot) \}_{l \rightarrow \infty}$ and denoted by $\lim_{l \rightarrow \infty} f_{l}(\cdot)$. 
\end{definition} 

A necessary condition for a sequence to be convergent is that it is a Cauchy sequence: 
\begin{definition}
A sequence $\{ f_{l} (\cdot) \in \mathcal{I} \}_{l \rightarrow \infty}$ of functions $f_{l}(\cdot)$ belonging to an inner-product function space $\mathcal{I}$ is called a Cauchy sequence if to any 
$\varepsilon > 0$ there exists a $n_{0} \in \mathbb{N}$ such that 
\begin{equation}
n_{1},n_{2} \geq n_{0} \,\, \Rightarrow \,\,  \| f_{n_{1}}(\cdot) - f_{n_{2}}(\cdot) \|_{\mathcal{I}} \leq \varepsilon.
\end{equation}
\end{definition} 

Now we are in the position to define the central mathematical structure for our thesis, i.e., a Hilbert space: 
\begin{definition}
A function Hilbert space $\mathcal{H}$ is an inner-product function space that has the additional property of being complete, i.e., 
every Cauchy sequence of functions belonging to $\mathcal{H}$ converges to a function belonging to $\mathcal{H}$. 
\end{definition} 
Every inner-product space $\mathcal{I}$ can be completed in a canonical way such that it becomes a Hilbert space \cite{aronszajn1950}.
In what follows, we will state the most important concepts and facts about Hilbert spaces, as required for our purposes.

\begin{definition}
Given a function Hilbert space $\mathcal{H}$, and an index set $\mathcal{T}$, we call a subset $\mathcal{S} = \{g_{l} (\cdot)  \}_{l \in \mathcal{T}} \subseteq \mathcal{H}$ dense, if any element of $\mathcal{H}$ can be approximated 
arbitrary well by an element of the set $\mathcal{S}$, i.e., for any $f(\cdot) \in \mathcal{H}$ and any $\varepsilon > 0$, there exists a function $g(\cdot) \in \mathcal{S}$ such that 
\begin{equation}
 \| g(\cdot) - f(\cdot) \|_{\mathcal{H}} \leq \varepsilon.
\end{equation}
\end{definition}

\begin{definition}
A function Hilbert space $\mathcal{H}$ is said to be separable if it contains a countable subset which is dense in $\mathcal{H}$. 
\end{definition}
Within this thesis, we exclusively work with Hilbert spaces that are separable. 

\begin{definition}
\label{def_complete_set}
Given a function Hilbert space $\mathcal{H}$, a set of functions $ \mathcal{S} = \{g_{l} (\cdot)  \}_{l \in \mathcal{T}}$ with an arbitrary index set $\mathcal{T}$ is called complete if the following implication holds:
\begin{equation}
\big\langle g_{l}(\cdot), f(\cdot)  \big\rangle_{\mathcal{H}} = 0 \,\,\, \forall l \in \mathcal{T} \quad \Rightarrow \quad f (\cdot) \equiv 0. 
\end{equation}
\end{definition} 

\begin{definition}
Consider a function Hilbert space $\mathcal{H}$ and a subset $\mathcal{U} \subseteq \mathcal{H}$. If 
the subset is such that any linear combination of two elements of $\mathcal{U}$ is still in $\mathcal{U}$, i.e., 
\begin{equation} 
\label{equ_def_subspace_linear_com}
f(\cdot) , g(\cdot) \in \mathcal{U}\,\mbox{ and }\, a, b \in \mathbb{R} \quad  \Rightarrow  \quad a f(\cdot) + b g(\cdot) \in \mathcal{U},
\end{equation} 
and if moreover every Cauchy sequence $\{ f_{l}(\cdot) \in \mathcal{U} \}_{l \rightarrow \infty}$ converges 
to a function belonging to $\mathcal{U}$, then we call $\,\,\mathcal{U}$ a subspace of the Hilbert space $\mathcal{H}$. 
\end{definition}
Note that according to this definition, we have that any subspace of a Hilbert space is itself a Hilbert space. Indeed, a subspace is a linear space since by \eqref{equ_def_subspace_linear_com} any 
linear combination of two elements are again in the subspace. A subspace is also an inner-product space since an inner product is given naturally by the inner product defined on the larger Hilbert space, in which 
the subspace is contained. Finally, by definition, a subspace is an inner-product space in which every Cauchy converges to a function of the same subspace, i.e., it is a complete inner-product space, i.e., a Hilbert space.

\begin{definition}
\label{def_relative_closure}
Given a function Hilbert space $\mathcal{H}$ that contains an inner-product space $\mathcal{I}$ (which is defined for the same inner product as that of $\mathcal{H}$), i.e., $\mathcal{I} \subseteq \mathcal{H}$, we define 
the closure of $\mathcal{I}$ relative to $\mathcal{H}$, denoted by $\closure \{ \mathcal{I} \}$, as the subspace of $\mathcal{H}$ that consists of all functions belonging to $\mathcal{I}$ as well as 
any limit of a Cauchy sequence of functions belonging to $\mathcal{I}$. 
\end{definition} 

\begin{definition}
Given an inner-product function space $\mathcal{I}$ (in particular, a Hilbert space), an arbitrary index set $\mathcal{T}$, and a set $\{ v_{l}(\cdot) \}_{l \in \mathcal{T}}$ of functions $v_{l}(\cdot)\in \mathcal{I}$, we define 
the linear span of  $\{ v_{l}(\cdot) \}_{l \in \mathcal{T}}$, denoted by $\linspan  \{ v_{l}(\cdot) \}_{l \in \mathcal{T}}$, as the set of all finite linear combinations of the functions $v_{l}(\cdot)$, i.e., 
\begin{equation}
\linspan  \{ v_{l}(\cdot) \}_{l \in \mathcal{T}} \triangleqÊ\bigg\{ v(\cdot) \in \mathcal{I} \bigg| v(\cdot) = \sum_{j \in [J]} a_{j} v_{l_j}(\cdot) \,\,\,  \mbox{ with }J \in \mathbb{N}, a_{j} \in \mathbb{R}, l_{j} \in \mathcal{T} \bigg\}.
\end{equation} 
\end{definition}

\begin{definition}
Given a function Hilbert space $\mathcal{H}$ and an arbitrary index set $\mathcal{T}$, we say that the set of functions $\{ f_{l} (\cdot) \in \mathcal{H} \}_{l \in \mathcal{T}}$ spans
$\mathcal{H}$ if the linear span $\linspan \{ f_{l}(\cdot) \}_{l \in \mathcal{T}}$ is dense in the Hilbert space $\mathcal{H}$. 
\end{definition}

\begin{definition}
\label{def_ONB}
Given a function Hilbert space $\mathcal{H}$, we call any set of functions $\{Êg_{l} (\cdot) \in \mathcal{H} \}_{l \in \mathcal{T}}$ an orthonormal basis (ONB) for the Hilbert space $\mathcal{H}$ if 
the functions  $\{Êg_{l} (\cdot) \}_{l \in \mathcal{T}}$ span the Hilbert space $\mathcal{H}$ and the functions $g_{l}(\cdot)$ are orthonormal, i.e., 
\begin{equation} 
\langle g_{k}(\cdot) , g_{l}(\cdot) \rangle_{\mathcal{H}} = \delta_{k,l}.
\end{equation} 
\end{definition} 


\begin{definition}
\label{def_orthog_projection}
Given a function Hilbert space $\mathcal{H}$ and a subspace $\mathcal{U} \subseteq \mathcal{H}$, we define the 
orthogonal projection of a function $f(\cdot) \in \mathcal{H}$ onto the subspace $\mathcal{U}$, denoted by $\mathbf{P}_{\mathcal{U}} f(\cdot)$, as the element of $\,\,\mathcal{U}$ given by 
\begin{equation}
\mathbf{P}_{\mathcal{U}} f(\cdot) \triangleq  \arg\!\min_{\hspace*{-3mm} g(\cdot) \in \mathcal{U}}    \| g(\cdot) - f(\cdot) \|_{\mathcal{H}}.
\end{equation}
\end{definition} 
From its definition, it is clear that the orthogonal projection $\mathbf{P}_{\mathcal{U}} f(\cdot)$ is the best approximation of $f(\cdot)$ by an element of $\mathcal{U}$. 

\begin{definition}
\label{def_isometry_congruence}
We call two Hilbert spaces $\mathcal{H}_{1}$ and $\mathcal{H}_{2}$ with associated inner products $\langle \, \cdot \, , \, \cdot \, \rangle_{\mathcal{H}_{1}}$ and $\langle \, \cdot \, , \, \cdot \, \rangle_{\mathcal{H}_{2}}$  
isometric if there exists a bijective linear map $\mathsf{J}[\cdot]: \mathcal{H}_{1} \rightarrow \mathcal{H}_{2}$, i.e., 
\begin{equation} 
f(\cdot),g(\cdot) \in \mathcal{H}_{1} \,\mbox{ and } \, a,b \in \mathbb{R} \,  \Rightarrow \, \mathsf{J}[af(\cdot) + bg(\cdot)]= a\mathsf{J}[f(\cdot)] + b\mathsf{J}[g(\cdot)],
\end{equation}
which preserves inner products, i.e., 
\begin{equation}
\label{equ_conser_inner_prod_def_isometry}
f(\cdot),g(\cdot) \in \mathcal{H}_{1} \, \Rightarrow \, \langle f(\cdot),g(\cdot) \rangle_{\mathcal{H}_{1}} = \langle \mathsf{J}[f(\cdot)], \mathsf{J}[g(\cdot)] \rangle_{\mathcal{H}_{2}}.
\end{equation}
Any bijective linear map $\mathsf{J}[\cdot]$ for which \eqref{equ_conser_inner_prod_def_isometry} holds is called a congruence (or congruence map) from $\mathcal{H}_{1}$ to $\mathcal{H}_{2}$. 
\end{definition}

\subsection{Some Basic Facts}
 
\begin{theorem}
\label{thm_approx_norm_dense_set}
Consider a function Hilbert space $\mathcal{H}$ and a subset $\mathcal{L} \subseteq \mathcal{H}$ which is dense in $\mathcal{H}$. Then 
we can express the norm $\| f(\cdot) \|_{\mathcal{H}}$ of any function $f(\cdot) \in \mathcal{H}$ as
\begin{equation}
\label{equ_approx_norm_HS_sup_dense_set}
\| f(\cdot) \|_{\mathcal{H}} = \sup_{\substack{g(\cdot) \in \mathcal{L} \\ \| g(\cdot) \|_{\mathcal{H}} > 0 }} \frac{\langle f(\cdot) , g(\cdot) \rangle_{\mathcal{H}}}{ \| g(\cdot) \|_{\mathcal{H}}}. 
\end{equation} 
\end{theorem}
\begin{proof}
By Theorem \ref{thm_cauchy_schwarz} we have that $|\langle f(\cdot) , g(\cdot) \rangle_{\mathcal{H}}| \leq \| f(\cdot) \|_{\mathcal{H}} \|g(\cdot) \|_{\mathcal{H}}$, 
which implies that 
\begin{equation}
\label{equ_proof_norm_approx_inner_prod_HS_dense_3} 
\sup_{\substack{g(\cdot) \in \mathcal{L}Ê\\ \| g(\cdot) \|_{\mathcal{H}} > 0}} \frac{\langle f(\cdot) , g(\cdot) \rangle_{\mathcal{H}}}{ \| g(\cdot) \|_{\mathcal{H}}} \leq \| f(\cdot) \|_{\mathcal{H}}.
\end{equation}
From this, the validity of \eqref{equ_approx_norm_HS_sup_dense_set} is obvious for $\| f(\cdot) \|_{\mathcal{H}} = 0$. 

It remains to prove \eqref{equ_approx_norm_HS_sup_dense_set} for the case $\| f(\cdot) \|_{\mathcal{H}} > 0$. 
However, in this case we can choose a $g(\cdot) \in \mathcal{L}$ with $\|g(\cdot)\|_{\mathcal{H}}>0$ such 
that $\frac{\langle f(\cdot) , g(\cdot) \rangle_{\mathcal{H}}}{ \| g(\cdot) \|_{\mathcal{H}}}$ is arbitrarily close to $ \| f(\cdot) \|_{\mathcal{H}}$. 
Indeed, since $\mathcal{L}$ is dense in $\mathcal{H}$ we can find for any sufficiently small $\varepsilon>0$ a function $g(\cdot) \in \mathcal{L}$ such that 
$\| f(\cdot) - g(\cdot) \|_{\mathcal{H}} \leq \varepsilon$ and $\| g(\cdot) \|_{\mathcal{H}} > 0$. 
This specific function $g(\cdot)$ yields
\begin{align} 
\frac{\langle f(\cdot) , g(\cdot) \rangle_{\mathcal{H}}}{ \| g(\cdot) \|_{\mathcal{H}}} & = \frac{\langle f(\cdot)- g(\cdot) + g(\cdot), g(\cdot)\rangle_{\mathcal{H}}}{ \| g(\cdot) \|_{\mathcal{H}}}
= \frac{\langle f(\cdot)-g(\cdot) , g(\cdot)\rangle_{\mathcal{H}}}{ \| g(\cdot) \|_{\mathcal{H}}} + \frac{\langle g(\cdot) , g(\cdot) \rangle_{\mathcal{H}}}{ \| g(\cdot) \|_{\mathcal{H}}}  \nonumber \\[3mm]
&  \stackrel{(a)}{\geq} -\frac{\|f(\cdot)-g(\cdot)\|_{\mathcal{H}} \|g(\cdot) \|_{\mathcal{H}} }{ \| g(\cdot) \|_{\mathcal{H}}} + \frac{\langle g(\cdot) , g(\cdot) \rangle_{\mathcal{H}}}{ \| g(\cdot) \|_{\mathcal{H}}}  \geq -\varepsilon +\| g(\cdot) \|_{\mathcal{H}} \nonumber \\[3mm]
& = -\varepsilon +\|  f(\cdot) + g(\cdot) - f(\cdot)  \|_{\mathcal{H}} \stackrel{\eqref{equ_triangle_inequality}}{\geq} -\varepsilon  + \| f(\cdot) \|_{\mathcal{H}} - \| g(\cdot) - f(\cdot) \|_{\mathcal{H}}  \nonumber \\[3mm]
& = -\varepsilon -  \| f(\cdot) - g(\cdot) \|_{\mathcal{H}} + \| f(\cdot) \|_{\mathcal{H}} \geq  - 2 \varepsilon+ \| f(\cdot) \|_{\mathcal{H}}    ,
\label{equ_proof_norm_approx_inner_prod_HS_dense_1}
\end{align}
where step $(a)$ is due to Theorem \ref{thm_cauchy_schwarz}.
From \eqref{equ_proof_norm_approx_inner_prod_HS_dense_1} we conclude that 
\begin{equation}
\label{equ_proof_norm_approx_inner_prod_HS_dense_4} 
\sup_{\substack{g(\cdot) \in \mathcal{L}Ê\\ \| g(\cdot) \|_{\mathcal{H}} > 0}} \frac{\langle f(\cdot) , g(\cdot) \rangle_{\mathcal{H}}}{ \| g(\cdot) \|_{\mathcal{H}}} \geq \| f(\cdot) \|_{\mathcal{H}}.
\end{equation} 
Combining \eqref{equ_proof_norm_approx_inner_prod_HS_dense_3} and \eqref{equ_proof_norm_approx_inner_prod_HS_dense_4} yields \eqref{equ_approx_norm_HS_sup_dense_set}.
\end{proof}
 
\begin{theorem}
\label{thm_linspan_finite_set_subspace}
Consider a function Hilbert space $\mathcal{H}$, a finite index set $\mathcal{T}$ (which can be assumed without loss of generality to be of the form $\mathcal{T} = [L]$ with some $L \in \mathbb{N}$), and the set of functions $\{ f_{l}(\cdot) \in \mathcal{H} \}_{l \in \mathcal{T}}$. 
We then have that the linear span $\linspan \{ f_{l}(\cdot)\}_{l \in \mathcal{T}}$ is a subspace of $\mathcal{H}$. 
\end{theorem} 
\begin{proof} 
The linear span $\linspan \{ f_{l}(\cdot)\}_{l \in \mathcal{T}}$, for a finite index set $\mathcal{T}$, together with the inner product $\langle \, \cdot \, , \, \cdot \, \rangle_{\mathcal{H}}$ is a finite dimensional inner-product function space, i.e.,
a special case of an abstract finite dimensional inner-product space. From this, the statement follows from the well-known fact, that any finite dimensional inner-product space is complete \cite{RudinBook,HalmosFiniteDimVecSp}. 
\end{proof} 

\begin{theorem}
\label{thm_fourier_onb}
Consider a function Hilbert space $\mathcal{H}$ and an ONB $\{ g_{l}(\cdot) \}_{l \in \mathcal{T}}$ for it, where $\mathcal{T}$ is an arbitrary index set. 
We then have for any $f(\cdot) \in \mathcal{H}$ 
\begin{equation}
\label{equ_fourier_series_squared_norm}
\| f(\cdot) \|_{\mathcal{H}}^{2} = \sum_{l \in \mathcal{T}} \big| \big\langle g_{l}(\cdot), f(\cdot) \big\rangle_{\mathcal{H}}\big|^{2},
\end{equation} 
and
\begin{equation} 
\label{equ_basis_expansion_onb_hilbert_space}
f(\cdot) = \sum_{l \in \mathcal{T}} \big \langle f(\cdot,  g_{l}(\cdot) \big \rangle_{\mathcal{H}} \,\, g_{l}(\cdot).
\end{equation}
\end{theorem}
\begin{proof}
\cite{RudinBook}
\end{proof}
The meaning of \eqref{equ_basis_expansion_onb_hilbert_space} in case of an infinite index set $\mathcal{T}$ is that the function $f(\cdot)$ can be 
approximated arbitrarily well by the sums $\sum_{l \in \mathcal{T}'} \big \langle f(\cdot,  g_{l}(\cdot) \big \rangle_{\mathcal{H}} g_{l}(\cdot)$, where $\mathcal{T}' \subseteq \mathcal{T}$ is an arbitrary finite subset of $\mathcal{T}$, i.e., 
\begin{equation}
\inf_{\mathcal{T}' \subseteq \mathcal{T}} \bigg\| f(\cdot) -  \sum_{l \in \mathcal{T}'} \big \langle f(\cdot,  g_{l}(\cdot) \big \rangle_{\mathcal{H}} \,\, g_{l}(\cdot) \bigg\|_{\mathcal{H}} = 0.
\end{equation}

\begin{theorem}
\label{thm_basic_congruence} 
Consider two Hilbert spaces $\mathcal{H}_{1}$ and $\mathcal{H}_{2}$ as well as two sets of functions $\mathcal{A} \triangleq \{ f_{j}(\cdot) \in \mathcal{H}_{1} \}_{j \in \mathcal{T}}$ and $ \mathcal{B} \triangleq \{ g_{j}(\cdot) \in \mathcal{H}_{2} \}_{j \in \mathcal{T}}$, indexed 
by the same set $\mathcal{T}$, that span $\mathcal{H}_{1}$ and $\mathcal{H}_{2}$ respectively, i.e., $\linspan \{ \mathcal{A} \}$ is dense in $\mathcal{H}_{1}$ and $\linspan \{ \mathcal{B} \}$ is dense in $\mathcal{H}_{2}$. Then, if
\begin{equation}
t,s \in \mathcal{T} \, \Rightarrow \, \langle f_{t}(\cdot) , f_{s}(\cdot)  \rangle_{\mathcal{H}_{1}} = \langle g_{t}(\cdot) , g_{s}(\cdot)  \rangle_{\mathcal{H}_{2}},
\end{equation} 
the two Hilbert spaces $\mathcal{H}_{1}$, $\mathcal{H}_{2}$ are isometric. 
A sufficient condition for a linear mapping $\mathsf{J}[\cdot]: \mathcal{H}_{1} \rightarrow \mathcal{H}_{2}$, which satisfies
\begin{equation}
\label{equ_basic_congruence_map_linear_space}
t \in \mathcal{T} \Rightarrow \mathsf{J}[f_{t}(\cdot)] = g_{t}(\cdot),
\end{equation}
to be a congruence from $\mathcal{H}_{1}$ to $\mathcal{H}_{2}$ is that 
\begin{equation} 
\label{equ_basic_congruence_map_continuous}
\mathsf{J}[f(\cdot)]  = \lim_{l \rightarrow \infty} \mathsf{J}[f_{l}(\cdot)],  
\end{equation}
for any Cauchy sequence $\{ f_{l} (\cdot) \in \linspan\{Ê\mathcal{A} \} \}_{l \rightarrow \infty}$ with limit $f(\cdot) \in \mathcal{H}_{1}$.

\end{theorem} 
\begin{proof}
\cite{Parzen59} 
\end{proof}

 \begin{theorem}
 \label{thm_orthog_proj_ineq}
 Consider a function Hilbert space $\mathcal{H}$, an arbitrary function $ f(\cdot) \in \mathcal{H}$, and a subspace $\mathcal{U} \subseteq \mathcal{H}$. We then have for any $ g(\cdot) \in \mathcal{U}$
 \begin{equation}
 \label{equ_error_orthog_projection}
 \langle f(\cdot) - \mathbf{P}_{\mathcal{U}} f(\cdot) , g(\cdot) \rangle_{\mathcal{H}} = 0, 
 \end{equation} 
 which implies that 
 \begin{equation}
 \label{equ_squared_norm_orthog_proj_pythag_thm}
 \| \mathbf{P}_{\mathcal{U}} f(\cdot) \|^{2}_{\mathcal{H}} =  \| f(\cdot) \|^{2}_{\mathcal{H}} -  \| f(\cdot)- \mathbf{P}_{\mathcal{U}} f(\cdot) \|^{2}_{\mathcal{H}} \leq \| f(\cdot) \|_{\mathcal{H}}^{2}.  
 \end{equation}
 \end{theorem} 
 \begin{proof}
 \cite{RudinBook}
 \end{proof}

\begin{theorem}
\label{thm_orth_proj_finite_dim_onb}
Consider a function Hilbert space $\mathcal{H}$ and a subspace $\mathcal{U} \subseteq \mathcal{H}$ for which an ONB $\{v_{l}(\cdot)\}_{l \in \mathcal{T}}$ with a finite index set $\mathcal{T}$ exists. 
Then the orthogonal projection and squared norm of a function $f(\cdot) \in \mathcal{H}$ onto the subspace $\mathcal{U}$ is given by 
\begin{equation}
\label{equ_orth_proj_finite_dim_onb}
\mathbf{P}_{\mathcal{U}} f(\cdot) = \sum_{l \in \mathcal{T}} \big\langle f(\cdot), v_{l}(\cdot) \big\rangle_{\mathcal{H}} \,\, v_{l}(\cdot), 
\end{equation}
and 
\begin{equation}
\label{equ_orth_proj_finite_dim_onb_squared_norm}
\| \mathbf{P}_{\mathcal{U}} f(\cdot) \|^{2}_{\mathcal{H}} = \sum_{l \in \mathcal{T}} \big| \big\langle f(\cdot), v_{l}(\cdot) \big\rangle_{\mathcal{H}} \big|^{2}, 
\end{equation}
respectively. 
\end{theorem} 
\begin{proof}  
\cite{RudinBook}
\end{proof}


\begin{theorem}
\label{thm_norm_projection_finite_dim_subspace} 
Consider a function Hilbert space $\mathcal{H}$ and a finite number of functions $v_{l} (\cdot) \in \mathcal{H}$, $l\in [L]$. 
Then the squared norm of the orthogonal projection of a function $f(\cdot) \in \mathcal{H}$ on the subspace $\mathcal{U} \triangleq \linspan \{ v_{l}(\cdot) \}_{l \in [L]}$ 
(cf.\ Theorem \ref{thm_linspan_finite_set_subspace}) is given by 
\begin{equation}
\label{equ_norm_projection_finite_dim_subspace}
\| \mathbf{P}_{\mathcal{U}} f(\cdot) \|_{\mathcal{H}}^{2} = \mathbf{c}^{T} \mathbf{G}^{\dagger} \mathbf{c}   
\end{equation}
where the vector $\mathbf{c} \in \mathbb{R}^{L}$ and matrix $\mathbf{G} \in \mathbb{R}^{L \times L}$ 
are defined elementwise by $c_{l} = \langle f(\cdot) , v_{l}(\cdot) \rangle_{\mathcal{H}}$ and $\left( \mathbf{G} \right)_{m,n} = \langle v_{m}(\cdot) , v_{n}(\cdot) \rangle_{\mathcal{H}}$ respectively. 
\end{theorem}

\begin{proof} 
According to Theorem \ref{thm_orthog_proj_ineq}, we can decompose any $f(\cdot) \in \mathcal{H}$ as a sum $f(\cdot) = \mathbf{P}_{\mathcal{U}} f(\cdot) + h(\cdot)$, where $\langle h(\cdot) , g (\cdot) \rangle_{\mathcal{H}}=0$ for 
every $g(\cdot) \in \mathcal{U}$. This implies that for any $l \in [L]$ we have
\begin{equation}
\label{equ_proof_squared_norm_projec_finite_subspace_inner_prod}
\langle  f(\cdot), v_{l}(\cdot) \rangle_{\mathcal{H}} = \langle  \mathbf{P}_{\mathcal{U}} f(\cdot) , v_{l}(\cdot) \rangle_{\mathcal{H}}.
\end{equation} 

Since, per definition, the orthogonal projection $\mathbf{P}_{\mathcal{U}} f(\cdot)$ belongs to the subspace $\mathcal{U} \triangleq \linspan \{ v_{l}(\cdot) \}_{l \in [L]}$, we 
can represent it as 
\begin{equation} 
\label{equ_proof_squared_norm_projec_finite_subspace_basis_expansion_projection}
\mathbf{P}_{\mathcal{U}} f(\cdot) = \sum_{l \in [L]}  d_{l} v_{l}(\cdot), 
\end{equation} 
with suitable coefficients $d_{l} \in \mathbb{R}$. The squared norm of the projection can then, due to the properties \eqref{equ_cond_inner_product}, 
be expressed as 
\begin{equation}
\label{equ_proof_squared_norm_projec_finite_subspace_squared_norm_projection}
\| \mathbf{P}_{\mathcal{U}} f(\cdot) \|_{\mathcal{H}}^{2} = \mathbf{d}^{T} \mathbf{G} \mathbf{d}, 
\end{equation} 
where $\mathbf{G}$ as defined above and the vector $\mathbf{d} \in \mathbb{R}^{L}$ is obtained by stacking the coefficients $d_{l}$. 
From \eqref{equ_proof_squared_norm_projec_finite_subspace_basis_expansion_projection} and the linearity of the inner product
one can verify using \eqref{equ_proof_squared_norm_projec_finite_subspace_inner_prod}, that the vector $\mathbf{d}$ and the above defined vector $\mathbf{c} \in \mathbb{R}^{L}$ are related by
\begin{equation}
\label{equ_proof_squared_norm_projec_finite_subspace_relation_coeff_vectors}
\mathbf{c} = \mathbf{G} \mathbf{d}.  
\end{equation} 

Using the identity $\mathbf{G} \mathbf{G}^{\dagger} \mathbf{G} = \mathbf{G}$ \cite{golub96}, we finally obtain
\begin{align} 
\mathbf{c}^{T} \mathbf{G}^{\dagger} \mathbf{c}  \stackrel{\eqref{equ_proof_squared_norm_projec_finite_subspace_relation_coeff_vectors}}{=} \mathbf{d}^{T} \mathbf{G} \mathbf{G}^{\dagger} \mathbf{G}\mathbf{d}
= \mathbf{d}^{T} \mathbf{G} \mathbf{d} \stackrel{\eqref{equ_proof_squared_norm_projec_finite_subspace_squared_norm_projection}}{=} \| \mathbf{P}_{\mathcal{U}} f(\cdot) \|_{\mathcal{H}}^{2}. \nonumber
\end{align}
\end{proof}

We will need in the following a slight variation of Theorem \ref{thm_norm_projection_finite_dim_subspace} stated in 
\begin{theorem} 
\label{thm_norm_projection_finite_dim_subspace_union_orthog_subspaces} 
Consider a function Hilbert space $\mathcal{H}$ and two sets of functions $\{v_{l} (\cdot) \in \mathcal{H}Ê\}_{l \in [L_{1}]}$ and $\{w_{l} (\cdot) \in \mathcal{H}Ê\}_{l \in [L_{2}]}$ 
that are mutually orthogonal, i.e., $\langle v_{l}(\cdot), w_{l'} (\cdot) \rangle_{\mathcal{H}}=0$ for any pair of indices $l \in [L_{1}]$, $l' \in [L_{2}]$. 
Then the squared norm of the orthogonal projection of a function $f(\cdot) \in \mathcal{H}$ on the subspace $\mathcal{U} \triangleq \linspanÊ\big\{ \{v_{l} (\cdot) Ê\}_{l \in [L_{1}]} \cup  \{ w_{l}(\cdot) \}_{l \in [L_2]} \big\}$ 
is given by 
\begin{equation}
\label{equ_idendity_projection_subspace_spanned_union_orthogonal_sets}
\| \mathbf{P}_{\mathcal{U}} f(\cdot) \|_{\mathcal{H}}^{2} = \mathbf{c}_{1}^{T} \mathbf{G}_{1}^{\dagger} \mathbf{c}_{1}   +  \mathbf{c}_{2}^{T} \mathbf{G}_{2}^{\dagger} \mathbf{c}_{2}, 
 \end{equation}
 where the vectors $\mathbf{c}_{1} \in \mathbb{R}^{L_{1}}$, $\mathbf{c}_{2} \in \mathbb{R}^{L_{2}}$ and matrices $\mathbf{G}_{1} \in \mathbb{R}^{L_{1} \times L_1}$, $\mathbf{G}_{2} \in \mathbb{R}^{L_2 \times L_2}$ are defined elementwise by 
 $c_{1,l} \triangleq \langle f(\cdot) , v_{l}(\cdot) \rangle_{\mathcal{H}}$,  $c_{2,l} \triangleq \langle f(\cdot) , w_{l}(\cdot) \rangle_{\mathcal{H}}$  and  $\left( \mathbf{G}_{1} \right)_{m,n} = \langle v_{m}(\cdot) , v_{n}(\cdot) \rangle_{\mathcal{H}}$, 
 $\left( \mathbf{G}_{2} \right)_{m,n} = \langle w_{m}(\cdot) , w_{n}(\cdot) \rangle_{\mathcal{H}}$, respectively.
\end{theorem} 

\begin{proof}
We merge the two sets $\{v_{l} (\cdot) \in \mathcal{H}Ê\}_{l \in [L_{1}]}$ and $\{w_{l} (\cdot) \in \mathcal{H}Ê\}_{l \in [L_{2}]}$ into one large set $\mathcal{A} \triangleq \{ u_{l}(\cdot) \}_{l \in [L_{1}+L_{2}] }$ such that 
$u_{l} (\cdot) = v_{l}(\cdot)$ for $l \in [L_{1}]$ and $u_{l}(\cdot) = w_{l-L_{1}}(\cdot)$ for $l \in [L_{1}+L_{2}] \setminus [L_{1}]$ and then define the matrix $\mathbf{G} \in \mathbb{R}^{(L_{1}+L_{2}) \times (L_{1} + L_{2})}$ elementwise 
by $\left( \mathbf{G} \right)_{k,l} \triangleq \langle u_{l}(\cdot) ,u_{k}(\cdot) \rangle_{\mathcal{H}}$.
Obviously, we have $\linspan \{ \mathcal{A} \} = \mathcal{U}$. 
Note that $\mathbf{G}$ is a block diagonal matrix, i.e., $\mathbf{G} = \begin{pmatrix} \mathbf{G}_{1} & \mathbf{0} \\ \mathbf{0} & \mathbf{G}_{2} \end{pmatrix}$, which implies that  
also its  pseudo-inverse is block diagonal, i.e., $\mathbf{G}^{\dagger} = \begin{pmatrix} \mathbf{G}_{1}^{\dagger} & \mathbf{0} \\ \mathbf{0} & \mathbf{G}_{2}^{\dagger} \end{pmatrix}$ (cf.\ \cite[Theorem 4.2.15]{Horn91}). 
Similarly, we define the vector $\mathbf{c} \in \mathbb{R}^{L_1 + L_2}$ via $c_{l} = \langle f(\cdot), u_{l}(\cdot) \rangle_{\mathcal{H}}$ and observe that 
$\mathbf{c} = \big( \mathbf{c}_{1}^{T}Ê\quad \mathbf{c}_{2}^{T} \big)^{T}$. 
The application of Theorem \ref{thm_norm_projection_finite_dim_subspace}, remember that $ \mathcal{U}=\linspan \{ \mathcal{A} \} $, yields 
\begin{equation}
\| \mathbf{P}_{\mathcal{U}} f(\cdot) \|_{\mathcal{H}}^{2}  = \mathbf{c}^{T} \mathbf{G}^{\dagger} \mathbf{c} =\big( \mathbf{c}_{1}^{T}Ê\quad \mathbf{c}_{2}^{T} \big) \begin{pmatrix} \mathbf{G}_{1}^{\dagger} & \mathbf{0} \\ \mathbf{0} & \mathbf{G}_{2}^{\dagger} \end{pmatrix}  \begin{pmatrix} \mathbf{c}_{1}Ê\\ \mathbf{c}_{2} \end{pmatrix} =  \mathbf{c}_{1}^{T} \mathbf{G}_{1}^{\dagger} \mathbf{c}_{1}   +  \mathbf{c}_{2}^{T} \mathbf{G}_{2}^{\dagger} \mathbf{c}_{2}. \nonumber
\end{equation} 
\end{proof}


\begin{theorem}
\label{thm_isometry_hilbert_space_coeffs_space}
Consider a function Hilbert space with an ONB $\{ g_{l}(\cdot) \in \mathcal{H} \}_{l \in \mathcal{T}}$ with an arbitrary index set $\mathcal{T}$. 
We then have that the Hilbert space $\mathcal{H}$ is isometric to the Hilbert space $\ell^{2}(\mathcal{T})$ with inner product $\big\langle f[\cdot], g[\cdot] \big\rangle_{\ell^{2}(\mathcal{T})} \triangleq \sum_{l \in \mathcal{T}} f[l]g[l]$. 
Furthermore, a congruence is given by $\mathsf{J}[\cdot]: \mathcal{H}Ê\rightarrow \ell^{2}(\mathcal{T}): f(\cdot) \in \mathcal{H} \mapsto h[\cdot] \in \ell^{2}(\mathcal{T})$ where $h[l] = \big\langle f(\cdot), g_{l} (\cdot) \big\rangle_{\mathcal{H}}$. 
Any element $f(\cdot) \in \mathcal{H}$ of the Hilbert space $\mathcal{H}$ can be written as 
\begin{equation}
\label{equ_basis_expansion_isometry_hilbert_space_coeffs}
f(\cdot) = \sum_{l \in \mathcal{T}} a[l] g_{l}(\cdot)
\end{equation}
with a coefficient sequence $a[l] \in \ell^{2}(\mathcal{T})$, and we have 
\begin{equation}
\label{equ_isometry_hilbert_space_coeffs_space_inner_prod_preserv}
\big\langle f(\cdot), g(\cdot) \big\rangle_{\mathcal{H}}  = \big\langle \mathsf{J}[f(\cdot)], \mathsf{J}[g(\cdot)] \big\rangle_{\ell^{2}(\mathcal{T})}.  
\end{equation}
Conversely, for any coefficient sequence $a[l] \in \ell^{2}(\mathcal{T})$ the function given by \eqref{equ_basis_expansion_isometry_hilbert_space_coeffs} 
is an element of $\mathcal{H}$.
\end{theorem}
\begin{proof}
\cite[p. 85]{RudinBook} 
\end{proof}
The meaning of \eqref{equ_basis_expansion_isometry_hilbert_space_coeffs} in case of an infinite index set $\mathcal{T}$ is that the function $f(\cdot)$ can be 
approximated arbitrarily well by the sums $\sum_{l \in \mathcal{T}'} a[l] g_{l}(\cdot)$, where $\mathcal{T}' \subseteq \mathcal{T}$ is an arbitrary finite subset of $\mathcal{T}$, i.e., 
\begin{equation}
\inf_{\substack{\mathcal{T}' \subseteq \mathcal{T}\\ a[l] \in \mathbb{R}}} \bigg\| f(\cdot) -  \sum_{l \in \mathcal{T}'} a[l] g_{l}(\cdot) \bigg\|_{\mathcal{H}} = 0.
\end{equation}



\section{Reproducing Kernel Hilbert Spaces (RKHS)}

A RKHS $\mathcal{H}$ is a special kind of Hilbert space consisting of real-valued functions $f(\cdot): \mathcal{D} \rightarrow \mathbb{R}$ 
defined over a specific domain $\mathcal{D}$ (which will always be a subset of $\mathbb{R}^{N}$ within this thesis). 
The special thing about a RKHS $\mathcal{H}$ is that the point evaluations $f_{\mathbf{y}}[\cdot]: \mathcal{H} \rightarrow \mathbb{R}: g(\cdot) \mapsto g(\mathbf{y})$ are continuous functionals \cite{aronszajn1950}. 
According to \cite{Parzen59}, we have 
\begin{definition} 
Given a function Hilbert space $\mathcal{H}$ over the domain $\mathcal{D}$ (i.e., the Hilbert space $\mathcal{H}$ consists of functions $f: \mathcal{D} \rightarrow \mathbb{R}$) 
and a function $R(\cdot,\cdot): \mathcal{D} \times \mathcal{D} \rightarrow \mathbb{R}$, we call $\mathcal{H}$ a RKHS with reproducing kernel $R(\cdot,\cdot)$
if the following two statements are valid: 
\begin{itemize}
\item The function $f(\cdot) \triangleq R(\cdot, \mathbf{x})$ obtained by fixing the second argument of the kernel to an arbitrary $\mathbf{x} \in \mathcal{D}$, i.e., $f(\mathbf{x}') = R(\mathbf{x}', \mathbf{x})$,
is an element of the RKHS $\mathcal{H}$, i.e., 
\begin{equation} 
\label{equ_kernel_func_in_RKHS}
\mathbf{x} \in \mathcal{D} \, \Rightarrow \, R(\cdot, \mathbf{x}) \in \mathcal{H}. 
\end{equation} 
\item The inner product $\big\langle \, \cdot \,  , \, \cdotÊ\, \big\rangle_{\mathcal{H}}$ on the RKHS $\mathcal{H}$ satisfies the following ``reproducing property'': 
\begin{equation} 
\label{equ_reproduction_property}
f(\cdot) \in \mathcal{H}, \mathbf{x} \in \mathcal{D} \,\,\, \Rightarrow \,\,\, \big\langle f( \cdot), R(\cdot, \mathbf{x}) \big\rangle_{\mathcal{H}} = f(\mathbf{x}). 
\end{equation}  
\end{itemize}
\end{definition}

We make the notion of a kernel function precise by 
\begin{definition}
\label{equ_def_kernel_function} 
Given an arbitrary set $\mathcal{D}$, we call a real-valued function $R(\cdot,\cdot): \mathcal{D} \times \mathcal{D} \rightarrow \mathbb{R}$ a ``kernel function over the domain $\mathcal{D}$'' if 
it is symmetric, i.e., 
\begin{equation}
\mathbf{x}_{1}, \mathbf{x}_{2} \in \mathcal{D} \,\,\, \Rightarrow \,\,\, R(\mathbf{x}_{1}, \mathbf{x}_{2}) =  R(\mathbf{x}_{2}, \mathbf{x}_{1}) 
\end{equation} 
and for any finite set $\{ \mathbf{x}_{k} \}_{k \in [K]}$ of points $\mathbf{x}_{k} \in \mathcal{D}$, the matrix $\mathbf{R} \in \mathbb{R}^{K \times K}$ defined elementwise 
by $\left( \mathbf{R} \right)_{k,l} = R(\mathbf{x}_{k}, \mathbf{x}_{l})$
is psd, i.e., 
\begin{equation}
\label{equ_finite_mtx_kernel_PSD}
\mathbf{R} \geq \mathbf{0}. 
\end{equation} 
\end{definition} 

As shown in \cite{aronszajn1950}, a necessary and sufficient condition for the existence of a (necessarily unique) reproducing kernel for a given Hilbert space $\mathcal{H}$ is 
that the point evaluations $f_{\mathbf{y}}[\cdot]: \mathcal{H} \rightarrow \mathbb{R}: g(\cdot)  \mapsto g(\mathbf{y})$ are continuous functionals. 
However, for our purposes, it is important to observe that in the opposite direction, i.e., given a kernel function $R(\cdot,\cdot): \mathcal{D} \times \mathcal{D} \rightarrow \mathbb{R}$, one can always associate 
a RKHS to it \cite{RKHSGaussprior, aronszajn1950,sun_jfaa,HeinRKHS2004}.  
\begin{definition}
\label{def_linear_space_kernel}
Given a kernel function $R(\cdot,\cdot): \mathcal{D} \times \mathcal{D} \rightarrow \mathbb{R}$ over a domain $\mathcal{D}$, we associate with 
it the linear function space denoted by $\mathcal{L}(R)$ and 
defined as 
\begin{equation}
\label{equ_def_linear_space}
\mathcal{L}(R) \triangleq \linspan \{ R(\cdot, \mathbf{x}) \}_{\mathbf{x} \in \mathcal{D}}.
\end{equation}
That is, any $f(\cdot) \in \mathcal{L}(R)$ can be written as a finite sum $f(\cdot) = \sum_{l \in [L]} a_{l} R(\cdot, \mathbf{x}_{l})$ with certain coefficients $a_{l} \in \mathbb{R}$ and points $\mathbf{x}_{l} \in \mathcal{D}$. 
\end{definition} 

\begin{theorem}
\label{thm_def_inner_prod_kernel}
Given an arbitrary kernel function $R(\cdot,\cdot): \mathcal{D} \times \mathcal{D} \rightarrow \mathbb{R}$, we have that on the linear space $\mathcal{L}(R)$, one can naturally define an inner product 
denoted by $\langle f(\cdot), g(\cdot) \rangle_{\mathcal{L}(R)}$. This inner product satisfies the basic properties \eqref{equ_cond_inner_product} and is completely specified by the requirement
\begin{equation}
\label{equ_def_inner_prod_kernel}
\mathbf{x}_{1}, \mathbf{x}_{2} \in \mathcal{D} \,\,\, \Rightarrow \,\,\, \langle R(\cdot, \mathbf{x}_{1}), R(\cdot, \mathbf{x}_{2}) \rangle_{\mathcal{L}(R)} = R(\mathbf{x}_{1}, \mathbf{x}_{2}).
\end{equation}  
Given two arbitrary functions $f(\cdot), g(\cdot) \in \mathcal{L}(R)$ and their representations $f(\cdot) = \sum_{l \in [L_{1}]} a_{l} R(\cdot,\mathbf{x}_{1,l})$, $g(\cdot) = \sum_{l \in [L_{2}]} b_{l} R(\cdot,\mathbf{x}_{2,l})$ with 
suitable coefficients $a_{l},b_{l} \in \mathbb{R}$ and points $\mathbf{x}_{1,l}, \mathbf{x}_{2,l} \in \mathcal{D}$, their inner product is given by 
\begin{equation} 
\label{equ_def_inner_prod_kernel_inner_prod_arb_functions}
\langle f(\cdot), g(\cdot) \rangle_{\mathcal{L}(R)} = \sum_{l \in [L_{1}]}\sum_{l' \in [L_{2}]} a_{l} b_{l'} R(\mathbf{x}_{2,l'},\mathbf{x}_{1,l}).  
\end{equation} 
Moreover, if two functions $f(\cdot), g(\cdot)$ are such that $\langle f(\cdot) - g(\cdot), f(\cdot) - g(\cdot) \rangle_{\mathcal{L}(R)} = 0$, we have that $f(\mathbf{x}) = g(\mathbf{x})$ for every $\mathbf{x} \in \mathcal{D}$. 
\end{theorem}
\begin{proof}
The properties $(a)$ and $(b)$ in \eqref{equ_cond_inner_product} are evident from \eqref{equ_def_inner_prod_kernel_inner_prod_arb_functions}. 
Note that \eqref{equ_def_inner_prod_kernel_inner_prod_arb_functions} contains as a special case the reproducing property \eqref{equ_reproduction_property} since for $g(\cdot) = R(\cdot,\mathbf{x}_{2})$, i.e., $L'=1$ and $b_{1}=1$, we obtain 
\begin{equation} 
\langle f(\cdot), g(\cdot) \rangle_{\mathcal{L}(R)} = \sum_{l \in [L_{1}]} a_{l} R(\mathbf{x}_{2},\mathbf{x}_{1,l}) = f(\mathbf{x}_{2}).
\end{equation} 
As already discussed, the property $(d)$ in \eqref{equ_cond_inner_product} is always satisfied if we 
define $f(\cdot)Ê\equiv 0$ to mean that $\langle f(\cdot) - n(\cdot), f(\cdot) -n(\cdot)\rangle_{\mathcal{L}(R)} = 0$, where $n(\cdot)$ denotes the zero function, i.e., $n(\mathbf{x}) = 0$ for every $\mathbf{x} \in \mathcal{D}$. 

The property $(c)$ in \eqref{equ_cond_inner_product}, i.e., $\langle f(\cdot), f(\cdot) \rangle_{\mathcal{L}(R)} \geq 0$ for any $f(\cdot) = \sum_{l \in [L_{1}]} a_{l} R(\cdot, \mathbf{x}_{l}) \in \mathcal{L}(R)$ 
follows straightforwardly from \eqref{equ_finite_mtx_kernel_PSD} and the formula \eqref{equ_def_inner_prod_kernel_inner_prod_arb_functions} evaluated for 
$f(\cdot)=g(\cdot)$. Indeed, 
\begin{equation} 
\langle f(\cdot), f(\cdot) \rangle_{\mathcal{L}(R)} = \sum_{l \in [L_{1}]}\sum_{l' \in [L_{1}]} a_{l} a_{l'} R(\mathbf{x}_{1,l'},\mathbf{x}_{1,l}) = \mathbf{a}^{T} \mathbf{R} \mathbf{a}Ê\geq 0,
\end{equation} 
with the vector $\mathbf{a} \in \mathbb{R}^{L_1}$ whose entries are the coefficients 
$a_{l}$ and psd matrix $\mathbf{R} \in \mathbb{R}^{L_1 \times L_1}$ whose elements are given by $\left( \mathbf{R} \right)_{m,n} \triangleq R(\mathbf{x}_{1,m},\mathbf{x}_{1,n})$. 

Finally, we have that $\langle f(\cdot) - g(\cdot), f(\cdot) - g(\cdot) \rangle_{\mathcal{L}(R)} = 0$ implies $f(\mathbf{x}) = g(\mathbf{x})$ for every $\mathbf{x} \in \mathcal{D}$. In particular, 
$f(\cdot) \equiv 0$ implies that $f(\mathbf{x}) = 0$ for every $\mathbf{x} \in \mathcal{D}$. 
Indeed, we have by \eqref{equ_reproduction_property} (which is a special case of \eqref{equ_def_inner_prod_kernel_inner_prod_arb_functions}) and Theorem \ref{thm_cauchy_schwarz} that for any 
$\mathbf{x} \in \mathcal{D}$
\begin{align} 
|f(\mathbf{x}) - g(\mathbf{x})| & \stackrel{\eqref{equ_reproduction_property}}{=} |\langle f(\cdot) - g(\cdot), R(\cdot, \mathbf{x}) \rangle_{\mathcal{L}(R)}| \nonumber \\[4mm]
& \stackrel{\eqref{equ_cauchy_schwarz}}{\leq} \sqrt{ \langle f(\cdot) - g(\cdot), f(\cdot) - g(\cdot) \rangle_{\mathcal{L}(R)}} \sqrt{ \langle R(\cdot,\mathbf{x}),  R(\cdot, \mathbf{x}) \rangle_{\mathcal{L}(R)}} \nonumber \\[4mm]
& \stackrel{\eqref{equ_def_inner_prod_kernel}}{=} \sqrt{ \langle f(\cdot) - g(\cdot), f(\cdot) - g(\cdot) \rangle_{\mathcal{L}(R)}} \sqrt{ R(\mathbf{x}, \mathbf{x})}=0, \nonumber
\end{align} 
i.e., 
$f(\mathbf{x}) = g(\mathbf{x})$ for any $\mathbf{x} \in \mathcal{D}$.  
\end{proof} 

With a slight abuse of notation, we will use from now on the same symbol $\mathcal{L}(R)$ for both the linear space given by \eqref{equ_def_linear_space} and the inner-product space that is given 
by this linear space together with the inner product $\langle f(\cdot), g(\cdot) \rangle_{\mathcal{L}(R)}$ defined by \eqref{equ_def_inner_prod_kernel}. It should be clear from the context 
what precise meaning of the symbol $\mathcal{L}(R)$ is used. 

Based on the inner-product space $\mathcal{L}(R)$, one can show
\begin{theorem} 
\label{thm_constr_RKHS_closure_linear_span}
Consider an arbitrary kernel function $R(\cdot,\cdot): \mathcal{D} \times \mathcal{D} \rightarrow \mathbb{R}$ over a domain $\mathcal{D}$, and let $\mathcal{C}$ denote 
the set of all functions $f(\cdot): \mathcal{D} \rightarrow \mathbb{R}$ that are the pointwise limits of a Cauchy sequence in the inner-product space $\mathcal{L}(R)$.
Then, the set of functions $f(\cdot): \mathcal{D} \rightarrow \mathbb{R}$ denoted by $\mathcal{H}(R)$ and defined as 
\begin{equation}
\mathcal{H}(R) \triangleq \mathcal{L}(R) \cup \mathcal{C}
\end{equation}
forms a RKHS with kernel $R(\cdot,\cdot)$. 
Moreover, the linear space $\mathcal{L}(R)$ is dense in the Hilbert space $\mathcal{H}(R)$, i.e., the set $\{ R(\cdot, \mathbf{x}) \in \mathcal{L}(R) \}_{\mathbf{x} \in \mathcal{D}}$ spans 
the RKHS $\mathcal{H}(R)$. 
The inner product on the Hilbert space $\mathcal{H}(R)$ is completely determined by the basic properties \eqref{equ_cond_inner_product} and by \eqref{equ_reproduction_property}. 

\end{theorem}
\begin{proof} 
\cite{aronszajn1950,Parzen59,RKHSGaussprior}
\end{proof} 
Given a kernel $R(\cdot, \cdot): \mathcal{D} \times \mathcal{D} \rightarrow \mathbb{R}$, we have that the RKHS $\mathcal{H}(R)$ constructed in Theorem \ref{thm_constr_RKHS_closure_linear_span} 
is the unique RKHS associated with this kernel, as stated in 
\begin{theorem}
\label{thm_exist_uniq_RKHS}
Given an abitrary kernel function $R(\cdot,\cdot): \mathcal{D} \times \mathcal{D} \rightarrow \mathbb{R}$ over a domain $\mathcal{D}$, there exists a unique RKHS which 
satisfies \eqref{equ_kernel_func_in_RKHS} and \eqref{equ_reproduction_property}. This unique RKHS coincides with $\mathcal{H}(R)$ constructed in Theorem \ref{thm_constr_RKHS_closure_linear_span}.
\end{theorem} 
\begin{proof}
The statement about the existence is proven trivially by Theorem \ref{thm_constr_RKHS_closure_linear_span}. 
A proof of the uniqueness can be found in 
\cite{Parzen59,aronszajn1950}. 
\end{proof}
\begin{example}
A well-known example of a RKHS is the Euclidean space $\mathbb{R}^{N}$ endowed with the inner product $\langle \mathbf{x}, \mathbf{y} \rangle = \sum_{k \in [N]} x_{k} y_{k}$. 
This is a RKHS of functions over the domain $\mathcal{D} = [N]$, which can be represented by length $N$ vectors. This RKHS is associated to the kernel $R(\cdot,\cdot): [N] \times [N] \rightarrow \mathbb{R}: R(k,l) = \delta_{k,l}$, i.e., 
the kernel can be represented by the identity matrix $\mathbf{I}_{N}$. For any $k \in [N]$ we have $\mathbf{e}_{k} = R(\cdot,k) \in \mathbb{R}^{N}$ which verifies \eqref{equ_kernel_func_in_RKHS}. 
Next, for any vector $\mathbf{a} \in \mathbb{R}^{N}$ it holds that $a_{k} = \mathbf{a}^{T} \mathbf{e}_{k} = \langle \mathbf{a}, R(\cdot,k) \rangle$ verifying \eqref{equ_reproduction_property}. 
Finally, one can verify \eqref{equ_def_inner_prod_kernel} since $\langle R(\cdot,k), R(\cdot,l) \rangle = \mathbf{e}_{k}^{T} \mathbf{e}_{l} = \delta_{k,l} = R(k,l)$. 
\end{example}

\section{Important Facts about RKHS} 

\begin{theorem}
\label{thm_densensess_lin_space_denseness_rkhs}
Consider a RKHS $\mathcal{H}(R)$ and a set $\{ v_{l}(\cdot) \in \mathcal{H}(R) \}_{l \in \mathcal{T}}$, where $\mathcal{T}$ is an arbitrary index set. 
If the set $\{ v_{l}(\cdot) \in \mathcal{H}(R) \}_{l \in \mathcal{T}}$ is dense in $\mathcal{L}(R)$, then it is also dense in the RKHS $\mathcal{H}(R)$. 
\end{theorem}
\begin{proof}
Consider an arbitrary $f(\cdot) \in \mathcal{H}(R)$. Since by Theorem \ref{thm_constr_RKHS_closure_linear_span} we have that the linear space $\mathcal{L}(R)$ is dense in $\mathcal{H}(R)$, 
we can find for any $\varepsilon >0$ a function $g(\cdot) \in \mathcal{L}(R)$ such that $\| g(\cdot) - f(\cdot) \|_{\mathcal{H}(R)} \leq \varepsilon$.  
Furthermore, since the set $\{ v_{l}(\cdot) \in \mathcal{H}(R) \}_{l \in \mathcal{T}}$ is dense in $\mathcal{L}(R)$, 
we can also find a function $v_{l}(\cdot)$ such that $\| g(\cdot) - v_{l} (\cdot) \|_{\mathcal{H}(R)} \leq \varepsilon$. By the triangle inequality \eqref{equ_triangle_inequality}, we have 
\begin{align}
\| f(\cdot) - v_{l} (\cdot) \|_{\mathcal{H}(R)} &= \| (f(\cdot)- g(\cdot)) + (g(\cdot) - v_{l} (\cdot) )\|_{\mathcal{H}(R)} \nonumber \\[4mm]
& \leq \| f(\cdot) - g(\cdot) \|_{\mathcal{H}(R)} +  \| g(\cdot) - v_{l} (\cdot) \|_{\mathcal{H}(R)}   \leq 2\varepsilon.
\end{align} 
Thus, we have shown that any $f(\cdot) \in \mathcal{H}(R)$ can be approximated arbitrary well by a function $v_{l}(\cdot)$. 
\end{proof}

\begin{theorem}
\label{thm_norm_impl_point_convergence} 
Consider a RKHS $\mathcal{H}(R)$ over the domain $\mathcal{D}$. If a sequence $\{ f_{l}(\cdot) \in \mathcal{H}(R) \}_{l \rightarrow \infty}$ converges to 
a function $f(\cdot) \in \mathcal{H}(R)$, i.e, $\lim_{l \rightarrow \infty} f_{l}(\cdot) = f(\cdot)$, then the sequence converges also pointwise to $f(\cdot)$, i.e., for every $\mathbf{x} \in \mathcal{D}$ 
we have that $\lim_{l \rightarrow \infty} f_{l}(\mathbf{x})   = f(\mathbf{x})$. 
Moreover, this point-wise convergence is uniform on every subset $\mathcal{D}' \subseteq \mathcal{D}$ on which $R(\mathbf{x}, \mathbf{x})$ is uniformly bounded, i.e.,
there exists a constant $C$ such that $R(\mathbf{x}, \mathbf{x}) \leq C$ for every $\mathbf{x} \in \mathcal{D}'$. 
\end{theorem} 
\begin{proof}
\cite{aronszajn1950}
\end{proof}

A direct consequence of Theorem \ref{thm_norm_impl_point_convergence} which is of major importance for our purposes is stated in 
\begin{theorem} 
\label{thm_point_convergent_impl_norm}
Consider a RKHS $\mathcal{H}(R)$ over the domain $\mathcal{D}$ and an arbitrary Cauchy sequence $\{ f_{l}(\cdot) \in \mathcal{H}(R) \}_{l \rightarrow \infty}$. 
If the sequence converges pointwise to the function $f(\cdot): \mathcal{D} \rightarrow \mathbb{R}$, i.e., for every $\mathbf{x} \in \mathcal{D}$ we have $\lim_{l \rightarrow \infty} f_{l}(\mathbf{x}) = f(\mathbf{x})$,  
then the sequence converges to $f(\cdot)$, i.e., $\lim_{l \rightarrow \infty} f_{l}(\cdot) = f(\cdot)$ and $f(\cdot) \in \mathcal{H}(R)$.  
\end{theorem}
\begin{proof}
This result is proven by contradiction: 
Denote the limit of the Cauchy sequence by $f'(\cdot)$, i.e., $f'(\cdot)=\lim_{l \rightarrow \infty} f_{l}(\cdot)$. 
Note that the limit exists since per definition any Cauchy sequence in a Hilbert space has a limit. 
Now let us assume that the limit $f'(\cdot)$ is different from $f(\cdot)$. Since $f'(\cdot)$ and $f(\cdot)$ are different, there must be some point $\mathbf{x}' \in \mathcal{D}$ for which $f'(\mathbf{x}') \neq f(\mathbf{x}')$. 
However, since the sequence $\{f_{l}(\cdot) \}_{l \rightarrow \infty}$ converges to $f'(\cdot)$ in RKHS norm, we have according to Theorem \ref{thm_norm_impl_point_convergence} 
that the sequence $\{f_{l}(\mathbf{x}') \}_{l \rightarrow \infty}$ must converge to $f'(\mathbf{x}')$. 
But this is a contradiction to the fact that $\{f_{l}(\mathbf{x}') \}_{l \rightarrow \infty}$ converges to $f(\mathbf{x}')$. 
\end{proof}

\begin{theorem}
\label{thm_suff_con_congruence_RKHS}
Consider two RKHSs $\mathcal{H}(R_1)$ and $\mathcal{H}(R_2)$ defined over the domains $\mathcal{D}_{1}$ and $\mathcal{D}_{2}$, respectively, 
as well as two sets $\mathcal{A} = \{ f_{j}(\cdot) \in \mathcal{H}(R_1) \}_{j \in \mathcal{T}}$ and $\mathcal{B} = \{ g_{j}(\cdot) \in \mathcal{H}(R_2) \}_{j \in \mathcal{T}}$, indexed 
by the same index set $\mathcal{T}$, that span $\mathcal{H}(R_1)$ and $\mathcal{H}(R_2)$, respectively. 
Assume that these two sets are such that 
\begin{equation} 
\label{equ_suff_con_congruence_inner_prod_conser}
t,s \in \mathcal{T} \,  \Rightarrow \, \big\langle f_{t} (\cdot), f_{s} (\cdot)  \big\rangle_{\mathcal{H}(R_1)} = \big\langle g_{t} (\cdot) , g_{s}(\cdot)  \big\rangle_{\mathcal{H}(R_2)}.
\end{equation} 
Then a sufficient condition for  a linear mapping $\mathsf{K}[\cdot]: \mathcal{H}(R_{1}) \rightarrow \mathcal{H}(R_{2})$, satisfying $\mathsf{K}[ f_{t}(\cdot)] = g_{t}(\cdot)$ for every $t \in \mathcal{T}$, 
to be a congruence from $\mathcal{H}(R_{1})$ to $\mathcal{H}(R_{2})$ is that 
for any sequence $\{ f_{l}(\cdot) \in \linspan \{ \mathcal{A}Ê\} \}_{l \rightarrow \infty }$ such that the pointwise limit function $f(\mathbf{x}) = \lim_{k \rightarrow \infty}f_{k}(\mathbf{x})$ exists and belongs to $\mathcal{H}(R_{1})$, we have that 
\begin{equation} 
\label{equ_suff_cond_cong_pointwise_converg}
\forall \mathbf{x} \in \mathcal{D}_{2}: \mathsf{K}[f(\cdot)] (\mathbf{x}) = \lim_{l \rightarrow \infty} \mathsf{K}[f_{l}(\cdot)] (\mathbf{x}).
\end{equation} 
\end{theorem}
\begin{proof} 
This result is proven by Theorem \ref{thm_basic_congruence}. Obviously we have that $\mathsf{K}[\cdot]$ satisfies \eqref{equ_basic_congruence_map_linear_space} and it only remains 
to verify \eqref{equ_basic_congruence_map_continuous}. 
To that end, consider an arbitrary Cauchy sequence $\{ f_{l}(\cdot) \in \linspan\{ \mathcal{A} \} \}_{l \rightarrow \infty}$ with limit $f(\cdot) \in \mathcal{H}(R_{1})$ and the associated sequence 
$\{ \mathsf{K}[f_{l}(\cdot)] \in \linspan\{ \mathcal{B}Ê\} \}_{l \rightarrow \infty}$ which due to \eqref{equ_suff_con_congruence_inner_prod_conser} is also a Cauchy sequence with some limit $g(\cdot) \in \mathcal{H}(R_{2})$. 
By Theorem \ref{thm_norm_impl_point_convergence}, we have that the sequence $\{ f_{l}(\cdot) \in \linspan \{ \mathcal{A}Ê\} \}_{l \rightarrow \infty}$ converges pointwise to $f(\cdot)$, 
which in turn implies via \eqref{equ_suff_cond_cong_pointwise_converg} that 
the Cauchy sequence $\{ \mathsf{K}[f_{l}(\cdot)](\cdot) \in \linspan \{ \mathcal{B} \} \}_{l \rightarrow \infty}$ converges pointwise to the function $\mathsf{K}[f(\cdot)](\cdot)$. 
Since by Theorem \ref{thm_point_convergent_impl_norm} and $\lim_{l \rightarrow \infty} \mathsf{K}[f_{l}(\cdot)] = g(\cdot)$, we have that $\mathsf{K}[f(\cdot)] = g(\cdot)$, we have 
verified \eqref{equ_basic_congruence_map_continuous}.  
\end{proof} 

A main theme of this thesis is how a certain RKHS $\mathcal{H}_{1}$ defined over domain $\mathcal{D}_{1}$ is related to another RKHS $\mathcal{H}_{2}$ that is obtained by restricting the functions 
$f(\cdot) \in \mathcal{H}_{1}$ to the subset $\mathcal{D}_{2} \subseteq \mathcal{D}_{1}$. 
The most important result concerning this question is stated in 
\begin{theorem}
\label{thm_reducing_domain_RKHS}
Consider a RKHS $\mathcal{H}(R_{1})$ associated with the kernel $R_{1}(\cdot,\cdot): \mathcal{D}_{1} \times \mathcal{D}_{1} \rightarrow \mathbb{R}$, whose elements are functions with domain $\mathcal{D}_{1}$. 
Then the set of functions that is obtained by restricting each function $f(\cdot) Ê\in \mathcal{H}(R_{1})$ to 
the subdomain $\mathcal{D}_{2} \subseteq \mathcal{D}_{1}$ coincides with the RKHS $\mathcal{H}(R_{2})$ associated with the kernel $R_{2}(\cdot,\cdot): \mathcal{D}_{2} \times \mathcal{D}_{2} \rightarrow \mathbb{R}$ 
that is the restriction of the kernel $R_{1}(\cdot, \cdot): \mathcal{D}_{1} \times \mathcal{D}_{1} \rightarrow \mathbb{R}$ to the subdomain $\mathcal{D}_{2} \times \mathcal{D}_{2}$, i.e., 
$R_{2} (\cdot, \cdot) = R_{1}(\cdot, \cdot)\big|_{\mathcal{D}_{2} \times \mathcal{D}_{2}}$. 
Furthermore, we have the relation 
\begin{equation} 
\label{equ_thm_reducing_domain_RKHS}
\| f(\cdot) \|_{\mathcal{H}(R_{2})}  = \min_{\substack{ g(\cdot) \in \mathcal{H}(R_{1}) \\ g(\cdot)\big|_{\mathcal{D}_{2}}\!=\!f(\cdot)}} \| g(\cdot) \|_{\mathcal{H}(R_{1})}.
\end{equation}
Thus, the norm of any element $f(\cdot) \in \mathcal{H}(R_{2})$ is equal to the minimal norm of a function $g(\cdot) \in \mathcal{H}(R_{1})$ that coincides with $f(\cdot)$ on the restricted domain $\mathcal{D}_{2}$. 
\end{theorem} 
\begin{proof}
\cite{aronszajn1950}
\end{proof}

In what follows, we will need
\begin{theorem}
\label{thm_pointwise_series_kernel_ONB}
Consider a RKHS $\mathcal{H}(R)$ whose kernel $R(\cdot,\cdot): \mathcal{D} \times \mathcal{D} \rightarrow \mathbb{R}$ is given pointwise by 
\begin{equation}
\label{equ_pointwise_series_kernel_ONB_kernel_sum}
R(\mathbf{x}_{1}, \mathbf{x}_{2}) = \sum\limits_{k \in \mathcal{T}} g_{k}(\mathbf{x}_{1})g_{k}(\mathbf{x}_{2}) 
\end{equation}
with an arbitrary index set $\mathcal{T}$.
The functions $g_{k}(\cdot) \in \mathcal{H}(R)$ are assumed to be orthonormal, i.e., 
\begin{equation} 
\big\langle g_{k}(\cdot) , g_{l}(\cdot) \big\rangle_{\mathcal{H}(R)} = \delta_{k,l}. 
\end{equation}
Then we have that the set $\{ g_{k}(\cdot)  \in \mathcal{H}(R) \}_{k \in \mathcal{T}}$ forms an ONB for the RKHS $\mathcal{H}(R)$. 
\end{theorem}
\begin{proof} 
The proof is a slight modification of the derivation of \cite[Theorem 3.12]{PaulsenRKHS}. 
First, we note that the sum $\sum\limits_{k \in \mathcal{T}} g_{k}(\mathbf{x}_{1})g_{k}(\mathbf{x}_{2})$ can be understood in the usual sense, i.e., as the supremum over all 
finite sums $\sum\limits_{k \in \mathcal{T}'} g_{k}(\mathbf{x}_{1})g_{k}(\mathbf{x}_{2})$ with a finite set $\mathcal{T}' \subseteq \mathcal{T}$. This is possible via the Cauchy-Schwarz inequality  
for the Hilbert space $\ell^{2}(\mathcal{T})$ (cf.\ Theorem \ref{thm_cauchy_schwarz} and \cite{RudinBook}) since the sequences $a[k] \triangleq g_{k}(\mathbf{x}_{1})$, $b[k] \triangleq g_{k}(\mathbf{x}_{2})$ 
both belong to $\ell^{2}(\mathcal{T})$ for any $\mathbf{x}_{1}, \mathbf{x}_{2} \in 
\mathcal{D}$. 

We now show that any function $f(\cdot) \in \mathcal{H}(R)$ can be approximated arbitrarily well by a finite linear combination of the functions $\{ g_{k}(\cdot)  \in \mathcal{H}(R) \}_{k \in \mathcal{T}}$, i.e., the set 
spans the RKHS $\mathcal{H}(R)$ which implies that $\{ g_{k}(\cdot)  \in \mathcal{H}(R) \}_{k \in \mathcal{T}}$ is an ONB for $\mathcal{H}(R)$. 
To this end, according to Theorem \ref{thm_densensess_lin_space_denseness_rkhs}, it is already sufficient 
to show that any function $f(\cdot) \in \mathcal{L}(R)$
in the linear span $\mathcal{L}(R)$ can be approximated arbitrarily well. 
 
Consider an arbitrary function $f(\cdot) \in \mathcal{L}(R)$ (cf. \eqref{equ_def_linear_space}), i.e., $f(\cdot) = \sum_{l \in [L]} a_{l} R(\cdot, \mathbf{x}_{l})$ with some $a_{l} \in \mathbb{R}$ and $\mathbf{x}_{l} \in \mathcal{D}$, $l \in [L]$. 
We then have that
\begin{align}
\label{equ_fourier_series_kernel_series_repr}
\| f(\cdot) \|^{2}_{\mathcal{H}(R)} & = \bigg\langle \sum_{l \in [L]} a_{l} R(\cdot, \mathbf{x}_{l}), \sum_{l \in [L]} a_{l} R(\cdot, \mathbf{x}_{l})\bigg\rangle_{\mathcal{H}(R)}  \nonumber \\[4mm]Ê
& =  \sum_{l,l' \in [L]} a_{l} a_{l'} R(\mathbf{x}_{l}, \mathbf{x}_{l'}) \stackrel{\eqref{equ_pointwise_series_kernel_ONB_kernel_sum}}{=} \sum_{l,l' \in [L]} a_{l} a_{l'} \sum\limits_{k \in \mathcal{T}} g_{k}(\mathbf{x}_{l})g_{k}(\mathbf{x}_{l'})  \nonumber \\[4mm]
& =  \sum_{l,l' \in [L]} a_{l} a_{l'} \sum\limits_{k \in \mathcal{T}} \bigg\langle g_{k}(\cdot), R(\cdot, \mathbf{x}_{l}) \bigg\rangle_{\mathcal{H}(R)} \bigg\langle g_{k}(\cdot), R(\cdot, \mathbf{x}_{l'}) \bigg\rangle_{\mathcal{H}(R)} \nonumber \\[4mm]
& \stackrel{(a)}{=} \sum\limits_{k \in \mathcal{T}} \bigg\langle g_{k}(\cdot),  \sum_{l \in [L]} a_{l} R(\cdot, \mathbf{x}_{l}) \bigg\rangle_{\mathcal{H}(R)} \bigg\langle g_{k}(\cdot),  \sum_{l' \in [L]} a_{l'} R(\cdot, \mathbf{x}_{l'}) \bigg\rangle_{\mathcal{H}(R)} \nonumber \\[4mm]
& = \sum\limits_{k \in \mathcal{T}} \bigg\langle g_{k}(\cdot),  f(\cdot) \bigg\rangle_{\mathcal{H}(R)} \bigg\langle g_{k}(\cdot),  f(\cdot) \bigg\rangle_{\mathcal{H}(R)} \nonumber \\[4mm]
& = \sum\limits_{k \in \mathcal{T}}  \big|\langle g_{k}(\cdot),  f(\cdot) \rangle_{\mathcal{H}(R)} \big|^{2},
\end{align}
where the change of summation order in $(a)$ can be validated by an argument based on the Cauchy-Schwarz inequality  
for the Hilbert space $\ell^{2}(\mathcal{T})$ (cf.\ Theorem \ref{thm_cauchy_schwarz} and \cite{RudinBook}).
Let us now introduce the function $s_{\mathcal{T}'}(\cdot) \triangleq \sum_{k \in \mathcal{T}'} \langle f(\cdot), g_{k}(\cdot) \rangle_{\mathcal{H}(R)} \, g_{k}(\cdot) \in \linspan\{ g_{k}(\cdot) \}_{k \in \mathcal{T}}$ for any finite set $\mathcal{T}' \subseteq \mathcal{T}$. 
It can be verified that according to Theorem \ref{thm_orth_proj_finite_dim_onb}, the function $s_{\mathcal{T}'}(\cdot)$ coincides with the orthogonal projection of $f(\cdot)$ on the subspace 
\begin{equation} 
\label{equ_proof_sum_kernel_ONB_subspaces_nested}
\mathcal{U}^{(\mathcal{T}')} \triangleq  \linspan\{ g_{k}(\cdot) \}_{k \in \mathcal{T}'} \subseteq \linspan\{ g_{k}(\cdot) \}_{k \in \mathcal{T}} \subseteq\mathcal{H}(R),
\end{equation}
i.e., 
$\mathbf{P}_{\mathcal{U}^{(\mathcal{T}')}} f(\cdot) = s_{\mathcal{T}'}(\cdot)$. Obviously, an ONB for the subspace $\mathcal{U}^{(\mathcal{T}')}$ is given by the 
function set $\linspan\{ g_{k}(\cdot) \}_{k \in \mathcal{T}'}$.
Therefore, by Theorem \ref{thm_orthog_proj_ineq} and Theorem \ref{thm_orth_proj_finite_dim_onb}, we have 
\begin{align}
\| f(\cdot) -s_{\mathcal{T}'}(\cdot) \|_{\mathcal{H}(R)}^{2} & = \| f(\cdot) -\mathbf{P}_{\mathcal{U}^{(\mathcal{T}')}} f(\cdot)  \|_{\mathcal{H}(R)}^{2} \stackrel{\eqref{equ_squared_norm_orthog_proj_pythag_thm}}{=} \| f(\cdot)  \|_{\mathcal{H}(R)}^{2} -\| \mathbf{P}_{\mathcal{U}^{(\mathcal{T}')}} f(\cdot)  \|_{\mathcal{H}(R)}^{2} \nonumber   \\[4mm]
& \stackrel{\eqref{equ_orth_proj_finite_dim_onb_squared_norm}}{=}  \| f(\cdot) \|^{2}_{\mathcal{H}} -   \sum\limits_{k \in \mathcal{T}'}  \big| \langle f(\cdot) , g_{k}(\cdot)  \rangle_{\mathcal{H}(R)}\big|^{2} \nonumber   \\[4mm]
& \stackrel{\eqref{equ_fourier_series_kernel_series_repr}}{=} \sum\limits_{k \in \mathcal{T}}  \big|\langle f(\cdot) , g_{k}(\cdot)  \rangle_{\mathcal{H}(R)} \big|^{2} -  \sum\limits_{k \in \mathcal{T}'}  
\big|\langle f(\cdot) , g_{k}(\cdot)  \rangle_{\mathcal{H}(R)}\big|^{2}  \nonumber \\[4mm]Ê
& = \sum\limits_{k \in \mathcal{T} \setminus \mathcal{T}'}  \big|\langle f(\cdot) , g_{k}(\cdot)  \rangle_{\mathcal{H}(R)}\big|^{2}. 
\end{align} 
From this it follows that 
\begin{equation}
\nonumber 
\inf_{\substack{|\mathcal{T}'| \inÊ\mathbb{N} \\ \mathcal{T}' \subseteq \mathcal{T}}} \| f(\cdot) - s_{\mathcal{T}'}(\cdot)\|_{\mathcal{H}(R)}=0,
\end{equation}
i.e., any function $f(\cdot) \in \mathcal{L}(R)$ can 
be 
approximated arbitrary well by an element of $\linspan \{ g_{k}(\cdot) \}_{k \in \mathcal{T}}$ due to \eqref{equ_proof_sum_kernel_ONB_subspaces_nested}. 
Thus, we have that $\linspan \{ g_{k}(\cdot) \}_{k \in \mathcal{T}}$ is dense in the linear space $\mathcal{L}(R)$ 
and by Theorem \ref{thm_densensess_lin_space_denseness_rkhs} it is also dense in $\mathcal{H}(R)$. This means in turn that the set $ \{ g_{k}(\cdot) \}_{k \in \mathcal{T}}$ of orthonormal functions 
spans the RKHS $\mathcal{H}(R)$, i.e., is an ONB for $\mathcal{H}(R)$.
\end{proof}



\section{RKHSs with a Differentiable Kernel} 
\label{sec_differentiable_kernel}

The main part of this thesis is concerned exclusively with RKHSs over a domain $\mathcal{D} \subseteq \mathbb{R}^{N}$ that are associated with a ``differentiable'' kernel \cite{sun_jfaa,zhou_jfaa}. 
It will turn out that under specific conditions, the functions of a RKHS associated with a differentiable kernel are characterized completely by their behavior in an arbitrarily small neighborhood around 
a single point belonging to the domain $\mathcal{D}$. Note that it will always be implicitly assumed that the standard topology, induced by the inner product $\mathbf{x}^{T} \mathbf{y}$, is used for $\mathbb{R}^{N}$. 

For a precise definition of a differentiable kernel, consider an index set $\mathcal{K} \subseteq [N]$ and 
a vector $\mathbf{x}_{c} \in \mathbb{R}^{N}$, and define the ``$\varepsilon$-$\mathcal{K}$-neighborhood'' of $\mathbf{x}_{c}$ by 
\begin{equation} 
\label{equ_def_epsilon_K_neighborhood}
\mathcal{N}^{\mathcal{K}}_{\mathbf{x}_{c}}(\varepsilon) \triangleq \{ \mathbf{x}_{c} + \mathbf{a} \big| \mathbf{a} \in \mathbb{R}_{+}^{N}, \supp(\mathbf{a}) = \mathcal{K}, \| \mathbf{a} \|_{\infty} \leq \varepsilon \}. 
\end{equation} 

\begin{definition}
\label{def_differentiable_kernel}
A kernel $R(\cdot,\cdot): \mathcal{D} \times \mathcal{D} \rightarrow \mathbb{R}$ over a domain $\mathcal{D} \subseteq \mathbb{R}^{N}$ is said to be differentiable up to order $m$
if for any $\mathbf{x}_{c}$ and $\mathcal{K}$ for which there exists an $\varepsilon >0$ such that $\mathcal{N}^{\mathcal{K}}_{\mathbf{x}_{c}}(\varepsilon) \subseteq \mathcal{D}$, 
the partial derivatives 
\begin{equation} 
\frac{\partial^{\mathbf{p}_{1}} \partial^{\mathbf{p}_{2}} R(\mathbf{x}_{1}, \mathbf{x}_{2})}{\partial \mathbf{x}_{1}^{\mathbf{p}_{1}} \partial \mathbf{x}_{2}^{\mathbf{p}_{2}} } \bigg|_{\mathbf{x}_{1}=\mathbf{x}_{2}= \mathbf{x}_{c}}
\end{equation}
exist for any orders $\mathbf{p}_{1},\mathbf{p}_{2} \in \mathbb{Z}_{+}^{N}$ with $\supp(\mathbf{p}_{1}), \supp(\mathbf{p}_{2})  \subseteq \mathcal{K}$ and $\|\mathbf{p}_{1}\|_{\infty},\|\mathbf{p}_{2}\|_{\infty} \leq m$, and moreover 
this partial derivatives are continuous functions of $(\mathbf{x}_{1},\mathbf{x}_{2})$ (viewed as a vector in $\mathbb{R}^{2N}$). 
We will call a RKHS $\mathcal{H}(R)$ differentiable, if it is associated to a differentiable kernel $R(\cdot,\cdot)$.  
\end{definition}

An important property of the RKHS associated to a differentiable kernel is stated in
\begin{theorem}
\label{thm_der_repr_prop}
Let $\mathcal{D} \subseteq \mathbb{R}^{N}$, and consider a RKHS $\mathcal{H}(R)$ associated to a kernel $R(\cdot,\cdot): \mathcal{D} \times \mathcal{D} \rightarrow \mathbb{R}$ which is differentiable up to order $m$. 
Then for any $\mathbf{p} \in \mathbb{Z}_{+}^{N}$ with $\|\mathbf{p}\|_{\infty}\leq m$, the function $g^{(\mathbf{p})}_{\mathbf{x}_{c}}(\cdot): \mathcal{D} \rightarrow \mathbb{R}$ defined by 
\begin{equation} 
\label{equ_def_part_der_func}
g^{(\mathbf{p})}_{\mathbf{x}_{c}} (\mathbf{x}) \triangleq \frac{\partial^{\mathbf{p}}  R(\mathbf{x}, \mathbf{x}_{2})}{\partial \mathbf{x}_{2}^{\mathbf{p}}} \bigg|_{\mathbf{x}_{2} = \mathbf{x}_{c}},
\end{equation}
where $\mathbf{x}_{c} \in \mathcal{D}$ is such that $\mathcal{N}^{\supp(\mathbf{p})}_{\mathbf{x}_{c}}(\varepsilon) \subseteq \mathcal{D}$
(for a suitable $\varepsilon >0$), is an element  of $\mathcal{H}(R)$, i.e., 
\begin{equation} 
g^{(\mathbf{p})}_{\mathbf{x}_{c}} (\cdot) \in \mathcal{H}(R). 
\end{equation} 
Furthermore, the inner product of the function $g^{(\mathbf{p})}_{\mathbf{x}_{c}} (\cdot)$ with an arbitrary function $f(\cdot) \in \mathcal{H}(R)$ is given by 
\begin{equation}
\label{equ_der_reproduction_prop} 
\bigg\langle g^{(\mathbf{p})}_{\mathbf{x}_{c}}(\cdot) , f(\cdot) \bigg\rangle_{\mathcal{H}(R)} = \frac{\partial^{\mathbf{p}}  f(\mathbf{x})}{\partial \mathbf{x}^{\mathbf{p}}} \bigg|_{\mathbf{x} = \mathbf{x}_{c}}.
\end{equation} 
\end{theorem} 
\begin{proof}
\cite{zhou_jfaa}
\end{proof}

We will also need the following slight variation of Theorem \ref{thm_der_repr_prop}. 
\begin{corollary}
\label{cor_der_repr_prop_mult_func}
Consider a RKHS $\mathcal{H}(R)$ associated to a kernel $R(\cdot,\cdot): \mathcal{D} \times \mathcal{D} \rightarrow \mathbb{R}$, $\mathcal{D} \subseteq \mathbb{R}^{N}$, which is differentiable up to order $m \in \mathbb{N}$. 
Furthermore, consider a matrix $\mathbf{A} \in \mathbb{R}^{N \times D}$, a vector $\mathbf{z}_{c} \in \mathbb{R}^{D}$, and a multi-index $\mathbf{p} \in \mathbb{Z}_{+}^{D}$ with $\|\mathbf{p}\|_{\infty}\leq m$ that 
satisfy $\mathbf{A}\mathcal{N}^{\supp(\mathbf{p})}_{\mathbf{z}_{c}}(\varepsilon) \subseteq \mathcal{D}$ for some small $\varepsilon >0$.
Then for any function $h(\cdot): \mathbb{R}^{D} \rightarrow \mathbb{R}$ 
for which the partial derivatives of any order $\mathbf{p}' \in \mathbb{Z}_{+}^{D}$ with $\|\mathbf{p}'\|_{\infty}\leq m$ exist
at the point $\mathbf{z}_{c} \in \mathbb{R}^{D}$, the function $f(\cdot): \mathcal{D} \rightarrow \mathbb{R}$ defined by 
\begin{equation} 
\label{equ_def_part_der_func_1}
f(\mathbf{x}) \triangleq \frac{\partial^{\mathbf{p}} \big( R(\mathbf{x}, \mathbf{A} \mathbf{z}_{2}) h(\mathbf{z}_{2}) \big)}{\partial \mathbf{z}_{2}^{\mathbf{p}}} \bigg|_{\mathbf{z}_{2} = \mathbf{z}_{c}},
\end{equation}
is an element  of $\mathcal{H}(R)$, i.e., 
\begin{equation} 
f(\cdot) \in \mathcal{H}(R). 
\end{equation} 
\end{corollary} 
\begin{proof}
The statement follows from Theorem \ref{thm_der_repr_prop} by induction on the order $\mathbf{p}$ by the chain rule and product rule of differentiation \cite{RudinBookPrinciplesMatheAnalysis}. 
We will restrict ourselves here to the begin of the induction, i.e., we consider $\mathbf{p} = \mathbf{e}_{k}$ with $k \in [D]$ which yields 
\begin{align}
f(\mathbf{x}) & =  \frac{\partial^{\mathbf{e}_{k}} \big( R(\mathbf{x}, \mathbf{A} \mathbf{z}_{2}) h(\mathbf{z}_{2}) \big)}{\partial \mathbf{z}_{2}^{\mathbf{e}_{k}}} \bigg|_{\mathbf{z}_{2} = \mathbf{z}_{c}}
 =  \frac{\partial^{\mathbf{e}_{k}}  R(\mathbf{x}, \mathbf{A} \mathbf{z}_{2}) }{\partial \mathbf{z}_{2}^{\mathbf{e}_{k}}} \bigg|_{\mathbf{z}_{2} = \mathbf{z}_{c}}  \frac{\partial^{\mathbf{e}_{k}}  h(\mathbf{z}_{2})}{\partial \mathbf{z}_{2}^{\mathbf{e}_{k}}} \bigg|_{\mathbf{z}_{2}= \mathbf{z}_{c}}
  \nonumber \\[3mm]
& = \Bigg[ \sum_{k' \in [D]} \frac{\partial^{\mathbf{e}_{k'}}R(\mathbf{x},x_{2})}{\partial \mathbf{x}_{2}^{\mathbf{e}_{k'}}}\bigg|_{\mathbf{x}_{2} = \mathbf{A} \mathbf{z}_{c}}  \left( \mathbf{A} \right)_{k',k} \Bigg]
 \frac{\partial^{\mathbf{e}_{k}}  h(\mathbf{z}_{2}) }{\partial \mathbf{z}_{2}^{\mathbf{e}_{k}}} \bigg|_{\mathbf{z}_{2}= \mathbf{z}_{c}}.
\end{align} 
By assumption, we have $\frac{\partial^{\mathbf{e}_{k'}}R(\cdot,x_{2})}{\partial \mathbf{x}_{2}^{\mathbf{e}_{k'}}}\bigg|_{\mathbf{x}_{2} = \mathbf{A} \mathbf{z}_{c}} \in \mathcal{H}(R)$ which implies, since $f(\cdot)$ is 
a linear combination of the functions $\frac{\partial^{\mathbf{e}_{k'}}R(\cdot,x_{2})}{\partial \mathbf{x}_{2}^{\mathbf{e}_{k'}}}\bigg|_{\mathbf{x}_{2} = \mathbf{A} \mathbf{z}_{c}}$ and $\mathcal{H}(R)$ is a Hilbert space (in particular a linear space), 
that $f(\cdot) \in \mathcal{H}(R)$. 
\end{proof}

In the next chapters, we will exclusively consider differentiable kernels that are such that the function $h(\mathbf{x}) \triangleq R(\mathbf{x}, \mathbf{x})$ is bounded on every bounded subset $\mathcal{D}' \subseteq \mathcal{D}$. 
Note that that a differentiable kernel up to an order $m \geq 1$ is always continuous, i.e., $\lim_{(\mathbf{x}, \mathbf{x}') \rightarrow (\mathbf{x}_{0},\mathbf{x}'_{0})} R(\mathbf{x}, \mathbf{x}') = R(\mathbf{x}_{0},\mathbf{x}'_{0})$, since 
the partial derivatives $\frac{\partial^{\mathbf{e}_{k}}R(\mathbf{x},\mathbf{x}_{2})}{\partial x_{2,k}}$ (and $\frac{\partial^{\mathbf{e}_{k}}R(\mathbf{x},\mathbf{x}_{2})}{\partial x_{k}}$ by symmetry of the kernel) are continuous and thus the kernel $R(\cdot,\cdot)$, viewed as a function 
with domain $\mathbb{R}^{2N}$ is differentiable in the usual sense of multivariable calculus \cite[Theorem 9.21]{RudinBookPrinciplesMatheAnalysis}. 
We have that the functions belonging to a RKHS whose kernel is continuous and bounded in the above sense must be continuous: 
\begin{theorem} 
\label{thm_continous_kernel_implies_cont_func}
Consider the RKHS $\mathcal{H}(R)$ associated with a continuous kernel $R(\cdot,\cdot): \mathcal{D} \times \mathcal{D} \rightarrow \mathbb{R}$, with $\mathcal{D} \subseteq \mathbb{R}^{N}$, for which 
$h(\mathbf{x}) \triangleq  R(\mathbf{x}, \mathbf{x})$ is bounded for $\mathbf{x} \in \mathcal{D}'$, where $\mathcal{D}'$ is any set of the form $\mathcal{D}' = \mathcal{D} \cap \mathcal{B}(\mathbf{x}_{c},r)$ 
with some center $\mathbf{x}_{c}Ê\in \mathbb{R}^{N}$ and radius $r >0$. 
Then every function $f(\cdot): \mathcal{D} \rightarrow \mathbb{R}$ belonging to $\mathcal{H}(R)$ is continuous.
\end{theorem}
\begin{proof}
Given an arbitrary but fixed function $f(\cdot) \in \mathcal{H}(R)$ and point $\mathbf{x}_{0} \in \mathcal{D}$, we have to show that for every $\varepsilon > 0$ we can find a $\delta>0$ such that 
$|f(\mathbf{x}_{0}) - f(\mathbf{x})| < \varepsilon$ for any $\mathbf{x} \in \mathcal{D}$ with $\| \mathbf{x} - \mathbf{x}_{0} \|_{2} < \delta$. 

First we observe that any function $f_{\mathbf{x}}(\cdot) = R(\cdot,\mathbf{x})$ with $\mathbf{x} \in \mathcal{D}$ is trivially continuous by our assumption that $R(\cdot, \cdot)$ in continuous. 
Since a finite linear combination of continuous functions 
is still continuous, we also have that any function $w(\cdot) \in \mathcal{L}(R)$ is continuous. 

Now consider the ball $\mathcal{B}(\mathbf{x}_{0},r)$ centered at $\mathbf{x}_{0}$ and with an arbitrary radius $r > 0$. According to our assumptions, 
there exists a constant $M$ such that $R(\mathbf{x},\mathbf{x}) < M$ for every $\mathbf{x} \in \mathcal{D} \cap \mathcal{B}(\mathbf{x}_{0},r)$. 
According to Theorem \ref{thm_constr_RKHS_closure_linear_span}, the linear space $\mathcal{L}(R)$ is dense in $\mathcal{H}(R)$, i.e., we can find a function $w(\cdot) \in \mathcal{L}(R)$ such that $\| w(\cdot) - f(\cdot) \|_{\mathcal{H}(R)} \leq \varepsilon/(4\sqrt{M})$ 
for an arbitrary but fixed $\varepsilon >0$. 
Since $w(\cdot)$ is continuous, we can find a radius $\delta\leq r$ such that $|w(\mathbf{x}_{0}) - w(\mathbf{x})| < \varepsilon/2$ for every $\mathbf{x} \in  \mathcal{D} \cap \mathcal{B}(\mathbf{x}_{0}, \delta)$. 
Now we have for any $\mathbf{x} \in \mathcal{D} \cap \mathcal{B}(\mathbf{x}_{0}, \delta)$ that
\begin{align} 
|f(\mathbf{x}_{0}) - f(\mathbf{x})| & = \big|w(\mathbf{x}_{0}) - w(\mathbf{x}) + [f(\mathbf{x}_{0}) - w(\mathbf{x}_{0})] + [w(\mathbf{x}) - f(\mathbf{x})] \big|  \nonumber \\[3mm]
& \stackrel{(a)}{=} \big|w(\mathbf{x}_{0}) - w(\mathbf{x}) +\langle f(\cdot)- w(\cdot), R(\cdot, \mathbf{x}_{0}) \rangle_{\mathcal{H}(R)}+\langle f(\cdot)- w(\cdot), R(\cdot, \mathbf{x}) \rangle_{\mathcal{H}(R)} \big|  \nonumber \\[3mm]
& \stackrel{(b)}{\leq} \big|w(\mathbf{x}_{0}) - w(\mathbf{x})\big| +\big|\langle f(\cdot)- w(\cdot), R(\cdot, \mathbf{x}_{0}) \rangle_{\mathcal{H}(R)}\big|+\big|\langle f(\cdot)- w(\cdot), R(\cdot, \mathbf{x}) \rangle_{\mathcal{H}(R)}\big|
\nonumber \\[3mm]
& \stackrel{(c)}{\leq} \varepsilon/2 + \| f(\cdot)- w(\cdot) \|_{\mathcal{H}(R)} \| R(\cdot, \mathbf{x}_{0}) \|_{\mathcal{H}(R)} + \| f(\cdot)- w(\cdot) \|_{\mathcal{H}(R)} \| R(\cdot, \mathbf{x}) \|_{\mathcal{H}(R)}
\nonumber \\[3mm]Ê
& \stackrel{(d)}{=} \varepsilon/2 + \| f(\cdot)- w(\cdot) \|_{\mathcal{H}(R)} \sqrt{R(\mathbf{x}_{0}, \mathbf{x}_{0})} + \| f(\cdot)- w(\cdot) \|_{\mathcal{H}(R)} \sqrt{R(\mathbf{x}, \mathbf{x})}
\nonumber \\[3mm]Ê
& \leq \varepsilon/2 +  \frac{\varepsilon}{4\sqrt{M}}  \sqrt{M} +  \frac{\varepsilon}{4\sqrt{M}}  \sqrt{M}
\nonumber \\[3mm]Ê
& = \varepsilon/2 +\varepsilon/4+\varepsilon/4 = \varepsilon,
\end{align} 
where $(a)$ and $(d)$ is due to the reproducing property \eqref{equ_reproduction_property}, $(b)$ follows from the inequality $|a+b| \leq |a| + |b|$ for $a,b \in \mathbb{R}$ and the step $(c)$ 
is obtained by an application of the Cauchy-Schwarz inequality (cf.\ Theorem \ref{thm_cauchy_schwarz}). 
\end{proof} 
While Theorem \ref{thm_continous_kernel_implies_cont_func} gives a sufficient condition on the kernel of a RKHS such that it consists only of continuous functions, a necessary and sufficient condition for this to be the case in a more 
general setting can be found in \cite[Theorem 5]{HeinRKHS2004}.

\chapter{The RKHS Approach to Minimum Variance Estimation}
\label{chap_RKHS_MVE}
\section{Introduction}

This chapter reviews the RKHS approach to minimum variance estimation as introduced in \cite{Duttweiler73b,Parzen59}. Using the 
RKHS approach, we will derive two fundamental results, which (to the best of the authors knowledge) seem to be novel. The first result makes a characterization of the minimum achievable variance $L_{\mathcal{M}}$ viewed as 
a function of the fixed parameter vector $\mathbf{x}_{0}$ of a minimum variance problem $\mathcal{M} = \minvarproblemscalar$, when $\mathbf{x}_{0}$ is varied. In particular, we show in Section \ref{sec_classes_min_var_problems_con_var_kernel} that 
$L_{\mathcal{M}}$, viewed as a function of $\mathbf{x}_{0}$, is lower semi-continuous. The second result is related to the concept of sufficient statistics. In particular, we prove in 
Section \ref{sec_suff_stat_RKHS}, that the RKHS associated to a minimum variance problem remains unchanged if the original observation is replaced by any sufficient statistic.
Within this chapter we will also use the RKHS approach to give a detailed geometric derivation and interpretation of some well-known lower bounds on the variance of estimators with a given bias. Although it has 
been already mentioned in \cite{Parzen59}, that virtually any known lower bound on the variance can be interpreted geometrically using RKHS theory, a detailed elaboration of this interpretation for some well-known bounds seems 
to be missing in the literature. 

From a mathematical viewpoint, the RKHS approach to minimum variance estimation 
deals with two specific Hilbert spaces which are naturally associated to a given minimum variance problem $\mathcal{M}=\minvarproblemscalar$ associated with the estimation problem $\scalarestproblem$. 

The first Hilbert space, which is a RKHS and denoted $\mathcal{H}(\mathcal{M})$, consists 
of all functions $\gamma(\cdot): \mathcal{X} \rightarrow \mathbb{R}$ that are obtained from a valid bias function $c(\cdot)$ for $\mathcal{M}$ by $\gamma(\cdot) = c(\cdot) + g(\cdot)$, where $g(\cdot): \mathcal{X} \rightarrow \mathbb{R}$ is the 
parameter function of $\mathcal{E}$. Thus, the RKHS $\mathcal{H}(\mathcal{M})$ consists of all functions $\gamma(\cdot)$ for which there exists an estimator for $\mathcal{M}$ with finite variance at $\mathbf{x}_{0}$ and whose mean function 
is $\gamma(\cdot)$. The function $\gamma(\cdot)$ will be referred to as the prescribed mean function of the minimum variance problem $\mathcal{M}$. 

The second Hilbert space, denoted $\mathcal{L}(\mathcal{M})$, which is naturally associated to the minimum variance problem $\mathcal{M}$ consists of estimators $\hat{g}(\cdot): \mathbb{R}^{M} \rightarrow \mathbb{R}$ 
with a finite variance at $\mathbf{x}_{0}$ but whose bias is not necessarily equal to the prescribed bias $c(\cdot)$ of $\mathcal{M}$. 

It will turn out that the Hilbert spaces $\mathcal{H}(\mathcal{M})$ and $\mathcal{L}(\mathcal{M})$ are isometric. Moreover, a specific congruence between these two Hilbert spaces will give us 
a tool for determining the LMV estimator (if it exists) for the minimum variance problem $\mathcal{M}$. 

In what follows, we will place an important constraint on the class of minimum variance problems that are considered. More precisely, 
we assume that for any considered minimum variance problem $\mathcal{M} = \minvarproblemscalar$ with associated estimation problem $\mathcal{E} = \scalarestproblem$, the following holds: 
\begin{postulate} 
\label{assumption_RKHS_minvarproblem} 
For any $\mathbf{x} \in \mathcal{X}$, 
\begin{equation}
\label{equ_finite_integral_RKHS_cond}
\mathsf{E}_{\mathbf{x}_{0}} \bigg \{ \left( \frac{ f(\mathbf{y};\mathbf{x})}{f(\mathbf{y}; \mathbf{x}_{0})} \right)^{2} \bigg \} < \infty.
\end{equation}
\end{postulate}
Unless otherwise stated, we will implicitly assume from now on that any minimum variance problem that is considered satisfies \eqref{equ_finite_integral_RKHS_cond}. We note 
that many works on minimum variance estimation introduce a variant of Postulate \ref{assumption_RKHS_minvarproblem} \cite{Duttweiler73b,Barankin49,Parzen59,stein50}. 

 
\section{A Hilbert Space of Estimators}
\label{sec_hilbert_space_est}

An estimator $\hat{g}(\cdot): \mathbb{R}^{M} \rightarrow \mathbb{R}$ for an estimation problem $\mathcal{E}=\scalarestproblem$ is nothing but a (possibly random) mapping that maps the observation space 
$\mathbb{R}^{M}$ to the range of the parameter function, i.e., $\mathbb{R}$.\footnote{We will assume a scalar parameter function $g(\mathbf{x})$, i.e., $P=1$ throughout this chapter (cf. Section \ref{sec_basic_concepts}).}  
However, for a given minimum variance problem $\mathcal{M}$, only the allowed estimators $\mathcal{F}(\mathcal{M})$ (see Definition \ref{def_est_finite_var_prescr_bias}) make sense. 
As mentioned in Section \ref{sec_esist_uniqu_mve}, one can easily verify that the set $\mathcal{F}(\mathcal{M})$ is an affine set of estimator functions, 
i.e., if $\hat{g}_{1}(\cdot)$, $\hat{g}_{2}(\cdot) \in \mathcal{F}(\mathcal{M})$ then also any affine combination (cf.\ \cite[p.\ 21]{BoydConvexBook}) $\hat{g}_{3}(\cdot) \triangleq a \hat{g}_{1}(\cdot) + (1-a) \hat{g}_{2}(\cdot)$, 
where $a \in \mathbb{R}$ is arbitrary, yields an allowed estimator, i.e., $\hat{g}_{3}(\cdot) \in \mathcal{F}(\mathcal{M})$. 
Moreover, every allowed estimator $\hat{g}(\cdot) \in \mathcal{F}(\mathcal{M})$ belongs to the Hilbert space $\mathcal{P}(\mathcal{M})$ associated with $\mathcal{M}$ that is defined as  
\begin{equation}
\label{equ_def_est_finite_power} 
\mathcal{P}(\mathcal{M}) \triangleq \big \{ \hat{g}(\cdot)  \big| P(\hat{g}(\cdot); \mathbf{x}_{0}) < \infty  \big \}.
\end{equation} 
(We recall that $P(\hat{g}(\cdot); \mathbf{x})  = v(\hat{g}(\cdot); \mathbf{x}) + \big| \mathsf{E}_{\mathbf{x}} \big\{ \hat{g}(\mathbf{y}) \} |^{2}$.) 
To put it formally, we have \cite{Parzen59}
\begin{theorem}
\label{thm_hilbert_space_est_functions}
Given a minimum variance problem $\mathcal{M}=\minvarproblemscalar$, the set of estimator functions $\mathcal{P}(\mathcal{M})$ defined in \eqref{equ_def_est_finite_power} forms 
a function Hilbert space together with the inner product $\langle \hat{g}_{1}(\cdot), \hat{g}_{2}(\cdot) \rangle_{\emph{\text{\tiny{RV}}}}$ defined by 
\begin{equation}
\label{equ_def_inner_prod_RV}
\langle \hat{g}_{1}(\cdot), \hat{g}_{2}(\cdot) \rangle_{\emph{\text{\tiny{RV}}}} \triangleq \mathsf{E}_{\mathbf{x}_{0}} \{ \hat{g}_{1}(\mathbf{y}) \hat{g}_{2}(\mathbf{y}) \}. 
\end{equation} 
The Hilbert space $\mathcal{P}(\mathcal{M})$ contains the set $\mathcal{F}(\mathcal{M})$ of allowed estimators for $\mathcal{M}$ as defined in \eqref{equ_est_finite_var_prescr_bias}, i.e., 
\begin{equation}
\label{equ_allowed_estimators_contained_reasonable_estimators}
\mathcal{F}(\mathcal{M}) \subseteq \mathcal{P}(\mathcal{M}). 
\end{equation} 
\end{theorem} 

\begin{proof} 
By definition, an (possibly random) estimator $\hat{g}_{1}(\cdot)$ with finite stochastic power at $\mathbf{x}_{0}$ is nothing but a random variable with a finite stochastic power at $\mathbf{x}_{0}$. Therefore, 
the Hilbert space structure of $\mathcal{P}(\mathcal{M})$ with the inner product  $\langle \hat{g}_{1}(\cdot), \hat{g}_{2}(\cdot) \rangle_{\text{\tiny{RV}}}$ follows from elementary probability- or measure theory (see, e.g., \cite{Parzen59, BillingsleyProbMeasure,AshProbMeasure,RudinBook}). 

Furthermore, by the definition of $\mathcal{F}(\mathcal{M})$ and \eqref{equ_rel_power_var}, we have that the stochastic power at $\mathbf{x}_{0}$ of 
any allowed estimator $\hat{g}(\cdot)\in \mathcal{F}(\mathcal{M})$ satisfies 
\begin{equation} 
P(\hat{g}(\cdot); \mathbf{x}_{0}) = v(\hat{g}(\cdot); \mathbf{x}_{0}) + \big[ b(\hat{g}(\cdot) ; \mathbf{x}_{0}) + g(\mathbf{x}_{0})\big]^2 
= v(\hat{g}(\cdot); \mathbf{x}_{0}) + \big[ c(\mathbf{x}_{0}) + g(\mathbf{x}_{0}) \big]^2  < \infty,  \nonumber
\end{equation}
verifying \eqref{equ_allowed_estimators_contained_reasonable_estimators}. 
\end{proof}

It will turn out that, in the context of minimum variance estimation, the set $\mathcal{P}(\mathcal{M})$ can be further reduced without affecting 
the minimum achievable variance (cf.\ \eqref{equ_def_min_ach_var}). 
In particular, we will show that for a given minimum variance problem $\mathcal{M}$ and associated estimation problem $\mathcal{E}$ 
that satisfy \eqref{equ_finite_integral_RKHS_cond}, 
it is sufficient to consider only estimators which belong to the subspace $\mathcal{L}(\mathcal{M}) \subseteq \mathcal{P}(\mathcal{M})$ defined by 
\begin{equation}
\label{equ_def_hilbert_space_rhos}
\mathcal{L}(\mathcal{M}) \triangleq \closure \{ \linspan \{ \rho_{\mathcal{M}}(\cdot, \mathbf{x})  \}_{\mathbf{x} \in \mathcal{X}} \} , 
\end{equation}
where $ \rho_{\mathcal{M}}(\cdot, \cdot) : \mathbb{R}^{M} \times \mathcal{X} \rightarrow \mathbb{R}$ is the likelihood ratio of the minimum variance problem $\mathcal{M}$ defined as 
\begin{equation} 
\label{equ_def_likelihood}
 \rho_{\mathcal{M}}(\mathbf{y}, \mathbf{x})  \triangleq \frac{f(\mathbf{y}; \mathbf{x})}{f(\mathbf{y}; \mathbf{x}_{0})}.
\end{equation}
Note that the condition \eqref{equ_finite_integral_RKHS_cond} of our Postulate \ref{assumption_RKHS_minvarproblem} can be written as 
\begin{equation} 
\mathsf{E}_{\mathbf{x}_{0}} \big \{ \big( \rho_{\mathcal{M}}(\mathbf{y}; \mathbf{x}) \big)^{2} \big \} < \infty.
\end{equation} 
The closure in \eqref{equ_def_hilbert_space_rhos} relative to the Hilbert space $\mathcal{P}(\mathcal{M})$ (see Definition \ref{def_relative_closure}) exists, since we have that 
\begin{equation}
\label{equ_span_rho_subspace}
 \linspan \{ \rho_{\mathcal{M}}(\cdot, \mathbf{x})  \}_{\mathbf{x} \in \mathcal{X}} \subseteq \mathcal{P}(\mathcal{M}).
\end{equation} 
This can be easily verified by noting that for any estimator of the form $\hat{g}(\cdot) \triangleq \rho_{\mathcal{M}}(\cdot, \mathbf{x})$, where $\mathbf{x} \in \mathcal{X}$, it holds that
\begin{align} 
P(\hat{g}(\cdot);\mathbf{x}_{0}) & = \mathsf{E}_{\mathbf{x}_{0}} \big \{Ê\left( \hat{g}(\mathbf{y}) \right)^{2} \big \} =  \mathsf{E}_{\mathbf{x}_{0}} \big \{Ê\left( \rho_{\mathcal{M}}(\mathbf{y}, \mathbf{x})\right)^{2} \big \}Ê < \infty, 
\end{align} 
where the inequality follows from \eqref{equ_finite_integral_RKHS_cond}.

We emphasize the fact that in contrast to the set of allowed estimators $\mathcal{F}(\mathcal{M})$, the set $\mathcal{L}(\mathcal{M})$ is by its very definition a function Hilbert space with associated 
inner product $\langle \hat{g}_{1}(\cdot), \hat{g}_{2}(\cdot) \rangle_{\text{\tiny{RV}}}$. 
This is a very important fact since, as we will see in Section \ref{sec_main_result_RKHS_MVE}, it guarantees that if the set $\mathcal{F}(\mathcal{M})$ is nonempty, the infimum in \eqref{equ_def_min_ach_var} can be achieved by a suitable estimator, 
i.e., there exists an LMV estimator if there is at least one estimator with the prescribed bias and finite variance at $\mathbf{x}_{0}$.

\section{The RKHS Associated to a Minimum Variance Problem}

In this section, we introduce the section Hilbert space, beside $\mathcal{L}(\mathcal{M})$ defined in \eqref{equ_def_hilbert_space_rhos}, which is naturally associated with a given minimum variance problem. 
This second Hilbert space has a particular structure, i.e., it is a RKHS. 

\subsection{Definition} 
\label{sec_def_RKHS_MVP}
Consider a minimum variance problem $\mathcal{M}=\minvarproblemscalar$, associated with the estimation problem $\mathcal{E}=\scalarestproblem$, 
that satisfies the condition of Postulate \ref{assumption_RKHS_minvarproblem}. 
We then have for any pair $\mathbf{x}_{1}, \mathbf{x}_{2} \in \mathcal{X}$ that the corresponding likelihood ratios $\rho_{\mathcal{M}}(\cdot, \mathbf{x}_{1})$, $\rho_{\mathcal{M}}(\cdot, \mathbf{x}_{2})$ 
belong to the Hilbert space $\mathcal{P}(\mathcal{M})$ with inner product $\langle \rho_{\mathcal{M}}(\cdot, \mathbf{x}_{1}), \rho_{\mathcal{M}}(\cdot, \mathbf{x}_{2}) \rangle_{\text{\tiny{RV}}}$. 

Therefore, one can naturally associate with the minimum variance problem $\mathcal{M}$ a real-valued function 
$R_{\mathcal{M}}(\cdot,\cdot): \mathcal{X} \times \mathcal{X} \rightarrow \mathbb{R}$ defined as \cite{Parzen59}
\vspace*{2mm}
\begin{align} 
\label{equ_def_kernel_M}
R_{\mathcal{M}} (\mathbf{x}_{1}, \mathbf{x}_{2}) & \triangleq \mathsf{E}_{\mathbf{x}_{0}} \big\{ \rho_{\mathcal{M}}(\mathbf{y}, \mathbf{x}_{1}) \rho_{\mathcal{M}}(\mathbf{y}, \mathbf{x}_{2}) \big\} = \big\langle\rho_{\mathcal{M}}(\cdot, \mathbf{x}_{1}), \rho_{\mathcal{M}}(\cdot, \mathbf{x}_{2}) \big\rangle_{\text{\tiny{RV}}} \nonumber \\[4mm]
& = \int_{\mathbf{y} \in \mathbb{R}^{M}}  \frac{ f(\mathbf{y};\mathbf{x}_{1})f(\mathbf{y}; \mathbf{x}_{2})}{f(\mathbf{y}; \mathbf{x}_{0})}  d \mathbf{y}.
\vspace*{2mm}
\end{align} 
\vspace*{2mm}
An easily verified and useful observation is that $R_{\mathcal{M}} (\mathbf{x}_{0}, \mathbf{x})  = 1$ for every $\mathbf{x} \in \mathcal{X}$.

The function $R_{\mathcal{M}}(\cdot,\cdot)$ has an important property as stated in
\begin{theorem}
The function $R_{\mathcal{M}}(\cdot,\cdot): \mathcal{X} \times \mathcal{X} \rightarrow \mathbb{R}$ 
is a kernel function with domain $\mathcal{X}$ in the sense of Definition \ref{equ_def_kernel_function}. Moreover, for every $\mathbf{x} \in \mathcal{X}$ it holds that 
\begin{equation}
\label{equ_kernel_one_arg_x_0_equal_1}
R_{\mathcal{M}}(\mathbf{x}, \mathbf{x}_{0}) =  R_{\mathcal{M}}(\mathbf{x}_{0}, \mathbf{x}) = 1, 
\end{equation} 
where $\mathbf{x}_{0}$ is the fixed parameter vector associated with the minimum variance problem $\mathcal{M}$. 
\end{theorem}

\begin{proof}
Consider an arbitrary finite set $\{ \mathbf{x}_{l} \}_{l \in [L]}$ consisting of $L \in \mathbb{N}$ points $\mathbf{x}_{l} \in \mathcal{X}$. The matrix $\mathbf{R} \in \mathbb{R}^{L \times L}$ defined elementwise by 
$\left( \mathbf{R}Ê\right)_{l,l'} \triangleq R_{\mathcal{M}}(\mathbf{x}_{l}, \mathbf{x}_{l'})$
can be verified to coincide with the correlation matrix $\mathsf{E}_{\mathbf{x}_{0}} \big\{ \mathbf{r} \mathbf{r}^{T} \big\}$ of the random vector $\mathbf{r} \in \mathbb{R}^{L}$ 
given elementwise by $r_l \triangleq \rho_{\mathcal{M}}(\mathbf{y}, \mathbf{x}_{l})$, i.e., 
$\mathbf{R} = \mathsf{E}_{\mathbf{x}_{0}} \big\{ \mathbf{r} \mathbf{r}^{T} \big\}$. 
Thus, since any correlation matrix of a real-valued random vector is psd \cite{papoulis}, we have that $\mathbf{R}$ is psd for any finite set $\{ \mathbf{x}_{l} \}_{lÊ\in [L]}$.
Furthermore, it is straightforward to verify that $R_{\mathcal{M}}(\mathbf{x}_{1}, \mathbf{x}_{2}) = R_{\mathcal{M}}(\mathbf{x}_{2}, \mathbf{x}_{1})$ 
for any pair $\mathbf{x}_{1}, \mathbf{x}_{2}Ê\inÊ\mathcal{X}$. Thus, the two axioms of Definition \ref{equ_def_kernel_function} are satisfied.

The identity \eqref{equ_kernel_one_arg_x_0_equal_1} follows from the symmetry of $R_{\mathcal{M}}(\cdot,\cdot)$ and 
\begin{align}
 R_{\mathcal{M}}(\mathbf{x}, \mathbf{x}_{0}) =  \int_{\mathbf{y}}  \frac{ f(\mathbf{y};\mathbf{x})f(\mathbf{y}; \mathbf{x}_{0})}{f(\mathbf{y}; \mathbf{x}_{0})}  d \mathbf{y} =
\int_{\mathbf{y}}  f(\mathbf{y};\mathbf{x})  d \mathbf{y} =1.
\end{align} 
\end{proof} 

We are now in the position to make 
\begin{definition}
\label{def_RKHS_min_var_problem}
Given a minimum variance problem $\mathcal{M}=\minvarproblemscalar$ associated with the estimation problem $\mathcal{E}=\scalarestproblem$, we associate with it the RKHS $\mathcal{H}(\mathcal{M})$ that is 
given via the kernel $R_{\mathcal{M}}(\cdot,\cdot): \mathcal{X} \times \mathcal{X} \rightarrow \mathbb{R}$ and Theorem \ref{thm_constr_RKHS_closure_linear_span}.
\end{definition} 
Thus, the RKHS $\mathcal{H}(\mathcal{M})$ is given by $\mathcal{H}(\mathcal{M}) = \mathcal{L}(R_{\mathcal{M}}) \cup \mathcal{C}$ where $\mathcal{L}(R_{\mathcal{M}})$ is the inner-product space given by Definition \ref{def_linear_space_kernel} and $\mathcal{C}$ denotes the set of functions $f(\cdot): \mathcal{X} \rightarrow \mathbb{R}$, that 
are pointwise limits of a Cauchy sequence in the inner-product space $\mathcal{L}(R_{\mathcal{M}})$. 

\subsection{Main Result}
\label{sec_main_result_RKHS_MVE}

Besides the RKHS $\mathcal{H}(\mathcal{M})$, a second Hilbert space naturally associated to a minimum variance problem $\mathcal{M}$ is 
the Hilbert space $\mathcal{L}(\mathcal{M})$ defined in \eqref{equ_def_hilbert_space_rhos}.
Since both Hilbert spaces, $\mathcal{H}(\mathcal{M})$ and $\mathcal{L}(\mathcal{M})$, are associated with the same minimum variance problem $\mathcal{M}$, 
it does not come as a big surprise that these two Hilbert spaces are closely related to each other. In particular, we have \cite{Parzen59,Duttweiler73b}
\begin{theorem}
\label{thm_isometry_RKHS_rhos}
Consider a minimum variance problem $\mathcal{M}$ with the associated Hilbert spaces $\mathcal{H}(\mathcal{M})$ and $\mathcal{L}(\mathcal{M})$ 
given by Definition \ref{def_RKHS_min_var_problem} and \eqref{equ_def_hilbert_space_rhos} respectively. We then have that these two Hilbert spaces are isometric.
Moreover, there exists a unique congruence $\mathsf{J}[\cdot]: \mathcal{H}(\mathcal{M}) \rightarrow \mathcal{L}(\mathcal{M})$ from $\mathcal{H}(\mathcal{M})$ to $\mathcal{L}(\mathcal{M})$ that satisfies 
\begin{equation} 
\label{equ_def_congruence_H_L_MVP}
\mathsf{J}[R_{\mathcal{M}}(\cdot, \mathbf{x})] = \rho_{\mathcal{M}}(\cdot, \mathbf{x})
\end{equation} 
for every $\mathbf{x} \in \mathcal{X}$. 
Finally, if two elements $f(\cdot) \in \mathcal{H}(\mathcal{M})$ and $g(\cdot) \in \mathcal{L}(\mathcal{M})$ satisfy $\mathsf{J}[f(\cdot)] = g(\cdot)$, we have the identity 
\begin{equation}
\label{equ_expect_inversion_J_isometry_RKHS_rhos}
f(\mathbf{x}) = \mathsf{E}_{\mathbf{x}} \{ g(\mathbf{y})Ê\} 
\end{equation}  
for every $\mathbf{x} \in \mathcal{X}$. 
\end{theorem} 
\begin{proof}
\cite{Parzen59}
\end{proof} 

If a minimum variance problem $\mathcal{M}$ is such that the associated kernel $R_{\mathcal{M}}(\cdot, \cdot)$ is differentiable, then one 
can show 
\begin{theorem}
\label{thm_isometry_RKHS_rhos_derivative_kernel}
Consider a minimum variance problem $\mathcal{M}$ with parameter set $\mathcal{X} \subseteq \mathbb{R}^{N}$ that is such 
that the associated kernel $R_{\mathcal{M}}(\cdot, \cdot): \mathcal{X} \times \mathcal{X} \rightarrow \mathbb{R}$ 
is differentiable up to order $m$ (see Definition \ref{def_differentiable_kernel}). 
We denote by $\mathsf{J}[\cdot]$ the congruence from $\mathcal{H}(\mathcal{M})$ to $\mathcal{L}(\mathcal{M})$ defined by \eqref{equ_def_congruence_H_L_MVP}. 
Then, given a matrix $\mathbf{A} \in \mathbb{R}^{N \times D}$, a vector $\mathbf{z}_{c} \in \mathbb{R}^{D}$, and a multi-index $\mathbf{p} \in \mathbb{Z}_{+}^{D}$ with $\|\mathbf{p}\|_{\infty}\leq m$ that 
satisfy $\mathbf{A}\mathcal{N}^{\supp(\mathbf{p})}_{\mathbf{x}_{c}}(\varepsilon) \subseteq \mathcal{X}$ for a sufficiently small $\varepsilon >0$, we have 
that the image $\mathsf{J}[f(\cdot)]$ of the function 
\begin{equation}
\label{equ_def_func_RKHS_min_var_problem_part_der_congruence} 
f(\cdot) \triangleq  \frac{\partial^{\mathbf{p}}  \big( R_{\mathcal{M}}(\cdot, \mathbf{A} \mathbf{z}_{2}) h(\mathbf{z}_{2}) \big) }{\partial \mathbf{z}_{2}^{\mathbf{p}}} \bigg|_{\mathbf{z}_{2} = \mathbf{z}_{c}} \in \mathcal{H}(\mathcal{M})
\end{equation}
is given by 
\begin{equation} 
\label{equ_isometry_RKHS_rhos_derivative_kernel}
\mathsf{J}[f(\cdot)] =\frac{\partial^{\mathbf{p}}  \big( \rho_{\mathcal{M}}(\cdot,\mathbf{A} \mathbf{z})h(\mathbf{z}) \big) }{\partial \mathbf{z}^{\mathbf{p}}} \bigg|_{\mathbf{z} = \mathbf{z}_{c}}.
\end{equation}
Here, the function $h(\cdot): \mathbb{R}^{D} \rightarrow \mathbb{R}$ is arbitrary except that it is assumed that its partial derivatives of any order $\mathbf{p} \in \mathbb{Z}^{D}_{+}$ with $\| \mathbf{p} \|_{\infty} \leq m$ exist at $\mathbf{z}_{c}$.  
\end{theorem} 

\begin{proof}
This statement can be shown by induction on $\mathbf{p}$. We restrict ourselves to the proof for the case $\mathbf{p} = \mathbf{e}_{k}$ where $k \in [D]$. 
We can then write the function $f(\cdot)$ defined in \eqref{equ_def_func_RKHS_min_var_problem_part_der_congruence} as the limit of the Cauchy sequence 
\begin{equation*}
\bigg\{ f_{l}(\cdot) \triangleq \frac{1}{l}  \bigg[ R_{\mathcal{M}} \bigg(\cdot, \mathbf{A} \bigg( \mathbf{z}_{c} + \frac{1}{l} \mathbf{e}_{k} \bigg) \bigg) h \bigg(\mathbf{z}_{c} + \frac{1}{l} \mathbf{e}_{k} \bigg)  - R_{\mathcal{M}}(\cdot, \mathbf{A} \mathbf{z}_{c}) h(\mathbf{z}_{c}) \bigg] \bigg\}_{l \rightarrow \infty},
\end{equation*} 
i.e., $f(\cdot) = \lim_{l \rightarrow \infty} f_{l}(\cdot)$.  
Since, almost trivially by definition (cf.\ Definition \ref{def_isometry_congruence}), any congruence is a continuous mapping (in particular, the limit operation commutes with the application of a congruence) we have 
that $\mathsf{J}[f(\cdot)] = \lim_{l \rightarrow \infty} \mathsf{J}[f_{l}(\cdot)]$, which by Theorem \ref{thm_isometry_RKHS_rhos} yields 
\begin{align} 
\mathsf{J}[f(\cdot)] & = \lim_{l \rightarrow \infty} \frac{1}{l}  \Bigg[  \rho_{\mathcal{M}}\bigg(\, \cdot \, ,\mathbf{A} \bigg( \mathbf{z}_{c} + \frac{1}{l} \mathbf{e}_{k} \bigg)Ê\bigg) h\bigg(\mathbf{z}_{c} + \frac{1}{l} \mathbf{e}_{k} \bigg)  - \rho_{\mathcal{M}}(\cdot,\mathbf{A} \mathbf{z}_{c})h(\mathbf{z}_{c}) \Bigg] \nonumber \\[3mm]
& =\frac{\partial^{\mathbf{e}_{k}}  \big( \rho_{\mathcal{M}} (\cdot,\mathbf{A} \mathbf{z}) h(\mathbf{z}) \big) }{\partial \mathbf{z}^{\mathbf{e}_{k}}} \bigg|_{\mathbf{z} = \mathbf{z}_{c}}.
\end{align}
\end{proof}

The connection between the theory of RKHS and minimum variance estimation is summarized in \cite{Parzen59,Duttweiler73b}
\begin{theorem} 
\label{thm_main_facts_RKHS_MVE} 
Given a minimum variance problem $\mathcal{M}=\minvarproblemscalar$ associated with the estimation problem $\mathcal{E}=\scalarestproblem$, 
satisfying \eqref{equ_finite_integral_RKHS_cond} such that the associated RKHS $\mathcal{H}(\mathcal{M})$ exists. Then, we have that:
\begin{itemize}
\item There exists at least one allowed estimator $\hat{g}(\cdot)$, i.e., the set $\mathcal{F}(\mathcal{M})$ defined in \eqref{equ_est_finite_var_prescr_bias} is nonempty, if and only if 
\begin{equation}
\label{equ_gamma_in_RKHS}
\gamma(\cdot) \in  \mathcal{H}(\mathcal{M})
\end{equation}  
where $\gamma(\cdot): \mathcal{X} \rightarrow \mathbb{R}: \gamma(\mathbf{x}) \triangleq c(\mathbf{x}) + g(\mathbf{x})$ denotes the prescribed mean function of the minimum variance problem $\mathcal{M}$. 
\item If \eqref{equ_gamma_in_RKHS} holds, then the minimum achievable variance $L_{\mathcal{M}}$ (see \eqref{equ_def_min_ach_var}) is given by 
\begin{equation}
\label{equ_squared_norm_min_achiev_var}
L_{\mathcal{M}} = \| \gamma(\cdot) \|_{\mathcal{H}(\mathcal{M})}^{2} - \big[ \gamma(\mathbf{x}_{0}) \big]^{2},
\end{equation} 
and the corresponding LMV estimator $\hat{g}^{(\mathbf{x}_{0})}(\cdot)$ (see Definition \ref{def_LMV}) is obtained as 
\begin{equation} 
\label{equ_lmv_estimator_general_congruence_L_M_RKHS}
\hat{g}^{(\mathbf{x}_{0})} (\cdot) = \mathsf{J}[ \gamma(\cdot) ],
\end{equation} 
where $\mathsf{J}[\cdot]: \mathcal{H}(\mathcal{M}) \rightarrow \mathcal{L}(\mathcal{M})$ denotes the congruence from $\mathcal{H}(\mathcal{M})$ to $\mathcal{L}(\mathcal{M})$ which is defined by \eqref{equ_def_congruence_H_L_MVP}.
\end{itemize}
\end{theorem} 
\begin{proof}
\cite{Parzen59} 
\end{proof} 

The relation \eqref{equ_squared_norm_min_achiev_var} together with the definition of the minimum achievable variance $L_{\mathcal{M}}$ (see \eqref{equ_def_min_ach_var}) 
gives us a powerful tool for the derivation of lower bounds on the 
estimator variance if the estimator's bias is known. Indeed, any lower bound on $\| \gamma(\cdot) \|_{\mathcal{H}(\mathcal{M})}^{2}$ induces via \eqref{equ_squared_norm_min_achiev_var} and \eqref{equ_def_min_ach_var} a 
lower bound on the variance $v(\hat{g}(\cdot); \mathbf{x}_{0})$ at $\mathbf{x}_{0}$ of any estimator $\hat{g}(\cdot)$ with bias equal to $c(\mathbf{x})$ 
(or equivalently with mean function equal to $\gamma(\mathbf{x}) = c(\mathbf{x}) + g(\mathbf{x})$). 

As already observed in \cite{Parzen59}, a specific class of lower bounds on the variance $v(\hat{g}(\cdot); \mathbf{x}_{0})$ can be readily obtained via Theorem \ref{thm_orthog_proj_ineq} by projecting $\gamma(\cdot) \in \mathcal{H}(\mathcal{M})$ 
on a subspace $\mathcal{U} \subseteq \mathcal{H}(R)$ of the RKHS $\mathcal{H}(R)$. 
Indeed, denoting by $\mathbf{P}_{\mathcal{U}} \gamma(\cdot)$ the orthogonal projection of the function $\gamma(\cdot) \in \mathcal{H}(\mathcal{M})$ on the subspace $\mathcal{U}$, 
we have by Theorem \ref{thm_orthog_proj_ineq} the lower bound 
\begin{equation} 
\label{equ_lower_bound_RKHS_projection}
\| \gamma(\cdot) \|_{\mathcal{H}(\mathcal{M})}^{2} \geq \| \mathbf{P}_{\mathcal{U}} \gamma(\cdot) \|_{\mathcal{H}(\mathcal{M})}^{2}. 
\end{equation}

\subsection{Transforming the Parameter Function}
A RKHS that is associated to a minimum variance problem has a specific property as stated in 
\begin{lemma}
\label{lem_aff_trafo_parameter_function_RKHS} 
Consider a minimum variance problem $\mathcal{M}=\minvarproblemscalar$ associated with the estimation problem $\mathcal{E}=\scalarestproblem$ and with prescribed bias $c(\cdot) \equiv 0$. 
Then, if the parameter function $g(\cdot)$ is estimable, i.e., $g(\cdot) \inÊ\mathcal{H}(\mathcal{M})$ (due to Theorem \ref{thm_main_facts_RKHS_MVE} and Definition \ref{def_estimable_par_function}), we have that any 
other parameter function $g'(\cdot) = g(\cdot) + a$ with $a \in \mathbb{R}$ is also estimable, i.e., $g'(\cdot) \in \mathcal{H}(\mathcal{M})$, and 
moreover we have that 
\begin{equation} 
\label{equ_aff_trafo_parameter_function_RKHS}
\| g' (\cdot) \|_{\mathcal{H}(R)}^{2} = \| g(\cdot) \|_{\mathcal{H}(R)}^{2} + 2a g(\mathbf{x}_{0})+a ^{2}.
\end{equation} 
The
minimum achievable variance remains unchanged after replacing the parameter function $g(\cdot)$ with $g'(\cdot)$, i.e., 
\begin{equation}
L_{\mathcal{M}} = L_{\mathcal{M}'}, 
\end{equation} 
where $\mathcal{M}'$ denotes the minimum variance problem that is obtained from $\mathcal{M}$ by changing the parameter function from $g(\cdot)$ to $g'(\cdot)$. 
\end{lemma}
 
\begin{proof} 
Consider the specific function $v(\cdot) = R_{\mathcal{M}}(\cdot, \mathbf{x}_{0})$ which obviously belongs to the RKHS $\mathcal{H}(\mathcal{M})$. 
By \eqref{equ_kernel_one_arg_x_0_equal_1}, we have that $v(\mathbf{x}) = 1$ for every $\mathbf{x} \in \mathcal{X}$, 
and in turn via the reproducing property \eqref{equ_reproduction_property} that $\langle v(\cdot), v(\cdot) \rangle_{\mathcal{H}(\mathcal{M})} = \langle v(\cdot), R(\cdot,\mathbf{x}_{0}) \rangle_{\mathcal{H}(\mathcal{M})} = v(\mathbf{x}_{0}) = 1$. 
Thus, we can represent $g'(\cdot) = g(\cdot)+a$ as $g'(\cdot) = g(\cdot) + a v(\cdot)$. 
The squared norm of $g'(\cdot)$ can then be developed as 
\begin{align} 
\label{equ_transforming_the_parameter_function_RKHS_proof_1}
\| g'(\cdot)  \|_{\mathcal{H}(M)}^{2} & = \| g(\cdot) + a v(\cdot) \|_{\mathcal{H}(R)}^{M} = \langle g(\cdot) + a v(\cdot),  g(\cdot) + a v(\cdot) \rangle_{\mathcal{H}(\mathcal{M})}  \nonumber \\[3mm]Ê
& = \| g(\cdot)  \|_{\mathcal{H}(M)}^{2} + 2a \langle g(\cdot) , v(\cdot) \rangle_{\mathcal{H}(\mathcal{M})}+ a^{2} \langle v(\cdot) , v(\cdot) \rangle_{\mathcal{H}(\mathcal{M})} \nonumber \\[3mm] 
& = \| g(\cdot)  \|_{\mathcal{H}(M)}^{2} + 2a \langle g(\cdot) , R(\cdot,\mathbf{x}_{0}) \rangle_{\mathcal{H}(\mathcal{M})}+ a^{2} \underbrace{\langle v(\cdot) , v(\cdot) \rangle_{\mathcal{H}(\mathcal{M})}}_{=1} \nonumber \\[-1mm] 
&   \stackrel{\eqref{equ_reproduction_property}}{=} \| g(\cdot)  \|_{\mathcal{H}(M)}^{2} + 2a g(\mathbf{x}_{0}) + a^{2}. 
\end{align}
Based on \eqref{equ_transforming_the_parameter_function_RKHS_proof_1}, and using \eqref{equ_squared_norm_min_achiev_var}, the minimum achievable variance for $\mathcal{M}'$ can be calculated as 
\begin{align}
L_{\mathcal{M}'} & = \| g'(\cdot) \|_{\mathcal{H}(\mathcal{M}')}^{2} - (g'(\mathbf{x}_{0}))^{2} \nonumber \\[3mm]
& =  \| g'(\cdot) \|_{\mathcal{H}(\mathcal{M})}^{2} - (g'(\mathbf{x}_{0}))^{2}  \nonumber \\[3mm]Ê
& \stackrel{\eqref{equ_transforming_the_parameter_function_RKHS_proof_1}}{=}   \| g(\cdot)  \|_{\mathcal{H}(\mathcal{M})}^{2} + 2a g(\mathbf{x}_{0}) + a^{2} - (g(\mathbf{x}_{0})+a)^{2}  \nonumber \\[3mm]
& =  \| g(\cdot)  \|_{\mathcal{H}(\mathcal{M})}^{2} + 2a g(\mathbf{x}_{0}) + a^{2} - (g(\mathbf{x}_{0}))^{2} - 2a g(\mathbf{x}_{0}) - a^{2} \nonumber \\[3mm] 
& =   \| g(\cdot)  \|_{\mathcal{H}(\mathcal{M})}^{2} - (g(\mathbf{x}_{0}))^{2} = L_{\mathcal{M}},
\end{align}
where we used the fact, that the two kernels $R_{\mathcal{M}}(\cdot,\cdot)$, $R_{\mathcal{M}'}(\cdot,\cdot)$ and therefore also the two RKHSs 
$\mathcal{H}(\mathcal{M})$ and $\mathcal{H}(\mathcal{M}')$ coincide. This follows from \eqref{equ_def_kernel_M}, since the parameter sets and statistical models of $\mathcal{M}$ and $\mathcal{M}'$ are identical.  
\end{proof} 
We note that Lemma \ref{lem_aff_trafo_parameter_function_RKHS} is the RKHS-based formulation of a special case of Theorem \ref{thm_aff_trafo_parameter_function}. 

\subsection{Reducing the Parameter Set} 
\label{sec_reducing_par_set_RKHS}
The effect of a reduction of the parameter set $\mathcal{X}$ of a minimum variance problem on its associated minimum achievable variance can be analyzed conveniently using the RKHS framework. Indeed, 
consider a minimum variance problem $\mathcal{M}=\minvarproblemscalar$ associated with the estimation problem $\mathcal{E} = \scalarestproblem$ and the modified minimum variance problem 
$\mathcal{M}'=\mathcal{M} \big|_{\mathcal{X}'}$ obtained for the new parameter set $\mathcal{X}' \subseteq \mathcal{X}$. 

We then have  that 
\begin{equation}
\vspace*{2mm}
\label{equ_min_var_problem_reduc_para_set_RKHS}
L_{\mathcal{M}} \stackrel{(a)}{=} \| \gamma(\cdot) \|^{2}_{\mathcal{H}(\mathcal{M})} - \left( \gamma(\mathbf{x}_{0}) \right)^{2} \stackrel{(b)}{\geq}  \| \gamma(\cdot)\big|_{\mathcal{X}'} \|^{2}_{\mathcal{H}(\mathcal{M}')} - \left( \gamma(\mathbf{x}_{0}) \right)^{2} = L_{\mathcal{M}'}, 
\end{equation} 
\vspace*{2mm}
where step $(a)$ follows from Theorem \ref{thm_main_facts_RKHS_MVE}, step $(b)$ is due to Theorem \ref{thm_reducing_domain_RKHS} and $\gamma(\cdot): \mathcal{X} \rightarrow \mathbb{R}: \gamma(\mathbf{x}) = c(\mathbf{x}) + g(\mathbf{x})$ is the prescribed mean function of $\mathcal{M}$. 
The main aspect of \eqref{equ_min_var_problem_reduc_para_set_RKHS} (which was already mentioned in Section \ref{sec_effect_red_par_set}) is that a reduction of the parameter set $\mathcal{X}$ 
can never result in a worse achievable performance, i.e., in a higher minimum achievable variance. 

In Chapter \ref{chap_SLM} and Chapter \ref{chap_SCM}, we will study two specific estimation problems and the associated minimum variance problems. 
These estimation problems are obtained from well-known estimation problems 
by a reduction of the original parameter set $\mathcal{X} = \mathbb{R}^{N}$. This reduction will be quite substantial, since the reduced parameter sets 
will have measure zero w.r.t.\ to the Lebesgue measure in $\mathbb{R}^{N}$. 
However, using RKHS theory, we will show that this large reduction is necessary in order to obtain strictly lower values for the minimum achievable variance (see Section \ref{sec_strict_sparstiy_SLM} and \ref{sec_unbiasecd_SDPCM}). 
This agrees with the observation that pure or strict inequality constraints on the parameter set $\mathcal{X}$ of an estimation problem have no influence on the \CRBfull when evaluated at a point in the interior of $\mathcal{X}$ \cite[p. 1292]{GormanHero}. 

\subsection{Classes of Minimum Variance Problems with Continuously Varying Kernels} 
\label{sec_classes_min_var_problems_con_var_kernel}

\begin{figure}
\vspace{-1mm}
\centering
\psfrag{x0}[c][c][.9]{\uput{1mm}[270]{0}{$\mathbf{x}_{0}$}}
\psfrag{x}[c][c][.9]{\uput{0mm}[270]{0}{$\mathbf{x}$}}
\psfrag{fx}[c][c][.9]{\uput{1mm}[0]{0}{$f(\mathbf{x})$}}
\centering
\hspace*{1.5mm}\includegraphics[height=5.4cm,width=8.7cm]{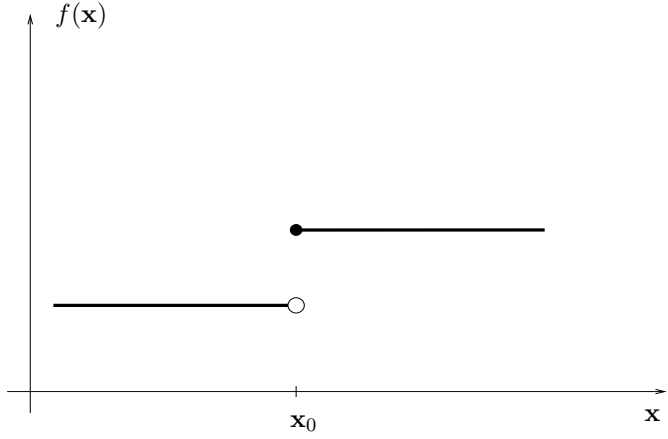}
\caption{Graph of a function that is lower semi-continuous at $\mathbf{x}_{0}$. The solid dot indicates the function value $f(\mathbf{x}_{0})$.} 
\label{fig_lower_upper_semi_cont}
\vspace*{-1.5mm}
\end{figure}

In some applications, it is of interest to characterize not only a single minimum variance problem $\mathcal{M}=\minvarproblemscalar$ associated with the estimation problem $\mathcal{E}=\scalarestproblem$, 
but rather to characterize the global properties of a whole class of minimum variance problems $\left\{ \mathcal{M}(\mathbf{x}_{0}) \right\}_{\mathbf{x}_{0} \in \mathcal{X}}$. This class is obtained 
by varying the parameter vector $\mathbf{x}_{0} \inÊ\mathcal{X}$ used for the definition of the minimum variance problem. 
Therefore, let us consider an estimation problem $\mathcal{E}=\scalarestproblem$ with a parameter set $\mathcal{X} \subseteq \mathbb{R}^{N}$ that is such that the kernel 
$R_{\mathcal{M}(\mathbf{x}_{0})}(\cdot,\cdot): \mathcal{X} \times \mathcal{X} \rightarrow \mathbb{R}$ 
is pointwise continuous w.r.t.\ the parameter $\mathbf{x}_{0}$, i.e., 
\begin{equation} 
\label{equ_def_cont_varying_kernel}
\lim\limits_{\mathbf{x}'_{0} \rightarrow \mathbf{x}_{0}} R_{\mathcal{M}(\mathbf{x}'_{0})}(\mathbf{x}_{1}, \mathbf{x}_{2}) = R_{\mathcal{M}(\mathbf{x}_{0})}(\mathbf{x}_{1}, \mathbf{x}_{2})\mbox{, } \quad  \forall \mathbf{x}_{0} ,\mathbf{x}_{1}, \mathbf{x}_{2} \in \mathcal{X}.
\end{equation}
In this case, the dependence of the minimum achievable variance $L_{\mathcal{M}(\mathbf{x}_{0})}$ 
on $\mathbf{x}_{0} \in \mathcal{X}$ is characterized by 
\begin{theorem} 
\label{thm_lower_semi_cont_varying_kernel}
Consider an estimation problem $\mathcal{E}= \scalarestproblem $ with parameter set $\mathcal{X} \subseteq \mathbb{R}^{N}$, and the class of minimum variance problems given by 
$\big \{ \mathcal{M}(\mathbf{x}_{0}) \triangleq \minvarproblemscalar \big \}_{\mathbf{x}_{0} \in \mathcal{X}}$ with a prescribed bias function $c(\cdot): \mathcal{X} \rightarrow \mathbb{R}$ 
that is valid for every $\mathcal{M}(\mathbf{x}_{0})$ with $\mathbf{x}_{0} \in \mathcal{X}$. 
We denote by $\gamma(\cdot): \mathcal{X} \rightarrow \mathbb{R} : \gamma(\mathbf{x}) = c(\mathbf{x}) + g(\mathbf{x})$ the prescribed mean function. 
Then, if $\gamma(\mathbf{x})$ is a continuous function of $\mathbf{x}$ and the kernel $R_{\mathcal{M}(\mathbf{x})}$ is pointwise continuous w.r.t.\ $\mathbf{x}$, i.e., 
it satisfies \eqref{equ_def_cont_varying_kernel}, the minimum achievable variance $L_{\mathcal{M} (\mathbf{x})}$ exists, i.e.\ is finite, for every $\mathbf{x}$.
Furthermore, viewed as a function of $\mathbf{x}$, $L_{\mathcal{M} (\mathbf{x})}$ is lower semi-continuous (see Figure \ref{fig_lower_upper_semi_cont}). 
\end{theorem}

\begin{proof} 
The finiteness of $L_{\mathcal{M} (\mathbf{x})}$ for every $\mathbf{x} \in \mathcal{X}$ follows trivially from the definition of a valid bias function (see Definition \ref{def_valid_bias_func_classic_est}) and Theorem \ref{thm_main_facts_RKHS_MVE}, 
since it is assumed that the prescribed bias function $c(\cdot)$ is valid for $\mathcal{M}(\mathbf{x})$ for every $\mathbf{x} \in \mathcal{X}$. 

Now, observe that 
\begin{align}
\label{equ_proof_cont_kernel_finite_approx_sup_1}
L_{\mathcal{M}(\mathbf{x})} & \stackrel{(a)}{=} \| \gamma(\cdot) \|_{\mathcal{H}(\mathcal{M})}^{2} - \left( \gamma(\mathbf{x}) \right)^{2} \nonumber \\[3mm]
   & \stackrel{(b)}{=}   \sup_{ \substack{f(\cdot) \in \mathcal{L}(R_{\mathcal{M}}) \\  \| f(\cdot) \|^{2}_{\mathcal{H}(\mathcal{M})} > 0}} 
\frac{\big( \langle  \gamma(\cdot), f(\cdot) \rangle_{\mathcal{H}(\mathcal{M})} \big)^{2}}{ \| f(\cdot) \|^{2}_{\mathcal{H}(\mathcal{M})} } - \left( \gamma(\mathbf{x}) \right)^{2}\nonumber \\[3mm]
   & \stackrel{(c)}{=}   \sup_{\substack{\mathcal{D}=\{\mathbf{x}_{1},\ldots,\mathbf{x}_{L}\} \\ L\in \mathbb{N}\mbox{, } \mathbf{x}_{l} \in \mathcal{X} \\ \mathbf{a} \in \mathcal{A}_{\mathcal{D}}}} 
\frac{\big( \langle  \gamma(\cdot), \sum_{l \in [L]} a_{l} R(\cdot,\mathbf{x}_{l}) \rangle_{\mathcal{H}(\mathcal{M})} \big)^{2}}{ \sum_{l,l'  \in [L]} a_{l} a_{l'} R_{\mathcal{M}(\mathbf{x})}(\mathbf{x}_{l}, \mathbf{x}_{l'})} - \left( \gamma(\mathbf{x}) \right)^{2}\nonumber \\[3mm]
   & =   \sup_{\substack{\mathcal{D}=\{\mathbf{x}_{1},\ldots,\mathbf{x}_{L}\} \\ L\in \mathbb{N}\mbox{, } \mathbf{x}_{l} \in \mathcal{X} \\ \mathbf{a} \in \mathcal{A}_{\mathcal{D}}}} 
\frac{\big( \sum_{l \in [L]} a_{l} \gamma(\mathbf{x}_{l}) \big)^{2}}{ \sum_{l,l'  \in [L]} a_{l} a_{l'} R_{\mathcal{M}(\mathbf{x})}(\mathbf{x}_{l}, \mathbf{x}_{l'})} - \left( \gamma(\mathbf{x}) \right)^{2}, 
\vspace*{4mm}
\end{align}
where we used $\mathcal{A}_{\mathcal{D}} \triangleq \big\{ \mathbf{a} \in \mathbb{R}^{L} \big|  \sum_{l,l'  \in [L]} a_{l} a_{l'} R_{\mathcal{M}(\mathbf{x})}(\mathbf{x}_{l}, \mathbf{x}_{l'}) >0 \big\}$.
The step $(a)$ follows from Theorem \ref{thm_main_facts_RKHS_MVE}, step $(b)$ is due to Theorem \ref{thm_approx_norm_dense_set} since the linear space $\mathcal{L}(R_{\mathcal{M}})$ 
(cf.\ Definition \ref{def_linear_space_kernel}) is dense in the RKHS $\mathcal{H}(\mathcal{M})$ by Theorem \ref{thm_constr_RKHS_closure_linear_span}. 
The step $(c)$ follows from the fact that any function $f(\cdot)$ belonging to the linear space $\mathcal{L}(R_{\mathcal{M}})$ can be written as a finite linear combination $f(\cdot) =  \sum_{l \in [L]} a_{l} R(\cdot,\mathbf{x}_{l})$ with 
suitable coefficients $a_{l} \in \mathbb{R}$ and points $\mathcal{D} = \{ \mathbf{x}_{1}, \ldots, \mathbf{x}_{L} \} \subseteq \mathcal{X}$, which can be used in combination with the reproducing property \eqref{equ_reproduction_property} to 
express the squared norm of $f(\cdot)$ as $ \| f(\cdot) \|^{2}_{\mathcal{H}(\mathcal{M})} =  \sum_{l,l'  \in [L]} a_{l} a_{l'} R_{\mathcal{M}(\mathbf{x})}(\mathbf{x}_{l}, \mathbf{x}_{l'})$. 
By \eqref{equ_proof_cont_kernel_finite_approx_sup_1}, we have 
\begin{equation} 
\label{equ_proof_cont_kernel_finite_approx_sup_4}
 L_{\mathcal{M}(\mathbf{x})} =  \sup_{\substack{\mathcal{D}=\{\mathbf{x}_{1},\ldots,\mathbf{x}_{L}\} \\ L\in \mathbb{N} \mbox{, } \mathbf{x}_{l} \in \mathcal{X} \\ 
 \mathbf{a} \in \mathcal{A}_{\mathcal{D}}}}h_{\mathcal{D},\mathbf{a}}(\mathbf{x}),
\end{equation} 
with the function $h_{\mathcal{D},\mathbf{a}}(\mathbf{x}): \mathcal{X} \rightarrow \mathbb{R}$ defined as 
\begin{equation} 
h_{\mathcal{D}, \mathbf{a}}(\mathbf{x}) \triangleq  \frac{\big( \sum_{l \in [L]} a_{l} \gamma(\mathbf{x}_{l}) \big)^{2}}{ \sum_{l,l'  \in [L]} a_{l} a_{l'} R_{\mathcal{M}(\mathbf{x})}(\mathbf{x}_{l}, \mathbf{x}_{l'})} 
- \left( \gamma(\mathbf{x}) \right)^{2},
\end{equation} 
where $L \triangleq | \mathcal{D}|$. 
For any finite set $\mathcal{D} = \{ \mathbf{x}_{l} \in \mathcal{X} \}_{l \in [L]}$ and $\mathbf{a} \in \mathcal{A}_{\mathcal{D}}$, it follows from the continuity of $R_{\mathcal{M}(\mathbf{x})}(\cdot,\cdot)$ and $\gamma(\mathbf{x})$ w.r.t.\ $\mathbf{x}$, that 
the function $h_{\mathcal{D},\mathbf{a}}(\mathbf{x})$ is continuous in a neighborhood around any point $\mathbf{x}_{0} \in \mathcal{X}$. 
Thus, for any point $\mathbf{x}_{0} \in \mathcal{X}$, there exists a radius $\delta_{0}>0$ such that the function $h_{\mathcal{D},\mathbf{a}}(\mathbf{x})$ is continuous on $\mathcal{B}(\mathbf{x}_{0},\delta_{0}) \cap \mathcal{X}$. 

We will show by contradiction that the function $L_{\mathcal{M}(\mathbf{x})}$ given by \eqref{equ_proof_cont_kernel_finite_approx_sup_4} must be lower semi-continuous. 
To that end, let us assume that the minimum achievable variance $L_{\mathcal{M}(\mathbf{x})}$ is not lower-semicontinuous at a specific parameter vector $\mathbf{x}_{0} \in \mathcal{X}$, i.e., we have 
that $\liminf\limits_{\mathbf{x} \rightarrow \mathbf{x}_{0}} L_{\mathcal{M}(\mathbf{x})} \leq  L_{\mathcal{M}(\mathbf{x}_{0})}- \varepsilon_{0}$ with a specific $\varepsilon_{0} > 0$. 
This implies that for any radius $r >0$, there exists at least one parameter vector $\mathbf{x}' \in \mathcal{X} \cap  \mathcal{B}(\mathbf{x}_{0}, r)$ such that $L_{\mathcal{M}(\mathbf{x}')} < L_{\mathcal{M}(\mathbf{x}_{0})} - \varepsilon_{0}/2$, i..e, 
\begin{equation}
\label{equ_proof_lower_semi_equ_1}
\forall r>0: \exists \mathbf{x}' \in \mathcal{X} \cap \mathcal{B}(\mathbf{x}_{0}, r)Ê\Rightarrow L_{\mathcal{M}(\mathbf{x}')} < L_{\mathcal{M}(\mathbf{x}_{0})} - \varepsilon_{0}/2. 
\end{equation} 
Now, we have that due to \eqref{equ_proof_cont_kernel_finite_approx_sup_4} there must be a finite subset $\mathcal{D}_{0} \subseteq \mathcal{X}$ and a vector $\mathbf{a}_{0} \in \mathcal{A}_{\mathcal{D}_{0}}$ such that
\begin{equation} 
\label{equ_proof_lower_semi_equ_2}
h_{\mathcal{D}_{0},\mathbf{a}_{0}}(\mathbf{x}_{0}) \geq  L_{\mathcal{M}(\mathbf{x}_{0})} - \varepsilon_{0}/4,  
\end{equation} 
since otherwise $h_{\mathcal{D},\mathbf{a}}(\mathbf{x}) \leq L_{\mathcal{M}(\mathbf{x}_{0})} - \varepsilon/4$ for every $\mathcal{D}$, $\mathbf{a}$, 
which would imply via \eqref{equ_proof_cont_kernel_finite_approx_sup_4} the contradiction 
$L_{\mathcal{M}(\mathbf{x}_{0})} \leq L_{\mathcal{M}(\mathbf{x}_{0})} - \varepsilon/4 < L_{\mathcal{M}(\mathbf{x}_{0})}$. 
Furthermore, $h_{\mathcal{D}_{0},\mathbf{a}_{0}}(\mathbf{x})$ is continuous at $\mathbf{x}_{0}$ within $\mathcal{B}(\mathbf{x}_{0},\delta_{0}) \cap \mathcal{X}$, implying the 
existence of a small radius $r_{0} >0$ (where $r_{0} < \delta_{0}$) such that for any $\mathbf{x} \in \mathcal{X} \cap \mathcal{B}(\mathbf{x}_{0},r_{0})$ we have 
\begin{equation} 
\label{equ_proof_lower_semi_equ_3}
h_{\mathcal{D}_{0},\mathbf{a}_{0}}(\mathbf{x}) \geq h_{\mathcal{D}_{0},\mathbf{a}_{0}}(\mathbf{x}_{0}) - \varepsilon_{0}/4.
\end{equation}  

By combining \eqref{equ_proof_lower_semi_equ_2} and \eqref{equ_proof_lower_semi_equ_3}, we see that there exists a radius $r_{0} >0$ (with $r_{0} < \delta_{0}$) such that for any $\mathbf{x} \in \mathcal{X} \cap \mathcal{B}(\mathbf{x}_{0},r_{0})$ we have
\begin{equation} 
\label{equ_proof_lower_semi_equ_4}
h_{\mathcal{D}_{0},\mathbf{a}_{0}}(\mathbf{x}) \stackrel{\eqref{equ_proof_lower_semi_equ_3}}{\geq} h_{\mathcal{D}_{0},\mathbf{a}_{0}}(\mathbf{x}_{0}) - \varepsilon_{0}/4 \stackrel{\eqref{equ_proof_lower_semi_equ_2}}{\geq} L_{\mathcal{M}(\mathbf{x}_{0})} - \varepsilon_{0}/4 - \varepsilon_{0}/4=L_{\mathcal{M}(\mathbf{x}_{0})} - \varepsilon_{0}/2.
\end{equation}  
This lower bound on the specific function $h_{\mathcal{D}_{0},\mathbf{a}_{0}}(\mathbf{x})$ implies also a lower bound on the supremum in \eqref{equ_proof_cont_kernel_finite_approx_sup_4}, i.e., for those parameter vectors $\mathbf{x}$ for which 
the bound in \eqref{equ_proof_lower_semi_equ_4} is in force, we have simultaneously the lower bound 
\begin{equation} 
\label{equ_proof_lower_semi_equ_6}
 L_{\mathcal{M}(\mathbf{x})} \stackrel{\eqref{equ_proof_cont_kernel_finite_approx_sup_4}}{=}  \sup_{\substack{\mathcal{D}=\{\mathbf{x}_{1},\ldots,\mathbf{x}_{L}\} \\ L\in \mathbb{N} \mbox{, } \mathbf{x}_{l} \in \mathcal{X} \\ 
 \mathbf{a} \in \mathcal{A}_{\mathcal{D}}}}h_{\mathcal{D},\mathbf{a}}(\mathbf{x}) \geq h_{\mathcal{D}_{0},\mathbf{a}_{0}}(\mathbf{x}) \stackrel{\eqref{equ_proof_lower_semi_equ_4}}{\geq} L_{\mathcal{M}(\mathbf{x}_{0})} - \varepsilon_{0}/2.
\end{equation}  

Putting together the pieces, we have that there exists a radius $r_{0} >0$ (with $r_{0} < \delta_{0}$) 
such that for any $\mathbf{x} \in \mathcal{X} \cap \mathcal{B}(\mathbf{x}_{0},r_{0})$ we have $L_{\mathcal{M}(\mathbf{x})}  \geq L_{\mathcal{M}(\mathbf{x}_{0})} - \varepsilon_{0}/2$, i.e., 
\begin{equation} 
\label{equ_proof_lower_semi_equ_5}
\forall \mathbf{x} \in  \mathcal{X} \cap \mathcal{B}(\mathbf{x}_{0},r_{0}) \quad \Rightarrow \quad L_{\mathcal{M}(\mathbf{x})}  \geq L_{\mathcal{M}(\mathbf{x}_{0})} - \varepsilon_{0}/2.
\end{equation}  
Since this fact contradicts \eqref{equ_proof_lower_semi_equ_1}, we have shown that $L_{\mathcal{M}(\mathbf{x})}$ must be lower-semicontinuous at every point $\mathbf{x} \in \mathcal{X}$. 
\end{proof}

\subsection{Sufficient Statistics from the RKHS Viewpoint} 
\label{sec_suff_stat_RKHS} 

Consider a minimum variance problem $\mathcal{M} = \minvarproblemscalar$ associated with the estimation problem $\mathcal{E} = \scalarestproblem$. 
In some cases, the observation $\mathbf{y} \in \mathbb{R}^{M}$, whose statistic is given by the statistical model $f(\mathbf{y};\mathbf{x})$ of the estimation problem, carries 
irrelevant information and can therefore be compressed in some sense. 
This compression or extraction of useful information (i.e., useful for the problem of estimating $g(\mathbf{x})$) can 
be performed by applying a (possibly randomized) mapping $\mathbf{T}(\cdot): \mathbb{R}^{M} \rightarrow \mathbb{R}^{K}$ to the observation $\mathbf{y}$, yielding a 
modified observation $\mathbf{z} \in \mathbb{R}^{K}$, $\mathbf{z} = \mathbf{T}(\mathbf{y})$. A compression is achieved if $K < M$. We will make the weak 
technical assumption that the mapping $\mathbf{T}(\cdot)$ is such that the modified observation $\mathbf{z}= \mathbf{T}(\mathbf{y})$, which is a random vector, 
possesses a pdf \cite{papoulis,BillingsleyProbMeasure,AshProbMeasure}. 
For clarity, we will denote the evaluation of the pdf of the random vector $\mathbf{z} = \mathbf{T}(\mathbf{y})$ at the specific vector $\mathbf{z}'  \in \mathbb{R}^{K}$ by $f_{\mathbf{z}}(\mathbf{z}'; \mathbf{x})$. 

\begin{definition}[\hspace*{-1mm}\cite{LC}]
A modified observation $\mathbf{z}=\mathbf{T}(\mathbf{y}) \in \mathbb{R}^{K}$ is called 
a sufficient statistic for the estimation problem $\mathcal{E}=\scalarestproblem$ if the conditional distribution $f(\mathbf{y} \big| \mathbf{z};\mathbf{x})$ of $\mathbf{y}$ 
given $\mathbf{z}$ does not depend on $\mathbf{x}$, i.e., 
if 
\begin{equation}
f(\mathbf{y} \big| \mathbf{z}; \mathbf{x}) = f(\mathbf{y} \big| \mathbf{z}),
\end{equation} 
with a suitable (conditional) pdf $f(\mathbf{y} \big| \mathbf{z})$. 
\end{definition}  

The verification if a given modified observation is a sufficient statistic can be based on 
\begin{theorem}
\label{thm_sufficient_stat_factor}
Consider an estimation problem $\mathcal{E}=\scalarestproblem$ and a modified observation $\mathbf{z}=\mathbf{T}(\mathbf{y})$ for $\mathcal{E}$ for which there exists a pdf. 
The modified observation $\mathbf{z}$ is a sufficient statistic for $\mathcal{E}$ if and only if we can factor the pdf $f(\mathbf{y}; \mathbf{x})$ of the observation $\mathbf{y}$ as 
\begin{equation} 
\label{equ_suff_stat_factor} 
f(\mathbf{y}; \mathbf{x}) =  f_{\mathbf{z}}(\mathbf{T}(\mathbf{y}); \mathbf{x}) h(\mathbf{y}),
\end{equation}  
where $h(\cdot)$ is a nonnegative function which does not depend on $\mathbf{x}$ and $f_{\mathbf{z}}(\cdot; \mathbf{x})$ 
denotes the pdf of the random vector $\mathbf{z} = \mathbf{T}(\mathbf{y})$. 
\end{theorem}
\begin{proof}
\cite{LC,kay,IbragimovBook} 
\end{proof} 
We note that Theorem \ref{thm_sufficient_stat_factor} is a variant of a famous result in classical estimation theory known as the ``Neyman-Fisher factorization theorem'' \cite{kay}. 
However, the Neyman-Fisher factorization theorem is more general than Theorem \ref{thm_sufficient_stat_factor} since it does not require the existence of a 
well-defined pdf of the modified observation $\mathbf{z} = \mathbf{T}(\mathbf{y})$. 

Consider an estimation problem $\mathcal{E}=\scalarestproblem$ and a sufficient statistic $\mathbf{z}=\mathbf{T}(\mathbf{y})$, i.e., 
we can factor the statistical model $f(\mathbf{y}; \mathbf{x})$ as in \eqref{equ_suff_stat_factor}. 
We can then define a new estimation problem $\mathcal{E}'$ by using the sufficient statistic as the observation, i.e., we define $\mathcal{E}' \triangleq \left( \mathcal{X}, f(\mathbf{z}; \mathbf{x}), g(\cdot) \right)$. 
If we consider the two minimum variance problems $\mathcal{M}$, $\mathcal{M}'$ that are obtained from $\mathcal{E}$ and $\mathcal{E}'$ and a common prescribed bias function 
$c(\cdot): \mathcal{X} \rightarrow \mathbb{R}$ and parameter vector $\mathbf{x}_{0} \in \mathcal{X}$, we have that the associated RKHSs are identical: 
\begin{theorem} 
\label{thm_suff_stat_RKHS_invar}
Consider an estimation problem $\mathcal{E}=\scalarestproblem$ with a sufficient statistic $\mathbf{z} = \mathbf{T}(\mathbf{y})$ and the modified estimation problem 
$\mathcal{E}'=\left(\mathcal{X},f(\mathbf{z}; \mathbf{x}),g(\cdot) \right)$ that is obtained from $\mathcal{E}$ by using 
the sufficient statistic $\mathbf{z}$ as the observation. 
Then we have for any fixed parameter vector $\mathbf{x}_{0} \in \mathcal{X}$ and any prescribed bias function $c(\cdot): \mathcal{X} \rightarrow \mathbb{R}$, 
that the RKHSs associated to the two minimum variance problems $\mathcal{M}=\left( \mathcal{E},c(\cdot), \mathbf{x}_{0} \right)$ and $\mathcal{M}'=\left( \mathcal{E}',c(\cdot), \mathbf{x}_{0} \right)$ are identical, i.e., 
\begin{equation} 
\mathcal{H}(\mathcal{M}) = \mathcal{H}(\mathcal{M}'). 
\end{equation}
In particular, this implies that the minimum achievable variances of 
$\mathcal{M}$ and $\mathcal{M}'$ are equal, i.e., 
\begin{equation}
\label{equ_suff_stat_RKHS_invar_min_var_same}
L_{\mathcal{M}} = L_{\mathcal{M}'}.  
\end{equation}
\end{theorem}

\begin{proof}
The statement is proven by showing that the kernels $R_{\mathcal{M}}(\cdot,\cdot): \mathcal{X} \times \mathcal{X} \rightarrow \mathbb{R}$ 
and $R_{\mathcal{M}'}(\cdot,\cdot): \mathcal{X} \times \mathcal{X} \rightarrow \mathbb{R}$ associated with $\mathcal{M}$ and $\mathcal{M}'$, respectively, are identical. 
Indeed, by denoting the value of the joint pdf and the conditional pdf of the 
two random vectors $\mathbf{y}$ and $\mathbf{z}$ evaluated at $\mathbf{z}=\mathbf{z}'$ by $f(\mathbf{y}, \mathbf{z}'; \mathbf{x})$ and $f(\mathbf{y} \big| \mathbf{z}'; \mathbf{x})$ respectively, 
we have by \eqref{equ_def_kernel_M} that
\vspace*{1mm}
\begin{align}
R_{\mathcal{M}}(\mathbf{x}_{1}, \mathbf{x}_{2}) & = 
\mathsf{E}_{\mathbf{x}_{0}} \bigg\{ \frac{ f(\mathbf{y}; \mathbf{x}_{1})}{ f (\mathbf{y}; \mathbf{x}_{0})} \frac{ f(\mathbf{y}; \mathbf{x}_{2})}{f (\mathbf{y}; \mathbf{x}_{0})}  \bigg \} 
 \stackrel{(a)}{=}  \mathsf{E}_{\mathbf{x}_{0}} \bigg\{  \frac{ f_{\mathbf{z}}(\mathbf{T}(\mathbf{y});  \mathbf{x}_{1})}{ f_{\mathbf{z}} (\mathbf{T}(\mathbf{y}); \mathbf{x}_{0})} \frac{ f_{\mathbf{z}}(\mathbf{T}(\mathbf{y}); \mathbf{x}_{2})}{f_{\mathbf{z}}(\mathbf{T}(\mathbf{y});\mathbf{x}_{0})}  \bigg \}\nonumber  \\[3mm] 
& \stackrel{(b)}{=}  \mathsf{E}_{\mathbf{x}_{0}} \bigg\{  \frac{ f_{\mathbf{z}}(\mathbf{z};  \mathbf{x}_{1})}{ f_{\mathbf{z}} (\mathbf{z}; \mathbf{x}_{0})} \frac{ f_{\mathbf{z}}(\mathbf{z}; \mathbf{x}_{2})}{f_{\mathbf{z}}(\mathbf{z};\mathbf{x}_{0})} \bigg \} 
= R_{\mathcal{M}'}(\mathbf{x}_{1}, \mathbf{x}_{2}). 
\end{align}
Here, the step $(a)$ follows from the factorization \eqref{equ_suff_stat_factor} and $(b)$ follows from a fundamental property of the expectation operation\footnote{For a random vector $\mathbf{y} \in \mathbb{R}^{M}$ 
and a function $\mathbf{t}(\cdot): \mathbb{R}^{M}Ê\rightarrow \mathbb{R}^{K}$, the expectation $\mathsf{E}_{\mathbf{x}_{0}} \{ \mathbf{t}(\mathbf{y}) \}$ is equal 
to the expectation $\mathsf{E}_{\mathbf{x}_{0}} \{\mathbf{b} \}$, where $\mathbf{b} \in \mathbb{R}^{K}$ is the random vector whose 
realization $\mathbf{b}'$ is obtained from $\mathbf{y}'$, the realization of the random vector $\mathbf{y}$, by $\mathbf{b}' = \mathbf{t}(\mathbf{y}')$.} \cite{papoulis,BillingsleyProbMeasure,AshProbMeasure}.
Since $R_{\mathcal{M}}(\cdot,\cdot) \!=\! R_{\mathcal{M}'}(\cdot,\cdot)$, we have by Theorem \ref{thm_exist_uniq_RKHS} that  $\mathcal{H}(\mathcal{M}) = \mathcal{H}(\mathcal{M}')$. 
Finally, the relation \eqref{equ_suff_stat_RKHS_invar_min_var_same} follows then by \eqref{equ_squared_norm_min_achiev_var} of Theorem \ref{thm_main_facts_RKHS_MVE}. 
\end{proof}

We already mentioned in Theorem \ref{thm_invariance_classical_est} of Section \ref{sec_invariance_mve} 
that any transformation of the observation, by an invertible map is irrelevant to minimum variance estimation. Indeed, as can be verified easily by Theorem \ref{thm_sufficient_stat_factor}, 
any modified observation that is obtained by an invertible map $\mathbf{T}(\cdot)$ is trivially a sufficient statistic, and hence Theorem \ref{thm_suff_stat_RKHS_invar} applies.

Furthermore, we note that Theorem \ref{thm_suff_stat_RKHS_invar} agrees with a famous result within classical (non-Bayesian) estimation theory, that is 
known as the ``Rao-Blackwell-Lehmann-Scheff{\'e}'' theorem \cite{kay}. This theorem states 
that given an estimator with a specific bias function, and given a sufficient statistic, one can always find another estimator, depending on the observation only via the sufficient statistic, 
with the same bias function and a variance that does not exceed the variance of the original estimator. 
This fact is consistent with Theorem \ref{thm_suff_stat_RKHS_invar}, in particular with \eqref{equ_suff_stat_RKHS_invar_min_var_same}.


\section{RKHS Interpretation of Known Variance Bounds}

In this section, we will interpret some well-known lower bounds on the estimator variance either as a specific evaluation of the squared norm $\| \gamma(\cdot) \|^{2}_{\mathcal{H}(\mathcal{M})}$ 
in \eqref{equ_squared_norm_min_achiev_var} or as specific lower bounds on $\| \gamma(\cdot) \|^{2}_{\mathcal{H}(\mathcal{M})}$. 

Throughout this section, we will 
consider only minimum variance problems with prescribed bias $c(\cdot) \equiv 0$, i..e, we consider unbiased estimation, which means that $\gamma(\cdot) = g(\cdot)$. 
However, by Section \ref{sec_equ_bias_param_function}, in particular by Theorem \ref{thm_equ_bias_param_function}, this is no real restriction since 
we do not constrain a priori the parameter functions $g(\cdot)$ of the considered estimation problems. 

\subsection{Barankin Bound} 
\label{sec_RKHS_interpretation_of_Barankin_bound}

The term ``Barankin bound'' refers to the following result \cite{Barankin49,mcaulay71}: 
\begin{theorem}[Barankin bound]
\label{thm_barankin_bound}
Consider a minimum variance problem $\mathcal{M}=\minvarproblemscalar$ associated with the estimation problem $\mathcal{E} = \scalarestproblem$, 
where $c(\mathbf{x}) \equiv 0$ and for which Postulate \ref{assumption_RKHS_minvarproblem} is fulfilled. 
Then, for any estimator $\hat{g}(\cdot): \mathbb{R}^{M} \rightarrow \mathbb{R}$ which is unbiased and has finite variance at $\mathbf{x}_{0}$, we have that its variance at $\mathbf{x}_{0}$ is lower bounded by 
\begin{equation}
\label{equ_barankin_bound_1}
v(\hat{g}(\cdot);  \mathbf{x}_{0}) 
\geq  \sup_{\substack{\mathcal{D}=\{\mathbf{x}_{1},\ldots,\mathbf{x}_{L}\} \\ L\in \mathbb{N}\mbox{,} \,\, \mathbf{x}_{l} \in \mathcal{X} \\ \mathbf{a} \in \mathcal{B}_{\mathcal{D}}}} 
\frac{ \bigg| \sum\limits_{l \in [L]} a_{l} h(\mathbf{x}_{l}) \bigg|^{2}}{ \mathsf{E}_{\mathbf{x}_{0}} \bigg\{ \bigg( \sum\limits_{l \in [L]} a_{l} \rho_{\mathcal{M}}(\mathbf{y},\mathbf{x}_{l}) \bigg)^{2} \bigg\}},
\end{equation}
where we used 
\begin{equation}
\label{equ_def_h_func_barankin_bound}
h(\cdot): \mathcal{X} \rightarrow \mathbb{R}: h(\mathbf{x}) \triangleq g(\mathbf{x}) - g(\mathbf{x}_{0}), 
\end{equation}
$\rho_{\mathcal{M}}(\mathbf{y},\mathbf{x}_{l})$ as defined in \eqref{equ_def_likelihood}, and $\mathcal{B}_{\mathcal{D}} \triangleq \big\{ \mathbf{a} \in \mathbb{R}^{L} | \mathsf{E}_{\mathbf{x}_{0}} \big\{ \big( \sum_{l \in [L]} a_{l} \rho_{\mathcal{M}}(\mathbf{y},\mathbf{x}_{l}) \big)^{2} \big\} > 0 \big\}$. 
Moreover, if there exists at least one unbiased estimator with a finite variance at $\mathbf{x}_{0}$, we have that the supremum in \eqref{equ_barankin_bound_1} is equal to the minimum achievable variance at $\mathbf{x}_{0}$, i.e., 
\begin{equation}
\label{equ_min_var_barankin_bound}
L_{\mathcal{M}} =  \sup_{\substack{\mathcal{D}=\{\mathbf{x}_{1},\ldots,\mathbf{x}_{L}\} \\ L\in \mathbb{N} \mbox{,} \,\, \mathbf{x}_{l} \in \mathcal{X} \\ \mathbf{a} \in \mathcal{B}_{\mathcal{D}}}} \frac{ \bigg| \sum\limits_{l \in [L]} a_{l} h(\mathbf{x}_{l}) \bigg|^{2}}{ \mathsf{E}_{\mathbf{x}_{0}} \bigg\{ \bigg( \sum\limits_{l \in [L]} a_{l} \rho_{\mathcal{M}}(\mathbf{y},\mathbf{x}_{l}) \bigg)^{2} \bigg\}},
\end{equation}
which is then called the Barankin bound for the minimum variance problem $\mathcal{M}$. 
\end{theorem} 

We note that the RKHS interpretation of the Barankin bound has been already discussed 
in \cite{GeomInterprBarankinBound}. However, the following derivation will be self-contained and more detailed than those presented in \cite{GeomInterprBarankinBound}. 

First, we have that by Theorem \ref{thm_main_facts_RKHS_MVE}, the existence of at least one unbiased estimator with finite variance at $\mathbf{x}_{0}$, i.e., the existence of at least one allowed estimator, is equivalent to 
the fact that the function $g(\cdot)$ (which 
coincides with the prescribed mean function $\gamma(\mathbf{x}) = g(\mathbf{x}) + c(\mathbf{x})$ for $\mathcal{M}$) belongs to the RKHS $\mathcal{H}(\mathcal{M})$, i.e., $g(\cdot) \in \mathcal{H}(\mathcal{M})$. 
However, since $h(\cdot) = g(\cdot) +a$ with $a = - g(\mathbf{x}_{0})$, we have by Lemma \ref{lem_aff_trafo_parameter_function_RKHS}, that $g(\cdot) \in \mathcal{H}(\mathcal{M})$ if and only if $h(\cdot) \in \mathcal{H}(\mathcal{M})$. 
Thus, there exists at least one unbiased estimator with finite variance at $\mathbf{x}_{0}$ if and only if $h(\cdot) \in \mathcal{H}(\mathcal{M})$. 
In what follows we will verify \eqref{equ_min_var_barankin_bound} in the case $h(\cdot) \in \mathcal{H}(\mathcal{M})$, which in turn shows that the supremum in \eqref{equ_barankin_bound_1} is finite 
if $h(\cdot) \in \mathcal{H}(\mathcal{M})$ (due to Theorem \ref{thm_main_facts_RKHS_MVE}), i.e., if there exists at least one unbiased estimator with finite variance at $\mathbf{x}_{0}$. 
The sufficiency of the condition $h(\cdot) \in \mathcal{H}(\mathcal{M})$ for the existence of a finite supremum in \eqref{equ_min_var_barankin_bound} is also mentioned in \cite[Lemma 4]{HeinRKHS2004}, where it is moreover shown that 
the condition $h(\cdot) \in \mathcal{H}(\mathcal{M})$ is also necessary for the existence of a finite supremum in \eqref{equ_min_var_barankin_bound}.

Let us now show that for the case $h(\cdot) \in \mathcal{H}(\mathcal{M})$ (which is equivalent to $g(\cdot) \in \mathcal{H}(\mathcal{M})$), the relation \eqref{equ_min_var_barankin_bound} of Theorem \ref{thm_barankin_bound} is
just a reformulation of \eqref{equ_squared_norm_min_achiev_var} in Theorem \ref{thm_main_facts_RKHS_MVE} which is based on Theorem \ref{thm_approx_norm_dense_set}. 
First, we note that by Theorem \ref{thm_main_facts_RKHS_MVE} and Lemma \ref{lem_aff_trafo_parameter_function_RKHS}, we have that 
\begin{align}
\label{equ_barankin_bound_affine_trafo}
L_{\mathcal{M}} & = \| g(\cdot) \|_{\mathcal{H}(\mathcal{M})}^{2} - ( g(\mathbf{x}_{0}))^{2} \stackrel{\eqref{equ_def_h_func_barankin_bound}}{=} \| h(\cdot)+ g(\mathbf{x}_{0}) \|_{\mathcal{H}(\mathcal{M})}^{2} -  ( g(\mathbf{x}_{0}))^{2}  \nonumber \\[4mm]
&\stackrel{\eqref{equ_aff_trafo_parameter_function_RKHS}}{=} \| h(\cdot) \|_{\mathcal{H}(\mathcal{M})}^{2} + 2 g(\mathbf{x}_{0}) \underbrace{h(\mathbf{x}_{0})}_{=0} +  ( g(\mathbf{x}_{0}))^{2}-  ( g(\mathbf{x}_{0}))^{2} =\| h(\cdot) \|_{\mathcal{H}(\mathcal{M})}^{2}  . 
\end{align}
By a similar derivation as used for \eqref{equ_proof_cont_kernel_finite_approx_sup_1}, we obtain 
\begin{align} 
\label{equ_barankin_bound_norm_as_sup}
\| h(\cdot) \|_{\mathcal{H}(\mathcal{M})}^{2} & =   \sup_{ \substack{f(\cdot) \in \mathcal{L}(R_{\mathcal{M}}) \\  \| f(\cdot) \|^{2}_{\mathcal{H}(\mathcal{M})} > 0}} 
\frac{\big( \langle  h(\cdot), f(\cdot) \rangle_{\mathcal{H}(\mathcal{M})} \big)^{2}}{ \| f(\cdot) \|^{2}_{\mathcal{H}(\mathcal{M})} } \nonumber \\[3mm]
   & =   \sup_{\substack{\mathcal{D}=\{\mathbf{x}_{1},\ldots,\mathbf{x}_{L}\} \\ L\in \mathbb{N}\mbox{, } \mathbf{x}_{l} \in \mathcal{X} \\ \mathbf{a} \in \mathcal{A}_{\mathcal{D}}}} 
\frac{\big( \langle  h(\cdot), \sum_{l \in [L]} a_{l} R(\cdot,\mathbf{x}_{l}) \rangle_{\mathcal{H}(\mathcal{M})} \big)^{2}}{ \sum_{l,l'  \in [L]} a_{l} a_{l'} R_{\mathcal{M}(\mathbf{x})}(\mathbf{x}_{l}, \mathbf{x}_{l'})} - \left( \gamma(\mathbf{x}) \right)^{2}\nonumber \\[3mm]
   & =   \sup_{\substack{\mathcal{D}=\{\mathbf{x}_{1},\ldots,\mathbf{x}_{L}\} \\ L\in \mathbb{N}\mbox{, } \mathbf{x}_{l} \in \mathcal{X} \\ \mathbf{a} \in \mathcal{A}_{\mathcal{D}}}} 
\frac{\big( \sum_{l \in [L]} a_{l} h(\mathbf{x}_{l}) \big)^{2}}{ \sum_{l,l'  \in [L]} a_{l} a_{l'} R_{\mathcal{M}(\mathbf{x})}(\mathbf{x}_{l}, \mathbf{x}_{l'})}, 
\end{align}
\vspace*{2mm}
where we used $\mathcal{A}_{\mathcal{D}} \triangleq \{ \mathbf{a} \in \mathbb{R}^{L} |  \sum_{l,l'  \in [L]} a_{l} a_{l'} R_{\mathcal{M}}(\mathbf{x}_{l}, \mathbf{x}_{l'}) >0 \}$. 
Combining \eqref{equ_barankin_bound_norm_as_sup} and \eqref{equ_barankin_bound_affine_trafo} yields 
\begin{equation}
L_{\mathcal{M}} =   \sup_{\substack{\mathcal{D}=\{\mathbf{x}_{1},\ldots,\mathbf{x}_{L}\} \\ L\in \mathbb{N} \mbox{,}\,\, \mathbf{x}_{l} \in \mathcal{X} \\ \mathbf{a} \in \mathcal{A}_{\mathcal{D}}}} 
\frac{\big| \sum_{l \in [L]} a_{l} h(\mathbf{x}_{l}) \big|^{2}}{ \sum_{l,l'  \in [L]} a_{l} a_{l'} R_{\mathcal{M}}(\mathbf{x}_{l}, \mathbf{x}_{l'})}.
\end{equation} 
But this is equal to \eqref{equ_min_var_barankin_bound} since for any finite $\mathcal{D} \subseteq \mathcal{X}$, of size $|\mathcal{D}| = L$, 
and arbitrary coefficient vector $\mathbf{a} \in \mathbb{R}^{L}$, we straightforwardly have that 
\begin{equation} 
 \sum_{k,k'  \in [L]} a_{k} a_{k'} R_{\mathcal{M}}(\mathbf{x}_{k}, \mathbf{x}_{k}') = 
\mathsf{E}_{\mathbf{x}_{0}} \bigg\{ \bigg( \sum\limits_{k \in [L]} a_{k} \rho_{\mathcal{M}}(\mathbf{y},\mathbf{x}_{k}) \bigg)^{2} \bigg\}.
\end{equation} 
This also implies that $\mathcal{A}_{\mathcal{D}} = \mathcal{B}_{\mathcal{D}}$. 
The lower bound in \eqref{equ_barankin_bound_1} follows then directly from  \eqref{equ_min_var_barankin_bound} and the definition of $L_{\mathcal{M}}$ in \eqref{equ_def_min_ach_var}.
 
\subsection{\CRBfull} 
\label{sec_CRB}
The \CRBfull (CRB) \cite{cramer45,rao45,kay} is probably the best known lower bound on the variance of estimators for a given estimation problem $\mathcal{E}=\scalarestproblem$ 
which satisfies some weak regularity conditions as discussed presently. 
Since the CRB applies to any unbiased estimator (the CRB applies also for biased estimators since it depends solely on the mean function of a given estimator) 
for $\mathcal{E}$, the CRB also yields a lower bound on the minimum achievable variance $L_{\mathcal{M}}$ 
for the minimum variance problem $\mathcal{M}=\left(\mathcal{E},c(\cdot)\equiv 0,\mathbf{x}_{0}\right)$ associated with the estimation problem $\mathcal{E}$.
We will consider the CRB for two different settings which are defined by the structure of the parameter set $\mathcal{X}$ 
associated with the estimation problem $\mathcal{E}$. These two instances of the CRB are termed the ``unconstrained CRB'' and the 
``constrained CRB'', respectively. 

The validity of the CRB requires some weak technical conditions to be placed on the estimation problem $\mathcal{E}=\scalarestproblem$ and associated 
minimum variance problem $\mathcal{M}=\left(\mathcal{E},c(\cdot)\equiv 0,\mathbf{x}_{0}\right)$, which are summarized in 
\begin{postulate}
\label{post_regular_cond_CRB}
The estimation problem $\mathcal{E}=\scalarestproblem$ and associated minimum variance problem  $\mathcal{M}=\left(\mathcal{E},c(\cdot) \equiv 0, \mathbf{x}_{0} \right)$ are such that:
\begin{itemize} 
\item The parameter set $\mathcal{X} \subseteq \mathbb{R}^{N}$ contains an $N$-dimensional open ball $\mathcal{B}(\mathbf{x}_{0},r)$ centered at $\mathbf{x}_{0} \in \mathcal{X}$ and 
with a positive radius $r >0$, i.e., 
\begin{equation}
\label{equ_unconstr_crb_param_set_contains_ball} 
\mathcal{B}(\mathbf{x}_{0},r) \subseteq \mathcal{X}. 
\end{equation} 
\item For every $\mathbf{x} \in \mathcal{X}$, the partial derivatives $ \frac{\partial^{\mathbf{p}} f(\mathbf{y}; \mathbf{x})}{\partial \mathbf{x}^{\mathbf{p}}} $ exist
and moreover 
\begin{equation} 
\label{equ_crb_regul_cond_finite_part_der}
\mathsf{E}_{\mathbf{x}_{0}} \bigg\{ \left( \frac{1}{f(\mathbf{y};\mathbf{x}_{0})} \frac{\partial^{\mathbf{p}} f(\mathbf{y}; \mathbf{x})}{\partial \mathbf{x}^{\mathbf{p}}}Ê\right)^{2} \bigg\} < \infty
\end{equation} 
for every multi-index $\mathbf{p}Ê \in \mathbb{Z}_{+}^{N}$ with $\| \mathbf{p} \|_{\infty} \leq m$, where $m \in \mathbb{N}$ is called the maximum differentiation order. 
\item For any function $h(\mathbf{y}): \mathbb{R}^{M} \rightarrow \mathbb{R}$ such that 
$\mathsf{E}_{\mathbf{x}} \{ h(\mathbf{y}) \}$ exists\footnote{Note that this implies that $\mathsf{E}_{\mathbf{x}} \{ | h(\mathbf{y})| \} < \infty$ \cite{AshProbMeasure,HalmosMeasure, BillingsleyProbMeasure}.} 
we have that for all $\mathbf{x} \in \mathcal{B}(\mathbf{x}_{0},r)$ the expectation operation commutes with partial differentiation, i.e.,  
\begin{equation} 
\label{equ_crb_regul_cond_intercah_diff_integr_intral_repr}
\frac{\partial^{\mathbf{p}} \int_{\mathbf{y}}  h(\mathbf{y}) f(\mathbf{y}; \mathbf{x}) d\mathbf{y}}{\partial \mathbf{x}^{\mathbf{p}}}  
= \int_{\mathbf{y}}  h(\mathbf{y}) \frac{\partial^{\mathbf{p}} f(\mathbf{y}; \mathbf{x})}{\partial \mathbf{x}^{\mathbf{p}}} d\mathbf{y}
\end{equation} 
or equivalently, expressed in terms of expectations, we have 
\begin{equation}
\label{equ_crb_regul_cond_intercah_diff_integr}
\frac{\partial^{\mathbf{p}} \mathsf{E}_{\mathbf{x}} \{  h(\mathbf{y}) \}}{\partial \mathbf{x}^{\mathbf{p}}}
= \mathsf{E}_{\mathbf{x}} \bigg\{  h(\mathbf{y})\frac{1}{f(\mathbf{y};\mathbf{x})} \frac{\partial^{\mathbf{p}} f(\mathbf{y}; \mathbf{x})}{\partial \mathbf{x}^{\mathbf{p}}}  \bigg\}
\end{equation}  
for every multi-index $\mathbf{p}Ê \in \mathbb{Z}_{+}^{N}$ with $\| \mathbf{p} \|_{\infty} \leq m$ provided that the right hand side of \eqref{equ_crb_regul_cond_intercah_diff_integr_intral_repr} and \eqref{equ_crb_regul_cond_intercah_diff_integr} is finite. 
\item The parameter function $g(\cdot): \mathcal{X} \rightarrow \mathbb{R}$ is such that the partial derivatives $\frac{\partial^{\mathbf{p}} g(\mathbf{x})}{\partial \mathbf{x}^{\mathbf{p}}} \big|_{\mathbf{x} = \mathbf{x}_{0}}$ exist 
at $\mathbf{x}_{0}$ for every multi-index $\mathbf{p}Ê \in \mathbb{Z}_{+}^{N}$ with $\| \mathbf{p} \|_{\infty} \leq m$.
\item The expectation 
\begin{equation}
\label{equ_CRB_regul_cond_inner_prod_cont_x_1_x_2}
 \mathsf{E}_{\mathbf{x}_{0}} \left\{ \frac{1}{f(\mathbf{y}; \mathbf{x}_{0})}  \frac{\partial^{\mathbf{p}_{1}}f(\mathbf{y};\mathbf{x}_{1})} {\partial \mathbf{x}_{1}^{\mathbf{p}_{1}} }   \frac{1}{f(\mathbf{y}; \mathbf{x}_{0})} \frac{\partial^{\mathbf{p}_{2}} f(\mathbf{y};\mathbf{x}_{2})} {\partial \mathbf{x}_{2}^{\mathbf{p}_{2}} }  \right\}
\end{equation} 
depends continuously on the parameter vectors 
$\mathbf{x}_1, \mathbf{x}_{2} \in \mathcal{X} \cap \mathcal{B}(\mathbf{x}_{0}, r)$ for every pair of multi-indices $\mathbf{p}_{1}, \mathbf{p}_{2}Ê \in \mathbb{Z}_{+}^{N}$ with $\| \mathbf{p}_{1} \|_{\infty},\| \mathbf{p}_{2} \|_{\infty} \leq m$.
\end{itemize}
\end{postulate} 
It can be shown that Postulate \ref{post_regular_cond_CRB} is satisfied in particular for estimation problems where the statistical model is given by an exponential family (see Section \ref{sec_exp_family}). 

In the following, we will need
\begin{theorem}
\label{thm_regul_cond_CRB_imply_diff}
Consider an estimation problem $\mathcal{E}=\scalarestproblem$ and associated minimum variance problem  $\mathcal{M}=\left(\mathcal{E},c(\cdot) \equiv 0, \mathbf{x}_{0} \right)$ 
for which Postulate \ref{assumption_RKHS_minvarproblem} and Postulate \ref{post_regular_cond_CRB} with maximum differentiation order $m$, and radius $r$ hold true.
Then, the kernel $R_{\mathcal{M}}(\cdot, \cdot): \mathcal{X} \times \mathcal{X} \rightarrow \mathbb{R}$ is differentiable up to order $m$. 
\end{theorem}

\begin{proof}
First, note that the condition \eqref{equ_crb_regul_cond_finite_part_der} implies that 
\begin{equation}
\label{equ_finite_value_expect_implied_by_postulate_CRB}
\mathsf{E}_{\mathbf{x}_{0}} \bigg\{ \frac{1}{f(\mathbf{y};\mathbf{x}_{0})} \frac{\partial^{\mathbf{p}} f(\mathbf{y}; \mathbf{x})}{\partial \mathbf{x}^{\mathbf{p}}} \bigg\} \leq 
\sqrt{\mathsf{E}_{\mathbf{x}_{0}} \bigg\{ \left( \frac{1}{f(\mathbf{y};\mathbf{x}_{0})} \frac{\partial^{\mathbf{p}} f(\mathbf{y}; \mathbf{x})}{\partial \mathbf{x}^{\mathbf{p}}}Ê\right)^{2} \bigg\}}< \infty
\end{equation}
due to the Cauchy-Schwarz inequality (cf.\ Theorem \ref{thm_cauchy_schwarz}) for the Hilbert space $\mathcal{P}(\mathcal{M})$ (see \eqref{equ_def_est_finite_power} and Theorem \ref{thm_hilbert_space_est_functions}) with inner product $\langle \, \cdot \, , \, \cdot \, \rangle_{\text{\tiny{RV}}}$ 
and the two specific estimators $\hat{g}_{1}(\mathbf{y}) = \frac{1}{f(\mathbf{y};\mathbf{x}_{0})} \frac{\partial^{\mathbf{p}} f(\mathbf{y}; \mathbf{x})}{\partial \mathbf{x}^{\mathbf{p}}}$ 
and $\hat{g}_{2}(\mathbf{y})=1$ where obviously $\langle \hat{g}_{2}(\cdot), \hat{g}_{2}(\cdot) \rangle_{\text{\tiny{RV}}} = 1$. 

Then, for any two multi-indices $\mathbf{p}_{1}, \mathbf{p}_{2} \in \mathbb{Z}_{+}^{N}$ with $\| \mathbf{p}_{1}\|_{\infty}, \| \mathbf{p}_{2}\|_{\infty} \leq m$ 
and parameter vectors $\mathbf{x}_{1}, \mathbf{x}_{2} \in \mathcal{X} \cap \mathcal{B}(\mathbf{x}_{0},r)$ we have  
\begin{align}
\label{equ_proof_crb_regul_impl_diff_1}
\infty & \stackrel{(a)}{>} \sqrt{\mathsf{E}_{\mathbf{x}_{0}} \bigg\{ \left( \frac{1}{f(\mathbf{y};\mathbf{x}_{0})} \frac{\partial^{\mathbf{p}_{1}} f(\mathbf{y}; \mathbf{x}_{1})}{\partial \mathbf{x}_{1}^{\mathbf{p}_{1}}}Ê\right)^{2} \bigg\} 
\mathsf{E}_{\mathbf{x}_{0}} \bigg\{ \left( \frac{1}{f(\mathbf{y};\mathbf{x}_{0})} \frac{\partial^{\mathbf{p}_{2}} f(\mathbf{y}; \mathbf{x}_{2})}{\partial \mathbf{x}_{2}^{\mathbf{p}_{2}}}Ê\right)^{2} \bigg\} }\nonumber \\[6mm] 
& \stackrel{(b)}{\geq} 
 \mathsf{E}_{\mathbf{x}_{0}} \left\{ \frac{1}{f(\mathbf{y}; \mathbf{x}_{0})} \frac{\partial^{\mathbf{p}_{1}} f(\mathbf{y};\mathbf{x}_{1})} {\partial \mathbf{x}_{1}^{\mathbf{p}_1} } 
 \frac{1}{f(\mathbf{y}; \mathbf{x}_{0})} \frac{\partial^{\mathbf{p}_2} f(\mathbf{y};\mathbf{x}_{2})} {\partial \mathbf{x}_{2}^{\mathbf{p}_2} }  \right\} \nonumber \\[6mm] 
& =
\int_{\mathbf{y}}  \frac{1}{f(\mathbf{y}; \mathbf{x}_{0})} \frac{\partial^{\mathbf{p}_{1}} f(\mathbf{y};\mathbf{x}_{1})} {\partial \mathbf{x}_{1}^{\mathbf{p}_1} } 
 \frac{1}{f(\mathbf{y}; \mathbf{x}_{0})} \frac{\partial^{\mathbf{p}_2} f(\mathbf{y};\mathbf{x}_{2})} {\partial \mathbf{x}_{2}^{\mathbf{p}_2} }  f(\mathbf{y}; \mathbf{x}_{0}) d \mathbf{y} \nonumber \\[6mm] 
& = 
\int_{\mathbf{y}}  \frac{1}{f(\mathbf{y}; \mathbf{x}_{0})} \frac{\partial^{\mathbf{p}_{1}} f(\mathbf{y};\mathbf{x}_{1})} {\partial \mathbf{x}_{1}^{\mathbf{p}_1} } 
 \frac{\partial^{\mathbf{p}_2} f(\mathbf{y};\mathbf{x}_{2})} {\partial \mathbf{x}_{2}^{\mathbf{p}_2} }  d \mathbf{y} \nonumber \\[6mm] 
& \stackrel{(c)}{=}
\frac{\partial^{\mathbf{p}_{2}} \int_{\mathbf{y}}  \frac{ f(\mathbf{y};\mathbf{x}_{2})}{f(\mathbf{y}; \mathbf{x}_{0})} 
\frac{\partial^{\mathbf{p}_{1}} f(\mathbf{y};\mathbf{x}_{1})} {\partial \mathbf{x}_{1}^{\mathbf{p}_1} } d \mathbf{y} } { \partial \mathbf{x}_{2}^{\mathbf{p}} }  \nonumber \\[3mm] 
& \stackrel{(d)}{=}
\frac{\partial^{\mathbf{p}_{2}} \partial^{\mathbf{p}_{1}} \int_{\mathbf{y}}  \frac{f(\mathbf{y};\mathbf{x}_{1})  f(\mathbf{y};\mathbf{x}_{2})}{f(\mathbf{y}; \mathbf{x}_{0})} 
d \mathbf{y} } { \partial \mathbf{x}_{1}^{\mathbf{p}_{1}} \partial \mathbf{x}_{2}^{\mathbf{p}_{2}} }  
  = \frac{\partial^{\mathbf{p}_{1}} \partial^{\mathbf{p}_{2}} R_{\mathcal{M}}(\mathbf{x}_1, \mathbf{x}_{2})}{ \partial \mathbf{x}_{1}^{\mathbf{p}_1} \partial \mathbf{x}_{2}^{\mathbf{p}_2} }, 
\end{align}
where $(a)$ is due to \eqref{equ_crb_regul_cond_finite_part_der}, 
step $(b)$ is due to the Cauchy-Schwarz inequality (cf.\ Theorem \ref{thm_cauchy_schwarz}) for the Hilbert space $\mathcal{P}(\mathcal{M})$ (see \eqref{equ_def_est_finite_power} and Theorem \ref{thm_hilbert_space_est_functions}), 
step $(c)$ follows from \eqref{equ_crb_regul_cond_intercah_diff_integr} with $h(\mathbf{y}) =  \frac{1}{f(\mathbf{y}; \mathbf{x}_{0})} \frac{\partial^{\mathbf{p}_{1}} f(\mathbf{y};\mathbf{x}_{1})} {\partial \mathbf{x}_{1}^{\mathbf{p}_1} }$ and $(d)$ follows by another application of \eqref{equ_crb_regul_cond_intercah_diff_integr} for the choice 
$h(\mathbf{y})= \frac{f(\mathbf{y}; \mathbf{x}_{2})}{f(\mathbf{y}; \mathbf{x}_{0})}$. The application of \eqref{equ_crb_regul_cond_intercah_diff_integr} is allowed 
in step $(c)$ and $(d)$ since a finite value of the right hand side of \eqref{equ_crb_regul_cond_intercah_diff_integr_intral_repr} and \eqref{equ_crb_regul_cond_intercah_diff_integr} is guaranteed by \eqref{equ_finite_value_expect_implied_by_postulate_CRB} and \eqref{equ_finite_integral_RKHS_cond}, respectively. 
Thus, we have shown that the partial derivative $\frac{\partial^{\mathbf{p}_{1}} \partial^{\mathbf{p}_{2}} R_{\mathcal{M}}(\mathbf{x}_1, \mathbf{x}_{2})}{ \partial \mathbf{x}_{1}^{\mathbf{p}_1} \partial \mathbf{x}_{2}^{\mathbf{p}_2} }$ exists for $\mathbf{p}_{1}, \mathbf{p}_{2} \in \mathbb{Z}_{+}^{N}$ with $\| \mathbf{p}_{1}\|_{\infty}, \| \mathbf{p}_{2}\|_{\infty} \leq m$ 
and parameter vectors $\mathbf{x}_{1}, \mathbf{x}_{2} \in \mathcal{X} \cap \mathcal{B}(\mathbf{x}_{0},r)$. Since moreover, these partial derivatives are 
continuous functions of $\mathbf{x}_{1}$ and $\mathbf{x}_{2}$, due the continuity of the function \eqref{equ_CRB_regul_cond_inner_prod_cont_x_1_x_2} in Postulate \ref{post_regular_cond_CRB}, we have that the kernel $R_{\mathcal{M}}(\cdot,\cdot):\mathcal{X} \times \mathcal{X} \rightarrow \mathbb{R}$ is differentiable up to order $m$.
\end{proof} 

\subsubsection{Unconstrained CRB}

The following form of the unconstrained CRB has been presented in \cite{HeroRecursiveCRB,HeroUniformCRB}: 
\begin{theorem}
\label{thm_unconstr_CR}
Consider an estimation problem $\mathcal{E}=\scalarestproblem$ and associated minimum variance problem $\mathcal{M}=\left(\mathcal{E},c(\cdot) \equiv 0, \mathbf{x}_{0} \right)$ that satisfy the conditions of Postulate \ref{assumption_RKHS_minvarproblem} and Postulate \ref{post_regular_cond_CRB} with $m=1$. Then the variance $v(\hat{g}(\cdot);\mathbf{x}_{0})$ at $\mathbf{x}_{0}$ 
of any unbiased estimator $\hat{g}(\cdot)$ with finite variance at $\mathbf{x}_{0}$ is 
lower bounded by 
\begin{equation}
\label{equ_unconstrained_CR}
v(\hat{g}(\cdot);\mathbf{x}_{0}) \geq  \left( \frac{\partial g(\mathbf{x})} {\partial \mathbf{x}}\bigg|_{\mathbf{x}_{0}} \right)^{T} \mathbf{J}_{\mathbf{x}_{0}}^{\dagger}  \,\, \frac{\partial g(\mathbf{x})} {\partial \mathbf{x}}\bigg|_{\mathbf{x}_{0}} 
\end{equation} 
where $\mathbf{J}_{\mathbf{x}_{0}} \in \mathbb{R}^{N \times N}$, known as the Fisher information matrix associated to the estimation problem $\mathcal{E}$, is given 
elementwise by 
\begin{equation} 
\label{equ_def_FIM_CRB}
\left( \mathbf{J}_{\mathbf{x}_{0}} \right)_{k,l} \triangleq \mathsf{E}_{\mathbf{x}_{0}} \bigg\{  \frac{\partial \log f(\mathbf{y}; \mathbf{x})}{\partial x_{k}} \bigg|_{\mathbf{x} = \mathbf{x}_{0}}  \frac{\partial \log f(\mathbf{y}; \mathbf{x})}{\partial x_{l}} \bigg|_{\mathbf{x} = \mathbf{x}_{0}} \bigg \}.
\end{equation} 
\end{theorem}
\begin{proof}
See \cite{HeroUniformCRB,kay}.  
\end{proof} 

We now provide an RKHS-based derivation of Theorem \ref{thm_unconstr_CR} by showing that under the validity of Postulate \ref{post_regular_cond_CRB}, the bound in \eqref{equ_unconstrained_CR} 
can be obtained via an orthogonal projection on a specific subspace of the RKHS $\mathcal{H}(\mathcal{M})$. 
To this end, we define the subspace $\mathcal{U}_{\text{\tiny{CR}}} \subseteq \mathcal{H}(\mathcal{M})$ as 
\begin{equation}
\label{equ_subspace_CR_RKHS}
\mathcal{U}_{\text{\tiny{CR}}} \triangleq  \linspan \big\{ \{ v_{0} (\cdot) \} \cup \{ v_{l}(\cdot) \}_{l \in [N]} \big\}, 
\end{equation} 
with the functions $v_{0} (\cdot ) \triangleq R_{\mathcal{M}}(\cdot, \mathbf{x}_{0}) \in \mathcal{H}(\mathcal{M})$ and 
\begin{equation}
\label{equ_subspace_CR_RKHS_vectors}
v_{l}(\cdot) \triangleq \frac{\partial R_{\mathcal{M}}(\cdot,\mathbf{x})}{\partial x_{l}} \bigg|_{\mathbf{x} = \mathbf{x}_{0}} \in \mathcal{H}(\mathcal{M}) \mbox{,} \quad l \in [N].  
\end{equation} 
By Theorem \ref{thm_regul_cond_CRB_imply_diff}, the kernel $R_{\mathcal{M}}(\cdot, \cdot)$ is differentiable up to order $m=1$. Therefore, 
the functions $\left\{ v_{l}(\cdot) \right\}_{l \in [N]}$ belong to $\mathcal{H}(R)$ due to Theorem \ref{thm_der_repr_prop}. In particular, 
the definition of the functions $v_{l}(\cdot)$ is of the same form as those in \eqref{equ_def_part_der_func}.
The inner product between any $v_{l}(\cdot)$ with $l \in [N]$ and $v_{0}(\cdot)$ is given by 
\begin{align} 
\label{equ_CR_orthog_vecs_1}
\big \langle v_{0}(\cdot) , v_{l}(\cdot) \big\rangle_{\mathcal{H}(\mathcal{M})} & \stackrel{(a)}{=}   \big\langle R_{\mathcal{M}}(\cdot, \mathbf{x}_{0}) , v_{l}(\cdot) \big\rangle_{\mathcal{H}(\mathcal{M})}
=  v_{l}(\mathbf{x}_{0}) \nonumber \\[6mm] 
&  \stackrel{(b)}{=}   \frac{\partial R_{\mathcal{M}}(\mathbf{x}_{0},\mathbf{x})}{\partial x_{l}} \bigg|_{\mathbf{x} = \mathbf{x}_{0}} 
 \stackrel{(c)}{=}   \frac{\partial}{\partial x_{l}}  \int_{\mathbf{y}} \frac{ f(\mathbf{y}; \mathbf{x}) f(\mathbf{y}; \mathbf{x}_{0})}{ f(\mathbf{y};\mathbf{x}_{0})} d\mathbf{y} \bigg|_{\mathbf{x} = \mathbf{x}_{0}}  \nonumber \\[6mm]
 & =   \frac{\partial}{\partial x_{l}}  \int_{\mathbf{y}} f(\mathbf{y}; \mathbf{x})d\mathbf{y} \bigg|_{\mathbf{x} = \mathbf{x}_{0}}  =   \frac{\partial}{\partial x_{l}} 1 \big|_{\mathbf{x} = \mathbf{x}_{0}}  = 0. 
\end{align} 
Here, we used for the step $(a)$ the reproducing property \eqref{equ_reproduction_property}, for $(b)$ the definition of the functions $v_{l}(\cdot)$ and for $(c)$ the definition of the kernel $R_{\mathcal{M}}(\cdot,\cdot)$ \eqref{equ_def_kernel_M}. 

Now we consider the projection $\mathbf{P}_{\mathcal{U}_{\text{\tiny{CR}}}} g(\cdot)$ of the parameter function $g(\cdot)$ onto the subspace $\mathcal{U}_{\text{\tiny{CR}}}$. 
Since Theorem \ref{thm_unconstr_CR} makes a statement about an unbiased estimator of $g(\cdot)$ with a finite variance at $\mathbf{x}_{0}$, we can assume that $g(\cdot)$ is estimable, i.e., we have $g(\cdot) \in \mathcal{H}(\mathcal{M})$ by Theorem \ref{thm_main_facts_RKHS_MVE}. (Note that $\gamma(\cdot)=g(\cdot)$ in Theorem \ref{thm_main_facts_RKHS_MVE}.)
Due to \eqref{equ_CR_orthog_vecs_1}, we have that subset $\mathcal{U}_{\text{\tiny{CR}}}$ (cf.\ \eqref{equ_subspace_CR_RKHS}) 
is spanned by the union of the two mutually orthogonal sets $\{ v_{0} (\cdot) \}$ and $\{ v_{l}(\cdot) \}_{l \in [N]}$. 
Therefore, we can apply Theorem \ref{thm_norm_projection_finite_dim_subspace_union_orthog_subspaces} to obtain the squared norm of the projection $\mathbf{P}_{\mathcal{U}_{\text{\tiny{CR}}}} g(\cdot)$ as 
\begin{align} 
\label{equ_projection_subspace_CR} 
\| \mathbf{P}_{\mathcal{U}_{\text{\tiny{CR}}}} g(\cdot) \|^{2}_{\mathcal{H}(\mathcal{M})} 
& =  \big\langle g(\cdot), v_{0}(\cdot) \big\rangle_{\mathcal{H}(\mathcal{M})}Ê\big( \big\langle v_{0}(\cdot), v_{0}(\cdot) \big\rangle_{\mathcal{H}(\mathcal{M})} \big)^{-1} \big\langle g(\cdot), v_{0}(\cdot) \big\rangle_{\mathcal{H}(\mathcal{M})}Ê + \mathbf{a}^{T} \mathbf{V}^{\dagger} \mathbf{a}  \nonumber \\[3mm]
& =  \left( g(\mathbf{x}_{0}) \right)^{2}+ \mathbf{a}^{T} \mathbf{V}^{\dagger} \mathbf{a},
\end{align} 
where, according to the reproducing property \eqref{equ_reproduction_property}, 
we used $\big\langle g(\cdot), v_{0}(\cdot) \big\rangle_{\mathcal{H}(\mathcal{M})}=g(\mathbf{x}_{0})$ and 
\begin{equation}
\big\langle v_{0}(\cdot), v_{0}(\cdot) \big\rangle_{\mathcal{H}(\mathcal{M})} = \big\langle R_{\mathcal{M}}(\cdot,\mathbf{x}_{0}), R_{\mathcal{M}}(\cdot,\mathbf{x}_{0}) \big\rangle_{\mathcal{H}(\mathcal{M})} = R(\mathbf{x}_{0}, \mathbf{x}_{0}) \stackrel{\eqref{equ_kernel_one_arg_x_0_equal_1}}{=}1.
\end{equation}
The vector $\mathbf{a} \in \mathbb{R}^{N}$ and the matrix $\mathbf{V} \in \mathbb{R}^{N \times N}$ 
in \eqref{equ_projection_subspace_CR} are defined elementwise as 
$a_{l} = \big\langle g(\cdot), v_{l}(\cdot)  \big\rangle_{\mathcal{H}(\mathcal{M})}$  
and $\left( \mathbf{V} \right)_{k,l}  = \big\langle v_{k}(\cdot), v_{l} (\cdot)  \big\rangle_{\mathcal{H}(\mathcal{M})}$, respectively. 
Note that due to Theorem \ref{thm_der_repr_prop} (using $g^{(\mathbf{p})}_{\mathbf{x}_{c}}(\cdot) = v_{l}(\cdot)$ for $\mathbf{p}=\mathbf{e}_{l}$, $\mathbf{x}_{c} = \mathbf{x}_{0}$ and $f(\cdot) = g(\cdot)$),
\begin{equation}
\label{equ_inner_prod_par_der_derivation_UCRB}
a_{l} = \big\langle  v_{l}(\cdot), g(\cdot)  \big\rangle_{\mathcal{H}(\mathcal{M})} = \frac{\partial g(\mathbf{x})}{\partial x_{l}} \bigg|_{\mathbf{x}=\mathbf{x}_{0}}, \quad l \in [N]. 
\end{equation} 
Furthermore,  
\begin{align} 
\left( \mathbf{V} \right)_{k,l}  & = \big\langle v_{k}(\cdot), v_{l}(\cdot)  \big\rangle_{\mathcal{H}(\mathcal{M})} \nonumber \\[4mm]
&  \stackrel{(a)}{=} \frac{\partial^{2} R_{\mathcal{M}}(\mathbf{x}_{1},\mathbf{x}_{2})}{\partial x_{1,k} x_{2,l}} \bigg|_{\mathbf{x}_{1} = \mathbf{x}_{2} = \mathbf{x}_{0}}   =   \frac{\partial^{2}}{\partial x_{1,k} x_{2,l}}  \int_{\mathbf{y}} \frac{ f(\mathbf{y}; \mathbf{x}_{1}) f(\mathbf{y}; \mathbf{x}_{2})}{ f(\mathbf{y};\mathbf{x}_{0})}  d\mathbf{y}  \bigg|_{\mathbf{x}_{1} = \mathbf{x}_{2} = \mathbf{x}_{0}} \nonumber \\[4mm]
   & \stackrel{(b)}{=}    \int_{\mathbf{y}} \frac{\partial^{2}}{\partial x_{1,k} x_{2,l}}  \frac{ f(\mathbf{y}; \mathbf{x}_{1}) f(\mathbf{y}; \mathbf{x}_{2})}{ f(\mathbf{y};\mathbf{x}_{0})} \bigg|_{\mathbf{x}_{1} = \mathbf{x}_{2} = \mathbf{x}_{0}} d\mathbf{y}  Ê\nonumber \\[4mm]
   & =  \int_{\mathbf{y}} \frac{\partial}{\partial x_{1,k}}  \frac{ f(\mathbf{y}; \mathbf{x}_{1})} { f(\mathbf{y};\mathbf{x}_{0})}  \bigg|_{\mathbf{x}_{1}  = \mathbf{x}_{0}}  \frac{\partial}{\partial x_{2,l}}  \frac{f(\mathbf{y}; \mathbf{x}_{2})}{ f(\mathbf{y};\mathbf{x}_{0})} \bigg|_{\mathbf{x}_{2} = \mathbf{x}_{0}}  f(\mathbf{y};\mathbf{x}_{0})d\mathbf{y} \nonumber   \\[4mm]Ê
   & = \int_{\mathbf{y}} \frac{\partial}{\partial x_{1,k}}  \log f(\mathbf{y}; \mathbf{x}_{1})\bigg|_{\mathbf{x}_{1}  = \mathbf{x}_{0}}  \frac{\partial}{\partial x_{2,l}}  \log f(\mathbf{y}; \mathbf{x}_{2}) \bigg|_{ \mathbf{x}_{2} = \mathbf{x}_{0}}  f(\mathbf{y};\mathbf{x}_{0})d\mathbf{y} \nonumber   \\[4mm]
   & = \mathsf{E}_{\mathbf{x}_{0}} \bigg\{ \frac{\partial}{\partial x_{1,k}}  \log f(\mathbf{y}; \mathbf{x}_{1})\bigg|_{\mathbf{x}_{1}  = \mathbf{x}_{0}}  \frac{\partial}{\partial x_{2,l}}  \log f(\mathbf{y}; \mathbf{x}_{2}) \bigg|_{ \mathbf{x}_{2} = \mathbf{x}_{0}}  \bigg\}, \label{equ_gramian_vectors_crb}
\end{align} 
where step $(a)$ follows from Theorem \ref{thm_der_repr_prop} (using $g^{(\mathbf{p})}_{\mathbf{x}_{c}}(\cdot)$ for the two choices $\mathbf{p}=\mathbf{e}_{k}$ and $\mathbf{p}=\mathbf{e}_{l}$, respectively, and $\mathbf{x}_{c} = \mathbf{x}_{0}$) and $(b)$ is obtained by two applications of \eqref{equ_crb_regul_cond_intercah_diff_integr_intral_repr} in Postulate \ref{post_regular_cond_CRB} (for the first application we use $h(\mathbf{y}) = \frac{f(\mathbf{y}; \mathbf{x}_{1})}{f(\mathbf{y}; \mathbf{x}_{0})}$ and differentiate w.r.t.\ 
$x_{2,l}$; and for the second application we set $h(\mathbf{y})=   \frac{\partial}{\partial x_{2,l}}  \frac{f(\mathbf{y}; \mathbf{x}_{2})}{ f(\mathbf{y};\mathbf{x}_{0})} \bigg|_{\mathbf{x}_{2} = \mathbf{x}_{0}}$ and differentiate w.r.t.\ $x_{1,k}$).
Thus we have $\mathbf{V} = \mathbf{J}_{\mathbf{x}_{0}}$, where $\mathbf{J}_{\mathbf{x}_{0}}$ is the Fisher information matrix as defined in \eqref{equ_def_FIM_CRB}.
As a consequence, by combining \eqref{equ_squared_norm_min_achiev_var} of Theorem \ref{thm_main_facts_RKHS_MVE}, \eqref{equ_lower_bound_RKHS_projection}, \eqref{equ_projection_subspace_CR} and \eqref{equ_inner_prod_par_der_derivation_UCRB} we obtain 
\begin{align}
\label{equ_proof_unconstrained_CRB_projection_subspace_main_facts}
L_{\mathcal{M}} & \stackrel{\eqref{equ_squared_norm_min_achiev_var}}{=} \| g(\cdot) \|_{\mathcal{H}(\mathcal{M})}^{2} - \big( g(\mathbf{x}_{0}) \big)^{2} 
\stackrel{\eqref{equ_lower_bound_RKHS_projection}}{\geq} \| \mathbf{P}_{\mathcal{U}_{\text{\tiny{CR}}}} g(\cdot) \|^{2}_{\mathcal{H}(\mathcal{M})} - \big( g(\mathbf{x}_{0}) \big)^{2}  \nonumber \\[3mm] 
& \stackrel{\eqref{equ_projection_subspace_CR}, \eqref{equ_inner_prod_par_der_derivation_UCRB}}{=} \big( g(\mathbf{x}_{0}) \big)^{2} +  \left( \frac{\partial g(\mathbf{x})} {\partial \mathbf{x}}\bigg|_{\mathbf{x}_{0}} \right)^{T} \mathbf{J}_{\mathbf{x}_{0}}^{\dagger}  \,\, \frac{\partial g(\mathbf{x})} {\partial \mathbf{x}}\bigg|_{\mathbf{x}_{0}} - \big( g(\mathbf{x}_{0}) \big)^{2} 
 =  \left( \frac{\partial g(\mathbf{x})} {\partial \mathbf{x}}\bigg|_{\mathbf{x}_{0}} \right)^{T} \mathbf{J}_{\mathbf{x}_{0}}^{\dagger}  \,\, \frac{\partial g(\mathbf{x})} {\partial \mathbf{x}}\bigg|_{\mathbf{x}_{0}}.
\end{align} 
This coincides with the unconstrained CRB in \eqref{equ_unconstrained_CR}, since we have $v(\hat{g}(\cdot); \mathbf{x}_{0}) \geq L_{\mathcal{M}}$ 
by the definition of $L_{\mathcal{M}}$ in \eqref{equ_def_min_ach_var}.

\subsubsection{Constrained CRB} 
Consider an estimation problem $\mathcal{E}=\scalarestproblem$ with associated minimum variance problem 
$\mathcal{M}=\left(\mathcal{E},c(\cdot) \equiv 0, \mathbf{x}_{0}\right)$. We asume that the parameter set $\mathcal{X}$ is defined via a set of equality constraints, i.e., 
\begin{equation} 
\label{equ_def_paramter_set_CCRB} 
\mathcal{X} = \{ÊÊ\mathbf{x} \in \mathbb{R}^{N} \big| \mathbf{f}(\mathbf{x})= \mathbf{0} \},
\end{equation} 
where $\mathbf{f}(\cdot): \mathbb{R}^{N} \rightarrow \mathbb{R}^{L}$ is a continuously differentiable vector-valued function. 
We also assume that $\mathbf{f}(\cdot)$ is such that the set $\mathcal{X}$ is nonempty and that its Jacobian $\mathbf{F}(\mathbf{x}) \triangleq \frac{\partial \mathbf{f}(\mathbf{x})}{ \partial \mathbf{x}} \in \mathbb{R}^{L \times N}$, 
given elementwise by $\left( \mathbf{F}(\mathbf{x}) \right)_{m,n} \triangleq \frac{\partial  f_{m}(\mathbf{x})}{\partial x_{n}}$, 
has rank $L$ whenever $\mathbf{f}(\mathbf{x})=\mathbf{0}$, i.e., for every $\mathbf{x} \in \mathcal{X}$.
Such parameter sets are considered e.g.\ in \cite{StoicaNgCCRB,ZvikaCCRB,MooreCCRB}. 

As discussed in Section \ref{sec_effect_red_par_set}, there are a two fundamentally different ways in how one can exploit the prior knowledge that is given by 
a specific form of the parameter set (e.g., such as \eqref{equ_def_paramter_set_CCRB}). 
For the problem of estimating the parameter vector itself, i.e., $\mathbf{g}(\mathbf{x}) = \mathbf{x}$, the authors of \cite{StoicaNgCCRB} require an estimator $\hat{\mathbf{g}}(\cdot): \mathbb{R}^{M} \rightarrow \mathbb{R}^{N}$ 
to take on values only in $\mathcal{X}$, i.e., $\hat{g}(\mathbf{y}) \in \mathcal{X}$ for all $\mathbf{y} \in \mathbb{R}^{M}$.This is reasonable since if it is known that the unknown 
parameter vector $\mathbf{x}$ is an element $\mathcal{X}$, then the value of an estimator of $\mathbf{x}$ should also be an element of $\mathcal{X}$. 

However, in this work, we will focus on another approach that is also considered by the authors of \cite{ZvikaCCRB}. Within this second approach, 
we do not constrain the values of the estimator itself to be an element of $\mathcal{X}$ but instead we only require the bias of the estimator to equal the 
prescribed bias for parameter vectors $\mathbf{x} \in \mathcal{X}$. 
In this approach, the a priori knowledge of a certain type of parameter set $\mathcal{X}$ helps us in that we must place fewer bias constraints on an estimator. 
Thus, the set of allowed estimators  $\mathcal{F}(\mathcal{M})$ in \eqref{equ_est_finite_var_prescr_bias} 
tends to become larger as the parameter set $\mathcal{X}$ is getting smaller, i.e., more ``informative''. 

The constrained CRB is a generalization of the unconstrained CRB to estimation problems with a parameter set $\mathcal{X}$ of the form \eqref{equ_def_paramter_set_CCRB}.
A specific formulation of the constrained CRB is given in \cite{ZvikaCCRB}:
\begin{theorem}
\label{thm_ccrb}
Consider an estimation problem $\mathcal{E}=\scalarestproblem$ and associated minimum variance problem $\mathcal{M}=\left(\mathcal{E},c(\cdot) \equiv 0, \mathbf{x}_{0} \right)$ that satisfy Postulate \ref{assumption_RKHS_minvarproblem} and Postulate \ref{post_regular_cond_CRB} with $m=1$ and whose parameter set $\mathcal{X}$
satisfies \eqref{equ_def_paramter_set_CCRB}. Then the variance $v(\hat{g}(\cdot); \mathbf{x}_{0})$ at $\mathbf{x}_{0} \in \mathcal{X}$ of any unbiased estimator $\hat{g}(\cdot)$ with a finite variance 
at $\mathbf{x}_{0}$ is lower bounded by 
\begin{equation} 
\label{equ_CCRB} 
v(\hat{g}(\cdot); \mathbf{x}_{0}) \geq  \left( \frac{\partial g(\mathbf{x})} {\partial \mathbf{x}}\bigg|_{\mathbf{x}_{0}} \right)^{T} \mathbf{U}(\mathbf{x}_{0}) \left( \mathbf{U}^{T}(\mathbf{x}_{0}) \mathbf{J}_{\mathbf{x}_{0}} \mathbf{U}(\mathbf{x}_{0}) \right)^{\dagger} \mathbf{U}^{T}(\mathbf{x}_{0})  \,\, \frac{\partial g(\mathbf{x})} {\partial \mathbf{x}}\bigg|_{\mathbf{x}_{0}}  
\end{equation} 
where $\mathbf{J}_{\mathbf{x}_{0}}$ is the FIM as defined in \eqref{equ_def_FIM_CRB} and $\mathbf{U}(\mathbf{x}_{0}) \in \mathbb{R}^{N \times (N-L)}$ is any matrix whose column vectors form an ONB for the null space $\mathcal{N}(\mathbf{F}(\mathbf{x}_{0}))$ 
of the Jacobian $\mathbf{F}(\mathbf{x}_{0})$, i.e., 
\begin{equation}
\mathbf{F}(\mathbf{x}_{0}) \mathbf{U}(\mathbf{x}_{0}) = \mathbf{0}, \quad \mathbf{U}^{T}(\mathbf{x}_{0}) \mathbf{U}(\mathbf{x}_{0}) = \mathbf{I}_{(N-L)}.  
\end{equation} 
\end{theorem} 
\begin{proof}
\cite{ZvikaCCRB} 
\end{proof}
Note that if the function $\mathbf{f}(\mathbf{x})$ in \eqref{equ_def_paramter_set_CCRB} is given as $\mathbf{f}(\mathbf{x}) \equiv \mathbf{0}$, 
one can verify that the matrix $\mathbf{U}(\mathbf{x}_{0})$ used in Theorem \ref{thm_ccrb} can be chosen as $\mathbf{U}(\mathbf{x}_{0}) = \mathbf{I}$ in which case 
the constrained CRB \eqref{equ_CCRB} reduces to the unconstrained CRB \eqref{equ_unconstrained_CR}. 
We now re-derive the bound \eqref{equ_CCRB} using the RKHS $\mathcal{H}(\mathcal{M})$ associated to the minimum variance problem $\mathcal{M}$ 
of Theorem \ref{thm_ccrb}. This derivation will be very similar to the RKHS-based derivation of the unconstrained CRB as 
performed above. However, there are some modifications necessary due to the additional assumption on the parameter set $\mathcal{X}$ given by \eqref{equ_def_paramter_set_CCRB}. These modifications are based mainly on the methods used in \cite{MooreCCRB}. 

In particular, will show that the bound \eqref{equ_CCRB} can be obtained via an orthogonal projection on a 
finite-dimensional subspace $\mathcal{U}_{\text{\tiny{CCRB}}} \subseteq \mathcal{H}(\mathcal{M})$ which is defined by 
\begin{equation}
\label{equ_def_U_CCRB}
\mathcal{U}_{\text{\tiny{CCRB}}} \triangleq  \linspan \big\{ \{ v_{0} (\cdot) \} \cup \{ v_{l}(\cdot) \}_{l \in [N-L]} \big\},
\end{equation} 
where $v_{0}(\cdot) \triangleq R_{\mathcal{M}}(\cdot, \mathbf{x}_{0}) \in \mathcal{H}(\mathcal{M})$ and the functions $\{ v_{l}(\cdot) \}_{l \in [N-L]} \in \mathcal{H}(\mathcal{M})$ will be specified in the following. 
We first note that given the fixed vector $\mathbf{x}_{0}$ associated to the minimum variance problem $\mathcal{M}$, we have according
to the implicit function theorem  (see \cite[Theorem 3.3]{MooreCCRB} or \cite{RudinBookPrinciplesMatheAnalysis})
that there exists a continuously differentiable and bijective map $\mathbf{r}(\cdot)$, which moreover has a continuously differentiable inverse\footnote{Such a function is called a diffeomorphism \cite{MooreCCRB}.} $\mathbf{r}^{-1}(\cdot)$, 
from an open set
$\mathbb{O} \subseteq \mathbb{R}^{N-L}$ into an open set $\mathbb{P} \subseteq \mathcal{X}$ which also contains $\mathbf{x}_{0}$, i.e., 
\begin{equation} 
\label{equ_ccrb_reparam} 
\mathbf{r}(\cdot): \mathbb{O} \subseteq \mathbb{R}^{N-L} \rightarrow \mathbb{P} \subseteq \mathcal{X}
\end{equation} 
where $\mathbf{x}_{0} \in \mathbb{P}$. 
Note that any function value $\mathbf{r}({\bm \theta})$ is an element of the parameter set, i.e., $\mathbf{r}({\bm \theta}) \in \mathcal{X}$. 
The Jacobian of $\mathbf{r}(\cdot)$ at an arbitrary ${\bm \theta} \in \mathbb{O}$ will be denoted by $\mathbf{G}({\bm \theta}) \in \mathbb{R}^{N \times (N-L)}$, i.e.,
\begin{equation}
\label{equ_def_jacobian_ccrb}
\left( \mathbf{G}({\bm \theta} ) \right)_{m,n} = \frac{\partial r_{m}({\bm \theta})}{ \partial \theta_{n}}. 
\end{equation}  
We now define the functions $v_{l}(\cdot) \in \mathcal{H}(\mathcal{M})$ in \eqref{equ_def_U_CCRB} as 
\begin{equation} 
\label{equ_def_vectors_ccrb}
v_{l}(\cdot) \triangleq \frac{\partial R_{\mathcal{M}}(\cdot,\mathbf{r}({\bm \theta}))}{\partial \theta_{l}} \big|_{{\bm \theta} = \mathbf{r}^{-1}(\mathbf{x}_{0})} \mbox{,} \quad l \in [N-L].  
\end{equation}
The fact that the so-defined functions  $v_{l}(\cdot)$ belong to $\mathcal{H}(\mathcal{M})$ can again be verified straightforwardly 
by Theorem \ref{thm_regul_cond_CRB_imply_diff} together with Theorem \ref{thm_der_repr_prop}. Indeed, using the chain rule for differentiation \cite[Theorem 9.15]{RudinBookPrinciplesMatheAnalysis}, we have 
\begin{align}
\label{equ_part_der_chain_rule_ccrb_v_l}
v_{l}(\cdot) & = \frac{\partial R_{\mathcal{M}}(\cdot,\mathbf{r}({\bm \theta}))}{\partial \theta_{l}} \bigg|_{{\bm \theta} = \mathbf{r}^{-1}(\mathbf{x}_{0})}  \nonumber \\[4mm] 
& = \sum_{l' \in [N]} \frac{\partial R_{\mathcal{M}}(\cdot,\mathbf{x})}{\partial x_{l'}} \bigg|_{\mathbf{x}=\mathbf{x}_{0}}   \frac{\partial r_{l'}({\bm \theta} )}{\partial \theta_{l}} \bigg|_{{\bm \theta} = \mathbf{r}^{-1}(\mathbf{x}_{0})} \mbox{,} \quad l \in [N-L],    
\end{align} 
which is a linear combination of the functions $\frac{\partial R_{\mathcal{M}}(\cdot,\mathbf{x})}{\partial x_{l'}} \bigg|_{\mathbf{x}=\mathbf{x}_{0}}$ 
which are of the same form (for the choice $\mathbf{p} = \mathbf{e}_{l'}$ and $\mathbf{x}_{c} = \mathbf{x}_{0}$) as the functions $g^{(\mathbf{p})}_{\mathbf{x}_{c}}(\cdot)$ defined in Theorem \ref{thm_der_repr_prop}, 

The inner product between any $v_{l}(\cdot)$ with $l \in [N-L]$ and $v_{0}(\cdot)$ is given by (cf.\ our derivation in \eqref{equ_CR_orthog_vecs_1}) 
\begin{align} 
\label{equ_CCR_orthog_vecs_1}
\langle v_{0}(\cdot) , v_{l}(\cdot) \rangle_{\mathcal{H}(\mathcal{M})} & = \frac{\partial R_{\mathcal{M}}(\mathbf{x}_{0},\mathbf{r} ( {\bm \theta}))}{\partial \theta_{l}} \bigg|_{{\bm \theta} = {\bm \theta}_{0}}
 =  \frac{\partial}{\partial \theta_{l}}  \int_{\mathbf{y}} \frac{ f(\mathbf{y};\mathbf{r} ( {\bm \theta}) ) f(\mathbf{y}; \mathbf{x}_{0})}{ f(\mathbf{y};\mathbf{x}_{0})} d\mathbf{y} \bigg|_{{\bm \theta} = {\bm \theta}_{0}}  \nonumber \\[4mm]Ê
& =   \frac{\partial}{\partial \theta_{l}}  \int_{\mathbf{y}} f(\mathbf{y}; \mathbf{r} ( {\bm \theta}))d\mathbf{y} \bigg|_{{\bm \theta} = {\bm \theta}_{0}}  =   \frac{\partial}{\partial \theta_{l}} 1 \big|_{{\bm \theta} = {\bm \theta}_{0}}   = 0,
\end{align} 
where we introduced the shorthand ${\bm \theta_{0}} \triangleq \mathbf{r}^{-1}(\mathbf{x}_{0})$. 

Since Theorem \ref{thm_ccrb} makes a statement about an unbiased estimator of the parameter function $g(\cdot)$ that has a finite variance at $\mathbf{x}_{0}$, we can assume 
that the parameter function $g(\cdot)$ is estimable, i.e., we have $g(\cdot) \in \mathcal{H}(\mathcal{M})$ by Theorem \ref{thm_main_facts_RKHS_MVE} (where again $\gamma(\cdot) = g(\cdot)$). 

According to Theorem \ref{thm_norm_projection_finite_dim_subspace} and \eqref{equ_CCR_orthog_vecs_1}, we can express the squared norm of the orthogonal 
projection $\mathbf{P}_{\mathcal{U}_{\text{\tiny{CCRB}}}} g(\cdot)$ as 
\begin{align} 
\label{equ_projection_subspace_CCR} 
\| \mathbf{P}_{\mathcal{U}_{\text{\tiny{CCRB}}}} g(\cdot) \|^{2}_{\mathcal{H}(\mathcal{M})}  = \left( \langle g(\cdot), v_{0}(\cdot) \rangle_{\mathcal{H}(\mathcal{M})}Ê\right)^{2} + \mathbf{a}^{T} \mathbf{V}^{\dagger} \mathbf{a}  = \left( g(\mathbf{x}_{0}) \right)^{2}+ \mathbf{a}^{T} \mathbf{V}^{\dagger} \mathbf{a}, 
\end{align} 
where the vector $\mathbf{a} \in \mathbb{R}^{N-L}$ and the matrix $\mathbf{V} \in \mathbb{R}^{(N-L) \times (N-L)}$ are defined elementwise as 
$a_{l}  \triangleq \langle g(\cdot), v_{l} (\cdot)  \rangle_{\mathcal{H}(\mathcal{M})}$ and $\left( \mathbf{V} \right)_{k,l}   \triangleq \langle v_{k}(\cdot), v_{l}(\cdot)  \rangle_{\mathcal{H}(\mathcal{M})}$, respectively. 
Due to the derivative-reproducing property \eqref{equ_der_reproduction_prop} and the chain rule of differentiation \cite[Theorem 9.15]{RudinBookPrinciplesMatheAnalysis}, we obtain
\begin{align} 
a_{l}  & = \langle g(\cdot), v_{l} (\cdot)  \rangle_{\mathcal{H}(\mathcal{M})} = 
\Bigg\langle g(\cdot), \frac{\partial R_{\mathcal{M}}(\cdot,\mathbf{r}({\bm \theta}))}{\partial \theta_{l}} \bigg|_{{\bm \theta} =  {\bm \theta_{0}}} \Bigg\rangle_{\mathcal{H}(\mathcal{M})} 
\nonumber \\[4mm] 
& \stackrel{\eqref{equ_part_der_chain_rule_ccrb_v_l}}{=}
 \Bigg\langle g(\cdot), \sum_{l' \in [N]} \frac{\partial R_{\mathcal{M}}(\cdot,\mathbf{x})}{\partial x_{l'}} \bigg|_{\mathbf{x}=\mathbf{x}_{0}}  
  \frac{\partial r_{l'}({\bm \theta} )}{\partial \theta_{l}} \bigg|_{{\bm \theta} = {\bm \theta_{0}} }  \Bigg\rangle_{\mathcal{H}(\mathcal{M})} \nonumber \\[4mm] 
  & =
\sum_{l' \in [N]} \frac{\partial r_{l'}({\bm \theta} )}{\partial \theta_{l}} \bigg|_{{\bm \theta} = {\bm \theta_{0}} }  \Bigg\langle g(\cdot),  \frac{\partial R_{\mathcal{M}}(\cdot,\mathbf{x})}{\partial x_{l'}} \bigg|_{\mathbf{x}=\mathbf{x}_{0}}  
  \Bigg\rangle_{\mathcal{H}(\mathcal{M})} \nonumber \\[4mm] 
& \stackrel{\eqref{equ_der_reproduction_prop}}{=}
 \sum_{l' \in [N]}  \frac{\partial r_{l'}({\bm \theta} )}{\partial \theta_{l}} \bigg|_{{\bm \theta} = {\bm \theta_{0}} } \frac{\partial g(\mathbf{x})}{\partial x_{l'}} \bigg|_{\mathbf{x}_{0}}
 \stackrel{\eqref{equ_def_jacobian_ccrb}}{=} \sum_{l' \in [N]} \left( \mathbf{G}({\bm \theta_{0}}) \right)_{l',l}  \frac{\partial g(\mathbf{x})}{\partial x_{l'}} \bigg|_{\mathbf{x}_{0}},
\end{align}
, i.e., in vector notation:
\begin{equation} 
\label{equ_inner_prod_ccrb}
\mathbf{a} =   \mathbf{G}^{T}({\bm \theta_{0}}) \frac{\partial g(\mathbf{x})}{\partial \mathbf{x}} \bigg|_{\mathbf{x}_{0}}.
\end{equation}

Furthermore, by using the chain rule for derivation once again, we obtain
\begin{align}
\left(\mathbf{V} \right)_{k,l} &  =   \langle v_{k}(\cdot), v_{l}(\cdot)  \rangle_{\mathcal{H}(\mathcal{M})} 
  = \Bigg\langle\frac{\partial R_{\mathcal{M}}(\cdot,\mathbf{r}({\bm \theta}))}{\partial \theta_{k}} \bigg|_{{\bm \theta} =  {\bm \theta_{0}}}, \frac{\partial R_{\mathcal{M}}(\cdot,\mathbf{r}({\bm \theta}))}{\partial \theta_{l}} \bigg|_{{\bm \theta} =  {\bm \theta_{0}}} \Bigg\rangle_{\mathcal{H}(\mathcal{M})} \nonumber \\[4mm]
 & = \Bigg\langle\sum_{l' \in [N]} \frac{\partial R_{\mathcal{M}}(\cdot,\mathbf{x})}{\partial x_{l'}} \bigg|_{\mathbf{x}=\mathbf{x}_{0}}  
  \frac{\partial r_{l'}({\bm \theta} )}{\partial \theta_{k}} \bigg|_{{\bm \theta} = {\bm \theta_{0}} },\sum_{l' \in [N]} \frac{\partial R_{\mathcal{M}}(\cdot,\mathbf{x})}{\partial x_{l'}} \bigg|_{\mathbf{x}=\mathbf{x}_{0}}  
  \frac{\partial r_{l'}({\bm \theta} )}{\partial \theta_{l}} \bigg|_{{\bm \theta} = {\bm \theta_{0}} } \Bigg\rangle_{\mathcal{H}(\mathcal{M})} \nonumber \\[4mm]
 & = \sum_{l',l'' \in [N]} \frac{\partial r_{l'}({\bm \theta} )}{\partial \theta_{k}} \bigg|_{{\bm \theta} = {\bm \theta_{0}} }\frac{\partial r_{l''}({\bm \theta} )}{\partial \theta_{l}} \bigg|_{{\bm \theta} = {\bm \theta_{0}} } \Bigg \langle\frac{\partial R_{\mathcal{M}}(\cdot,\mathbf{x})}{\partial x_{l'}} \bigg|_{\mathbf{x}=\mathbf{x}_{0}}  
 , \frac{\partial R_{\mathcal{M}}(\cdot,\mathbf{x})}{\partial x_{l''}} \bigg|_{\mathbf{x}=\mathbf{x}_{0}}  
   \Bigg\rangle_{\mathcal{H}(\mathcal{M})} \nonumber \\[4mm]
 & \stackrel{(a)}{=} \sum_{l',l'' \in [N]} \left( \mathbf{G} ({\bm \theta}_{0})\right)_{l',k} \left( \mathbf{G} ({\bm \theta}_{0})\right)_{l'',l}   \frac{\partial^{2} R_{\mathcal{M}}(\mathbf{x}_{1},\mathbf{x}_{2})}{ \partial x_{1,l'} \partial x_{2,l''}} \bigg|_{ \mathbf{x}_{1} = \mathbf{x}_{2} = \mathbf{x}_{0}} 
 \nonumber \\[4mm]
& = \left( \mathbf{G}({\bm \theta}_{0})^{T}  \mathbf{W} \mathbf{G}({\bm \theta}_{0}) \right)_{k,l},  \label{equ_gramian_vectors_ccrb}
\end{align}
where step $(a)$ is due to Theorem \ref{thm_der_repr_prop} (using $g^{(\mathbf{p})}_{\mathbf{x}_{c}}(\cdot)$ for the choices $\mathbf{p}=\mathbf{e}_{l'}$, $\mathbf{p} = \mathbf{e}_{l''}$ and $\mathbf{x}_{c} = \mathbf{x}_{0}$), and the matrix $\mathbf{W} \in \mathbb{R}^{N \times N}$ is given elementwise via 
\begin{equation} 
\label{equ_intermed_matrix_ccrb}
\left( \mathbf{W} \right)_{k,l} \triangleq \frac{\partial^{2} R_{\mathcal{M}}(\mathbf{x}_{1},\mathbf{x}_{2})}{\partial x_{1,k} \partial x_{2,l}} \bigg|_{\mathbf{x}_{1} = \mathbf{x}_{2} = \mathbf{x}_{0}}  = \mathsf{E}_{\mathbf{x}_{0}} \bigg\{ \frac{\partial}{\partial x_{1,k}}  \log f(\mathbf{y}; \mathbf{x}_{1})\bigg|_{\mathbf{x}_{1}  = \mathbf{x}_{0}}  \frac{\partial}{\partial x_{2,l}}  \log f(\mathbf{y}; \mathbf{x}_{2}) \bigg|_{ \mathbf{x}_{2} = \mathbf{x}_{0}}  \bigg\}.
\end{equation} 
Here, we reused for \eqref{equ_intermed_matrix_ccrb} the derivation of \eqref{equ_gramian_vectors_crb}. According to \eqref{equ_intermed_matrix_ccrb}, the matrix $\mathbf{W}$ 
coincides with the FIM $\mathbf{J}_{\mathbf{x}_{0}}$ as defined in \eqref{equ_def_FIM_CRB}, i.e., $\mathbf{W} = \mathbf{J}_{\mathbf{x}_{0}}$. 
Putting together the pieces, we obtain 
\begin{align}
\label{equ_lower_bound_L_ccrb}
L_{\mathcal{M}} & \stackrel{\eqref{equ_squared_norm_min_achiev_var}}{=} \| g(\cdot) \|^{2}_{2} - \big( g(\mathbf{x}_{0}) \big)^{2} \stackrel{\eqref{equ_lower_bound_RKHS_projection}}{\geq}  \| \mathbf{P}_{\mathcal{U}_{\text{\tiny{CCRB}}}} g(\cdot) \|^{2}_{\mathcal{H}(\mathcal{M})}- \big( g(\mathbf{x}_{0}) \big)^{2}  \nonumber \\[4mm] 
&  \stackrel{\eqref{equ_projection_subspace_CCR}}{=}\mathbf{a}^{T} \mathbf{V}^{\dagger} \mathbf{a} \stackrel{\eqref{equ_gramian_vectors_ccrb},\eqref{equ_inner_prod_ccrb}}{=} \left( \frac{\partial g(\mathbf{x})} {\partial \mathbf{x}}\bigg|_{\mathbf{x}_{0}} \right)^{T} \mathbf{G}({\bm \theta}_{0}) \left( \mathbf{G}^{T}({\bm \theta}_{0}) \mathbf{J}_{\mathbf{x}_{0}} \mathbf{G}({\bm \theta}_{0}) \right)^{\dagger} \mathbf{G}^{T}({\bm \theta}_{0})  \,\, \frac{\partial g(\mathbf{x})} {\partial \mathbf{x}}\bigg|_{\mathbf{x}_{0}}. 
\end{align} 

As shown in \cite[p.\ 29]{MooreCCRB}, the Jacobian $\mathbf{G}({\bm \theta}_{0})$ in \eqref{equ_def_jacobian_ccrb} and the matrix $\mathbf{U}(\mathbf{x}_{0})$ appearing in Theorem \ref{thm_ccrb} are related by an invertible 
linear transformation, i.e., there 
exists a nonsingular matrix $\mathbf{S} \in \mathbb{R}^{(N-L) \times (N-L)}$ (wich may depend on $\mathbf{x}_{0}$) such that 
\begin{equation}
\label{equ_ccrb_relation_onb_null_space_G_matrix}
\mathbf{G}({\bm \theta}_{0}) \mathbf{S} = \mathbf{U}(\mathbf{x}_{0}). 
\end{equation}  
However, one can always find a map $\mathbf{r}(\cdot)$ such that $\mathbf{S} = \mathbf{I}$, in which case $\mathbf{G}({\bm \theta}_{0}) = \mathbf{U}(\mathbf{x}_{0})$ 
\cite[p.\ 50] {MooreCCRB}. Indeed, suppose that the map $\mathbf{r}(\cdot)$ is such that 
$\mathbf{S} \neq \mathbf{I}$ and therefore $\mathbf{G}({\bm \theta}_{0}) \neq \mathbf{U}(\mathbf{x}_{0})$. We can then 
use for the definition of the functions $v_{l}(\cdot)$ in \eqref{equ_def_vectors_ccrb} instead of $\mathbf{r}(\cdot)$ a new map $\mathbf{r}'(\cdot)$ 
defined as $\mathbf{r}'({\bm \theta}) = \mathbf{r} ( \mathbf{S} {\bm \theta})$, whose 
Jacobian $\mathbf{G}'({\bm \theta}_{0})$ satisfies  
\begin{equation}
\mathbf{G}' ({\bm \theta}_{0}) \stackrel{(a)}{=} \mathbf{G}({\bm \theta}_{0}) \mathbf{S} \stackrel{\eqref{equ_ccrb_relation_onb_null_space_G_matrix}}{=} \mathbf{U}(\mathbf{x}_{0}),
\end{equation} 
where step $(a)$ is due to the chain rule for differentiation \cite[Theorem 9.15]{RudinBookPrinciplesMatheAnalysis}. 

Thus, since we can assume that $\mathbf{G}({\bm \theta}_{0}) = \mathbf{U}(\mathbf{x}_{0})$, we have by the definition of $L_{\mathcal{M}}$ in \eqref{equ_def_min_ach_var} that the bound \eqref{equ_lower_bound_L_ccrb} coincides with the bound \eqref{equ_CCRB}.

\subsection{Bhattacharyya Bound} 
The CRB of a given minimum variance problem depends only on the first-order partial derivatives of the statistical model $f(\mathbf{y}, \mathbf{x})$ w.r.t.\ to the parameter $\mathbf{x}$. 
By incorporating also higher-order partial derivatives, one obtains the Bhattacharyya bound \cite{LowerBoundAbel,bhattacharyya}. Like the CRB the Bhattacharyya bound, is a lower bound on the variance 
of any unbiased estimator for a given estimation problem (which is required to satisfy some regularity conditions). The Bhattacharyya bound may be formulated in our setting as follows.
\begin{theorem} 
\label{thm_bhattacharyya}
Consider an estimation problem $\mathcal{E}=\scalarestproblem$ with associated minimum variance problem $\mathcal{M}=\left(\mathcal{E},c(\cdot) \equiv 0, \mathbf{x}_{0} \right)$ 
that satisfy the conditions of Postulate \ref{assumption_RKHS_minvarproblem} and Postulate \ref{post_regular_cond_CRB} with maximum order of differentiation equal to $m$. 
For any set of $L \in \mathbb{N}$ multi-indices $\{ \mathbf{p}_{l} \in \mathbb{Z}_{+}^{N} \}_{l \in [L]}$ with $\| \mathbf{p}_{l} \|_{\infty} \leq m$, we have that 
the variance $v(\hat{g}(\cdot);\mathbf{x}_{0})$ at $\mathbf{x}_{0}$ 
of any unbiased estimator $\hat{g}(\cdot)$ with finite variance at $\mathbf{x}_{0}$ is lower bounded by 
\begin{equation}
\label{equ_bhattacharyya}
v(\hat{g}(\cdot);\mathbf{x}_{0}) \geq  \mathbf{b}^{T} \mathbf{B}^{\dagger} \mathbf{b},
\end{equation} 
where the vector $\mathbf{b} \in \mathbb{R}^{L}$ and the matrix $\mathbf{B} \in \mathbb{R}^{L \times L}$ are given elementwise by 
\begin{equation} 
b_{l} \triangleq \frac{\partial^{\mathbf{p}_{l}} g(\mathbf{x})} {\partial \mathbf{x}^{\mathbf{p}_{l}} }\bigg|_{\mathbf{x}_{0}}
\end{equation}
and 
\begin{equation} 
\label{equ_def_bhattacharyya_matrix}
\left( \mathbf{B} \right)_{k,l} \triangleq \mathsf{E}_{\mathbf{x}_{0}}
 \bigg\{  \frac{1}{f^{2}(\mathbf{y}; \mathbf{x}_{0})} \frac{\partial^{\mathbf{p}_{k}}  f(\mathbf{y}; \mathbf{x})}{\partial \mathbf{x}^{\mathbf{p}_{k}} } \bigg|_{\mathbf{x} = \mathbf{x}_{0}}  \frac{\partial^{\mathbf{p}_{l}}  f(\mathbf{y}; \mathbf{x})}{\partial \mathbf{x}^{\mathbf{p}_{l}} }  \bigg|_{\mathbf{x} = \mathbf{x}_{0}} \bigg \},
\end{equation} 
respectively. 
\end{theorem}
A derivation of Theorem \ref{thm_bhattacharyya} can be based on the RKHS $\mathcal{H}(\mathcal{M})$. 
This derivation is completely analogous to that of Theorem \ref{thm_unconstr_CR}. The only 
difference is that instead of the subspace $\mathcal{U}_{\text{\tiny{CR}}}$ defined via \eqref{equ_subspace_CR_RKHS} and \eqref{equ_subspace_CR_RKHS_vectors}, 
which is used for the derivation of Theorem \ref{thm_unconstr_CR}, one has to use 
the subspace $\mathcal{U}_{\text{\tiny{BHATT}}} \triangleq  \linspan \{ \{ v_{0} (\cdot) \} \cup \{ v_{l}(\cdot) \}_{l \in [L]} \}$. 
Here, we used the functions $v_{0} (\cdot ) \triangleq R_{\mathcal{M}}(\cdot, \mathbf{x}_{0}) \in \mathcal{H}(\mathcal{M})$ and $v_{l}(\cdot) \triangleq \frac{\partial^{\mathbf{p}_{l}} R_{\mathcal{M}}(\cdot,\mathbf{x})}{\partial \mathbf{x}^{\mathbf{p}_{l}}} \big|_{\mathbf{x} = \mathbf{x}_{0}} \in \mathcal{H}(\mathcal{M})$ for $l \in [L]$. The fact that the functions $ \{ v_{l}(\cdot) \}_{l \in [L]}$ belong to the RKHS $\mathcal{H}(\mathcal{M})$ may be verified by Theorem \ref{thm_regul_cond_CRB_imply_diff} together with Theorem \ref{thm_der_repr_prop}.

While the RKHS interpretation of the Bhattacharyya bound has already been presented in \cite{Duttweiler73b} for a specific estimation problem, the above derivation holds for general estimation problems.
Finally, note that the unconstrained CRB in Theorem \ref{thm_unconstr_CR} is obtained as a special case of the Bhattacharyya bound for the choices $L=N$, $m=1$ and multi-indices $\{ \mathbf{p}_{l} = \mathbf{e}_{l} \}_{lÊ\in [L]}$ in Theorem \ref{thm_bhattacharyya}. 

\subsection{Hammersley-Chapman-Robbins Bound} 
\label{sec_HCRB}
Despite its popularity, a drawback of the CRB (and also of the Bhattacharyya bound) is that it only exploits the local structure of an estimation problem $\mathcal{E}$ around a specific point $\mathbf{x}_{0} \in \mathcal{X}$ \cite{LowerBoundAbel}. 
As an illustrative example, consider two different estimation problems $\mathcal{E}_{1} = \left(\mathcal{X}_{1},f(\mathbf{y}; \mathbf{x}),g(\cdot)\right)$ and $\mathcal{E}_{2} = \left(\mathcal{X}_{2},f(\mathbf{y}; \mathbf{x}),g(\cdot)\right)$ 
that share the same statistical model $f(\mathbf{y}; \mathbf{x})$ and parameter function $g(\cdot)$ but are defined for two different parameter sets $\mathcal{X}_{1}$ and $\mathcal{X}_{2}$. Let us assume that both parameter sets are open balls 
centered at $\mathbf{x}_{0}$ but with different radii $r_{1}$ and $r_{2}$, i.e., $\mathcal{X}_{1} = \mathcal{B}(\mathbf{x}_{0},r_{1})$ and $\mathcal{X}_{2} = \mathcal{B}(\mathbf{x}_{0},r_{2})$ with $r_{1} \neq r_{2}$. Then the CRB at $\mathbf{x}_{0}$ will be the same for both estimation problems irrespective of the precise values of $r_{1}$, $r_{2}$. 
Thus, the CRB ignores in some sense a part of the information that is contained in the parameter set $\mathcal{X}$. 

On the other hand, the Barankin bound exploits the full information carried by the parameter set $\mathcal{X}$ since it is the tightest possible lower bound on the estimator variance. 
However, in the general case, it is difficult to evaluate the exact Barankin bound in \eqref{equ_min_var_barankin_bound}. 

The Hammersley-Chapman-Robbins bound (HCRB) \cite{GormanHero,ChapmanRobbins51,Hammersley50} is a lower bound on the estimator variance which also incorporates the 
``global'' structure of the estimation problem, i.e., it depends on the parameter set $\mathcal{X}$ not only through its local structure as does the CRB. 
However, it can be evaluated much more easily than the Barankin bound.
In our context, the HCRB can be stated as follows \cite{GormanHero}: 
\begin{theorem} 
\label{thm_HCRB} 
Consider an estimation problem $\mathcal{E}=\scalarestproblem$ and associated minimum variance problem  $\mathcal{M}=\left(\mathcal{E},c(\cdot) \equiv 0, \mathbf{x}_{0} \right)$ 
that satisfy Postulate \ref{assumption_RKHS_minvarproblem}. Furthermore, consider a finite set of ``test-points'' $\{ \mathbf{x}_{l} \in \mathcal{X} \}_{l \in [L]}$. 
Then the variance $v(\hat{g}(\cdot);\mathbf{x}_{0})$ at $\mathbf{x}_{0} \in \mathcal{X}$ of any unbiased estimator $\hat{g}(\cdot)$ with a finite variance 
at $\mathbf{x}_{0}$ is lower bounded by 
\begin{equation}
\label{equ_HCRB} 
v(\hat{g}(\cdot);\mathbf{x}_{0}) \geq \mathbf{m}^{T} \mathbf{V}^{\dagger} \mathbf{m},
\end{equation}  
where the vector $\mathbf{m} \in \mathbb{R}^{L}$ and the matrix $\mathbf{V} \in \mathbb{R}^{L \times L}$ are defined elementwise by 
\begin{equation}
\label{equ_def_vector_HCRB_thm}
m_{l} \triangleq g(\mathbf{x}_{l}) - g(\mathbf{x}_{0})
\end{equation} 
and 
\begin{equation} 
\left( \mathbf{V} \right)_{l,l'} \triangleq \mathsf{E}_{\mathbf{x}_{0}} \bigg\{  \frac{ \big[f(\mathbf{y};\mathbf{x}_{l}) - f(\mathbf{y};\mathbf{x}_{0}) \big]\big[f(\mathbf{y};\mathbf{x}_{l'}) - f(\mathbf{y};\mathbf{x}_{0}) \big]}{ f^{2}(\mathbf{y}; \mathbf{x}_{0})} \bigg\},
\end{equation} 
respectively.
\end{theorem} 

\begin{proof}
\cite{GormanHero}
\end{proof}
We now show how the bound in \eqref{equ_HCRB} can be obtained via the RKHS $\mathcal{H}(\mathcal{M})$. 
First, note that since Theorem \ref{thm_HCRB} makes a statement about an existing unbiased estimator of $g(\mathbf{x})$ with finite variance at $\mathbf{x}_{0}$, we can assume that $g(\cdot)$ is estimable, i.e., 
$g(\cdot) \in \mathcal{H}(\mathcal{M})$ according to Theorem \ref{thm_main_facts_RKHS_MVE}. 
Then, let us define the subspace 
\begin{equation} 
\label{equ_def_subspace_HCRB}
\mathcal{U}_{\text{\tiny{HCRB}}} \triangleq \linspan \big\{ \{ v_{0} (\cdot) \} \cup \{ v_{l} (\cdot) \}_{l \in [L]} \big\},
\end{equation} 
where $L$ is the number of test-points $\big\{ \mathbf{x}_{l} \big\}_{l \in [L]}$ used in Theorem \ref{thm_HCRB}. The subspace $\mathcal{U}_{\text{\tiny{HCRB}}}$ is spanned by the functions  
$v_{0}(\cdot) \triangleq R_{\mathcal{M}}(\cdot, \mathbf{x}_{0})$ and $\big\{ v_{l} (\cdot) \triangleq R_{\mathcal{M}}(\cdot, \mathbf{x}_{l}) - R_{\mathcal{M}}(\cdot, \mathbf{x}_{0})\big\}_{l \in [L]}$. 
Note that we obviously have by Theorem \ref{thm_constr_RKHS_closure_linear_span} that these functions belong to $\mathcal{H}(\mathcal{M})$. 

The function $v_{0}(\cdot) \in \mathcal{H}(\mathcal{M})$ is orthogonal to any $v_{l} (\cdot) \in \mathcal{H}(\mathcal{M})$ with $l \in [L]$, i.e., 
the two sets $\{ v_{0} (\cdot) \}$ and $\{ v_{l} (\cdot) \}_{l \in [L]}$ are mutually 
orthogonal, since by the reproducing property \eqref{equ_reproduction_property} we have
\begin{align}
\label{equ_orthog_proof_HCRB_1} 
\big\langle v_{0}(\cdot), v_{l}(\cdot) \big\rangle_{\mathcal{H}(\mathcal{M})} 
& =  \big\langle R_{\mathcal{M}}(\cdot, \mathbf{x}_{0}),  R_{\mathcal{M}}(\cdot, \mathbf{x}_{l}) - R_{\mathcal{M}}(\cdot, \mathbf{x}_{0})\big\rangle_{\mathcal{H}(\mathcal{M})} \nonumber \\[4mm]
& =  \big\langle R_{\mathcal{M}}(\cdot, \mathbf{x}_{0}),  R_{\mathcal{M}}(\cdot, \mathbf{x}_{l}) \big\rangle_{\mathcal{H}(\mathcal{M})} -  
\big\langle R_{\mathcal{M}}(\cdot, \mathbf{x}_{0}),  R_{\mathcal{M}}(\cdot, \mathbf{x}_{0})\big\rangle_{\mathcal{H}(\mathcal{M})} \nonumber \\[4mm]
& = R_{\mathcal{M}}(\mathbf{x}_{l}, \mathbf{x}_{0}) - R_{\mathcal{M}}(\mathbf{x}_{0}, \mathbf{x}_{0}) \stackrel{\eqref{equ_kernel_one_arg_x_0_equal_1}}{=} 1 -1 = 0.
\end{align}
Now we define the matrix $\mathbf{V} \in \mathbb{R}^{L \times L}$ elementwise by
\begin{align} 
\label{equ_def_Grammian_proof_RKHS_HCRB}
\left( \mathbf{V} \right)_{l,l'} \triangleq \big\langle v_{l}(\cdot), v_{l'}(\cdot) \big\rangle_{\mathcal{H}(\mathcal{M}}  
& = \bigg \langle R_{\mathcal{M}}(\cdot, \mathbf{x}_{l}) - R_{\mathcal{M}}(\cdot, \mathbf{x}_{0}),  R_{\mathcal{M}}(\cdot, \mathbf{x}_{l'}) - R_{\mathcal{M}}(\cdot, \mathbf{x}_{0}) \bigg\rangle_{\mathcal{H}(\mathcal{M})} \nonumber \\[4mm]
& \stackrel{\eqref{equ_reproduction_property}}{=}R_{\mathcal{M}}(\mathbf{x}_{l'}, \mathbf{x}_{l}) - R_{\mathcal{M}}(\mathbf{x}_{l'}, \mathbf{x}_{0}) - R_{\mathcal{M}}(\mathbf{x}_{l}, \mathbf{x}_{0})+R_{\mathcal{M}}(\mathbf{x}_{0}, \mathbf{x}_{0}).
\end{align} 
The inner products of $g(\cdot)$ with the vectors $v_{l}(\cdot)$ are calculated (also using the reproducing property \eqref{equ_reproduction_property}) as 
\begin{equation} 
\label{equ_inner_prod_HCRB}
\big\langle g(\cdot), v_{0}(\cdot) \big\rangle_{\mathcal{H}(\mathcal{M})} = g(\mathbf{x}_{0})\mbox{, } \,\,\, \mbox{and } \,\,\, \langle g(\cdot), v_{l}(\cdot) \rangle_{\mathcal{H}(\mathcal{M})} = g(\mathbf{x}_{l}) - g(\mathbf{x}_{0}) = m_{l} \,\, \mbox{ for } l \in [L]. 
\end{equation}
By projecting the prescribed mean function $\gamma(\cdot) = g(\cdot) + c(\cdot) = g(\cdot)$ of the minimum variance problem $\mathcal{M}$ on the subspace $\mathcal{U}_{\text{\tiny{HCRB}}}$, we obtain 
\begin{align}
\label{equ_min_var_bound_HCRB}
L_{\mathcal{M}} & \stackrel{\eqref{equ_squared_norm_min_achiev_var}}{=} \| g(\cdot) \|^{2}_{2} - \big( g(\mathbf{x}_{0}) \big)^{2} \stackrel{\eqref{equ_lower_bound_RKHS_projection}}{\geq}  \| \mathbf{P}_{\mathcal{U}_{\text{\tiny{HCRB}}}} g(\cdot) \|^{2}_{\mathcal{H}(\mathcal{M})}- \big( g(\mathbf{x}_{0}) \big)^{2}  \nonumber \\[4mm] 
&  \stackrel{(a)}{=} \big\langle g(\cdot), v_{0}(\cdot) \big\rangle_{\mathcal{H}(\mathcal{M})} \big( \big\langle v_{0}(\cdot), v_{0}(\cdot) \big\rangle_{\mathcal{H}(\mathcal{M})}\big)^{-1} \big\langle g(\cdot), v_{0}(\cdot) \big\rangle_{\mathcal{H}(\mathcal{M})} + \mathbf{m}^{T} \mathbf{V}^{\dagger} \mathbf{m} - \big( g(\mathbf{x}_{0}) \big)^{2} \nonumber \\[4mm]
&  \stackrel{(b)}{=} \big( g(\mathbf{x}_{0}) \big)^{2} + \mathbf{m}^{T} \mathbf{V}^{\dagger} \mathbf{m} - \big( g(\mathbf{x}_{0}) \big)^{2} = \mathbf{m}^{T} \mathbf{V}^{\dagger} \mathbf{m} ,
\end{align} 
where for step $(a)$ we used \eqref{equ_inner_prod_HCRB} and Theorem \ref{thm_norm_projection_finite_dim_subspace_union_orthog_subspaces}, which applies since $\mathcal{U}_{\text{\tiny{HCRB}}}$ is defined as the linear span of the union of the two 
mutually orthogonal function sets $\{ v_{0} (\cdot) \}$ and $\{ v_{l} (\cdot) \}_{l \in [L]}$. The step $(b)$ follows from 
\begin{equation}
\big\langle v_{0}(\cdot), v_{0}(\cdot) \big\rangle_{\mathcal{H}(\mathcal{M})} = \big\langle R_{\mathcal{M}}(\cdot,\mathbf{x}_{0}), R_{\mathcal{M}}(\cdot,\mathbf{x}_{0}) \big\rangle_{\mathcal{H}(\mathcal{M})} \stackrel{\eqref{equ_reproduction_property}}{=} R_{\mathcal{M}}(\mathbf{x}_{0}, \mathbf{x}_{0}) \stackrel{\eqref{equ_kernel_one_arg_x_0_equal_1}}{=}1.
\end{equation}
and \eqref{equ_inner_prod_HCRB}. 

The bound in \eqref{equ_HCRB} follows from \eqref{equ_min_var_bound_HCRB}, the definition of $L_{\mathcal{M}}$ in \eqref{equ_def_min_ach_var} and the fact that 
\begin{align} 
 \left( \mathbf{V} \right)_{l,l'} & \stackrel{\eqref{equ_def_Grammian_proof_RKHS_HCRB}}{=} 
 R_{\mathcal{M}}(\mathbf{x}_{l'}, \mathbf{x}_{l}) - R_{\mathcal{M}}(\mathbf{x}_{l'}, \mathbf{x}_{0}) - R_{\mathcal{M}}(\mathbf{x}_{l}, \mathbf{x}_{0})+R_{\mathcal{M}}(\mathbf{x}_{0}, \mathbf{x}_{0})  \nonumber \\[4mm]
 & = \mathsf{E}_{\mathbf{x}_{0}}\bigg\{  \frac{ f(\mathbf{y};\mathbf{x}_{l'})} { f(\mathbf{y}; \mathbf{x}_{0})} \frac{ f(\mathbf{y};\mathbf{x}_{l})}{f(\mathbf{y}; \mathbf{x}_{0})} \bigg\} 
- \mathsf{E}_{\mathbf{x}_{0}}\bigg\{   \frac{ f(\mathbf{y};\mathbf{x}_{l'})}{f(\mathbf{y}; \mathbf{x}_{0})}  \frac{ f(\mathbf{y};\mathbf{x}_{0})} { f(\mathbf{y}; \mathbf{x}_{0})}\bigg\} 
-\mathsf{E}_{\mathbf{x}_{0}}\bigg\{  \frac{ f(\mathbf{y};\mathbf{x}_{l})} { f(\mathbf{y}; \mathbf{x}_{0})} \frac{ f(\mathbf{y};\mathbf{x}_{0})}{f(\mathbf{y}; \mathbf{x}_{0})} \bigg\} \nonumber \\
& +\mathsf{E}_{\mathbf{x}_{0}}\bigg\{  \frac{ f(\mathbf{y};\mathbf{x}_{0})} { f(\mathbf{y}; \mathbf{x}_{0})} \frac{ f(\mathbf{y};\mathbf{x}_{0})}{f(\mathbf{y}; \mathbf{x}_{0})} \bigg\} \nonumber \\[4mm]
& = \mathsf{E}_{\mathbf{x}_{0}} \bigg\{  \frac{ f(\mathbf{y};\mathbf{x}_{l'}) - f(\mathbf{y};\mathbf{x}_{0}) }{ f(\mathbf{y}; \mathbf{x}_{0})} \frac{ f(\mathbf{y};\mathbf{x}_{l}) - f(\mathbf{y};\mathbf{x}_{0}) }{ f(\mathbf{y}; \mathbf{x}_{0})} \bigg\} \nonumber \\[4mm]
& =\mathsf{E}_{\mathbf{x}_{0}} \bigg\{  \frac{ \big[f(\mathbf{y};\mathbf{x}_{l}) - f(\mathbf{y};\mathbf{x}_{0}) \big]\big[f(\mathbf{y};\mathbf{x}_{l'}) - f(\mathbf{y};\mathbf{x}_{0}) \big]}{ f^{2}(\mathbf{y}; \mathbf{x}_{0})} \bigg\}.
\end{align}

\chapter{The Sparse Linear Model}
\label{chap_SLM}
\section{Introduction}
In this chapter, we study the problem of estimating a nonrandom parameter vector $\mathbf{x} \!\in\! \mathbb{R}^{N}\!$
which is known to be strictly
$S$-sparse, i.e., at most $S$ of its entries are nonzero, where $S \in [N]$ (typically $S \!\ll\! N$). We thus 
\vspace{.5mm}
have
\be
\label{equ_SLM_parameter}
\mathbf{x} \!\in\!  \mathcal{X}_{S} \,, \quad\; \text{with} \;\; \mathcal{X}_{S} = \big\{\mathbf{x}' \rmv\!\in\! \mathbb{R}^{N} \big| \ist {\|\mathbf{x}' \|}_{0} \rmv\leq\rmv S \big\} \,.
\vspace{.5mm}
\ee 
While the sparsity degree $S \in \mathbb{N}$ is assumed to be known, the set of positions of the nonzero entries of $\mathbf{x}$, i.e., the support $\supp(\mathbf{x})$, is unknown. 

The estimation of $\mathbf{x}$ is based on the observed vector $\mathbf{y} \!\in\! \mathbb{R}^{M}\!$ 
given 
\vspace{-.5mm}
by
\be
\label{equ_linear_observation_model}
\mathbf{y} \ist=\ist \mathbf{H} \mathbf{x} + \mathbf{n} \,,
\vspace{.2mm}
\ee
with a known 
system matrix $\mathbf{H} \!\in\! \mathbb{R}^{M \times N}\!$ and additive white Gaussian noise (AWGN) $\mathbf{n} \sim \mathcal{N}(\mathbf{0},  \sigma^{2} \mathbf{I})$ with 
a known positive (i.e., nonzero) noise variance $\sigma^{2} >0$. 
The matrix $\mathbf{H}$ is arbitrary, typically, it will be assumed to satisfy the 
\vspace{-.5mm}
requirement 
\begin{equation} 
\spark(\mathbf{H}) \rmv>\rmv   \, S \,, 
\label{equ_spark_cond}
\end{equation} 
where $\rm{spark}(\mathbf{H})$ denotes the minimum number of linearly dependent columns of $\mathbf{H}$ \cite{GreedisGood}.\footnote{For a matrix $\mathbf{H} \in \mathbb{R}^{M \times N}$ with full column rank, i.e., $\rank(\mathbf{H}) = N$, we define $\rm{spark}(\mathbf{H}) = N$.}
Note that we also allow $M \!<\! N$ 
(this case is relevant to compressed sensing methods as discussed in Section \ref{sec_slm_viewpoint_CS});
however, condition \eqref{equ_spark_cond} implies that $M \!\ge\! S$. 
Furthermore, note that the requirement \eqref{equ_spark_cond} is weaker than the standard requirement \cite{ZvikaCRB} 
\begin{equation} 
\spark(\mathbf{H}) \rmv>\rmv   \, 2S \,, 
\end{equation} 
which is reasonable since otherwise we can find two different parameter vectors $\mathbf{x}_{1}, \mathbf{x}_{2} \in \mathcal{X}_{S}$ for 
which the pdf of the observation is identical, i.e., $f(\mathbf{y}; \mathbf{x}_{1}) = f(\mathbf{y}; \mathbf{x}_{2})$ for every $\mathbf{y}$, 
which means that it is completely impossible to distinguish 
between $\mathbf{x}_{1}$ and $\mathbf{x}_{2}$ only by looking at the observation $\mathbf{y}$.

More formally, by combining the observation model \eqref{equ_linear_observation_model} with the requirement \eqref{equ_SLM_parameter}, 
we define the \emph{sparse linear model} (SLM) as the specific estimation problem denoted $\mathcal{E}_{\text{SLM}}$ and given as
\begin{equation} 
\label{equ_def_SLM_est_problem}
\mathcal{E}_{\text{SLM}} \triangleq \big( \mathcal{X}_{S}, f_{\mathbf{H}}(\mathbf{y}; \mathbf{x}), g(\mathbf{x}) = x_{k} \big), 
\end{equation} 
where $k \in [N]$ is an arbitrary but fixed index (the precise choice of $k$ will be made explicit whenever necessary) and with the statistical model 
\begin{equation} 
\label{equ_def_LGM_stat_model}
 f_{\mathbf{H}}(\mathbf{y}; \mathbf{x}) = \frac{1}{(2 \pi \sigma^{2})^{M/2} } \exp \left(- \frac{1}{2 \sigma^{2}} \| \mathbf{y} - \mathbf{H} \mathbf{x} \|^{2}_{2} \right).
\end{equation}
One may argue that considering the parameter function $g(\mathbf{x}) = x_{k}$ is to restrictive, since in practice one may be intersted in the whose parameter 
vector $\mathbf{x}$ rather than only a single entry $x_{k}$. However, as discussed in Section \ref{sec_sep_vector_scalar_param_function}, 
in the context of minimum variance estimation, 
we have that estimation of the parameter vector 
$\mathbf{x}$ is completely equivalent to the $N$ separate estimation problems 
\begin{equation}
\bigg\{\big( \mathcal{X}_{S}, f_{\mathbf{H}}(\mathbf{y}; \mathbf{x}), g(\mathbf{x}) = x_{k} \big) \bigg\}_{k \in [N]}.
\end{equation}
By comparing \eqref{equ_def_LGM_stat_model} with \eqref{equ_def_pdf_exp_fam}, we see that the statistical model of the SLM is an instance of the 
standard exponential family, with sufficient statistic ${\bm \Phi}(\mathbf{y}) = \frac{1}{\sigma^{2}} \mathbf{H}^{T} \mathbf{y}$, parameter function $\mathbf{u}(\mathbf{x}) = \mathbf{x}$, cumulant function $ A^{( {\bm \Phi})}(\mathbf{x}) =   \frac{1}{2 \sigma^{2}} \mathbf{x}^{T} \mathbf{H}^{T} \mathbf{H} \mathbf{x}$, 
and weight $h(\mathbf{y}) = \exp \left( - \frac{1}{2Ê\sigma^{2}} \mathbf{y}^{T} \mathbf{y} \right)$. 

The SLM is very similar to a well-known estimation problem, i.e., the \emph{linear Gaussian model} (LGM) \cite{poorsspbook,kay,scharf91,HeroUniformCRB}, which is given by 
\begin{equation}
\label{equ_def_LGM_est_problem}
\mathcal{E}_{\text{LGM}} \triangleq \big( \mathbb{R}^{N}, f_{\mathbf{H}}(\mathbf{y}; \mathbf{x}), g(\mathbf{x}) = x_{k} \big).
\end{equation} 
Indeed, the only difference between the LGM $\mathcal{E}_{\text{LGM}}$ and the SLM $\mathcal{E}_{\text{SLM}}$ is the parameter set. 
The SLM is obtained from the LGM by reducing the parameter set from $\mathcal{X} = \mathbb{R}^{N}$ to 
the set of $S$-sparse vectors, i.e., $\mathcal{X} = \mathcal{X}_{S}$. If the sparsity degree $S$ is equal to the dimension of the observation, i.e., $S=N$, then the SLM coincides 
with the LGM. 

The SLM is relevant, e.g., for sparse channel estimation \cite{Carbonelli06} where the sparse parameter vector $\mathbf{x}$ represents the taps of a linear time-invariant 
channel and the system matrix $\mathbf{H}$ represents the training signal. Generally, the SLM can be used for any type of sparse deconvolution \cite{MallatBook} like in optical coherence tomography (OCT) \cite{KailOCTIcassp09}, where the sparse vector $\mathbf{x}$ models the layer structure of the human eye and the matrix 
$\mathbf{H}$ represents the properties of the measurement device. 

An important special case of the SLM is given by $\mathbf{H} \!=\! \mathbf{I}$ (so that $M \!=\! N$), 
\vspace{-1mm}
i.e., 
\be
\label{equ_observation_model_SSNM}
\mathbf{y} \ist=\ist \mathbf{x} + \mathbf{n} \,,
\ee
where again $\mathbf{x} \!\in\!  \mathcal{X}_{S}$ and $\mathbf{n} \in \mathbb{R}^{N}$ represents AWGN, i.e., $\mathbf{n} \sim \mathcal{N}(\mathbf{0}, \sigma^{2}  \mathbf{I})$ with known variance $\sigma^{2} > 0$. This will be referred to as the \emph{sparse signal in noise model} (SSNM). 
Formally, the SSNM is defined as the estimation problem
\begin{equation} 
\label{equ_def_SSNM_est_problem}
\mathcal{E}_{\text{SSNM}} \triangleq \big( \mathcal{X}_{S}, f_{\mathbf{H}=\mathbf{I}}(\mathbf{y}; \mathbf{x}), g(\mathbf{x}) = x_{k} \big), 
\end{equation} 
i.e., with the statistical model \eqref{equ_def_LGM_stat_model} for the specific choice $\mathbf{H} = \mathbf{I}$ for the system matrix. 
The SSNM can be used, e.g., for channel estimation \cite{Carbonelli06} when 
the 
channel consists only of 
few significant taps and an orthogonal training signal is used \cite{OptimalTrainingDong}. 
Another potential application for the SSNM is image denoising using an orthonormal wavelet basis \cite{Donoho94idealspatial}. 
Note that due to Theorem \ref{thm_data_proc_inequ_classical_est}, any lower bound on the minimum achievable variance of 
a minimum variance problem $\mathcal{M}_{\text{SSNM}}=\left( \mathcal{E}_{\text{SSNM}}, c(\cdot), \mathbf{x}_{0} \right)$ entails a 
lower bound on the minimum achievable variance for a minimum variance problem for any estimation problem with parameter set $\mathcal{X}_{S}$ and 
observation model 
\begin{equation}
\mathbf{y} = \mathbf{H}(\mathbf{x} + \mathbf{n}), 
\end{equation}
where $\mathbf{n}$ denotes AWGN, i.e., $\mathbf{n} \sim \mathcal{N}(\mathbf{0}, \sigma^{2}\mathbf{I})$ and 
$\mathbf{H} \in \mathbb{R}^{M \times N}$ is an arbitrary (possibly random) matrix. 
However, if the system matrix $\mathbf{H} \in \mathbb{R}^{M \times N}$ in \eqref{equ_linear_observation_model} is a deterministic 
orthonormal matrix, i.e., $\mathbf{H}^{T} \mathbf{H} = \mathbf{I}$ (which implies 
that $M \geq N$), then we can transform the observation $\mathbf{y}$ in an invertible manner to obtain the modified observation $\mathbf{y}' = \mathbf{H}^{T} \mathbf{y} = \mathbf{x} + \mathbf{n}'$, where 
$\mathbf{n}' \sim \mathcal{N}(\mathbf{0}, \sigma^{2} \mathbf{I})$. The modified observation corresponds then to an SSNM and moreover, by Theorem \ref{thm_invariance_classical_est}, is equivalent to the original SLM since the transformation of the observation is invertible. Thus, any SLM with an orthonormal system matrix is, from the viewpoint of minimum variance estimation, completely equivalent to the SSNM 
(cf. \cite{AlexZvikaICASSP,AlexZvikaJournal}). 

The SLM has been considered in \cite{JustRelax}, where the authors analyze the performance of specific estimation schemes that are based on convex optimization problems. 
In a minimax estimation context, the SLM has been studied in \cite{RaskuttiMinmaxSLM,VerzelenMinmaxSLM}, where the authors derive bounds on the minimax risk (see Section \ref{sec_minimax_est}) and also 
discuss estimators that come close to these bounds. An asymptotic analysis of minimax estimation for the SSNM is given in the seminal work \cite{DonohoJohnstone94,Donoho94idealspatial}.
In our context, i.e., minimum variance estimation, lower bounds on the minimum achievable estimation variance for the SLM have been studied recently. 
In particular, the CRB for the SLM was derived and analyzed in \cite{ZvikaCRB,ZvikaSSP}.
In \cite{AlexZvikaICASSP}, lower and upper bounds on the minimum achievable variance of unbiased estimators were derived for the SSNM by using techniques that differ from the RKHS approach that will 
be used in this thesis. A remarkable property of the lower bounds presented in \cite{ZvikaCRB} and \cite{AlexZvikaICASSP} is the fact that they 
exhibit a discontinuity when passing from the case ${\| \mathbf{x} \|}_{0} \!=\! S$ to the case ${\| \mathbf{x} \|}_{0} \!<\! S$. 

In this chapter, we use the mathematical framework of RKHS 
to derive novel lower bounds on the variance of estimators for the SLM. 
These bounds hold for estimators with a given differentiable bias function.
For the special case of the SSNM, we obtain a lower bound for unbiased estimators 
which is tighter
than the bounds in \cite{ZvikaCRB,ZvikaSSP,AlexZvikaICASSP} and, moreover, is a continuous function of $\mathbf{x}$. As we will show, RKHS theory relates the 
lower bound for the SLM to that obtained for a 
specific LGM. The RKHS framework has been previously applied to estimation \cite{Parzen59,Duttweiler73b} but, 
to the best of our knowledge, not to the SLM.

This chapter is organized as follows:
In Section \ref{sec_RKHS_SLM}, we introduce and discuss the RKHS associated to minimum variance problems that arise from the SLM. 
Based on this RKHS, we discuss in Section \ref{sec_RKHS_basic_facts_SLM} some fundamental facts about minimum variance estimation for the SLM. 
In Section \ref{sec_RKHS_Lower_Bounds_SLM}, we use RKHS theory first to reinterpret existing lower bounds from the geometric RKHS perspective and then to 
derive novel lower bounds on the variance of estimators for the SLM which have a prescribed bias, or, equivalently, a prescribed mean function. 
The important special case of the SLM given by the SSNM is discussed in Section \ref{sec_SSNM}, where we derive closed-form expressions 
for the minimum achievable variance, i.e., the Barankin bound and for the corresponding LMV estimator. 
A discussion on the necessity of strict sparsity requirements is presented in Section \ref{sec_strict_sparstiy_SLM}. The SLM viewpoint on 
the CS recovery problem is presented in Section \ref{sec_slm_viewpoint_CS}. 
Finally, we compare the theoretical bounds with the actual variance behavior of some popular estimation schemes in Section \ref{sec_numerical_SLM}.

The key results of this chapter have been presented in part in \cite{AlexZvikaICASSP,RKHSAsilomar2010,AlexZvikaJournal,RKHSISIT2011}. 


\section{RKHS Associated with the SLM}
\label{sec_RKHS_SLM}

We will denote any minimum variance problem that 
is obtained from $\mathcal{E}_{\text{SLM}}$ \eqref{equ_def_SLM_est_problem} by additionally fixing a prescribed bias function 
$c(\cdot): \mathcal{X}_{S} \rightarrow \mathbb{R}$ and parameter vector $\mathbf{x}_{0} \in \mathcal{X}_{S}$, by 
\begin{equation}
\label{equ_def_min_var_problem_SLM}
\mathcal{M}_{\text{SLM}} \triangleq \left( \mathcal{E}_{\text{SLM}}, c(\cdot), \mathbf{x}_{0} \right).
\end{equation} 
It should be clear from the context which bias function $c(\cdot)$ and fixed parameter vector $\mathbf{x}_{0}$ is used within the minimum variance problem $\mathcal{M}_{\text{SLM}}$. 
Note that the restriction to the specific parameter function $g(\mathbf{x}) = x_{k}$ in the SLM \eqref{equ_def_SLM_est_problem} is no real restriction, 
since we do not constrain the prescribed bias $c(\cdot)$ used for the minimum variance problem $\mathcal{M}_{\text{SLM}}$ 
and, according to Theorem \ref{thm_equ_bias_param_function}, it is no loss of generality if one fixes either the prescribed bias or the parameter function of a minimum variance problem, but not both.

The estimation problem $\mathcal{E}_{\text{SLM}}$ satisfies the requirements of Postulate \ref{assumption_RKHS_minvarproblem}, 
since it can be easily verified that $\mathsf{E}_{\mathbf{x}_{0}} \bigg \{ \bigg[ \frac{ f_{\mathbf{H}}(\mathbf{y}; \mathbf{x}) }{  f_{\mathbf{H}}(\mathbf{y}; \mathbf{x}_{0})} \bigg]^{2} \bigg\} < \infty$ for any $\mathbf{x} \inÊ\mathcal{X}_{S}$.
Therefore, according to Section \ref{sec_def_RKHS_MVP}, we can associate with the minimum variance problem $\mathcal{M}_{\text{SLM}}$ 
an RKHS $\mathcal{H}(\mathcal{M}_{\text{SLM}})$ whose kernel $R_{\mathcal{M}_{\text{SLM}}}(\cdot,\cdot): \mathcal{X}_{S} \times \mathcal{X}_{S} \rightarrow \mathbb{R}$ is given 
as (cf. \eqref{equ_def_kernel_M}): 
\begin{align}
\label{equ_kernel_SLM}
R_{\mathcal{M}_{\text{SLM}}}(\mathbf{x}_{1}, \mathbf{x}_{2})  & =  \mathsf{E}_{\mathbf{x}_{0}} \{ \rho_{\mathcal{M}_{\text{SLM}}}(\mathbf{y},\mathbf{x}_{1}) \rho_{\mathcal{M}_{\text{SLM}}}(\mathbf{y},\mathbf{x}_{2}) \}  
=  \mathsf{E}_{\mathbf{x}_{0}} \left \{ \frac{ f_{\mathbf{H}}(\mathbf{y}; \mathbf{x}_{1}) }{  f_{\mathbf{H}}(\mathbf{y}; \mathbf{x}_{0})}   \frac{ f_{\mathbf{H}}(\mathbf{y}; \mathbf{x}_{2}) }{  f_{\mathbf{H}}(\mathbf{y}; \mathbf{x}_{0})} \right \} Ê\nonumberÊ\\[4mm]
&  \hspace*{-20mm} = \int_{\mathbf{y}} \frac{ f_{\mathbf{H}}(\mathbf{y}; \mathbf{x}_{1})  f_{\mathbf{H}}(\mathbf{y}; \mathbf{x}_{2}) }{  f_{\mathbf{H}}(\mathbf{y}; \mathbf{x}_{0})}d \mathbf{y}  \nonumberÊ\\[4mm]
&  \hspace*{-20mm} = \int_{\mathbf{y}} \frac{1}{(2 \pi \sigma^{2})^{M/2} }  
\frac{ \exp  \left(-\frac{1}{2\sigma^{2}} \| \mathbf{y} - \mathbf{H} \mathbf{x}_{1} \|^{2}_{2} \right)  \exp \left(-\frac{1}{2 \sigma^{2}}  \| \mathbf{y} - \mathbf{H} \mathbf{x}_{2} \|^{2}_{2} \right) }
{   \exp \left( -\frac{1}{2 \sigma^{2}} \| \mathbf{y} - \mathbf{H} \mathbf{x}_{0} \|^{2}_{2} \right)}d \mathbf{y}  \nonumber \\[4mm]
& \hspace*{-20mm} = \int_{\mathbf{y}} \frac{1}{(2 \pi \sigma^{2} )^{M/2} }   \exp \bigg( -\frac{1}{2 \sigma^{2}} \bigg[ \| \mathbf{y} \|^{2}_{2}- 2\mathbf{y}^{T} \mathbf{H} (\mathbf{x}_{1}+\mathbf{x}_{2}-\mathbf{x}_{0}) +  \|  \mathbf{H} \mathbf{x}_{1} \|^{2}_{2} + \|  \mathbf{H} \mathbf{x}_{2} \|^{2}_{2} -  \| \mathbf{H} \mathbf{x}_{0} \|^{2}_{2} \bigg] \bigg) d \mathbf{y} \nonumber \\[4mm]
& \hspace*{-20mm} = \int_{\mathbf{y}} \frac{1}{(2 \pi \sigma^{2} )^{M/2} }  \exp \bigg(-\frac{1}{2\sigma^{2}}  \bigg[ \| \mathbf{y} - \mathbf{H}(\mathbf{x}_{1}+\mathbf{x}_{2}-\mathbf{x}_{0})  \|^{2} + 2 \mathbf{x}_{0}^{T}\mathbf{H}^{T} \mathbf{H} (\mathbf{x}_{1} + \mathbf{x}_{2})  \nonumber \\[4mm]
& - 2 \mathbf{x}_{1}^{T}\mathbf{H}^{T} \mathbf{H} \mathbf{x}_{2} - 2 \mathbf{x}_{0}^{T} \mathbf{H}^{T} \mathbf{H} \mathbf{x}_{0} \bigg] \bigg) d \mathbf{y} \nonumber \\[4mm] 
& \hspace*{-20mm} =  \exp \bigg(-\frac{1}{2\sigma^{2}}  \bigg[+ 2 \mathbf{x}_{0}^{T}\mathbf{H}^{T} \mathbf{H} (\mathbf{x}_{1} + \mathbf{x}_{2})  
 - 2 \mathbf{x}_{1}^{T}\mathbf{H}^{T} \mathbf{H} \mathbf{x}_{2} - 2 \mathbf{x}_{0}^{T} \mathbf{H}^{T} \mathbf{H} \mathbf{x}_{0} \bigg] \bigg) \times \nonumber \\[4mm] 
& \underbrace{\int_{\mathbf{y}} \frac{1}{(2 \pi \sigma^{2} )^{M/2} }  \exp \bigg(-\frac{1}{2\sigma^{2}}  \bigg[ \| \mathbf{y} - \mathbf{H}(\mathbf{x}_{1}+\mathbf{x}_{2}-\mathbf{x}_{0})  \|^{2} \bigg] \bigg ) d \mathbf{y}}_{= 1} \nonumber \\[4mm] 
& \hspace*{-20mm}  =  \exp   \left( \frac{1}{\sigma^{2}} (\mathbf{x}_{1} - \mathbf{x}_{0})^{T}\mathbf{H}^{T} \mathbf{H} (\mathbf{x}_{2} - \mathbf{x}_{0}) \right). 
\end{align} 
Obviously, the kernel $R_{\mathcal{M}_{\text{SLM}}}(\cdot,\cdot)$ is differentiable in the sense of Definition \ref{def_differentiable_kernel} 
up to any order $m \in \mathbb{N}$, and therefore Theorem \ref{thm_der_repr_prop} applies 
to the RKHS $\mathcal{H}(\mathcal{M}_{\text{SLM}})$ (note that this RKHS consists of functions $f(\cdot): \mathcal{X}_{S} \rightarrow \mathbb{R}$). 

Let us now consider the SLM for the specific case $S=N$, for which the SLM coincides with the LGM. 
We will denote by 
\begin{equation} 
\label{equ_def_min_var_problem_LGM}
\mathcal{M}_{\text{LGM}} \triangleq \left( \mathcal{E}_{\text{LGM}}, \tilde{c}(\cdot), \mathbf{x}_{0} \right)
\end{equation} 
any minimum variance problem arising from the LGM $\mathcal{E}_{\text{LGM}}$ by fixing a prescribed bias function $\tilde{c}(\cdot): \mathbb{R}^{N} \rightarrow \mathbb{R}$ 
and parameter vector $\mathbf{x}_{0} \in \mathcal{X}_{S}$.\footnote{Note that we require $\mathbf{x}_{0} \in \mathcal{X}_{S}$ even if the parameter set of the LGM is $\mathbb{R}^{N} \supseteq \mathcal{X}_{S}$.}
The corresponding RKHS will be denoted by $\mathcal{H}(\mathcal{M}_{\text{LGM}})$, i.e., $\mathcal{H}(\mathcal{M}_{\text{LGM}})$ is the 
RKHS associated to the kernel $R_{\mathcal{M}_{\text{LGM}}}(\cdot,\cdot): \mathbb{R}^{N} \times \mathbb{R}^{N} \rightarrow \mathbb{R}$: 
\begin{equation} 
\label{equ_kernel_LGM}
R_{\mathcal{M}_{\text{LGM}}}(\mathbf{x}_{1}, \mathbf{x}_{2}) =  \exp   \left(\frac{1}{\sigma^{2}}  (\mathbf{x}_{1} - \mathbf{x}_{0})^{T} \mathbf{H}^{T} \mathbf{H} (\mathbf{x}_{2} - \mathbf{x}_{0}) \right). 
\end{equation} 
Note that the functional form of $R_{\mathcal{M}_{\text{LGM}}}(\mathbf{x}_{1}, \mathbf{x}_{2})$ is identical with that of $R_{\mathcal{M}_{\text{SLM}}}(\mathbf{x}_{1}, \mathbf{x}_{2}) $ (see \eqref{equ_kernel_SLM}). However, these two kernel functions differ in their domain, which is $\mathbb{R}^{N}Ê\times \mathbb{R}^{N}$ and $\mathcal{X}_{S} \times \mathcal{X}_{S}$, respectively. 

As can be verified easily, the restriction $\mathcal{M}_{\text{LGM}}\big|_{\mathcal{X}_{S}}$ of the minimum variance problem $\mathcal{M}_{\text{LGM}}$ given by \eqref{equ_def_min_var_problem_LGM} coincides with the minimum variance problem $\mathcal{M}_{\text{SLM}}$ given by \eqref{equ_def_min_var_problem_SLM} with prescribed bias $c(\cdot) = \tilde{c}(\cdot) \big|_{\mathcal{X}_{S}}$. We highlight the fact that we prescribe the estimator bias within minimum variance estimation for the SLM 
only for the smaller set $\mathcal{X}_{S}$, i.e., we require $b(\hat{x}_{k}(\cdot); \mathbf{x}) = c(\mathbf{x})$ only for $\mathbf{x} \in \mathcal{X}_{S} \subseteq \mathbb{R}^{N}$, whereas for the LGM we prescribe the bias for every $\mathbf{x} \in \mathbb{R}^{N}$.
\begin{theorem}
\label{thm_relation_RKHS_LGM_SLM}
Consider the minimum variance problem $\mathcal{M}_{\emph{LGM}}= \left( \mathcal{E}_{\emph{LGM}},\tilde{c}(\cdot),\mathbf{x}_{0} \right)$ with an arbitrary but fixed choice for 
$\sigma$, $M$, $N$, $\mathbf{H}$, prescribed bias $\tilde{c}(\cdot): \mathbb{R}^{N} \rightarrow \mathbb{R}$, and parameter vector $\mathbf{x}_{0} \in \mathcal{X}_{S}$. We then have that 
for any sparsity degree $S \in [N]$, the RKHS $\mathcal{H}(\mathcal{M}_{\emph{SLM}})$ associated with $\mathcal{M}_{\emph{SLM}}=\mathcal{M}_{\emph{LGM}}\big|_{\mathcal{X}_{S}}$
consists of the restrictions of all functions $f(\cdot): \mathbb{R}^{N} \rightarrow \mathbb{R}$ belonging to the RKHS $\mathcal{H}(\mathcal{M}_{\emph{LGM}})$ to the subdomain $\mathcal{X}_{S} \subseteq \mathbb{R}^{N}$, i.e., 
\begin{equation}
 \mathcal{H}(\mathcal{M}_{\emph{SLM}}) = \bigg\{ f(\cdot) \big|_{\mathcal{X}_{S}} \,\, \bigg| \,\,  f(\cdot) \in \mathcal{H}(\mathcal{M}_{\emph{LGM}}) \bigg\}. 
\end{equation} 
Moreover, the norm of any element $f(\cdot) \in \mathcal{H}(\mathcal{M}_{\emph{SLM}})$ is given by 
\begin{equation}
\label{equ_relation_norm_SLM_norm_LGM} 
\| f(\cdot) \|_{\mathcal{H}(\mathcal{M}_{\emph{SLM}})}   = 
\min_{\substack{g(\cdot) \in \mathcal{H}(\mathcal{M}_{\emph{LGM}}) \\  g(\cdot)\big|_{\mathcal{X}_{S}} = f(\cdot)} } \| g(\cdot) \|_{\mathcal{H}(\mathcal{M}_{\emph{LGM}})}
\end{equation}
\end{theorem} 
\begin{proof}
This result is obtained from the application of Theorem \ref{thm_reducing_domain_RKHS} to the specific RKHSs $\mathcal{H}(\mathcal{M}_{\text{SLM}})$ and $\mathcal{H}(\mathcal{M}_{\text{LGM}})$
since the associated kernels obviously satisfy $R_{\mathcal{M}_{\text{SLM}}}(\cdot,\cdot) =R_{\mathcal{M}_{\text{LGM}}}(\cdot,\cdot) \big|_{\mathcal{X}_{S} \times \mathcal{X}_{S}}$. 
\end{proof}

The RKHS $\mathcal{H}(\mathcal{M}_{\text{LGM}})$ associated with the minimum variance problem $\mathcal{M}_{\text{LGM}}$ has already been analyzed in \cite{Duttweiler73b}. 
However, we will now present an alternative characterization of the RKHS $\mathcal{H}(\mathcal{M}_{\text{LGM}})$ that is based on a congruence with 
a specific RKHS that has some pleasing properties. 
\begin{theorem}
\label{thm_isometry_LGM}
Consider the RKHS $\mathcal{H}(\mathcal{M}_{\emph{LGM}})$ associated with  $\mathcal{M}_{\emph{LGM}} \triangleq \left( \mathcal{E}_{\emph{LGM}}, \tilde{c}(\cdot), \mathbf{x}_{0} \right)$ 
with an arbitrary system matrix $\mathbf{H} \in \mathbb{R}^{M \times N}$, not necessarily satisfying \eqref{equ_spark_cond}, 
and denote its thin SVD by $\mathbf{H} = \mathbf{U} \mathbf{\Sigma} \mathbf{V}^{T}$. 
Then we have that any function $f(\cdot) \in \mathcal{H}(\mathcal{M}_{\emph{LGM}})$ is invariant w.r.t.\ translations by vectors $\mathbf{x}' \in \mathbb{R}^{N}$ 
that belong to the null-space of $\mathbf{H}$, i.e., 
\begin{equation} 
\label{equ_invar_prop_LGM_kernel_H}
f(\mathbf{x}) = f(\mathbf{x} + \mathbf{x}')
\end{equation}
for any $\mathbf{x}' \in \mathcal{N}(\mathbf{H})$, $\mathbf{x}\in \mathbb{R}^{N}$.  
Furthermore, the RKHS $\mathcal{H}(\mathcal{M}_{\emph{LGM}})$ is isometric to the RKHS $\mathcal{H}(R_{g}^{(D)})$ for $D=\rank( \mathbf{H} )$, 
which is defined via the kernel $R_{g}^{(D)}(\cdot,\cdot): \mathbb{R}^{D} \times \mathbb{R}^{D} \rightarrow \mathbb{R}$: 
\begin{equation}
\label{equ_def_kernel_R_g}
R_{g}^{(D)}(\mathbf{x}_{1}, \mathbf{x}_{2}) \triangleq \exp \big( \mathbf{x}_{1}^{T} \mathbf{x}_{2} \big).
\end{equation}
The mapping $\mathsf{K}_{g}[\cdot]: \mathcal{H}(R_{g}^{(D)}) \rightarrow \mathcal{H}(\mathcal{M}_{\emph{LGM}})$ given by 
\begin{equation} 
\label{equ_def_isometry_LGM}
f(\cdot) \mapsto \widetilde{f}(\cdot) =  \mathsf{K}_{g}[f(\cdot)] :  \widetilde{f}(\mathbf{x}') 
= f\bigg(\frac{1}{\sigma} \widetilde{\mathbf{H}}^{\dagger} \mathbf{x}' \bigg) \exp\left( \frac{1}{2 \sigma^{2}} \| \mathbf{H} \mathbf{x}_{0} \|^{2}_{2} - \frac{1}{\sigma^{2}} (\mathbf{x}')^{T}  \mathbf{H}^{T}\mathbf{H} \mathbf{x}_{0}\right),
\end{equation} 
where $\widetilde{\mathbf{H}} \triangleq \mathbf{V} \mathbf{ \Sigma}^{-1}$, 
is a congruence from $\mathcal{H}(R_{g}^{(D)})$ to $\mathcal{H}(\mathcal{M}_{\emph{LGM}})$. 
The inverse mapping $\mathsf{K}^{-1}_{g}[\cdot]:  \mathcal{H}(\mathcal{M}_{\emph{LGM}}) \rightarrow \mathcal{H}(R_{g}^{(D)})$  given by 
\begin{equation} 
\label{equ_def_isometry_LGM_inverse}
f(\cdot) \mapsto \widetilde{f}(\cdot) =  \mathsf{K}^{-1}_{g}[f(\cdot)] :  \widetilde{f}(\mathbf{x}') 
= f\big(\sigma \widetilde{\mathbf{H}} \mathbf{x}' \big) \exp\left( -\frac{1}{2 \sigma^{2}} \| \mathbf{H} \mathbf{x}_{0} \|^{2}_{2} + \frac{1}{\sigma} (\mathbf{x}')^{T}  \widetilde{\mathbf{H}}^{\dagger} \mathbf{x}_{0}\right),
\end{equation}  
is a congruence from $\mathcal{H}(\mathcal{M}_{\emph{LGM}})$ to $\mathcal{H}(R_{g}^{(D)})$. 
\end{theorem}

\begin{proof}
Consider the function $f_{\mathbf{x}_{1}}(\cdot) = R_{\mathcal{M}_{\text{LGM}}}(\cdot, \mathbf{x}_{1} )$, where $\mathbf{x}_{1} \in \mathbb{R}^{N}$ is arbitrary but fixed.
Then, for any $\mathbf{x}' \in \mathcal{N}(\mathbf{H})$ and $\mathbf{x} \inÊ\mathbb{R}^{N}$,
\begin{align}
f_{\mathbf{x}_{1}}(\mathbf{x} + \mathbf{x}') & = \exp   \left( \frac{1}{\sigma^{2}}  (\mathbf{x} + \mathbf{x}' - \mathbf{x}_{0})^{T}\mathbf{H}^{T} \mathbf{H} (\mathbf{x}_{1} - \mathbf{x}_{0}) \right) \nonumber \\[4mm]
& = \exp   \left( \frac{1}{\sigma^{2}}  (\mathbf{x} - \mathbf{x}_{0})^{T}\mathbf{H}^{T} \mathbf{H} (\mathbf{x}_{1} - \mathbf{x}_{0})
+ \underbrace{\big( (\mathbf{x}_{1} - \mathbf{x}_{0})^{T} \mathbf{H}^{T} \mathbf{H} \mathbf{x}' \big)^{T}}_{=0} \right) \nonumber \\[4mm]
& = \exp\left( \frac{1}{\sigma^{2}} (\mathbf{x} - \mathbf{x}_{0})^{T}\mathbf{H}^{T} \mathbf{H} (\mathbf{x}_{1} - \mathbf{x}_{0}) \right) = f_{\mathbf{x}_{1}}(\mathbf{x}).
\end{align} 
From this, it follows straightforwardly that any function $f(\cdot) \in \linspan \left \{ R_{\mathcal{M}_{\text{LGM}}}(\cdot, \mathbf{x} ) \right \}_{\mathbf{x} \in \mathbb{R}^{N}}$ satisfies \eqref{equ_invar_prop_LGM_kernel_H}. 
Consider now a general $f(\cdot) \in  \mathcal{H}(\mathcal{M}_{\text{LGM}})$ and fix an arbitrary $\mathbf{x} \in \mathbb{R}^{N}$ and $\mathbf{x}' \in \mathcal{N}(\mathbf{H})$. 
Then, for any $f'(\cdot) \in \linspan \left \{ R_{\mathcal{M}_{\text{LGM}}}(\cdot, \mathbf{x} ) \right \}_{\mathbf{x} \in \mathbb{R}^{N}}$, we have 
\begin{align} 
\label{equ_inver_prop_LGM_arb_f_kernel_H}
f(\mathbf{x}+ \mathbf{x}') - f(\mathbf{x}) & = \big[f(\mathbf{x}+ \mathbf{x}')-f'(\mathbf{x}+\mathbf{x}')\big] - \big[f(\mathbf{x})-f'(\mathbf{x})\big]+ f'(\mathbf{x}+\mathbf{x}')-f'(\mathbf{x}) \nonumber \\[4mm]Ê
 &\stackrel{(a)}{=} \big[f(\mathbf{x}+ \mathbf{x}')-f'(\mathbf{x}+\mathbf{x}')\big] - \big[f(\mathbf{x})-f'(\mathbf{x})\big] 
\end{align} 
where $(a)$ follows from the fact that $f'(\cdot)$ satisfies \eqref{equ_invar_prop_LGM_kernel_H}. 
Since according to Theorem \ref{thm_constr_RKHS_closure_linear_span}, the linear space $\linspan \left \{ R_{\mathcal{M}_{\text{LGM}}}(\cdot, \mathbf{x} ) \right \}_{\mathbf{x} \in \mathbb{R}^{N}}$ is dense 
in $\mathcal{H}(\mathcal{M}_{\text{LGM}})$ we can choose for any $\varepsilon>0$ a function $f'(\cdot) \in \linspan \left \{ R_{\mathcal{M}_{\text{LGM}}}(\cdot, \mathbf{x} ) \right \}_{\mathbf{x} \in \mathbb{R}^{N}}$ such that $\| f'(\cdot) - f(\cdot) \|_{\mathcal{H}(\mathcal{M}_{\text{LGM}})} \leq \varepsilon$ which implies via Theorem \ref{thm_cauchy_schwarz} and the 
reproducing property \eqref{equ_reproduction_property} that 
\begin{align}
\big|f(\mathbf{x}+ \mathbf{x}')-f'(\mathbf{x}+\mathbf{x}')\big| & \stackrel{\eqref{equ_reproduction_property}}{=}
 \bigg| \big\langle f(\cdot) - f'(\cdot), R_{\mathcal{M}_{\text{LGM}}}(\cdot,\mathbf{x}+ \mathbf{x}') \big\rangle_{\mathcal{H}(\mathcal{M}_{\text{LGM}})} \bigg| \nonumber \\[4mm]
& \hspace*{-30mm} \stackrel{\eqref{equ_cauchy_schwarz}}{\leq} \| f'(\cdot) - f(\cdot) \|_{\mathcal{M}_{\text{LGM}}} \| R_{\mathcal{M}_{\text{LGM}}}(\cdot,\mathbf{x}+\mathbf{x}') \|_{\mathcal{H}(\mathcal{M}_{\text{LGM}})} \leq \varepsilon  \| R_{\mathcal{M}_{\text{LGM}}}(\cdot,\mathbf{x}+\mathbf{x}') \|_{\mathcal{H}(\mathcal{M}_{\text{LGM}})} \nonumber \\[4mm]
  & \hspace*{-30mm} \stackrel{\eqref{equ_reproduction_property}}{=} \varepsilon \sqrt{ R_{\mathcal{M}_{\text{LGM}}}(\mathbf{x}+\mathbf{x}',\mathbf{x}+\mathbf{x}')}, 
\end{align}
and similarly that 
\begin{align}
\big|f(\mathbf{x})-f'(\mathbf{x})\big| \leq  \varepsilon \sqrt{ R_{\mathcal{M}_{\text{LGM}}}(\mathbf{x},\mathbf{x})}. 
\end{align}
Thus, we can make the last expression of \eqref{equ_inver_prop_LGM_arb_f_kernel_H} arbitrarily small by a suitable choice for $f'(\cdot)$. 
This concludes the proof of \eqref{equ_invar_prop_LGM_kernel_H}. 

For the following, we note that based on the thin SVD $\mathbf{H} = \mathbf{U} \mathbf{\Sigma} \mathbf{V}^{T}$ of $\mathbf{H}$, we have 
\begin{equation}
\label{equ_proof_isometr_R_g_quadr_form} 
\mathbf{H}^{T}\mathbf{H} = \mathbf{V} \mathbf{\Sigma}^{2} \mathbf{V}^{T} =  \big( \widetilde{\mathbf{H}}^{\dagger} \big)^{T}\widetilde{\mathbf{H}}^{\dagger}.
\end{equation} 
Consider the two sets of functions given by
\begin{equation} 
\mathcal{A} \triangleq \bigg \{ g_{\mathbf{x}} (\cdot) \triangleq \exp\left( \frac{1}{2 \sigma^{2}} \| \mathbf{H} \mathbf{x}_{0} \|^{2}_{2} - 
 \frac{1}{\sigma}\mathbf{x}^{T} \mathbf{H}^{T} \mathbf{H} \mathbf{x}_{0} \right)R_{g}^{(D)}\big(\cdot, \widetilde{\mathbf{H}}^{\dagger} \mathbf{x}\big) \bigg \}_{\mathbf{x} \in \mathbb{R}^{N}}
 \label{equ_proof_isometr_R_g_set_A_in_R_LGM}
\end{equation} 
and 
\begin{equation} 
 \label{equ_proof_isometr_R_g_set_b_in_R_LGM}
\mathcal{B} \triangleq \big \{ f_{\mathbf{x}}(\cdot) \triangleq  R_{\mathcal{M}_{\text{LGM}}}(\cdot, \sigma \mathbf{x} ) \big \}_{\mathbf{x} \in \mathbb{R}^{N}}.
\end{equation} 
We have obviously that 
\begin{equation}
\mathcal{B} = \big\{ R_{\mathcal{M}_{\text{LGM}}}(\cdot, \mathbf{x} ) \big \}_{\mathbf{x} \in \mathbb{R}^{N}}
\end{equation} 
and in turn 
\begin{equation} 
\linspan \big\{\mathcal{B} \big\}=  \linspan  \big \{ R_{\mathcal{M}_{\text{LGM}}}(\cdot, \mathbf{x} ) \big \}_{\mathbf{x} \in \mathbb{R}^{N}}.
\end{equation} 
Furthermore, since for any $\mathbf{x} \in \mathbb{R}^{D}$ and $\mathbf{x}' \triangleq \mathbf{V} \mathbf{\Sigma} \mathbf{x} \in \mathbb{R}^{N}$ 
(note that $\widetilde{\mathbf{H}}^{\dagger}  \mathbf{x}' = \mathbf{x}$)  we have 
\begin{equation}
\label{equ_proof_charac_RKHS_LGM_equl_funcs_R_N_R_D}
R_{g}^{(D)}(\cdot, \widetilde{\mathbf{H}}^{\dagger} \mathbf{x}') = R_{g}^{(D)}(\cdot, \mathbf{x}),
\end{equation}
it holds that 
\begin{equation}
\label{equ_proof_charac_RKHS_LGM_equl_funcs_R_N_R_D_subset_1}
\bigg\{ R_{g}^{(D)}(\cdot, \mathbf{x})  \bigg\}_{\mathbf{x} \in \mathbb{R}^{D}} \subseteq \bigg\{ R_{g}^{(D)}(\cdot, \widetilde{\mathbf{H}}^{\dagger} \mathbf{x}')   \bigg\}_{\mathbf{x}' \in \mathbb{R}^{N}}.
\end{equation} 
Similarly, since any $\mathbf{x}' \in \mathbb{R}^{N}$ and $\mathbf{x} \triangleq \widetilde{\mathbf{H}}^{\dagger}  \mathbf{x}' \in \mathbb{R}^{D}$ 
also satisfy \eqref{equ_proof_charac_RKHS_LGM_equl_funcs_R_N_R_D} it holds that 
\begin{equation}
\label{equ_proof_charac_RKHS_LGM_equl_funcs_R_N_R_D_subset_2}
\bigg\{ R_{g}^{(D)}(\cdot, \widetilde{\mathbf{H}}^{\dagger} \mathbf{x}')   \bigg\}_{\mathbf{x}' \in \mathbb{R}^{N}}  \subseteq \bigg\{ R_{g}^{(D)}(\cdot, \mathbf{x})  \bigg\}_{\mathbf{x} \in \mathbb{R}^{D}}.
\end{equation}
Combining \eqref{equ_proof_charac_RKHS_LGM_equl_funcs_R_N_R_D_subset_1} and \eqref{equ_proof_charac_RKHS_LGM_equl_funcs_R_N_R_D_subset_2} yields
\begin{equation}
\bigg\{ R_{g}^{(D)}(\cdot, \mathbf{x})  \bigg\}_{\mathbf{x} \in \mathbb{R}^{D}}  =  \bigg\{ R_{g}^{(D)}\big(\cdot, \widetilde{\mathbf{H}}^{\dagger} \mathbf{x}' \big)   \bigg\}_{\mathbf{x}' \in \mathbb{R}^{N}},
\end{equation} 
and in turn 
\begin{equation}
\linspan \{\mathcal{A} \}=  \linspan  \left \{R_{g}^{(D)}(\cdot, \mathbf{x} ) \right \}_{\mathbf{x} \in \mathbb{R}^{D}},
\end{equation} 
since the weights $\exp\bigg( \frac{1}{2 \sigma^{2}} \| \mathbf{H} \mathbf{x}_{0} \|^{2}_{2} -  \frac{1}{\sigma}\mathbf{x}^{T} (\mathbf{H}^{T} \mathbf{H}) \mathbf{x}_{0} \bigg)$ 
appearing in the definition of $\mathcal{A}$ have no influence on the set $\linspan \{\mathcal{A} \}$. 
Thus, we have shown that $\mathcal{A}$ and $\mathcal{B}$ span the two RKHSs $\mathcal{H}(R_{g}^{(D)})$ and $\mathcal{H}(\mathcal{M}_{\text{LGM}})$, respectively.
 
For any two vectors $\mathbf{x}_{1}, \mathbf{x}_{2} \in \mathbb{R}^{N}$:
\begin{align}
\label{equ_proof_LGM_isom_Rg_inner_prod_preserv}
 \big\langle f_{\mathbf{x}_{1}}(\cdot) ,  f_{\mathbf{x}_{2}}(\cdot) \big\rangle_{\mathcal{H}(\mathcal{M}_{\text{LGM}})} & =
 \big\langle R_{\mathcal{M}_{\text{LGM}}}(\cdot,\sigma \mathbf{x}_{1}), R_{\mathcal{M}_{\text{LGM}}}(\cdot,Ê\sigma \mathbf{x}_{2}) \big\rangle_{\mathcal{H}(\mathcal{M}_{\text{LGM}})} \nonumber  \\[4mm]
 & \hspace*{-30mm} \stackrel{\eqref{equ_reproduction_property}}{=} R_{\mathcal{M}_{\text{LGM}}}(\sigma \mathbf{x}_{1},\sigma \mathbf{x}_{2}) \nonumber \\[4mm]
& \hspace*{-30mm} = \exp   \left( \frac{1}{\sigma^{2}} ( \sigma \mathbf{x}_{1} - \mathbf{x}_{0})^{T} \mathbf{H}^{T} \mathbf{H} (\sigma \mathbf{x}_{2} - \mathbf{x}_{0}) \right) \nonumber \\[4mm]Ê
&\hspace*{-30mm} = \exp\left( \frac{1}{2 \sigma^{2}} \| \mathbf{H} \mathbf{x}_{0} \|^{2}_{2} -  \frac{1}{\sigma} \mathbf{x}_{1}^{T} \mathbf{H}^{T} \mathbf{H}\mathbf{x}_{0}+ \mathbf{x}_{1}^{T} \mathbf{H}^{T} \mathbf{H} \mathbf{x}_{2} + \frac{1}{2 \sigma^{2}} \| \mathbf{H} \mathbf{x}_{0} \|^{2}_{2} - \frac{1}{\sigma} \mathbf{x}_{2}^{T} \mathbf{H}^{T} \mathbf{H} \mathbf{x}_{0} \right) \nonumber \\[4mm]
&\hspace*{-30mm} = \exp\left( \frac{1}{2 \sigma^{2}} \| \mathbf{H} \mathbf{x}_{0} \|^{2}_{2} -  \frac{1}{\sigma} \mathbf{x}_{1}^{T} \mathbf{H}^{T} \mathbf{H}\mathbf{x}_{0} \right)  \exp   \left(  \mathbf{x}_{1}^{T} \mathbf{H}^{T} \mathbf{H} \mathbf{x}_{2}  \right)\exp\left( \frac{1}{2 \sigma^{2}} \| \mathbf{H} \mathbf{x}_{0} \|^{2}_{2} - \frac{1}{\sigma} \mathbf{x}_{2}^{T} \mathbf{H}^{T} \mathbf{H} \mathbf{x}_{0} \right) \nonumber \\[4mm]
&\hspace*{-30mm} \stackrel{\eqref{equ_proof_isometr_R_g_quadr_form}}{=} \exp\left( \frac{1}{2\sigma^{2}} \| \mathbf{H} \mathbf{x}_{0} \|^{2}_{2} -  \frac{1}{\sigma} \mathbf{x}_{1}^{T} \mathbf{H}^{T} \mathbf{H} \mathbf{x}_{0} \right)  \exp   \left(  \mathbf{x}_{1}^{T}( \widetilde{\mathbf{H}}^{\dagger})^{T}\widetilde{\mathbf{H}}^{\dagger} \mathbf{x}_{2}  \right)\exp\left( \frac{1}{2 \sigma^{2}} \| \mathbf{H} \mathbf{x}_{0} \|^{2}_{2} - \frac{1}{\sigma} \mathbf{x}_{2}^{T} \mathbf{H}^{T} \mathbf{H} \mathbf{x}_{0} \right) \nonumber \\[4mm]
& \hspace*{-30mm} \stackrel{\eqref{equ_def_kernel_R_g}}{=}   \exp\left( \frac{1}{2 \sigma^{2}} \| \mathbf{H} \mathbf{x}_{0} \|^{2}_{2} - \frac{1}{\sigma} \mathbf{x}_{1}^{T} \mathbf{H}^{T} \mathbf{H} \mathbf{x}_{0} \right)R_{g}^{(D)} \big( \widetilde{\mathbf{H}}^{\dagger} \mathbf{x}_{1}, \widetilde{\mathbf{H}}^{\dagger} \mathbf{x}_{2}\big) 
\exp\left( \frac{1}{2 \sigma^{2}} \| \mathbf{H} \mathbf{x}_{0} \|^{2}_{2} -  \frac{1}{\sigma}\mathbf{x}_{2}^{T} \mathbf{H}^{T} \mathbf{H} \mathbf{x}_{0} \right)  \nonumber \\[4mm]Ê
& \hspace*{-30mm} \stackrel{\eqref{equ_reproduction_property}}{=}  \bigg\langle \exp\left( \frac{1}{2 \sigma^{2}} \| \mathbf{H} \mathbf{x}_{0} \|^{2}_{2} - \frac{1}{\sigma} \mathbf{x}_{1}^{T} \mathbf{H}^{T} \mathbf{H} \mathbf{x}_{0} \right)R_{g}^{(D)}\big(\cdot, \widetilde{\mathbf{H}}^{\dagger} \mathbf{x}_{1}\big), \nonumber \\[4mm]
& \exp\left( \frac{1}{2 \sigma^{2}} \| \mathbf{H} \mathbf{x}_{0} \|^{2}_{2} -  \frac{1}{\sigma}\mathbf{x}_{2}^{T} \mathbf{H}^{T} \mathbf{H} \mathbf{x}_{0} \right)R_{g}^{(D)}\big(\cdot,\widetilde{\mathbf{H}}^{\dagger} \mathbf{x}_{2}\big) \bigg\rangle_{\mathcal{H}(R_{g}^{(D)})} \nonumber \\[4mm]Ê
&\hspace*{-30mm}  = \big\langle g_{\mathbf{x}_{1}}(\cdot) , g_{\mathbf{x}_{2}}(\cdot) \big\rangle_{\mathcal{H}(R_{g}^{(D)})}.
\end{align}

Now, for an arbitrary $\mathbf{x} \in \mathbb{R}^{N}$, let us denote by $h_{\mathbf{x}}(\cdot) \triangleq \mathsf{K}_{g}[g_{\mathbf{x}}(\cdot)] \in  \mathcal{H}(\mathcal{M}_{\text{LGM}})$ the image of the 
function $g_{\mathbf{x}}(\cdot) \triangleq \exp\left( \frac{1}{2 \sigma^{2}} \| \mathbf{H} \mathbf{x}_{0} \|^{2}_{2} -  \frac{1}{\sigma}\mathbf{x}^{T} (\mathbf{H}^{T} \mathbf{H}) \mathbf{x}_{0} \right)R_{g}^{(D)}\big(\cdot, \widetilde{\mathbf{H}}^{\dagger} \mathbf{x}\big)   \in \mathcal{H}(R_{g}^{(D)})$ (see \eqref{equ_proof_isometr_R_g_set_A_in_R_LGM}) under the mapping $\mathsf{K}_{g}[\cdot]$ defined in \eqref{equ_def_isometry_LGM}. 
We have
\begin{align}
\label{equ_proof_LGM_isom_Rg_mapping_is_congruence_f_g}
h_{\mathbf{x}}(\mathbf{x}') & =  g_{\mathbf{x}}\bigg(\frac{1}{\sigma} \widetilde{\mathbf{H}}^{\dagger} \mathbf{x}' \bigg) \exp\left( \frac{1}{2 \sigma^{2}} \| \mathbf{H} \mathbf{x}_{0} \|^{2}_{2} - \frac{1}{\sigma^{2}} (\mathbf{x}')^{T}  \mathbf{H}^{T}\mathbf{H} \mathbf{x}_{0}\right)
 \nonumber \\[4mm]Ê
& \hspace*{-10mm} \stackrel{\eqref{equ_proof_isometr_R_g_set_A_in_R_LGM}}{=}   \exp\left( \frac{1}{2 \sigma^{2}} \| \mathbf{H} \mathbf{x}_{0} \|^{2}_{2} - 
 \frac{1}{\sigma}\mathbf{x}^{T} \mathbf{H}^{T} \mathbf{H} \mathbf{x}_{0} \right)R_{g}^{(D)}\bigg(\frac{1}{\sigma} \widetilde{\mathbf{H}}^{\dagger} \mathbf{x}', \widetilde{\mathbf{H}}^{\dagger} \mathbf{x}\bigg) \exp\left( \frac{1}{2 \sigma^{2}} \| \mathbf{H} \mathbf{x}_{0} \|^{2}_{2} - \frac{1}{\sigma^{2}} (\mathbf{x}')^{T}  \mathbf{H}^{T}\mathbf{H} \mathbf{x}_{0}\right)
 \nonumber  \\[4mm]
 & \hspace*{-10mm} \stackrel{\eqref{equ_def_kernel_R_g}}{=}   \exp\left( \frac{1}{ \sigma^{2}} \| \mathbf{H} \mathbf{x}_{0} \|^{2}_{2} - 
 \frac{1}{\sigma}\mathbf{x}^{T} \mathbf{H}^{T} \mathbf{H} \mathbf{x}_{0} \right) \exp\left( \frac{1}{\sigma} \big(\mathbf{x}'\big)^{T} \big(\widetilde{\mathbf{H}}^{\dagger}\big)^{T} \widetilde{\mathbf{H}}^{\dagger} \mathbf{x} \right) \exp\left( - \frac{1}{\sigma^{2}} (\mathbf{x}')^{T}  \mathbf{H}^{T}\mathbf{H} \mathbf{x}_{0}\right)
 \nonumber  \\[4mm]
  & \hspace*{-10mm} \stackrel{\eqref{equ_proof_isometr_R_g_quadr_form}}{=}   \exp\left( \frac{1}{ \sigma^{2}} \| \mathbf{H} \mathbf{x}_{0} \|^{2}_{2} - 
 \frac{1}{\sigma}\mathbf{x}^{T} \mathbf{H}^{T} \mathbf{H} \mathbf{x}_{0} +\frac{1}{\sigma} \big(\mathbf{x}'\big)^{T} \mathbf{H}^{T} \mathbf{H} \mathbf{x}- \frac{1}{\sigma^{2}} (\mathbf{x}')^{T}  \mathbf{H}^{T}\mathbf{H} \mathbf{x}_{0}\right)
 \nonumber  \\[4mm]
   & \hspace*{-10mm} =   \exp\left( \frac{1}{\sigma^{2}} (\mathbf{x}'- \mathbf{x}_{0})^{T} \mathbf{H}^{T} \mathbf{H} (\sigma \mathbf{x} - \mathbf{x}_{0}) \right) 
\stackrel{\eqref{equ_kernel_SLM}}{=}  R_{\mathcal{M}_{\text{LGM}}}(\mathbf{x}', \sigma \mathbf{x} ) 
= f_{\mathbf{x}}(\mathbf{x}').
\end{align} 
The fact that the mapping $\mathsf{K}_{g}[\cdot]$ defined in \eqref{equ_def_isometry_LGM} is a congruence from $\mathcal{H}(R_{g}^{(D)})$ to $\mathcal{H}(\mathcal{M}_{\text{LGM}})$ follows then 
from Theorem \ref{thm_suff_con_congruence_RKHS}, since for every argument $\mathbf{x} \in \mathbb{R}^{N}$ the 
function value $\mathsf{K}_{g}[g(\cdot)](\mathbf{x})$ depends continuously on the function value $g\bigg(\frac{1}{\sigma} \widetilde{\mathbf{H}}^{\dagger} \mathbf{x} \bigg) $, which implies that the image $\mathsf{K}_{g}[g(\cdot)]$ of a function $g(\cdot) \in \mathcal{H}(R_{g}^{(D)})$ which is the pointwise limit of a sequence $\big\{ g_{l}(\cdot) \in \linspan \{ \mathcal{A} \} \big\}_{l \rightarrow \infty}$ is the pointwise limit of the functions $\mathsf{K}_{g}[g_{l}(\cdot)] \in \mathcal{H}(\mathcal{M}_{\text{LGM}})$.
\end{proof}

In the following, we will make use of the obvious facts summarized in 
\begin{lemma}
\label{lem_R_g_is_differentiable} 
The kernel $R_{g}^{(D)}(\cdot,\cdot): \mathbb{R}^{D} \times \mathbb{R}^{D} \rightarrow \mathbb{R}$ defined in \eqref{equ_def_kernel_R_g} is differentiable up to any order $m \in \mathbb{N}$ (see Definition \ref{def_differentiable_kernel}).
We have that for any vector $\mathbf{x}_{c} \in \mathbb{R}^{D}$, index set $\mathcal{K} \subseteq [D]$, and $\varepsilon >0$, the $\varepsilon$-$\mathcal{K}$-neighborhood of $\mathbf{x}_{c}$ (see \eqref{equ_def_epsilon_K_neighborhood}) belongs
to the domain of $R_{g}^{(D)}(\cdot,\cdot)$, i.e., $\mathcal{N}^{\mathcal{K}}_{\mathbf{x}_{c}}(\varepsilon) \subseteq \mathbb{R}^{D}$.
\end{lemma}

The RKHS $\mathcal{H}(R_{g}^{(D)})$ is entirely characterized by
\begin{theorem}
\label{thm_entire_character_H_R_g}
The RKHS $\mathcal{H}(R_{g}^{(D)})$ is separable and it contains the 
functions $g^{(\mathbf{p})}(\cdot): \mathbb{R}^{D} \rightarrow \mathbb{R}$ given by 
\begin{equation} 
g^{(\mathbf{p})} (\mathbf{x}) \triangleq  \frac{1}{\sqrt{\mathbf{p}!}} \frac{ \partial^{\mathbf{p}} R_{g}^{(D)}(\mathbf{x}, \mathbf{x}_{2})}{\partial \mathbf{x}_{2}^{\mathbf{p}}}\bigg|_{\mathbf{x}_{2} = \mathbf{0}} = \frac{1}{\sqrt{\mathbf{p}!}} \mathbf{x}^{\mathbf{p}},
\end{equation}
where $\mathbf{p} \in \mathbb{Z}_{+}^{D}$ is an arbitrary multi-index, i.e., 
\begin{equation}
\label{equ_part_der_element_RKHS_Rg}
\mathbf{p} \in \mathbb{Z}_{+}^{D} \Rightarrow g^{(\mathbf{p})}(\cdot) \in \mathcal{H}(R_{g}^{(D)}),
\end{equation}
with $\mathbf{p}! \triangleq \prod_{l \in [D]} p_{l}!$ and $\mathbf{x}^{\mathbf{p}} \triangleq \prod_{l \in [D]} (x_{l})^{p_l}$. 
The inner product of an arbitrary function $f(\cdot) \in \mathcal{H}(R_{g}^{(D)})$ with $g^{(\mathbf{p})}(\cdot)$ is given by 
\begin{equation} 
\label{equ_inner_prod_part_der_RKHS_Rg}
\big\langle f (\cdot), g^{(\mathbf{p})}(\cdot) \big\rangle_{\mathcal{H}(R_{g}^{(D)})} = \frac{1}{\sqrt{\mathbf{p}!}} \frac{ \partial^{\mathbf{p}} f(\mathbf{x})} {\partial \mathbf{x}^{\mathbf{p}}} \bigg|_{\mathbf{x} = 0}. 
\end{equation} 
Moreover, the set $\big\{ g^{(\mathbf{p})}(\cdot) \in \mathcal{H}(R_{g}^{(D)}) \big\}_{\mathbf{p} \in \mathbb{Z}_{+}^{D}}$ forms an ONB for $\mathcal{H}(R_{g}^{(D)})$. Thus, 
a function $f(\cdot): \mathbb{R}^{D} \rightarrow \mathbb{R}$ belongs to $\mathcal{H}(R_{g}^{(D)})$, i.e., $f(\cdot)  \in \mathcal{H}(R_{g}^{(D)})$ if and only if it can be written pointwise as 
\begin{equation}
\label{equ_series_repr_RKHS_R_g}
f(\mathbf{x}) = \sum_{\mathbf{p} \in \mathbb{Z}_{+}^{D}} a[\mathbf{p}] g^{(\mathbf{p})}(\mathbf{x}) =   \sum_{\mathbf{p} \in \mathbb{Z}_{+}^{D}} a[\mathbf{p}] \frac{1}{\sqrt{\mathbf{p}!}} \mathbf{x}^{\mathbf{p}},
\end{equation} 
with  a unique coefficient sequence $a[\mathbf{p}] \in \ell^{2}(\mathbb{Z}_{+}^{D})$. 
\end{theorem} 
\begin{proof}
The separability of $\mathcal{H}(R_{g}^{(D)})$ follows from \cite[Theorem 7]{HeinRKHS2004}, since obviously the domain $\mathbb{R}^{D}$ is separable and the kernel $R_{g}^{(D)}(\cdot,\cdot)$ is continuous. 
The validity of \eqref{equ_part_der_element_RKHS_Rg} and \eqref{equ_inner_prod_part_der_RKHS_Rg} follows from Lemma \ref{lem_R_g_is_differentiable} and Theorem \ref{thm_der_repr_prop}. 
By the power (Taylor) series representation of the exponential function \cite{RudinBook}, we have that the kernel $R_{g}^{(D)}(\cdot,\cdot)$ can be written pointwise as 
\begin{equation} 
R_{g}^{(D)}(\mathbf{x}_{1}, \mathbf{x}_{2}) = \exp \big( \mathbf{x}_{1}^{T} \mathbf{x}_{2} \big) = 
\sum_{\mathbf{p} \in \mathbb{Z}_{+}^{D}} \frac{1}{\mathbf{p}!} \mathbf{x}_{1}^{\mathbf{p}} \mathbf{x}_{2}^{\mathbf{p}} = 
\sum_{\mathbf{p} \in \mathbb{Z}_{+}^{D}} g^{(\mathbf{p})}(\mathbf{x}_{1}) g^{(\mathbf{p})}(\mathbf{x}_{2}).
\end{equation} 
This implies via Theorem \ref{thm_pointwise_series_kernel_ONB} that the set $\big\{ g^{(\mathbf{p})}(\cdot) \in \mathcal{H}(R_{g}^{(D)}) \big\}_{\mathbf{p} \in \mathbb{Z}_{+}^{D}}$ forms an ONB for $\mathcal{H}(R_{g}^{(D)})$, since it can be easily verified by \eqref{equ_inner_prod_part_der_RKHS_Rg} that $\big\langle g^{(\mathbf{p})}(\cdot),g^{(\mathbf{p}')}(\cdot) \big\rangle_{\mathcal{H}(R_{g}^{(D)})} = \delta_{\mathbf{p}, \mathbf{p}'}$.
The series representation in \eqref{equ_series_repr_RKHS_R_g} finally follows from Theorem \ref{thm_isometry_hilbert_space_coeffs_space}.  
\end{proof}

A useful consequence of Theorem \ref{thm_entire_character_H_R_g} is obtained by Theorem \ref{thm_isometry_hilbert_space_coeffs_space}: 
\begin{corollary}
\label{cor_R_g_congr_sequence_space} 
The mapping $\mathsf{K}_{a}[\cdot]:  \ell^{2}\big(\mathbb{Z}_{+}^{D}\big) \rightarrow \mathcal{H}(R_{g}^{(D)}):$ 
\begin{equation} 
\label{equ_cor_R_g_congr_sequence_space_a_p}
\mathsf{K}_{a}\big[a[\mathbf{p}]\big](\mathbf{x}) = \sum_{\mathbf{p} \in \mathbb{Z}_{+}^{D}} a[\mathbf{p}] g^{(\mathbf{p})}(\mathbf{x})  = \sum_{\mathbf{p} \in \mathbb{Z}_{+}^{D}} a[\mathbf{p}] \frac{1}{\sqrt{\mathbf{p}!}} \mathbf{x}^{\mathbf{p}} \quad \mbox{, $\mathbf{x} \in \mathbb{R}^{D}$,}
\end{equation}
is a congruence from $\ell^{2}\big(\mathbb{Z}_{+}^{D}\big)$ to $\mathcal{H}(R_{g}^{(D)})$. 
\end{corollary} 

\begin{proof}
The linearity of the mapping $\mathsf{K}_{a}[\cdot]$ is obvious. 
Since the functions $\big\{ g^{(\mathbf{p})}(\cdot) \in \mathcal{H}(R_{g}^{(D)}) \big\}_{\mathbf{p} \in \mathbb{Z}_{+}^{D}}$ form an ONB for $\mathcal{H}(R_{g}^{(D)})$, 
we have by Theorem \ref{thm_isometry_hilbert_space_coeffs_space} that any function $f(\cdot) \in \mathcal{H}(R_{g}^{(D)})$ 
can be written as a sum $\sum_{\mathbf{p} \in \mathbb{Z}_{+}^{D}} a[\mathbf{p}] \frac{1}{\sqrt{\mathbf{p}!}} \mathbf{x}^{\mathbf{p}}$, i.e., 
as the image $\mathsf{K}_{a}\big[a[\mathbf{p}]\big]$ of some coefficient sequence $a[\mathbf{p}] \in \ell^{2}\big(\mathbb{Z}_{+}^{D}\big)$. 
We also have that the mapping preserves inner products, i.e., given two functions $f(\cdot) =  \sum_{\mathbf{p} \in \mathbb{Z}_{+}^{D}} a[\mathbf{p}] g^{(\mathbf{p})}(\mathbf{x})$, 
$f'(\cdot) =  \sum_{\mathbf{p} \in \mathbb{Z}_{+}^{D}} a'[\mathbf{p}] g^{(\mathbf{p})}(\mathbf{x})$ we have 
\begin{equation}
\big\langle f(\cdot), f'(\cdot) \big \rangle_{ \mathcal{H}(R_{g}^{(D)})} \stackrel{\eqref{equ_isometry_hilbert_space_coeffs_space_inner_prod_preserv}}{=} \big\langle a[\mathbf{p}], a'[\mathbf{p}] \big\rangle_{ \ell^{2}\big(\mathbb{Z}_{+}^{D}\big)} =\sum_{\mathbf{p} \in \mathbb{Z}_{+}^{D}} a[\mathbf{p}] a'[\mathbf{p}].
\end{equation}
This relation implies that two different coefficient sequences $a[\mathbf{p}]$ and $a'[\mathbf{p}]$, i.e., 
$\| a[\mathbf{p}]- a'[\mathbf{p}] \|_{\ell^{2}\big(\mathbb{Z}_{+}^{D}\big)} \neq 0$ 
cannot yield the same image under the mapping $\mathsf{K}_{a}[\cdot]$. 
To summarize, we have that the mapping $\mathsf{K}_{a}[\cdot]$ is linear, inner-product preserving, and bijective, i.e., a congruence (see Definition \ref{def_isometry_congruence}). 
\end{proof} 

A specific class of functions that belong to $\mathcal{H}(R_{g}^{(D)})$ is presented in   
\begin{theorem} 
\label{thm_exp_times_finite_order_polynom_Rg}
Consider a coefficient sequence $a[\mathbf{p}]: \mathbb{Z}_{+}^{D} \rightarrow \mathbb{R}$ which satisfies the inequality
\begin{equation} 
\label{thm_specific_class_R_g_bound_coeffs}
| a[\mathbf{p}] | \leq C^{|\mathbf{p}|},
\end{equation} 
where $C \in \mathbb{R}_{+}$ is an arbitrary constant. 
The RKHS $\mathcal{H}(R_{g}^{(D)})$ contains any function $f(\cdot): \mathbb{R}^{D} \rightarrow \mathbb{R}$ given by 
\begin{equation}
f(\mathbf{x}) = \exp \big( \mathbf{x}^{T} \mathbf{x}_{1} \big) \sum_{\mathbf{p} \in \mathbb{Z}_{+}^{D}} \frac{a[\mathbf{p}]}{\mathbf{p}!} \mathbf{x}^{\mathbf{p}},
\end{equation} 
where $\mathbf{x}_{1} \in \mathbb{R}^{D}$ is arbitrary. 
\end{theorem}
\begin{proof}
Since $\exp\big( \mathbf{x}^{T} \mathbf{x}_{1}Ê\big) =\sum_{\mathbf{p} \in \mathbb{Z}_{+}^{D}} \frac{\mathbf{x}_{1}^{\mathbf{p}}}{\mathbf{p}!} \mathbf{x}^{\mathbf{p}}$ and $\sum_{\mathbf{p} \in \mathbb{Z}_{+}^{D}} \frac{a[\mathbf{p}]}{\mathbf{p}!} \mathbf{x}^{\mathbf{p}}$ are power series that 
converge everywhere in $\mathbb{R}^{D}$, we can by \cite{KranzPrimerAnalytic}  write their product $f(\cdot)$ also as a power series in the form
\begin{equation}
\label{equ_proof_prod_analy_exp_in_R_g_4}
f(\mathbf{x}) = \sum_{\mathbf{p} \in \mathbb{Z}_{+}^{D}} c[\mathbf{p}] \mathbf{x}^{\mathbf{p}},
\end{equation} 
where the coefficients $c[\mathbf{p}]$ are obtained as \cite{KranzPrimerAnalytic} 
\begin{equation}
\label{equ_proof_prod_analy_exp_in_R_g_133}
c[\mathbf{p}] = \sum_{ \mathbf{n} \leq \mathbf{p}}   \frac{1}{\mathbf{n}!  (\mathbf{p}- \mathbf{n})!} a[\mathbf{p}- \mathbf{n}] \mathbf{x}_{1}^{\mathbf{n}}. 
\end{equation} 
Now we show that the coefficients $d[\mathbf{p}]\triangleq c[\mathbf{p}] \sqrt{\mathbf{p}!}$ are square summable, i.e., $d[\mathbf{p}] \in \ell^{2}(\mathbb{Z}_{+}^{D})$ which 
then implies via \eqref{equ_proof_prod_analy_exp_in_R_g_4} and Theorem \ref{thm_entire_character_H_R_g} that $f(\cdot) \in \mathcal{H}(R_{g}^{(D)})$.

Indeed, using \eqref{thm_specific_class_R_g_bound_coeffs}, \eqref{equ_proof_prod_analy_exp_in_R_g_133} and the fact that for any $\mathbf{x} \in \mathbb{R}^{D}$ and 
multi-index $\mathbf{p} \in \mathbb{Z}_{+}^{D}$ it holds that 
\begin{equation}
\label{equ_bound_power_x_p_infty_norm}
\mathbf{x}^{\mathbf{p}} = \prod_{l \in [D]} x_{l}^{p_{l}} \leq \prod_{l \in [D]}  \| \mathbf{x}\|_{\infty}^{p_{l}}  = \| \mathbf{x}\|_{\infty}^{|\mathbf{p}|}, 
\end{equation}
we obtain the bound 
\begin{align} 
 \big|d[\mathbf{p}] \big| &= \bigg| \sum_{ \mathbf{n} \leq \mathbf{p}}   \frac{\sqrt{\mathbf{p}!}}{\mathbf{n}!  (\mathbf{p}- \mathbf{n})!}  a[\mathbf{p}-\mathbf{n}]    \mathbf{x}_{1}^{\mathbf{n}}\bigg| \leq  \sum_{ \mathbf{n} \leq \mathbf{p}}   \frac{\sqrt{\mathbf{p}!}}{\mathbf{n}!  (\mathbf{p}- \mathbf{n})!} \bigg| a[\mathbf{p}- \mathbf{n}]    \mathbf{x}_{1}^{\mathbf{n}}\bigg|  \nonumber \\[4mm]
& = \frac{1}{\sqrt{\mathbf{p}!}}  \sum_{ \mathbf{n} \leq \mathbf{p}}   \frac{\mathbf{p}!}{\mathbf{n}!  (\mathbf{p}- \mathbf{n})!} \bigg| a[\mathbf{p}- \mathbf{n}]    \mathbf{x}_{1}^{\mathbf{n}}\bigg|  \stackrel{\eqref{equ_bound_power_x_p_infty_norm}}{\leq} 
 \frac{1}{\sqrt{\mathbf{p}!}} \sum_{ \mathbf{n} \leq \mathbf{p}}   \frac{\mathbf{p}!}{\mathbf{n}!  (\mathbf{p}- \mathbf{n})!} \big| a[\mathbf{p}- \mathbf{n}] \big|  \| \mathbf{x}_{1}\|_{\infty}^{|\mathbf{n}|}  \nonumber \\[4mm]
 &Ê\leq  \frac{1}{\sqrt{\mathbf{p}!}} \sum_{ \mathbf{n} \leq \mathbf{p}}   \frac{\mathbf{p}!}{\mathbf{n}!  (\mathbf{p}- \mathbf{n})!}  C^{|\mathbf{p}- \mathbf{n}|}  \| \mathbf{x}_{1}\|_{\infty}^{|\mathbf{n}|} =
  \frac{1}{\sqrt{\mathbf{p}!}} \sum_{ \mathbf{n} \leq \mathbf{p}}   \prod_{l \in [D]} \frac{p_{l}!}{n_{l}!  (p_{l} - n_{l})!}  C^{p_{l}- n_{l}}  \| \mathbf{x}_{1}\|_{\infty}^{n_{l}} \nonumber \\[4mm]
& = \frac{1}{\sqrt{\mathbf{p}!}} \sum_{ \mathbf{n} \leq \mathbf{p}}   \prod_{l \in [D]} \binom{p_{l}}{n_{l}}  C^{p_{l}- n_{l}}  \| \mathbf{x}_{1}\|_{\infty}^{n_{l}} =  \frac{1}{\sqrt{\mathbf{p}!}}  \prod_{l \in [D]}  (C + \| \mathbf{x}_{1}\|_{\infty})^{p_{l}}  \nonumber \\[4mm]
& = \frac{1}{\sqrt{\mathbf{p}!}} (C + \| \mathbf{x}_{1}\|_{\infty})^{|\mathbf{p}|}, 
\end{align} 
which implies that 
\begin{equation} 
\label{equ_proof_prod_analy_exp_upper_bound_squared_d}
 \big( d[\mathbf{p}] \big)^{2} \leq   \frac{1}{\mathbf{p}!} (C + \| \mathbf{x}_{1}\|_{\infty})^{2|\mathbf{p}|}. 
\end{equation} 
By the ratio test \cite[p.\ 66]{RudinBookPrinciplesMatheAnalysis}) we have that 
\begin{equation} 
\label{equ_proof_prod_analy_exp_in_marginal_series_converges_ratio_test}
A \triangleq \sum_{p_l  \in \mathbb{Z}_{+}}  \frac{1}{p_{l}!} (C + \| \mathbf{x}_{1}\|_{\infty})^{2p_{l}} < \infty, 
\end{equation}
since 
\begin{equation} 
\lim_{p_{l} \rightarrow \infty} \frac{\frac{1}{(p_{l}+1)!} (C + \| \mathbf{x}_{1}\|_{\infty})^{2(p_{l}+1)}}{\frac{1}{p_{l}!} (C + \| \mathbf{x}_{1}\|_{\infty})^{2p_{l}}} = \lim_{p_{l} \rightarrow \infty}  \frac{1}{p_{l}+1} (C + \| \mathbf{x}_{1}\|_{\infty})^{2} = 0.
\end{equation}

Now, consider an arbitrary finite index set $\mathcal{T} \subseteq \mathbb{Z}_{+}^{D}$ and denote the largest entry of any element of $\mathcal{T}$ by 
$k(\mathcal{T}) \triangleq \max_{\mathbf{p}' \in \mathcal{T}} \| \mathbf{p}' \|_{\infty}$.
Given the index set $\mathcal{T}$, we define a larger index set $\mathcal{T}'$, i.e., 
\begin{equation}
\mathcal{T}' \triangleq \big\{ \mathbf{p} \in \mathbb{Z}_{+}^{D} \big| p_{l} \leq k(\mathcal{T}) \mbox{ for every } lÊ\in [D] \big\},  
\end{equation} 
which obviously contains $\mathcal{T}$, i.e., $\mathcal{T} \subseteq \mathcal{T}'$. 
Then, we can bound the sum $\sum_{\mathbf{p} \in \mathcal{T}}  \big( d[\mathbf{p}] \big)^{2}$ as 
\begin{align}
\sum_{\mathbf{p}\in \mathcal{T}}  \big( d[\mathbf{p}] \big)^{2} & \leq \sum_{\mathbf{p}\in \mathcal{T}'} \big( d[\mathbf{p}] \big)^{2} \stackrel{\eqref{equ_proof_prod_analy_exp_upper_bound_squared_d}}{\leq}   \sum_{\mathbf{p}\in \mathcal{T}'}   \frac{1}{\mathbf{p}!} (C + \| \mathbf{x}_{1}\|_{\infty})^{2|\mathbf{p}|}
=  \prod_{l \in [D]}  \sum_{p_{l} \in [ k(\mathcal{T})]}   \frac{1}{p_{l}!} (C + \| \mathbf{x}_{1}\|_{\infty})^{2p_{l}} \nonumber \\[4mm]
& \leq \prod_{l \in [D]}  \sum_{p_{l} \in \mathbb{Z}_{+}}   \frac{1}{p_{l}!} (C + \| \mathbf{x}_{1}\|_{\infty})^{2p_{l}} \stackrel{\eqref{equ_proof_prod_analy_exp_in_marginal_series_converges_ratio_test}}{=} A^{D}, 
\end{align}
i.e., we can upper bound any finite sum $\sum_{\mathbf{p} \in \mathcal{T}}  \big( d[\mathbf{p}] \big)^{2}$ by the finite quantity $A^{D}$ (which does not depend on the index set $\mathcal{T}$). Thus, we have 
verified that  $d[\mathbf{p}] \in \ell^{2}(\mathbb{Z}_{+}^{D})$.
\end{proof} 

Note that Theorem \ref{thm_exp_times_finite_order_polynom_Rg} trivially implies that $\mathcal{H}(R_{g}^{(D)})$ contains any function $f(\cdot): \mathbb{R}^{D} \rightarrow \mathbb{R}$ that can be written as
\begin{equation}
f(\mathbf{x}) = \exp\big( \mathbf{x}^{T} \mathbf{x}_{1}\big) \sum_{l \in [L]} a_{l} \mathbf{x}^{\mathbf{p}_{l}},
\end{equation}
with arbitrary $\mathbf{x}_{1} \in \mathbb{R}^{D}$, $L \in \mathbb{N}$, $a_{l} \in \mathbb{R}$, and $\mathbf{p}_{l} \in \mathbb{Z}_{+}^{D}$. 

We will also make use of 
\begin{theorem}
\label{thm_part_der_any_origin_complete_R_g}
For an arbitrary parameter vector $\mathbf{x}_{c} \in \mathbb{R}^{D}$ and multi-index $\mathbf{p} \in \mathbb{Z}_{+}^{D}$, consider the function $g^{(\mathbf{p})}_{\mathbf{x}_{c}}(\cdot): \mathbb{R}^{D} \rightarrow \mathbb{R}$ given by
\begin{equation}
g^{(\mathbf{p})}_{\mathbf{x}_{c}}(\mathbf{x}) \triangleq  \frac{\partial^{\mathbf{p}}  R_{g}^{(D)}(\mathbf{x}, \mathbf{x}_{2})}{\partial \mathbf{x}_{2}^{\mathbf{p}}} \bigg|_{\mathbf{x}_{2} = \mathbf{x}_{c}}.
\end{equation}
Then $g^{(\mathbf{p})}_{\mathbf{x}_{c}}(\cdot)$ belongs to the RKHS $\mathcal{H}(R_{g}^{(D)})$, i.e., 
\begin{equation} 
\label{equ_part_der_arb_point_eval_in_Rg}
\mathbf{p} \in \mathbb{Z}_{+}^{D},\mathbf{x}_{c} \in \mathbb{R}^{D} \quad  \Rightarrow \quad g^{(\mathbf{p})}_{\mathbf{x}_{c}}(\cdot) \in \mathcal{H}(R_{g}^{(D)}). 
\end{equation}
The inner product of an arbitrary function $f(\cdot) \in \mathcal{H}(R_{g}^{(D)})$ with $g^{(\mathbf{p})}_{\mathbf{x}_{c}}(\mathbf{x})$ is given by  
\begin{equation}
\label{equ_part_der_inner_prod_R_g}
\big\langle f(\cdot) , g^{(\mathbf{p})}_{\mathbf{x}_{c}}(\cdot) \big\rangle_{\mathcal{H}(R_{g}^{(D)})} = \frac{\partial^{\mathbf{p}}  f(\mathbf{x})}{\partial \mathbf{x}^{\mathbf{p}}} \bigg|_{\mathbf{x} = \mathbf{x}_{c}}. 
\end{equation} 
Moreover, the set $\big\{ g^{(\mathbf{p})}_{\mathbf{x}_{c}}(\cdot) \big\}_{\mathbf{p} \in \mathbb{Z}_{+}^{D}}$ is complete for $\mathcal{H}(R_{g}^{(D)})$, which implies 
in particular that if two functions $f(\cdot), f'(\cdot) \in \mathcal{H}(R_{g}^{(D)})$ take on the same values on $\mathcal{B}(\mathbf{x}_{c},r) \subseteq \mathbb{R}^{D}$ with an arbitrary radius $r>0$, then 
they are identical, i.e., $f(\cdot) = f'(\cdot)$. 
\end{theorem}
\begin{proof} 
The validity of \eqref{equ_part_der_arb_point_eval_in_Rg} and \eqref{equ_part_der_inner_prod_R_g} follows from Theorem \ref{thm_der_repr_prop} by the fact that $R_{g}^{(D)}(\cdot,\cdot)$ is a differentiable kernel. 
To prove the completeness (see Definition \ref{def_complete_set}) of the set $\big\{ g^{(\mathbf{p})}_{\mathbf{x}_{c}}(\cdot) \big\}_{\mathbf{p} \in \mathbb{Z}_{+}^{D}}$, we consider an arbitrary function 
$f(\cdot) \in \mathcal{H}(R_{g}^{(D)})$ and note that by \eqref{equ_series_repr_RKHS_R_g} in Theorem \ref{thm_entire_character_H_R_g}, it can be represented as a 
power series which converges everywhere in $\mathbb{R}^{D}$. 
According to \cite[Proposition 1.2.3]{KranzPrimerAnalytic}, we can then represent the function as
\begin{equation} 
\label{equ_proof_part_der_span_R_g_LGM_taylor}
f(\mathbf{x}) = \sum_{\mathbf{p}: |\mathbf{p}| \leq L}  \frac{\partial^{\mathbf{p}}  f(\mathbf{x})}{\partial \mathbf{x}^{\mathbf{p}}} \bigg|_{\mathbf{x} = \mathbf{x}_{c}}  \frac{1}{\mathbf{p}!} ( \mathbf{x}- \mathbf{x}_c ) ^{\mathbf{p}}  + R_{L}(\mathbf{x}),
\end{equation} 
where the remainder term $R_{L}(\mathbf{x})$ is such that $\lim_{L \rightarrow \infty} |R_{L}(\mathbf{x})| = 0$ for any $\mathbf{x} \in \mathbb{R}^{D}$. 
By using \eqref{equ_part_der_inner_prod_R_g}, we can rewrite \eqref{equ_proof_part_der_span_R_g_LGM_taylor} as 
\begin{equation} 
\label{equ_proof_part_der_span_R_g_LGM_taylor_inner_prod}
f(\mathbf{x})= \sum_{\mathbf{p}: |\mathbf{p}| \leq L}  \frac{\big\langle f (\cdot), g^{(\mathbf{p})}_{\mathbf{x}_{c}}(\cdot) \big\rangle_{\mathcal{H}(R_{g}^{(D)})} }{\mathbf{p}!} (\mathbf{x}-\mathbf{x}_{c})^{\mathbf{p}}    + R_{L}(\mathbf{x}).
\end{equation} 
From this, we conclude that if $\big\langle f (\cdot), g^{(\mathbf{p})}_{\mathbf{x}_{c}}(\cdot) \big\rangle_{\mathcal{H}(R_{g}^{(D)})}=0$ for all $\mathbf{p} \in \mathbb{Z}_{+}^{D}$, the function $f(\cdot)$ necessarily has to be zero for all arguments, i.e., $f(\cdot) \equiv 0$. 
Thus, according to Definition \ref{def_complete_set}, $\big\{ g^{(\mathbf{p})}_{\mathbf{x}_{c}}(\cdot) \big\}_{\mathbf{p} \in \mathbb{Z}_{+}^{D}}$ is complete for $\mathcal{H}(R_{g}^{(D)})$. 

Finally, consider two functions $f(\cdot), f'(\cdot) \in \mathcal{H}(R_{g}^{(D)})$ that take on the same values on $\mathcal{B}(\mathbf{x}_{c},r) \subseteq \mathbb{R}^{D}$ with an arbitrary radius $r>0$, i.e., 
the difference $f(\cdot) - f'(\cdot)$ is zero for all $\mathbf{x} \in \mathcal{B}(\mathbf{x}_{c},r)$ and in turn the partial derivative $\frac{\partial^{\mathbf{p}}  \big[ f(\mathbf{x})-f'(\mathbf{x})\big]}{\partial \mathbf{x}^{\mathbf{p}}} \bigg|_{\mathbf{x} = \mathbf{x}_{c}}$ is 
zero for every multi-index $\mathbf{p}$. However, by \eqref{equ_part_der_inner_prod_R_g}, this means that the inner product $\big\langle f(\cdot)-f'(\cdot) , g^{(\mathbf{p})}_{\mathbf{x}_{c}}(\cdot) \big\rangle_{\mathcal{H}(R_{g}^{(D)})}$ is equal to zero 
for every $\mathbf{p} \in \mathbb{Z}_{+}^{D}$. Since the set  $\big\{ g^{(\mathbf{p})}_{\mathbf{x}_{c}}(\cdot) \big\}_{\mathbf{p} \in \mathbb{Z}_{+}^{D}}$ is complete for $\mathcal{H}(R_{g}^{(D)})$, this implies that $f(\cdot)-f'(\cdot) \equiv 0$, i.e., $f(\cdot) = f'(\cdot)$. 

\end{proof}

\section{Minimum Variance Estimation for the SLM}
\label{sec_RKHS_basic_facts_SLM}
In the next section, we will derive lower bounds on the minimum achievable variance $L_{\mathcal{M}_{\text{SLM}}(\mathbf{x}_{0})}$ for the class of minimum variance problems 
$\{ \mathcal{M}_{\text{SLM}}(\mathbf{x}_{0}) \}_{\mathbf{x}_{0} \in \mathcal{X}_{S}}$. Here, 
\begin{equation}
\label{equ_def_SLM_func_x_0}
\mathcal{M}_{\text{SLM}}(\mathbf{x}_{0}) \triangleq \left( \mathcal{E}_{\text{SLM}}, c(\cdot), \mathbf{x}_{0} \right),
\end{equation}
for a fixed choice of the SLM parameters $\sigma$, $S$, $M$, $N$, $\mathbf{H}$ and prescribed bias function $c(\cdot): \mathcal{X}_{S} \rightarrow \mathbb{R}$. 
For the special case $S=N$, i.e., where the SLM coincides with the LGM, 
\begin{equation}
\mathcal{M}_{\text{LGM}}(\mathbf{x}_{0}) \triangleq \left( \mathcal{E}_{\text{LGM}}, \tilde{c}(\cdot), \mathbf{x}_{0} \right).
\end{equation}
In this section, we present some fundamental properties of $L_{\mathcal{M}_{\text{SLM}}(\mathbf{x}_{0})}$. 

The following statement gives a detailed characterization of the class of valid bias functions for $\mathcal{M}_{\text{SLM}}(\mathbf{x}_{0})$ and moreover 
states explicit expressions of the minimum achievable variance and corresponding LMV estimator for a valid bias function. 

\begin{theorem}
\label{thm_condition_valid_bias_SLM}
Consider the SLM $\mathcal{E}_{\emph{SLM}}= \big( \mathcal{X}_{S}, f_{\mathbf{H}}(\mathbf{y}; \mathbf{x}), g(\mathbf{x}) = x_{k} \big)$ with specific values of $k$, $S$, $M$, $N$, $\mathbf{H} \in \mathbb{R}^{M \times N}$, and 
denote the thin SVD of $\mathbf{H}$ by $\mathbf{H} = \mathbf{U} \mathbf{\Sigma} \mathbf{V}^{T}$.
\begin{enumerate}
\item
For a specific choice for $\mathbf{x}_{0} \in \mathcal{X}_{S}$, a prescribed bias function $c(\cdot): \mathcal{X}_{S} \rightarrow \mathbb{R}$ 
is valid for $\mathcal{M}_{\emph{SLM}}(\mathbf{x}_{0}) \triangleq \left(\mathcal{E}_{\emph{SLM}},c(\cdot), \mathbf{x}_{0}\right)$, i.e., $L_{\mathcal{M}_{\emph{SLM}}(\mathbf{x}_{0})} < \infty$, 
if and only if it can be written as 
\begin{equation} 
\label{equ_condition_gamma_valid_SLM_fourier_series_R_g}
c(\mathbf{x}) = \exp \left(  \frac{1}{2 \sigma^{2}} \| \mathbf{H} \mathbf{x}_{0} \|^{2}_{2} - \frac{1}{\sigma^{2}} \mathbf{x}^{T}\mathbf{H}^{T} \mathbf{H} \mathbf{x}_{0} \right) \sum_{\mathbf{p} \in \mathbb{Z}_{+}^{D}} \frac{1}{\sqrt{\mathbf{p}!}} a[\mathbf{p}] \bigg( \frac{1}{\sigma} \widetilde{\mathbf{H}}^{\dagger} \mathbf{x}\bigg)^{\mathbf{p}} - x_{k} \mbox{, $\mathbf{x} \in \mathcal{X}_{S}$},
\end{equation} 
with a suitable coefficient sequence $a[\mathbf{p}] \in \ell^{2}(\mathbb{Z}_{+}^{D})$,  $D = \rank(\mathbf{H})$, and $\widetilde{\mathbf{H}} \triangleq \mathbf{V} \mathbf{ \Sigma}^{-1}$. 
If $c(\cdot)$ is valid for $\mathcal{M}_{\emph{SLM}}(\mathbf{x}_{0})$, then it is necessarily continuous. 

\item The minimum achievable variance for $\mathcal{M}_{\emph{SLM}}(\mathbf{x}_{0})$ is obtained in terms of the minimum squared 
$\ell^{2}(\mathbb{Z}_{+}^{D})$ norm $\| a[\cdot] \|_{\ell^{2}(\mathbb{Z}_{+}^{D})} \triangleq \sqrt{\sum_{\mathbf{p} \in \mathbb{Z}_{+}^{D}} ( a[\mathbf{p}] )^{2}}$ 
among all coefficient sequences $a[\mathbf{p}] \in \ell^{2}(\mathbb{Z}_{+}^{D})$ that are consistent with \eqref{equ_condition_gamma_valid_SLM_fourier_series_R_g}, i.e., 
\begin{equation}
\label{equ_def_min_var_SLM_opt_coeff_sequences}
L_{\mathcal{M}_{\emph{SLM}}(\mathbf{x}_{0})} = 
\min_{ \substack{ a[\mathbf{p}] \in \ell^{2}(\mathbb{Z}_{+}^{D}) \\ a[\mathbf{p}] \emph{\tiny{ consistent with }} \eqref{equ_condition_gamma_valid_SLM_fourier_series_R_g}}} 
\| a[\cdot] \|^{2}_{\ell^{2}(\mathbb{Z}_{+}^{D})} - \big[c(\mathbf{x}_{0}) + x_{0,k}\big]^{2}.
\end{equation} 

\item Furthermore, given any valid bias function $c(\cdot): \mathcal{X}_{S} \rightarrow \mathbb{R}$, we have that for every coefficient sequence $a[\mathbf{p}] \in \ell^{2}(\mathbb{Z}_{+}^{D})$ 
which is consistent with \eqref{equ_condition_gamma_valid_SLM_fourier_series_R_g}, the estimator $\hat{g}(\cdot): \mathbb{R}^{M} \rightarrow \mathbb{R}$ given by 
\begin{equation} 
\label{equ_est_arb_coeffs_SLM}
\hat{g}(\cdot) \triangleq \exp \bigg( - \frac{1}{2\sigma^{2}} \| \mathbf{H} \mathbf{x}_{0} \|^{2}_{2} \bigg) \sum_{\mathbf{p} \in \mathbb{Z}^{D}_{+}} \frac{a[\mathbf{p}]}{\sqrt{\mathbf{p}!}}
 \frac{ \partial^{\mathbf{p}} \big[ \rho_{\mathcal{M}_{\emph{LGM}}(\mathbf{x}_{0})}(\cdot,\sigma \widetilde{\mathbf{H}} \mathbf{x}) \exp \big( \frac{1}{\sigma} \mathbf{x}_{0}^{T} \mathbf{H}^{T} \mathbf{H}
  \widetilde{\mathbf{H}} \mathbf{x} \big) \big]}{ \partial \mathbf{x}^{\mathbf{p}}} \Bigg|_{\mathbf{x} = \mathbf{0}} 
\end{equation}
is an allowed estimator for $\mathcal{M}_{\emph{SLM}}(\mathbf{x}_{0})$, i.e., $\hat{g}(\cdot) \in \mathcal{F}(\mathcal{M}_{\emph{SLM}}(\mathbf{x}_{0}))$. 

\item Finally, the estimator obtained from \eqref{equ_est_arb_coeffs_SLM} by the specific coefficient sequence $a_{0}[\mathbf{p}] \in  \ell^{2}(\mathbb{Z}_{+}^{D})$ that has 
the minimum $\ell^{2}(\mathbb{Z}_{+}^{D})$ norm $\| a[\cdot] \|_{\ell^{2}(\mathbb{Z}_{+}^{D})}$ among all coefficient 
sequences that are consistent with \eqref{equ_condition_gamma_valid_SLM_fourier_series_R_g}, is the LMV estimator for $\mathcal{M}_{\emph{SLM}}(\mathbf{x}_{0})$.
\end{enumerate} 
\end{theorem} 

\begin{proof}
\begin{enumerate}
\item 
By Theorem \ref{thm_main_facts_RKHS_MVE}, we have that the prescribed bias function $c(\cdot): \mathcal{X}_{S} \rightarrow \mathbb{R}$ is valid for $\mathcal{M}_{\text{SLM}}(\mathbf{x}_{0})$ 
if and only if $c(\mathbf{x}) + x_{k} \in \mathcal{H}(\mathcal{M}_{\text{SLM}}(\mathbf{x}_{0}))$ which by Theorem 
\ref{thm_relation_RKHS_LGM_SLM} is the case if and only if there exists a function $\gamma(\cdot) \in \mathcal{H}(\mathcal{M}_{\text{LGM}}(\mathbf{x}_{0}))$ such that $\gamma(\mathbf{x}) = c(\mathbf{x}) + x_{k}$ for every $\mathbf{x} \in \mathcal{X}_{S}$. 
Furthermore, by Theorem \ref{thm_isometry_LGM}, $\gamma(\cdot) \in \mathcal{H}(\mathcal{M}_{\text{LGM}}(\mathbf{x}_{0}))$ if and only if 
the function $\gamma(\cdot): \mathbb{R}^{N} \rightarrow \mathbb{R}$ is the image under the congruence $\mathsf{K}_{g}[f(\cdot)]$(see \eqref{equ_def_isometry_LGM}) 
of a function $f(\cdot) \in \mathcal{H}(R_{g}^{(D)})$ (with $D=\rank(\mathbf{H})$), i.e., 
\begin{equation}
\label{equ_image_cond_valid_bias_SLM}
\gamma(\mathbf{x}) =   f\bigg(\frac{1}{\sigma} \widetilde{\mathbf{H}}^{\dagger} \mathbf{x}' \bigg) \exp\left( \frac{1}{2 \sigma^{2}} \| \mathbf{H} \mathbf{x}_{0} \|^{2}_{2} - \frac{1}{\sigma^{2}} (\mathbf{x}')^{T}  \mathbf{H}^{T}\mathbf{H} \mathbf{x}_{0}\right), 
\end{equation} 
with a suitable function $f(\cdot) \in \mathcal{H}(R_{g}^{(D)})$.
By Theorem \ref{thm_entire_character_H_R_g}, a function $f(\cdot): \mathbb{R}^{D} \rightarrow \mathbb{R}$ belongs to $ \mathcal{H}(R_{g}^{(D)})$ if and only 
if it can be written as \eqref{equ_series_repr_RKHS_R_g} with  a unique coefficient sequence $a[\mathbf{p}] \in \ell^{2}(\mathbb{Z}_{+}^{D})$. 
The condition \eqref{equ_condition_gamma_valid_SLM_fourier_series_R_g} for a bias function to be valid follows then by combining \eqref{equ_image_cond_valid_bias_SLM} with \eqref{equ_series_repr_RKHS_R_g}.

The continuity of any valid bias function $c(\cdot)$ can be verified by Theorem \ref{thm_continous_kernel_implies_cont_func}, since the kernel $R_{\mathcal{M}_{\text{SLM}}}$ (see \eqref{equ_kernel_SLM}) is obviously continuous. 

\item Given a valid bias function $c(\mathbf{x})$, and using the shorthand $\tilde{\gamma}(\cdot): \mathcal{X}_{S} \rightarrow \mathbb{R}: \tilde{\gamma}(\mathbf{x}) = c(\mathbf{x})+x_{k}$, the minimum achievable variance is then obtained by Theorem \ref{thm_main_facts_RKHS_MVE} and \eqref{equ_relation_norm_SLM_norm_LGM} in Theorem \ref{thm_relation_RKHS_LGM_SLM} as 
\begin{align}
\label{equ_expr_min_achiev_var_general_SLM_min}
L_{\mathcal{M}_{\text{SLM}}(\mathbf{x}_{0})}  & \stackrel{\eqref{equ_squared_norm_min_achiev_var}}{=} 
\| \tilde{\gamma}(\cdot) \|^{2}_{\mathcal{H}(\mathcal{M}_{\text{SLM}}(\mathbf{x}_{0}))} - \big[ \tilde{\gamma}(\mathbf{x}_{0}) \big]^{2}    \nonumber \\[4mm]
& \stackrel{\eqref{equ_relation_norm_SLM_norm_LGM}}{=} 
\min_{\substack{g(\cdot) \in \mathcal{H}(\mathcal{M}_{\text{LGM}}(\mathbf{x}_{0})) \\  g(\cdot)\big|_{\mathcal{X}_{S}} = \tilde{\gamma}(\cdot)} } \| g(\cdot) \|^{2}_{\mathcal{H}(\mathcal{M}_{\text{LGM}}(\mathbf{x}_{0}))} -\big[ \tilde{\gamma}(\mathbf{x}_{0}) \big]^{2} 
\end{align}  
Using Theorem \ref{thm_isometry_LGM} and Corollary \ref{cor_R_g_congr_sequence_space}, we can rewrite \eqref{equ_expr_min_achiev_var_general_SLM_min} as 
\begin{align}
\label{equ_SLM_norm_squared_Norm_coeffs_restr_constr_proof}
L_{\mathcal{M}_{\text{SLM}}(\mathbf{x}_{0})}  & =
\min_{\substack{g(\cdot) \in \mathcal{H}(\mathcal{M}_{\text{LGM}}(\mathbf{x}_{0})) \\  g(\cdot)\big|_{\mathcal{X}_{S}} = \tilde{\gamma}(\cdot)} } \| g(\cdot) \|^{2}_{\mathcal{H}(\mathcal{M}_{\text{LGM}}(\mathbf{x}_{0}))} - \big[ \tilde{\gamma}(\mathbf{x}_{0}) \big]^{2} \nonumber \\[4mm]
& = 
\min_{\substack{f(\cdot) \in  \mathcal{H}(R_{g}^{(D)}) \\  \mathsf{K}_{g}[f(\cdot)]\big|_{\mathcal{X}_{S}} = \tilde{\gamma}(\cdot)} } \bigg\|  \mathsf{K}_{g}[f(\cdot)] \bigg \|^{2}_{\mathcal{H}(\mathcal{M}_{\text{LGM}}(\mathbf{x}_{0}))} - \big[ \tilde{\gamma}(\mathbf{x}_{0}) \big]^{2} \nonumber \\[4mm]
& = 
\min_{\substack{a[\mathbf{p}]  \in  \ell^{2}(\mathbb{Z}_{+}^{D}) ) \\  \mathsf{K}_{g}\big[\mathsf{K}_{a}\big[a[\mathbf{p}]\big]\big] \big|_{\mathcal{X}_{S}} 
= \tilde{\gamma}(\cdot)} } \bigg\|   \mathsf{K}_{g}\big[\mathsf{K}_{a}\big[a[\mathbf{p}]\big]\big] \bigg \|^{2}_{\mathcal{H}(\mathcal{M}_{\text{LGM}}(\mathbf{x}_{0}))} - \big[ \tilde{\gamma}(\mathbf{x}_{0}) \big]^{2} \nonumber \\[4mm]
& = 
\min_{\substack{a[\mathbf{p}]  \in  \ell^{2}(\mathbb{Z}_{+}^{D}) ) \\  \mathsf{K}_{g}\big[\mathsf{K}_{a}\big[a[\mathbf{p}]\big]\big] \big|_{\mathcal{X}_{S}} 
= \tilde{\gamma}(\cdot)} } \big\|   a[\mathbf{p}] \big \|^{2}_{\ell^{2}(\mathbb{Z}_{+}^{D})} - \big[ c(\mathbf{x}_{0})+x_{0,k} \big]^{2}.
\end{align} 
Since 
\begin{align}
\label{equ_proof_general_SLM_cond_restriction_congruences}
\mathsf{K}_{g}\big[\mathsf{K}_{a}\big[a[\mathbf{p}]\big]\big] & \stackrel{\eqref{equ_cor_R_g_congr_sequence_space_a_p}}{=} \mathsf{K}_{g}\Bigg[ \sum_{\mathbf{p} \in \mathbb{Z}_{+}^{N}} \frac{1}{\sqrt{\mathbf{p}!}} a[\mathbf{p}] \mathbf{x}^{\mathbf{p}} \Bigg] \nonumber \\[4mm]
& \stackrel{\eqref{equ_def_isometry_LGM}}{=} \exp \left(  \frac{1}{2 \sigma^{2}} \| \mathbf{H} \mathbf{x}_{0} \|^{2}_{2} - \frac{1}{\sigma^{2}} \mathbf{x}^{T}\mathbf{H}^{T} \mathbf{H} \mathbf{x}_{0} \right) \sum_{\mathbf{p} \in \mathbb{Z}_{+}^{D}} \frac{1}{\sqrt{\mathbf{p}!}} a[\mathbf{p}] \bigg( \frac{1}{\sigma} \widetilde{\mathbf{H}}^{\dagger} \mathbf{x}\bigg)^{\mathbf{p}}, 
\end{align}
we have that a coefficient sequence $a[\mathbf{p}]  \in  \ell^{2}(\mathbb{Z}_{+}^{D})$ satisfies 
\begin{equation} 
\mathsf{K}_{g}\big[\mathsf{K}_{a}\big[a[\mathbf{p}]\big]\big] \big|_{\mathcal{X}_{S}} 
= \tilde{\gamma}(\mathbf{x})=c(\mathbf{x}) + x_{k}
\end{equation}
if and only if it is consistent with \eqref{equ_condition_gamma_valid_SLM_fourier_series_R_g}. 
Based on this observation, the relation \eqref{equ_def_min_var_SLM_opt_coeff_sequences} follows then from \eqref{equ_SLM_norm_squared_Norm_coeffs_restr_constr_proof}. 

\item
Now consider 
a coefficient sequence $a[\mathbf{p}] \in \ell^{2}(\mathbb{Z}_{+}^{D})$ which is consistent \eqref{equ_condition_gamma_valid_SLM_fourier_series_R_g}. 
This implies that the function 
\begin{equation}
\label{equ_proof_number_3_bar_gamma}
\bar{\gamma}(\cdot): \mathbb{R}^{N} \rightarrow \mathbb{R}: \bar{\gamma}(\mathbf{x}) \triangleq  \exp \left(  \frac{1}{2 \sigma^{2}} \| \mathbf{H} \mathbf{x}_{0} \|^{2}_{2} - \frac{1}{\sigma^{2}} \mathbf{x}^{T}\mathbf{H}^{T} \mathbf{H} \mathbf{x}_{0} \right) \sum_{\mathbf{p} \in \mathbb{Z}_{+}^{D}} \frac{1}{\sqrt{\mathbf{p}!}} a[\mathbf{p}] \bigg( \frac{1}{\sigma} \widetilde{\mathbf{H}}^{\dagger} \mathbf{x}\bigg)^{\mathbf{p}}
\end{equation} 
satisfies 
\begin{equation}
\bar{\gamma}(\cdot) \big|_{\mathcal{X}_{S}} = c(\mathbf{x}) + x_{k}.  
\end{equation} 
Since by \eqref{equ_proof_general_SLM_cond_restriction_congruences}   
\begin{equation} 
\bar{\gamma}(\cdot) =  \mathsf{K}_{g}\big[\mathsf{K}_{a}\big[a[\mathbf{p}]\big]\big], 
\end{equation}
we have that the function $\bar{\gamma}(\cdot)$ belongs to $\mathcal{H}(\mathcal{M}_{\text{LGM}}(\mathbf{x}_{0}))$, i.e., we can define 
the estimator 
\begin{equation} 
\hat{g}^{a[\mathbf{p}]}(\cdot) \triangleq \mathsf{J}[ \bar{\gamma}(\cdot)] =   \mathsf{J} \big[\mathsf{K}_{g}\big[\mathsf{K}_{a}\big[a[\mathbf{p}]\big]\big] \big]
\end{equation} 
as the image of $\bar{\gamma}(\cdot)$ under the congruence $\mathsf{J}[\cdot]: \mathcal{H}(\mathcal{M}_{\text{LGM}}(\mathbf{x}_{0})) \rightarrow 
\mathcal{L}(\mathcal{H}(\mathcal{M}_{\text{LGM}}(\mathbf{x}_{0})))$ defined in \eqref{equ_def_congruence_H_L_MVP}. 
By Theorem \ref{thm_isometry_RKHS_rhos}  and Theorem \ref{thm_main_facts_RKHS_MVE}, the estimator $\hat{g}^{a[\mathbf{p}]}(\cdot)$ has variance at $\mathbf{x}_{0}$ equal to 
\begin{align}
\label{equ_var_def_g_a_p_congruence}
 v(\hat{g}^{a[\mathbf{p}]}(\cdot); \mathbf{x}_{0}) & =   \|  \bar{\gamma}(\cdot)  \|^{2}_{\mathcal{H}(\mathcal{M}_{\text{LGM}}(\mathbf{x}_{0}))} - \big[ \bar{\gamma}(\mathbf{x}_{0}) \big]^{2}   \nonumber \\[4mm]
   &=   \|   \mathsf{K}_{g}\big[\mathsf{K}_{a}\big[a[\mathbf{p}]\big]\big] \|^{2}_{\mathcal{H}(\mathcal{M}_{\text{LGM}}(\mathbf{x}_{0}))} - \big[ \bar{\gamma}(\mathbf{x}_{0}) \big]^{2}
=  \big\|   a[\mathbf{p}] \big \|^{2}_{\ell^{2}(\mathbb{Z}_{+}^{D})} - \big[ \bar{\gamma}(\mathbf{x}_{0}) \big]^{2}  < \infty
\end{align} 
and mean function equal to $\bar{\gamma}(\mathbf{x})$. In particular, the bias of $\hat{g}^{a[\mathbf{p}]}(\cdot)$ is equal to $c(\mathbf{x})$ for every 
$\mathbf{x} \in \mathcal{X}_{S}$. This implies that the estimator $\hat{g}^{a[\mathbf{p}]}(\cdot)$ is an allowed estimator for $\mathcal{M}_{\text{SLM}}(\mathbf{x}_{0})$, i.e., 
$\hat{g}^{a[\mathbf{p}]}(\cdot) \in \mathcal{F}(\mathcal{M}_{\text{SLM}}(\mathbf{x}_{0}))$

Let us now show that the estimator $\hat{g}^{a[\mathbf{p}]}(\cdot)$ coincides with the estimator \eqref{equ_est_arb_coeffs_SLM}, when using the same coefficient sequence $a[\mathbf{p}]$. 
Based on the identity
\begin{equation}
\label{equ_proof_isometr_R_g_pseudo_inver_tilde_H} 
\widetilde{\mathbf{H}}^{T}\mathbf{H}^{T}\mathbf{H} =\mathbf{\Sigma}^{-1} \mathbf{V}^{T}  \mathbf{V} \mathbf{\Sigma}^{2} \mathbf{V}^{T}  =\mathbf{\Sigma} \mathbf{V}^{T} =\widetilde{\mathbf{H}}^{\dagger},
\end{equation} 
we observe that 
\begin{align} 
&  \exp \left(  \frac{1}{ \sigma^{2}} \| \mathbf{H} \mathbf{x}_{0} \|^{2}_{2} - \frac{1}{\sigma^{2}} \mathbf{x}^{T}\mathbf{H}^{T} \mathbf{H} \mathbf{x}_{0} \right)\bigg( \frac{1}{\sigma} \widetilde{\mathbf{H}}^{\dagger} \mathbf{x}\bigg)^{\mathbf{p}} \nonumber \\[4mm]
&  = \exp \left(  \frac{1}{ \sigma^{2}} \| \mathbf{H} \mathbf{x}_{0} \|^{2}_{2} - \frac{1}{\sigma^{2}} \mathbf{x}^{T}\mathbf{H}^{T} \mathbf{H} \mathbf{x}_{0} \right)\frac{ \partial^{\mathbf{p}} \exp \big( \frac{1}{\sigma} (\mathbf{x}')^{T} \widetilde{\mathbf{H}}^{\dagger} \mathbf{x} \big) }{\partial \mathbf{x}'^{\mathbf{p}}}\Bigg|_{\mathbf{x}'=\mathbf{0}} \nonumber \\[4mm]
& \stackrel{\eqref{equ_proof_isometr_R_g_pseudo_inver_tilde_H}}{=} \exp \left(  \frac{1}{ \sigma^{2}} \| \mathbf{H} \mathbf{x}_{0} \|^{2}_{2} - \frac{1}{\sigma^{2}} \mathbf{x}^{T}\mathbf{H}^{T} \mathbf{H} \mathbf{x}_{0} \right)\frac{ \partial^{\mathbf{p}} \exp \big( \frac{1}{\sigma} \mathbf{x}^{T} \mathbf{H}^{T} \mathbf{H} 
\widetilde{\mathbf{H}} \mathbf{x}' \big) }{\partial \mathbf{x}'^{\mathbf{p}}}\Bigg|_{\mathbf{x}'=\mathbf{0}} \nonumber \\[4mm]
& = \frac{ \partial^{\mathbf{p}} \exp \big(  \frac{1}{ \sigma^{2}} \| \mathbf{H} \mathbf{x}_{0} \|^{2}_{2} - \frac{1}{\sigma^{2}} \mathbf{x}^{T}\mathbf{H}^{T} \mathbf{H} \mathbf{x}_{0} + \frac{1}{\sigma} \mathbf{x}^{T} \mathbf{H}^{T} \mathbf{H} 
\widetilde{\mathbf{H}} \mathbf{x}' \big) }{\partial \mathbf{x}'^{\mathbf{p}}}\Bigg|_{\mathbf{x}'=\mathbf{0}} \nonumber \\[4mm]
& = \frac{ \partial^{\mathbf{p}} \exp \big(  \frac{1}{ \sigma^{2}} (\mathbf{x} - \mathbf{x}_{0})^{T} \mathbf{H}^{T} \mathbf{H} \big(\sigma \widetilde{\mathbf{H}} \mathbf{x}' -\mathbf{x}_{0} \big) + \frac{1}{\sigma} \mathbf{x}_{0}^{T} \mathbf{H}^{T} \mathbf{H} \widetilde{\mathbf{H}} \mathbf{x}' \big) }{\partial \mathbf{x}'^{\mathbf{p}}}\Bigg|_{\mathbf{x}'=\mathbf{0}} \nonumber \\[4mm]
& \stackrel{\eqref{equ_kernel_LGM}}{=} \frac{ \partial^{\mathbf{p}} \big[ R_{\mathcal{M}_{\text{LGM}}(\mathbf{x}_{0})}\big(\mathbf{x},\sigma \widetilde{\mathbf{H}} \mathbf{x}'\big) \exp \big( \frac{1}{\sigma} \mathbf{x}_{0}^{T} \mathbf{H}^{T} \mathbf{H} \widetilde{\mathbf{H}} \mathbf{x}' \big) \big] }{\partial \mathbf{x}'^{\mathbf{p}}}\Bigg|_{\mathbf{x}'=\mathbf{0}}.
\end{align}  
This allows us to rewrite \eqref{equ_proof_number_3_bar_gamma} as 
\begin{equation} 
\bar{\gamma}(\cdot) = 
 \exp \big( - \frac{1}{2\sigma^{2}} \| \mathbf{H} \mathbf{x}_{0} \|^{2}_{2} \big)  \sum_{\mathbf{p} \in \mathbb{Z}_{+}^{D}} \frac{a[\mathbf{p}]}{\sqrt{\mathbf{p}!}}   
\frac{ \partial^{\mathbf{p}} \bigg( R_{\text{LGM}}(\mathbf{x}, \sigma \widetilde{\mathbf{H}} \mathbf{x}') \exp \big( \frac{1}{\sigma} \mathbf{x}_{0}^{T} \mathbf{H}^{T} \mathbf{H} 
\widetilde{\mathbf{H}} \mathbf{x}' \big) \bigg)}{\partial \mathbf{x}'^{\mathbf{p}}}\Bigg|_{\mathbf{x}'=\mathbf{0}}, 
\end{equation} 
from which it follows by an application of Theorem \ref{thm_isometry_RKHS_rhos_derivative_kernel} that 
\begin{align} 
\hat{g}^{a[\mathbf{p}]}(\cdot) & = \mathsf{J}[ \bar{\gamma}(\cdot)]  Ê\nonumber \\[4mm]
& \hspace*{-10mm}=  \mathsf{J}\Bigg[  \exp \bigg( - \frac{1}{2\sigma^{2}} \| \mathbf{H} \mathbf{x}_{0} \|^{2}_{2} \bigg)  \sum_{\mathbf{p} \in \mathbb{Z}_{+}^{D}} \frac{a[\mathbf{p}]}{\sqrt{\mathbf{p}!}}   
\frac{ \partial^{\mathbf{p}} \bigg( R_{\text{LGM}}(\mathbf{x}, \sigma \widetilde{\mathbf{H}} \mathbf{x}') \exp \big( \frac{1}{\sigma} \mathbf{x}_{0}^{T} \mathbf{H}^{T} \mathbf{H} 
\widetilde{\mathbf{H}} \mathbf{x}' \big) \bigg)}{\partial \mathbf{x}'^{\mathbf{p}}}\Bigg|_{\mathbf{x}'=\mathbf{0}} \Bigg]  \nonumber \\[4mm]
& \hspace*{-10mm}= \exp \bigg( - \frac{1}{2\sigma^{2}} \| \mathbf{H} \mathbf{x}_{0} \|^{2}_{2} \bigg)  \sum_{\mathbf{p} \in \mathbb{Z}_{+}^{D}} \frac{a[\mathbf{p}]}{\sqrt{\mathbf{p}!}}   
\mathsf{J}\Bigg[  \frac{ \partial^{\mathbf{p}} \bigg( R_{\text{LGM}}(\mathbf{x}, \sigma \widetilde{\mathbf{H}} \mathbf{x}') \exp \big( \frac{1}{\sigma} \mathbf{x}_{0}^{T} \mathbf{H}^{T} \mathbf{H} 
\widetilde{\mathbf{H}} \mathbf{x}' \big) \bigg)}{\partial \mathbf{x}'^{\mathbf{p}}}\Bigg|_{\mathbf{x}'=\mathbf{0}} \Bigg]  \nonumber \\[4mm]
& \hspace*{-10mm}= \bigg( - \frac{1}{2\sigma^{2}} \| \mathbf{H} \mathbf{x}_{0} \|^{2}_{2} \bigg) \sum_{\mathbf{p} \in \mathbb{Z}^{D}_{+}} \frac{a[\mathbf{p}]}{\sqrt{\mathbf{p}!}}
 \frac{ \partial^{\mathbf{p}} \big[ \rho_{\mathcal{M}_{\text{LGM}}(\mathbf{x}_{0})}(\cdot,\sigma \widetilde{\mathbf{H}} \mathbf{x}) \exp \big( \frac{1}{\sigma} \mathbf{x}_{0}^{T} \mathbf{H}^{T} \mathbf{H}
  \widetilde{\mathbf{H}} \mathbf{x} \big) \big]}{ \partial \mathbf{x}^{\mathbf{p}}} \Bigg|_{\mathbf{x} = \mathbf{0}}, 
\end{align} 
which coincides with \eqref{equ_est_arb_coeffs_SLM}.

\item  Note that by definition, the coefficient sequence $a_{0}[\mathbf{p}] \in  \ell^{2}(\mathbb{Z}_{+}^{D})$ is a minimizer 
of \eqref{equ_SLM_norm_squared_Norm_coeffs_restr_constr_proof}.
Since the estimator \eqref{equ_est_arb_coeffs_SLM} obtained for the coefficients $a_{0}[\mathbf{p}]$ coincides with the estimator 
$\hat{g}^{a_{0}[\mathbf{p}]}(\cdot)=  \mathsf{J} \big[\mathsf{K}_{g}\big[\mathsf{K}_{a}\big[a_{0}[\mathbf{p}]\big]\big] \big]
$, which is an allowed estimator for  $\mathcal{M}_{\text{SLM}}(\mathbf{x}_{0})$ whose variance 
at $\mathbf{x}_{0}$ equals (cf.\ \eqref{equ_var_def_g_a_p_congruence})
\begin{align} 
\big\|   a_{0}[\mathbf{p}] \big \|^{2}_{\ell^{2}(\mathbb{Z}_{+}^{D})} - \big[ \bar{\gamma}(\mathbf{x}_{0}) \big]^{2}& = 
\big\|   a_{0}[\mathbf{p}] \big \|^{2}_{\ell^{2}(\mathbb{Z}_{+}^{D})} - \big[ c(\mathbf{x}_{0}) + x_{0,k}) \big]^{2}  \nonumber \\[4mm]
& \stackrel{\eqref{equ_SLM_norm_squared_Norm_coeffs_restr_constr_proof}}{=} L_{\mathcal{M}_{\text{SLM}}(\mathbf{x}_{0})}.  
\end{align}
\end{enumerate}
\end{proof}

One specific implication of Theorem \ref{thm_condition_valid_bias_SLM} concerning the problem of sparsity pattern detection \cite{EstThSparsPattRec,KarbasiSuppRecCS,WainInfThBoundsSparsPattRec} 
is stated in 
\begin{theorem}
\label{thm_non_esit_finite_var_unbiased_est_supp}
There exists no estimator $\hat{g}(\mathbf{y})$ of an injective real-valued function $g(\supp(\mathbf{x}))$ of the support of $\mathbf{x} \in \mathcal{X}_{S}$ which uses only the observation \eqref{equ_linear_observation_model} of the SLM $\mathcal{E}_{\emph{SLM}}$, is unbiased for every $\mathbf{x} \in \mathcal{X}_{S}$, and has 
a finite variance at any $\mathbf{x}_{0} \in \mathcal{X}_{S}$.
\end{theorem}
\begin{proof}
Consider the SLM $\mathcal{E}_{\text{SLM}}$, but instead of estimating the parameter function $g(\mathbf{x}) = x_{k}$ (for some fixed index $k \in [N]$), we are 
interested in unbiased estimation of an injective real-valued function $g(\supp(\mathbf{x}))$ of the support $\supp(\mathbf{x})$. 
According to Section \ref{sec_equ_bias_param_function}, this is equivalent to a minimum variance problem associated with the ordinary 
SLM $\mathcal{E}_{\text{SLM}}$, i.e., with the parameter function $g(\mathbf{x}) = x_{k}$, and using the prescribed bias function 
$c'(\cdot): \mathcal{X}_{S} \rightarrow \mathbb{R}: c'(\mathbf{x}) = g(\supp(\mathbf{x})) - x_{k}$, which is obviously discontinuous. 
Indeed, consider the set $\{ \mathbf{x}^{(a)} \triangleq a \mathbf{e}_{1} \}_{a \in \mathbb{R}} \subseteq \mathcal{X}_{S}$ consisting of parameter vectors of the 
form $\mathbf{x}^{(a)} = \big(a,0,\ldots,0\big)^{T}$, i.e., with exactly one nonzero entry whose value is equal to the number $a \in \mathbb{R}$. 
We then have $\lim_{a \rightarrow 0} c'(\mathbf{x}^{(a)}) = g(Ê\{1\})$ but $c'(\mathbf{x}^{(0)}) = g( \emptyset )$, where $g(\emptyset) \neq g(\{1 \})$ since 
$g(\supp(\mathbf{x}))$ is assumed to be injective. 

However, since the prescribed bias function $c'(\cdot)$ is discontinuous, we have by Theorem \ref{thm_condition_valid_bias_SLM} that this bias function is not valid for any minimum variance problem  
$\mathcal{M}_{\text{SLM}}(\mathbf{x}_{0})$ with $\mathbf{x}_{0} \in \mathcal{X}_{S}$, i.e., there exists no estimator with the bias $c'(\cdot)$ and finite variance at any point $\mathbf{x}_{0} \in \mathcal{X}_{S}$. 
Thus, by Theorem \ref{cor_estimable_equiv_valid}, we have that the parameter function $g(\supp(\mathbf{x}))$ is not estimable for the 
SLM $\mathcal{E}_{\text{SLM}}$ at any $\mathbf{x}_{0} \in \mathcal{X}_{S}$, which is just an equivalent formulation of the statement to be proved. 
\end{proof}

It is important to note that Theorem \ref{thm_non_esit_finite_var_unbiased_est_supp} asserts the impossibility of a finite variance estimator of the support 
of the parameter vector $\mathbf{x}$ which is unbiased for all $\mathbf{x} \in \mathcal{X}_{S}$. 
However, even if it might seem at first sight, Theorem \ref{thm_non_esit_finite_var_unbiased_est_supp} does not contradict \cite{EstThSparsPattRec,KarbasiSuppRecCS} 
since the authors of \cite{EstThSparsPattRec,KarbasiSuppRecCS} intentionally reduce the parameter set $\mathcal{X}_{S}$ of the SLM to a smaller parameter set $\mathcal{X}'$ which includes only vectors $\mathbf{x}$ 
with exactly $S$ non-zeros, i.e., $\| \mathbf{x} \|_{0} = S$, and whose smallest nonzero magnitude is not smaller than a given threshold $\theta_{\text{min}}$, i.e., $\min_{k \in \supp(\mathbf{x})} |x_{k}| \geq \theta_{\text{min}}$.

The shape of $L_{\mathcal{M}_{\text{SLM}}(\mathbf{x}_{0})}$, viewed as a function of $\mathbf{x}_{0} \in \mathcal{X}_{S}$, is characterized by 
\begin{theorem} 
\label{thm_SLM_min_achiev_var_lower_semi_cont}
Consider the SLM $\mathcal{E}_{\emph{SLM}}$ with arbitrary values of $S$, $M$, $N$, and $\mathbf{H}$ and a prescribed bias function $c(\cdot): \mathcal{X}_{S} \rightarrow \mathbb{R}$ that is valid for 
any $\mathcal{M}_{\emph{SLM}}(\mathbf{x}_{0})$ (cf. \eqref{equ_def_SLM_func_x_0}) with $\mathbf{x}_{0} \in \mathcal{X}_{S}$. 
Then the minimum achievable variance $L_{\mathcal{M}_{\emph{SLM}}(\mathbf{x}_{0})}$ exists, i.e., is finite for every $\mathbf{x}_{0}$, and viewed as a function of $\mathbf{x}_{0}$ is lower semi-continuous. 
\end{theorem} 
\begin{proof} 
This fact follows from Theorem \ref{thm_lower_semi_cont_varying_kernel} by noting that the kernel $R_{\mathcal{M}_{\text{SLM}}(\mathbf{x}_{0})}(\cdot,\cdot)$ 
(see \eqref{equ_kernel_SLM}) obviously satisfies \eqref{equ_def_cont_varying_kernel}.
\end{proof} 

Finally, we give two sufficient conditions on the prescribed bias $c(\cdot)$ to be valid for $\mathcal{M}_{\text{SLM}}(\mathbf{x}_{0})$ with arbitrary $\mathbf{x}_{0} \in \mathcal{X}_{S}$:  
\begin{theorem} 
\label{thm_suff_cond_valid_erervwhere_SLM_power_series}
Consider the minimum variance problem $\mathcal{M}_{\emph{SLM}}(\mathbf{x}_{0})$ with system matrix $\mathbf{H} \in \mathbb{R}^{M \times N}$ and a prescribed bias function $c(\mathbf{x}) : \mathcal{X}_{S} \rightarrow \mathbb{R}$ which is given as
\begin{equation}
\label{equ_prescr_bias_finite_order_polynom_SLM}
c(\mathbf{x}) = \exp \big( \mathbf{x}_{1}^{T}  \widetilde{\mathbf{H}}^{\dagger} \mathbf{x} \big) \sum_{\mathbf{p} \in \mathbb{Z}_{+}^{D}} \frac{a[\mathbf{p}]}{\mathbf{p}!}  \bigg( \frac{1}{\sigma} \widetilde{\mathbf{H}}^{\dagger} \mathbf{x}\bigg)^{\mathbf{p}} - x_{k} , 
\end{equation} 
with an arbitrary vector $\mathbf{x}_{1} \in \mathbb{R}^{D}$. 
Then, if the coefficients satisfy $| a[\mathbf{p}] | \leq C^{|\mathbf{p}|}$ with an arbitrary constant $C \in \mathbb{R}_{+}$, the bias $c(\cdot)$ is valid for any $\mathcal{M}_{\emph{SLM}}(\mathbf{x}_{0})$ with $\mathbf{x}_{0} \in \mathcal{X}_{S}$. 
\end{theorem}

\begin{proof} 
By Theorem \ref{thm_main_facts_RKHS_MVE}, we have that a prescribed bias function $c(\mathbf{x})$ is valid for $\mathcal{M}_{\text{SLM}}(\mathbf{x}_{0})$, 
if (and only if) the prescribed mean function $\gamma(\cdot): \mathcal{X}_{S} \rightarrow \mathbb{R}:$ 
\begin{equation} 
\label{equ_suff_cond_valid_erervwhere_SLM_power_series_proof_gamma}
\gamma(\mathbf{x}) = c(\mathbf{x}) + x_{k} = \exp \big( \mathbf{x}_{1}^{T}  \widetilde{\mathbf{H}}^{\dagger} \mathbf{x} \big) \sum_{\mathbf{p} \in \mathbb{Z}_{+}^{D}} \frac{a[\mathbf{p}]}{\mathbf{p}!}  \bigg( \frac{1}{\sigma} \widetilde{\mathbf{H}}^{\dagger} \mathbf{x}\bigg)^{\mathbf{p}}
\end{equation}
belongs to the RKHS 
$\mathcal{H}(\mathcal{M}_{\text{SLM}}(\mathbf{x}_{0}))$. 
According to Theorem \ref{thm_exp_times_finite_order_polynom_Rg} the function $f(\mathbf{x}): \mathbb{R}^{D} \rightarrow \mathbb{R}:$ 
\begin{equation}
f(\mathbf{x}) =  \exp\left(- \frac{1}{2 \sigma^{2}} \| \mathbf{H} \mathbf{x}_{0} \|^{2}_{2} \right)  \exp \big( \big(\sigma \mathbf{x}_{1}+ \frac{1}{\sigma}  \widetilde{\mathbf{H}}^{\dagger}\mathbf{x}_{0}\big)^{T} \mathbf{x} \big) \sum_{\mathbf{p} \in \mathbb{Z}_{+}^{D}} \frac{a[\mathbf{p}]}{\mathbf{p}!}  \mathbf{x}^{\mathbf{p}}, 
\end{equation}
with the same $\mathbf{x}_{1} \in \mathbb{R}^{D}$ and coefficients $a[\mathbf{p}]$ as in \eqref{equ_prescr_bias_finite_order_polynom_SLM}, belongs to the RKHS $\mathcal{H}(R_{g}^{(D)})$ since $| a[\mathbf{p}] | \leq C^{|\mathbf{p}|}$.
The image $\mathsf{K}_{g}[f(\cdot)] \in \mathcal{H}(\mathcal{M}_{\text{LGM}})$ under the congruence \eqref{equ_def_isometry_LGM} is then obtained as 
\begin{align}
\label{equ_suff_cond_valid_erervwhere_SLM_power_series_proof_gamma_image_congruence}
\mathsf{K}_{g}[f(\cdot)](\mathbf{x}) & = f\bigg(\frac{1}{\sigma} \widetilde{\mathbf{H}}^{\dagger} \mathbf{x} \bigg) \exp\left( \frac{1}{2 \sigma^{2}} \| \mathbf{H} \mathbf{x}_{0} \|^{2}_{2} - \frac{1}{\sigma^{2}} \mathbf{x}^{T}  \mathbf{H}^{T}\mathbf{H} \mathbf{x}_{0}\right) \nonumber \\[4mm]Ê
& \hspace*{-20mm} = \exp\left( \sigma\big(\mathbf{x}_{1}+ \frac{1}{\sigma^{2}}  \widetilde{\mathbf{H}}^{\dagger}\mathbf{x}_{0}\big)^{T} \frac{1}{\sigma} \widetilde{\mathbf{H}}^{\dagger} \mathbf{x} \right) \sum_{\mathbf{p} \in \mathbb{Z}_{+}^{D}} \frac{a[\mathbf{p}]}{\mathbf{p}!} \bigg(\frac{1}{\sigma} \widetilde{\mathbf{H}}^{\dagger} \mathbf{x}\bigg)^{\mathbf{p}}  \exp\left(- \frac{1}{\sigma^{2}} \mathbf{x}^{T}  \mathbf{H}^{T}\mathbf{H} \mathbf{x}_{0}\right)
\nonumber \\[4mm]
& \hspace*{-20mm} = \exp\left(- \frac{1}{\sigma^{2}} \mathbf{x}^{T}  \mathbf{H}^{T}\mathbf{H} \mathbf{x}_{0} + \mathbf{x}_{1}\widetilde{\mathbf{H}}^{\dagger} \mathbf{x}+   \frac{1}{\sigma^{2}} \mathbf{x}_{0}^{T} \big( \widetilde{\mathbf{H}}^{\dagger}\big)^{T}\widetilde{\mathbf{H}}^{\dagger} \mathbf{x} \right) \sum_{\mathbf{p} \in \mathbb{Z}_{+}^{D}} \frac{a[\mathbf{p}]}{\mathbf{p}!} \bigg(\frac{1}{\sigma} \widetilde{\mathbf{H}}^{\dagger} \mathbf{x}\bigg)^{\mathbf{p}} \nonumber \\[4mm]
& \hspace*{-20mm} \stackrel{\eqref{equ_proof_isometr_R_g_quadr_form}}{=} \exp\left(- \frac{1}{\sigma^{2}} \mathbf{x}^{T}  \mathbf{H}^{T}\mathbf{H} \mathbf{x}_{0} + \mathbf{x}_{1}\widetilde{\mathbf{H}}^{\dagger} \mathbf{x}+   \frac{1}{\sigma^{2}} \mathbf{x}_{0}^{T} \mathbf{H}^{T}\mathbf{H} \mathbf{x} \right) \sum_{\mathbf{p} \in \mathbb{Z}_{+}^{D}} \frac{a[\mathbf{p}]}{\mathbf{p}!} \bigg(\frac{1}{\sigma} \widetilde{\mathbf{H}}^{\dagger} \mathbf{x}\bigg)^{\mathbf{p}} \nonumber \\[4mm]
& \hspace*{-20mm} = \exp\left(\mathbf{x}_{1}\widetilde{\mathbf{H}}^{\dagger} \mathbf{x}\right) \sum_{\mathbf{p} \in \mathbb{Z}_{+}^{D}} \frac{a[\mathbf{p}]}{\mathbf{p}!} \bigg(\frac{1}{\sigma} \widetilde{\mathbf{H}}^{\dagger} \mathbf{x}\bigg)^{\mathbf{p}}.
\end{align} 
Comparing \eqref{equ_suff_cond_valid_erervwhere_SLM_power_series_proof_gamma} with \eqref{equ_suff_cond_valid_erervwhere_SLM_power_series_proof_gamma_image_congruence} reveals that the prescribed mean function $\gamma(\cdot)$, which corresponds to the prescribed bias in \eqref{equ_prescr_bias_finite_order_polynom_SLM}, 
is the restriction of $\mathsf{K}_{g}[f(\cdot)](\mathbf{x}) \in  \mathcal{H}(\mathcal{M}_{\text{LGM}})$ to the subdomain $\mathcal{X}_{S}$, i.e., 
$\gamma(\cdot) = \mathsf{K}_{g}[f(\cdot)]\big|_{\mathcal{X}_{S}}$, which implies by Theorem \ref{thm_relation_RKHS_LGM_SLM} that 
$\gamma(\cdot) \in \mathcal{H}(\mathcal{M}_{\text{SLM}}(\mathbf{x}_{0}))$ and in turn that the prescribed bias 
$c(\mathbf{x})$ is valid for $\mathcal{M}_{\text{SLM}}(\mathbf{x}_{0})$ with $\mathbf{x}_{0} \in \mathcal{X}_{S}$.
\end{proof}

By the definition of a valid bias function in Definition \ref{def_valid_bias_func_classic_est}, we have trivially that any given estimator with finite variance for 
all $\mathbf{x} \in \mathbb{R}^{N}$ induces a valid bias function: 
\begin{lemma}
Consider an estimator $\hat{x}_{k}(\mathbf{y})$ for the LGM $\mathcal{E}_{\emph{LGM}}$ that has finite variance everywhere, i.e., $v(\hat{x}_{k}(\cdot); \mathbf{x}) < \infty$ for every $\mathbf{x} \in \mathbb{R}^{N}$. 
Then the bias function $c(\cdot): \mathcal{X}_{S} \rightarrow \mathbb{R}: c(\mathbf{x}) = b(\hat{x}_{k}(\cdot); \mathbf{x})$ is valid for $\mathcal{M}_{\emph{SLM}}(\mathbf{x}_{0})$ (where $\sigma$, $M$, $N$ and $\mathbf{H}$ are the same as for $\mathcal{E}_{\emph{LGM}}$) with arbitrary $\mathbf{x}_{0} \in \mathcal{X}_{S}$.
\end{lemma}

\section{Lower Bounds on the Estimator Variance for the SLM}
\label{sec_RKHS_Lower_Bounds_SLM}
In this section we will derive lower bounds on the minimum achievable variance $L_{\mathcal{M}_{\text{SLM}}}$ where 
the parameters $\sigma$, $S$, $M$, $N$, $\mathbf{H}$, $c(\cdot): \mathcal{X}_{S} \rightarrow \mathbb{R}$, and 
$\mathbf{x}_{0} \in \mathcal{X}_{S}$ of the minimum variance problem $\mathcal{M}_{\text{SLM}}=\left( \mathcal{E}_{\text{SLM}},c(\cdot),\mathbf{x}_{0}\right)$ are assumed arbitrary but fixed. 
Furthermore, we will assume from now on that the prescribed bias function $c(\cdot)$ is valid. This is no real restriction since the lower bounds that will be derived are 
finite and therefore trivially apply also formally if $c(\cdot)$ is not valid, in which case $L_{\mathcal{M}_{\text{SLM}}}=\infty$ by definition (see Definition \ref{def_min_ach_var}). 

The various bounds that will be presented are similar in that they are based on a projection of the 
prescribed mean function $\gamma(\cdot): \mathcal{X}_{S} \rightarrow \mathbb{R}: \gamma(\mathbf{x}) \triangleq c(\mathbf{x}) + x_{k}$ 
(which satisfies $\gamma(\cdot) \in \mathcal{H}(\mathcal{M}_{\text{SLM}})$ if $c(\cdot)$ is valid) onto some 
subspace $\mathcal{U}$ of the RKHS $\mathcal{H}(\mathcal{M}_{\text{SLM}})$. Indeed (as already discussed after Theorem \ref{thm_main_facts_RKHS_MVE}), by Theorem \ref{thm_main_facts_RKHS_MVE} and Theorem \ref{thm_orthog_proj_ineq}, we have 
for an arbitrary subspace $\mathcal{U} \subseteq \mathcal{H}(\mathcal{M}_{\text{SLM}})$ the bound 
\begin{equation}
L_{\mathcal{M}_{\text{SLM}}} = \| \gamma(\cdot) \|^{2}_{\mathcal{H}(\mathcal{M}_{\text{SLM}})} - \big[ \gamma(\mathbf{x}_{0}) \big]^{2} \geq 
\| \mathbf{P}_{\mathcal{U}} \gamma(\cdot) \|^{2}_{\mathcal{H}(\mathcal{M}_{\text{SLM}})} - \big[ \gamma(\mathbf{x}_{0}) \big]^{2}.
\end{equation}

The first bound is based on a generalization of the CRB and has been previously presented in \cite{ZvikaCRB}:
\begin{theorem}
\label{thm_CRB_SLM}
If the prescribed bias $c(\cdot): \mathcal{X}_{S} \rightarrow \mathbb{R}$ is such that the partial derivatives $\frac{\partial c(\mathbf{x})}{\partial x_{l}} \big|_{\mathbf{x} = \mathbf{x}_{0}}$ exist for $l \in [N]$, then 
\begin{align} 
L_{\mathcal{M}_{\emph{SLM}}} & \geq \sigma^{2} \mathbf{b}^{T} \left( \mathbf{H}^{T} \mathbf{H} \right)^{\dagger} \mathbf{b}   \hspace*{19mm} \mbox{when}\,\, \| \mathbf{x}_{0} \|_{0} < S \label{equ_CRB_SLM_not_full_sparsity} \\Ê
 L_{\mathcal{M}_{\emph{SLM}}} & \geq  \sigma^{2}\mathbf{b}_{\mathbf{x}_{0}}^{T} \left( \mathbf{H}_{\mathbf{x}_{0}}^{T} \mathbf{H}_{\mathbf{x}_{0}} \right)^{\dagger} \mathbf{b}_{\mathbf{x}_{0}} \hspace*{10mm} \mbox{when}\,\, \| \mathbf{x}_{0} \|_{0} = S, \label{equ_CRB_SLM_full_sparsity} 
\end{align} 
where $\mathbf{b} \in \mathbb{R}^{N}$ is defined elementwise by $b_{l} \triangleq  \delta_{k,l} + \frac{\partial c(\mathbf{x})}{\partial x_{l}} \big|_{\mathbf{x} = \mathbf{x}_{0}}$ and 
$\mathbf{b}_{\mathbf{x}_{0}} \in \mathbb{R}^{S}$, $\mathbf{H}_{\mathbf{x}_{0}} \in \mathbb{R}^{M \times S}$ denote the restrictions to the entries and columns of $\mathbf{b}$ and $\mathbf{H}$, respectively, which are indexed by $\supp(\mathbf{x}_{0})=(i_{1},\ldots,i_{S})$, i.e., $\big(\mathbf{b}_{\mathbf{x}_{0}}\big)_{j} = b_{i_{j}}$. 
\end{theorem} 
\begin{proof}
Since $c(\cdot)$ is assumed valid, we have by Theorem \ref{thm_main_facts_RKHS_MVE} that the prescribed mean function $\gamma(\cdot): \mathcal{X}_{S} \rightarrow \mathbb{R}: \gamma(\mathbf{x}) = c(\mathbf{x}) + x_{k}$ 
belongs to the RKHS $\mathcal{H}(\mathcal{M}_{\text{SLM}})$. 

For the case $\| \mathbf{x}_{0} \|_{0} < S$, consider the subspace $\mathcal{U}_{1} \triangleq \linspan \big\{ \{v_{0}(\cdot) \} \cup \{ v_{l}(\cdot) \big\}_{l \in [N]}Ê\big\}$ 
spanned by the functions $v_{0} (\cdot) \triangleq R_{\mathcal{M}_{\text{SLM}}} (\cdot, \mathbf{x}_{0})$ and 
\begin{equation}
v_{l}(\cdot) \triangleq \frac{ \partial^{\mathbf{e}_{l}} R_{\mathcal{M}_{\text{SLM}}} (\cdot, \mathbf{x}_{2})}  
{\partial \mathbf{x}_{2}^{\mathbf{e}_{l}}} \bigg|_{\mathbf{x}_{2} = \mathbf{x}_{0}} \quad \mbox{, }  l \in [N].
\end{equation} 
We have trivially $v_{0}(\cdot) \in \mathcal{H}(\mathcal{M}_{\text{SLM}})$ and by Theorem \ref{thm_der_repr_prop} we have also that $v_{l} (\cdot)  \in \mathcal{H}(\mathcal{M}_{\text{SLM}})$ for $l \in [N]$. 
In a completely analogous manner as the RKHS-based derivation of Theorem \ref{thm_unconstr_CR} in Section \ref{sec_CRB}, one can show that 
$\big\langle v_{0}(\cdot), v_{l} (\cdot) \big\rangle_{ \mathcal{H}(\mathcal{M}_{\text{SLM}})} = 0$ for $l \in [N]$ (see \eqref{equ_CR_orthog_vecs_1}),  
$\big\langle v_{l}(\cdot), v_{l'}(\cdot) \big\rangle_{ \mathcal{H}(\mathcal{M}_{\text{SLM}})} = \sigma^{2} \left( \mathbf{H}^{T} \mathbf{H} \right)_{l,l'}$ for $l,l' \in [N]$ (see \eqref{equ_gramian_vectors_crb}), 
and $\big\langle  v_{l}(\cdot), \gamma(\cdot) \big\rangle_{ \mathcal{H}(\mathcal{M}_{\text{SLM}})} = b_{l}$ for $l \in [N]$ (see \eqref{equ_inner_prod_par_der_derivation_UCRB}). 
The bound \eqref{equ_CRB_SLM_not_full_sparsity} is then obtained via Theorem \ref{thm_main_facts_RKHS_MVE}, Theorem \ref{thm_orthog_proj_ineq}, and Theorem \ref{thm_norm_projection_finite_dim_subspace} 
by projecting $\gamma(\cdot)$ on the subspace $\mathcal{U}_{1}$ (see \eqref{equ_proof_unconstrained_CRB_projection_subspace_main_facts}). 

In an almost identical manner one can prove the bound \eqref{equ_CRB_SLM_full_sparsity} for the case $\| \mathbf{x}_{0} \|_{0} = S$. Indeed, the only difference is that 
instead of using the subspace $\mathcal{U}_{1}$ one has to use the subspace $\mathcal{U}_{2} \triangleq \linspan \big\{ \{ v_{0}(\cdot) \} \cup \{ v_{l}(\cdot) \}_{l \in \supp(\mathbf{x}_{0})}Ê\big\}$.
\end{proof} 
The formulation of the bound in \cite{ZvikaCRB} slightly differs from Theorem \ref{thm_CRB_SLM} in that the authors of \cite{ZvikaCRB} present 
a bound on the variance $v(\hat{\mathbf{x}}(\cdot); \mathbf{x}_{0})$ of a vector-valued estimator that estimates $\mathbf{x}$ itself and not only the $k$th coefficient $x_{k}$ as we consider. 
However, according to Section \ref{sec_sep_vector_scalar_param_function}, by summing the bound in Theorem \ref{thm_CRB_SLM} over all individual indices $k \in [N]$, 
one obtains the bound presented in \cite{ZvikaCRB}. This is just another consequence of the fact that minimum variance estimation of the parameter vector $\mathbf{x}$ 
is equivalent to separate minimum variance estimation of the coefficients $x_{k}$.

We note two important aspects of Theorem \ref{thm_CRB_SLM}: 
First, if the matrix $\mathbf{H}$ has full column rank, i.e, $\rank(\mathbf{H}) = N$, and we consider unbiased estimation, i.e., $c(\cdot) \equiv 0$, then the bound in \eqref{equ_CRB_SLM_not_full_sparsity} 
coincides with the variance of the well-known least-squares (LS) estimator \cite{kay,scharf91} given as 
\begin{equation}
\label{equ_def_ord_LS_est}
\hat{x}_{k,\text{LS}} (\mathbf{y}) \triangleq \mathbf{e}_{k}^{T} (\mathbf{H}^{T} \mathbf{H})^{-1} \mathbf{H}^{T} \mathbf{y}, 
\end{equation}
which implies that in this case the inequality in \eqref{equ_CRB_SLM_not_full_sparsity} becomes an equality, i.e., $L_{\mathcal{M}_{\text{SLM}}} = \sigma^{2} \mathbf{b}^{T} \left( \mathbf{H}^{T} \mathbf{H} \right)^{\dagger} \mathbf{b}$.
Indeed, a straightforward calculation reveals that $\hat{x}_{k,\text{LS}}(\cdot)$ is unbiased and has variance 
\begin{equation} 
\label{equ_var_ord_LS_est}
v(\hat{x}_{k,\text{LS}}(\cdot); \mathbf{x}_{0})= \sigma^{2} \mathbf{e}_{k}^{T} \left( \mathbf{H}^{T} \mathbf{H} \right)^{\dagger} \mathbf{e}_{k} =\sigma^{2} \mathbf{e}_{k}^{T} \left( \mathbf{H}^{T} \mathbf{H} \right)^{-1} \mathbf{e}_{k},
\end{equation} 
where the last equality follows because $\rank(\mathbf{H}) = N$. Thus, under these conditions, the LS estimator is the LMVU estimator for the SLM at 
every parameter vector $\mathbf{x}_{0} \in \mathcal{X}_{S}$ with 
$\| \mathbf{x}_{0} \|_{0} < S$.
Since (by Theorem \ref{thm_uniqueness_LMV}) it is the unique LMVU estimator for $\mathbf{x}_{0}$ with $\|\mathbf{x}_{0}\|_{0} < S$, we have that if there existed a UMVU estimator, it necessarily would be the ordinary LS estimator $\hat{x}_{k,\text{LS}}(\cdot)$. Note that the LS estimator does 
not exploit the sparsity assumption expressed by the parameter set $\mathcal{X}_{S}$, and it has the constant variance \eqref{equ_var_ord_LS_est} for every $\mathbf{x} \in \mathcal{X}_{S}$. Thus, even if there existed a UMVU estimator it would not be able to take advantage of the sparsity assumption.
However, as already shown in \cite{AlexZvikaICASSP,AlexZvikaJournal} and discussed in Section \ref{sec_SSNM}, there does not exist a UMVU estimator for the SLM in 
general. 

A second important aspect of Theorem \ref{thm_CRB_SLM} is that the lower bound given by \eqref{equ_CRB_SLM_not_full_sparsity} and \eqref{equ_CRB_SLM_full_sparsity} is not a continuous function of $\mathbf{x}_{0}$ in general.
Indeed, for the case $\mathbf{H}=\mathbf{I}$ and $c(\cdot) Ê\equiv 0$ considered in \cite{AlexZvikaICASSP,AlexZvikaJournal}, it can be verified that the bound 
is a strictly upper semi-continuous function of $\mathbf{x}_{0}$. 
Consider e.g. $\mathcal{M}_{\text{SLM}}=\left( \mathcal{E}_{\text{SLM}},c(\cdot),\mathbf{x}_{0}\right)$ with $\mathbf{H} = \mathbf{I}$, $S=1$, $M=N=2$, $k=2$, $c(\cdot)Ê\equiv 0$, and 
a parameter vector $\mathbf{x}_{0}$ of the form $\mathbf{x}_{0} = a (1,0)^{T}$ where $a \in \mathbb{R}_{+}$. The corresponding bound 
given by Theorem \ref{thm_CRB_SLM} is equal to $0$ for all $a>0$ but equals $1$ for $a=0$.
However, since by Theorem \ref{thm_SLM_min_achiev_var_lower_semi_cont} the minimum achievable variance $L_{\mathcal{M}_{\text{SLM}}}$ has to be a \emph{lower} semi-continuous function of $\mathbf{x}_{0}$, we have 
that the bound in Theorem \ref{thm_CRB_SLM} cannot be tight, i.e., we have a strict inequality in \eqref{equ_CRB_SLM_not_full_sparsity} or \eqref{equ_CRB_SLM_full_sparsity} in general. 

While the bound in Theorem \ref{thm_CRB_SLM} was based on a generalization of the CRB (see Section \ref{sec_CRB}), the authors of \cite{AlexZvikaICASSP} present a lower bound on 
$L_{\mathcal{M}_{\text{SLM}}}$ for $\mathbf{H}=\mathbf{I}$ and $c(\cdot) \equiv 0$ that is 
based on a generalization of the HCRB (see Section \ref{sec_HCRB}):
\begin{theorem}
\label{thm_HCRB_SLM}
Consider the minimum variance problem $\mathcal{M}_{\emph{SLM}}=\left( \mathcal{E}_{\emph{SLM}},c(\cdot),\mathbf{x}_{0}\right)$ with $\mathbf{H} = \mathbf{I}$ and $c(\cdot) \equiv 0$. 
Then we have the bound  
\begin{align} 
\label{equ_HCRB_SLM_not_full_sparsity}  & L_{\mathcal{M}_{\emph{SLM}}} \geq \sigma^{2}   \hspace*{50mm} \mbox{when}\,\, |\supp(\mathbf{x}_{0}) \cup \{k\} | < S+1  \\Ê
\label{equ_HCRB_SLM_full_sparsity} & L_{\mathcal{M}_{\emph{SLM}}} \geq \sigma^{2} \frac{N-S-1}{N-S}\exp \left(-\xi_{0}^{2}/ \sigma^{2} \right) \hspace*{8mm}  \mbox{when}\,\, |\supp(\mathbf{x}_{0}) \cup \{k\} | = S+1,
\end{align} 
where $\xi_{0}$ denotes the value of the $S$-largest (in magnitude) entry of $\mathbf{x}_{0}$. 
\end{theorem} 
\begin{proof}
The bias function $c(\cdot) \equiv 0$ is assumed to be valid, i.e., the mean function $\gamma(\cdot): \mathcal{X}_{S} \rightarrow \mathbb{R}: \gamma(\mathbf{x}) = x_{k}$ belongs to $\mathcal{H}(\mathcal{M}_{\text{SLM}})$. 

First consider the case where $|\supp(\mathbf{x}_{0}) \cup \{k\} | = S+1$ (which can only occur if $N > S$). 
In order to show \eqref{equ_HCRB_SLM_full_sparsity} we denote by $j_{0}$ and $\xi_{0}$ respectively the index and value of the $S$-largest (in magnitude) entry of $\mathbf{x}_{0}$. 
We then introduce the subspace 
$\mathcal{U}_{1}^{(t)} \triangleq \linspan \big\{ v_{0}(\cdot) \cup \{ v^{(t)}_{l}(\cdot) \}_{l \in [N]}Ê\big\}$ which is parametrized by $t \in \mathbb{R}$ and spanned by the vectors 
$v_{0}(\cdot) \triangleq R_{\mathcal{M}_{\text{SLM}}} (\cdot, \mathbf{x}_{0})$, $v_{l}^{(t)}(\cdot) \triangleq R_{\mathcal{M}_{\text{SLM}}} (\cdot, \mathbf{x}_{0}-\mathbf{e}_{j_{0}}\xi_{0}+t\mathbf{e}_{l})-R_{\mathcal{M}_{\text{SLM}}} (\cdot, \mathbf{x}_{0})$ for $l \in [N] \setminus \supp(\mathbf{x}_{0})$ and 
$v_{l}^{(t)}(\cdot) \triangleq R_{\mathcal{M}_{\text{SLM}}} (\cdot, \mathbf{x}_{0}+t\mathbf{e}_{l})-R_{\mathcal{M}_{\text{SLM}}} (\cdot, \mathbf{x}_{0})$ for $l \in \supp(\mathbf{x}_{0})$. 
Similarly to the RKHS-based derivation of the HCRB in Theorem \ref{thm_HCRB}, one can show by a projection of the mean function $\gamma(\cdot) \in \mathcal{H}(\mathcal{M}_{\text{SLM}})$ on the subspace $\mathcal{U}^{(t)}$ that (see \eqref{equ_min_var_bound_HCRB})
\begin{equation}
\label{equ_derivation_SLM_HCRB_bound_1}
L_{\mathcal{M}_{\text{SLM}}} \geq  t^2 \mathbf{e}_{k}^{T} \mathbf{V}^{\dagger} \mathbf{e}_{k},
\end{equation} 
where $\mathbf{V} \in \mathbb{R}^{N \times N}$ is defined elementwise as 
\begin{align} 
\left(\mathbf{V} \right)_{m,n} & \triangleq \big\langle v^{(t)}_{m}(\cdot), v^{(t)}_{n}(\cdot) \big\rangle_{\mathcal{H}(\mathcal{M}_{\text{SLM}})}   \nonumber \\[4mm]
& = \big\langle R_{\mathcal{M}_{\text{SLM}}} (\cdot, \mathbf{x}_{0}+t\mathbf{e}_{m})-R_{\mathcal{M}_{\text{SLM}}} (\cdot, \mathbf{x}_{0}), R_{\mathcal{M}_{\text{SLM}}} (\cdot, \mathbf{x}_{0}+t\mathbf{e}_{n})-R_{\mathcal{M}_{\text{SLM}}} (\cdot, \mathbf{x}_{0}) \big\rangle_{\mathcal{H}(\mathcal{M}_{\text{SLM}})}   \nonumber \\[4mm]
& \stackrel{\eqref{equ_reproduction_property}}{=}Ê R_{\mathcal{M}_{\text{SLM}}} (\mathbf{x}_{0}+t\mathbf{e}_{m}, \mathbf{x}_{0}+t\mathbf{e}_{n}) 
- R_{\mathcal{M}_{\text{SLM}}} (\mathbf{x}_{0}, \mathbf{x}_{0}+t\mathbf{e}_{m}) \nonumber \\[4mm]
& - R_{\mathcal{M}_{\text{SLM}}} (\mathbf{x}_{0}+t\mathbf{e}_{n}, \mathbf{x}_{0}) + 
R_{\mathcal{M}_{\text{SLM}}} (\mathbf{x}_{0}, \mathbf{x}_{0})  \nonumber \\[4mm]
& \stackrel{\eqref{equ_kernel_one_arg_x_0_equal_1}}{=} \exp( t^{2} \mathbf{e}_{n}^{T}\mathbf{e}_{m}/ \sigma^{2}) - 1 -1 + 1  \nonumber\\[4mm]
& = \delta_{m,n} \exp(t^{2}/\sigma^{2}) -1 .
\end{align}  
As shown in \cite{AlexZvikaJournal}, the bound in \eqref{equ_HCRB_SLM_full_sparsity} is obtained from \eqref{equ_derivation_SLM_HCRB_bound_1} as the limit $\lim_{t \rightarrow 0} t^2 \mathbf{e}_{k}^{T} \mathbf{V}^{\dagger} \mathbf{e}_{k}$. 

Now we consider the remaining case where $|\supp(\mathbf{x}_{0}) \cup \{k\} | < S+1$. Instead of $\mathcal{U}_{1}^{(t)}$, we will now use 
the subspace $\mathcal{U}_{2}^{(t)} \triangleq \{ v_{0}(\cdot), v^{(t)}(\cdot) Ê\}$, where $v^{(t)}(\cdot) \triangleq R_{\mathcal{M}_{\text{SLM}}} (\cdot, \mathbf{x}_{0}+t\mathbf{e}_{k})-R_{\mathcal{M}_{\text{SLM}}} (\cdot, \mathbf{x}_{0})$.
In an analogous manner as before, one obtains 
\begin{equation}
\label{equ_derivation_SLM_HCRB_bound_2}
L_{\mathcal{M}_{\text{SLM}}} \geq \frac{ t^2}{\|v^{(t)}(\cdot) \|^{2}_{\mathcal{H}(\mathcal{M}_{\text{SLM}})}}.
\end{equation} 
The squared norm of $v^{(t)}(\cdot)$ can be calculated as 
\begin{align}
 \|v^{(t)}(\cdot) \|^{2}_{\mathcal{H}(\mathcal{M}_{\text{SLM}})} & = 
\big\langle v^{(t)}(\cdot), v^{(t)}(\cdot) \big\rangle_{\mathcal{H}(\mathcal{M}_{\text{SLM}})}   \nonumber \\[4mm]Ê
& \stackrel{\eqref{equ_reproduction_property}}{=}Ê R_{\mathcal{M}_{\text{SLM}}} (\mathbf{x}_{0}+t\mathbf{e}_{k}, \mathbf{x}_{0}+t\mathbf{e}_{k}) 
- R_{\mathcal{M}_{\text{SLM}}} (\mathbf{x}_{0}, \mathbf{x}_{0}+t\mathbf{e}_{k})  \nonumber \\[4mm]
& - R_{\mathcal{M}_{\text{SLM}}} (\mathbf{x}_{0}+t\mathbf{e}_{k}, \mathbf{x}_{0}) + 
R_{\mathcal{M}_{\text{SLM}}} (\mathbf{x}_{0}, \mathbf{x}_{0})  \nonumber \\[4mm]Ê
& \stackrel{\eqref{equ_kernel_one_arg_x_0_equal_1}}{=} \exp( t^{2} \mathbf{e}_{k}^{T}\mathbf{e}_{k}/ \sigma^{2}) - 1 - 1 + 1  \nonumber \\[4mm]
& =  \exp(t^{2}/\sigma^{2}) -1.
\end{align} 
Thus, one straightforwardly obtains the bound $L_{\mathcal{M}_{\text{SLM}}} \geq \sigma^{2}$ in \eqref{equ_HCRB_SLM_not_full_sparsity} by the limit 
\begin{equation} 
\lim_{t \rightarrow 0}   t^2/[\exp(t^{2}/ \sigma^{2}) -1].
\end{equation}
Since the variance of the specific unbiased estimator $\hat{g}(\mathbf{y}) = y_{k}$ for the SLM is equal to $\sigma^{2}$, we have that $L_{\mathcal{M}_{\text{SLM}}} = \sigma^{2}$. 
\end{proof} 
The bound given in Theorem \ref{thm_HCRB_SLM} has been presented in a slightly different form in \cite{AlexZvikaICASSP,AlexZvikaJournal}, where 
a lower bound on the variance $v(\hat{\mathbf{x}}(\cdot), \mathbf{x}_{0})$ of an unbiased estimator of the parameter vector $\mathbf{x}$ itself, instead of the single coefficient $x_{k}$, has been derived.\footnote{More precisely, 
a lower bound on the MSE of any unbiased estimator of $\mathbf{x}$ was presented in \cite{AlexZvikaICASSP,AlexZvikaJournal}. However, the variance of an unbiased estimator is equal to its MSE.}
As discussed in Section \ref{sec_sep_vector_scalar_param_function}, minimum variance estimation of the parameter vector is equivalent to separate minimum variance estimation of the coefficients $x_{k}$.
In particular, the bound given in \cite{AlexZvikaICASSP,AlexZvikaJournal} is obtained by summing the bound in Theorem \ref{thm_HCRB_SLM} over all indices $k \in [N]$.

While in general the bound in Theorem \ref{thm_HCRB_SLM} is tighter, i.e., higher than the bound of Theorem \ref{thm_CRB_SLM} (when specialized to $c(\cdot) \equiv 0$ and $\mathbf{H} = \mathbf{I}$), it 
is again upper semi-continuous. This again implies by Theorem \ref{thm_SLM_min_achiev_var_lower_semi_cont} that the bound cannot be tight, i.e., we have a strict inequality in \eqref{equ_HCRB_SLM_full_sparsity} in general. 
However, for the case 
$|\supp(\mathbf{x}_{0}) \cup \{k\} | < S+1$,\footnote{Note that the case $|\supp(\mathbf{x}_{0}) \cup \{k\} | < S+1$ can only occur if either $\| \mathbf{x} \|_{0} < S$ or $\| \mathbf{x} \|_{0} = S$ and additionally $k \in \supp(\mathbf{x}_{0})$.}
the bound in \eqref{equ_HCRB_SLM_not_full_sparsity} is tight since, as already mentioned above, 
it is achieved by the ordinary LS estimator $\hat{x}_{k,\text{LS}}(\mathbf{y}) = y_{k}$.

An improved lower bound on $L_{\mathcal{M}_{\text{SLM}}}$ has been presented in \cite{RKHSAsilomar2010}: 
\begin{theorem} 
\label{thm_bound_asilomar}
Consider the minimum variance problem $\mathcal{M}_{\emph{SLM}}=\left( \mathcal{E}_{\emph{SLM}},c(\cdot),\mathbf{x}_{0}\right)$ with a system matrix $\mathbf{H} \in \mathbb{R}^{M \times N}$ satisfying \eqref{equ_spark_cond} 
and with prescribed mean function $\gamma(\cdot): \mathcal{X}_{S} \rightarrow \mathbb{R}: \gamma(\mathbf{x}) = c(\mathbf{x}) + x_{k}$. 
We assume that the prescribed bias function $c(\cdot): \mathcal{X}_{S} \rightarrow \mathbb{R}$ is such that the partial derivatives $\frac{ \partial^{\mathbf{e}_{l}} c(\mathbf{x})}{\partial \mathbf{x}^{\mathbf{e}_{l}}}$ exist for all $l \in [N]$. 
Then, for an arbitrary set $\mathcal{K} =\{i_1,\ldots,i_{|\mathcal{K}|} \} \subseteq [N]$ consisting of no more than $S$ different indices, i.e., $|\mathcal{K}| \leq S$, we have the bound
\vspace*{2mm}
\begin{equation} 
\label{equ_bound_asilomar_1}
L_{\mathcal{M}_{\emph{SLM}}} \geq  \exp\left( - \frac{1}{\sigma^{2}}\| (\mathbf{I} - \mathbf{P}_{\mathcal{K}}) \mathbf{H} \mathbf{x}_{0} \|^{2}_{2} \right) \left[ \sigma^{2} \mathbf{r}_{\mathbf{x}_{0}}^{T} \left( \mathbf{H}_{\mathcal{K}}^{T} \mathbf{H}_{\mathcal{K}} \right)^{-1} \mathbf{r}_{\mathbf{x}_{0}} + \gamma^{2}(\widetilde{\mathbf{x}}_{0}) \right] - \gamma^{2}(\mathbf{x}_{0}), 
\vspace*{2mm}
\end{equation} 
where $\mathbf{P}_{\mathcal{K}} \triangleq \mathbf{H}_{\mathcal{K}} \mathbf{H}_{\mathcal{K}}^{\dagger} \in \mathbb{R}^{M \times M}$, $\mathbf{r}_{\mathbf{x}_{0}} \in \mathbb{R}^{|\mathcal{K}|}$ is defined elementwise as $r_{\mathbf{x}_{0},l} \triangleq \frac{ \partial^{\mathbf{e}_{i_{l}}} \gamma(\mathbf{x})}{\partial \mathbf{x}^{\mathbf{e}_{i_{l}}}} \big|_{\mathbf{x} = \widetilde{\mathbf{x}}_{0}}$, and $\widetilde{\mathbf{x}}_{0} \in \mathbb{R}^{N}$ is defined as the unique (due to \eqref{equ_spark_cond}) 
vector with $\supp(\widetilde{\mathbf{x}}_{0})=\mathcal{K}$ that solves 
\begin{equation} 
\label{equ_def_x_0_tilde_asiloamr_bound_1}
\mathbf{H} \widetilde{\mathbf{x}}_{0} = \mathbf{H}_{\mathcal{K}} \mathbf{H}_{\mathcal{K}}^{\dagger} \mathbf{H} \mathbf{x}_{0}. 
\end{equation} 
\end{theorem}

\begin{proof}
First we note that since $\mathbf{H}$ is assumed to satisfy \eqref{equ_spark_cond}, the matrix $\mathbf{H}_{\mathcal{K}} \in \mathbb{R}^{M \times |\mathcal{K}|}$ has full column rank, which implies that 
$\mathbf{H}_{\mathcal{K}}^{\dagger} = (\mathbf{H}_{\mathcal{K}}^{T} \mathbf{H}_{\mathcal{K}})^{-1} \mathbf{H}_{\mathcal{K}}^{T}$ \cite{golub96}. Thus, 
\begin{align}
\label{equ_proof_lower_bound_asilomar_relation_1}
\mathbf{H}_{\mathcal{K}}^{T} \mathbf{H} (\widetilde{\mathbf{x}}_{0} - \mathbf{x}_{0}) & = \mathbf{H}_{\mathcal{K}}^{T}( \mathbf{H}_{\mathcal{K}} \mathbf{H}_{\mathcal{K}}^{\dagger} \mathbf{H} \mathbf{x}_{0} - \mathbf{H}\mathbf{x}_{0})   = \mathbf{H}_{\mathcal{K}}^{T}\mathbf{H}_{\mathcal{K}} \mathbf{H}_{\mathcal{K}}^{\dagger} \mathbf{H} \mathbf{x}_{0} - \mathbf{H}_{\mathcal{K}}^{T}\mathbf{H} \mathbf{x}_{0} \nonumber \\[4mm]
&=  \mathbf{H}_{\mathcal{K}}^{T}\mathbf{H}_{\mathcal{K}}(\mathbf{H}_{\mathcal{K}}^{T} \mathbf{H}_{\mathcal{K}})^{-1} \mathbf{H}_{\mathcal{K}}^{T} \mathbf{H} \mathbf{x}_{0} - \mathbf{H}_{\mathcal{K}}^{T}\mathbf{H} \mathbf{x}_{0} =   \mathbf{H}_{\mathcal{K}}^{T} \mathbf{H} \mathbf{x}_{0} - \mathbf{H}_{\mathcal{K}}^{T}\mathbf{H} \mathbf{x}_{0} = \mathbf{0}, 
\end{align} 
which implies in turn that $\linspan(\mathbf{H}_{\mathcal{K}})$ is orthogonal to $\mathbf{H} (\widetilde{\mathbf{x}}_{0} - \mathbf{x}_{0})$. 
Let us then consider the subspace 
\begin{equation} 
\label{equ_proof_lower_bound_asilomar_def_subspace}
\mathcal{U}_{3} \triangleq \linspan \big\{ v_{0} (\cdot) \cup \{ v_{l}(\cdot) \}_{l \in \mathcal{K}}Ê\big\}
\end{equation}
spanned by the functions $v_{0}(\cdot) \triangleq R_{\mathcal{M}_{\text{SLM}}}(\cdot, \widetilde{\mathbf{x}}_{0})$ and 
$v_{l}(\cdot) \triangleq \frac{\partial^{\mathbf{e}_{l}}R_{\mathcal{M}_{\text{SLM}}}(\cdot, \mathbf{x}_{2})}{\partial \mathbf{x}_{2}^{\mathbf{e}_{l}}}Ê\big|_{\mathbf{x}_{2} = \widetilde{\mathbf{x}}_{0}}$, which by Theorem \ref{thm_der_repr_prop} belong to the RKHS $\mathcal{H}(\mathcal{M}_{\text{SLM}})$.
The inner product between any $v_{l}(\cdot)$ with $l \in \mathcal{K}$ and $v_{0}(\cdot)$ is given due to Theorem \ref{thm_der_repr_prop} by 
\begin{align} 
\label{equ_proof_lower_bound_asilomar_orthog_vecs_1}
\big\langle v_{0}(\cdot) , v_{l}(\cdot) \big\rangle_{\mathcal{H}(\mathcal{M}_{\text{SLM}})} & = \frac{\partial R_{\mathcal{M}_{\text{SLM}}}(\mathbf{x}_{1},\mathbf{x}_{2})}{\partial x_{2,l}} \bigg|_{\mathbf{x}_{1} = \mathbf{x}_{2} = \widetilde{\mathbf{x}}_{0}}  
\nonumber \\[4mm]
& = \frac{1}{\sigma^{2}} ( \widetilde{\mathbf{x}}_{0} - \mathbf{x}_{0})^{T} \mathbf{H}^{T} \mathbf{H} \mathbf{e}_{l} \exp\left( \frac{1}{\sigma^{2}} \big\| \mathbf{H} ( \widetilde{\mathbf{x}}_{0} - \mathbf{x}_{0})\big\|^{2}_{2} \right)\stackrel{(a)}{=} 0, 
\end{align} 
where $(a)$ follows from the fact that $\mathbf{H} \mathbf{e}_{l}  \in \linspan(\mathbf{H}_{\mathcal{K}})$ and therefore, as shown by \eqref{equ_proof_lower_bound_asilomar_relation_1}, has to be orthogonal to $\mathbf{H} (\widetilde{\mathbf{x}}_{0} - \mathbf{x}_{0})$, i.e., 
$( \widetilde{\mathbf{x}}_{0} - \mathbf{x}_{0})^{T} \mathbf{H}^{T} \mathbf{H} \mathbf{e}_{l} =0$.
Now, we define the matrix $\mathbf{V}Ê\in \mathbb{R}^{|\mathcal{K}| \times |\mathcal{K}|}$ elementwise via
\begin{align} 
\label{equ_proof_lower_bound_asilomar_grammian}
\left( \mathbf{V} \right)_{m,n} & \triangleq \big\langle v_{i_m}(\cdot) , v_{i_n}(\cdot) \big\rangle_{\mathcal{H}(\mathcal{M}_{\text{SLM}})}  \nonumber \\[4mm]
& \stackrel{(a)}{=} \frac{\partial R_{\mathcal{M}_{\text{SLM}}}(\mathbf{x}_{1},\mathbf{x}_{2})}{\partial x_{1,i_n} \partial x_{2,i_m}} \bigg|_{\mathbf{x}_{1} = \mathbf{x}_{2} = \widetilde{\mathbf{x}}_{0}} \nonumber \\[4mm]
&= \frac{1}{\sigma^{2}} \frac{\partial \big[ (\mathbf{x}_{1}- \mathbf{x}_{0})^{T} \mathbf{H}^{T} \mathbf{H} \mathbf{e}_{i_{m}} \exp \left( \frac{1}{\sigma^{2}} (\mathbf{x}_{1} - \mathbf{x}_{0})^{T} \mathbf{H}^{T} \mathbf{H} (\widetilde{\mathbf{x}}_{0} - \mathbf{x}_{0}) \right)\big]}{\partial x_{1,i_n}} \bigg|_{\mathbf{x}_{1} = \widetilde{\mathbf{x}}_{0}} \nonumber \\[4mm] 
& =  \frac{1}{\sigma^{2}} \mathbf{e}_{i_n}^{T} \mathbf{H}^{T} \mathbf{H} \mathbf{e}_{i_{m}} \exp \left( \frac{1}{\sigma^{2}} (\widetilde{\mathbf{x}}_{0} - \mathbf{x}_{0})^{T} \mathbf{H}^{T} \mathbf{H} (\widetilde{\mathbf{x}}_{0} - \mathbf{x}_{0}) \right) \nonumber \\[4mm] 
& \hspace*{3mm} + \frac{1}{\sigma^{4}}  \bigg( (\widetilde{\mathbf{x}}_{0}- \mathbf{x}_{0})^{T} \mathbf{H}^{T} \mathbf{H} \mathbf{e}_{i_{m}} \bigg) \mathbf{e}^{T}_{i_{n}}\mathbf{H}^{T} \mathbf{H} (\widetilde{\mathbf{x}}_{0} - \mathbf{x}_{0}) \exp \left( \frac{1}{\sigma^{2}} (\widetilde{\mathbf{x}}_{0}- \mathbf{x}_{0})^{T} \mathbf{H}^{T} \mathbf{H} (\widetilde{\mathbf{x}}_{0} - \mathbf{x}_{0}) \right)
 \nonumber \\[4mm] 
& \stackrel{(b)}{=} \frac{1}{\sigma^{2}} \mathbf{e}_{i_n}^{T} \mathbf{H}^{T} \mathbf{H}  \mathbf{e}_{i_{m}} \exp \left( \frac{1}{\sigma^{2}} (\widetilde{\mathbf{x}}_{0} - \mathbf{x}_{0})^{T} \mathbf{H}^{T} \mathbf{H} (\widetilde{\mathbf{x}}_{0} - \mathbf{x}_{0}) \right),
\end{align} 
where $(a)$ follows again from Theorem \ref{thm_der_repr_prop} and $(b)$ follows from the fact that $\mathbf{H} \mathbf{e}_{i_{m}} \! \in \! \linspan(\mathbf{H}_{\mathcal{K}})$, implying via \eqref{equ_proof_lower_bound_asilomar_relation_1} that $(\widetilde{\mathbf{x}}_{0}- \mathbf{x}_{0})^{T} \mathbf{H}^{T} \mathbf{H} \mathbf{e}_{i_{m}}=0$. 
From \eqref{equ_proof_lower_bound_asilomar_grammian}, we have by elementwise comparison that 
\begin{equation} 
\label{equ_proof_lower_bound_asilomar_grammian_matrix_complete_notation}
\mathbf{V} =  \frac{1}{\sigma^{2}} \mathbf{H}_{\mathcal{K}}^{T} \mathbf{H}_{\mathcal{K}} \exp \left( \frac{1}{\sigma^{2}}\big\| \mathbf{H} ( \widetilde{\mathbf{x}}_{0} - \mathbf{x}_{0})\big\|^{2}_{2} \right).
\end{equation} 
Since furthermore 
\begin{align} 
\label{equ_proof_lower_bound_asilomar_grammian_exponente_projection_identity}
\| \mathbf{H} (\widetilde{\mathbf{x}}_{0} - \mathbf{x}_{0})\|^{2}_{2}  
 \stackrel{\eqref{equ_def_x_0_tilde_asiloamr_bound_1}}{=} \| \mathbf{H} \mathbf{x}_{0} - \mathbf{H}_{\mathcal{K}}\mathbf{H}_{\mathcal{K}}^{\dagger} \mathbf{H} \mathbf{x}_{0} \|^{2}_{2} 
   = \| (\mathbf{I} - \mathbf{H}_{\mathcal{K}}\mathbf{H}_{\mathcal{K}}^{\dagger}) \mathbf{H} \mathbf{x}_{0} \|^{2}_{2}  = \| (\mathbf{I} -\mathbf{P}_{\mathcal{K}}) \mathbf{H} \mathbf{x}_{0} \|^{2}_{2},
\end{align} 
the bound in \eqref{equ_bound_asilomar_1} follows from Theorem \ref{thm_main_facts_RKHS_MVE} and Theorem \ref{thm_norm_projection_finite_dim_subspace_union_orthog_subspaces} 
(note that due to \eqref{equ_spark_cond} we have that $\left( \mathbf{H}_{\mathcal{K}}^{T} \mathbf{H}_{\mathcal{K}} \right)^{-1} =\left( \mathbf{H}_{\mathcal{K}}^{T} \mathbf{H}_{\mathcal{K}} \right)^{\dagger}$) 
via \eqref{equ_proof_lower_bound_asilomar_grammian}, \eqref{equ_proof_lower_bound_asilomar_orthog_vecs_1} and 
the inner products
\begin{equation} 
\label{equ_proof_lower_bound_asilomar_inner_prod_v_0_gamma}
\big\langle v_{0}(\cdot) , \gamma(\cdot) \big\rangle_{\mathcal{H}(\mathcal{M}_{\text{SLM}})}  = \gamma(\widetilde{\mathbf{x}}_{0})
\end{equation} 
and 
\begin{equation} 
\label{equ_proof_lower_bound_asilomar_inner_prod_v_l_gamma}
\big\langle v_{l}(\cdot) , \gamma(\cdot) \big\rangle_{\mathcal{H}(\mathcal{M}_{\text{SLM}})} = \frac{ \partial^{\mathbf{e}_{i_{l}}} \gamma(\mathbf{x})}{\partial \mathbf{x}^{\mathbf{e}_{i_{l}}}} \bigg|_{\mathbf{x} = \widetilde{\mathbf{x}}_{0}} =  r_{\mathbf{x}_{0},l}
\end{equation}  
(cf.\ Theorem \ref{thm_der_repr_prop}): 
\begin{align}
\label{equ_proof_lower_bound_asilomar_step_by_step}
L_{\mathcal{M}_{\text{SLM}}} & \stackrel{\eqref{equ_squared_norm_min_achiev_var}}{=} \| \gamma(\cdot) \|^{2}_{\mathcal{H}(\mathcal{M}_{\text{SLM}})} - \big[ \gamma(\mathbf{x}_{0}) \big]^{2} \stackrel{\eqref{equ_lower_bound_RKHS_projection}}{\geq}  \| \mathbf{P}_{\mathcal{U}_{3}} \gamma(\cdot) \|^{2}_{\mathcal{H}(\mathcal{M}_{\text{SLM}})}-  \big[ \gamma(\mathbf{x}_{0}) \big]^{2}  \nonumber \\[4mm] 
&  \hspace*{-5mm} \stackrel{\eqref{equ_idendity_projection_subspace_spanned_union_orthogonal_sets},\eqref{equ_proof_lower_bound_asilomar_inner_prod_v_l_gamma}}{=} \big\langle \gamma(\cdot), v_{0}(\cdot) \big\rangle_{\mathcal{H}(\mathcal{M}_{\text{SLM}})} \big( 
\big\langle v_{0}(\cdot), v_{0}(\cdot) \big\rangle_{\mathcal{H}(\mathcal{M}_{\text{SLM}})} \big)^{-1} 
\underbrace{\big\langle \gamma(\cdot), v_{0}(\cdot) \big\rangle_{\mathcal{H}(\mathcal{M}_{\text{SLM}})}}_{\stackrel{\eqref{equ_proof_lower_bound_asilomar_inner_prod_v_0_gamma}}{=} \gamma(\widetilde{\mathbf{x}}_{0})} \nonumber \\[4mm]
& + \mathbf{r}_{\mathbf{x}_{0}}^{T} \mathbf{V}^{\dagger} \mathbf{r}_{\mathbf{x}_{0}} - \big[ \gamma(\mathbf{x}_{0}) \big]^{2} \nonumber \\[4mm]
& \hspace*{-5mm}  \stackrel{(a)}{=} \exp\left( - \frac{1}{\sigma^{2}}\| (\mathbf{I} - \mathbf{P}_{\mathcal{K}}) \mathbf{H} \mathbf{x}_{0} \|^{2}_{2} \right)\big[  \gamma(\widetilde{\mathbf{x}}_{0}) \big]^{2} +  \mathbf{r}_{\mathbf{x}_{0}}^{T} \mathbf{V}^{\dagger} \mathbf{r}_{\mathbf{x}_{0}} -\big[ \gamma(\mathbf{x}_{0}) \big]^{2} \nonumber \\[4mm]
&\hspace*{-5mm}  \stackrel{\eqref{equ_proof_lower_bound_asilomar_grammian_matrix_complete_notation}}{=} \exp\left( - \frac{1}{\sigma^{2}}\| (\mathbf{I} - \mathbf{P}_{\mathcal{K}}) \mathbf{H} \mathbf{x}_{0} \|^{2}_{2} \right) \left[ \sigma^{2} \mathbf{r}_{\mathbf{x}_{0}}^{T} \left( \mathbf{H}_{\mathcal{K}}^{T} \mathbf{H}_{\mathcal{K}} \right)^{-1} \mathbf{r}_{\mathbf{x}_{0}} + \gamma^{2}(\widetilde{\mathbf{x}}_{0}) \right] - \big[ \gamma(\mathbf{x}_{0}) \big]^{2},
\end{align}
where step $(a)$ is due to 
\begin{align} 
\big\langle v_{0}(\cdot), v_{0}(\cdot) \big\rangle_{\mathcal{H}(\mathcal{M}_{\text{SLM}})}& = \big\langle R_{\mathcal{M}_{\text{SLM}}}(\cdot, \widetilde{\mathbf{x}}_{0}), R_{\mathcal{M}_{\text{SLM}}}(\cdot, \widetilde{\mathbf{x}}_{0}) \big\rangle_{\mathcal{H}(\mathcal{M}_{\text{SLM}})}
 \stackrel{\eqref{equ_reproduction_property}}{=} R_{\mathcal{M}_{\text{SLM}}}(\widetilde{\mathbf{x}}_{0}, \widetilde{\mathbf{x}}_{0}) \nonumber \\[4mm]
 & =  \exp \left( \frac{1}{\sigma^{2}} (\widetilde{\mathbf{x}}_{0} - \mathbf{x}_{0})^{T} \mathbf{H}^{T} \mathbf{H} (\widetilde{\mathbf{x}}_{0} - \mathbf{x}_{0}) \right)  = \exp \left( \frac{1}{\sigma^{2}}\big\| \mathbf{H} ( \widetilde{\mathbf{x}}_{0} - \mathbf{x}_{0})\big\|^{2}_{2} \right) \nonumber \\[4mm]
& \stackrel{\eqref{equ_proof_lower_bound_asilomar_grammian_exponente_projection_identity}}{=}  \exp\left( \frac{1}{\sigma^{2}}\| (\mathbf{I} - \mathbf{P}_{\mathcal{K}}) \mathbf{H} \mathbf{x}_{0} \|^{2}_{2} \right). \nonumber
\end{align}
\end{proof}

Note that the matrix $\mathbf{P}_{\mathcal{K}}= \mathbf{H}_{\mathcal{K}} \mathbf{H}_{\mathcal{K}}^{\dagger}$ defined in Theorem \ref{thm_bound_asilomar} 
is an orthogonal projection matrix \cite{golub96} on the subspace $\mathcal{U}_{\mathcal{K}} \triangleq \linspan(\mathbf{H}_{\mathcal{K}}) \subseteq \mathbb{R}^{M}$; 
consequently, $(\mathbf{I} - \mathbf{P}_{\mathcal{K}})$ is an orthogonal projection matrix on the orthogonal complement \cite{HalmosFiniteDimVecSpace,golub96} of $\mathcal{U}_{\mathcal{K}}$ and the 
norm $\|(\mathbf{I} - \mathbf{P}_{\mathcal{K}}) \mathbf{H} \mathbf{x}_{0} \|$ thus represents the distance between the point $\mathbf{H} \mathbf{x}_{0}$ and the subspace $\mathcal{U}_{\mathcal{K}}$ \cite{RudinBook}. 
Therefore, the factor $\exp\left( - \frac{1}{\sigma^{2}} \| (\mathbf{I} - \mathbf{P}_{\mathcal{K}}) \mathbf{H} \mathbf{x}_{0} \|^{2}_{2} \right)$ appearing in the  bound \eqref{equ_bound_asilomar_1} can be interpreted as a measure 
of the distance between the point $\mathbf{H} \mathbf{x}_{0}$ and the subspace $\mathcal{U}_{\mathcal{K}}$. In general, the bound \eqref{equ_bound_asilomar_1} is tighter, i.e. higher, if $\mathcal{K}$ is chosen such that this distance is small, i.e., $\| (\mathbf{I} - \mathbf{P}_{\mathcal{K}}) \mathbf{H} \mathbf{x}_{0} \|^{2}_{2}$ 
is small. 

One of the appealing properties of the bound \eqref{equ_bound_asilomar_1} is that it is continuous for the special case $\mathbf{H}=\mathbf{I}$ and $c(\cdot) \equiv 0$, 
irrespective of what index $k$ is chosen for the parameter function in the SLM (see \eqref{equ_def_SLM_est_problem}). 
We also have for this specific instance of $\mathcal{M}_{\text{SLM}}$ that the 
bound \eqref{equ_bound_asilomar_1} is tighter, i.e., higher, than the bounds given in Theorem \ref{thm_CRB_SLM} and Theorem \ref{thm_HCRB_SLM}. 

The quantity $\sigma^{2} \mathbf{r}_{\mathbf{x}_{0}}^{T} \left( \mathbf{H}_{\mathcal{K}}^{T} \mathbf{H}_{\mathcal{K}} \right)^{-1} \mathbf{r}_{\mathbf{x}_{0}}$ appearing in the bound \eqref{equ_bound_asilomar_1} 
has a remarkable interpretation as the unconstrained CRB (see Section \ref{sec_CRB}) for a minimum variance problem associated with the LGM with $N = |\mathcal{K}|$. 
In fact, for $k \in \supp(\mathbf{x}_{0})$, the factor $\exp\left( - \frac{1}{\sigma^{2}} \| (\mathbf{I} - \mathbf{P}_{\mathcal{K}}) \mathbf{H} \mathbf{x}_{0} \|^{2}_{2} \right)$ can be made equal to $1$ by choosing $\mathcal{K} = \supp(\mathbf{x}_{0})$. 
On the other hand, consider $k \notin \supp(\mathbf{x}_{0})$, and chose $\mathcal{K} =\{ \{k\} \cup \mathcal{L} \}$, where $\mathcal{L}$ denotes the indices of the $S-1$ largest (in magnitude) entries of $\mathbf{x}_{0}$. 
The bound \eqref{equ_bound_asilomar_1} then indicates a transition from a ``low''-SNR regime, where 
\begin{equation} 
\exp\left( - \frac{1}{\sigma^{2}} \| (\mathbf{I} - \mathbf{P}_{\mathcal{K}}) \mathbf{H} \mathbf{x}_{0} \|^{2}_{2} \right) \approx 1, \nonumber 
\end{equation} 
to a ``high''-SNR regime, where 
\begin{equation}
\exp\left( - \frac{1}{\sigma^{2}} \| (\mathbf{I} - \mathbf{P}_{\mathcal{K}}) \mathbf{H} \mathbf{x}_{0} \|^{2}_{2} \right) \approx 0. \nonumber
\end{equation}
In the low-SNR regime, the bound \eqref{equ_bound_asilomar_1} is equal to the CRB of a LGM with $N = |\mathcal{K}|$. In the 
high-SNR regime, the bound becomes approximately equal to $0$; this means that the nonzero entries $x_{k}$ with $k \notin \supp(\mathbf{x})$ can be estimated with small variance. 

We note that in \cite{RKHSAsilomar2010}, we presented a version of the bound \eqref{equ_bound_asilomar_1}, 
that is obtained by choosing the index set $\mathcal{K}$ for which the right hand side of \eqref{equ_bound_asilomar_1} is largest.
We note also that the proof of the bound \eqref{equ_bound_asilomar_1} used here is slightly more direct than that sketched in \cite{RKHSAsilomar2010}. 

By a slight modification of the proof of Theorem \ref{thm_bound_asilomar}, we obtain 
\begin{theorem} 
\label{thm_bound_asilomar_2}
Consider the minimum variance problem $\mathcal{M}_{\emph{SLM}}=\left( \mathcal{E}_{\emph{SLM}},c(\cdot),\mathbf{x}_{0}\right)$ with a system matrix $\mathbf{H}$ 
satisfying \eqref{equ_spark_cond} and with prescribed mean function $\gamma(\cdot): \mathcal{X}_{S} \rightarrow \mathbb{R}: \gamma(\mathbf{x}) = c(\mathbf{x}) + x_{k}$.
We assume that $c(\cdot)$ is valid and such that the partial derivatives $\frac{ \partial^{\mathbf{e}_{l}} c(\mathbf{x})}{\partial \mathbf{x}^{\mathbf{e}_{l}}}$ exist for all $l \in [N]$. 
Then, for an arbitrary set $\mathcal{K} =\{i_1,\ldots,i_{|\mathcal{K}|} \} \subseteq [N]$ consisting of no more than $S$ different indices, i.e., $|\mathcal{K}| \leq S$, we have the bound
\begin{equation} 
\label{equ_bound_asilomar_2}
L_{\mathcal{M}_{\emph{SLM}}} \geq  \exp\left( - \frac{1}{\sigma^{2}} \| (\mathbf{I} - \mathbf{P}_{\mathcal{K}}) \mathbf{H} \mathbf{x}_{0} \|^{2}_{2} \right) \sigma^{2} \mathbf{r}_{\mathbf{x}_{0}}^{T} \left( \mathbf{H}_{\mathcal{K}}^{T} \mathbf{H}_{\mathcal{K}} \right)^{-1} \mathbf{r}_{\mathbf{x}_{0}}, 
\end{equation} 
where $\mathbf{P}_{\mathcal{K}} \in \mathbb{R}^{M \times M}$ and $\mathbf{r}_{\mathbf{x}_{0}} \in \mathbb{R}^{|\mathcal{K}|}$ are defined as in Theorem \ref{thm_bound_asilomar}.
\end{theorem}

\begin{proof}
We use the same derivation as in the proof of Theorem \ref{thm_bound_asilomar}, but with the subspace $\mathcal{U}_{3}$ in \eqref{equ_proof_lower_bound_asilomar_def_subspace} replaced by 
the subspace 
\begin{equation} 
\mathcal{U}_{4} \triangleq \linspan \big\{ v_{0} (\cdot) \cup \{ v_{l}(\cdot) \}_{l \in \mathcal{K}}Ê\big\}
\end{equation}
which is spanned by the functions $v_{0}(\cdot) \triangleq R_{\mathcal{M}_{\text{SLM}}}(\cdot,\mathbf{x}_{0})$ and 
$v_{l}(\cdot) \triangleq \frac{\partial^{\mathbf{e}_{l}}R_{\mathcal{M}_{\text{SLM}}}(\cdot, \mathbf{x}_{2})}{\partial \mathbf{x}_{2}^{\mathbf{e}_{l}}}Ê\big|_{\mathbf{x}_{2} = \widetilde{\mathbf{x}}_{0}}$.
The difference between $\mathcal{U}_{3}$ and $\mathcal{U}_{4}$ is the definition of the function $v_{0}(\cdot)$. Therefore, only 
the equations in the proof of Theorem \ref{thm_bound_asilomar} which involve $v_{0}(\cdot)$ change. 
However, we still have that $v_{l}(\cdot)$ with $l \in \mathcal{K}$ is orthogonal to $v_{0}(\cdot)$ (see \eqref{equ_proof_lower_bound_asilomar_orthog_vecs_1}), since 
\begin{align} 
\label{equ_proof_lower_bound_asilomar_2_orthog_vecs_1}
\big\langle v_{0}(\cdot) , v_{l}(\cdot) \big\rangle_{\mathcal{H}(\mathcal{M}_{\text{SLM}})} & =
 \bigg\langle R_{\mathcal{M}_{\text{SLM}}}(\cdot,\mathbf{x}_{0}) , \frac{\partial^{\mathbf{e}_{l}}R_{\mathcal{M}_{\text{SLM}}}(\cdot, \mathbf{x}_{2})}{\partial \mathbf{x}_{2}^{\mathbf{e}_{l}}}Ê\bigg|_{\mathbf{x}_{2}Ê  = \widetilde{\mathbf{x}}_{0}} \bigg\rangle_{\mathcal{H}(\mathcal{M}_{\text{SLM}})} \nonumber \\[4mm]
&  \stackrel{\eqref{equ_reproduction_property}}{=} \frac{\partial^{\mathbf{e}_{l}}R_{\mathcal{M}_{\text{SLM}}}(\mathbf{x}_{0}, \mathbf{x}_{2})}{\partial \mathbf{x}_{2}^{\mathbf{e}_{l}}}Ê\bigg|_{\mathbf{x}_{2} = \widetilde{\mathbf{x}}_{0}} \nonumber \\[4mm]
& \stackrel{\eqref{equ_kernel_one_arg_x_0_equal_1}}{=} \frac{\partial^{\mathbf{e}_{l}}1}{\partial \mathbf{x}_{2}^{\mathbf{e}_{l}}}Ê\bigg|_{\mathbf{x}_{2} = \widetilde{\mathbf{x}}_{0}} = 0.
\end{align}
We also obtain the inner products
\begin{align} 
\label{equ_proof_lower_bound_asilomar_2_inner_prod_v_0}
\big\langle v_{0}(\cdot), v_{0}(\cdot) \big\rangle_{\mathcal{H}(\mathcal{M}_{\text{SLM}})}& = \big\langle R_{\mathcal{M}_{\text{SLM}}}(\cdot,\mathbf{x}_{0}), R_{\mathcal{M}_{\text{SLM}}}(\cdot,\mathbf{x}_{0}) \big\rangle_{\mathcal{H}(\mathcal{M}_{\text{SLM}})}
\nonumber \\[4mm]
& \stackrel{\eqref{equ_reproduction_property}}{=}  R_{\mathcal{M}_{\text{SLM}}}(\mathbf{x}_{0},\mathbf{x}_{0}) \stackrel{\eqref{equ_kernel_one_arg_x_0_equal_1}} = 1, 
\end{align} 
and 
\begin{equation} 
\label{equ_proof_lower_bound_asilomar_2_inner_prod_v_0_gamma}
\big\langle  \gamma(\cdot), v_{0}(\cdot) \big\rangle_{\mathcal{H}(\mathcal{M}_{\text{SLM}})} = \big\langle \gamma(\cdot), R_{\mathcal{M}_{\text{SLM}}}(\cdot,\mathbf{x}_{0}) \big\rangle_{\mathcal{H}(\mathcal{M}_{\text{SLM}})}
 \stackrel{\eqref{equ_reproduction_property}}{=}  \gamma(\mathbf{x}_{0}). 
\end{equation} 

Therefore, the relations of \eqref{equ_proof_lower_bound_asilomar_step_by_step} up to the last expression before step $(a)$ are still valid, i.e., 
\begin{align}
\label{equ_proof_lower_bound_asilomar_2_step_by_step}
L_{\mathcal{M}_{\text{SLM}}} & \stackrel{\eqref{equ_squared_norm_min_achiev_var}}{=} \| \gamma(\cdot) \|^{2}_{\mathcal{H}(\mathcal{M}_{\text{SLM}})} - \big[ \gamma(\mathbf{x}_{0}) \big]^{2} \stackrel{\eqref{equ_lower_bound_RKHS_projection}}{\geq}  \| \mathbf{P}_{\mathcal{U}_{4}} \gamma(\cdot) \|^{2}_{\mathcal{H}(\mathcal{M}_{\text{SLM}})}-  \big[ \gamma(\mathbf{x}_{0}) \big]^{2}  \nonumber \\[4mm] 
&  \hspace*{-5mm} \stackrel{\eqref{equ_idendity_projection_subspace_spanned_union_orthogonal_sets},\eqref{equ_proof_lower_bound_asilomar_inner_prod_v_l_gamma}}{=} \big\langle \gamma(\cdot), v_{0}(\cdot) \big\rangle_{\mathcal{H}(\mathcal{M}_{\text{SLM}})} \big( 
\underbrace{\big\langle v_{0}(\cdot), v_{0}(\cdot) \big\rangle_{\mathcal{H}(\mathcal{M}_{\text{SLM}})}}_{\stackrel{\eqref{equ_proof_lower_bound_asilomar_2_inner_prod_v_0}}{=} 1} \big)^{-1} 
\underbrace{\big\langle \gamma(\cdot), v_{0}(\cdot) \big\rangle_{\mathcal{H}(\mathcal{M}_{\text{SLM}})}}_{\stackrel{\eqref{equ_proof_lower_bound_asilomar_2_inner_prod_v_0_gamma}}{=} \gamma(\mathbf{x}_{0})}  \nonumber \\[4mm]
& \hspace*{10mm}+ \mathbf{r}_{\mathbf{x}_{0}}^{T} \mathbf{V}^{\dagger} \mathbf{r}_{\mathbf{x}_{0}} - \big[ \gamma(\mathbf{x}_{0}) \big]^{2} \nonumber \\[4mm]
& \hspace*{-5mm}  = \big[  \gamma(\mathbf{x}_{0}) \big]^{2} +  \mathbf{r}_{\mathbf{x}_{0}}^{T} \mathbf{V}^{\dagger} \mathbf{r}_{\mathbf{x}_{0}} -\big[ \gamma(\mathbf{x}_{0}) \big]^{2} \nonumber \\[4mm]
&\hspace*{-5mm}  \stackrel{\eqref{equ_proof_lower_bound_asilomar_grammian_matrix_complete_notation}}{=} \exp\left( - \frac{1}{\sigma^{2}}\| (\mathbf{I} - \mathbf{P}_{\mathcal{K}}) \mathbf{H} \mathbf{x}_{0} \|^{2}_{2} \right) \sigma^{2} \mathbf{r}_{\mathbf{x}_{0}}^{T} \left( \mathbf{H}_{\mathcal{K}}^{T} \mathbf{H}_{\mathcal{K}} \right)^{-1} \mathbf{r}_{\mathbf{x}_{0}}.
\end{align}
\end{proof}

The difference between the bound \eqref{equ_bound_asilomar_2} and the bound \eqref{equ_bound_asilomar_1} is equal to 
\begin{equation} 
\label{equ_diff_asilomar_bounds}
 \gamma^{2}(\mathbf{x}_{0}) - \exp\left( - \frac{1}{\sigma^{2}}\| (\mathbf{I} - \mathbf{P}_{\mathcal{K}}) \mathbf{H} \mathbf{x}_{0} \|^{2}_{2} \right)  \gamma^{2}(\widetilde{\mathbf{x}}_{0}), 
\end{equation} 
which depends on the choice of the index set $\mathcal{K}$. 
If for some index set $\mathcal{K}$ (note that $\mathbf{r}_{\mathbf{x}_{0}}$ depends on $\mathcal{K}$), the prescribed bias is such that 
$\gamma^{2}(\widetilde{\mathbf{x}}_{0}) Ê\approx \gamma^{2}(\mathbf{x}_{0})$,\footnote{This is the case e.g. if $c(\cdot)Ê\equiv 0$, i.e. for unbiased estimation, and the columns of $\mathbf{H}_{\mathcal{K}}$ are 
nearly orthonormal.}
then the difference \eqref{equ_diff_asilomar_bounds} is approximately nonnegative since 
\begin{equation} 
\label{equ_SLM_diff_nonnegative_asilomar_bounds}
\exp\left( - \frac{1}{\sigma^{2}}\| (\mathbf{I} - \mathbf{P}_{\mathcal{K}}) \mathbf{H} \mathbf{x}_{0} \|^{2}_{2} \right) \leq 1. \nonumber 
\end{equation} 
This in turn, implies that the bound \eqref{equ_bound_asilomar_2} is tighter, i.e. higher, in this case. 


%

\section{The Sparse Signal in Noise Model}
\label{sec_SSNM} 
We now specialize the results obtained for minimum variance estimation associated with the general SLM  $\mathcal{E}_{\text{SLM}}$ 
to the special case of the SLM given by the SSNM $\mathcal{E}_{\text{SSNM}}$ \eqref{equ_def_SSNM_est_problem}. 
The SSNM is obtained from the SLM by choosing $\mathbf{H} = \mathbf{I}$ (which implies that $M=N$). 
We denote by $\mathcal{M}_{\text{SSNM}} \triangleq \left( \mathcal{E}_{\text{SSNM}},c(\cdot),\mathbf{x}_{0} \right)$ any minimum variance problem that is associated with 
$\mathcal{E}_{\text{SSNM}}$, a prescribed bias function $c(\cdot):\mathcal{X}_{S} \rightarrow \mathbb{R}$, and a parameter vector $\mathbf{x}_{0}Ê\in \mathcal{X}_{S}$.
The specific choices for $c(\cdot)$ and $\mathbf{x}_{0}$ should be clear from the context. 
The following discussions and results are presented in part in \cite{RKHSISIT2011}. 

In principle, we can characterize the RKHS $\mathcal{H}(\mathcal{M}_{\text{SSNM}})$, which is associated with the kernel 
\begin{equation} 
\label{equ_def_kernel_SSNM}
R_{\mathcal{M}_{\text{SSNM}}}(\cdot,\cdot): \mathcal{X}_{S} \times \mathcal{X}_{S} \rightarrow \mathbb{R}: R_{\mathcal{M}_{\text{SSNM}}}(\mathbf{x}_{1},\mathbf{x}_{2})\triangleq \exp \left( \frac{1}{\sigma^{2}} (\mathbf{x}_{1} - \mathbf{x}_{0})^{T} (\mathbf{x}_{2} - \mathbf{x}_{0}) \right),
\end{equation} 
using Theorem \ref{thm_isometry_LGM} and Theorem \ref{thm_relation_RKHS_LGM_SLM} specialized to $\mathbf{H}=\mathbf{I}$. 
However, it will turn out that a more convenient characterization is possible by exploiting the specific structure of 
$\mathcal{H}(\mathcal{M}_{\text{SLM}})$ induced by the choice $\mathbf{H} = \mathbf{I}$:

\begin{theorem} 
\label{thm_entire_charact_RKHS_SSNM}
Consider the minimum variance problem $\mathcal{M}_{\emph{SSNM}} = \left( \mathcal{E}_{\emph{SSNM}},c(\cdot),\mathbf{x}_{0} \right)$ with 
some $\sigma$, $S$, $N$, $\mathbf{x}_{0} \in \mathcal{X}_{S}$, and $c(\cdot): \mathcal{X}_{S} \rightarrow \mathbb{R}$. 
The RKHS $\mathcal{H}(\mathcal{M}_{\emph{SSNM}})$ is isometric to the RKHS $\mathcal{H}(R_{e})$ associated with 
the kernel $R_{e}(\cdot,\cdot): \mathcal{X}_{S} \times \mathcal{X}_{S} \rightarrow \mathbb{R}$: 
\begin{equation} 
\label{equ_def_R_e}
R_{e}(\mathbf{x}_{1}, \mathbf{x}_{2}) = \exp \big(\mathbf{x}_{1}^{T} \mathbf{x}_{2} \big). 
\end{equation}
A congruence $\mathsf{K}_{e}[\cdot]: \mathcal{H}(\mathcal{M}_{\emph{SSNM}}) \rightarrow \mathcal{H}(R_{e})$ is given by
\begin{equation} 
\label{equ_def_congruence_SSNM}
f(\cdot) \mapsto \widetilde{f}(\cdot) =  \mathsf{K}_{e}[f(\cdot)] :  \widetilde{f}(\mathbf{x}') 
= f(\sigma \mathbf{x}') \nu_{\mathbf{x}_{0}}(\mathbf{x}')
\end{equation} 
with  $\nu_{\mathbf{x}_{0}}(\mathbf{x}) \triangleq \exp\left( - \frac{1}{2 \sigma^{2}} \|\mathbf{x}_{0} \|^{2}_{2} + Ê\frac{1}{\sigma} \mathbf{x}^{T} \mathbf{x}_{0}\right)$.
The RKHS $\mathcal{H}(R_{e})$ is differentiable up to any order and contains 
the functions 
\begin{equation} 
\label{equ_def_func_g_p_part_der_R_e}
g^{(\mathbf{p})}(\mathbf{x}) \triangleq \frac{1}{\sqrt{\mathbf{p}!}} \frac{ \partial^{\mathbf{p}} R_{e}(\mathbf{x}, \mathbf{x}_{2})}{\partial \mathbf{x}_{2}^{\mathbf{p}}}\bigg|_{\mathbf{x}_{2} = \mathbf{0}} = \frac{1}{\sqrt{\mathbf{p}!}} \mathbf{x}^{\mathbf{p}}, 
\end{equation}
i.e., 
\begin{equation} 
\label{equ_def_SSNM_par_der_R_e}
g^{(\mathbf{p})}(\cdot) \in \mathcal{H}(R_{e}), 
\end{equation}
where $\mathbf{p} \in \mathbb{Z}_{+}^{N}$ is an arbitrary multi-index with $\| \mathbf{p} \|_{0} \leq S$, i.e., $\mathbf{p} \in \mathbb{Z}_{+}^{N} \cap \mathcal{X}_{S}$.
The inner product of an arbitrary function $f(\cdot) \in \mathcal{H}(R_{e})$ with any function $g^{(\mathbf{p})}(\cdot)$ is given by 
\begin{equation} 
\label{equ_inner_prod_part_der_RKHS_Re}
\big\langle f (\cdot), g^{(\mathbf{p})}(\cdot) \big\rangle_{\mathcal{H}(R_{e})} = \frac{1}{\sqrt{\mathbf{p}!}} \frac{ \partial^{\mathbf{p}} f(\mathbf{x})} { \mathbf{x}^{\mathbf{p}}} \bigg|_{\mathbf{x} = 0}. 
\end{equation} 
Moreover, the set $\{g^{(\mathbf{p})}(\cdot) \}_{\mathbf{p} \in \mathbb{Z}_{+}^{N} \cap \mathcal{X}_{S}}$ is an ONB for $\mathcal{H}(R_{e})$ and therefore
any function $f(\cdot) \in \mathcal{H}(R_{e})$ can be written pointwise as 
\begin{equation}
\label{equ_series_repr_RKHS_R_e}
f(\mathbf{x}) = \sum_{\mathbf{p} \in \mathbb{Z}_{+}^{N} \cap \mathcal{X}_{S}} a[\mathbf{p}] g^{(\mathbf{p})}(\mathbf{x}) =   \sum_{\mathbf{p} \in \mathbb{Z}_{+}^{N} \cap \mathcal{X}_{S}} a[\mathbf{p}] \frac{1}{\sqrt{\mathbf{p}!}} \mathbf{x}^{\mathbf{p}},
\end{equation} 
with  a coefficient sequence $a[\mathbf{p}] \in \ell^{2}(\mathbb{Z}_{+}^{N} \cap \mathcal{X}_{S})$. Conversely, for 
any coefficient sequence $a[\mathbf{p}] \in \ell^{2}(\mathbb{Z}_{+}^{N} \cap \mathcal{X}_{S})$ the function 
given by \eqref{equ_series_repr_RKHS_R_e} belongs to $\mathcal{H}(R_{e})$.
\end{theorem}

\begin{proof} 
The derivation of the congruence in \eqref{equ_def_congruence_SSNM} is completely analogous to the proof of Theorem \ref{thm_isometry_LGM} specialized to $\mathbf{H} = \mathbf{I}$: 

Let us define the function sets (similar to \eqref{equ_proof_isometr_R_g_set_A_in_R_LGM} and \eqref{equ_proof_isometr_R_g_set_b_in_R_LGM})
\begin{equation} 
\mathcal{A}_{1} \triangleq \bigg \{ g_{\mathbf{x}} (\cdot) \triangleq R_{\mathcal{M}_{\text{SSNM}}}(\cdot, \sigma \mathbf{x} ) \bigg \}_{\mathbf{x} \in \mathcal{X}_{S}}
 \label{equ_proof_isometr_R_e_set_A_in_R_SSNM}
\end{equation} 
and 
\begin{equation} 
 \label{equ_proof_isometr_R_e_set_B_in_R_SSNM}
\mathcal{B}_{1} \triangleq \bigg \{ f_{\mathbf{x}}(\cdot) \triangleq  
\exp\left( \frac{1}{2 \sigma^{2}} \| \mathbf{x}_{0} \|^{2}_{2} - 
 \frac{1}{\sigma}\mathbf{x}^{T} \mathbf{x}_{0} \right)R_{e}\big(\cdot, \mathbf{x}\big) \bigg \}_{\mathbf{x} \in \mathcal{X}_{S}}.
\end{equation} 
Obviously, we have that the sets $\mathcal{A}_{1}$ and $\mathcal{B}_{1}$ span $\mathcal{H}(\mathcal{M}_{\text{SSNM}})$ and $\mathcal{H}(R_{e})$, respectively, i.e., 
$\linspan \{\mathcal{A}_{1} \} =\mathcal{H}(\mathcal{M}_{\text{SSNM}}) $ and $\linspan \{ \mathcal{B}_{1} \} = \mathcal{H}(R_{e})$. 
Similar to \eqref{equ_proof_LGM_isom_Rg_inner_prod_preserv} one can verify that 
\begin{align}
  \label{equ_proof_SSNM_isom_Re_inner_prod_preserv}
 \big\langle g_{\mathbf{x}_{1}}(\cdot) ,  g_{\mathbf{x}_{2}}(\cdot) \big\rangle_{\mathcal{H}(\mathcal{M}_{\text{SSNM}})} & =
 \big\langle R_{\mathcal{M}_{\text{SSNM}}}(\cdot,\sigma \mathbf{x}_{1}), R_{\mathcal{M}_{\text{SSNM}}}(\cdot,Ê\sigma \mathbf{x}_{2}) \big\rangle_{\mathcal{H}(\mathcal{M}_{\text{SSNM}})} \nonumber  \\[4mm]
& \hspace*{-30mm} \stackrel{\eqref{equ_reproduction_property}}{=} R_{\mathcal{M}_{\text{SSNM}}}(\sigma \mathbf{x}_{1},\sigma \mathbf{x}_{2}) \nonumber \\[4mm]
& \hspace*{-30mm} = \exp   \left( \frac{1}{\sigma^{2}} ( \sigma \mathbf{x}_{1} - \mathbf{x}_{0})^{T}  (\sigma \mathbf{x}_{2} - \mathbf{x}_{0}) \right) \nonumber \\[4mm]Ê
&\hspace*{-30mm} = \exp\left( \frac{1}{2 \sigma^{2}} \|  \mathbf{x}_{0} \|^{2}_{2} -  \frac{1}{\sigma} \mathbf{x}_{1}^{T} \mathbf{x}_{0}+ \mathbf{x}_{1}^{T}  \mathbf{x}_{2} + \frac{1}{2 \sigma^{2}} \|  \mathbf{x}_{0} \|^{2}_{2} - \frac{1}{\sigma} \mathbf{x}_{2}^{T}  \mathbf{x}_{0} \right) \nonumber \\[4mm]
&\hspace*{-30mm} = \exp\left( \frac{1}{2 \sigma^{2}} \| \mathbf{x}_{0} \|^{2}_{2} -  \frac{1}{\sigma} \mathbf{x}_{1}^{T} \mathbf{x}_{0} \right)  \exp   \left(  \mathbf{x}_{1}^{T}  \mathbf{x}_{2}  \right)\exp\left( \frac{1}{2 \sigma^{2}} \|\mathbf{x}_{0} \|^{2}_{2} - \frac{1}{\sigma} \mathbf{x}_{2}^{T}  \mathbf{x}_{0} \right) \nonumber \\[4mm]
& \hspace*{-30mm} \stackrel{\eqref{equ_def_R_e}}{=}   \exp\left( \frac{1}{2 \sigma^{2}} \|  \mathbf{x}_{0} \|^{2}_{2} - \frac{1}{\sigma} \mathbf{x}_{1}^{T} \mathbf{x}_{0} \right)R_{e} \big(\mathbf{x}_{1}, \mathbf{x}_{2}\big) 
\exp\left( \frac{1}{2 \sigma^{2}} \|  \mathbf{x}_{0} \|^{2}_{2} -  \frac{1}{\sigma}\mathbf{x}_{2}^{T} \mathbf{x}_{0} \right)  \nonumber \\[4mm]Ê
& \hspace*{-30mm} \stackrel{\eqref{equ_reproduction_property}}{=}  \bigg\langle \exp\left( \frac{1}{2 \sigma^{2}} \| \mathbf{x}_{0} \|^{2}_{2} - \frac{1}{\sigma} \mathbf{x}_{1}^{T}  \mathbf{x}_{0} \right)R_{e}\big(\cdot, \mathbf{x}_{1}\big), \nonumber \\[4mm]
& \exp\left( \frac{1}{2 \sigma^{2}} \| \mathbf{x}_{0} \|^{2}_{2} -  \frac{1}{\sigma}\mathbf{x}_{2}^{T} \mathbf{x}_{0} \right)R_{e}\big(\cdot, \mathbf{x}_{2}\big) \bigg\rangle_{\mathcal{H}(R_{e})} \nonumber \\[4mm]Ê
&\hspace*{-30mm}  = \big\langle f_{\mathbf{x}_{1}}(\cdot) , f_{\mathbf{x}_{2}}(\cdot) \big\rangle_{\mathcal{H}(R_{e})}.
  \end{align}
for any two vectors $\mathbf{x}_{1}, \mathbf{x}_{2} \in \mathcal{X}_{S}$. 
Similar to \eqref{equ_proof_LGM_isom_Rg_mapping_is_congruence_f_g}, let us, for an arbitrary $\mathbf{x} \in \mathcal{X}_{S}$, denote by $h_{\mathbf{x}}(\cdot) \triangleq \mathsf{K}_{e}[g_{\mathbf{x}}(\cdot)] \in  \mathcal{H}(R_{e})$ the image of the 
function $g_{\mathbf{x}} (\cdot) = R_{\mathcal{M}_{\text{SSNM}}}(\cdot, \sigma \mathbf{x} )$ (see \eqref{equ_proof_isometr_R_e_set_A_in_R_SSNM}) 
under the mapping $\mathsf{K}_{e}[\cdot]$ defined in \eqref{equ_def_congruence_SSNM}. 
We have
\begin{align}
\label{equ_proof_SSNM_isom_Re_mapping_is_congruence_f_g}
h_{\mathbf{x}}(\mathbf{x}') & =  g_{\mathbf{x}}\big(\sigma \mathbf{x}' \big) \exp\left(- \frac{1}{2 \sigma^{2}} \|  \mathbf{x}_{0} \|^{2}_{2} + \frac{1}{\sigma} (\mathbf{x}')^{T}   \mathbf{x}_{0}\right)
 \nonumber \\[4mm]Ê
& \hspace*{-10mm} Ê\stackrel{\eqref{equ_proof_isometr_R_e_set_A_in_R_SSNM}}{=}   \exp\left( -\frac{1}{2 \sigma^{2}} \| \mathbf{x}_{0} \|^{2}_{2} + 
 \frac{1}{\sigma} (\mathbf{x}')^{T}   \mathbf{x}_{0} \right) R_{\mathcal{M}_{\text{SSNM}}}(\sigma \mathbf{x}', \sigma \mathbf{x} ) 
 \nonumber  \\[4mm]
 & \hspace*{-10mm} \stackrel{\eqref{equ_def_kernel_SSNM}}{=}   \exp\left(- \frac{1}{2 \sigma^{2}} \| \mathbf{x}_{0} \|^{2}_{2} + 
 \frac{1}{\sigma} (\mathbf{x}')^{T}   \mathbf{x}_{0} \right) \exp\left( \frac{1}{\sigma^{2}} \big( \sigma \mathbf{x}' - \mathbf{x}_{0}\big)^{T} \big( \sigma \mathbf{x}- \mathbf{x}_{0} \big) \right)
 \nonumber  \\[4mm]
  & \hspace*{-10mm} =  \exp\left(- \frac{1}{ 2\sigma^{2}} \| \mathbf{x}_{0} \|^{2}_{2} + 
 \frac{1}{\sigma} (\mathbf{x}')^{T}   \mathbf{x}_{0} + (\mathbf{x}')^{T} \mathbf{x} +  \frac{1}{ \sigma^{2}} \| \mathbf{x}_{0} \|^{2}_{2} - \frac{1}{\sigma} (\mathbf{x}')^{T}   \mathbf{x}_{0} - \frac{1}{\sigma} \mathbf{x}^{T}   \mathbf{x}_{0}    \right)
 \nonumber  \\[4mm]
   & \hspace*{-10mm} =  \exp\left( \frac{1}{ 2\sigma^{2}} \| \mathbf{x}_{0} \|^{2}_{2} - 
 \frac{1}{\sigma} \mathbf{x}^{T}   \mathbf{x}_{0} \right)   \exp\left((\mathbf{x}')^{T} \mathbf{x}  \right) Ê\nonumber \\[4mm]
& \hspace*{-10mm} \stackrel{\eqref{equ_def_R_e}}{=} \exp\left( \frac{1}{ 2\sigma^{2}} \| \mathbf{x}_{0} \|^{2}_{2} - 
 \frac{1}{\sigma} \mathbf{x}^{T}   \mathbf{x}_{0} \right) R_{e}(\mathbf{x}', \mathbf{x} ) 
 \stackrel{\eqref{equ_proof_isometr_R_e_set_B_in_R_SSNM}}{=} f_{\mathbf{x}}(\mathbf{x}'),
\end{align} 
i.e., 
\begin{align}
 \mathsf{K}_{e}[g_{\mathbf{x}}(\cdot)] = f_{\mathbf{x}}(\cdot).
\end{align} 
for any $\mathbf{x}Ê\in \mathcal{X}_{S}$. 
The fact that the mapping $\mathsf{K}_{e}[\cdot]$ defined in \eqref{equ_def_congruence_SSNM} is a congruence from $\mathcal{H}(\mathcal{M}_{\text{SSNM}})$ to $\mathcal{H}(R_{e})$ follows then from Theorem \ref{thm_suff_con_congruence_RKHS}, since for every argument $\mathbf{x} \in \mathcal{X}_{S}$ the 
function value $\mathsf{K}_{e}[g(\cdot)](\mathbf{x})$ depends continuously on the function value $g\big(\sigma  \mathbf{x} \big) $, which implies that the image $\mathsf{K}_{e}[g(\cdot)]$ of a function $g(\cdot) \in \mathcal{H}(R_{e})$ which is the pointwise limit of a sequence $\big\{ g_{l}(\cdot) \in \linspan \{ \mathcal{A} \} \big\}_{l \rightarrow \infty}$, is the pointwise limit of the functions $\mathsf{K}_{e}[g_{l}(\cdot)] \in \mathcal{H}(\mathcal{M}_{\text{SSNM}})$.

From Theorem \ref{thm_der_repr_prop}, we obtain \eqref{equ_def_SSNM_par_der_R_e} and \eqref{equ_inner_prod_part_der_RKHS_Re}. 

The functions $\{g^{(\mathbf{p})}(\cdot) \}_{\mathbf{p} \in \mathbb{Z}_{+}^{N} \cap \mathcal{X}_{S}}$ are orthogonal since 
for two difference indices $\mathbf{p}_{1}, \mathbf{p}_{2} \mathbb{Z}_{+}^{N}\cap \mathcal{X}_{S}$ there must be at least one index $l \in [N]$ such that $p_{1,l}> p_{2,l}$ or $p_{1,l} < p_{2,l}$. By the symmetry 
of the inner product we can restrict ourselves without loss of generality to the case $p_{1,l}> p_{2,l}$, and compute the inner product
\begin{align} 
\big\langle g^{(\mathbf{p}_{2})}(\cdot) , g^{(\mathbf{p}_{1})}(\cdot) \big\rangle_{\mathcal{H}(R_{e})} 
&= \frac{1}{\sqrt{\mathbf{p}_{1}!}\sqrt{\mathbf{p}_{2}!}} \frac{ \partial^{\mathbf{p}_{1}}  \mathbf{x}^{\mathbf{p}_{2}}} { \mathbf{x}^{\mathbf{p}_{1}}} \bigg|_{\mathbf{x} = 0} \nonumber \\[4mm]
& = \frac{1}{\sqrt{\mathbf{p}_{1}!}\sqrt{\mathbf{p}_{2}!}} \Bigg[ \prod_{l' \in [N] \setminus \{l\}}\frac{ \partial^{{p}_{1,l'}}  \mathbf{x}^{{p}_{2,l'}}} { \mathbf{x}^{{p}_{1,l'}}} \bigg|_{\mathbf{x} = 0}\Bigg] \underbrace{\frac{ \partial^{{p}_{1,l}}  \mathbf{x}^{{p}_{2,l}}} { \mathbf{x}^{{p}_{1,l}}} \bigg|_{\mathbf{x} = 0}}_{=0} \nonumber \\[4mm]
& =0. 
\end{align} 
Since moreover
\begin{align}  
\big\langle g^{(\mathbf{p})}(\cdot) , g^{(\mathbf{p})}(\cdot) \big\rangle_{\mathcal{H}(R_{e})} 
&= \frac{1}{\mathbf{p}!} \frac{ \partial^{\mathbf{p}}  \mathbf{x}^{\mathbf{p}}} { \mathbf{x}^{\mathbf{p}}} \bigg|_{\mathbf{x} = 0}  = \frac{\mathbf{p}!}{\mathbf{p}!} =1, 
\end{align} 
we have that $\{g^{(\mathbf{p})}(\cdot) \}_{\mathbf{p} \in \mathbb{Z}_{+}^{N} \cap \mathcal{X}_{S}}$ is an orthonormal set, i.e., 
\begin{equation} 
\label{equ_proof_isom_R_e_set_g_p_orthonormal}
\big\langle g^{(\mathbf{p}_{2})}(\cdot) , g^{(\mathbf{p}_{1})}(\cdot) \big\rangle_{\mathcal{H}(R_{e})}  = \delta_{\mathbf{p}_{1}, \mathbf{p}_{2}}. 
\end{equation} 
The fact that the set $\{g^{(\mathbf{p})}(\cdot) \}_{\mathbf{p} \in \mathbb{Z}_{+}^{N} \cap \mathcal{X}_{S}}$ is an ONB follows then from Theorem \ref{thm_pointwise_series_kernel_ONB} and the identity 
\begin{align}
R_{e}(\mathbf{x}_{1}, \mathbf{x}_{2}) & \stackrel{\eqref{equ_def_R_e}}{=} \exp \big(\mathbf{x}_{1}^{T} \mathbf{x}_{2} \big) = \prod_{l \in [N]} \exp( x_{1,l} x_{2,l}) = \prod_{l \in [N]}  \sum_{p_{l} \in \mathbb{Z}_{+}} \frac{1}{p_{l}!}(x_{1,l}x_{2,l})^{p_l} \nonumber \\[4mm]
& = \sum_{\mathbf{p} \in \mathbb{Z}_{+}^{N}} \frac{1}{\mathbf{p}!} \mathbf{x}_{1}^{\mathbf{p}}Ê\mathbf{x}_{2}^{\mathbf{p}}Ê\stackrel{(a)}{=} 
 \sum_{\mathbf{p} \in\mathbb{Z}_{+}^{N} \cap \mathcal{X}_{S}} \frac{ \mathbf{x}_{1}^{\mathbf{p}}}{\sqrt{\mathbf{p}!}}\frac{ \mathbf{x}_{2}^{\mathbf{p}}}{\sqrt{\mathbf{p}!}} 
 \stackrel{\eqref{equ_def_func_g_p_part_der_R_e}}{=} \sum_{\mathbf{p} \in\mathbb{Z}_{+}^{N} \cap \mathcal{X}_{S}}g^{(\mathbf{p})}(\mathbf{x}_{1})g^{(\mathbf{p})}(\mathbf{x}_{2}),
\end{align}
which is valid for every $\mathbf{x}_{1}, \mathbf{x}_{2} \in \mathcal{X}_{S}$ (for $\mathbf{x}_{1}, \mathbf{x}_{2} \notin \mathcal{X}_{S}$ the step $(a)$ is not valid).

Finally, the series representation in \eqref{equ_series_repr_RKHS_R_e} follows from Theorem \ref{thm_isometry_hilbert_space_coeffs_space}.  
\end{proof} 

Note that while it might seem at first sight that \eqref{equ_def_congruence_SSNM} 
is just the specialization of \eqref{equ_def_isometry_LGM} to the case $\mathbf{H}=\mathbf{I}$, their meanings are rather different, since 
the domain of $R_{\mathcal{M}_{\text{SSNM}}}(\cdot,\cdot)$ and $R_{e}(\cdot, \cdot)$ is $\mathcal{X}_{S} \times \mathcal{X}_{S}$ 
and not $\mathbb{R}^{N} \times \mathbb{R}^{N}$, which 
is the domain of $R_{\mathcal{M}_{\text{LGM}}}(\cdot,\cdot)$ and $R_{g}^{(N)}(\cdot,\cdot)$ in the special case $\mathbf{H}=\mathbf{I}$. 
Also note, that the kernel $R_{e}(\cdot,\cdot)$ in Theorem \ref{thm_entire_charact_RKHS_SSNM} is the restriction to the sub-domain $\mathcal{X}_{S} \times \mathcal{X}_{S}$ 
of the kernel $R_{g}^{(N)}(\cdot,\cdot): \mathbb{R}^{D} \times \mathbb{R}^{D} \rightarrow \mathbb{R}$ defined in Theorem \ref{thm_isometry_LGM}, i.e., 
$R_{e}(\cdot,\cdot) =  R_{g}^{(N)}(\cdot,\cdot)\big|_{\mathcal{X}_{S} \times \mathcal{X}_{S}}$. From this we have that the RKHSs $\mathcal{H}(R_{e})$ and $\mathcal{H}(R_{g}^{(N)})$ are related via Theorem \ref{thm_reducing_domain_RKHS}. 

Based on Theorem \ref{thm_entire_charact_RKHS_SSNM}, we have the following general characterization of the minimum achievable variance and the corresponding LMV estimator 
for the minimum variance problem $\mathcal{M}_{\text{SSNM}} = \left( \mathcal{E}_{\text{SSNM}},c(\cdot),\mathbf{x}_{0} \right)$:
\begin{theorem} 
\label{thm_general_min_achiev_var_LMV_SSNM}
For the minimum variance problem $\mathcal{M}_{\emph{SSNM}} = \left( \mathcal{E}_{\emph{SSNM}},c(\cdot),\mathbf{x}_{0} \right)$, the prescribed bias function 
$c(\cdot): \mathcal{X}_{S} \rightarrow \mathbb{R}$ is valid if and only if 
there exists a (necessarily unique) coefficient sequence $a[\mathbf{p}] \in \ell^{2}(\mathbb{Z}_{+}^{N} \cap \mathcal{X}_{S})$ such that 
\begin{equation}
\label{equ_condition_gamma_valid_SSNM_fourier_series_R_e}
\gamma(\sigma \mathbf{x}) \nu_{\mathbf{x}_{0}}(\mathbf{x}) = \sum_{\mathbf{p} \in \mathbb{Z}_{+}^{N} \cap \mathcal{X}_{S}} a[\mathbf{p}] \frac{\mathbf{x}^{\mathbf{p}}}{\sqrt{\mathbf{p}!}}
\end{equation}
for every $\mathbf{x} \in \mathcal{X}_{S}$, where $\nu_{\mathbf{x}_{0}}(\mathbf{x})$ as defined in Theorem \ref{thm_entire_charact_RKHS_SSNM} and $\gamma(\mathbf{x}) = c(\mathbf{x}) + x_{k}$ denotes the prescribed mean function. 
If the prescribed bias $c(\cdot)$ is valid for $\mathcal{M}_{\emph{SSNM}}$, then the minimum achievable variance is given by 
\begin{equation}
\label{equ_ssnm_expr_min_ach_var}
L_{\mathcal{M}_{\emph{SSNM}}}  = \sum_{\mathbf{p} \in \mathbb{Z}_{+}^{N} \cap \mathcal{X}_{S}}   \left( a_{\mathbf{x}_{0}}[\mathbf{p}] \right)^{2} - (\gamma(\mathbf{x}_{0}))^{2}
=\|   a_{\mathbf{x}_{0}}[\mathbf{p}] \|^{2}_{\ell^{2}(\mathbb{Z}_{+}^{N} \cap \mathcal{X}_{S})} - (\gamma(\mathbf{x}_{0}))^{2}
\end{equation}
with coefficients 
\begin{equation} 
\label{equ_coeffs_general_min_achiev_var_LMV_SSNM}
a_{\mathbf{x}_{0}}[\mathbf{p}] Ê\triangleq  \big\langle \gamma(\sigma \mathbf{x}) \nu_{\mathbf{x}_{0}}(\mathbf{x}), g^{(\mathbf{p})}(\mathbf{x}) \big\rangle_{\mathcal{H}(R_{e})} = 
\frac{1}{\sqrt{\mathbf{p}!}} \frac{ \partial^{\mathbf{p}} \left( \gamma( \sigma \mathbf{x}) \nu_{\mathbf{x}_{0}}(\mathbf{x}) \right)}{\partial \mathbf{x}^{\mathbf{p}}} \bigg|_{\mathbf{x}=\mathbf{0}}.
\end{equation} 
The corresponding LMV estimator $\hat{x}_{k}^{(\mathbf{x}_{0})}(\cdot)$ is obtained as 
\begin{equation}
\label{equ_ssnm_general_expr_LMV}
\hat{x}_{k}^{(\mathbf{x}_{0})}(\mathbf{y}) = \sum_{\mathbf{p} \in \mathbb{Z}_{+}^{N} \cap \mathcal{X}_{S}} \frac{a_{\mathbf{x}_{0}}[\mathbf{p}]}{\sqrt{\mathbf{p}!}} 
\frac{\partial^{\mathbf{p}} \psi_{\mathbf{x}_{0}}(\mathbf{x}; \mathbf{y})}{ \partial \mathbf{x}^{\mathbf{p}}} \bigg|_{\mathbf{x}= \mathbf{0}}
\end{equation}
where $\psi_{\mathbf{x}_{0}}(\mathbf{x};\mathbf{y}) \triangleq  
\exp\big(Ê\! \frac{1}{\sigma^{2}} \mathbf{y}^{T}( \sigma \mathbf{x} - \mathbf{x}_{0}) + \frac{1}{\sigma} \mathbf{x}_{0}^{T} \mathbf{x} \!-\! \frac{1}{2} \| \mathbf{x} \|^{2}_{2}  \big)$,
i.e., $b(\hat{x}_{k}^{(\mathbf{x}_{0})}(\cdot); \mathbf{x})  = c(\mathbf{x})$ for every $\mathbf{x} \in \mathcal{X}_{S}$ and $v(\hat{x}_{k}^{(\mathbf{x}_{0})}(\cdot); \mathbf{x}_{0}) = L_{\mathcal{M}_{\emph{SSNM}}}$. 
\end{theorem}

\begin{proof}
By Theorem \ref{thm_main_facts_RKHS_MVE}, we have that the prescribed bias function $c(\cdot)$ is valid for $\mathcal{M}_{\text{SSNM}}$ 
if and only if the prescribed mean $\gamma(\cdot)$ belongs to 
$\mathcal{H}(\mathcal{M}_{\text{SSNM}})$. 
By Theorem \ref{thm_entire_charact_RKHS_SSNM}, this is the case if and only if 
\begin{equation} 
\label{equ_cond_valid_gamma_element_RKHS}
\gamma( \sigma \mathbf{x}) \nu_{\mathbf{x}_{0}}(\mathbf{x})  \in \mathcal{H}(R_{e}).
\end{equation} 
which, again by Theorem \ref{thm_entire_charact_RKHS_SSNM}, is the case if and only if there exists a coefficient sequence $a[\mathbf{p}] \in \ell^{2}(\mathbb{Z}_{+}^{N} \cap \mathcal{X}_{S})$ that satisfies \eqref{equ_condition_gamma_valid_SSNM_fourier_series_R_e} for every $\mathbf{x} \in \mathcal{X}_{S}$.

Let us assume from now on that the prescribed bias function $c(\cdot): \mathcal{X}_{S} \rightarrow \mathbb{R}$ is valid for $\mathcal{M}_{\text{SSNM}}$, i.e., 
\eqref{equ_condition_gamma_valid_SSNM_fourier_series_R_e} and \eqref{equ_cond_valid_gamma_element_RKHS} hold. 
By Theorem \ref{thm_main_facts_RKHS_MVE}, Theorem \ref{thm_entire_charact_RKHS_SSNM}, and Theorem \ref{thm_fourier_onb}, we obtain 
\begin{align}
L_{\mathcal{M}_{\text{SSNM}}}  
& \stackrel{\eqref{equ_squared_norm_min_achiev_var}}{=} \|Ê\gamma(\cdot) \|_{\mathcal{H}(\mathcal{M}_{\text{SSNM}})}^{2} -  \big[ \gamma(\mathbf{x}_{0})Ê\big]^{2} \nonumber \\[4mm]
& \stackrel{\eqref{equ_def_congruence_SSNM}}{=}\| \gamma(\sigma \mathbf{x}) \nu_{\mathbf{x}_{0}}(\mathbf{x}) \|_{\mathcal{H}(R_{e})}^{2} - \big[ \gamma(\mathbf{x}_{0})Ê\big]^{2}  \nonumber \\[4mm]
& \stackrel{\eqref{equ_fourier_series_squared_norm}}{=} \sum_{\mathbf{p} \in \mathbb{Z}_{+}^{N} \cap \mathcal{X}_{S}}  \left( \big\langle \gamma(\sigma \mathbf{x}) \nu_{\mathbf{x}_{0}}(\mathbf{x}), g^{(\mathbf{p})}(\mathbf{x}) \big\rangle_{\mathcal{H}(R_{e})} \right)^{2}   -  \big[ \gamma(\mathbf{x}_{0})Ê\big]^{2} \nonumber \\[4mm]
& \stackrel{\eqref{equ_inner_prod_part_der_RKHS_Re}}{=} \sum_{\mathbf{p} \in \mathbb{Z}_{+}^{N} \cap \mathcal{X}_{S}} \frac{1}{\mathbf{p}!}  \left( \frac{ \partial^{\mathbf{p}} ( \gamma(\sigma\mathbf{x}) \nu_{\mathbf{x}_{0}}(\mathbf{x}))}{\partial \mathbf{x}^{\mathbf{p}}} \bigg|_{\mathbf{x}=\mathbf{0}} \right)^{2} - \big[ \gamma(\mathbf{x}_{0})Ê\big]^{2},
\end{align} 
which verifies \eqref{equ_ssnm_expr_min_ach_var}. 
We then have by  
\eqref{equ_condition_gamma_valid_SSNM_fourier_series_R_e} that 
\begin{equation} 
\label{equ_ssnm_expr_gamma_valid}
\gamma(\mathbf{x}) =  \exp\left(  \frac{1}{2 \sigma^{2}} \|\mathbf{x}_{0} \|^{2}_{2} - Ê
\frac{1}{\sigma^2 } \mathbf{x}^{T} \mathbf{x}_{0}\right) \sum_{\mathbf{p} \in \mathbb{Z}_{+}^{N} \cap \mathcal{X}_{S}} \frac{a_{0}[\mathbf{p}]}{\sqrt{\mathbf{p}!}}  \bigg( \frac{\mathbf{x}}{\sigma}\bigg)^{\mathbf{p}},
\end{equation}
with a suitable coefficient sequence $a_{0}[\mathbf{p}] \in \ell^{2}(\mathbb{Z}_{+}^{N} \cap \mathcal{X}_{S})$. 
The image $\mathsf{K}_{e}[\gamma(\cdot)] \in \mathcal{H}(R_{e})$ of the function $\gamma(\cdot)$ given by \eqref{equ_ssnm_expr_gamma_valid} is 
obtained as 
\begin{equation}
\mathsf{K}_{e}[\gamma(\cdot)] = \gamma(\sigma \mathbf{x}) \nu_{\mathbf{x}_{0}}(\mathbf{x}) = \sum_{\mathbf{p} \in \mathbb{Z}_{+}^{N} \cap \mathcal{X}_{S}} \frac{a_{0}[\mathbf{p}]}{\sqrt{\mathbf{p}!}} \mathbf{x}^{\mathbf{p}},
\end{equation}
i.e.,
\begin{equation} 
\label{equ_image_valid_bias_function_R_e}
\mathsf{K}_{e}[\gamma(\cdot)]= \sum_{\mathbf{p} \in \mathbb{Z}_{+}^{N} \cap \mathcal{X}_{S}} \frac{a_{0}[\mathbf{p}]}{\sqrt{\mathbf{p}!}}  \mathbf{x}^{\mathbf{p}} 
\stackrel{\eqref{equ_def_func_g_p_part_der_R_e}}{=} 
\sum_{\mathbf{p} \in \mathbb{Z}_{+}^{N} \cap \mathcal{X}_{S}} a_{0}[\mathbf{p}]g^{(\mathbf{p})}(\cdot).
\end{equation}
Since the functions $\big\{ g^{(\mathbf{p})}(\cdot) \big\}_{\mathbf{p} \in \mathbb{Z}_{+}^{N} \cap \mathcal{X}_{S}}$ are orthonormal (cf.\ \eqref{equ_proof_isom_R_e_set_g_p_orthonormal}), we 
have 
\begin{align}
a_{0}[\mathbf{p}] & \stackrel{\eqref{equ_image_valid_bias_function_R_e}}{=}  \big\langle\mathsf{K}_{e}[\gamma(\cdot)], g^{(\mathbf{p})}(\mathbf{x}) \big\rangle_{\mathcal{H}(R_{e})}  \nonumber \\[4mm]
& \stackrel{\eqref{equ_def_congruence_SSNM}}{=} \big\langle \gamma(\sigma \mathbf{x}) \nu_{\mathbf{x}_{0}}(\mathbf{x}), g^{(\mathbf{p})}(\mathbf{x}) \big\rangle_{\mathcal{H}(R_{e})}  \stackrel{\eqref{equ_inner_prod_part_der_RKHS_Re}}{=} \frac{1}{\sqrt{\mathbf{p}!}} \frac{ \partial^{\mathbf{p}} \left( \gamma( \sigma \mathbf{x}) \nu_{\mathbf{x}_{0}}(\mathbf{x}) \right)}{\partial \mathbf{x}^{\mathbf{p}}} \bigg|_{\mathbf{x}=\mathbf{0}},
\end{align} 
implying also their uniqueness. 

Using 
\begin{align} 
&  \exp \left(  \frac{1}{ \sigma^{2}} \| \mathbf{x}_{0} \|^{2}_{2} - \frac{1}{\sigma^{2}} \mathbf{x}^{T} \mathbf{x}_{0} \right)\bigg( \frac{1}{\sigma} \mathbf{x}\bigg)^{\mathbf{p}} \nonumber \\[4mm]
&  = \exp \left(  \frac{1}{ \sigma^{2}} \|  \mathbf{x}_{0} \|^{2}_{2} - \frac{1}{\sigma^{2}} \mathbf{x}^{T}\mathbf{x}_{0} \right)\frac{ \partial^{\mathbf{p}} \exp \big( \frac{1}{\sigma} (\mathbf{x}')^{T}  \mathbf{x} \big) }{\partial \mathbf{x}'^{\mathbf{p}}}\Bigg|_{\mathbf{x}'=\mathbf{0}} \nonumber \\[4mm]
& = \frac{ \partial^{\mathbf{p}} \exp \big(  \frac{1}{ \sigma^{2}} \|  \mathbf{x}_{0} \|^{2}_{2} - \frac{1}{\sigma^{2}} \mathbf{x}^{T}\mathbf{x}_{0} + \frac{1}{\sigma} \mathbf{x}^{T} \mathbf{x}' \big) }{\partial \mathbf{x}'^{\mathbf{p}}}\Bigg|_{\mathbf{x}'=\mathbf{0}} \nonumber \\[4mm]
& = \frac{ \partial^{\mathbf{p}} \exp \big(  \frac{1}{ \sigma^{2}} (\mathbf{x} - \mathbf{x}_{0})^{T} \big(\sigma \mathbf{x}' -\mathbf{x}_{0} \big) + \frac{1}{\sigma} \mathbf{x}_{0}^{T}  \mathbf{x}' \big) }{\partial \mathbf{x}'^{\mathbf{p}}}\Bigg|_{\mathbf{x}'=\mathbf{0}} \nonumber \\[4mm]
& \stackrel{\eqref{equ_def_kernel_SSNM}}{=} \frac{ \partial^{\mathbf{p}} \big[ R_{\mathcal{M}_{\text{SSNM}}}\big(\mathbf{x},\sigma \mathbf{x}'\big) \exp \big( \frac{1}{\sigma} \mathbf{x}_{0}^{T} \mathbf{x}' \big) \big] }{\partial \mathbf{x}'^{\mathbf{p}}}\Bigg|_{\mathbf{x}'=\mathbf{0}},
\end{align}  
we can further develop \eqref{equ_ssnm_expr_gamma_valid} as 
\begin{align}
\label{equ_ssnm_gamma_sum_part_derivatives_kernel} 
\gamma(\mathbf{x}) & =   \exp\left(  \frac{1}{2 \sigma^{2}} \|\mathbf{x}_{0} \|^{2}_{2} - Ê\frac{1}{\sigma^2 } \mathbf{x}^{T} \mathbf{x}_{0}\right) 
\sum_{\mathbf{p} \in \mathbb{Z}_{+}^{N} \cap \mathcal{X}_{S}} \frac{a_{0}[\mathbf{p}]}{\sqrt{\mathbf{p}!}}  \bigg( \frac{\mathbf{x} }{\sigma} \bigg)^{\mathbf{p}} \nonumber \\[4mm]
& = \exp \bigg( - \frac{1}{2\sigma^{2}} \|  \mathbf{x}_{0} \|^{2}_{2} \bigg)  \sum_{\mathbf{p} \in \mathbb{Z}_{+}^{N}} \frac{a_{0}[\mathbf{p}]}{\sqrt{\mathbf{p}!}}   
\frac{ \partial^{\mathbf{p}} \big[ R_{\mathcal{M}_{\text{SSNM}}}(\mathbf{x}, \sigma \mathbf{x}') \exp \big( \frac{1}{\sigma} \mathbf{x}_{0}^{T} \mathbf{x}' \big) \big]}{\partial \mathbf{x}'^{\mathbf{p}}}\Bigg|_{\mathbf{x}'=\mathbf{0}}.
\end{align} 
This implies via Corollary \ref{cor_der_repr_prop_mult_func}, Theorem \ref{thm_main_facts_RKHS_MVE}, Theorem \ref{thm_isometry_RKHS_rhos_derivative_kernel} and the identity 
\begin{equation}
\label{equ_rho_M_SSNM}
\rho_{\mathcal{M}_{\text{SSNM}}}(\mathbf{y}, \mathbf{x}) = \exp \left( \frac{1}{\sigma^{2}} \mathbf{y}^{T}(\mathbf{x}-\mathbf{x}_{0}) - \frac{1}{2 \sigma^{2}}(\| \mathbf{x} \|^{2}_{2} - \|Ê\mathbf{x}_{0} \|^{2}_{2}) \right), 
\end{equation}
that the estimator defined by \eqref{equ_ssnm_general_expr_LMV} is the 
LMV estimator for $\mathcal{M}_{\text{SSNM}}$: 
\begin{align}
\hat{x}_{k}^{(\mathbf{x}_{0})}(\mathbf{y}) \triangleq \hat{g}^{(\mathbf{x}_{0})}(\mathbf{y}) & \stackrel{\eqref{equ_lmv_estimator_general_congruence_L_M_RKHS}}{=} \mathsf{J}[\gamma(\cdot)]  \nonumber \\[4mm] 
& \stackrel{\eqref{equ_ssnm_gamma_sum_part_derivatives_kernel}}{=} \mathsf{J} \Bigg[  \exp \bigg( - \frac{1}{2\sigma^{2}} \|  \mathbf{x}_{0} \|^{2}_{2} \bigg)  \sum_{\mathbf{p} \in \mathbb{Z}_{+}^{N}} \frac{a_{0}[\mathbf{p}]}{\sqrt{\mathbf{p}!}}   
\frac{ \partial^{\mathbf{p}} \big[ R_{\mathcal{M}_{\text{SSNM}}}(\mathbf{x}, \sigma \mathbf{x}') \exp \big( \frac{1}{\sigma} \mathbf{x}_{0}^{T} \mathbf{x}' \big) \big]}{\partial \mathbf{x}'^{\mathbf{p}}}\Bigg|_{\mathbf{x}'=\mathbf{0}} \Bigg]  \nonumber \\[4mm] 
& =  \exp \bigg( - \frac{1}{2\sigma^{2}} \|  \mathbf{x}_{0} \|^{2}_{2} \bigg)   \sum_{\mathbf{p} \in \mathbb{Z}_{+}^{N}} \frac{a_{0}[\mathbf{p}]}{\sqrt{\mathbf{p}!}}   
\mathsf{J} \Bigg[ \frac{ \partial^{\mathbf{p}} \big[ R_{\mathcal{M}_{\text{SSNM}}}(\mathbf{x}, \sigma \mathbf{x}') \exp \big( \frac{1}{\sigma} \mathbf{x}_{0}^{T} \mathbf{x}' \big) \big]}{\partial \mathbf{x}'^{\mathbf{p}}}\Bigg|_{\mathbf{x}'=\mathbf{0}} \Bigg]  
\nonumber \\[4mm] 
& \stackrel{\eqref{equ_isometry_RKHS_rhos_derivative_kernel}}{=}  \exp \bigg( - \frac{1}{2\sigma^{2}} \|  \mathbf{x}_{0} \|^{2}_{2} \bigg)   \sum_{\mathbf{p} \in \mathbb{Z}_{+}^{N}} \frac{a_{0}[\mathbf{p}]}{\sqrt{\mathbf{p}!}}   
 \frac{ \partial^{\mathbf{p}} \big[ \rho_{\mathcal{M}_{\text{SSNM}}}(\cdot, \sigma \mathbf{x}') \exp \big( \frac{1}{\sigma} \mathbf{x}_{0}^{T} \mathbf{x}' \big) \big]}{\partial \mathbf{x}'^{\mathbf{p}}}\Bigg|_{\mathbf{x}'=\mathbf{0}}  \nonumber \\[4mm] 
& \stackrel{\eqref{equ_rho_M_SSNM}}{=}   \sum_{\mathbf{p} \in \mathbb{Z}_{+}^{N}} \frac{a_{0}[\mathbf{p}]}{\sqrt{\mathbf{p}!}}   
 \frac{ \partial^{\mathbf{p}} \big[  \exp \big( \frac{1}{\sigma^{2}} \mathbf{y}^{T}(\sigma \mathbf{x}' - \mathbf{x}_{0})- \frac{1}{2}\|\mathbf{x}'\|^{2}_{2}+ \frac{1}{\sigma} \mathbf{x}_{0}^{T} \mathbf{x}' \big) \big]}{\partial \mathbf{x}'^{\mathbf{p}}}\Bigg|_{\mathbf{x}'=\mathbf{0}}  \nonumber \\[4mm]
 & = \sum_{\mathbf{p} \in \mathbb{Z}_{+}^{N} \cap \mathcal{X}_{S}} \frac{a_{\mathbf{x}_{0}}[\mathbf{p}]}{\sqrt{\mathbf{p}!}} 
\frac{\partial^{\mathbf{p}} \psi_{\mathbf{x}_{0}}(\mathbf{x}'; \mathbf{y})}{ \partial \mathbf{x}'^{\mathbf{p}}} \bigg|_{\mathbf{x}'= \mathbf{0}}.
\end{align} 
\end{proof}

The statement of Theorem \ref{thm_general_min_achiev_var_LMV_SSNM} is stronger than that of Theorem \ref{thm_condition_valid_bias_SLM} 
when specialized to $\mathbf{H}=\mathbf{I}$, because it contains explicit expressions for the minimum achievable variance $L_{\mathcal{M}_{\text{SSNM}}}$ and the corresponding LMV estimator $\hat{x}_{k}^{(\mathbf{x}_{0})}(\mathbf{y})$. This is possible because of the 
isometry between the RKHS $\mathcal{H}(\mathcal{M}_{\text{SSNM}})$ and the RKHS $\mathcal{H}(R_{e})$, i.e., without the intermediate step via the 
LGM and Theorem \ref{thm_reducing_domain_RKHS}, as used (indirectly through Theorem \ref{thm_relation_RKHS_LGM_SLM}) in Theorem \ref{thm_condition_valid_bias_SLM}. 

The expression \eqref{equ_ssnm_expr_min_ach_var} nicely demonstrates the influence of the sparsity constraints on minimum variance estimation for the SSNM. 
Indeed, let us consider the SSNM without any sparsity constraints, i.e., where $S=N$ and the SSNM reduces to the LGM with $\mathbf{H}=\mathbf{I}$. 
For a fixed prescribed bias $c_{0}(\cdot): \mathbb{R}^{N} \rightarrow \mathbb{R}$ and parameter vector $\mathbf{x}_{0} \in \mathcal{X}_{S'}$ with some sparsity degree 
$S' < N$, denote the minimum achievable variance for $\mathcal{M}_{\text{non-sparse}}= \left( \mathcal{E}_{\text{SSNM}},c_{0}(\cdot),\mathbf{x}_{0} \right)$ 
by $L_{\text{non-sparse}}$. Then consider the SSNM with sparsity degree $S' < N$, i.e., where sparsity constraints are present, and denote 
the minimum achievable variance for $\mathcal{M}_{\text{SSNM}}= \left( \mathcal{E}_{\text{SSNM}},c_{0}(\cdot)\big|_{\mathcal{X}_{S'}},\mathbf{x}_{0} \right)=\mathcal{M}_{\text{non-sparse}}\big|_{\mathcal{X}_{S'}}$ by $L_{S'}$. 
Then we have 
\begin{equation}
\label{equ_diff_min_achiev_var_SSNM}
L_{\text{non-sparse}}-L_{S'}  =
 \sum_{\mathbf{p} \in    \mathbb{Z}_{+}^{N} \setminus \left(\mathbb{Z}_{+}^{N} \cap \mathcal{X}_{S'} \right)} \frac{1}{\mathbf{p}!}  
 \left( \frac{ \partial^{\mathbf{p}} \big[ (c_{0}(\sigma \mathbf{x}) + \sigma x_{k})\nu_{\mathbf{x}_{0}}(\mathbf{x})\big]}{\partial \mathbf{x}^{\mathbf{p}}} \bigg|_{\mathbf{x}=\mathbf{0}} \right)^{2}.
\end{equation}
By reducing the sparsity degree $S'$, the set $ \left(\mathbb{Z}_{+}^{N} \cap \mathcal{X}_{S'} \right)$ becomes smaller. 
This in turn implies that number of summands in \eqref{equ_diff_min_achiev_var_SSNM} 
is increased and therefore, since the summands in \eqref{equ_diff_min_achiev_var_SSNM} are non-negative, the minimum achievable variance $L_{S'}$ becomes smaller. 

Based on Theorem \ref{thm_general_min_achiev_var_LMV_SSNM}, we have the following obvious result concerning the existence of a UMV estimator for the SSNM:  
\begin{lemma} 
Consider the SSNM $\mathcal{E}_{\emph{SSNM}}=\left(\mathcal{X}_{S},f_{\mathbf{H}=\mathbf{I}}(\mathbf{y};\mathbf{x}),g(\mathbf{x}) = x_{k} \right)$ 
and a prescribed bias function $c(\cdot): \mathcal{X}_{S} \rightarrow \mathbb{R}$. Then there exists a UMV estimator for $\mathcal{E}_{\emph{SSNM}}$ and $c(\cdot)$ if 
and only if the LMV $\hat{x}_{k}^{(\mathbf{x}_{0})}(\mathbf{y})$ for the minimum variance problem $\mathcal{M}_{\emph{SSNM}} = \left(\mathcal{E}_{\emph{SSNM}},c(\cdot), \mathbf{x}_{0} \right)$ given 
explicitly by \eqref{equ_ssnm_general_expr_LMV} does not depend on $\mathbf{x}_{0}$. 
\end{lemma}

In what follows, we will need the $l$th order (probabilists') Hermite polynomial $\mathsf{H}_{l}(\cdot): \mathbb{R} \rightarrow \mathbb{R}$ defined as \cite{AbramowitzStegun}
\begin{equation}
\label{equ_def_hermite_polynomial}
\mathsf{H}_{l}(x) = (-1)^l e^{x^2/2}\frac{d^l}{dx^l}e^{-x^2/2}. 
\end{equation} 

We now specialize Theorem \ref{thm_general_min_achiev_var_LMV_SSNM} to a specific class of bias functions $c(\cdot): \mathcal{X}_{S} \rightarrow \mathbb{R}$ that we call ``diagonal'' bias functions. 
A diagonal bias function $c(\mathbf{x})$ depends on the parameter $\mathbf{x}$ only through the $k$th component (where $k$ is the same index as used for the parameter 
function of the SSNM $\mathcal{E}_{\text{SSNM}} = \big( \mathcal{X}_{S}, f_{\mathbf{H}=\mathbf{I}}(\mathbf{y}; \mathbf{x}), g(\mathbf{x}) = x_{k} \big)$), i.e., $c(\mathbf{x}) = \tilde{c}(x_{k})$ with some function $\tilde{c}(\cdot) : \mathbb{R} \rightarrow \mathbb{R}$ that may vary with the index $k$ used for the definition 
of the parameter function of $\mathcal{E}_{\text{SSNM}}$. 
Similarly, we say that an estimator $\hat{x}_{k}(\mathbf{y})$ for the SSNM is 
diagonal if it depends on the observation $\mathbf{y}$ only via the $k$th entry, i.e., $\hat{x}_{k}(\mathbf{y}) = \hat{x}_{k}(y_{k})$. 
Obviously, the bias function of a diagonal estimator for the SSNM is diagonal. 
The well-known hard- and soft-thresholding estimators are diagonal estimators. However, the maximum likelihood (ML) estimator for the SSNM is not diagonal. 

Making the weak assumption that $\tilde{c}(x_{k}): \mathbb{R} \rightarrow \mathbb{R}$ can be represented by a power-series which converges everywhere on $\mathbb{R}$, we obtain 
\begin{theorem} 
\label{thm_diag_bias_min_achiev_var_LMV}
Consider the minimum variance problem $\mathcal{M}_{\emph{SSNM}}=\left(\mathcal{E}_{\emph{SSNM}},c(\cdot),\mathbf{x}_{0}\right)$ for some choice of 
$\sigma$, $S$, $M$, $N$, $\mathbf{x}_{0}$ and for a diagonal prescribed bias function $c(\cdot): \mathcal{X}_{S} \rightarrow \mathbb{R}$, i.e., $c(\mathbf{x}) = \tilde{c}(x_{k})$. 
Moreover, we assume that the function $\tilde{c}(\cdot): \mathbb{R} \rightarrow \mathbb{R}$ is such that the prescribed mean function $\gamma(\cdot): \mathcal{X}_{S} \rightarrow \mathbb{R}: \gamma(\mathbf{x}) = \tilde{c}(x_{k}) + x_{k}$ 
can be written as a power series centered at $\mathbf{x}_{0}$ and which converges everywhere, i.e., 
\begin{equation}
\label{equ_prescr_mean_diag_bias}
\gamma(\mathbf{x}) = \sum\limits_{l \in \mathbb{Z}_{+}} \frac{m_{l}}{l!} (x_{k} - x_{0,k})^{l}
\end{equation}
with suitable coefficients $m_{l} \in \mathbb{R}$. 
We then have that the bias $c(\cdot)$ is valid for $\mathcal{M}_{\emph{SSNM}}$ if and only if 
\begin{equation} 
\label{equ_cond_valid_bias_diagon_bias_SSNM}
\sum\limits_{l \in \mathbb{Z}_{+}} \frac{(m_{l})^{2} \sigma^{2l}}{l!}  < \infty. 
\end{equation}
If the coefficients $m_{l}$ satisfy \eqref{equ_cond_valid_bias_diagon_bias_SSNM}, then we have for the case where $| \supp(\mathbf{x}_{0}) \cup \{k\}| < S+1$ that
\begin{equation}
\label{equ_expr_min_achiev_var_case_1_SSNM}
L_{\mathcal{M}_{\emph{SSNM}}}  = \sum\limits_{lÊ\in \mathbb{N}} \frac{(m_{l})^{2}}{l!} \sigma^{2l},
\end{equation} 
with the corresponding LMV estimator 
\begin{equation}
\label{equ_expr_LMV_diag_bias_case_1_SSNM}
\hat{x}_{k}^{(\mathbf{x}_{0})}(\mathbf{y}) = \sum_{l \in \mathbb{Z}_{+}} \frac{ m_{l}}{l!} \sigma^{l} \mathsf{H}_{l} \left(\frac{ y_k-x_{0,k} }{\sigma} \right), 
\end{equation} 
i.e., $b(\hat{x}_{k}^{(\mathbf{x}_{0})}(\cdot);\mathbf{x})  = \gamma(\mathbf{x}) - x_{k}=\tilde{c}(x_{k})$ for every $\mathbf{x} \in \mathcal{X}_{S}$ and 
$v(\hat{x}_{k}^{(\mathbf{x}_{0})}(\cdot); \mathbf{x}_{0}) = L_{\mathcal{M}_{\emph{SSNM}}}$. 
For the complementary case where $| \supp(\mathbf{x}_{0}) \cup \{k\}| = S+1$, we have, by denoting the indices of the support of $\mathbf{x}_{0}$ as 
$\supp(\mathbf{x}_{0})=\big\{ i_{1},\ldots,i_{S}\big\}$, that
\begin{equation}
\label{equ_min_achiev_var_case_2_ssnm_main}
L_{\mathcal{M}_{\emph{SSNM}}}  = \left(\sum_{l \in \mathbb{Z}_{+}} \frac{m_{l}^{2} \sigma^{2l}}{l!} \right) \Bigg[ \sum_{j \in [S]} \exp\bigg( - \frac{x_{0,i_{j}}^{2}}{\sigma^{2}} \bigg) 
\prod_{j' \in [j-1]}  \bigg[Ê\! 1 - \exp\bigg( - \frac{x_{0,i_{j'}}^{2}}{\sigma^{2}} \bigg) \! \bigg]Ê\Bigg]- \big[\gamma(\mathbf{x}_{0})Ê\big]^{2},
\end{equation} 
with the corresponding LMV estimator 
\begin{align}
\label{equ_lmv_case_2_ssnm_main}
\hat{x}_{k}^{(\mathbf{x}_{0})}(\mathbf{y}) & = \nonumber \\ 
& \hspace*{-10mm} \! \Bigg[ \sum_{l\in \mathbb{Z}_{+}} \frac{ m_{l}}{l!} \sigma^{l} \mathsf{H}_{l} \left( \frac{y_k}{\sigma} \right)\! \Bigg] \sum_{j \in [S]} \exp\bigg( -  \frac{x_{0,i_{j}}^{2}+2y_{i_{j}}x_{0,i_{j}}}{2 \sigma^{2}} \! \bigg)
 \prod_{j' \in [j-1]}  \bigg[Ê\! 1  -  \exp\bigg( - \frac{x_{0,i_{j'}}^{2}+2y_{i_{j'}}x_{0,i_{j'}}}{2 \sigma^{2}} \! \bigg) \! \bigg].
\end{align} 
\end{theorem} 
\begin{proof}
Appendix \ref{chap_appendix_A}. 
\end{proof}

A slight reformulation of Theorem \ref{thm_diag_bias_min_achiev_var_LMV} yields 
\begin{corollary} 
\label{cor_thm_min_var_lmv_ssnm_diag_bias_ISIT}
Consider the minimum variance problem $\mathcal{M}_{\emph{SSNM}}=\left(\mathcal{E}_{\emph{SSNM}},c(\cdot),\mathbf{x}_{0}\right)$ for some choice of 
$\sigma$, $S$, $N$, $\mathbf{x}_{0}\in \mathcal{X}_{S}$, and for a diagonal prescribed bias function $c(\cdot): \mathcal{X}_{S} \rightarrow \mathbb{R}$, i.e., $c(\mathbf{x}) = \tilde{c}(x_{k})$. 
Moreover we assume that the bias function $c(\cdot)$ is valid and such that the prescribed mean function $\gamma(\cdot): \mathcal{X}_{S} \rightarrow \mathbb{R}: \gamma(\mathbf{x}) = \tilde{c}(x_k) + x_{k}$ can be written 
as a power series centered at $x_{0,k}$ and which converges everywhere, i.e., we have $\gamma(\mathbf{x}) = \sum\limits_{l \in \mathbb{Z}_{+}} \frac{m_{l}}{l!} (x_{k} - x_{0,k})^{l}$ with coefficients $m_{l}$ 
satisfying \eqref{equ_cond_valid_bias_diagon_bias_SSNM}. 
Then we obtain for the minimum achievable variance 
\begin{equation}
\label{equ_expr_min_achiev_var_ISIT_SSNM}
L_{\mathcal{M}_{\emph{SSNM}}}  = t_{\mathcal{M}_{\emph{SSNM}}}(\mathbf{x}_{0})   \sum_{l \in \mathbb{Z}_{+}} \frac{m_{l}^{2} \sigma^{2l}}{l!} - \big[\gamma(\mathbf{x}_{0})Ê\big]^{2},
\end{equation} 
and for the corresponding LMV estimator 
\begin{equation}
\label{equ_expr_LMV_diag_bias_ISIT_SSNM}
\hat{x}_{k}^{(\mathbf{x}_{0})}(\mathbf{y})  = h_{\mathcal{M}_{\emph{SSNM}}}(\mathbf{y}, \mathbf{x}_{0}) \sum_{l\in \mathbb{Z}_{+}} \frac{ m_{l}}{l!} \sigma^{l} \mathsf{H}_{l} \left( \frac{y_k - x_{0,k}}{\sigma} \right),
\end{equation} 
where 
\begin{align} 
\label{equ_factor_ISIT_SSNM_min_achiev_var}
 t_{\mathcal{M}_{\emph{SSNM}}}(\mathbf{x}) & \triangleqÊ\nonumber \\ 
 & \hspace*{-20mm} \begin{cases} \sum_{j \in [S]} \exp\bigg( - \frac{x_{i_{j}}^{2}}{\sigma^{2}} \bigg) \prod_{j' \in [j-1]}  \bigg[Ê\! 1 - \exp\bigg( - \frac{x_{i_{j'}}^{2}}{\sigma^{2}} \bigg) \bigg] \leq 1, & \mbox{when }| \supp(\mathbf{x}_{0}) \cup \{k\}| = S+1,\\Ê
1, & \mbox{when } | \supp(\mathbf{x}_{0}) \cup \{k\}| < S+1\end{cases}
 \\[4mm]
\label{equ_factor_ISIT_SSNM_LMV}
h_{\mathcal{M}_{\emph{SSNM}}}(\mathbf{y}, \mathbf{x}) & \triangleq \nonumber \\Ê
& \hspace*{-25mm} \begin{cases} \sum\limits_{j \in [S]} \exp\bigg( \!-  \frac{x_{i_{j}}^{2}+2y_{i_{j}}x_{i_{j}}}{2 \sigma^{2}} \! \bigg)   \prod\limits_{j' \in [j-1]}  \bigg[Ê\! 1  -  \exp\bigg( \!- \frac{x_{i_{j'}}^{2}+2y_{i_{j'}}x_{i_{j'}}}{2 \sigma^{2}} \! \bigg) \bigg] , & \mbox{when }| \supp(\mathbf{x}_{0}) \cup \{k\}| = S+1,\\
1, & \mbox{when } | \supp(\mathbf{x}_{0}) \cup \{k\}| < S+1\end{cases}
\end{align}
with an arbitrary index set $\mathcal{I} = \{ i_{1}, \ldots, i_{S} \}$ such that $|\mathcal{I}| = S$, $k \notin \mathcal{I} $ and $\supp(\mathbf{x}_{0}) \setminus \{k\} \subseteq \mathcal{I}$. 
\end{corollary}

\begin{proof}
The relation \eqref{equ_expr_min_achiev_var_ISIT_SSNM} follows from \eqref{equ_expr_min_achiev_var_case_1_SSNM} and \eqref{equ_min_achiev_var_case_2_ssnm_main} since $\gamma(\mathbf{x}_{0}) = m_{0}$, which implies that 
\begin{equation} 
\sum_{l \in \mathbb{Z}_{+}} \frac{m_{l}^{2} \sigma^{2l}}{l!} - \big[\gamma(\mathbf{x}_{0})Ê\big]^{2} =  \sum_{l \in \mathbb{N}} \frac{m_{l}^{2} \sigma^{2l}}{l!}.
\end{equation}  
The relation \eqref{equ_expr_LMV_diag_bias_ISIT_SSNM} follows from \eqref{equ_expr_LMV_diag_bias_case_1_SSNM} and \eqref{equ_lmv_case_2_ssnm_main}, since for the case $| \supp(\mathbf{x}_{0}) \cup \{k\}| = S+1$ we must have $k \notin \supp(\mathbf{x}_{0})$, i.e., $x_{0,k} = 0$. 
Indeed, if $k \in \supp(\mathbf{x})$ we would have $\supp(\mathbf{x}_{0}) \cup \{k\} = \supp(\mathbf{x}_{0})$, and since the cardinality of the support $\supp(\mathbf{x}_{0})$ cannot be larger than $S$ ($\mathbf{x}_{0} \in \mathcal{X}_{S}$), 
this would mean that $| \supp(\mathbf{x}_{0}) \cup \{k\}| =| \supp(\mathbf{x}_{0})| \leq S <S+1$.

The inequality $t_{\mathcal{M}_{\text{SSNM}}}(\mathbf{x}) \leq 1$ in \eqref{equ_factor_ISIT_SSNM_min_achiev_var} can be verified by observing that
\begin{align} 
& \sum_{j \in [S]} \exp\bigg( - \frac{x_{i_{j}}^{2}}{\sigma^{2}} \bigg) \prod_{j' \in [j-1]}  \bigg(Ê\! 1 - \exp\bigg( - \frac{x_{i_{j'}}^{2}}{\sigma^{2}} \bigg) \bigg)  \nonumber \\[4mm]Ê
& = \Bigg[\prod_{j'' \in [S]} \exp\bigg( - \frac{x_{i_{j''}}^{2}}{\sigma^{2}} \bigg)Ê\Bigg] \sum_{j \in [S]}  \Bigg[\prod_{j''' \in [S]} \exp\bigg(  \frac{x_{i_{j'''}}^{2}}{\sigma^{2}} \bigg)Ê\Bigg]  \exp\bigg( - \frac{x_{i_{j}}^{2}}{\sigma^{2}} \bigg) \Bigg[\prod_{j' \in [j-1]}  \bigg(Ê\! 1 - \exp\bigg( - \frac{x_{i_{j'}}^{2}}{\sigma^{2}} \bigg) \bigg) \Bigg]  \nonumber \\[4mm]Ê
& =\Bigg[\prod_{j'' \in [S]} \exp\bigg( - \frac{x_{i_{j''}}^{2}}{\sigma^{2}} \bigg)Ê\Bigg] \sum_{j \in [S]}  \Bigg[\prod_{j''' \in [S]\setminus [j]} \exp\bigg(  \frac{x_{i_{j'''}}^{2}}{\sigma^{2}} \bigg)Ê\Bigg]   \Bigg[\prod_{j' \in [j-1]}  \bigg(Ê\!  \exp\bigg(  \frac{x_{i_{j'}}^{2}}{\sigma^{2}}-1  \bigg) \bigg) \Bigg]  \nonumber \\[4mm]Ê
&   \stackrel{(a)}{=} \Bigg[\prod_{j'' \in [S]} \exp\bigg( - \frac{x_{i_{j''}}^{2}}{\sigma^{2}} \bigg)Ê\Bigg] \sum_{j \in [S]}
\sum_{ \substack{\mathbf{k} \in \mathbb{Z}_{+}^{S} \\ k_j = 0 \\ k_{j'''} >0\mbox{, } j''' \in [S] \setminus [j]}} \prod_{j' \in [S]} \bigg( \frac{x_{i_{j'}}^{2}}{\sigma^{2}} \bigg)^{k_{j'}}   \nonumber \\[4mm]
&   \leq \Bigg[\prod_{j'' \in [S]} \exp\bigg( - \frac{x_{i_{j''}}^{2}}{\sigma^{2}} \bigg)Ê\Bigg] 
\sum_{ \mathbf{k} \in \mathbb{Z}_{+}^{S}} \prod_{j' \in [S]} \bigg( \frac{x_{i_{j'}}^{2}}{\sigma^{2}} \bigg)^{k_{j'}}   \nonumber \\[4mm]
& \stackrel{(b)}{=} \Bigg[\prod_{j'' \in [S]} \exp\bigg( - \frac{x_{i_{j''}}^{2}}{\sigma^{2}} \bigg)Ê\Bigg] \Bigg[ \prod_{j \in [S]}  \exp\bigg( \frac{x_{i_{j}}^{2}}{\sigma^{2}} \bigg) \Bigg]   =1, \nonumberÊ
\end{align}  
where the step $(a)$ and $(b)$ can be verified by the distributive law and the series representation $\exp\bigg( \frac{x_{i_{j}}^{2}}{\sigma^{2}} \bigg) = \sum_{l \in \mathbb{Z}_{+}} \frac{1}{l!} \bigg( \frac{x_{i_{j}}^{2}}{\sigma^{2}} \bigg)^{l}$. 
\end{proof}

Note that the coefficients $m_{l}$ appearing in the representation \eqref{equ_prescr_mean_diag_bias} of the prescribed mean function 
depend on the parameter vector $\mathbf{x}_{0} \in \mathcal{X}_{S}$ associated with the minimum variance problem $\mathcal{M}_{\text{SSNM}}$, because 
the power series representing $\gamma(\mathbf{x})$ is centered at $\mathbf{x}_{0}$. 

If the prescribed bias function $c(\cdot)$ is the actual bias function $b(\hat{x}'_{k}(\cdot); \mathbf{x})$ of a given diagonal estimator $\hat{x}'_{k}(\mathbf{y}) = \hat{x}'_{k}(y_{k})$ with finite variance at $\mathbf{x}_{0}$, i.e., $\hat{x}'_{k}(\cdot) \in \mathcal{P}(\mathcal{M}_{\text{SSNM}})$ (see \eqref{equ_def_est_finite_power}), the coefficients $m_{l}$ appearing in Corollary \ref{cor_thm_min_var_lmv_ssnm_diag_bias_ISIT} and Theorem \ref{thm_diag_bias_min_achiev_var_LMV} 
have a particular interpretation. 

In what follows, we need the following result:
\begin{lemma}
\label{lem_def_diag_finite_var_est_hilbert_space}
Consider the minimum variance problem $\mathcal{M}_{\emph{SSNM}}(\mathbf{x})=\left(\mathcal{E}_{\emph{SSNM}},c(\cdot),\mathbf{x}\right)$ and the associated function Hilbert space $\mathcal{P}(\mathcal{M}_{\emph{SSNM}})$ of 
finite-variance estimators with inner product 
\begin{align} 
\big\langle \hat{g}_{1}(\cdot), \hat{g}_{2}(\cdot) \big\rangle_{\emph{RV}} & \triangleq  \mathsf{E}_{\mathbf{x}} \big\{  \hat{g}_{1}(\mathbf{y})\hat{g}_{2}(\mathbf{y}) \big\} \nonumber \\[4mm]
& =
 \frac{1}{(2 \pi \sigma^{2})^{N/2}} \int_{\mathbf{y} \in \mathbb{R}^{N}} \hat{g}_{1}(\mathbf{y}) \hat{g}_{2}(\mathbf{y}) \exp \left( - \frac{1}{2 \sigma^{2}} \| \mathbf{y} - \mathbf{x} \|^{2}_{2} \right) dÊ\mathbf{y}. 
\end{align}
Then the subset $\mathcal{D}(\mathbf{x}) \subseteq \mathcal{P}(\mathcal{M}_{\emph{SSNM}})$ which consists 
of all diagonal estimators $\hat{g}(\mathbf{y}) = \hat{g}(y_{k})$ is a Hilbert subspace with induced inner product 
\begin{align} 
\label{equ_def_inner_prod_D_SSNM}
\big\langle \hat{g}_{1}(\cdot), \hat{g}_{2}(\cdot) \big\rangle_{\mathcal{D}(\mathbf{x})}& \triangleq \mathsf{E}_{\mathbf{x}_{0}} \big\{  \hat{g}_{1}(y_{k}) \hat{g}_{2}(y_{k}) \big\} \nonumber \\[4mm]
& = \frac{1}{ \sqrt{2 \pi} \sigma} \int_{y_{k} \in \mathbb{R}} \hat{g}_{1}(y_{k}) \hat{g}_{2}(y_{k}) \exp \left( - \frac{1}{2 \sigma^{2}} (y_k - x_{k} )^{2} \right) dÊy_{k}. 
\end{align}
A specific ONB for $\mathcal{D}(\mathbf{x})$ is given by $\left\{ h^{(l)}(\cdot) \right\}_{l \in \mathbb{Z}_{+}}$ with functions $h^{(l)}(\cdot): \mathbb{R}^{N} \rightarrow \mathbb{R}$ defined as 
\begin{equation}
\label{equ_def_ONB_h_l_D_M_SSNM}
h^{(l)}(\mathbf{y}) \triangleq \frac{1}{\sqrt{l!}} \mathsf{H}_{l}\left(\frac{y_{k} - x_{k}}{\sigma}\right).
\end{equation}
\end{lemma}
\begin{proof}
\cite{SzegoOrthogPoly}
\end{proof}

We have the following straightforward consequence of Corollary \ref{cor_thm_min_var_lmv_ssnm_diag_bias_ISIT}:
\begin{corollary} 
\label{cor_diag_est_min_achiev_var_LMV_ISIT}
Consider the minimum variance problem $\mathcal{M}_{\emph{SSNM}}=\left(\mathcal{E}_{\emph{SSNM}},c(\cdot),\mathbf{x}_{0}\right)$ for some choice of 
$\sigma$, $S$, $N$, $\mathbf{x}_{0} \in \mathcal{X}_{S}$, and assume that the prescribed bias function 
$c(\cdot): \mathcal{X}_{S} \rightarrow \mathbb{R}$ is given as the actual 
bias function of a diagonal estimator $\hat{x}_{k}(\mathbf{y}) = \hat{x}_{k}(y_{k})$, i.e., $c(\mathbf{x}) = b(\hat{x}_{k}(\cdot); \mathbf{x})$. The estimator $\hat{x}_{k}(\mathbf{y})$ is 
assumed to have finite variance and stochastic power at all parameter vectors $\mathbf{x} \in \mathcal{X}_{S}$, i.e., $\hat{x}_{k}(\cdot) \in \mathcal{D}(\mathbf{x})$ for every $\mathbf{x} \in \mathcal{X}_{S}$. 
Then, the prescribed mean function $\gamma(\cdot): \mathcal{X}_{S} \rightarrow \mathbb{R}: \gamma(\mathbf{x}) = c(\mathbf{x}) + x_{k}= \mathsf{E}_{\mathbf{x}} \{ \hat{x}_{k}(\mathbf{y}) \}$ can be written 
as a power series \eqref{equ_prescr_mean_diag_bias} which converges everywhere. 
The coefficients $m_{l}$ in \eqref{equ_prescr_mean_diag_bias} are given by scaled expansion coefficients of $\hat{x}_{k}(\mathbf{y})$ w.r.t.\ the ONB $\left\{ h^{(l)}(\cdot) \right\}_{l \in \mathbb{Z}_{+}}$ defined in \eqref{equ_def_ONB_h_l_D_M_SSNM}, i.e., 
\begin{align} 
m_{l}& = \frac{\sqrt{l!}}{\sigma^{l}} \big\langle \hat{x}_{k}(\cdot), h^{(l)}(\cdot) \big\rangle_{\mathcal{D}(\mathcal{M}_{\emph{SSNM}})}.
\end{align} 
The prescribed bias $c(\cdot)$ is valid by construction. The minimum achievable variance and corresponding LMV estimator are obtained as 
\begin{align}
\label{equ_expr_min_achiev_var_diag_est_ISIT_SSNM}
L_{\mathcal{M}_{\emph{SSNM}}}  & = v(\hat{x}_{k}(\cdot); \mathbf{x}_{0}) t_{\mathcal{M}_{\emph{SSNM}}}(\mathbf{x}_{0}) + \big[ t_{\mathcal{M}_{\emph{SSNM}}}(\mathbf{x}_{0})  - 1 \big]\big[\gamma(\mathbf{x}_{0})Ê\big]^{2},
 \\[4mm]
\hat{x}_{k}^{(\mathbf{x}_{0})}(\mathbf{y})  & = \hat{x}_{k}(y_{k}) h_{\mathcal{M}_{\emph{SSNM}}}(\mathbf{y}, \mathbf{x}_{0}), \label{equ_expr_LMV_diag_est_ISIT_SSNM}
\end{align} 
respectively, with the functions $t_{\mathcal{M}_{\emph{SSNM}}}(\mathbf{x}) $, $h_{\mathcal{M}_{\emph{SSNM}}}(\mathbf{y}, \mathbf{x})$ as defined in \eqref{equ_factor_ISIT_SSNM_min_achiev_var}, and \eqref{equ_factor_ISIT_SSNM_LMV}, respectively.
\end{corollary} 

\begin{proof} 
Since the functions $\left\{ h^{(l)}(\cdot) \right\}_{l \in \mathbb{Z}_{+}}$ form an ONB for $\mathcal{D}(\mathbf{x})$ for every $\mathbf{x} \in \mathcal{X}_{S}$, we have by Theorem \ref{thm_fourier_onb}: 
\begin{align}
\label{equ_bessel_series_power_diag_est_SSNM_ISIT}
\sum_{l \in \mathbb{Z}_{+}} \bigg[ \big\langle \hat{x}_{k}(\cdot), h^{(l)}(\cdot) \big\rangle_{\mathcal{D}(\mathbf{x}))} \bigg]^{2} \stackrel{\eqref{equ_fourier_series_squared_norm}}{=} \| \hat{x}_{k}(\cdot) \|^{2}_{\mathcal{D}(\mathbf{x})} & \stackrel{\eqref{equ_def_inner_prod_D_SSNM}}{=}   P(\hat{x}_{k}(\cdot); \mathbf{x}) Ê = v(\hat{x}_{k}(\cdot); \mathbf{x}) + \big[ \mathsf{E}_{\mathbf{x}} \{ \hat{x}_{k}(\mathbf{y}) \} \big]^{2} \nonumber \\[4mm]
& = v(\hat{x}_{k}(\cdot); \mathbf{x}) + \big[ \gamma(\mathbf{x}) \big]^{2}.
\end{align} 
and
\begin{equation} 
\label{equ_fourier_series_power_diag_est_SSNM_ISIT}
\sum_{l \in \mathbb{Z}_{+}} \big\langle \hat{x}_{k}(\cdot), h^{(l)}(\cdot) \big\rangle_{\mathcal{D}(\mathbf{x})}h^{(l)}(\cdot) = \hat{x}_{k}(\cdot). 
\end{equation} 
%

The prescribed mean function can be expressed as 
\begin{align} 
\label{equ_expr_gamma_diagonal_est_mean_function_ISIT}
\gamma(\mathbf{x}) & = \mathsf{E}_{\mathbf{x}} \{ \hat{x}_{k}(\mathbf{y}) \} = \frac{1}{\sqrt{2 \pi} \sigma}  \int_{y_{k} \in \mathbb{R}} \hat{x}_{k}(y_{k}) \exp \left( - \frac{1}{2 \sigma^{2}} (y_k - x_{k} )^{2} \right) d y_{k}. 
\end{align}
Note that $\gamma(\mathbf{x})$ depends only on $x_{k}$, i.e., $\gamma(\mathbf{x}) = \gamma(x_{k})$. 
Differentiating the mean function $l$ times w.r.t.\ $x_{k}$ yields 
\begin{align} 
\label{equ_part_der_inner_prod_diag_est_SSNM_ISIT}
\frac{\partial^{l \mathbf{e}_{k}} \gamma(\mathbf{x})}{\partial \mathbf{x}^{l \mathbf{e}_{k}}} 
& \stackrel{\eqref{equ_expr_gamma_diagonal_est_mean_function_ISIT}}{=}  \frac{1}{\sqrt{2 \pi} \sigma} \frac{\partial^{l}  \int_{y_{k} \in \mathbb{R}} \hat{x}_{k}(y_{k}) \exp \left( - \frac{1}{2 \sigma^{2}} (y_k - x_{k} )^{2} \right)d y_{k} }{\partial x_{k}^{l}}  \nonumber \\[4mm]
& \stackrel{(a)}{=}  \frac{1}{\sqrt{2 \pi} \sigma}  \int_{y_{k} \in \mathbb{R}} \hat{x}_{k}(y_{k}) \frac{\partial^{l} \exp \left( - \frac{1}{2 \sigma^{2}} (y_k - x_{k} )^{2} \right)}{\partial x_{k}^{l}} d y_{k} \nonumber \\[4mm]Ê
& \stackrel{(b)}{=}  \frac{1}{\sqrt{2 \pi} \sigma^{l+1}}  \int_{y_{k} \in \mathbb{R}} \hat{x}_{k}(y_{k})  \mathsf{H}_{l}\left( \frac{y_{k}- x_{k}}{\sigma} \right) \exp \left( - \frac{1}{2 \sigma^{2}} (y_k - x_{k} )^{2} \right) d y_{k}Ê\nonumber \\[4mm]
& =   \frac{\sqrt{l!}}{\sigma^{l}} \big\langle \hat{x}_{k}(\cdot), h^{(l)}(\cdot) \big\rangle_{\mathcal{D}(\mathbf{x})},
\end{align} 
where $(a)$ follows from changing the order of differentiation and integration (cf.\ \cite[Theorem 1.5.8]{LC}), which be verified, e.g., by a standard argument using the dominated convergence theorem \cite{FundmentExpFamBrown,RudinBookPrinciplesMatheAnalysis,HalmosMeasure}. The step 
$(b)$ can be verified by a calculation similar to \eqref{equ_proof_main_thm_diag_bias_SSNM_cases_identiy_case_1_hermite}. 

Since (by the assumption $\hat{x}_{k}(\cdot) \in \mathcal{D}(\mathbf{x})$ for every $\mathbf{x} \in \mathcal{X}_{S}$) it holds that $P(\hat{x}_{k}(\cdot); \mathbf{x}) < \infty$ for every $\mathbf{x} \in \mathcal{X}_{S}$, we have from \eqref{equ_bessel_series_power_diag_est_SSNM_ISIT} hat the inner products $\big\langle \hat{x}_{k}(\cdot), h^{(l)}(\cdot) \big\rangle_{\mathcal{D}(\mathbf{x})}$ must be bounded by some constant $C$ (strictly speaking, the sequence of the inner products has 
to decay faster than $1/l$). This in turn yields by \eqref{equ_part_der_inner_prod_diag_est_SSNM_ISIT}, that the $l$th order partial derivative 
of the mean function $\gamma(\mathbf{x})$ is bounded by $C\frac{\sqrt{l!}}{\sigma^{l}}$, which implies by Taylor's theorem \cite{RudinBookPrinciplesMatheAnalysis} that the mean function can be written as a power series in the form \eqref{equ_prescr_mean_diag_bias} with coefficients 
\begin{equation} 
\label{equ_coeffs_power_series_inner_prod_hermit} 
m_{l} =  \frac{\sqrt{l!}}{\sigma^{l}} \big\langle \hat{x}_{k}(\cdot), h^{(l)}(\cdot) \big\rangle_{\mathcal{D}(\mathbf{x}_{0})}. 
\end{equation}
The power series is centered at $\mathbf{x}_{0}$ and converges everywhere on $\mathcal{X}_{S}$. 
 
The relation \eqref{equ_expr_min_achiev_var_diag_est_ISIT_SSNM} follows then from \eqref{equ_expr_min_achiev_var_ISIT_SSNM} via \eqref{equ_bessel_series_power_diag_est_SSNM_ISIT} and \eqref{equ_coeffs_power_series_inner_prod_hermit} as 
\begin{align}
\label{equ_expr_min_achiev_var_ISIT_SSNM_mean_func_diag_est}
L_{\mathcal{M}_{\text{SSNM}}} & \stackrel{\eqref{equ_expr_min_achiev_var_ISIT_SSNM}}{=} t_{\mathcal{M}_{\text{SSNM}}}(\mathbf{x}_{0})   \sum_{l \in \mathbb{Z}_{+}} \frac{m_{l}^{2} \sigma^{2l}}{l!} - \big[\gamma(\mathbf{x}_{0})Ê\big]^{2} \nonumber \\[4mm] 
& \stackrel{\eqref{equ_coeffs_power_series_inner_prod_hermit}}{=}  t_{\mathcal{M}_{\text{SSNM}}}(\mathbf{x}_{0})   \sum_{l \in \mathbb{Z}_{+}}  \bigg[\frac{\sqrt{l!}}{\sigma^{l}} \big\langle \hat{x}_{k}(\cdot), h^{(l)}(\cdot) \big\rangle_{\mathcal{D}(\mathbf{x}_{0})}\bigg]^{2} \frac{ \sigma^{2l}}{l!} - \big[\gamma(\mathbf{x}_{0})Ê\big]^{2} \nonumber \\[4mm]
& =  t_{\mathcal{M}_{\text{SSNM}}}(\mathbf{x}_{0})   \sum_{l \in \mathbb{Z}_{+}}  \big[ \big\langle \hat{x}_{k}(\cdot), h^{(l)}(\cdot) \big\rangle_{\mathcal{D}(\mathbf{x}_{0})}\big]^{2}
 - \big[\gamma(\mathbf{x}_{0})Ê\big]^{2} \nonumber \\[4mm]Ê
 & \stackrel{\eqref{equ_bessel_series_power_diag_est_SSNM_ISIT}}{=}  t_{\mathcal{M}_{\text{SSNM}}}(\mathbf{x}_{0})   v(\hat{x}_{k}(\cdot); \mathbf{x}_{0}) +   t_{\mathcal{M}_{\text{SSNM}}}(\mathbf{x}_{0})   \big[ \gamma(\mathbf{x}_{0}) \big]^{2}
 - \big[\gamma(\mathbf{x}_{0})Ê\big]^{2} \nonumber \\[4mm]Ê
  & \stackrel{\eqref{equ_bessel_series_power_diag_est_SSNM_ISIT}}{=} v(\hat{x}_{k}(\cdot); \mathbf{x}_{0}) t_{\mathcal{M}_{\text{SSNM}}}(\mathbf{x}_{0}) + \big[ t_{\mathcal{M}_{\text{SSNM}}}(\mathbf{x}_{0})  - 1 \big]\big[\gamma(\mathbf{x}_{0})Ê\big]^{2}.
  \end{align} 

Finally, the relation \eqref{equ_expr_LMV_diag_est_ISIT_SSNM} follows from \eqref{equ_expr_LMV_diag_bias_ISIT_SSNM} by \eqref{equ_fourier_series_power_diag_est_SSNM_ISIT} as  
\begin{align}
 \hat{x}_{k}^{(\mathbf{x}_{0})}(\mathbf{y})  & \stackrel{\eqref{equ_expr_LMV_diag_bias_ISIT_SSNM}}{=} h_{\mathcal{M}_{\text{SSNM}}}(\mathbf{y}, \mathbf{x}_{0}) \sum_{l\in \mathbb{Z}_{+}} \frac{ m_{l}}{l!} \sigma^{l} \mathsf{H}_{l} \left( \frac{y_k - x_{0,k}}{\sigma} \right) \nonumber \\[4mm]
& \stackrel{\eqref{equ_def_ONB_h_l_D_M_SSNM}}{=} h_{\mathcal{M}_{\text{SSNM}}}(\mathbf{y}, \mathbf{x}_{0}) \sum_{l\in \mathbb{Z}_{+}} \frac{ m_{l}}{\sqrt{l!}} \sigma^{l}h^{(l)}(\mathbf{y})  \nonumber \\[4mm]
& \stackrel{\eqref{equ_coeffs_power_series_inner_prod_hermit}}{=} h_{\mathcal{M}_{\text{SSNM}}}(\mathbf{y}, \mathbf{x}_{0}) \sum_{l\in \mathbb{Z}_{+}}  \big\langle \hat{x}_{k}(\cdot), h^{(l)}(\cdot) \big\rangle_{\mathcal{D}(\mathbf{x}_{0})}  h^{(l)}(\mathbf{y})  \nonumber \\[4mm] 
& \stackrel{\eqref{equ_fourier_series_power_diag_est_SSNM_ISIT}}{=} h_{\mathcal{M}_{\text{SSNM}}}(\mathbf{y}, \mathbf{x}_{0}) \hat{x}_{k}(\mathbf{y}).  
\end{align} 
\end{proof} 

Note that Corollary \ref{cor_diag_est_min_achiev_var_LMV_ISIT} also applies to unbiased estimation where $c(\cdot) \equiv 0$ and $\gamma(\mathbf{x}) = x_{k}$, since these are the bias and mean functions of the 
ordinary LS estimator $\hat{x}_{\text{LS},k}(\mathbf{y}) = y_{k}$, which obviously is diagonal and has finite variance at every $\mathbf{x} \in \mathcal{X}_{S}$. 

Remarkably, according to Corollary \ref{cor_diag_est_min_achiev_var_LMV_ISIT} (in particuar, due to 
\eqref{equ_factor_ISIT_SSNM_min_achiev_var} and \eqref{equ_factor_ISIT_SSNM_LMV}), every diagonal estimator $\hat{x}_{k}(\cdot): \mathbb{R}^{N} \rightarrow \mathbb{R}$ for the SSNM with finite variance at $\mathbf{x}_{0} \in \mathcal{X}_{S}$ is 
necessarily the LMV estimator for the minimum variance problem $\mathcal{M}_{\text{SSNM}}=\left(\mathcal{E}_{\text{SSNM}},c(\mathbf{x})=b(\hat{x}_{k}(\cdot); \mathbf{x}),\mathbf{x}_{0}\right)$ whenever $|\supp(\mathbf{x}_{0}) \cup \{ k \}| < S +1$, 
irrespective of the sparsity degree $S$.
Thus, in this case, the sparsity information does not help anything with regard to minimum variance estimation, i.e., the estimator $\hat{x}_{k}(\cdot)$ is the LMV estimator for any parameter set $\mathcal{X}_{S}$ including the non-sparse case $\mathcal{X} = \mathbb{R}^{N}$ 
and the variance $v(\hat{x}_{k}(\cdot);\mathbf{x}_{0})$ coincides with the minimum achievable variance $L_{\mathcal{M}_{\text{SSNM}}}$, i.e., the Barankin bound. 

However, whenever $|\supp(\mathbf{x}_{0}) \cup \{ k \}| =S+1$, we have by Corollary \ref{cor_diag_est_min_achiev_var_LMV_ISIT} and the inequality in \eqref{equ_factor_ISIT_SSNM_min_achiev_var} 
that there are estimators with the same bias and mean as $\hat{x}_{k}(\cdot)$ 
but with a smaller variance at $\mathbf{x}_{0}$. The LMV estimator $\hat{x}_{k}^{(\mathbf{x}_{0})}(\cdot)$ that achieves the minimum variance at $\mathbf{x}_{0}$ is obtained by a multiplication 
of the estimator $\hat{x}_{k}(\mathbf{y})=\hat{x}_{k}(y_{k})$ with a ``correction'' factor $h_{\mathcal{M}_{\text{SSNM}}}(\mathbf{y}, \mathbf{x}_{0})$ that 
does not depend on $y_{k}$. 
Note that by a straightforward calculation one can show that the LMV estimator in \eqref{equ_expr_LMV_diag_est_ISIT_SSNM} 
is robust against deviations from the nominal parameter $\mathbf{x}_{0}$. In particular, that estimator has always finite mean 
and variance for any parameter vector $\mathbf{x} \in \mathbb{R}^{N}$ and an observation $\mathbf{y} = \mathbf{x} + \mathbf{n}$ that follows the statistical model of the LGM with $\mathbf{H} = \mathbf{I}$ and 
an arbitrary noise variance $\sigma^{2}$. 
We verify this claim separately for the two complementary cases $|\supp(\mathbf{x}_{0}) \cup \{ k \}| =S+1$ and $|\supp(\mathbf{x}_{0}) \cup \{ k \}|  < S+1$. First consider 
the case $|\supp(\mathbf{x}_{0}) \cup \{ k \}| =S+1$, where, as already observed above, necessarily 
\begin{equation} 
\label{equ_k_notin_supp_x_0}
\{k\} \notin \supp(\mathbf{x}_{0})=\{i_{1},\ldots,i_{S}\},
\end{equation}
i.e., $x_{0,k} = 0$ and moreover $|\supp(\mathbf{x}_{0})| = \|Ê\mathbf{x}_{0} \|_{0} =S$. 
Given an arbitrary vector $\mathbf{x} \in \mathbb{R}^{N}$, the stochastic power $P(\hat{x}_{k}^{(\mathbf{x}_{0})}(\cdot); \mathbf{x})$ at $\mathbf{x}$ 
of the LMV estimator in \eqref{equ_expr_LMV_diag_est_ISIT_SSNM} is then obtained as 
\begin{align}
P(\hat{x}_{k}^{(\mathbf{x}_{0})}(\cdot); \mathbf{x}) & = \mathsf{E}_{\mathbf{x}} \bigg\{ \bigg[ \hat{x}_{k}^{(\mathbf{x}_{0})}(\cdot) \bigg]^{2} \bigg\} \nonumber \\[4mm]
& \stackrel{\eqref{equ_expr_LMV_diag_est_ISIT_SSNM}}{=}  \mathsf{E}_{\mathbf{x}} \bigg\{ \bigg[  \hat{x}_{k}(y_{k}) h_{\mathcal{M}_{\text{SSNM}}}(\mathbf{y}, \mathbf{x}_{0}) \bigg]^{2} \bigg\}\nonumber \\[4mm]
& \stackrel{\eqref{equ_factor_ISIT_SSNM_LMV}}{=}  \mathsf{E}_{\mathbf{x}} \bigg\{ \bigg[  \hat{x}_{k}(y_{k})  \sum\limits_{j \in [S]}\alpha_{i_{j}}   \prod\limits_{j' \in [j-1]}  (\! 1  -  \alpha_{i_{j'}} ) \bigg]^{2}  \bigg \} \\[4mm]
& \stackrel{\eqref{equ_factor_ISIT_SSNM_LMV}}{=}  \mathsf{E}_{\mathbf{x}} \bigg\{ \bigg[  \hat{x}_{k}(y_{k})  \sum\limits_{j \in [S]}\alpha_{i_{j}}   \prod\limits_{j' \in [j-1]}  (Ê\! 1  -  \alpha_{i_{j'}} )  \bigg]^{2}  \bigg \} \\[4mm]
& \stackrel{\eqref{equ_k_notin_supp_x_0}}{=}  \underbrace{\mathsf{E}_{\mathbf{x}} \bigg\{ \big[  \hat{x}_{k}(y_{k}) \big]^{2} \bigg\}}_{< \infty} \mathsf{E}_{\mathbf{x}}Ê
\bigg\{ \bigg[  \sum\limits_{j \in [S]}\alpha_{i_{j}}   \prod\limits_{j' \in [j-1]}  (Ê\! 1  -  \alpha_{i_{j'}} )  \bigg]^{2}  \bigg \} \nonumber \\[4mm]
& =  \mathsf{E}_{\mathbf{x}} \bigg\{ \big[  \hat{x}_{k}(y_{k}) \big]^{2} \bigg\} \mathsf{E}_{\mathbf{x}}Ê
\bigg\{ \sum\limits_{j,j'' \in [S]} \alpha_{i_{j}}  \alpha_{i_{j''}}  \bigg[  \prod\limits_{j' \in [j-1]}  (Ê\! 1  -  \alpha_{i_{j'}} ) \bigg]
\bigg[\prod\limits_{j''' \in [j''-1]}  (Ê\! 1  -  \alpha_{i_{j'''}} )  \bigg] \bigg \} < \infty
\end{align}
where we introduced the shorthand $\alpha_{i} \triangleq  \exp\bigg( \!- \frac{x_{0,i}^{2}+2y_{i}x_{0,i}}{2 \sigma^{2}} \! \bigg)$ and the last 
step follows by collecting terms, and the fact that 
\begin{equation}
\mathsf{E}_{\mathbf{x}} \bigg\{ \prod_{i \in \supp(\mathbf{x}_{0})} \big[\alpha_{i}\big]^{n_{i}} \bigg\}  = \prod_{i \in \supp(\mathbf{x}_{0})}\mathsf{E}_{\mathbf{x}} \big\{ \big[\alpha_{i}\big]^{n_{i}} \big\} < \infty, 
\end{equation}
for any nonnegative exponents $n_{i} \in \mathbb{Z}_{+}$. 
Similarly, one can also show that the mean of $\hat{x}_{k}^{(\mathbf{x}_{0})}(\cdot)$ is finite at every $\mathbf{x} \in \mathbb{R}^{N}$, i.e., 
$\mathsf{E}_{\mathbf{x}} \bigg\{ \hat{x}_{k}^{(\mathbf{x}_{0})}(\cdot)  \bigg\}< \infty$.

\section{Necessity of Strict Sparsity} 
\label{sec_strict_sparstiy_SLM}

In what follows, we consider the LGM $\mathcal{E}_{\text{LGM}}$ with a fixed choice of $\sigma$, $M$, $N$, $\mathbf{H}$, and a fixed choice of 
the sparsity degree $S$, and the associated minimum variance problem $\mathcal{M}_{\text{LGM}}= \left( \mathcal{E}_{\text{LGM}},\tilde{c}(\cdot), \mathbf{x}_{0} \right)$ 
with parameter vector $\mathbf{x}_{0} \in \mathcal{X}_{S}$ and prescribed bias $\tilde{c}(\cdot): \mathbb{R}^{N} \rightarrow \mathbb{R}$. 
Since the SLM $\mathcal{E}_{\text{SLM}}$ for $S <N$ is obtained by reducing the parameter set $\mathcal{X} = \mathbb{R}^{N}$ of the LGM $\mathcal{E}_{\text{LGM}}$ to the 
set of $S$-sparse vectors $\mathcal{X}_{S}$, the set of allowed estimators $\mathcal{F}(\mathcal{M}_{\text{SLM}})$ (see \eqref{equ_est_finite_var_prescr_bias}) 
for the minimum variance problem $\mathcal{M}_{\text{SLM}} = \mathcal{M}_{\text{LGM}}\big|_{\mathcal{X}_{S}}$ is in general larger than 
the set of allowed estimators $\mathcal{F}(\mathcal{M}_{\text{LGM}})$ for $\mathcal{M}_{\text{LGM}}$, i.e., 
$\mathcal{F}(\mathcal{M}_{\text{SLM}}) \supseteq \mathcal{F}(\mathcal{M}_{\text{LGM}})$. 

It is this increase of the set of allowed estimators that causes the minimum achievable variance $L_{\mathcal{M}_{\text{SLM}}}$ (for $S < N$) to be strictly smaller than $L_{\mathcal{M}_{\text{LGM}}}$ in general. This fact agrees with the discussion of Section \ref{sec_reducing_par_set_RKHS} according to which we have $L_{\mathcal{M}_{\text{SLM}}} \leq L_{\mathcal{M}_{\text{LGM}}}$. 
An important exception to this rule is the special case of an affine prescribed bias, i.e., $c(\mathbf{x}) = \mathbf{a}^{T} \mathbf{x} + b$ with a fixed vector $\mathbf{a} \in \mathbb{R}^{N}$ and number $b \in \mathbb{R}$, and a full column rank system matrix $\mathbf{H}$. In this case, the lower bound in Theorem \ref{thm_CRB_SLM} for the case $\| \mathbf{x}_{0} \|_{0} < S$ coincides with $L_{\mathcal{M}_{\text{LGM}}}$. Therefore, the minimum achievable variance $L_{\mathcal{M}_{\text{SLM}}}$ for the SLM coincides with the 
minimum achievable variance $L_{\mathcal{M}_{\text{LGM}}}$ for the LGM (since $L_{\mathcal{M}_{\text{SLM}}} \leq L_{\mathcal{M}_{\text{LGM}}}$).

Note that the reduction from the parameter set $\mathcal{X} = \mathbb{R}^{N}$ of the LGM $\mathcal{E}_{\text{LGM}}$ to the parameter set $\mathcal{X} = \mathcal{X}_{S}$ of the SLM $\mathcal{E}_{\text{SLM}}$ with $S < N$ is rather large, 
since the Lebesgue measure of $\mathcal{X}_{S}$ (with $S < N$) is zero w.r.t.\ $\mathbb{R}^{N}$ \cite{HalmosMeasure}. 

Intuitively, it should make not a big difference if we relax the constraint $\mathbf{x} \in \mathcal{X}_{S}$ and include 
also \emph{approximately} $S$-sparse vectors in the parameter set. Thus, we prescribe the estimator bias not only for strictly but also for approximately $S$-sparse vectors.  
A specific definition of approximate $S$-sparsity is used, e.g., in \cite{RaskuttiMinmaxSLM, DonohoJohnstone94,Donoho94idealspatial}, where the authors model the set of approximately 
$S$-sparse vectors by an $\ell_{q}$ ball of radius $S$, where $0 \leq q \leq 1$. The $\ell_{q}$ ball $\mathcal{B}_{q}(S)$ of radius $S$ is defined by 
\begin{equation}
\label{equ_def_approximate_S_sparse_vecs_q_ball} 
\mathcal{B}_{q}(S)Ê\triangleq \bigg \{ \mathbf{x}' \in \mathbb{R}^{N} \big|  \| \mathbf{x}' \|^{q}_{q} \leq S \bigg\}.
\end{equation} 
Note that the set $\mathcal{X}_{S}$ of strictly sparse vectors  is obtained as the limiting case when $q=0$, i.e., $\mathcal{X}_{S} = \mathcal{B}_{0}(S)$. 
\vspace*{10mm}
\begin{figure}[htbp]
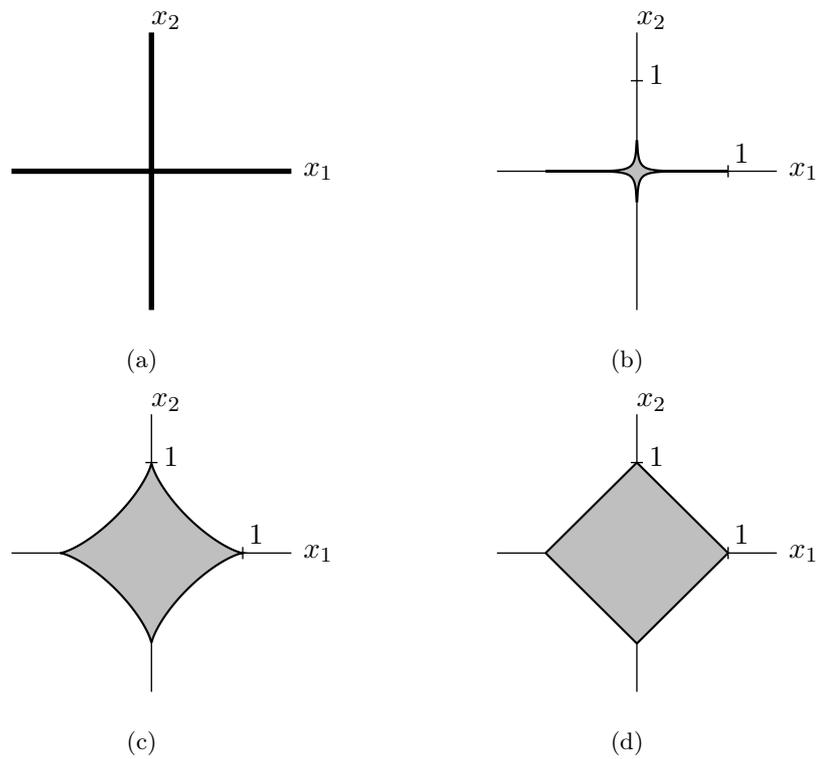

\psset{xunit=0.8cm,yunit=0.8cm}
\begin{center}
\subfigure[]
{\label{fig_B_q_0}
\pspicture(-2.5,-2.5)(2.5,2.5) 
\psline[linewidth=2pt](-2.3,0)(2.3,0)
\psline[linewidth=2pt](0,-2.3)(0,2.3) 
\rput[l](2.5,0){$x_{1}$}
\rput[l](0,2.5){$x_{2}$}
\endpspicture
}
\hspace*{2cm}
\subfigure[]
{\label{fig_B_q_0_25}
\pspicture(-2.5,-2.5)(2.5,2.5) 
\psline[linewidth=0.5pt](-2.3,0)(2.3,0)
\psline[linewidth=0.5pt](0,-2.3)(0,2.3) 
\pscustom[fillcolor=lightgray, fillstyle=solid]{
\psplot[linewidth=0.1pt,plotpoints=300]{-1.5}{1.5} {1.5 0.25 exp x abs 0.25 exp  sub 4 exp}
\psplot[linewidth=0.1pt,plotpoints=300]{1.5}{-1.5} {1.5 0.25 exp x abs 0.25 exp  sub 4 exp -1 mul}
}
\psline[linewidth=0.5pt](1.5,-0.1)(1.5,0.1)
\psline[linewidth=0.5pt](-0.1,1.5)(0.1,1.5)
\rput[l](1.6,0.3){$1$}
\rput[l](0.2,1.6){$1$}
\rput[l](2.5,0){$x_{1}$}
\rput[l](0,2.5){$x_{2}$}
\endpspicture
} \\
\subfigure[]
{
\label{fig_B_q_0_75}
\pspicture(-2.5,-2.5)(2.5,2.5) 
\psline[linewidth=0.5pt](-2.3,0)(2.3,0)
\psline[linewidth=0.5pt](0,-2.3)(0,2.3) 
\pscustom[fillcolor=lightgray, fillstyle=solid]{
\psplot[linewidth=0.1pt,plotpoints=300]{-1.5}{1.5} {1.5 0.75 exp x abs 0.75 exp  sub 1.333 exp}
\psplot[linewidth=0.1pt,plotpoints=300]{1.5}{-1.5} {1.5 0.75 exp x abs 0.75 exp  sub 1.333 exp -1 mul}
}
\psline[linewidth=0.5pt](1.5,-0.1)(1.5,0.1)
\psline[linewidth=0.5pt](-0.1,1.5)(0.1,1.5)
\rput[l](1.6,0.3){$1$}
\rput[l](0.2,1.6){$1$}
\rput[l](2.5,0){$x_{1}$}
\rput[l](0,2.5){$x_{2}$}
\endpspicture
}
\hspace*{2cm}
\subfigure[]
{
\label{fig_B_q_1}
\pspicture(-2.5,-2.5)(2.5,2.5) 
\psline[linewidth=0.5pt](-2.3,0)(2.3,0)
\psline[linewidth=0.5pt](0,-2.3)(0,2.3) 
\pscustom[fillcolor=lightgray, fillstyle=solid]{
\psplot[linewidth=0.1pt,plotpoints=300]{-1.5}{1.5} {1.5 1 exp x abs 1 exp  sub 1 exp}
\psplot[linewidth=0.1pt,plotpoints=300]{1.5}{-1.5} {1.5 1 exp x abs 1 exp  sub 1 exp -1 mul}
}
\psline[linewidth=0.5pt](1.5,-0.1)(1.5,0.1)
\psline[linewidth=0.5pt](-0.1,1.5)(0.1,1.5)
\rput[l](1.6,0.3){$1$}
\rput[l](0.2,1.6){$1$}
\rput[l](2.5,0){$x_{1}$}
\rput[l](0,2.5){$x_{2}$}
\endpspicture
}
\end{center}
\caption{Different parameter sets consisting of strictly or approximately $1$-sparse 
vectors $\mathbf{x} =Ê\big( x_{1}, x_{2} \big)^{T} \in \mathbb{R}^{2}$.  \subref{fig_B_q_0} $ \mathcal{X}= \mathcal{B}_{0}(1) = \mathcal{X}_{1}$. 
\subref{fig_B_q_0_25} $ \mathcal{X} = \mathcal{B}_{0.25}(1)$. \subref{fig_B_q_0_75} $ \mathcal{X} = \mathcal{B}_{0.75}(1)$. \subref{fig_B_q_1} $ \mathcal{X} = \mathcal{B}_{1}(1).$}
\label{fig_approx_sparse_par_set}
\end{figure}
\vspace*{10mm}

In Figure \ref{fig_approx_sparse_par_set}, we illustrate $\mathcal{B}_{q}(S)$ for $N=2$, $S=1$ and various $q$. 
Note also that in contrast to $\mathcal{X}_{S}$, the parameter sets $\mathcal{B}_{q}(S)$ are bounded for $q>0$, i.e., for every $q >0$ and $S$ there exists a radius $r$ such that $\mathcal{B}_{q}(S) \subseteq \mathcal{B}(\mathbf{0}, r) \subseteq \mathbb{R}^{N}$.

Let us now define the minimum variance problem $\mathcal{M}^{(q)} \triangleq \mathcal{M}_{\text{LGM}}\big|_{\mathcal{B}_{q}(S)}$, which is obtained from $\mathcal{M}_{\text{LGM}}$ by reducing 
the parameter set from $\mathcal{X} = \mathbb{R}^{N}$ to $\mathcal{X} =\mathcal{B}_{q}(S)$. Note that $\mathcal{M}_{\text{SLM}} = \mathcal{M}^{(0)}$. 
Then, one might expect that for $q \rightarrow 0$, i.e., when $\mathcal{B}_{q}(S) \rightarrow \mathcal{X}_{S}$, we should have 
$L_{\mathcal{M}^{(q)}} \rightarrow L_{\mathcal{M}_{\text{SLM}}}$. 

Despite the intuition that a change of the parameter set from $\mathcal{X}_{S}$ to $\mathcal{B}_{q}(S)$ with small $q>0$, i.e., considering $\mathcal{M}^{(q)}$ instead of $\mathcal{M}_{\text{SLM}}$, should make not a big difference, the next result 
tells us that for minimum variance estimation it makes a very big difference if we use $\mathcal{X}_{S}$ or $\mathcal{B}_{q}(S)$, no matter how small $q$ is.  
\begin{theorem}
\label{thm_strict_sparsity_nec_SLM} 
Consider a parameter set $\mathcal{X} \subseteq \mathbb{R}^{N}$ and the minimum variance problem $\mathcal{M}'=\mathcal{M}_{\emph{LGM}}\big|_{\mathcal{X}} =
 \left(\mathcal{E}', \tilde{c}(\cdot)\big|_{\mathcal{X}}, \mathbf{x}_{0}\in \mathcal{X}Ê\right)$ associated with the estimation problem $\mathcal{E}'=\left(\mathcal{X},f_{\mathbf{H}} (\mathbf{y}; \mathbf{x}), g(\mathbf{x}) = x_{k} \right)$. Note that the minimum variance problem $\mathcal{M}'$ 
 is identical to $\mathcal{M}_{\emph{LGM}} =\left( \mathcal{E}_{\emph{LGM}},\tilde{c}(\cdot), \mathbf{x}_{0} \right)$ except for the parameter set. 
 We assume that the prescribed bias function $\tilde{c}(\cdot): \mathbb{R}^{N} \rightarrow \mathbb{R}$ is valid for $\mathcal{M}_{\emph{LGM}}$.
%
Then, if the parameter set $\mathcal{X}$ contains an open ball $\mathcal{B}(\mathbf{x}_{c},r)$ with some radius $r>0$ and center $\mathbf{x}_{c} \in \mathcal{X}$, i.e., 
\begin{equation} 
\label{equ_thm_strict_spar_nec_contains_ball}
\mathcal{B}(\mathbf{x}_{c},r) \subseteq \mathcal{X},
\end{equation}
we have that 
\begin{equation}
L_{\mathcal{M}'} = L_{\mathcal{M}_{\emph{LGM}}}.
\end{equation}
\end{theorem} 
\begin{proof} 
We will denote the prescribed mean functions of the minimum variance problems $\mathcal{M}'$ and $\mathcal{M}_{\text{LGM}}$ by $\gamma'(\cdot): \mathcal{X} \rightarrow \mathbb{R}: \gamma'(\mathbf{x}) \triangleq \tilde{c}(\mathbf{x}) + x_{k}$ and $\gamma(\cdot): \mathbb{R}^{N} \rightarrow \mathbb{R}: \gamma(\mathbf{x}) \triangleq \tilde{c}(\mathbf{x}) + x_{k}$, respectively. 
These two mean functions are related by 
\begin{equation}
\label{equ_proof_necessity_mean_func_restriction}
\gamma'(\cdot) = \gamma(\cdot)\big|_{\mathcal{X}}.
\end{equation} 
Note that the assumption that $\tilde{c}(\cdot)$ is valid for $\mathcal{M}_{\text{LGM}}$ implies that the bias function $\tilde{c}(\cdot)\big|_{\mathcal{X}}$ is valid for $\mathcal{M}'$ (due to Theorem \ref{thm_reducing_domain_RKHS}). 
We therefore have that $\gamma'(\cdot) \in \mathcal{H}(\mathcal{M}')$ and $\gamma(\cdot) \in  \mathcal{H}(\mathcal{M}_{\text{LGM}})$ (due to Theorem \ref{thm_main_facts_RKHS_MVE}). 

Since the kernel $R_{\mathcal{M}'}(\cdot,\cdot)$ is the restriction of the kernel $R_{\mathcal{M}_{\text{LGM}}}(\cdot,\cdot)$ to the subdomain $\mathcal{X} \times \mathcal{X}$, i.e., 
$R_{\mathcal{M}'}(\cdot,\cdot) = R_{\mathcal{M}_{\text{LGM}}}(\cdot,\cdot)\big|_{\mathcal{X} \times \mathcal{X}}$, we have by
Theorem \ref{thm_main_facts_RKHS_MVE} and Theorem \ref{thm_reducing_domain_RKHS} that
\begin{align} 
\label{equ_min_var_strict_spar_nec}
L_{\mathcal{M}'}  & \stackrel{\eqref{equ_squared_norm_min_achiev_var}}{=} \| \gamma'(\cdot) \|_{\mathcal{H}(\mathcal{M}')}^{2} - \big[ \gamma'(\mathbf{x}_{0}) \big]^{2}  \nonumber \\[4mm] 
& \stackrel{\eqref{equ_thm_reducing_domain_RKHS}}{=} \min_{\substack{\gamma'' (\cdot) \in \mathcal{H}(\mathcal{M}_{\text{LGM}}) \\ \gamma''(\cdot)\big|_{\mathcal{X}} = \gamma'(\cdot) }}\| \gamma'' (\cdot) \|^{2}_{\mathcal{H}(\mathcal{M}_{\text{LGM}})}- \big[ \gamma'(\mathbf{x}_{0}) \big]^{2} \nonumber \\[4mm]
& \stackrel{\eqref{equ_proof_necessity_mean_func_restriction}}{=} \min_{\substack{\gamma'' (\cdot) \in \mathcal{H}(\mathcal{M}_{\text{LGM}}) \\ \gamma''(\cdot)\big|_{\mathcal{X}} \equiv 0 }}\| \gamma'' (\cdot)+ \gamma(\cdot) \|^{2}_{\mathcal{H}(\mathcal{M}_{\text{LGM}})}- \big[ \gamma'(\mathbf{x}_{0}) \big]^{2}. 
\end{align}
Using Theorem \ref{thm_isometry_LGM}, we can reformulate \eqref{equ_min_var_strict_spar_nec} as 
\begin{equation} 
L_{\mathcal{M}'}  = \min_{\substack{\gamma''(\cdot) \in \mathcal{H}(\mathcal{M}_{\text{LGM}}) \\ \gamma''(\cdot)\big|_{\mathcal{X}} \equiv 0}} \big\| \mathsf{K}^{-1}_{g}[\gamma''(\cdot) + \gamma(\cdot)] \big\|^{2}_{\mathcal{H}(R_{g}^{(D)})}-\big[ \gamma'(\mathbf{x}_{0}) \big]^{2} 
\end{equation}
with $D = \rank(\mathbf{H})$. 

Consider then an arbitrary function $\gamma''(\cdot) \in \mathcal{H}(\mathcal{M}_{\text{LGM}})$ which vanishes on $\mathcal{X}$, i.e., $\gamma''(\mathbf{x}) = 0$ 
for every $\mathbf{x} \in \mathcal{X}$.
Due to \eqref{equ_thm_strict_spar_nec_contains_ball}, the function $\mathsf{K}^{-1}_{g}[\gamma''(\cdot)](\mathbf{x}')$ vanishes for every 
$\mathbf{x}' \in \mathcal{B}(\widetilde{\mathbf{H}} \mathbf{x}_{c},r_{1}) \subseteq \mathbb{R}^{D}$ 
with a suitable radius $r_{1} >0$ and with center $\widetilde{\mathbf{H}} \mathbf{x}_{c}$, where $\widetilde{\mathbf{H}}$ is defined as in Theorem \ref{thm_isometry_LGM}. 
From this, it follows via Theorem \ref{thm_part_der_any_origin_complete_R_g} that for every $\gamma''(\cdot) \in \mathcal{H}(\mathcal{M}_{\text{LGM}})$ that 
satisfies $\gamma''(\cdot)\big|_{\mathcal{X}} \equiv 0$ we have
\begin{equation} 
\mathsf{K}^{-1}_{g}[\gamma''(\cdot)+ \gamma(\cdot)]=\mathsf{K}^{-1}_{g}[\gamma''(\cdot)]+\mathsf{K}^{-1}_{g}[\gamma(\cdot)] =\mathsf{K}^{-1}_{g}[\gamma(\cdot)].
\end{equation} 
This, in turn, implies that 
\begin{align}
L_{\mathcal{M}'}  & = \min_{\substack{\gamma''(\cdot) \in \mathcal{H}(\mathcal{M}_{\text{LGM}}) \\ \gamma''(\cdot)\big|_{\mathcal{X}} \equiv 0}} \big\| \mathsf{K}^{-1}_{g}[\gamma''(\cdot) + \gamma(\cdot)] \big\|^{2}_{\mathcal{H}(R_{g}^{(D)})}-\big[ \gamma'(\mathbf{x}_{0}) \big]^{2}  \nonumber \\[4mm]
 & =  \| \mathsf{K}^{-1}_{g}[\gamma(\cdot)] \|^{2}_{\mathcal{H}(R_{g}^{(D)})} -\big[ \gamma'(\mathbf{x}_{0}) \big]^{2}  \nonumber \\[4mm]
 &  \stackrel{(a)}{=} \|\gamma(\cdot) \|^{2}_{ \mathcal{H}(\mathcal{M}_{\text{LGM}})} -\big[ \gamma(\mathbf{x}_{0}) \big]^{2}  \stackrel{\eqref{equ_squared_norm_min_achiev_var}}{=} L_{\mathcal{M}_{\text{LGM}}},
\end{align}
where for the step $(a)$ we used again Theorem \ref{thm_isometry_LGM}. 
\end{proof}

Since we have obviously that the parameter set $\mathcal{X} = \mathcal{B}_{q}(S)$ satisfies \eqref{equ_thm_strict_spar_nec_contains_ball} for any $q>0$, Theorem \ref{thm_strict_sparsity_nec_SLM} implies that 
$L_{\mathcal{M}^{(q)}} = L_{\mathcal{M}_{\text{LGM}}}$ for any $q>0$. Thus, the minimum achievable variance for $\mathcal{M}^{(q)}$ with $q >0$ is always equal to that of the minimum variance problem $\mathcal{M}_{\text{LGM}}$, which 
does not include any sparsity constraints since its parameter set is given by $\mathcal{X} = \mathbb{R}^{N}$. 
Furthermore, since in general (see \eqref{equ_diff_min_achiev_var_SSNM}) $L_{\mathcal{M}_{\text{LGM}}} >  L_{\mathcal{M}_{\text{SLM}}}$, we have also 
that $L_{\mathcal{M}^{(q)}}$ does not converge to $L_{\mathcal{M}_{\text{SLM}}}$ as $q$ goes to zero in general.

\section{The SLM Viewpoint on Compressed Sensing}
\label{sec_slm_viewpoint_CS}
The compressive measurement process of a compressed sensing (CS) application is often modeled by a linear measurement equation of the form 
\cite{JustRelax,GreedisGood,Can06a,ZvikaCoherenceTSP,MallatBook}
\begin{equation}
\label{equ_measurment_equ_CS}
\mathbf{y} = \mathbf{H} \mathbf{x} + \mathbf{n}, 
\end{equation} 
where $\mathbf{y} \in \mathbb{R}^{M}$ denotes the compressive measurements; $\mathbf{H} \in \mathbb{R}^{M \times N}$ (where typically $M \ll N$) denotes the CS measurement matrix; the signal or parameter vector $\mathbf{x} \in \mathbb{R}^{N}$ is 
assumed to be sparse, i.e., $\mathbf{x} \in \mathcal{X}_{S}$ with known sparsity degree $S$ (where typically $S \ll N$); and $\mathbf{n}$ represents some additive measurement noise. 
If the vector $\mathbf{n} \in \mathbb{R}^{M}$ in \eqref{equ_measurment_equ_CS} is AWGN, i.e., $\mathbf{n} \sim \mathcal{N}(\mathbf{0}, \sigma^{2} \mathbf{I})$, then the CS measurement equation \eqref{equ_measurment_equ_CS} is identical to the observation model of the SLM (see \eqref{equ_linear_observation_model}).  

Any CS recovery method, e.g., the ``Basis Pursuit'' (BP) \cite{JustRelax,Chen98atomicdecomposition} or the ``Orthogonal Matching Pursuit'' (OMP) 
\cite{GreedisGood,TroppGilbertOMP} to name two well-known instances,\footnote{A comprehensive overview is provided at \htmladdnormallink{http://dsp.rice.edu/cs}{http://dsp.rice.edu/cs}.}
can be interpreted as an estimator $\hat{\mathbf{x}}(\mathbf{y})$ that aims to estimate the sparse vector $\mathbf{x}$ using 
the observation $\mathbf{y}$ given by \eqref{equ_measurment_equ_CS}. 

In general, it is infeasible to characterize the analytical properties of the measurement matrix $\mathbf{H}$ in an exact manner (e.g., to compute its thin SVD), due to the large dimension (typically $M \geq 10$ and $NÊ\geq 100$). Instead, one can use an incomplete characterization via the concept of (mutual) coherence and the restricted isometry property \cite{JustRelax,GreedisGood,Can06a,ZvikaCoherenceTSP}, which we will now define.
\begin{definition}
The coherence of a matrix $\mathbf{H} \in \mathbb{R}^{M \times N}$ is defined as 
\begin{equation} 
\mu(\mathbf{H}) \triangleq \max\limits_{i \neq j} | \mathbf{h}_{j}^{T} \mathbf{h}_{i} |, 
\end{equation} 
where $\mathbf{h}_{i} \in \mathbb{R}^{M}$ denotes the $i$th column of $\mathbf{H}$. 
\end{definition} 

\begin{definition} 
A matrix $\mathbf{H} \in \mathbb{R}^{M \times N}$ is said to satisfy the restricted isometry property (RIP) of order $K$ 
with RIP constant $\delta_{K} \in \mathbb{R}_{+}$ if for every index set $\mathcal{I}Ê\subseteq [N]$ of size $K$, i.e., $|\mathcal{I}| = K$, we have 
\begin{equation} 
\label{equ_RIP_condition}
(1-\delta_{K}) \| \mathbf{z} \|^{2}_{2} \leq \| \mathbf{H}_{\mathcal{I}} \mathbf{z} \|^{2}_{2} \leq (1+\delta_{K}) \| \mathbf{z} \|^{2}_{2}
\end{equation}
for every $\mathbf{z} \in \mathbb{R}^{K}$. 
\end{definition}
Note that it can be straightforwardly shown that $\delta_{K} \leq \delta_{K'}$ if $K' \geq K$.

While the coherence of a matrix can be calculated more efficiently than the RIP constants $\delta_{K}$ \cite{ZvikaCoherenceTSP}, it is 
a coarser description of the matrix $\mathbf{H}$ than the RIP constants. However, as stated in \cite{ZvikaCoherenceTSP}, we have the following relation between these two concepts: 
\begin{theorem}
For a matrix $\mathbf{H}$ with coherence $\mu(\mathbf{H})$, the RIP constant $\delta_{K}$ of order $K$ satisfies 
\begin{equation} 
\label{equ_bound_coherence_RIP}
\delta_{K} \leq (K-1) \mu(\mathbf{H}). 
\end{equation} 
\end{theorem} 
\begin{proof}
\cite{ZvikaCoherenceTSP}
\end{proof}

We now specialize Theorem \ref{thm_bound_asilomar_2}  to the CS scenario of the SLM, i.e., where the system matrix $\mathbf{H}$ of the SLM is a CS measurement matrix 
with known RIP constants.   
\begin{theorem}
\label{thm_CS_measur_matrix_asilomar_bound}
Consider the minimum variance problem $\mathcal{M}_{\emph{SLM}}=\left( \mathcal{E}_{\emph{SLM}}, c(\cdot), \mathbf{x}_{0} \right)$ 
with some choice of $\sigma$, $S$, $M$, $N$, $\mathbf{x}_{0} \in \mathcal{X}_{S}$ and with system matrix $\mathbf{H}$, where $\mathbf{H} \in \mathbb{R}^{M \times N}$ is a CS measurement matrix satisfying the RIP of order $S$ with RIP constant  $\delta_{S} < 1$. We denote the prescribed mean by $\gamma(\cdot): \mathcal{X}_{S} \rightarrow \mathbb{R}: \gamma(\mathbf{x}) = c(\mathbf{x}) + x_{k}$.
Then, for an arbitrary set $\mathcal{K} =\{i_1,\ldots,i_{|\mathcal{K}|} \} \subseteq [N]$ consisting of no more than $S$ different indices, i.e., $|\mathcal{K}| \leq S$, 
for which the partial derivatives $\frac{ \partial^{\mathbf{e}_{i_{l}}} \gamma(\mathbf{x})}{\partial \mathbf{x}^{\mathbf{e}_{i_{l}}}}$ exist, we have that
\begin{equation}
\label{equ_CS_measur_matrix_asilomar_bound}
L_{\mathcal{M}_{\emph{SLM}}} \geq  \expÊ\bigg( - \frac{1+\delta_{S}}{\sigma^{2}} \bigg\| \mathbf{x}_{0}^{\supp(\mathbf{x}_{0}) \setminus \mathcal{K}} \bigg \|^{2}_{2}    \bigg)\sigma^{2}  \mathbf{r}_{\mathbf{x}_{0}}^{T} \left( \mathbf{H}_{\mathcal{K}}^{T} 
\mathbf{H}_{\mathcal{K}} \right)^{-1} \mathbf{r}_{\mathbf{x}_{0}}, 
\end{equation} 
where the vector $\mathbf{r}_{\mathbf{x}_{0}} \in \mathbb{R}^{|\mathcal{K}|}$ 
is defined elementwise as $r_{\mathbf{x}_{0},l} \triangleq \frac{ \partial^{\mathbf{e}_{i_{l}}} \gamma(\mathbf{x})}{\partial \mathbf{x}^{\mathbf{e}_{i_{l}}}} \big|_{\mathbf{x} = \widetilde{\mathbf{x}}_{0}}$ with $\widetilde{\mathbf{x}}_{0}$ as defined by \eqref{equ_def_x_0_tilde_asiloamr_bound_1} in Theorem \ref{thm_bound_asilomar}. 
\end{theorem}

\begin{proof}
Let $\mathbf{P}_{\mathcal{K}} \in \mathbb{R}^{M \times |\mathcal{K}|}$ denote the orthogonal projection matrix on the 
subspace $\linspan( \mathbf{H}_{\mathcal{K}} ) \subseteq \mathbb{R}^{M}$. We have 
that $\mathbf{P}_{\mathcal{K}} = \mathbf{H}_{\mathcal{K}} \left(\mathbf{H}_{\mathcal{K}}^{T}\mathbf{H}_{\mathcal{K}}\right)^{-1} \mathbf{H}_{\mathcal{K}}^{T}$ since $\mathbf{H}$ satisfies the 
RIP or order $S$, which implies that $\mathbf{H}_{\mathcal{K}}$ has full column rank \cite{golub96}. 
We also have
\begin{equation} 
(\mathbf{I} - \mathbf{P}_{\mathcal{K}})^{T} \mathbf{H}_{\mathcal{K}} = (\mathbf{I} - \mathbf{P}_{\mathcal{K}}) \mathbf{H}_{\mathcal{K}} = \mathbf{0} \nonumber,
\end{equation} 
i.e., for every vector $\mathbf{x}' \in \linspan(\mathbf{H}_{\mathcal{K}})$ we have that $ (\mathbf{I} - \mathbf{P}_{\mathcal{K}}) \mathbf{x}' = \mathbf{0}$. 
This implies that 
\begin{equation} 
\label{equ_slm_CS_proof_reformulation_factor}
(\mathbf{I} - \mathbf{P}_{\mathcal{K}})^{T} \mathbf{H} \mathbf{x}_{0} = 
(\mathbf{I} - \mathbf{P}_{\mathcal{K}})\mathbf{H} ( \mathbf{x}_{0}^{\supp(\mathbf{x}_{0}) \setminus \mathcal{K}}  + \mathbf{x}_{0}^{\mathcal{K}}) =  
(\mathbf{I} - \mathbf{P}_{\mathcal{K}})\mathbf{H} \mathbf{x}_{0}^{\supp(\mathbf{x}_{0}) \setminus \mathcal{K}},
\end{equation} 
since $\mathbf{H} \mathbf{x}_{0}^{\mathcal{K}} \in  \linspan(\mathbf{H}_{\mathcal{K}})$. 
Based on \eqref{equ_slm_CS_proof_reformulation_factor} and using the shorthand $\mathbf{x}'_{0} \triangleq  \mathbf{x}_{0}^{\supp(\mathbf{x}_{0}) \setminus \mathcal{K}}$, we have 
\begin{align}
\label{equ_slm_CS_inequ_rip}
\big \| (\mathbf{I} - \mathbf{P}_{\mathcal{K}} ) \mathbf{H} \mathbf{x}_{0} \big \|^{2}_{2} & = \big \| (\mathbf{I} - \mathbf{P}_{\mathcal{K}} ) \mathbf{H}\mathbf{x}_{0}' \big \|^{2}_{2}  \nonumber  \\[4mm] 
& =  \big \| \mathbf{H} \mathbf{x}'_{0} \big \|^{2}_{2} - 2 \left( \mathbf{x}'_{0} \right)^{T} \mathbf{H}^{T} \mathbf{P}_{\mathcal{K}} \mathbf{H}\mathbf{x}'_{0} + 
\left(\mathbf{x}'_{0} \right)^{T} \mathbf{H}^{T} \mathbf{P}^{2}_{\mathcal{K}} \mathbf{H}\mathbf{x}'_{0} \nonumber \\[4mm] 
& \stackrel{(a)}{=}  \big \| \mathbf{H}\mathbf{x}'_{0} \big \|^{2}_{2} - \left(\mathbf{x}'_{0}\right)^{T} \mathbf{H}^{T} \mathbf{P}_{\mathcal{K}} \mathbf{H}\mathbf{x}'_{0} \nonumber \\[4mm] 
& \stackrel{(b)}{\leq}  \big \| \mathbf{H}\mathbf{x}'_{0} \big \|^{2}_{2}  \nonumber \\[4mm] 
&\stackrel{(c)}{\leq}  (1+\delta_{S}) \big \| \mathbf{x}'_{0} \big \|^{2}_{2},
\end{align} 
where the step $(a)$ follows from $\left( \mathbf{P}_{\mathcal{K}} \right)^{2} = \mathbf{P}_{\mathcal{K}}$ since $\mathbf{P}_{\mathcal{K}}$ is an orthogonal projection matrix \cite{golub96}, the 
step $(b)$ follows from the fact that $\mathbf{P}_{\mathcal{K}}$ is psd, and $(c)$ is due to $\| \mathbf{x}_{0}' \|_{0} \leq S$ and the assumption that $\mathbf{H}$ 
satisfies the RIP of order $S$ with RIP constant $\delta_{S} < 1$.  
The statement follows then as a reformulation of \eqref{equ_bound_asilomar_2} from Theorem \ref{thm_bound_asilomar_2} using \eqref{equ_slm_CS_inequ_rip}, i.e., 
\begin{align}
L_{\mathcal{M}_{\text{SLM}}} & \stackrel{\eqref{equ_bound_asilomar_2}}{\geq}  \exp\left( - \frac{1}{\sigma^{2}} \| (\mathbf{I} - \mathbf{P}_{\mathcal{K}}) \mathbf{H} \mathbf{x}_{0} \|^{2}_{2} \right) \sigma^{2} \mathbf{r}_{\mathbf{x}_{0}}^{T} \left( \mathbf{H}_{\mathcal{K}}^{T} \mathbf{H}_{\mathcal{K}} \right)^{-1} \mathbf{r}_{\mathbf{x}_{0}} \nonumber \\[4mm]
& \stackrel{\eqref{equ_slm_CS_inequ_rip}}{\geq}   \expÊ\bigg( - \frac{1+\delta_{S}}{\sigma^{2}} \big\| \mathbf{x}_{0}^{\supp(\mathbf{x}_{0}) \setminus \mathcal{K}} \big \|^{2}_{2}    \bigg)\sigma^{2}  \mathbf{r}_{\mathbf{x}_{0}}^{T} \left( \mathbf{H}_{\mathcal{K}}^{T} 
\mathbf{H}_{\mathcal{K}} \right)^{-1} \mathbf{r}_{\mathbf{x}_{0}}.
\end{align}  
\end{proof}

If we want to use Theorem \ref{thm_CS_measur_matrix_asilomar_bound} for comparing the actual variance behavior of a given CS recovery scheme (using the CS measurement matrix $\mathbf{H}$), which is 
an estimator $\hat{x}_{k}(\cdot)$ for the SLM with system matrix $\mathbf{H}$, with the theoretically minimum variance achievable for the bias of the estimator $\hat{x}_{k}(\cdot)$, we 
have to ensure that the first-order partial derivatives of the mean function $\mathsf{E}_{\mathbf{x}} \{ \hat{x}_{k}(\mathbf{y}) \}$ of the given estimator $\hat{x}_{k}(\cdot)$ exist. 
The next lemma states that this is indeed the case for a very broad class of estimators $\hat{x}_{k}(\cdot)$ resulting from a CS recovery scheme. It is 
essentially a more rigorous formulation of the relation $(19)$ in \cite{HeroUniformCRB}.
\begin{lemma}
\label{lem_cond_exist_CS_recovery_partial_der}
Consider the SLM $\mathcal{E}_{\emph{SLM}}$ with system matrix $\mathbf{H} \in \mathbb{R}^{M \times N}$ 
and an estimator $\hat{x}_{k}(\cdot): \mathbb{R}^{M} \rightarrow \mathbb{R}$ which 
may arise from a CS recovery scheme like BP or OMP. We denote the mean function of $\hat{x}_{k}(\cdot)$ by 
$\gamma(\cdot): \mathcal{X}_{S} \rightarrow \mathbb{R}: \gamma(\mathbf{x}) = \mathsf{E}_{\mathbf{x}} \{ \hat{x}_{k}(\mathbf{y}) \}$. 
If the estimator $\hat{x}_{k}(\mathbf{y})$ is a continuous function of $\mathbf{y}$, and moreover for every $\mathbf{y} \in \mathbb{R}^{M}$ we have 
\begin{equation}
\label{equ_bound_suff_cond_exist_par_der_CS_recovery}
|\hat{x}_{k}(\mathbf{y})| \leq C \| \mathbf{y} \|_{2}^{L}
\end{equation}
with some constants $C,L \in \mathbb{R}_{+}$, then 
the partial derivatives $\frac{ \partial^{\mathbf{e}_{l}} \gamma(\mathbf{x})}{\partial \mathbf{x}^{\mathbf{e}_{l}}}$ exist for every $l \in [N]$ and are given by
\begin{equation}
\label{equ_expr_par_der_CS_recovery_scheme}
\frac{ \partial^{\mathbf{e}_{l}} \gamma(\mathbf{x})}{\partial \mathbf{x}^{\mathbf{e}_{l}}}
= \frac{1}{\sigma^{2}} \mathsf{E}_{\mathbf{x}} \big\{ \hat{x}_{k}(\mathbf{y})   (\mathbf{y}-\mathbf{H}\mathbf{x})^{T} \mathbf{H} \mathbf{e}_{l}   \big\}. 
\end{equation} 
\end{lemma} 
\begin{proof}
We have 
\begin{align} 
\label{equ_expr_par_der_CS_reocevery_schme_mean_integral}
\gamma(\mathbf{x})  = \mathsf{E}_{\mathbf{x}} \big\{ \hat{x}_{k}(\mathbf{y}) \big\} 
= \frac{1}{(2 \pi \sigma^{2})^{M/2}} \int_{\mathbf{y} \in \mathbb{R}^{M}} \hat{x}_{k}(\mathbf{y}) \exp \left( - \frac{1}{2 \sigma^{2}} \| \mathbf{y} - \mathbf{H} \mathbf{x} \|^{2}_{2} \right) d \mathbf{y}, 
\end{align} 
where the existence of the integral follows from the dominated convergence theorem \cite{HalmosMeasure,RudinBookPrinciplesMatheAnalysis,RudinBook} since (i) the function 
\begin{equation*}
\hat{x}_{k}(\mathbf{y}) \exp \left( - \frac{1}{2 \sigma^{2}} \| \mathbf{y} - \mathbf{H} \mathbf{x} \|^{2}_{2}Ê\right)
\end{equation*}
is continuous and therefore measurable and (ii) it is upper bounded in magnitude (dominated) by the function 
$C \| \mathbf{y} \|_{2}^{L} \exp \left( - \frac{1}{2 \sigma^{2}} \| \mathbf{y} - \mathbf{H} \mathbf{x} \|^{2}_{2}\right)$, which is obviously integrable. 
We further have
\begin{align}
\label{equ_expr_par_der_CS_reocevery_schme_mean_integral_proof_part_der}
\frac{ \partial^{\mathbf{e}_{l}} \gamma(\mathbf{x})}{\partial \mathbf{x}^{\mathbf{e}_{l}}} & \stackrel{\eqref{equ_expr_par_der_CS_reocevery_schme_mean_integral}}{=}
\frac{ \partial^{\mathbf{e}_{l}}} {\partial \mathbf{x}^{\mathbf{e}_{l}}} \frac{1}{(2 \pi \sigma^{2})^{M/2}} \int_{\mathbf{y} \in \mathbb{R}^{M}} 
\hat{x}_{k}(\mathbf{y}) \exp \left( - \frac{1}{2 \sigma^{2}} \| \mathbf{y} - \mathbf{H} \mathbf{x} \|^{2}_{2} \right) d \mathbf{y} \nonumber \\[4mm]
& \stackrel{(a)}{=} \frac{1}{(2 \pi \sigma^{2})^{M/2}} \int_{\mathbf{y} \in \mathbb{R}^{M}} \hat{x}_{k}(\mathbf{y}) \frac{ \partial^{\mathbf{e}_{l}} \exp \left( - \frac{1}{2 \sigma^{2}} \| \mathbf{y} - \mathbf{H} \mathbf{x} \|^{2}_{2} \right)}{\partial \mathbf{x}^{\mathbf{e}_{l}}}  d \mathbf{y} \nonumber \\[4mm]
& \stackrel{(b)}{=}  \frac{1}{(2 \pi \sigma^{2})^{M/2}} \int_{\mathbf{y} \in \mathbb{R}^{M}} \hat{x}_{k}(\mathbf{y})  
\frac{1}{\sigma^{2}} (\mathbf{y}-\mathbf{H}\mathbf{x})^{T} \mathbf{H} \mathbf{e}_{l} \exp \left( - \frac{1}{2 \sigma^{2}} \| \mathbf{y} - \mathbf{H} \mathbf{x} \|^{2}_{2} \right) d \mathbf{y} \nonumber \\[4mm]
&=  \frac{1}{\sigma^{2}} \mathsf{E}_{\mathbf{x}} \big\{ \hat{x}_{k}(\mathbf{y})   (\mathbf{y}-\mathbf{H}\mathbf{x})^{T} \mathbf{H} \mathbf{e}_{l}   \big\} ,
\end{align} 
where step $(a)$ follows from a change of the order of differentiation and integration (see \cite[Theorem 1.5.8]{LC}), and step $(b)$ is due to the chain rule for differentiation 
\cite{RudinBookPrinciplesMatheAnalysis}. The existence of the last integral in \eqref{equ_expr_par_der_CS_reocevery_schme_mean_integral_proof_part_der} can again be verified by the continuity of the integrand and the validity of the bound \eqref{equ_bound_suff_cond_exist_par_der_CS_recovery} 
on the estimator $\hat{x}_{k}(\mathbf{y})$. 
\end{proof}
For the application of Lemma \ref{lem_cond_exist_CS_recovery_partial_der}, we have to verify for a given CS recovery scheme, i.e., for a given estimator $\hat{\mathbf{x}}(\mathbf{y})$, 
that its entries $\hat{x}_{k}(\mathbf{y})$ satisfy the bound \eqref{equ_bound_suff_cond_exist_par_der_CS_recovery} and are continuous functions of $\mathbf{y}$. While showing 
\eqref{equ_bound_suff_cond_exist_par_der_CS_recovery} is rather straightforward for the BP and OMP, the rigorous proof of continuity is subtle. Consider, e.g.,  
the BP recovery method \cite{Chen98atomicdecomposition}, which is defined as the solution to the problem 
\begin{equation}
\label{equ_def_BP}
\min_{\mathbf{x} \in \mathbb{R}^{N}} \| \mathbf{x} \|_{1} \quad \mbox{s.t.} \quad \mathbf{H} \mathbf{x} = \mathbf{y}.  
\end{equation}
Since there may be multiple solutions to this problem, there are many possible estimator functions $\hat{\mathbf{x}}_{\text{BP}}(\cdot)$ corresponding to the BP. However, 
some popular implementations (e.g., the primal dual algorithm presented in \cite{BoydConvexBook,ell1magicUsersGuide}) 
of the estimator $\hat{\mathbf{x}}_{\text{BP}}(\mathbf{y}): \mathbb{R}^{M} \rightarrow \mathbb{R}^{N}$ yield a continuous function of the observation $\mathbf{y}$. 
The resulting approximate BP estimator is then continuous ``by design.'' 
As can be verified easily, the statement and proof of Lemma \ref{lem_cond_exist_CS_recovery_partial_der} remain unchanged if one replaces the requirement of a continuous estimator function 
with that of a Lebesgue measurable estimator function. 
It seems then fairly unlikely that any reasonable and practically implemented CS recovery scheme does not yield an 
estimator that is both measurable and bounded as in \eqref{equ_bound_suff_cond_exist_par_der_CS_recovery}.
 
We note that condition \eqref{equ_spark_cond} is necessary for a matrix $\mathbf{H}$ to have the RIP of order $S$ with a constant $\delta_{S}<1$.\footnote{Indeed, assume that $\spark(\mathbf{H}) \leq S$ which 
means that there exists an index set $\mathcal{I} \subseteq [N]$ consisting of $S$ indices, i.e., $| \mathcal{I}| =S$ such that the columns of $\mathbf{H}_{\mathcal{I}}$ are linearly dependent. This implies that 
there is a coefficient vector $\mathbf{z}_{0} \in \mathbb{R}^{S}$ with $\mathbf{z}_{0} \neq \mathbf{0}$ (i.e., $\| \mathbf{z}_{0}\|^{2}_{2} > 0$) 
such that $ \mathbf{H}_{\mathcal{I}} \mathbf{z}_{0} = \mathbf{0}$ and in turn $\|\mathbf{H}_{\mathcal{I}} \mathbf{z}_{0}Ê\|^{2}_{2} = 0$. Therefore, 
there cannot exist a constant $\delta_{K} < 1$ that satisfies the RIP condition \eqref{equ_RIP_condition} for every $\mathbf{z} \in \mathbb{R}^{S}$.}
Moreover, for CS applications, one favors in general measurement matrices $\mathbf{H}$ that have RIP constants close to zero, i.e., $\delta_{S} \approx 0$ \cite{Can06a,ROMP_stab,AnalysisOMPRIP,CoSAMP,ZvikaCoherenceTSP}. However, in this 
case the bound in \eqref{equ_CS_measur_matrix_asilomar_bound} is nearly identical to the bound \eqref{equ_bound_asilomar_2} for the SLM with $\mathbf{H} = \mathbf{I}$, given explicitly as 
\begin{equation}
L_{\mathcal{M}_{\text{SLM}}} \geq  \exp\left( - \frac{1}{\sigma^{2}} \| (\mathbf{I} - \mathbf{P}_{\mathcal{K}})  \mathbf{x}_{0} \|^{2}_{2} \right) \sigma^{2} \mathbf{r}_{\mathbf{x}_{0}}^{T}  \mathbf{r}_{\mathbf{x}_{0}}.
\end{equation}  
Note that for $\mathbf{H}= \mathbf{I}$, implying that 
\begin{equation} 
\mathbf{P}_{\mathcal{K}} = \mathbf{H}_{\mathcal{K}} \mathbf{H}_{\mathcal{K}}^{\dagger} = \sum_{l \in \mathcal{K}} \mathbf{e}_{l} \mathbf{e}_{l}^{T}, 
\end{equation} 
the effect of the multiplication $(\mathbf{I} - \mathbf{P}_{\mathcal{K}})\mathbf{x}_{0}$ is a zeroing of all entries of the vector $\mathbf{x}_{0}$ whose indices belong to $\mathcal{K}$, i.e., 
\begin{equation}
(\mathbf{I} - \mathbf{P}_{\mathcal{K}})  \mathbf{x}_{0}= \mathbf{x}_{0}^{\supp(\mathbf{x}_{0}) \setminus \mathcal{K}}.
\end{equation}
Thus, if one uses a ``good'' CS measurement matrix, i.e., with RIP constant $\delta_{S}$ close to zero, then Theorem \ref{thm_CS_measur_matrix_asilomar_bound} 
suggests that the minimum achievable variance 
is close to the minimum achievable variance of another instance of the SLM with the same $S$, and $N$ but with system the matrix being the identity matrix $\mathbf{I}$ 
instead of the measurement matrix $\mathbf{H}$. This means that in terms of achievable estimation accuracy, there is no loss of information incurred by applying the CS measurement matrix 
$\mathbf{H} \in \mathbb{R}^{M \times N}$ and thereby reducing the signal dimension from $N$ to $M$, where $M \ll N$ in general. 
This agrees with the fact that if a CS measurement matrix satisfies the RIP with RIP constants that are sufficiently small, 
then one can recover - e.g., by using the BP - the sparse parameter vector $\mathbf{x} \in \mathcal{X}_{S}$ 
from the compressed observation $\mathbf{y}$ given by \eqref{equ_measurment_equ_CS} up to a reconstruction error that is solely determined by the noise in \eqref{equ_measurment_equ_CS} \cite{Can06a,DantzigCandes}. 

The CS measurement model in \eqref{equ_measurment_equ_CS} assumes that the noise component is added after the measurement process. However, 
in some applications it may be more accurate to model the effect of the noise before the measurement process takes place, i.e., to use the following measurement model: 
\begin{equation}
\label{equ_CS_model_noise_before_measurement}
\mathbf{y}' = \mathbf{H} (\mathbf{x} + \mathbf{n}),
\end{equation}
where the meanings of $\mathbf{H}$, $\mathbf{x}$, and $\mathbf{n}$ are the same as for \eqref{equ_measurment_equ_CS}. 
Note that the observation $\mathbf{y}'$ is obtained by applying a mapping, i.e., the matrix multiplication by $\mathbf{H}$, to the 
observation $\mathbf{y} = \mathbf{x} + \mathbf{n}$ of the SSNM.  
Therefore, if the CS measurement process is modeled by \eqref{equ_CS_model_noise_before_measurement}, then we can analyze the resulting minimum variance problem of estimating the sparse vector $\mathbf{x} \in \mathcal{X}_{S}$ from the observation $\mathbf{y}'$ by using the results on the SSNM in Section \ref{sec_SSNM}, in particular Theorem \ref{thm_general_min_achiev_var_LMV_SSNM}, together with Theorem \ref{thm_data_proc_inequ_classical_est}. 
We thus obtain the following result: 

\begin{theorem}
\label{thm_CS_measur_matrix_after_noise}
Consider the observation model \eqref{equ_CS_model_noise_before_measurement} and an estimator $\hat{x}_{k}(\mathbf{y}'): \mathbb{R}^{M} \rightarrow \mathbb{R}$, which 
may arise from a CS recovery scheme, 
of $x_{k}$. The estimator $\hat{x}_{k}(\mathbf{y}')$ can use only the observation 
$\mathbf{y}'$ and is assumed to have finite variance at $\mathbf{x}_{0} \in \mathcal{X}_{S}$, i.e., $v(\hat{x}_{k}(\cdot); \mathbf{x}_{0}) < \infty$. 
 We assume that the mean of this estimator exists for all parameter vectors $\mathbf{x} \in \mathcal{X}_{S}$ and denote its mean function by 
$\gamma(\cdot): \mathcal{X}_{S} \rightarrow \mathbb{R}: \gamma(\mathbf{x}) = \mathsf{E}_{\mathbf{x}} \{ \hat{x}_{k}(\mathbf{y}') \}$.
Then, we have that
\begin{equation}
\label{equ_CS_measur_matrix_after_noise}
v(\hat{x}_{k}(\cdot); \mathbf{x}_{0}) \geq  \|   a_{\mathbf{x}_{0}}[\mathbf{p}] \|^{2}_{\ell^{2}(\mathbb{Z}_{+}^{N} \cap \mathcal{X}_{S})} - \big[ \gamma(\mathbf{x}_{0}) \big]^{2}, 
\end{equation}  
where the coefficient sequence $a_{\mathbf{x}_{0}}[\mathbf{p}] \in \ell^{2}(\mathbb{Z}_{+}^{N} \cap \mathcal{X}_{S})$ are given by \eqref{equ_coeffs_general_min_achiev_var_LMV_SSNM}.
\end{theorem}

\begin{proof}
Let us denote by $\mathcal{M}'$ the minimum variance problem obtained for the observation model \eqref{equ_CS_model_noise_before_measurement},the parameter set $\mathcal{X}_{S}$, the fixed parameter vector $\mathbf{x}_{0} \in \mathcal{X}_{S}$, the parameter function 
$g(\mathbf{x}) = x_{k}$, and prescribed bias function $c(\cdot) \triangleq \gamma(\cdot) - x_{k} = \mathsf{E}_{\mathbf{x}} \big\{ \hat{x}_{k}(\mathbf{y}') \big \} - x_{k}$. 
Since $c(\cdot)$ is identical with the bias function $b(\hat{x}_{k}(\cdot), \mathbf{x})$ of the estimator $\hat{x}_{k}(\cdot)$ which is assumed to have finite variance at $\mathbf{x}_{0}$, we have that $c(\cdot)$ is valid for $\mathcal{M}'$ and $v(\hat{x}_{k}(\cdot); \mathbf{x}_{0}) \geq L_{\mathcal{M}'}$. 

The minimum variance problem $\mathcal{M}'$ is obtained from the minimum variance problem $\mathcal{M}_{\text{SSNM}}=\left( \mathcal{E}_{\text{SSNM}}, c(\cdot), \mathbf{x}_{0} \right)$ (which uses the same index $k \in [N]$ in its parameter function (see \eqref{equ_def_SSNM_est_problem}) as used within $\mathcal{M}'$) 
by a matrix multiplication of the observation of the SSNM with the matrix $\mathbf{H}$. Therefore, by
Theorem \ref{thm_data_proc_inequ_classical_est}, we have that $c(\cdot)$ is valid for $\mathcal{M}_{\text{SSNM}}$. Indeed, since $c(\cdot)$ is valid for $\mathcal{M}'$ we have that $L_{\mathcal{M}'}$ is finite and this implies via \eqref{equ_data_proc_inequ_classic_est_min_var} that $L_{\mathcal{M}_{\text{SSNM}}}$ is finite, i.e., $c(\cdot)$ must be valid for $\mathcal{M}_{\text{SSNM}}$.  
Since $c(\cdot)$ is valid for $\mathcal{M}_{\text{SSNM}}$ we have that $L_{\mathcal{M}_{\text{SSNM}}}$ is given by \eqref{equ_ssnm_expr_min_ach_var} in Theorem \ref{thm_general_min_achiev_var_LMV_SSNM}. 
The bound \eqref{equ_CS_measur_matrix_after_noise} follows then as 
\begin{equation}
 v(\hat{x}_{k}(\cdot); \mathbf{x}_{0}) \geq L_{\mathcal{M}'} \stackrel{\eqref{equ_data_proc_inequ_classic_est_min_var}}{\geq} L_{\mathcal{M}_{\text{SSNM}}}
 \stackrel{\eqref{equ_ssnm_expr_min_ach_var}}{=} \|   a_{\mathbf{x}_{0}}[\mathbf{p}] \|^{2}_{\ell^{2}(\mathbb{Z}_{+}^{N} \cap \mathcal{X}_{S})} - \big[\gamma(\mathbf{x}_{0})\big]^{2}.
\end{equation}
\end{proof} 
It is important to note that the statement and proof of Theorem \ref{thm_CS_measur_matrix_after_noise} are also valid for the case where $\mathbf{H}$ in \eqref{equ_CS_model_noise_before_measurement} is a random matrix, which is independent of the noise $\mathbf{n}$. This is practically relevant since it is often the case that CS measurement matrices are modeled as random matrices.


\section{Comparison of the Bounds with Existing Estimators} 
\label{sec_numerical_SLM}
Let us now compare the lower bounds derived so far with the actual variance behavior of some well-known estimators. 
In what follows, we will denote by $L^{(c(\cdot),k)}_{\mathcal{K},a}(\mathbf{x}_{0})$, $L^{(c(\cdot),k)}_{\mathcal{K},b}(\mathbf{x}_{0})$ the bounds \eqref{equ_bound_asilomar_1} and \eqref{equ_bound_asilomar_2} obtained by 
Theorem \ref{thm_bound_asilomar} and Theorem \ref{thm_bound_asilomar_2}, respectively, using the index set $\mathcal{K}$ and prescribed bias function $c(\cdot)$.

Given an estimator $\hat{\mathbf{x}}(\cdot)$ whose bias is equal to $\mathbf{c}(\mathbf{x})$ and mean equal to ${\bm \gamma}(\mathbf{x}) \triangleq \mathbf{c}(\mathbf{x}) + \mathbf{x}$, we obtain due to \eqref{equ_sum_min_var_scalar_min_var},
a lower bound on the estimator variance $v(\hat{\mathbf{x}}(\cdot);\mathbf{x}_{0})$ by summing the quantities $L^{(c_{k}(\cdot),k)}_{\mathcal{K},a/b}(\mathbf{x}_{0})$, i.e., 
\begin{equation}
\label{equ_lower_bound_var_vec_est_numerical_SLM}
v(\hat{\mathbf{x}}(\cdot);\mathbf{x}_{0}) \geq L^{a/b}_{{\bm \gamma}(\cdot),\mathbf{x}_{0}} \triangleq \sum_{k \in [N]} \rmv\rmv L^{(c_{k}(\cdot),k)}_{\mathcal{K}_{k},a/b} (\mathbf{x}_{0}).
\end{equation}
Note that in \eqref{equ_lower_bound_var_vec_est_numerical_SLM}. we allow the index $\mathcal{K} = \mathcal{K}_{k}$ in \eqref{equ_bound_asilomar_1}, \eqref{equ_bound_asilomar_2} 
to vary with the index $k$. 

\subsection{SLM Viewpoint on Fourier Analysis}

The first experiment is inspired by \cite[Example 4.2 -- `Fourier Analysis']{kay} and considers the general SLM with $\sigma^{2}=1$ and system matrix 
\begin{equation} 
\mathbf{H} = \begin{pmatrix} \cos( \theta_{1} 0) & \ldots & \cos( \theta_{L} 0) &  \sin( \theta_{1} 0) & \ldots & \sin( \theta_{L} 0)  \\Ê
					     \cos( \theta_{1} 1) & \ldots & \cos( \theta_{L} 1) &  \sin( \theta_{1} 1) & \ldots & \sin( \theta_{L} 1)  \\Ê
					     \vdots & \ddots & \vdots & \vdots & \ddots & \vdots  \\Ê
					        \cos( \theta_{1} (M-1)) & \ldots & \cos( \theta_{L} (M-1)) &  \sin( \theta_{1} (M-1)) & \ldots & \sin( \theta_{L} (M-1))  \end{pmatrix}, 
\end{equation} 
i.e., $\mathbf{H} \in \mathbb{R}^{M \times N}$ with $N = 2L$. Note that this system matrix corresponds to a finite-length discrete Fourier transform evaluated at the frequencies\footnote{Here, we consider the normalized frequency $\theta \in [0,1]$ of a discrete-time 
harmonic signal $x[n] = \sin(2 \pi \theta n)$.} $\{ \theta_{l} \}_{l \in [L]}$. 
In our simulation we choose $M=128$, $L=8$, and $\theta_{l} = 0.2 + 3.9 \cdot 10^{-3} (l-1)$. Our choice of $\theta_{l}$ corresponds to a frequency resolution of $\Delta \theta = 3.9Ê\times 10^{-3}$; this is about half of the DFT frequency resolution, which is given by $\Delta \theta = \frac{1}{128} \approx 7.8 \times 10^{-3}$. We assume that 
the parameter vector can have at maximum $4$ nonzero entries, i.e., we consider the SLM with sparsity degree $S=4$.

\begin{figure}
\vspace{-1mm}
\centering
\psfrag{SNR}[c][c][.9]{\uput{3.4mm}[270]{0}{\hspace{3mm}SNR [dB]}}
\psfrag{title}[c][c][.9]{\uput{2.5mm}[270]{0}{}}
\psfrag{x_0}[c][c][.9]{\uput{0.3mm}[270]{0}{$0$}}
\psfrag{x_0_001}[c][c][.9]{\uput{0.3mm}[270]{0}{$-60$}}
\psfrag{x_0_01}[c][c][.9]{\uput{0.3mm}[270]{0}{$-40$}}
\psfrag{x_0_1}[c][c][.9]{\uput{0.3mm}[270]{0}{$-20$}}
\psfrag{x_1}[c][c][.9]{\uput{0.3mm}[270]{0}{$0$}}
\psfrag{x_10}[c][c][.9]{\uput{0.3mm}[270]{0}{$20$}}
\psfrag{x_100}[c][c][.9]{\uput{0.3mm}[270]{0}{$40$}}
\psfrag{y_0}[c][c][.9]{\uput{0.1mm}[180]{0}{$0$}}
\psfrag{y_4}[c][c][.9]{\uput{0.1mm}[180]{0}{$4$}}
\psfrag{y_8}[c][c][.9]{\uput{0.1mm}[180]{0}{$8$}}
\psfrag{y_12}[c][c][.9]{\uput{0.1mm}[180]{0}{$12$}}
\psfrag{y_16}[c][c][.9]{\uput{0.1mm}[180]{0}{$16$}}
\psfrag{y_20}[c][c][.9]{\uput{0.1mm}[180]{0}{$20$}}
\psfrag{variance}[c][c][.9]{\uput{1mm}[90]{0}{variance/bound}}
\psfrag{bML}[l][l][0.8]{bound on $v(\ML_est(\cdot);\mathbf{x}_0)$}
\psfrag{data1}[l][l][0.8]{$v(\hat{\mathbf{x}}_{\text{OMP}}(\cdot); \mathbf{x}_{0})$}
\psfrag{data2}[l][l][0.8]{$L_{\text{OMP},b}(\mathbf{x}_{0})$}
\psfrag{data3}[l][l][0.8]{$L_{\text{OMP},a}(\mathbf{x}_{0})$}
\psfrag{data4}[l][l][0.8]{$L_{\text{OMP},c}(\mathbf{x}_{0})$}
\psfrag{CRB}[l][l][0.8]{Oracle bound}
\psfrag{MLarrow}[l][l][0.8]{ML}
\psfrag{T3}[l][l][0.8]{HT \!($T \!\rmv=\! 3$)}
\psfrag{T4}[l][l][0.8]{HT \!($T \!\rmv=\! 4$)}
\psfrag{T5}[l][l][0.8]{HT \!($T \!\rmv=\! 5$)}
\centering
\hspace*{-0mm}\includegraphics[height=8cm,width=15cm]{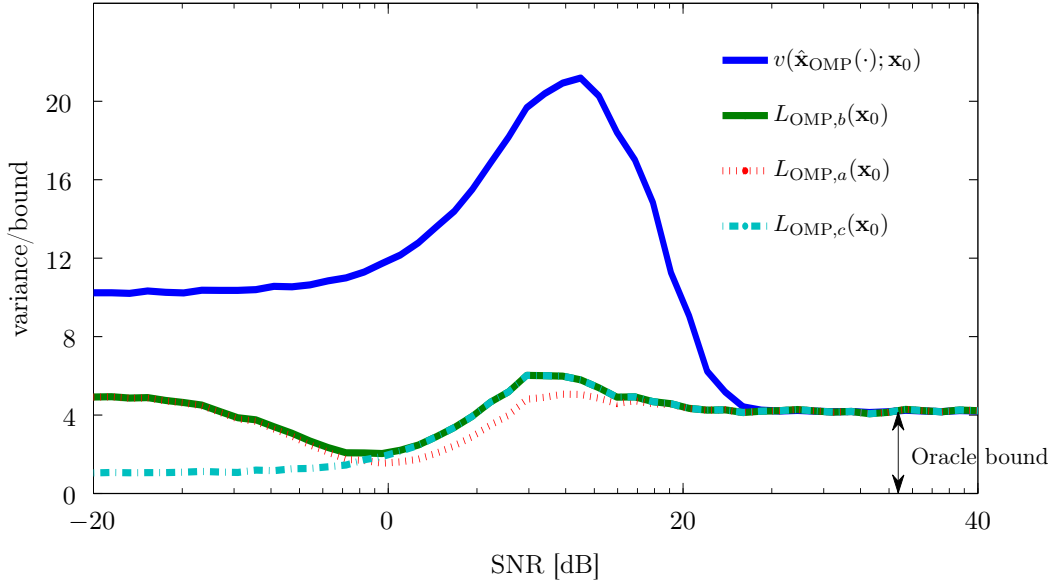}
\vspace*{1mm}
  \caption{Variance of the OMP estimator and corresponding lower bounds versus the SNR,
  for the SLM with $N\!=\!16$, $M = 128$, $\sigma^{2}=1$ and $S\!=\!4$.} 
\label{fig_SLM_high_res}
\vspace*{3.5mm}
\end{figure}

In Fig. \ref{fig_SLM_high_res}, we compare the variance of the OMP estimator $\hat{\mathbf{x}}_{\text{OMP}}(\mathbf{y})$, which is obtained by applying OMP (see \cite{GreedisGood,TroppGilbertOMP}) with $S$ iterations to the observation, with the lower bounds $L_{\text{OMP},a/b}(\mathbf{x}_{0})$ defined as (cf.\ \eqref{equ_lower_bound_var_vec_est_numerical_SLM})
\begin{align}
\label{equ_def_lower_bound_OMP}
L_{\text{OMP},a/b}(\mathbf{x}_{0}) & =\sum_{k \in [N]} \rmv\rmv L^{(c_{k}(\cdot),k)}_{\mathcal{K}_{k},a/b} (\mathbf{x}_{0}), 
\end{align} 
where $c_{k}(\cdot)$ is chosen as $c_{k}(\cdot) = \gamma_{k,\text{OMP}}(\mathbf{x}) - x_k$, with $\gamma_{k,\text{OMP}}(\mathbf{x}) \triangleq \mathsf{E}_{\mathbf{x}} \big\{ \hat{x}_{k,\text{OMP}}(\mathbf{y}) \big \}$ 
being the actual mean function of the OMP estimator. 
We also plotted the lower bound $L_{\text{OMP},c}(\mathbf{x}_{0})$ obtained by summing the lower bound \eqref{equ_CRB_SLM_not_full_sparsity}, 
\eqref{equ_CRB_SLM_full_sparsity} (using $\gamma_{k,\text{OMP}}-x_{k}$ as prescribed bias function)  in Theorem \ref{thm_CRB_SLM} over all indices $k \in [N]$. 
By \eqref{equ_sum_min_var_scalar_min_var}, we have that $v(\hat{\mathbf{x}}_{\text{OMP}}(\cdot); \mathbf{x}_{0})Ê\geq L_{\text{OMP},c}(\mathbf{x}_{0})$. 

The variance $v(\hat{\mathbf{x}}_{\text{OMP}}(\cdot); \mathbf{x}_{0})$ and the lower bounds $L_{\text{OMP},a/b/c}(\mathbf{x}_{0})$ are 
evaluated for parameter vectors $\mathbf{x}_{0} = \sqrt{\text{SNR}} \sigma \mathbf{x}_{1}$, where $\mathbf{x}_{1} \in \{0,1\}^{16}$, 
$\supp(\mathbf{x}_{1}) = \{3,6,67,70\}$, and $\sqrt{\text{SNR}}$ 
varies between $0.1$ and $100$.
The variance $v(\hat{\mathbf{x}}_{\text{OMP}}(\cdot); \mathbf{x}_{0})$ of the OMP estimator is estimated by means of Monte Carlo simulation. 
For the evaluation of the bound $L_{\text{OMP},a/b/c}(\mathbf{x}_{0})$ (cf.\ \eqref{equ_bound_asilomar_1}, \eqref{equ_bound_asilomar_2}, \eqref{equ_CRB_SLM_not_full_sparsity}, 
\eqref{equ_CRB_SLM_full_sparsity}), we have to compute the first-order partial derivatives of the mean 
function $\gamma_{k,\text{OMP}}(\mathbf{x})$.  
We estimated these partial derivatives by means of Lemma \ref{lem_cond_exist_CS_recovery_partial_der} and a Monte Carlo simulation. 
This method of estimating the derivative of an estimator bias or mean function is also used and explained in detail 
in \cite{HeroUniformCRB}.
For simplicity, the index sets $\mathcal{K}_{k}$ in \eqref{equ_def_lower_bound_OMP} are chosen as $\mathcal{K}_{k} = \supp(\mathbf{x}_{0})$ for $k \in \supp(\mathbf{x}_{0})$ and $\mathcal{K}_{k} = \{k\}$ for $k \notin \supp(\mathbf{x}_{0})$. This is the simplest non-trivial choice that yields a bound $L_{\text{OMP},b}(\mathbf{x}_{0})$ that is tighter than the bound $L_{\text{OMP},c}(\mathbf{x}_{0})$, which is based on the sparse CRB in Theorem \ref{thm_CRB_SLM} (cf.\ \cite{ZvikaCRB}). 
In Fig. \ref{fig_SLM_high_res}, we also indicated the ``oracle CRB,'' which is defined as the CRB under the assumption that one knows 
the support of $\mathbf{x}_{0}$, i.e., the oracle bound is the CRB of a linear {G}aussian model with system matrix $\mathbf{H}_{\supp(\mathbf{x}_{0})}$ and is thus 
given by \cite{kay} 
\begin{equation}
\sigma^2 \trace \big\{ \big[ \mathbf{H}_{\supp(\mathbf{x}_{0})}^{T} \mathbf{H}_{\supp(\mathbf{x}_{0})} \big]^{-1} \big\} \approx 4.19 \times \sigma^{2}.
\end{equation}
As can be seen from Fig. \ref{fig_SLM_high_res}, there are two SNR regimes regarding the variance of the OMP estimator: Below $20$ dB, i.e., in the low-SNR 
regime the variance of $\hat{\mathbf{x}}_{\text{OMP}}(\cdot)$ is significantly higher than the oracle CRB and also significantly higher than the lower bounds 
$L_{\text{OMP},a/b/c}(\mathbf{x}_{0})$. This suggests that there might exist estimators with a reduced variance but the same bias and mean as the OMP estimator. 
Above $20$ dB, i.e., in the high-SNR regime the variance $\hat{\mathbf{x}}_{\text{OMP}}(\cdot)$ shows a fast convergence towards the oracle CRB as well as to the bounds $L_{\text{OMP},a/b/c}(\mathbf{x}_{0})$, which is due to the fact that for SNR values above $20$ dB the OMP estimator is able to detect the support of $\mathbf{x}_{0}$ with very high probability. Note also that the curves of $L_{\text{OMP},a/b}(\mathbf{x}_{0})$ agrees with the discussion around \eqref{equ_SLM_diff_nonnegative_asilomar_bounds}, i.e., that the bound $L_{\text{OMP},b}(\mathbf{x}_{0})$ tends 
to be higher than $L_{\text{OMP},a}(\mathbf{x}_{0})$ in general.

\subsection{Minimum Variance Analysis of Estimators for the SSNM}

Let us now consider the SSNM in \eqref{equ_observation_model_SSNM} for $N \!=\! 50$, $S \!=\! 5$, and $\sigma^2 = 1$.
We will compute the lower bound on the estimator variance given in \eqref{equ_lower_bound_var_vec_est_numerical_SLM} and compare it with the 
variance of two established estimators, namely, the maximum likelihood (ML) estimator and the hard-thresholding (HT) estimator.
The ML estimator is given by 
\[
\ML_est(\mathbf{y}) \,\triangleq\, \argmax_{\mathbf{x}' \in \mathcal{X}_{S}} f ( \mathbf{y};  \mathbf{x}' ) \ist=\, {\mathsf P}_{\!S} ( \mathbf{y}  ) \,,
\]
where the operator $\mathsf{P}_{\! S}$ retains the $S$ largest (in magnitude) entries and zeros out
all others. The HT estimator $\hat{\mathbf{x}}_{\text{HT}}  (\mathbf{y})$ is given 
\vspace{-.5mm}
by
\begin{equation}
\label{equ_def_thr_func}
\hat{x}_{\text{HT},k}  (\mathbf{y}) = \hat{x}_{\text{HT},k}  (y_k) = \begin{cases} y_k \,,  & |y_k| \geq T\\[.5mm]
  0 \,, & \text{else} \ist,
\end{cases} 
\end{equation} 
where $T$ is a fixed threshold. Note that in the limiting case where $T=0$, the HT estimator coincides with the least square (LS) estimator \cite{kay,scharf91,LC} for the SSNM $\hat{\mathbf{x}}_{\text{LS}}(\mathbf{y}) = \mathbf{y}$. The prescribed bias $c_{k}(\cdot)$ in \eqref{equ_lower_bound_var_vec_est_numerical_SLM} is chosen as $k$th entry of the bias function 
of the ML and HT estimator, respectively.

The mean and variance of the HT estimator are given by 
\begin{equation}
\label{equ_mean_hard_th_est_ssnm}
\mathsf{E}_{\mathbf{x}} \big\{ \hat{x}_{\text{HT},k}(\mathbf{y})  \big\}  = \frac{1}{\sqrt{2 \pi \sigma^{2}}}\int_{y_{k} \in \mathbb{R}^{N} \setminus [-T,T]} y_{k} \exp \bigg( - \frac{1}{2 \sigma^{2}} (y_{k} - x_{k})^{2} \bigg) d y_{k}
\end{equation} 
and 
\begin{align}
\label{equ_var_hard_th_est_ssnm}
v(\hat{x}_{\text{HT}(\cdot),k}; \mathbf{x}) & \stackrel{\eqref{equ_rel_power_var}}{=} P(\hat{x}_{\text{HT},k}(\cdot); \mathbf{x}) - \big[ \mathsf{E}_{\mathbf{x}} \big\{ \hat{x}_{\text{HT},k}(\mathbf{y})  \big\}  \big]^{2}    \nonumber \\[4mm]
& = \frac{1}{\sqrt{2 \pi \sigma^{2}}}\int_{y_{k} \in \mathbb{R} \setminus [-T,T]} y_{k}^{2} \exp \bigg( - \frac{1}{2 \sigma^{2}} (y_{k} - x_{k})^{2} \bigg) d y_{k}-\big[ \mathsf{E}_{\mathbf{x}} \big\{ \hat{x}_{\text{HT},k}(\mathbf{y})  \big\}  \big]^{2} ,
\end{align} 
respectively. The integrals are calculated using numerical integration. 
For the evaluation of the mean and variance of the ML estimator, we used the (complicated) closed-form expressions derived in \cite{AlexZvikaJournal}. 

Let us denote by $c_{\text{ML},k}(\mathbf{x}) \triangleq b(\hat{x}_{\text{ML},k}(\cdot), \mathbf{x})$ and 
$c_{\text{HT},k}(\mathbf{x}) \triangleq b(\hat{x}_{\text{ML},k}(\cdot), \mathbf{x})$ the bias function of the $k$th entry of the ML and HT estimator, respectively.
Using these bias functions, we compare the variance of the HT and ML estimators with the lower bounds 
\begin{equation}
\label{equ_def_lound_asilomar_ML}
L_{\text{ML},a/b}(\mathbf{x}_{0}) \triangleq \sum_{k \in [N]} \rmv\rmv L^{(c_{\text{ML},k}(\cdot),k)}_{\mathcal{K}_{k},a/b} (\mathbf{x}_{0}) 
\end{equation} 
and 
\begin{equation}
\label{equ_def_lound_asilomar_HT}
L_{\text{HT},a/b}(\mathbf{x}_{0}) \triangleq \sum_{k \in [N]}  \rmv\rmv L^{(c_{\text{HT},k}(\cdot),k)}_{\mathcal{K}_{k},a/b} (\mathbf{x}_{0}), 
\end{equation} 
respectively. 
\begin{figure}
\vspace{-1mm}
\centering
\psfrag{SNR}[c][c][.9]{\uput{3.4mm}[270]{0}{\hspace{3mm}SNR [dB]}}
\psfrag{title}[c][c][.9]{\uput{2.5mm}[270]{0}{}}
\psfrag{x_0}[c][c][.9]{\uput{0.3mm}[270]{0}{$0$}}
\psfrag{x_0_001}[c][c][.9]{\uput{0.3mm}[270]{0}{$-30$}}
\psfrag{x_0_01}[c][c][.9]{\uput{0.3mm}[270]{0}{$-20$}}
\psfrag{x_0_1}[c][c][.9]{\uput{0.3mm}[270]{0}{$-10$}}
\psfrag{x_1}[c][c][.9]{\uput{0.3mm}[270]{0}{$0$}}
\psfrag{x_10}[c][c][.9]{\uput{0.3mm}[270]{0}{$10$}}
\psfrag{x_100}[c][c][.9]{\uput{0.3mm}[270]{0}{$20$}}
\psfrag{y_0}[c][c][.9]{\uput{0.1mm}[180]{0}{$0$}}
\psfrag{y_5}[c][c][.9]{\uput{0.1mm}[180]{0}{$5$}}
\psfrag{y_10}[c][c][.9]{\uput{0.1mm}[180]{0}{$10$}}
\psfrag{y_15}[c][c][.9]{\uput{0.1mm}[180]{0}{$15$}}
\psfrag{y_20}[c][c][.9]{\uput{0.1mm}[180]{0}{$20$}}
\psfrag{y_25}[c][c][.9]{\uput{0.1mm}[180]{0}{$25$}}
\psfrag{y_30}[c][c][.9]{\uput{0.1mm}[180]{0}{$30$}}
\psfrag{y_35}[c][c][.9]{\uput{0.1mm}[180]{0}{$35$}}
\psfrag{y_40}[c][c][.9]{\uput{0.1mm}[180]{0}{$40$}}
\psfrag{y_45}[c][c][.9]{\uput{0.1mm}[180]{0}{$45$}}
\psfrag{y_50}[c][c][.9]{\uput{0.1mm}[180]{0}{$50$}}
\psfrag{y_55}[c][c][.9]{\uput{0.1mm}[180]{0}{$55$}}
\psfrag{variance}[c][c][.9]{\uput{3mm}[90]{0}{variance/bound}}
\psfrag{bML}[l][l][0.8]{bound on $v(\ML_est(\cdot);\mathbf{x}_0)$}
\psfrag{ML}[l][l][0.8]{ML}
\psfrag{data3}[l][l][0.8]{$v(\hat{\mathbf{x}}_{\text{HT}}(\cdot); \mathbf{x}_0)$, $T \!\rmv=\! 2$}
\psfrag{data5}[l][l][0.8]{$v(\hat{\mathbf{x}}_{\text{HT}}(\cdot); \mathbf{x}_0)$, $T \!\rmv=\! 3$}
\psfrag{data7}[l][l][0.8]{$v(\hat{\mathbf{x}}_{\text{HT}}(\cdot); \mathbf{x}_0)$, $T \!\rmv=\! 4$}
\psfrag{data4}[l][l][0.8]{$L_{\text{HT},a/b}(\mathbf{x}_{0}) $, $T \!\rmv=\! 2$}
\psfrag{data6}[l][l][0.8]{$L_{\text{HT},a/b}(\mathbf{x}_{0}) $, $T \!\rmv=\! 3$}
\psfrag{data8}[l][l][0.8]{$L_{\text{HT},a/b}(\mathbf{x}_{0}) $, $T \!\rmv=\! 4$}
\psfrag{data1}[l][l][0.8]{$v(\hat{\mathbf{x}}_{\text{LS}}(\cdot); \mathbf{x}_0)$}
\psfrag{data2}[l][l][0.8]{$L_{\text{HT},a/b}(\mathbf{x}_{0}) $, $T \!\rmv=\! 0$}
\psfrag{data9}[l][l][0.8]{$v(\hat{\mathbf{x}}_{\text{ML}}(\cdot); \mathbf{x}_0)$}
\psfrag{data10}[l][l][0.8]{$L_{\text{ML},a/b}(\mathbf{x}_{0})$}
\psfrag{MLarrow}[l][l][0.8]{ML}
\psfrag{T=3}[l][l][0.8]{\uput{0mm}[90]{0}{\hspace{0mm}$T \!\rmv=\! 3$}}
\psfrag{T=4}[l][l][0.8]{$T \!\rmv=\! 4$}
\psfrag{T=2}[l][l][0.8]{\uput{.5mm}[0]{0}{\vspace*{3mm} \!\!\!\!$T \!\rmv=\! 2$}}
\psfrag{T=0}[l][l][0.8]{\uput{0mm}[0]{0}{\hspace{-10mm}$T \!\rmv=\! 0\,\ist$(LS)}}
\psfrag{oracle}[l][l][0.8]{\uput{0mm}[0]{0}{$S \sigma^{2}$}}
\centering
\hspace*{-0mm}\includegraphics[height=8cm,width=15cm]{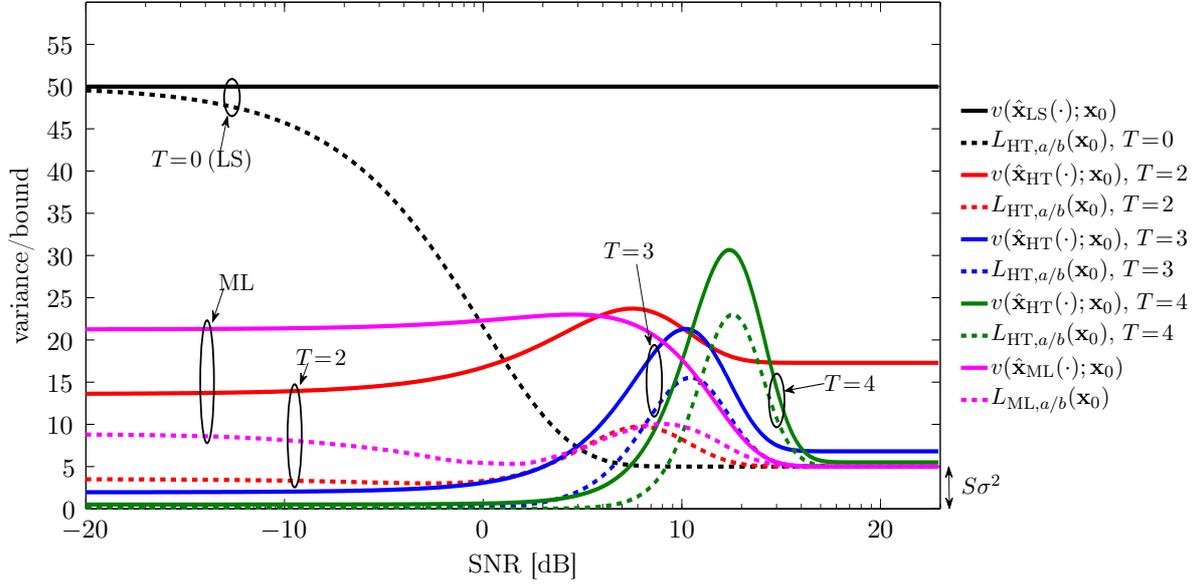}
\vspace*{1mm}
  \caption{Variance of the ML and HT estimators and corresponding lower bounds versus the SNR,
  for the SSNM
  with $N\!=\!50$, $S\!=\!5$ and $\sigma^{2}=1$.} 
\label{fig_bounds_1}
\vspace*{3.5mm}
\end{figure}
As before, the index set $\mathcal{K}_{k} \subseteq [N]$ is allowed to vary with the index $k$. In particular, we choose $\mathcal{K}_{k} = \supp(\mathbf{x})$ 
for $k \in \supp(\mathbf{x})$ and $\mathcal{K}_{k} = \{k\} \cup \{ \supp(\mathbf{x}_{0}) \setminus \{j_{0} \} \}$ for $k \notin \supp(\mathbf{x}_{0})$, where $j_{0}$ denotes the index of the $S$-largest (in magnitude) entry of $\mathbf{x}_{0}$. It can be verified easily that for these choices of the index sets $\mathcal{K}_{k}$, the 
bounds for the HT estimator, $L^{(c_{\text{HT},k}(\cdot),k)}_{\mathcal{K}_{k},a} (\mathbf{x}_{0})$ and $L^{(c_{\text{HT},k}(\cdot),k)}_{\mathcal{K}_{k},b} (\mathbf{x}_{0})$, as well as the bounds for the ML estimator, 
$L^{(c_{\text{ML},k}(\cdot),k)}_{\mathcal{K}_{k},a} (\mathbf{x}_{0})$ and $L^{(c_{\text{ML},k}(\cdot),k)}_{\mathcal{K}_{k},b} (\mathbf{x}_{0})$
coincide for the SSNM (i.e., when $\mathbf{H} = \mathbf{I}$). Thus, we have in this case that $L_{\text{HT},a}(\mathbf{x}_{0}) = L_{\text{HT},b}(\mathbf{x}_{0})$ and $L_{\text{ML},a}(\mathbf{x}_{0}) = L_{\text{ML},b}(\mathbf{x}_{0})$.

In order to evaluate the bounds $L^{(c_{\text{ML},k}(\cdot),k)}_{\mathcal{K}_{k},a/b} (\mathbf{x}_{0})$ and $L^{(c_{\text{HT},k}(\cdot),k)}_{\mathcal{K}_{k},a/b} (\mathbf{x}_{0})$, we need to compute the first-order partial derivatives of the mean function.
This is accomplished by using Lemma \ref{lem_cond_exist_CS_recovery_partial_der} for the HT estimator and by using a finite-difference quotient approximation \cite{HeroUniformCRB} for the ML estimator, i.e., 
\begin{equation}
\frac{\partial \mathsf{E}_{\mathbf{x}} \big\{\hat{x}_{\text{ML},k}(\mathbf{y})\big\} }{\partial x_{l}} \approx \frac{ \mathsf{E}_{\mathbf{x}+ \Delta \mathbf{e}_{l}} \big\{\hat{x}_{\text{ML},k}(\mathbf{y})\big\} -\mathsf{E}_{\mathbf{x}} \big\{\hat{x}_{\text{ML},k}(\mathbf{y})\big\} }{\Delta}, 
\end{equation}  
where $\Delta \in \mathbb{R}_{+}$ is a small stepsize and the expectations are calculated using numerical integration.


We generated parameter vectors $\mathbf{x}_{0} = \sqrt{\text{SNR}} \sigma \mathbf{x}_{1}$, where $\mathbf{x}_{1} \in \{0,1\}^{50}$, 
$\supp(\mathbf{x}_{1}) = [S]$, and $\text{SNR}$ 
varies between $0.01$ and $100$.
(The fixed choice $\supp(\mathbf{x}_{0}) = [S]$ is justified by the fact that neither the variances of the ML and HT estimators nor the corresponding variance bounds 
depend on the location of $\supp(\mathbf{x}_{0})$.)
In Fig.~\ref{fig_bounds_1}, we plot the variances $v(\ML_est(\cdot); \mathbf{x}_{0})$ and $v( \hat{\mathbf{x}}_{\text{HT}} (\cdot); \mathbf{x}_0)$ 
(the latter for four different choices of 
$T$ in \eqref{equ_def_thr_func}) along with the corresponding bounds $L_{\text{ML},a/b}(\mathbf{x}_{0})$ and $L_{\text{HT},a/b}(\mathbf{x}_{0})$ 
(cf.\ \eqref{equ_def_lound_asilomar_ML} and \eqref{equ_def_lound_asilomar_ML}), as a function of SNR. 
The variances are obtained by numerical integration from \eqref{equ_var_hard_th_est_ssnm} for the HT estimator and 
from the closed-form expression for the variance of the ML estimator presented in \cite{AlexZvikaJournal}.
It is seen that for SNR larger than about 18 dB, 
all variances and bounds are effectively equal (for the HT estimator, 
this is true if $T$ is not too small). However, in the medium-SNR
range, the variances of the ML and HT estimators are significantly higher than
the corresponding lower bounds. 
We can conclude that there \emph{might} exist estimators with the same mean as that of the ML or HT estimator but a smaller variance.
On the other hand, a positive statement regarding the existence of an estimator that has the same mean as the ML or HT estimator but a \emph{uniformly} lower variance cannot be 
based on our analysis. 


However, for the special case of diagonal estimators for the SSNM, such as the HT estimator, Theorem \ref{thm_diag_bias_min_achiev_var_LMV} and 
Corollary \ref{cor_diag_est_min_achiev_var_LMV_ISIT} make a positive statement about the existence of estimators that 
have \emph{locally} a smaller variance than the HT estimator. In particular, we can use Corollary \ref{cor_diag_est_min_achiev_var_LMV_ISIT} to obtain the LMV 
estimator and corresponding minimum achievable variance at a parameter vector $\mathbf{x}_{0} \in \mathcal{X}_{S}$ for the given bias function of the HT estimator. 

\begin{figure}
\vspace{-1mm}
\centering
\psfrag{SNR}[c][c][.85]{\uput{4mm}[270]{0}{\hspace{0mm}SNR\,[dB]}}
\psfrag{title}[c][c][.9]{\uput{2.5mm}[270]{0}{}}
\psfrag{x_0}[c][c][.9]{\uput{0.3mm}[270]{0}{$0$}}
\psfrag{x_0_001}[c][c][.9]{\uput{0.3mm}[270]{0}{$-30$}}
\psfrag{x_0_01}[c][c][.9]{\uput{0.3mm}[270]{0}{$-20$}}
\psfrag{x_0_1}[c][c][.8]{\uput{0.3mm}[270]{0}{\hspace{-2mm}$-10$}}
\psfrag{x_1}[c][c][.8]{\uput{0.3mm}[270]{0}{\hspace{.3mm}$0$}}
\psfrag{x_10}[c][c][.8]{\uput{0.3mm}[270]{0}{\hspace{.2mm}$10$}}
\psfrag{x_100}[c][c][.8]{\uput{0.3mm}[270]{0}{\hspace{.2mm}$20$}}
\psfrag{y_0}[c][c][.9]{\uput{0.1mm}[180]{0}{$0$}}
\psfrag{y_5}[c][c][.9]{\uput{0.1mm}[180]{0}{$5$}}
\psfrag{y_10}[c][c][.9]{\uput{0.1mm}[180]{0}{$10$}}
\psfrag{y_15}[c][c][.9]{\uput{0.1mm}[180]{0}{$15$}}
\psfrag{y_20}[c][c][.9]{\uput{0.1mm}[180]{0}{$20$}}
\psfrag{y_25}[c][c][.9]{\uput{0.1mm}[180]{0}{$25$}}
\psfrag{y_30}[c][c][.9]{\uput{0.1mm}[180]{0}{$30$}}
\psfrag{y_35}[c][c][.9]{\uput{0.1mm}[180]{0}{$35$}}
\psfrag{y_40}[c][c][.9]{\uput{0.1mm}[180]{0}{$40$}}
\psfrag{y_45}[c][c][.9]{\uput{0.1mm}[180]{0}{$45$}}
\psfrag{y_50}[c][c][.9]{\uput{0.1mm}[180]{0}{$50$}}
\psfrag{y_55}[c][c][.9]{\uput{0.1mm}[180]{0}{$55$}}
\psfrag{variance}[c][c][.9]{\uput{3mm}[90]{0}{variance/Barankin bound}}
\psfrag{bML}[l][l][0.8]{bound on $v((\cdot);\mathbf{x}_0)$}
\psfrag{ML}[l][l][0.8]{$v((\cdot);\mathbf{x}_0)$}
\psfrag{data1}[l][l][0.8]{$v(\hat{\mathbf{x}}_{\text{LS}}(\cdot); \mathbf{x}_0)$}
\psfrag{data2}[l][l][0.8]{$L_{\text{HT}}(\mathbf{x}_{0})$, $T \!\rmv=\! 0$}
\psfrag{data3}[l][l][0.8]{$v(\hat{\mathbf{x}}_{\text{HT}}(\cdot); \mathbf{x}_0)$, $T \!\rmv=\! 2$}
\psfrag{data5}[l][l][0.8]{$v(\hat{\mathbf{x}}_{\text{HT}}(\cdot); \mathbf{x}_0)$, $T \!\rmv=\! 3$}
\psfrag{data7}[l][l][0.8]{$v(\hat{\mathbf{x}}_{\text{HT}}(\cdot); \mathbf{x}_0)$, $T \!\rmv=\! 4$}
\psfrag{data4}[l][l][0.8]{$L_{\text{HT}}(\mathbf{x}_{0})$, $T \!\rmv=\! 2$}
\psfrag{data6}[l][l][0.8]{$L_{\text{HT}}(\mathbf{x}_{0})$, $T \!\rmv=\! 3$}
\psfrag{data8}[l][l][0.8]{$L_{\text{HT}}(\mathbf{x}_{0})$, $T \!\rmv=\! 4$}
\psfrag{data9}[l][l][0.8]{$v(\hat{\mathbf{x}}_{\text{ML}}(\cdot); \mathbf{x}_0)$}
\psfrag{MLarrow}[l][l][0.8]{ML}
\psfrag{T=3}[l][l][0.8]{\uput{0mm}[-10]{0}{\hspace{1mm}$T \!\rmv=\! 4$}}
\psfrag{T=4}[l][l][0.8]{$T \!\rmv=\! 3$}
\psfrag{T=2}[l][l][0.8]{\uput{.5mm}[0]{0}{\vspace*{3mm} \!\!\!\!$T \!\rmv=\! 2$}}
\psfrag{T=0}[l][l][0.8]{\uput{0mm}[0]{0}{\hspace{-10mm}$T \!\rmv=\! 0\,\ist$(LS)}}
\psfrag{oracle}[l][l][0.8]{\uput{0mm}[0]{0}{$S \sigma^{2}$}}
\centering
\includegraphics[height=8cm,width=15cm]{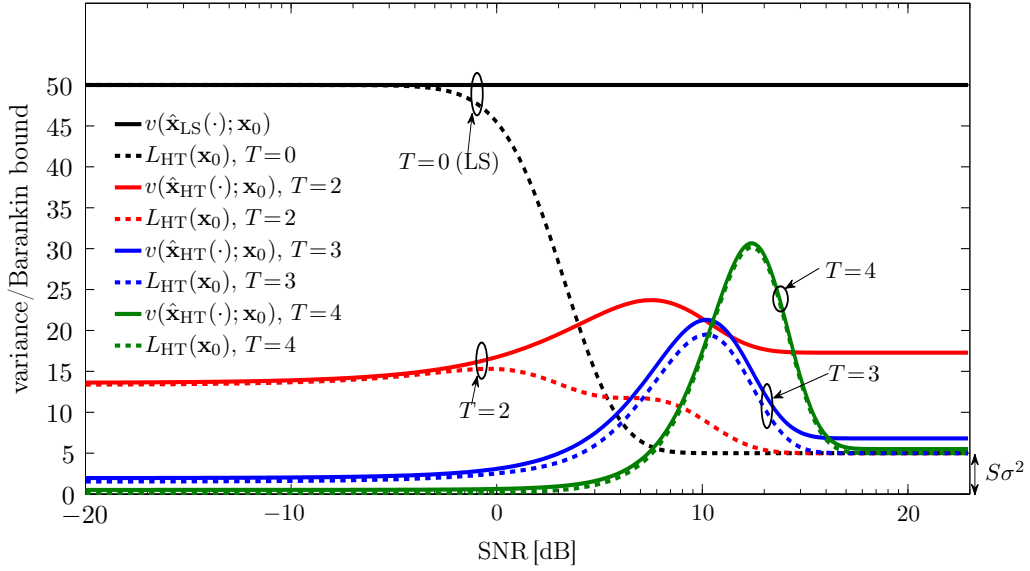}
\vspace*{1mm}
  \caption{Variance of the HT estimator, $v(\hat{\mathbf{x}}_{\text{HT}}(\cdot); \mathbf{x}_{0})$, for different $T$ (solid lines) 
  and corresponding minimum achievable variance (Barankin bound) $L_{\text{HT}}(\mathbf{x}_{0})$ (dashed lines)
  versus the SNR, for the SSNM with $N\!=\!50$, $S\!=\!5$, and $\sigma^2\!=\!1$.} 
\label{fig_barankin_bound_HT_ISIT}
\vspace*{3.5mm}
\end{figure}

As before, we generated parameter vectors $\mathbf{x}_{0} = \sqrt{\text{SNR}} \sigma \mathbf{x}_{1}$, where $\mathbf{x}_{1} \in \{0,1\}^{50}$, 
$\supp(\mathbf{x}_{1}) = [S]$, and $\text{SNR}$ 
varies between $0.01$ and $100$.
(The fixed choice $\supp(\mathbf{x}_{0}) = [S]$ is justified by the fact that neither the variance of the HT estimators nor the corresponding minimum achievable variance 
depend on the location of $\supp(\mathbf{x}_{0})$.)
In Fig.~\ref{fig_barankin_bound_HT_ISIT}, we plot the variances  $v(\hat{\mathbf{x}}_{\text{HT}}(\cdot); \mathbf{x}_{0})$ (obtained by numerical integration using \eqref{equ_var_hard_th_est_ssnm}) 
for four different choices of $T$ in \eqref{equ_def_thr_func}, as a function of SNR. We also plot the
corresponding minimum achievable variance (Barankin bound) $L_{\text{HT}}(\mathbf{x}_{0}) \triangleq \sum_{k \in [N]} L_{\mathcal{M}_{k}}$ applying to estimators for the 
SSNM that have the same bias (and mean function) as the HT estimator. Note that, by definition, $L_{\text{HT}}(\mathbf{x}_{0}) \triangleq \sum_{k \in [N]} L_{\mathcal{M}_{k}}$ is the maximally 
tight lower bound on the variance of any estimator that has the same bias (and mean function) as the HT estimator. This bound is achieved by the LMV estimator given component-wise by \eqref{equ_expr_LMV_diag_est_ISIT_SSNM} with $\hat{x}_{k}(y_{k}) = \hat{x}_{\text{HT},k}(y_{k})$. 
The minimum achievable variance $L_{\text{HT}}(\mathbf{x}_{0})$ is obtained by summing the minimum achievable variances $L_{\mathcal{M}_{k}}$ of the scalar minimum variance problems $\mathcal{M}_{k}$ which coincide with $\mathcal{M}_{\text{SSNM}}$ using the parameter function $g(\mathbf{x}) = x_{k}$ and 
prescribed bias $c(\mathbf{x}) = b(\hat{x}_{\text{HT},k}(\mathbf{y}),\mathbf{x})$. 
We obtained $L_{\mathcal{M}_{k}}$ by an application of Corollary \ref{cor_diag_est_min_achiev_var_LMV_ISIT} since the estimator $\hat{x}_{\text{HT},k}(\mathbf{y})$ is diagonal and has finite variance at every $\mathbf{x}_{0} \in \mathcal{X}_{S}$. 

It is seen that for small $T$ (in particular, for $T\!=\!0$ where the HT estimator reduces to the LS estimator), the Barankin bound $L_{\text{HT}}(\mathbf{x}_{0})$ is 
significantly below the corresponding variance curve. However, as $T$ increases, the gap between variance and Barankin bound $L_{\text{HT}}(\mathbf{x}_{0})$ becomes smaller;
in particular, the two curves are already indistinguishable for $T\!=\rmv 4$. 
For high SNR, the Barankin bound $L_{\text{HT}}(\mathbf{x}_{0})$ converges to $S \sigma^{2}$ for any value of $T$; this equals 
the variance of an oracle estimator that knows the support 
of $\mathbf{x}$.

\chapter{The Sparse Parametric Covariance Model}
\label{chap_SCM}
\section{Introduction}
\label{sec_intro_SPCM}

In the previous chapter, we considered the problem of estimating a sparse parameter vector using a noisy and linearly distorted (by the system matrix $\mathbf{H}$) 
observation. Now, we consider the fundamentally different problem of estimating a sparse parameter vector which determines the covariance matrix of 
a Gaussian signal vector of which a noisy version is observed. 

More specifically, we consider a Gaussian signal vector $\mathbf{s} \in \mathbb{R}^{M}$, i.e., $\mathbf{s} \sim \mathcal{N}({\bm \mu}, \mathbf{C})$ 
with the psd covariance matrix $\mathbf{C} \in \mathbb{R}^{M \times M}$, embedded in white Gaussian 
noise $\mathbf{n} \sim \mathcal{N}(\mathbf{0}, \sigma^{2} \mathbf{I}_{M})$. The observed vector is 
\begin{equation}
\label{equ_obs_model_SPCM}
\mathbf{y} = \mathbf{s} + \mathbf{n},
\end{equation}
where $\mathbf{s}$ and $\mathbf{n}$ are independent and the signal mean ${\bm \mu} \in \mathbb{R}^{M}$ and noise variance $\sigma^{2} > 0$ are known. 
These assumptions imply that the observation $\mathbf{y}$ is a Gaussian random vector with mean ${\bm \mu}$ and covariance matrix $\mathbf{C} + \sigma^{2} \mathbf{I}$, 
i.e., $\mathbf{y} \sim \mathcal{N}({\bm \mu}, \mathbf{C}+ \sigma^{2} \mathbf{I})$. 
In what follows, we assume that ${\bm \mu} = \mathbf{0}$ since a nonzero ${\bm \mu}$ can always be removed from $\mathbf{s}$ by subtracting it from $\mathbf{y}$. 

The signal covariance matrix $\mathbf{C} \in \mathbb{R}^{M \times M}$ is unknown. We will parametrize it in a linear manner by 
\begin{equation} 
\label{equ_cov_param_SPCM}
\mathbf{C} = \mathbf{C}(\mathbf{x}) = \sum_{k \in [N]} x_{k} \mathbf{C}_{k}, 
\end{equation} 
with unknown nonrandom coefficients $x_{k} \in \mathbb{R}_{+}$ and known psd ``basis matrices'' $\mathbf{C}_{k} \in \mathbb{R}^{M \times M}$, i.e., $\mathbf{C}_{k} \geq \mathbf{0}$ with 
ranks $r_{k} \triangleq \rank(\mathbf{C}_{k})$. It what follows, we will use
\begin{equation}
\label{equ_def_signal_space_dimension_SPCM}
R \triangleq \sum_{k \in [N]} r_{k}
\end{equation} 
to denote the maximum overall dimension of the signal space. We use the term ``maximum overall dimension'' since the dimension of the signal space 
is not equal to $R$ in general. It is only equal to $R$  if the basis matrices $\mathbf{C}_{k}$ fulfill additional requirements. One such additional requirement would be that the basis matrices span disjoint orthogonal subspaces of $\mathbb{R}^{M}$.
Due to the parameterization \eqref{equ_cov_param_SPCM}, the estimation of the covariance matrix $\mathbf{C}$ reduces to the estimation of the coefficient vector $\mathbf{x} \triangleq \big(x_{1},\ldots,x_{N}\big)^{T} \in \mathbb{R}_{+}^{N}$. 

Let
\begin{equation}
\label{equ_def_cov_matrix_observation_SPCM}
\widetilde{\mathbf{C}}(\mathbf{x})\triangleq \mathbf{C}(\mathbf{x}) + \sigma^{2} \mathbf{I}
\end{equation} 
denote the covariance matrix of the Gaussian observation vector $\mathbf{y}$ (cf.\ \eqref{equ_obs_model_SPCM}), i.e., $\mbox{cov}\{\mathbf{y}\} = \widetilde{\mathbf{C}}(\mathbf{x})$. Due to the assumption 
that $\sigma^{2} > 0$, one can verify by Lemma \ref{lem_eigvalues_sum_symmetric_matrices} that $\| \widetilde{\mathbf{C}}(\mathbf{x}) \|_{2} \geq \sigma^{2}$, i.e., the matrix $\widetilde{\mathbf{C}}(\mathbf{x})$ is positive definite 
($\widetilde{\mathbf{C}}(\mathbf{x})>\mathbf{0}$) for every parameter vector $\mathbf{x} \in \mathbb{R}_{+}^{N}$. Another useful fact which can be verified easily (see, e.g., \cite{golub96}) is that it holds 
\begin{equation} 
\label{equ_inverse_norm_leq_sigma_2}
\| \widetilde{\mathbf{C}}^{-1}(\mathbf{x}) \|_{2} \leq \sigma^{-2}
\end{equation} 
for every $\mathbf{x} \in \mathbb{R}_{+}^{N}$.

The parameter vector $\mathbf{x}$ is assumed to be $S$-sparse, i.e.,
\begin{equation}
\mathbf{x} \in \mathcal{X}_{S,+} \triangleq \mathcal{X}_{S} \cap \mathbb{R}_{+}^{N} = \left\{ \mathbf{x}' \in \mathbb{R}_{+}^{N} \big| \| \mathbf{x}' \|_{0} \leq S \right\}.  
\end{equation} 
While the sparsity degree $S \in [N]$ is assumed known, the support of $\mathbf{x}$, i.e., $\supp(\mathbf{x})$ is 
unknown. Typically, $S \ll N$. 
 
We now define the sparse parametric covariance model (SPCM) $\mathcal{E}_{\text{SPCM}}$ as the estimation problem given as 
\begin{equation}
\label{equ_def_SPCM}
\mathcal{E}_{\text{SPCM}} = \left( \mathcal{X}_{S,+},f_{\text{SPCM}}(\mathbf{y}; \mathbf{x}), g(\mathbf{x})  \right), 
\end{equation} 
with an arbitrary parameter function $g(\cdot): \mathcal{X}_{S,+} \rightarrow \mathbb{R}$ and the statistical model 
\begin{equation}
\label{equ_def_SPCM_stat_model}
f_{\text{SPCM}}(\mathbf{y}; \mathbf{x}) = \frac{1}{(2 \pi)^{M/2} [\detm{\widetilde{\mathbf{C}}(\mathbf{x})}]^{1/2}}\exp \left(- \frac{1}{2} \mathbf{y}^{T}  \widetilde{\mathbf{C}}^{-1}(\mathbf{x}) \mathbf{y} \right). 
\end{equation} 
The statistical model of the SPCM belongs to an exponential family (see Section \ref{sec_exp_family}). In particular, it is obtained for the sufficient statistic 
${\bm \Phi}(\mathbf{y}) = -\frac{1}{2} \vectorize{\mathbf{y}\mathbf{y}^{T}}$, 
parameter function $\mathbf{u}(\mathbf{x}) = \vectorizesymb\big\{ ( \mathbf{C}(\mathbf{x}) + \sigma^{2} \mathbf{I} )^{-1}\big\}$, weight function $h(\mathbf{y}) = \frac{1}{\sqrt{(2 \pi)^{M}}}$, and cumulant function 
$A^{( {\bm \Phi})}(\mathbf{x}) = \frac{1}{2} \log \big(\detm{\mathbf{C}(\mathbf{x})+\sigma^{2} \mathbf{I}}\big)$. 

Note that by contrast to the SLM, the entries of the parameter vector $\mathbf{x} \in \mathcal{X}_{S,+}$ are non-negative. Another difference from the SLM is that we allow for a general parameter function 
$g(\cdot): \mathcal{X}_{S,+} \rightarrow \mathbb{R}$ instead of fixing the parameter function to $g(\mathbf{x}) = x_{k}$ as for the SLM. 
However, according to Theorem \ref{thm_equ_bias_param_function}, 
this latter difference is only a notational one as long as we allow for arbitrary prescribed bias functions $c(\cdot): \mathcal{X}_{S,+} \rightarrow \mathbb{R}$ 
for the minimum variance problems arising from the SLM $\mathcal{E}_{\text{SLM}}$.

Since the basis matrices $\mathbf{C}_{k}$ are assumed to be psd, we can factor them in the form $\mathbf{C}_{k} = \mathbf{H}_{k} \mathbf{H}_{k}^{T}$,
where the factors $\mathbf{H}_{k} \in \mathbb{R}^{M \times r_{k}}$ are not unique \cite{golub96,Horn91}. 
One possible choice of $\mathbf{H}_{k}$ is obtained via the thin EVD $\mathbf{C}_{k} = \mathbf{U} \mathbf{\Sigma} \mathbf{U}^{T}$ by setting $\mathbf{H}_{k} = \mathbf{U} \sqrt{ \mathbf{\Sigma}}$. 
If we consider a specific factorization of the basis matrices and stack the resulting factors $\mathbf{H}_{k}$ into the matrix 
\begin{equation} 
\label{eq_def_stacked_system_matrix_SPCM}
\mathbf{H} \triangleq \begin{pmatrix} \mathbf{H}_{1} & \mathbf{H}_{2} & \ldots & \mathbf{H}_{N} \end{pmatrix}  \in \mathbb{R}^{M \times  R},
\end{equation}
we can rewrite the parameterization \eqref{equ_cov_param_SPCM} as 
\begin{equation} 
\label{equ_cov_param_SPCM_1}
\mathbf{C} = \mathbf{C}(\mathbf{x}) = \mathbf{H} \mathbf{D}(\mathbf{x}) \mathbf{H}^{T}.
\end{equation}
Here, $\mathbf{D}(\mathbf{x}) \in \mathbb{R}^{R \times R}$ is a diagonal matrix which has on its main diagonal the entries of $\mathbf{x} \in \mathbb{R}^{N}$, where 
each entry $x_{k}$ appears $r_{k}$ times, i.e., 
\begin{equation}
\label{equ_diag_matrix_param_vec_SPCM}
\mathbf{D}(\mathbf{x}) = \begin{pmatrix} x_{1} \mathbf{I}_{r_{1}} & \mathbf{0} & \ldots & \mathbf{0} \\ 
\mathbf{0} & x_{2} \mathbf{I}_{r_{2}} & \ddots  &  \vdots \\ 
\vdots & \ddots& \ddots &   \mathbf{0} \\ 
\mathbf{0} & \ldots&  \mathbf{0} &  x_{N} \mathbf{I}_{r_{N}} \\ 
\end{pmatrix}.
\end{equation} 
A partial characterization of the set of basis matrices can be made based on  
\begin{definition} 
\label{def_spcm_basis_matrices_coherence}
Given the SPCM $\mathcal{E}_{\emph{SPCM}}$, we say that the set of basis matrices $\{ \mathbf{C}_{k} \in \mathbb{R}^{M \times M} \}_{k \in [N]}$ has RIP of order $K$ with constant $\delta_{K} \in [0,1]$ if 
there exist factorizations $\mathbf{C}_{k} = \mathbf{H}_{k} \mathbf{H}_{k}^{T}$ 
such that the singular values of any matrix $\mathbf{H}^{(\mathcal{I})} = \begin{pmatrix} \mathbf{H}_{i_{1}} & \mathbf{H}_{i_{2}} & \ldots & \mathbf{H}_{i_{K}} \end{pmatrix} 
\in \mathbb{R}^{M \times \sum_{k \in \mathcal{I}} r_{k}}$, where 
$\mathcal{I} = \{i_{1}, \ldots, i_{K}Ê\} \subseteq [N]$ 
is an arbitrary set of $K$ different indices, are located in the interval $[1- \delta_{K}, 1+ \delta_{K}]$ and moreover any matrix $\mathbf{H}^{(\mathcal{I})}$ has full column rank, i.e., 
$\rank(\mathbf{H}^{(\mathcal{I})})=\sum_{k \in \mathcal{I}} r_{k}$. 
\end{definition}
 
The SPCM and estimation of $\mathbf{x}$ are relevant, e.g., in time-frequency analysis \cite{fla-book2,hla-book}, where the basis matrices $\mathbf{C}_{k}$ 
correspond to disjoint time-frequency regions and $x_{k}$ represents the mean signal energy in the $k$th time-frequency region. 
Such a sparsity-based time-frequency analysis may be used, e.g., for cognitive radio scene analysis \cite{CognitiveRadio2005}.
Another application that fits our scope would be compressive spectrum estimation \cite{jung-ssp09,jung-icassp09}, where classical (non-stationary) spectrum 
estimation techniques are combined with novel CS methods. 

An important special case of the SPCM is obtained if the basis matrices $\mathbf{C}_{k}$ appearing in the parameterization \eqref{equ_cov_param_SPCM} are orthogonal projection matrices on orthogonal subspaces of 
$\mathbb{R}^{M}$, i.e., 
\begin{equation} 
\label{equ_conds_basis_matrices_SDPCM}
\mathbf{C}_{k} \mathbf{C}_{k'} = \delta_{k,k'} \mathbf{C}_{k}\mbox{ , } \quad \forall k,k' \in [N]. 
\end{equation} 
Note that condition \eqref{equ_conds_basis_matrices_SDPCM} implies that the basis matrices satisfy the RIP of any order $K \in [N]$ with constant $\delta_{K}=0$. 
The so-obtained special case of the SPCM will be termed the \emph{sparse diagonalizable parametric covariance model} (SDPCM). This naming is due to the 
fact that for the SDPCM, i.e., when the basis matrices $\mathbf{C}_{k}$ in \eqref{equ_cov_param_SPCM} satisfy \eqref{equ_conds_basis_matrices_SDPCM}, 
the covariance matrix $\mathbf{C}(\mathbf{x})$ given by \eqref{equ_cov_param_SPCM} can be diagonalized by a signal transformation $\mathbf{s}' = \mathbf{U} \mathbf{s}$, with 
an orthonormal matrix $\mathbf{U} \in \mathbb{R}^{M \times M}$, i.e., $\mathbf{U}^{T} \mathbf{U} = \mathbf{I}$, that does not depend on the parameter vector $\mathbf{x}$. 
We will denote by $\mathcal{E}_{\text{SDPCM}}$ the special case of the estimation problem $\mathcal{E}_{\text{SPCM}}$ in \eqref{equ_def_SPCM} 
with basis matrices satisfying \eqref{equ_conds_basis_matrices_SDPCM}. More specifically, 
\begin{equation}
\label{equ_def_SDPCM}
\mathcal{E}_{\text{SDPCM}} = \left( \mathcal{X}_{S,+},f_{\text{SDPCM}}(\mathbf{y}; \mathbf{x}), g(\mathbf{x})  \right), 
\end{equation} 
with an arbitrary parameter function $g(\cdot): \mathcal{X}_{S,+} \rightarrow \mathbb{R}$ and the statistical model 
\begin{equation}
\label{equ_def_SDPCM_stat_model}
f_{\text{SDPCM}}(\mathbf{y}; \mathbf{x}) = \frac{1}{(2 \pi)^{M/2} [\detm{\widetilde{\mathbf{C}}(\mathbf{x})}]^{1/2}}\exp \left(- \frac{1}{2} \mathbf{y}^{T}  \widetilde{\mathbf{C}}^{-1}(\mathbf{x}) \mathbf{y} \right),  
\end{equation} 
where $\widetilde{\mathbf{C}}(\mathbf{x})$ is defined as in \eqref{equ_def_cov_matrix_observation_SPCM} but with basis matrices $\mathbf{C}_{k}$ satisfying \eqref{equ_conds_basis_matrices_SDPCM}. 
 
The key results of this chapter have been presented in part in \cite{AlexSebICASSP2011}. 

\section{The $\mathcal{D}$-restricted SPCM} 
\label{sec_D_restricted_SPCM}

As already observed in \cite{Duttweiler73b}, the RKHS-based analysis of minimum variance estimation for estimation problems where the parameter determines 
the covariance matrix of the observation (e.g., the SPCM) is more involved than for those problems where the parameter determines the mean (e.g., the SLM). 
In particular, any minimum variance problem $\mathcal{M}_{\text{SPCM}}=\left(\mathcal{E}_{\text{SPCM}},c(\cdot) \equiv 0,\mathbf{x}_{0}Ê\right)$ 
associated with the SPCM and parameter vector $\mathbf{x}_{0} \in \mathcal{X}_{S,+}$ does not fulfill Postulate \ref{assumption_RKHS_minvarproblem}. 
This precludes a direct RKHS approach to the SPCM.

One specific reason for the minimum variance problem $\mathcal{M}_{\text{SPCM}}$ to fail in satisfying Postulate \ref{assumption_RKHS_minvarproblem} is 
that the space $\mathcal{L}(\mathcal{M}_{\text{SPCM}})$ (cf. \eqref{equ_def_hilbert_space_rhos}) 
with the inner product $\langleÊ\,  \cdot \,  , \,  \cdot \, \rangle_{\text{RV}}$ (see \eqref{equ_def_inner_prod_RV}) fails to be a Hilbert space. 
An ad-hoc approach in order to facilitate an RKHS-based analysis of the minimum variance problem $\mathcal{M}_{\text{SPCM}}$ could be to reduce the space 
$\mathcal{L}(\mathcal{M}_{\text{SPCM}})$ (consisting of linear combinations of the elementary estimator functions $\hat{g}_{\mathbf{x}}(\cdot) = \rho_{\mathcal{M}}(\cdot,\mathbf{x})$) 
to obtain a Hilbert space 
$\mathcal{L}'(\mathcal{M}_{\text{SPCM}}) \subseteq \mathcal{L}(\mathcal{M}_{\text{SPCM}})$. 
However, the resulting RKHS associated with $\mathcal{L}'(\mathcal{M}_{\text{SPCM}})$ depends on the 
completely arbitrary reduction of $\mathcal{L}(\mathcal{M}_{\text{SPCM}})$ to $\mathcal{L}'(\mathcal{M}_{SPCM})$, and a convenient analytical description of the 
resulting RKHS can be hardly expected in general. Therefore, as discussed presently, we will pursue a different RKHS approach to the minimum variance problem 
$\mathcal{M}_{\text{SPCM}}=\left(\mathcal{E}_{\text{SPCM}},c(\cdot) \equiv 0,\mathbf{x}_{0}Ê\right)$.

Specifically, in order to derive lower bounds on the minimum achievable variance $L_{\mathcal{M}_{\text{SPCM}}}$, we can further restrict the parameter set of $\mathcal{E}_{\text{SPCM}}$ 
to obtain the so called ``$\mathcal{D}$-restricted'' SPCM 
\begin{equation}
\mathcal{E}^{(\mathcal{D})}_{\text{SPCM}} \triangleq \left( \mathcal{D},f_{\text{SPCM}}(\mathbf{y}; \mathbf{x}), g(\mathbf{x})  \right)
\end{equation}
with a certain set $\mathcal{D} \subseteq \mathcal{X}_{S,+}$. 
Indeed, consider the minimum variance problem $\mathcal{M}_{\text{SPCM}}=\left(\mathcal{E}_{\text{SPCM}},c(\cdot) \equiv 0, \mathbf{x}_{0}Ê\right)$ associated with the SPCM $\mathcal{E}_{\text{SPCM}}$ and 
the ``$\mathcal{D}$-restricted'' minimum variance problem 
\begin{equation} 
\mathcal{M}_{\mathcal{D},\text{SPCM}}\triangleq \mathcal{M}_{\text{SPCM}}\big|_{\mathcal{D}} = \big( \mathcal{E}^{(\mathcal{D})}_{\text{SPCM}} , c(\cdot) \equiv 0, \mathbf{x}_{0} \big)
\end{equation}
associated with the $\mathcal{D}$-restricted SPCM $\mathcal{E}^{(\mathcal{D})}_{\text{SPCM}}$ and the same parameter vector $\mathbf{x}_{0}$ as used for $\mathcal{M}_{\text{SPCM}}$. 
We do not require the set $\mathcal{D}$ to contain $\mathbf{x}_{0}$, i.e., we may have $\mathbf{x}_{0} \notin \mathcal{D}$. 
According to Theorem \ref{thm_par_set_reduction_classic_est_mve} (cf.\ also Section \ref{sec_reducing_par_set_RKHS}), we have then the inequality 
\begin{equation} 
\label{equ_lower_bound_spcm_restr_spcm}
L_{\mathcal{M}_{\text{SPCM}}} \geq L_{\mathcal{M}_{\mathcal{D}, \text{SPCM}}}.
\end{equation} 

The idea is now to choose the set $\mathcal{D}\subseteq \mathcal{X}_{S,+}$ such that $\mathcal{M}_{\mathcal{D},\text{SPCM}}$ fulfills Postulate \ref{assumption_RKHS_minvarproblem}, so that
 the RKHS $\mathcal{H}(\mathcal{M}_{\mathcal{D},\text{SPCM}})$ as in Definition \ref{def_RKHS_min_var_problem} exists for the minimum variance problem $\mathcal{M}_{\mathcal{D},\text{SPCM}}$. 
Based on \eqref{equ_lower_bound_spcm_restr_spcm}, we can then find 
lower bounds for $L_{\mathcal{M}_{\text{SPCM}}}$ by deriving lower bounds on $L_{\mathcal{M}_{\mathcal{D}, \text{SPCM}}}$ using the RKHS theory for $\mathcal{H}(\mathcal{M}_{\mathcal{D},\text{SPCM}})$. 

Postulate \ref{assumption_RKHS_minvarproblem} can be written for the SPCM as (see \eqref{equ_def_SPCM_stat_model})
\vspace*{4mm}
 \begin{align}
 \label{equ_cond_min_var_RKHS_SPCM}
& \mathsf{E}_{\mathbf{x}_{0}} \Bigg \{  \Bigg( \frac{f_{\text{SPCM}}(\mathbf{y}; \mathbf{x})}{f_{\text{SPCM}}(\mathbf{y}; \mathbf{x}_{0})} \Bigg)^{2} \Bigg\}  =   \int_{\mathbf{y}} \frac{f_{\text{SPCM}}(\mathbf{y}; \mathbf{x}) f_{\text{SPCM}}(\mathbf{y}; \mathbf{x})}{f_{\text{SPCM}}(\mathbf{y}; \mathbf{x}_{0})} d\mathbf{y}\nonumber \\[5mm]
&  \hspace*{10mm}=  \frac{1}{(2 \pi)^{M/2}} \frac{ \big[ \detmb{\widetilde{\mathbf{C}}(\mathbf{x}_{0})}Ê\big]^{1/2}}{\detmb{\widetilde{\mathbf{C}}(\mathbf{x})}} \int_{\mathbf{y}} 
\exp \left( - \frac{1}{2} \mathbf{y}^{T} \big[ 2 \widetilde{\mathbf{C}}^{-1}(\mathbf{x}) - \widetilde{\mathbf{C}}^{-1}(\mathbf{x}_{0}) \big]Ê\mathbf{y} \right)  d\mathbf{y} < \infty. \\[-4mm]
\nonumber 
\end{align}
Since by \eqref{equ_def_cov_matrix_observation_SPCM} we have $\widetilde{\mathbf{C}}(\mathbf{x}) > \mathbf{0}$ and therefore $\detmb{\widetilde{\mathbf{C}}(\mathbf{x})}> 0$ 
for every $\mathbf{x} \in \mathcal{X}_{S,+}$, we can conclude from \eqref{equ_cond_min_var_RKHS_SPCM} that in order to satisfy Postulate \ref{assumption_RKHS_minvarproblem}, it is necessary and sufficient to choose the set $\mathcal{D}$ such that $\big(2 \widetilde{\mathbf{C}}^{-1}(\mathbf{x}) - \widetilde{\mathbf{C}}^{-1}(\mathbf{x}_{0})\big)$ is positive definite, i.e.,
\vspace*{3mm}
\begin{align}
\label{equ_nec_suff_cond_D_SPCM_kernel_exists}
 2 \widetilde{\mathbf{C}}^{-1}(\mathbf{x}) - \widetilde{\mathbf{C}}^{-1}(\mathbf{x}_{0}) > \mathbf{0} \quad\quad \forall \mathbf{x} \in \mathcal{D}. \\[-4mm]
\nonumber
\end{align}

We will now give a sufficient condition for the set $\mathcal{D}$ to fulfill \eqref{equ_nec_suff_cond_D_SPCM_kernel_exists}.  
To that end, the following two technical results will prove handy. 
Remember that given a square matrix $\mathbf{H} \in \mathbb{R}^{M \times M}$, we denote by $\| \mathbf{H} \|_{2}$ its matrix $2$-norm defined as $\| \mathbf{H} \|_{2} \triangleq \sup\limits_{\| \mathbf{x} \|_{p}=1} \| \mathbf{H} \mathbf{x} \|_{2}$.
\begin{lemma}
\label{lem_matrix_inversion_perturb}
Consider an invertible matrix $\mathbf{A} \in \mathbb{R}^{M \times M}$ and a matrix $\mathbf{E} \in \mathbb{R}^{MÊ\times M}$. 
Then if $\| \mathbf{E} \|_{2} < \frac{1}{\| \mathbf{A}^{-1} \|_{2}}$, we have that $\mathbf{A} + \mathbf{E}$ is invertible and
\vspace*{3mm}
\begin{align} 
\label{equ_matrix_inversion_perturb}
\big\| \left( \mathbf{A} + \mathbf{E} \right)^{-1} - \mathbf{A}^{-1} \big\|_{2} \leq \frac{\| \mathbf{A}^{-1}Ê\|^{2}_{2} \| \mathbf{E} \|_{2}}{1 - \| \mathbf{A}^{-1} \|_{2} \| \mathbf{E} \|_{2}}. \\[-4mm] 
\nonumber
\end{align} 
\end{lemma} 

\begin{proof}
The statement follows by combining \cite[Theorem 2.3.4]{golub96} with the fact that for two matrices $\mathbf{N}, \mathbf{M} \in \mathbb{R}^{M \times M}$ we have 
$\| \mathbf{N} \mathbf{M} \|_{2} \leq \| \mathbf{N} \|_{2} \| \mathbf{M} \|_{2}$ \cite[p.\ 55]{golub96}.  
\end{proof}

In what follows, we denote the $k$th largest eigenvalue of a symmetric matrix $\mathbf{C} \in \mathbb{R}^{L \times L}$, with $k \in [L]$, by $\lambda_{k}(\mathbf{C})$, i.e., 
\vspace*{1mm}
\begin{align}
\lambda_{L}(\mathbf{C})\leq \ldots \leq \lambda_{2}(\mathbf{C}) \leq \lambda_{1}(\mathbf{C}). \\[-5mm]
\nonumber 
\end{align}  
The eigenvalues of a symmetric matrix are nonnegative (positive) if and only if the matrix is positive semi-definite (positive definite) \cite{golub96,Horn85}. 
It can be verified that $\| \mathbf{C} \|_{2} = \max \big\{ \big| \lambda_{1}(\mathbf{C})\big| , \big| \lambda_{L}(\mathbf{C})\big| \big\}$ \cite[p.\ 394]{golub96} and 
in turn 
\begin{equation} 
\label{equ_any_eigval_lower_mag_norm_2}
| \lambda_{k}(\mathbf{C})| \leq \| \mathbf{C} \|_{2}.
\end{equation}
Moreover, for any matrix $\mathbf{C} \in \mathbb{R}^{L \times L}$ it holds that \cite[p.\ 310]{golub96} 
\begin{equation}
\label{equ_rel_eigvals_det}
\detmb{ \mathbf{C} } = \prod_{k \in [L]} \lambda_{k}(\mathbf{C}).    
\end{equation}

\begin{lemma}
\emph{(\hspace*{-1.5mm}\cite[Theorem 8.1.5]{golub96})}
\label{lem_eigvalues_sum_symmetric_matrices} 
Consider two symmetric matrices $\mathbf{A},\mathbf{E} \in \mathbb{R}^{L \times L}$, i.e., $\mathbf{A}^{T} = \mathbf{A}$ and $\mathbf{E}^{T} =\mathbf{E}$.
We then have the double inequality 
\vspace*{2mm}
\begin{align} 
\label{equ_double_ineq_eigvals_sum}
\lambda_{k}(\mathbf{A}) + \lambda_{L}(\mathbf{E}) \leq \lambda_{k}(\mathbf{A}+ \mathbf{E}) \leq \lambda_{k}(\mathbf{A}) + \lambda_{1}(\mathbf{E})\mbox{,} \quad \mbox{for every $k \in [L]$.} \\[-5mm]
\nonumber
\end{align} 
\end{lemma} 

\begin{proof}
\cite{golub96} 
\end{proof} 

The next result shows that $\mathcal{D}$ can be chosen as the intersection of $\mathcal{X}_{S,+}$ and a ball $\mathcal{B}(\mathbf{x}_{c},r)$, 
with a suitable center $\mathbf{x}_{c} \in \mathcal{X}_{S,+}$ and a sufficiently small radius $r$:
\begin{theorem}
\label{thm_cond_d_restr_SPCM_ball}
Consider the minimum variance problem $\mathcal{M}_{\emph{SPCM}}=\left(\mathcal{E}_{\emph{SPCM}},c(\cdot) \equiv 0, \mathbf{x}_{0}Ê\right)$. 
Any set $\mathcal{D}$ given as  
\vspace*{2mm}
\begin{align}
\label{equ_def_D_thm_cond_d_restr_SPCM_ball}
\mathcal{D} = \mathcal{X}_{S,+} \cap \mathcal{B}(\mathbf{x}_{c},r), \\[-5mm]
\nonumber
\end{align} 
with $\mathbf{x}_{c} \in \mathcal{X}_{S,+}$ such that 
\vspace*{2mm}
\begin{align} 
\label{equ_def_D_thm_cond_d_restr_SPCM_ball_valid_x_c}
 2\widetilde{\mathbf{C}}^{-1}(\mathbf{x}_{c}) - \widetilde{\mathbf{C}}^{-1}(\mathbf{x}_{0})  > \mathbf{0}, \\[-5mm]
\nonumber
\end{align} 
and a sufficiently small radius $r >0$, satisfies \eqref{equ_nec_suff_cond_D_SPCM_kernel_exists}. Two specific choices for $\mathbf{x}_{c}$ which satisfy \eqref{equ_def_D_thm_cond_d_restr_SPCM_ball_valid_x_c} 
are $\mathbf{x}_{c} = \mathbf{x}_{0}$ and $\mathbf{x}_{c} = \mathbf{0}$. 
\end{theorem}

\begin{proof}
Let us consider a parameter vector $\mathbf{x}_{c} \in \mathcal{X}_{S,+}$ such that \eqref{equ_def_D_thm_cond_d_restr_SPCM_ball_valid_x_c} is satisfied.
We can then write $\widetilde{\mathbf{C}}(\mathbf{x}) =\widetilde{\mathbf{C}}(\mathbf{x}_{c}) + \mathbf{E}(\mathbf{x})$, where 
\begin{equation} 
\label{equ_proof_suff_cond_D_E_x_SPCM}
\mathbf{E}(\mathbf{x}) \triangleq \widetilde{\mathbf{C}}(\mathbf{x})  - \widetilde{\mathbf{C}}(\mathbf{x}_{c})  =\sum_{k \in [N]} \mathbf{C}_{k} \bar{x}_{k}
\end{equation} 
is a symmetric matrix, i.e, $\mathbf{E}^{T}(\mathbf{x}) =\mathbf{E}(\mathbf{x})$, with $\bar{x}_{k} \triangleq x_{k} - x_{c,k}$. 
For every $\mathbf{x} \in \mathcal{D}= \mathcal{X}_{S,+} \cap \mathcal{B}(\mathbf{x}_{c},r)$ (implying that $\| \mathbf{x} - \mathbf{x}_{c} \|_{2} \leq r$), we have that 
\begin{equation} 
\label{equ_proof_suff_cond_D_x_k_SPCM}
|\bar{x}_{k}|  =  |x_{k} - x_{c,k}| = \sqrt{|x_{k} - x_{c,k}|^{2}} \leq \sqrt{\sum_{k \in [N]}  |x_{k} - x_{c,k}|^{2}} = \| \mathbf{x} - \mathbf{x}_{c} \|_{2}  \leq r. 
\end{equation} 
This implies, 
via the triangle inequality for the matrix $2$-norm, that 
\begin{equation} 
\label{equ_proof_suff_cond_norm_E_x_bound_SPCM}
\| \mathbf{E}(\mathbf{x}) \|_{2} \leq N \, q \, r,
\end{equation}
with $q \triangleq \max_{k \in [N]} \| \mathbf{C}_{k} \|_{2}$.
Indeed, we have 
\begin{align}
 \| \mathbf{E}(\mathbf{x}) \|_{2}  & \stackrel{\eqref{equ_proof_suff_cond_D_E_x_SPCM}}{=} \bigg\|\sum_{k \in [N]} \mathbf{C}_{k} \bar{x}_{k} \bigg\|_{2} \leq \sum_{k \in [N]}   \|\mathbf{C}_{k} \bar{x}_{k} \|_{2}  \nonumber \\[4mm]
 & =  \sum_{k \in [N]}  |\bar{x}_{k}|  \|\mathbf{C}_{k} \|_{2}
 \stackrel{\eqref{equ_proof_suff_cond_D_x_k_SPCM}}{\leq}  r \sum_{k \in [N]}   \|\mathbf{C}_{k} \|_{2} \leq   r  N  \max_{k \in [N]} \| \mathbf{C}_{k} \|_{2}. 
\end{align} 

Let us now choose a radius $r_{0}$ such that  $N \, q \, r_{0} < \frac{1}{\big\|  \widetilde{\mathbf{C}}^{-1}(\mathbf{x}_{c}) \big\|}$ which implies via \eqref{equ_proof_suff_cond_norm_E_x_bound_SPCM} that 
\begin{equation} 
\label{equ_suff_cond_D_SPCM_E_x_bound}
\| \mathbf{E}(\mathbf{x}) \|_{2} \leq  N \, q \, r_{0} <\frac{1}{\big\|  \widetilde{\mathbf{C}}^{-1}(\mathbf{x}_{c}) \big\|}
\end{equation} 
for every $\mathbf{x} \in \mathcal{X}_{S,+} \cap \mathcal{B}(\mathbf{x}_{c},r_{0})$.
The fulfillment of condition \eqref{equ_suff_cond_D_SPCM_E_x_bound} allows us then to invoke 
Lemma \ref{lem_matrix_inversion_perturb}, using $\mathbf{A} = \widetilde{\mathbf{C}}(\mathbf{x}_{c})$ and $\mathbf{E} = \mathbf{E}(\mathbf{x})$, which 
yields that for every $\mathbf{x} \in \mathcal{X}_{S,+} \cap \mathcal{B}(\mathbf{x}_{c},r_{0})$ we have
\begin{equation} 
\label{equ_proof_suff_cond_set_D_SPCM_diff_E_0}
\widetilde{\mathbf{C}}^{-1}(\mathbf{x}) = \widetilde{\mathbf{C}}^{-1}(\mathbf{x}_{c}) + \mathbf{E}_{0},
\end{equation} 
with a symmetric matrix $\mathbf{E}_{0} \in \mathbb{R}^{M \times M}$ 
satisfying 
\begin{align} 
\label{equ_E_0_bound_N_q_r}
\| \mathbf{E}_{0} \|_{2} & = \big\| \widetilde{\mathbf{C}}^{-1}(\mathbf{x})  - \widetilde{\mathbf{C}}^{-1}(\mathbf{x}_{c})  \big\|_{2}
 \stackrel{\eqref{equ_matrix_inversion_perturb},\eqref{equ_proof_suff_cond_D_E_x_SPCM}}{\leq}  \frac{\| \widetilde{\mathbf{C}}^{-1}(\mathbf{x}_{c})Ê\|^{2}_{2} \|\mathbf{E}(\mathbf{x}) \|_{2}}
 {1 - \| \widetilde{\mathbf{C}}^{-1}(\mathbf{x}_{c}) \|_{2} \|\mathbf{E}(\mathbf{x}) \|_{2}} \nonumber \\[4mm]
 &  \stackrel{\eqref{equ_suff_cond_D_SPCM_E_x_bound},\eqref{equ_inverse_norm_leq_sigma_2}}{\leq} \frac{N q r_{0} \sigma^{-4}}{1- \sigma^{-2} N q r_{0}}. 
\end{align} 
By choosing $r_{0}$ small enough we can, according to \eqref{equ_E_0_bound_N_q_r}, make $\| \mathbf{E}_{0} \|_{2}$ arbitrarily small. 
Observe that by assumption the matrix $\mathbf{F} \triangleq 2\widetilde{\mathbf{C}}^{-1}(\mathbf{x}_{c})-\widetilde{\mathbf{C}}^{-1}(\mathbf{x}_{0}) \in \mathbb{R}^{M \times M}$ is positive definite, which implies that its eigenvalues $\lambda_{k}(\mathbf{F})$ are positive, i.e., they satisfy $\lambda_{k}(\mathbf{F}) \geq \varepsilon_{0}$ with some 
positive constant $\varepsilon_{0} > 0$. 
We can then choose $r_{0}$ such that $\| \mathbf{E}_{0} \|_{2} < \varepsilon_{0}/2$ (implying that $|\lambda_{M}(\mathbf{E}_{0})| < \varepsilon_{0}/2$) for every $\mathbf{x} \in \mathcal{X}_{S,+} \cap \mathcal{B}(\mathbf{x}_{c},r_{0})$, which implies by Lemma \ref{lem_eigvalues_sum_symmetric_matrices} that 
the eigenvalues $\lambda_{k}(\mathbf{G})$ of the symmetric matrix 
\begin{equation} 
\mathbf{G} \triangleq 2\widetilde{\mathbf{C}}^{-1}(\mathbf{x})  - \widetilde{\mathbf{C}}^{-1}(\mathbf{x}_{0}) \stackrel{\eqref{equ_proof_suff_cond_set_D_SPCM_diff_E_0}}{=} 2\widetilde{\mathbf{C}}^{-1}(\mathbf{x}_{c})-\widetilde{\mathbf{C}}^{-1}(\mathbf{x}_{0})  + 2\mathbf{E}_{0} = \mathbf{F} +  2\mathbf{E}_{0}.
\end{equation} 
satisfy 
\begin{equation}
\lambda_{k}(\mathbf{G}) = \lambda_{k}( \mathbf{F} +  2\mathbf{E}_{0}) \stackrel{\eqref{equ_double_ineq_eigvals_sum}}{\geq} \lambda_{k}(\mathbf{F}) + \lambda_{M}(2\mathbf{E}_{0})=  \lambda_{k}(\mathbf{F}) + 2\lambda_{M}(\mathbf{E}_{0})
\geq \varepsilon_{0} - 2\| \mathbf{E}_{0} \|_{2} >	0, 
\end{equation} 
i.e., all eigenvalues of the matrix $\mathbf{G} =  2\widetilde{\mathbf{C}}^{-1}(\mathbf{x})  - \widetilde{\mathbf{C}}^{-1}(\mathbf{x}_{0})$ are 
strictly positive. This in turn means that the matrix $2\widetilde{\mathbf{C}}^{-1}(\mathbf{x})  - \widetilde{\mathbf{C}}^{-1}(\mathbf{x}_{0})$ is positive definite for every $\mathbf{x} \in \mathcal{X}_{S,+} \cap \mathcal{B}(\mathbf{x}_{c},r_{0})$.

Finally, for $\mathbf{x}_{c} = \mathbf{x}_{0}$, the matrix $2\widetilde{\mathbf{C}}^{-1}(\mathbf{x}_{c})-\widetilde{\mathbf{C}}^{-1}(\mathbf{x}_{0}) =\widetilde{\mathbf{C}}^{-1}(\mathbf{x}_{0})$ is 
obviously positive definite since the matrix $\widetilde{\mathbf{C}}(\mathbf{x}_{0})$ is positive definite. Also, for $\mathbf{x}_{c} = \mathbf{0}$, 
the matrix 
\begin{align} 
\label{equ_proof_suff_cond_D_SPCM_widetilde_C_widetilde_C_1_x_c_x_0}
2\widetilde{\mathbf{C}}^{-1}(\mathbf{x}_{c})-\widetilde{\mathbf{C}}^{-1}(\mathbf{x}_{0}) \stackrel{\eqref{equ_def_cov_matrix_observation_SPCM}}{=} 2\big[\mathbf{C}(\mathbf{0})+ \sigma^{2} \mathbf{I} \big]^{-1} -\widetilde{\mathbf{C}}^{-1}(\mathbf{x}_{0}) \stackrel{\eqref{equ_cov_param_SPCM}}{=} 2 \sigma^{-2} \mathbf{I} - \widetilde{\mathbf{C}}^{-1}(\mathbf{x}_{0})
\end{align}
is positive definite since its eigenvalues obey 
\begin{align}
\lambda_{k}( 2\widetilde{\mathbf{C}}^{-1}(\mathbf{x}_{c})-\widetilde{\mathbf{C}}^{-1}(\mathbf{x}_{0}) ) & \stackrel{\eqref{equ_proof_suff_cond_D_SPCM_widetilde_C_widetilde_C_1_x_c_x_0}}{=} 
\lambda_{k}(2 \sigma^{-2} \mathbf{I} - \widetilde{\mathbf{C}}^{-1}(\mathbf{x}_{0})) \stackrel{\eqref{equ_double_ineq_eigvals_sum}}{\geq} \lambda_{k}(2 \sigma^{-2} \mathbf{I})  + \lambda_{M}(- \widetilde{\mathbf{C}}^{-1}(\mathbf{x}_{0}))  \nonumber \\[4mm]
& =  2 \sigma^{-2} + \lambda_{M}(- \widetilde{\mathbf{C}}^{-1}(\mathbf{x}_{0})) \stackrel{\eqref{equ_inverse_norm_leq_sigma_2}}{\geq}  2 \sigma^{-2}  - \sigma^{-2} = \sigma^{-2} > 0, 
\end{align}
i.e., are all strictly positive.
\end{proof} 

In what follows, we will always assume that the set $\mathcal{D}$ is of the form \eqref{equ_def_D_thm_cond_d_restr_SPCM_ball}. Note 
that any set $\mathcal{D} \subseteq \mathbb{R}^{N}$ of the form \eqref{equ_def_D_thm_cond_d_restr_SPCM_ball} is closed and bounded and 
therefore \emph{compact} \cite[Theorem 2.41]{RudinBookPrinciplesMatheAnalysis}.
If the set $\mathcal{D}$ is chosen such that \eqref{equ_nec_suff_cond_D_SPCM_kernel_exists} is fulfilled, we  
obtain the kernel $R_{\mathcal{M}_{\mathcal{D},\text{SPCM}}}(\cdot, \cdot): \mathcal{D} \times \mathcal{D} \rightarrow \mathbb{R}$ (see \eqref{equ_def_kernel_M} and \eqref{equ_def_SPCM_stat_model}): 
\begin{align}
\label{equ_kernel_D_restr_SPCM}
R_{\mathcal{M}_{\mathcal{D},\text{SPCM}}}(\mathbf{x}_{1}, \mathbf{x}_{2}) 
& \triangleq \mathsf{E}_{\mathbf{x}_{0}}\left\{ \rho_{\mathcal{M}_{\text{SPCM}}}(\mathbf{y} , \mathbf{x}_{1}) \rho_{\mathcal{M}_{\text{SPCM}}}(\mathbf{y},\mathbf{x}_{2}) \right\} 
= \int_{\mathbf{y} \in \mathbb{R}^{M}} \frac{f_{\text{SPCM}}(\mathbf{y}; \mathbf{x}_{1}) f_{\text{SPCM}}(\mathbf{y}; \mathbf{x_{2}})}{f_{\text{SPCM}}(\mathbf{y}; \mathbf{x}_{0})} d\mathbf{y} \nonumber \\[4mm]
& \hspace*{-10mm} \hspace*{-20mm}=  \frac{(2 \pi)^{-M/2} \big[ \detmb{\widetilde{\mathbf{C}}(\mathbf{x}_{0})} \big]^{1/2}}
{\big[\detmb{\widetilde{\mathbf{C}}(\mathbf{x}_{1})} \detmb{\widetilde{\mathbf{C}}(\mathbf{x}_{2})} \big]^{1/2}} \int_{\mathbf{y} \in \mathbb{R}^{M}} 
\exp \left( - \frac{1}{2} \mathbf{y}^{T} \left( \widetilde{\mathbf{C}}^{-1}(\mathbf{x}_{1}) +  \widetilde{\mathbf{C}}^{-1}(\mathbf{x}_{2}) - \widetilde{\mathbf{C}}^{-1}(\mathbf{x}_{0}) \right)Ê\mathbf{y} \right)  
d\mathbf{y} \nonumber \\[4mm]
& \hspace*{-30mm}=  \big[ \detmb{ \widetilde{\mathbf{C}}(\mathbf{x}_{0}) } \big]^{1/2} \ist 
  \big[ \detmb{ \widetilde{\mathbf{C}}(\mathbf{x}_{1}) \ist \widetilde{\mathbf{C}}(\mathbf{x}_{2}) \big( 
  \widetilde{\mathbf{C}}^{-1}(\mathbf{x}_{1}) + \widetilde{\mathbf{C}}^{-1}(\mathbf{x}_{2}) - \widetilde{\mathbf{C}}^{-1}(\mathbf{x}_{0}) \big)} \big]^{-1/2} \nonumber \\[4mm]
& \hspace*{-30mm}=  \big[ \detmb{ \widetilde{\mathbf{C}}(\mathbf{x}_{0}) } \big]^{1/2} \ist 
  \big[ \detmb{ \widetilde{\mathbf{C}}(\mathbf{x}_{1}) + \widetilde{\mathbf{C}}(\mathbf{x}_{2}) - \widetilde{\mathbf{C}}(\mathbf{x}_{1}) \widetilde{\mathbf{C}}^{-1}(\mathbf{x}_{0})\widetilde{\mathbf{C}}(\mathbf{x}_{2}) } \big]^{-1/2}  ,
\end{align}
with $\widetilde{\mathbf{C}}(\mathbf{x})$ as defined in \eqref{equ_def_cov_matrix_observation_SPCM}. 
Obviously, according to elementary linear algebra \cite{golub96}, the kernel $R_{\mathcal{M}_{\mathcal{D},\text{SPCM}}}(\mathbf{x}_{1}, \mathbf{x}_{2})$ is differentiable in the sense of Definition \ref{def_differentiable_kernel} up to any order $m$ and therefore also continuous. 
Moreover, since we assume that the set $	\mathcal{D}$ is of the form \eqref{equ_def_D_thm_cond_d_restr_SPCM_ball}, the function 
\begin{equation} 
k(\cdot): \mathbb{R}^{N}Ê\rightarrow \mathbb{R}: k(\mathbf{x}) = R_{\mathcal{M}_{\mathcal{D},\text{SPCM}}}(\mathbf{x}, \mathbf{x}) = \big[ \detmb{ \widetilde{\mathbf{C}}(\mathbf{x}_{0}) } \big]^{1/2} \ist 
  \big[ \detmb{ 2\widetilde{\mathbf{C}}(\mathbf{x}) - \widetilde{\mathbf{C}}(\mathbf{x}) \widetilde{\mathbf{C}}^{-1}(\mathbf{x}_{0})\widetilde{\mathbf{C}}(\mathbf{x}) } \big]^{-1/2} 
\end{equation}
is bounded over its domain $\mathcal{D}$. 
In order to verify the boundedness, we observe that (i) every set of the form \eqref{equ_def_D_thm_cond_d_restr_SPCM_ball} is 
closed\footnote{We implicitly assume the Hilbert space structure in $\mathbb{R}^{N}$, which is induced by the inner product $\langle \mathbf{x} , \mathbf{y} \rangle \triangle \mathbf{x}^{T} \mathbf{y}$} and (ii) 
the fact that $\detmb{  2 \widetilde{\mathbf{C}}^{-1}(\mathbf{x}) - \widetilde{\mathbf{C}}^{-1}(\mathbf{x}_{0}) }$ is a continuous function of $\mathbf{x}$. 
It follows then from the extreme value theorem \cite[Theorem 4.16]{RudinBookPrinciplesMatheAnalysis} and the compactness, i.e., boundedness and closedness, of $\mathcal{D}$ that there must a exist a positive constant $\varepsilon >0$ such that 
\begin{equation}
\label{equ_nec_suff_cond_D_SPCM_kernel_exists_closed_subset_larger_epsilon}
\detmb{  2 \widetilde{\mathbf{C}}^{-1}(\mathbf{x}) - \widetilde{\mathbf{C}}^{-1}(\mathbf{x}_{0}) } \geq \varepsilon
\end{equation} 
for every vector $\mathbf{x} \in \mathcal{D}$. This implies that 
\begin{align}
k(\mathbf{x}) & =\big[ \detmb{ \widetilde{\mathbf{C}}(\mathbf{x}_{0}) } \big]^{1/2} \ist 
  \big[ \detmb{ 2\widetilde{\mathbf{C}}(\mathbf{x}) - \widetilde{\mathbf{C}}(\mathbf{x}) \widetilde{\mathbf{C}}^{-1}(\mathbf{x}_{0})\widetilde{\mathbf{C}}(\mathbf{x}) } \big]^{-1/2}  \nonumber \\[4mm]
  & =\big[ \detmb{ \widetilde{\mathbf{C}}(\mathbf{x}_{0}) } \big]^{1/2} \ist 
  \bigg[ \big[ \detmb{ \widetilde{\mathbf{C}}(\mathbf{x})} \big]^{2} \detmb{ 2\widetilde{\mathbf{C}}^{-1}(\mathbf{x}) - \widetilde{\mathbf{C}}^{-1}(\mathbf{x}_{0})} \bigg]^{-1/2}  \nonumber \\[4mm]
  & \stackrel{\deto\{\mathbf{A}^{-1}\} = \deto^{-1}\{\mathbf{A}\}}{=}  \big[ \detmb{ \widetilde{\mathbf{C}}(\mathbf{x}_{0}) } \big]^{1/2} \ist  \detmb{ \widetilde{\mathbf{C}}^{-1}(\mathbf{x}) } 
  \big[ \detmb{ 2\widetilde{\mathbf{C}}^{-1}(\mathbf{x}) -  \widetilde{\mathbf{C}}^{-1}(\mathbf{x}_{0}) } \big]^{-1/2}  \nonumber \\[4mm]
& \stackrel{\eqref{equ_inverse_norm_leq_sigma_2},\eqref{equ_nec_suff_cond_D_SPCM_kernel_exists_closed_subset_larger_epsilon}}{\leq}
 \big[ \detmb{ \widetilde{\mathbf{C}}(\mathbf{x}_{0}) } \big]^{1/2} \sigma^{-2M} \varepsilon^{-1/2},
\end{align} 
 i.e., the function $k(\mathbf{x})$ is bounded by the constant (w.r.t.\ $\mathbf{x}$) $\big[ \detmb{ \widetilde{\mathbf{C}}(\mathbf{x}_{0}) } \big]^{1/2} \sigma^{-2M} \varepsilon^{-1/2}$.
Therefore, we can apply Theorem \ref{thm_continous_kernel_implies_cont_func} to derive 
\begin{theorem} 
\label{thm_cont_requirement_estimable_func_SPCM}
Consider the minimum variance problem $\mathcal{M}_{\emph{SPCM}}=\left(\mathcal{E}_{\emph{SPCM}},c(\cdot) \equiv 0, \mathbf{x}_{0}Ê\in \mathcal{X}_{S,+} \right)$ and a set $\mathcal{D}= \mathcal{X}_{S,+} \cap \mathcal{B}(\mathbf{x}_{c},r)$ satisfying 
\eqref{equ_nec_suff_cond_D_SPCM_kernel_exists} so that $R_{\mathcal{M}_{\mathcal{D},\emph{SPCM}}}(\cdot, \cdot): \mathcal{D} \times \mathcal{D} \rightarrow \mathbb{R}$
exists. Then, if the parameter function $g(\cdot): \mathcal{X}_{S,+} \rightarrow \mathbb{R}$ of the associated SPCM $\mathcal{E}_{\emph{SPCM}}$ has a discontinuity on the set $\mathcal{D}$, i.e., there exists a point $\mathbf{x} \in \mathcal{D}$ at which the function $g(\cdot)$ is discontinuous, there exists no allowed 
estimator (i.e., unbiased and with finite variance at $\mathbf{x}_{0}$) for $\mathcal{M}_{\emph{SPCM}}$. 
\end{theorem}

\begin{proof}
Consider the $\mathcal{D}$-restricted minimum variance problem $\mathcal{M}' = \mathcal{M}_{\text{SPCM}} \big|_{\mathcal{D}}$ with the set $\mathcal{D}$ of the form 
\eqref{equ_def_D_thm_cond_d_restr_SPCM_ball}. The RKHS $\mathcal{H}(\mathcal{M}')$ associated to the kernel 
$R_{\mathcal{M}'}$ (see \eqref{equ_kernel_D_restr_SPCM}) exists since the set $\mathcal{D}$ is assumed to satisfy \eqref{equ_nec_suff_cond_D_SPCM_kernel_exists}. 
Moreover, according to our discussion above, Theorem \ref{thm_continous_kernel_implies_cont_func} applies to the RKHS $\mathcal{H}(\mathcal{M}')$.
We have by Theorem \ref{thm_continous_kernel_implies_cont_func} that $\mathcal{H}(\mathcal{M}')$ consists exclusively of continuous functions. Therefore, it cannot contain 
any function with a discontinuity within the set $\mathcal{D}$. By Theorem \ref{thm_main_facts_RKHS_MVE}, this implies that 
any discontinuous function is not estimable for $\mathcal{M}'$. However, since $\mathcal{M}'$ is obtained by a reduction of the parameter set underlying 
$\mathcal{M}_{\text{SPCM}}$, it follows from Theorem \ref{thm_par_set_reduction_classic_est_mve} that also for the minimum variance problem 
$\mathcal{M}_{\text{SPCM}}$ every 
discontinuous function is not estimable, i.e., there exists no unbiased estimator of a discontinuous (on $\mathcal{D}$) function with finite variance at 
$\mathbf{x}_{0} \in \mathcal{X}_{S,+}$. 
\end{proof}

As for the SLM (cf.\ Theorem \ref{thm_non_esit_finite_var_unbiased_est_supp}), we have also for the SPCM that there exists no finite-variance unbiased estimator of an injective function of the support $\supp(\mathbf{x})$: 
\begin{corollary}
\label{cor_non_esit_finite_var_unbiased_est_supp_SPCM}
There exists no estimator $\hat{g}(\mathbf{y})$ of an injective real-valued function $g(\supp(\mathbf{x}))$ of the support of $\mathbf{x} \in \mathcal{X}_{S,+}$ which uses only the observation \eqref{equ_obs_model_SPCM} of the SPCM $\mathcal{E}_{\emph{SPCM}}$, is unbiased for every $\mathbf{x} \in \mathcal{X}_{S}$, and has 
a finite variance at any $\mathbf{x}_{0} \in \mathcal{X}_{S,+}$.
\end{corollary}

\begin{proof}
Consider a set $\mathcal{D} = \mathcal{X}_{S,+} \cap \mathcal{B}(\mathbf{0},r)$ with a sufficiently small radius $r>0$ such that, according to Theorem \ref{thm_cond_d_restr_SPCM_ball}, the set $\mathcal{D}$ fulfills \eqref{equ_nec_suff_cond_D_SPCM_kernel_exists}, and with an arbitrary parameter vector $\mathbf{x}_{0}\in \mathcal{X}_{S,+}$. 
As verified in the proof of Theorem \ref{thm_non_esit_finite_var_unbiased_est_supp}, any injective function $g(\supp(\mathbf{x}))$ of the support must have at least one discontinuity on the set $\mathcal{D}$. 
Therefore, by Theorem \ref{thm_cont_requirement_estimable_func_SPCM}, there does not exist an allowed estimator for the minimum variance problem 
$\mathcal{M}_{\text{SPCM}}=\left(\mathcal{E}_{\text{SPCM}},c(\cdot) \equiv 0, \mathbf{x}_{0}Ê\in \mathcal{X}_{S,+} \right)$ associated with the SPCM with parameter function $g(\supp(\mathbf{x}))$.
By the definition of an allowed estimator, this means that there exists no estimator of $g(\supp(\mathbf{x}))$ that is unbiased for every $\mathbf{x} \in \mathcal{X}_{S,+}$ and 
has finite variance at $\mathbf{x}_{0}$ . 
\end{proof}

\section{Lower Bounds on the Estimator Variance for the SPCM} 
\label{sec_lower_bounds_SPCM}

Let us consider the minimum variance problem $\mathcal{M}_{\text{SPCM}}=\left(\mathcal{E}_{\text{SPCM}},c(\cdot) \equiv 0, \mathbf{x}_{0}Ê\right)$ associated 
with the SPCM, for a parameter vector $\mathbf{x}_{0} \in \mathcal{X}_{S,+}$. Note again that the prescribed bias $c(\cdot)$ is assumed to 
vanish on $\mathcal{X}_{S,+}$, i.e., 
we consider only unbiased estimation. However, by Theorem \ref{thm_equ_bias_param_function}, this
is no real restriction since we allow for general parameter functions $g(\mathbf{x})$ in the SPCM. 
We then define a second minimum variance problem 
$\mathcal{M}_{\mathcal{D},\text{SPCM}}Ê\triangleq \mathcal{M}_{\text{SPCM}}Ê\big|_{\mathcal{D}} = \big(\mathcal{E}^{(\mathcal{D})}_{\text{SPCM}},c(\cdot)\equiv 0,\mathbf{x}_{0}Ê\big)$
associated to the $\mathcal{D}$-restricted SPCM, which is identical to $\mathcal{M}_{\text{SPCM}}$ except for the parameter set $\mathcal{D} \subseteq \mathcal{X}_{S,+}$. 

We will now derive lower bounds on the minimum achievable variance $L_{\mathcal{M}_{\text{SPCM}}}$ via \eqref{equ_lower_bound_spcm_restr_spcm}. 
We choose $\mathcal{D}$ such that Theorem \ref{thm_cond_d_restr_SPCM_ball} applies 
and therefore the RKHS $\mathcal{H}(\mathcal{M}_{\mathcal{D},\text{SPCM}})$, associated with the kernel 
$R_{\mathcal{M}_{\mathcal{D},\text{SPCM}}}(\cdot,\cdot)$ given in \eqref{equ_kernel_D_restr_SPCM}, exists. 
We will assume that the parameter function $g(\cdot)$ is estimable for $\mathcal{M}_{\text{SPCM}}$, which implies via Theorem \ref{thm_par_set_reduction_classic_est_mve} that the restricted parameter function $g(\cdot)\big|_{\mathcal{D}}$ is estimable for $\mathcal{M}_{\mathcal{D},\text{SPCM}}$. 
This assumption is justified by the fact that the bounds that we will derive are always finite, i.e., they trivially apply also in the case when $g(\cdot)\big|_{\mathcal{D}}$ is not estimable for $\mathcal{M}_{\mathcal{D},\text{SPCM}}$, for which $L_{\mathcal{M}_{\text{SPCM}}}=\infty$. 
These bounds will be obtained by projecting the restricted parameter function $g\big|_{\mathcal{D}}(\cdot): \mathcal{D} \rightarrow \mathbb{R}$, which is equal to the prescribed mean function $\gamma(\mathbf{x}) = c(\mathbf{x}) + g(\mathbf{x})$ of $\mathcal{M}_{\mathcal{D},\text{SPCM}}$ since $c(\mathbf{x})=0$, on suitable subspaces 
$\mathcal{U}$ of the RKHS $\mathcal{H}(\mathcal{M}_{\mathcal{D},\text{SPCM}})$. 
Our choice of $\mathcal{U}$ will be inspired by the derivation of the lower bounds for the SLM in Chapter \ref{chap_SLM}. Hence, there will 
be a natural correspondence between the bounds for the SPCM and the bounds for the SLM presented in Chapter \ref{chap_SLM}. 

For the derivation of the bounds, we will need an identity which is stated in  
\begin{lemma}
\label{lem_der_det_function}
Consider the function $f(\cdot): \mathbb{R}^{N} \rightarrow \mathbb{R}: f(\mathbf{x}) = \detmb{\mathbf{H}(\mathbf{x})}$, where the 
matrix-valued function $\mathbf{H}(\cdot): \mathbb{R}^{N} \rightarrow \mathbb{R}^{M \times M}$ is such that there exist the partial derivatives 
$\frac{\partial^{\mathbf{p}} \left( \mathbf{H}(\mathbf{x}) \right)_{m,n}}{\partial \mathbf{x}^{\mathbf{p}}}\big|_{\mathbf{x} = \mathbf{x}_{0}}$ 
for any order $\mathbf{p} \in \mathbb{Z}_{+}^{N}$ with $\| \mathbf{p} \|_{\infty} \leq 1$ as well as the inverse $\mathbf{H}^{-1}(\mathbf{x}_{0})$. Then we have 
\vspace*{4mm}
\begin{align}
\label{equ_lem_der_det_function}
\frac{ \partial  f(\mathbf{x})}{\partial x_{k}} \bigg|_{\mathbf{x} = \mathbf{x}_{0}} = \detm{\mathbf{H}(\mathbf{x}_{0})} \trace \big\{ \mathbf{H}^{-1}(\mathbf{x}_{0}) \left( \partial \mathbf{H}(\mathbf{x}_{0}) \right)^{T} \big\},Ê\\[-4mm]
\nonumber
\end{align}  
where $\partial_{k} \mathbf{H}(\mathbf{x}_{0})$ is the matrix whose entry in the $m$th row and $n$th column is given by $\frac{\partial \left( \mathbf{H}(\mathbf{x}) \right)_{m,n}}{\partial x_{k}}\bigg|_{\mathbf{x} = \mathbf{x}_{0}}$.
\end{lemma} 
\begin{proof}
\cite[p.\ 73]{kay}  
\end{proof} 

The first lower bound on $L_{\mathcal{M}_{\text{SPCM}}}$ will be the ``SPCM analogue'' of the SLM bound given by Theorem \ref{thm_CRB_SLM}. 
\begin{theorem}
\label{thm_CRB_SPCM}
Consider the minimum variance problem $\mathcal{M}_{\emph{SPCM}}=\left(\mathcal{E}_{\emph{SPCM}},c(\cdot) \equiv 0 , \mathbf{x}_{0}Ê\right)$ associated with the SPCM with arbitrary 
psd basis matrices $\big\{ \mathbf{C}_{k} \in \mathbb{R}^{M \times M} \big\}_{k \in [N]}$. 
If the parameter function $g(\cdot):\mathcal{X}_{S,+} \rightarrow \mathbb{R}$ is such that the partial derivatives $\frac{\partial g(\mathbf{x})}{\partial x_{l}} \big|_{\mathbf{x} = \mathbf{x}_{0}}$ exist for $l \in [N]$, 
we have 
\vspace*{2mm}
\begin{align} 
& L_{\mathcal{M}_{\emph{SPCM}}} \geq  \mathbf{b}^{T}  \mathbf{J}^{\dagger} \mathbf{b}  \hspace*{17mm} \mbox{when}\,\, \| \mathbf{x}_{0} \|_{0} < S \label{equ_CRB_SPCM_not_full_sparsity} \\[4mm]
& L_{\mathcal{M}_{\emph{SPCM}}} \geq  \mathbf{b}_{\mathbf{x}_{0}}^{T} \mathbf{J}_{\mathbf{x}_{0}}^{\dagger} \mathbf{b}_{\mathbf{x}_{0}}  \hspace*{10.6mm} \mbox{when}\,\, \| \mathbf{x}_{0} \|_{0} = S. \label{equ_CRB_SPCM_full_sparsity}  \\[-5mm]
\nonumber 
\end{align} 
Here, $\mathbf{b} \in \mathbb{R}^{N}$ is defined elementwise by $b_{l} \triangleq  \frac{\partial g(\mathbf{x})}{\partial x_{l}} \big|_{\mathbf{x} = \mathbf{x}_{0}}$, and 
$\mathbf{b}_{\mathbf{x}_{0}} \in \mathbb{R}^{S}$ denotes the restriction to the entries of $\mathbf{b}$ indexed by $\supp(\mathbf{x}_{0})=(i_{1},\ldots,i_{S})$, i.e., $\big( \mathbf{b}_{\mathbf{x}_{0}}\big)_j  = b_{i_{j}}$.
Furthermore, the matrix $\mathbf{J} \in \mathbb{R}^{N \times N}$ is given elementwise by 
\begin{equation} 
\label{equ_FIM_CRB_SPCM_dim_N}
\big( \mathbf{J} \big)_{m,n} =Ê\frac{1}{2} \trace\big\{ \widetilde{\mathbf{C}}^{-1}(\mathbf{x}_{0}) \mathbf{C}_{m} \widetilde{\mathbf{C}}^{-1}(\mathbf{x}_{0}) \mathbf{C}_{n} \big\},
\end{equation} 
and the matrix $\mathbf{J}_{\mathbf{x}_{0}} \in \mathbb{R}^{S \times S}$ is given elementwise by 
\begin{equation}
\big( \mathbf{J}_{\mathbf{x}_{0}} \big)_{m,n} = \frac{1}{2} \trace\big\{ \widetilde{\mathbf{C}}^{-1}(\mathbf{x}_{0})\mathbf{C}_{i_{m}} \widetilde{\mathbf{C}}^{-1}(\mathbf{x}_{0}) \mathbf{C}_{i_{n}}\big\}.
\end{equation}  
\end{theorem} 

\begin{proof}
Consider the minimum variance problem $\mathcal{M}_{\mathcal{D},\text{SPCM}}=\mathcal{M}_{\text{SPCM}}\big|_{\mathcal{D}}$, 
where the set $\mathcal{D}$ is chosen as in \eqref{equ_def_D_thm_cond_d_restr_SPCM_ball} with $\mathbf{x}_{c} = \mathbf{x}_{0}$, 
so that Theorem \ref{thm_cond_d_restr_SPCM_ball} applies, i.e., 
the RKHS $\mathcal{H}(\mathcal{M}_{\mathcal{D},\text{SPCM}})$ exists. 
We can assume that $g(\cdot)$ is estimable for $\mathcal{M}_{\text{SDPCM}}$, which implies via Theorem \ref{thm_par_set_reduction_classic_est_mve} that 
the restriction $g(\cdot)\big|_{\mathcal{D}}$ is estimable for $\mathcal{M}_{\mathcal{D}, \text{SDPCM}}$. Thus, according to Theorem \ref{thm_main_facts_RKHS_MVE}, 
the prescribed mean function $\gamma(\cdot): \mathcal{D} \rightarrow \mathbb{R}: \gamma(\mathbf{x}) = c(\mathbf{x}) + g(\mathbf{x})= g(\mathbf{x})$ belongs to the RKHS $\mathcal{H}(\mathcal{M}_{\mathcal{D},\text{SPCM}})$. 

For the case $\| \mathbf{x}_{0} \|_{0} < S$, in order to prove \eqref{equ_CRB_SPCM_not_full_sparsity}, consider the subspace $\mathcal{U}_{1} \triangleq \linspan \big\{ w_{0}(\cdot) \cup \{ w_{l}(\cdot) \}_{l \in [N]}Ê\big\}$ 
spanned by the functions $w_{0} (\cdot) \triangleq R_{\mathcal{M}_{\mathcal{D},\text{SPCM}}} (\cdot, \mathbf{x}_{0})$ and 
\begin{equation}
\label{equ_def_w_l_proof_SPCM_CRB_part_der}
w_{l}(\cdot) \triangleq \frac{ \partial^{\mathbf{e}_{l}} R_{\mathcal{M}_{\mathcal{D},\text{SPCM}}} (\cdot, \mathbf{x}_{2})}  {\partial \mathbf{x}_{2}^{\mathbf{e}_{l}}} \bigg|_{\mathbf{x}_{2} = \mathbf{x}_{0}} 
\end{equation} 
for $l \! \in \! [N]$. 
We have trivially $w_{0}(\cdot) \in \mathcal{H}(\mathcal{M}_{\mathcal{D},\text{SPCM}})$, and by Theorem \ref{thm_der_repr_prop} under the condition 
$\| \mathbf{x}_{0} \|_{0} < S$, we have also that 
$w_{l} (\cdot)  \in \mathcal{H}(\mathcal{M}_{\mathcal{D},\text{SPCM}})$ for $l \in [N]$. 
In a completely analogous manner as the RKHS-based derivation of Theorem \ref{thm_unconstr_CR} (see Section \ref{sec_CRB}), one can show that 
$\big\langle w_{0}(\cdot), w_{l} (\cdot) \big\rangle_{ \mathcal{H}(\mathcal{M}_{\mathcal{D},\text{SPCM}})} = 0$ (cf.\ \eqref{equ_CR_orthog_vecs_1}), 
$\big\langle w_{0}(\cdot), \gamma (\cdot) \big\rangle_{ \mathcal{H}(\mathcal{M}_{\mathcal{D},\text{SPCM}})} = \gamma (\mathbf{x}_{0})$, 
and $\big\langle  w_{l}(\cdot), \gamma(\cdot) \big\rangle_{ \mathcal{H}(\mathcal{M}_{\mathcal{D},\text{SPCM}})} = b_{l}$ (cf.\ \eqref{equ_inner_prod_par_der_derivation_UCRB}) for $l \in [N]$ (note that $\gamma(\cdot)=g(\cdot)$). 
For $l,l' \in [N]$, we obtain by \eqref{equ_kernel_D_restr_SPCM}, Lemma \ref{lem_der_det_function}, and the obvious fact $\frac{\partial^{\mathbf{e}_{l}}}{\partial \mathbf{x}^{\mathbf{e}_{l}}} \widetilde{\mathbf{C}}(\mathbf{x}) = \mathbf{C}_{l}$ that 
\vspace*{2mm}
\begin{align}
\big\langle w_{l}(\cdot), w_{l'}(\cdot) \big\rangle_{ \mathcal{H}(\mathcal{M}_{\mathcal{D},\text{SPCM}})} & = \bigg\langle \frac{ \partial^{\mathbf{e}_{l}} R_{\mathcal{M}_{\mathcal{D},\text{SPCM}}} (\cdot, \mathbf{x}_{2})}  {\partial \mathbf{x}_{2}^{\mathbf{e}_{l}}} \bigg|_{\mathbf{x}_{2} = \mathbf{x}_{0}} , \frac{ \partial^{\mathbf{e}_{l'}} R_{\mathcal{M}_{\mathcal{D},\text{SPCM}}} (\cdot, \mathbf{x}_{2})}  {\partial \mathbf{x}_{2}^{\mathbf{e}_{l'}}} \bigg|_{\mathbf{x}_{2} = \mathbf{x}_{0}}  \bigg\rangle_{ \mathcal{H}(\mathcal{M}_{\mathcal{D},\text{SPCM}})} \nonumber \\[4mm]
& \hspace*{-30mm} \stackrel{\eqref{equ_der_reproduction_prop}}{=}
 \frac{ \partial^{\mathbf{e}_{l}}\partial^{\mathbf{e}_{l'}} R_{\mathcal{M}_{\mathcal{D},\text{SPCM}}} (\mathbf{x}_{1}, \mathbf{x}_{2})}  
 {\partial  \mathbf{x}_{1}^{\mathbf{e}_{l}} \partial\mathbf{x}_{2}^{\mathbf{e}_{l'}}} \bigg|_{\mathbf{x}_{1}=\mathbf{x}_{2} = \mathbf{x}_{0}} Ê\nonumber \\[4mm]
 & \hspace*{-30mm} \stackrel{\eqref{equ_kernel_D_restr_SPCM}}{=}  \frac{ \partial^{\mathbf{e}_{l}}\partial^{\mathbf{e}_{l'}} \big[ \detmb{ \widetilde{\mathbf{C}}(\mathbf{x}_{0}) } \big]^{1/2} \ist 
  \big[ \detmb{\widetilde{\mathbf{C}}(\mathbf{x}_{1}) + \widetilde{\mathbf{C}}(\mathbf{x}_{2}) - \widetilde{\mathbf{C}}(\mathbf{x}_{1}) \widetilde{\mathbf{C}}^{-1}(\mathbf{x}_{0})\widetilde{\mathbf{C}}(\mathbf{x}_{2})} \big]^{-1/2} }  
  {\partial  \mathbf{x}_{1}^{\mathbf{e}_{l}}\mathbf{x}_{2}^{\mathbf{e}_{l'}}} \bigg|_{\mathbf{x}_{1}=\mathbf{x}_{2} = \mathbf{x}_{0}} \nonumber \\[4mm]
  & \hspace*{-30mm} \stackrel{(a)}{=}  - \frac{1}{2} 
  \big[ \detmb{ \widetilde{\mathbf{C}}(\mathbf{x}_{0}) } \big]^{1/2} \times  \nonumber \\[2mm]
  & \hspace*{-20mm}
 \frac{ \partial^{\mathbf{e}_{l}}  \tracem{\widetilde{\mathbf{C}}^{-1}(\mathbf{x}_{0})(\mathbf{I} - \widetilde{\mathbf{C}}(\mathbf{x}_{1})\widetilde{\mathbf{C}}^{-1}(\mathbf{x}_{0})) \mathbf{C}_{l'}}
  \big[ \detmb{ \widetilde{\mathbf{C}}(\mathbf{x}_{1}) + \widetilde{\mathbf{C}}(\mathbf{x}_{0}) - \widetilde{\mathbf{C}}(\mathbf{x}_{1}) } \big]^{-1/2}} 
  {\partial  \mathbf{x}_{1}^{\mathbf{e}_{l}}}  \bigg|_{\mathbf{x}_{1}=\mathbf{x}_{0}}    \nonumber \\[4mm]
  & \hspace*{-30mm} =  - \frac{1}{2} \frac{ \partial^{\mathbf{e}_{l}}  
  \tracem{\widetilde{\mathbf{C}}^{-1}(\mathbf{x}_{0})(\mathbf{I} - \widetilde{\mathbf{C}}(\mathbf{x}_{1})\widetilde{\mathbf{C}}^{-1}(\mathbf{x}_{0}))\mathbf{C}_{l'}}} 
  {\partial  \mathbf{x}_{1}^{\mathbf{e}_{l}}} \bigg|_{\mathbf{x}_{1} = \mathbf{x}_{0}} \nonumber \\[4mm]
  & \hspace*{-30mm} =  \frac{1}{2}  
  \tracem{\widetilde{\mathbf{C}}^{-1}(\mathbf{x}_{0}) \mathbf{C}_{l} \widetilde{\mathbf{C}}^{-1}(\mathbf{x}_{0})\mathbf{C}_{l'}} \stackrel{\eqref{equ_FIM_CRB_SPCM_dim_N}}{=} \big( \mathbf{J} \big)_{l,l'}.
\end{align}
Here, the step $(a)$ is due to Lemma \ref{lem_der_det_function} using $\mathbf{H}(\mathbf{x}_{2}) = \widetilde{\mathbf{C}}(\mathbf{x}_{1}) + \widetilde{\mathbf{C}}(\mathbf{x}_{2}) - \widetilde{\mathbf{C}}(\mathbf{x}_{1}) \widetilde{\mathbf{C}}^{-1}(\mathbf{x}_{0})\widetilde{\mathbf{C}}(\mathbf{x}_{2})$ and the chain rule for differentiation \cite{RudinBookPrinciplesMatheAnalysis}. 
The bound \eqref{equ_CRB_SPCM_not_full_sparsity} is then obtained by projecting $\gamma(\cdot)$ on the subspace $\mathcal{U}_{1}$ as derived in the following. 
Using Theorem \ref{thm_main_facts_RKHS_MVE}, and Theorem \ref{thm_orthog_proj_ineq},  
\begin{align}
L_{\mathcal{M}_{\text{SPCM}}}  & \stackrel{\eqref{equ_squared_norm_min_achiev_var}}{=} \| \gamma(\cdot) \|^{2}_{\mathcal{H}(\mathcal{M}_{\mathcal{D},\text{SPCM}})} 
- \big[ \underbrace{\gamma(\mathbf{x}_{0})}_{=g(\mathbf{x}_{0})} \big]^{2}  \stackrel{\eqref{equ_squared_norm_orthog_proj_pythag_thm}}{\geq} \| \mathbf{P}_{\mathcal{U}_{1}} \gamma(\cdot) \|^{2} -
 \big( g(\mathbf{x}_{0}) \big)^{2}.
 \end{align}  
 According to Theorem \ref{thm_norm_projection_finite_dim_subspace_union_orthog_subspaces}, this becomes further
 \begin{align}
 \label{equ_proof_derivation_CRB_SPCM_not_full_sparsity}
L_{\mathcal{M}_{\text{SPCM}}} & \stackrel{\eqref{equ_idendity_projection_subspace_spanned_union_orthogonal_sets}}{\geq}  \frac{\big\langle w_{0}(\cdot), \gamma (\cdot) \big\rangle_{ \mathcal{H}(\mathcal{M}_{\mathcal{D},\text{SPCM}})}^{2}}{
\big\langle w_{0}(\cdot), w_{0}(\cdot) \big\rangle_{ \mathcal{H}(\mathcal{M}_{\mathcal{D},\text{SPCM}})}}  + \mathbf{b}^{T} \mathbf{J}^{\dagger} \mathbf{b} - \big[ g(\mathbf{x}_{0}) \big]^{2}  \nonumber \\[4mm]
& =  \frac{\big\langle R_{\mathcal{M}_{\mathcal{D},\text{SPCM}}} (\cdot, \mathbf{x}_{0}), \gamma (\cdot) \big\rangle_{ \mathcal{H}(\mathcal{M}_{\mathcal{D},\text{SPCM}})}^{2}}{
\big\| R_{\mathcal{M}_{\mathcal{D},\text{SPCM}}} (\cdot, \mathbf{x}_{0}) \big\|^{2}_{ \mathcal{H}(\mathcal{M}_{\mathcal{D},\text{SPCM}})}}  + \mathbf{b}^{T} \mathbf{J}^{\dagger} \mathbf{b} - \big[ g(\mathbf{x}_{0}) \big]^{2}  \nonumber \\[4mm]
&  \stackrel{(a)}{=}  \frac{ \big[ g(\mathbf{x}_{0}) \big]^{2}}{1}  + \mathbf{b}^{T} \mathbf{J}^{\dagger} \mathbf{b} - \big[ g(\mathbf{x}_{0}) \big]^{2}  \nonumber \\[4mm]
&  =  \mathbf{b}^{T} \mathbf{J}^{\dagger} \mathbf{b},
\end{align} 
which is \eqref{equ_CRB_SPCM_not_full_sparsity}. 
In \eqref{equ_proof_derivation_CRB_SPCM_not_full_sparsity}, the step $(a)$ follows from the facts that 
\begin{equation}
\big\langle R_{\mathcal{M}_{\mathcal{D},\text{SPCM}}} (\cdot, \mathbf{x}_{0}), \gamma (\cdot) \big\rangle_{ \mathcal{H}(\mathcal{M}_{\mathcal{D},\text{SPCM}})} \stackrel{\eqref{equ_reproduction_property}}{=} \gamma(\mathbf{x}_{0}) = g(\mathbf{x}_{0}) 
\end{equation} 
and 
\begin{align}
 \big\| R_{\mathcal{M}_{\mathcal{D},\text{SPCM}}} (\cdot, \mathbf{x}_{0}) \big\|^{2}_{ \mathcal{H}(\mathcal{M}_{\mathcal{D},\text{SPCM}})} & =
 \big\langle R_{\mathcal{M}_{\mathcal{D},\text{SPCM}}} (\cdot, \mathbf{x}_{0}),R_{\mathcal{M}_{\mathcal{D},\text{SPCM}}} (\cdot, \mathbf{x}_{0}) \big\rangle_{ \mathcal{H}(\mathcal{M}_{\mathcal{D},\text{SPCM}})} \nonumber \\[4mm]
 & \stackrel{\eqref{equ_reproduction_property}}{=}R_{\mathcal{M}_{\mathcal{D},\text{SPCM}}} (\mathbf{x}_{0}, \mathbf{x}_{0}) \stackrel{\eqref{equ_kernel_one_arg_x_0_equal_1}}{=} 1.
\end{align} 

In an almost identical manner, one can prove the bound \eqref{equ_CRB_SPCM_full_sparsity} for the case $\| \mathbf{x}_{0} \|_{0} = S$. The only difference is that 
instead of using the subspace $\mathcal{U}_{1}$ (which does not exist in this case), one has to use the subspace 
$\mathcal{U}_{2} \triangleq \linspan \big\{ w_{0}(\cdot) \cup \{ w_{l}(\cdot) \}_{l \in \supp(\mathbf{x}_{0})}Ê\big\}$, 
with $w_{l}(\cdot)$ as defined in \eqref{equ_def_w_l_proof_SPCM_CRB_part_der}.
\end{proof} 

The bounds \eqref{equ_CRB_SPCM_not_full_sparsity}, \eqref{equ_CRB_SPCM_full_sparsity} have a particular interpretation as the CRB of a non-sparse estimation problem which is 
identical to the SPCM except for the sparsity constraint. More precisely, the bounds \eqref{equ_CRB_SPCM_not_full_sparsity}, \eqref{equ_CRB_SPCM_full_sparsity} are  
the CRB for a general non-sparse estimation problem with a Gaussian observation $\mathbf{y} \sim \mathcal{N}(\mathbf{0},\mathbf{C}(\mathbf{x}))$. Here, the 
covariance matrix $\mathbf{C}(\mathbf{x})$ depends arbitrarily on the parameter vector $\mathbf{x}$, which is not required to be sparse \cite[p.\ 47]{kay}.  

A specific similarity of Theorem \ref{thm_CRB_SPCM} to Theorem \ref{thm_CRB_SLM} becomes apparent when the corresponding lower bounds are 
evaluated for the special case of the SPCM with $g(\mathbf{x}) = x_{k}$ where $k \notin \supp(\mathbf{x}_{0})$ and the basis 
matrices $\mathbf{C}_{l} =Ê\mathbf{e}_{l}\mathbf{e}^{T}_{l}$. This choice yields an instance of the SDPCM (that will 
be discussed presently). For this special case, we observe a discontinuity between the two domains of $\mathcal{X}$ corresponding to $\| \mathbf{x} \|_{0} = S$ and $\| \mathbf{x} \|_{0} < S$. 
In fact, while the bound in \eqref{equ_CRB_SPCM_not_full_sparsity} for $\| \mathbf{x}_{0} \|_{0} < S$ is then equal to $2 \sigma^4$, the bound \eqref{equ_CRB_SPCM_full_sparsity} for $\| \mathbf{x}_{0} \|_{0} = S$ is equal to $0$.

In what follows, we will need a standard result in matrix calculus: 
\begin{lemma}
\label{lem_derivation_inverse_matrix} 
Consider a matrix-valued function $\mathbf{H}(\cdot): \mathbb{R}Ê\rightarrow \mathbb{R}^{M \times M}$ such 
that the derivatives of any component function $f_{m,n}(x) \triangleq \left( \mathbf{H}(x) \right)_{m,n}$ exist to a sufficient order. 
We denote by  $\partial \mathbf{H}(x) \in \mathbb{R}^{M \times M}$ the specific matrix whose entries are given by $\frac{\partial f_{m,n}(x)}{\partial x}$. 
We then have the identity 
\begin{equation} 
\label{equ_derivation_inverse_matrix}
\partial \mathbf{H}^{-1}(x)= -\mathbf{H}^{-1}(x) \partial \mathbf{H}(x) \mathbf{H}^{-1}(x). 
\end{equation}
\end{lemma}
\begin{proof} 
\cite{golub96}
\end{proof}

Whereas Theorem \ref{thm_CRB_SPCM} places no assumption on the basis matrices $\big\{ \mathbf{C}_{k} \big\}_{k \in [N]}$ of the SPCM, 
the next result requires the basis matrices to satisfy the RIP of order $S+1$ with a sufficiently small constant $\delta_{S+1}$:
\begin{theorem} 
\label{thm_sparse_lower_bound_SPCM}
Consider the minimum variance problem $\mathcal{M}_{\emph{SPCM}}=\left(\mathcal{E}_{\emph{SPCM}},c(\cdot) \equiv 0 , \mathbf{x}_{0}Ê\right)$ 
associated with the SPCM with sparsity degree $S$ and basis matrices $\big\{ \mathbf{C}_{k} \big\}_{k \in [N]}$, with rank $r_{k} \triangleq \rank(\mathbf{C}_{k})$, that satisfy the RIP of order $S+1$ with RIP constant $\delta_{S+1}< 1/32$. 
If the parameter function $g(\cdot):\mathcal{X}_{S,+} \rightarrow \mathbb{R}$ is such that the partial derivatives $b_{l} \triangleq \frac{\partial g(\mathbf{x})}{\partial x_{l}} \big|_{\mathbf{x} = \mathbf{0}}$ exist 
for $l \in [N]\setminus \supp(\mathbf{x}_{0})$, 
we have 
\begin{align} 
\label{equ_lower_bound_SPCM_RIP}
L_{\mathcal{M}_{\emph{SPCM}}} &\geq  \nonumber \\Ê
&\hspace*{-5mm} \frac{2\sigma^{4}  b_{l}^{2}}{r_{l}} \frac{ (\sigma^2 \beta)^{Q/2}(\sigma^{2} (1+5\delta_{S+1}))^{-r_{l}/2}}{ \prod_{k \in \supp(\mathbf{x}_{0})} (x_{0,k}+ \sigma^{2} (1+5\delta_{S+1}))^{r_{k}/2}}
  \beta^{2} \Bigg[ r_{l}\frac{\left( 12 \delta_{S+1}\right)^{2}}{2} +   (12 \delta_{S+1} )^{2} +  \beta (6 \delta_{S+1}+ 1) \Bigg]^{-1} 
\end{align} 
for any $l \in [N] \setminus \supp(\mathbf{x}_{0})$. Here,   
\begin{equation} 
\beta \triangleq \frac{1}{(1+\delta_{S+1})^{2}} \bigg( 2 - \frac{(1+\delta_{S+1})^{4}}{(1-\delta_{S+1})^{4}} \bigg), 
\end{equation}  
and $Q \triangleq \sum_{k \in \supp(\mathbf{x}_{0})} r_{k} + r_{l}$. 
\end{theorem} 
\begin{proof} 
Appendix \ref{chap_appendix_D}.
\end{proof}

Note that Theorem \ref{thm_sparse_lower_bound_SPCM} yields a lower bound on the minimum achievable variance $L_{\mathcal{M}_{\text{SPCM}}}$ 
that is continuous with respect to the parameter vector $\mathbf{x}_{0}$. 
This is in contrast to the bound given by Theorem \ref{thm_CRB_SPCM} which yields a discontinuous bound in general. 

However, as can be verified easily, for the special case of the SPCM given by the SDPCM with parameter function $g(\mathbf{x}) = x_k$ where $k \notin \supp(\mathbf{x}_{0})$, the lower bound \eqref{equ_lower_bound_SPCM_RIP} of Theorem \ref{thm_sparse_lower_bound_SPCM} is in general looser, i.e., lower, than the bound \eqref{equ_CRB_SPCM_not_full_sparsity}, \eqref{equ_CRB_SPCM_full_sparsity} of Theorem \ref{thm_CRB_SPCM}. 

In the following, we specialize our RKHS approach for the general SPCM to the SDPCM and derive a lower bound that is both continuous with respect to $\mathbf{x}_{0}$ and, for the parameter function $g(\mathbf{x})= x_k$ with $k \notin \supp(\mathbf{x}_{0})$, at least as tight as the bound given by 
\eqref{equ_CRB_SPCM_not_full_sparsity} and \eqref{equ_CRB_SPCM_full_sparsity}. 

\section{The Sparse Diagonalizable Parametric Covariance Model}
\label{sec_SDPCM}

We now specialize and extend the results found for the general SPCM to the special case termed the \emph{sparse diagonalizable parametric covariance model} (SDPCM). 
The SDPCM is obtained from the SPCM by restricting to basis matrices that are orthogonal projection matrices on orthogonal subspaces and, thus, satisfy \eqref{equ_conds_basis_matrices_SDPCM}, i.e., $\mathbf{C}_{k} \mathbf{C}_{k'} = \delta_{k,k'} \mathbf{C}_{k}$ for all $k,k' \in [N]$. 
This special case, which can be seen as the SPCM analogue of the SSNM,\footnote{Indeed, by a suitable transformation of the observation vector $\mathbf{y}$ of the SDPCM, 
we have that within the SDPCM and the SSNM, each entry $x_{k}$ of the parameter vector determines the statistics of an individual subset of the entries of the observation vector. The difference between the SDPCM and 
the SSNM is then how the parameter vector entry $x_{k}$ determines the statistics of the corresponding observation vector entries. In the SSNM $x_{k}$ 
determines the mean, while in the SDPCM $x_{k}$ determines the variance. Another difference is that for the the SSNM, the subset of 
observation vector entries that correspond to $x_{k}$ consist only of the single entry $y_{k}$, while for the SDPCM these subsets may contain more than one entry of the observation vector.} will allow for a slightly more 
precise analysis via RKHS theory than was possible for the general SPCM. 
In particular, we will derive lower bounds that are tighter than the bound given by Theorem \ref{thm_CRB_SPCM} when specialized to 
the SDPCM. 
We will use the generic notation $\mathcal{M}_{\text{SDPCM}} = \left( \mathcal{E}_{\text{SDPCM}},c(\cdot) \equiv 0, \mathbf{x}_{0} \right)$ 
for a minimum variance problem that is obtained from the SDPCM $\mathcal{E}_{\text{SDPCM}} = \left( \mathcal{X}_{S,+},f_{\text{SDPCM}}(\mathbf{y}; \mathbf{x}), g(\mathbf{x})  \right)$, the identically vanishing prescribed bias (i.e., we consider unbiased estimation) 
and a parameter vector $\mathbf{x}_{0} \in \mathcal{X}_{S,+}$. 

Since the basis matrices $\mathbf{C}_{k} \in \mathbb{R}^{MÊ\times M}$ are assumed psd, the condition \eqref{equ_conds_basis_matrices_SDPCM} is equivalent to requiring 
the basis matrices to be orthogonal projection matrices on orthogonal subspaces of $\mathbb{R}^{M}$. 
The condition \eqref{equ_conds_basis_matrices_SDPCM} evaluated for $k'=k$, i.e., $\mathbf{C}_{k}^{2} = \mathbf{C}_{k}$, is the defining requirement for the psd matrix $\mathbf{C}_{k}$ 
to be an orthogonal projection matrix (cf.\ \cite{golub96}), i.e., it can be written as 
\begin{equation} 
\label{equ_basis_matrix_sum_rank_1_projection_SDPCM}
\mathbf{C}_{k} = \sum_{i \in [r_{k}]} \mathbf{u}_{m_{k,i}} \mathbf{u}_{m_{k,i}}^{T}.  
\end{equation}
Here, $\{ \mathbf{u}_{m} \in \mathbb{R}^{M} \}_{m \in [R] }$ is an orthonormal set of vectors, i.e., $\mathbf{u}_{m}^{T}  \mathbf{u}_{m'} = \delta_{m,m'}$, and the 
sets $\mathcal{U}_{k} \triangleq \{ \mathbf{u}_{m_{k,i}} \}_{i \in [r_{k}]}$ are disjoint, so that they span orthogonal subspaces of $\mathbb{R}^{M}$.
Indeed, assume that for two different basis matrices $\mathbf{C}_{k}$, $\mathbf{C}_{k'}$ ($k \neq k'$), there would be at least two vectors 
$\mathbf{a} \in \mathcal{U}_{k}$, $\mathbf{b} \in \mathcal{U}_{k'}$ which are nonorthogonal, i.e., 
$\mathbf{a}^{T} \mathbf{b} \neq 0$. Then, using \eqref{equ_basis_matrix_sum_rank_1_projection_SDPCM}, we would obtain 
$\mathbf{C}_{k} \mathbf{C}_{k'} \mathbf{b} = \mathbf{C}_{k} \mathbf{b} \neq \mathbf{0}$, which is a contradiction to \eqref{equ_conds_basis_matrices_SDPCM}. 
From the representation  \eqref{equ_basis_matrix_sum_rank_1_projection_SDPCM}, it also follows that for the SDPCM, there must be $M \geq R= \sum_{k \in [N]} r_{k}$ (in particular, $M \geq N$). 

Using the representation \eqref{equ_basis_matrix_sum_rank_1_projection_SDPCM} of the basis matrices, we can interpret the observation model of the SDPCM in \eqref{equ_obs_model_SPCM} as a 
latent variable model 
\begin{equation} 
\mathbf{y} = \sum_{k \in [N]} \mathbf{s}_{k}  + \mathbf{n},
\end{equation} 
with $\mathbf{s}_{k} = \sum_{i \in [r_{k}]} \xi_{m_{k,i}} \mathbf{u}_{m_{k,i}}$, where the $\xi_{m_{k,i}}$ are independent zero-mean Gaussian random variables with variance $x_{k}$ 
for all $i$, i.e., $\xi_{m_{k,i}} \sim \mathcal{N}(0,x_{k})$. This is similar to the latent variable model used in probabilistic principal component analysis \cite{TippingProbPCA} 
except that our ``factors'' $\mathbf{u}_{m}$ are fixed. 

Another useful implication of \eqref{equ_basis_matrix_sum_rank_1_projection_SDPCM} is based on a fundamental result in linear algebra \cite{HalmosFiniteDimVecSpace} which states 
that we can add vectors $\{ \mathbf{v}_{l} \in \mathbb{R}^{M} \}_{l \in [L]}$, with $L = M-R$, to the set $\{ \mathbf{u}_{m} \in \mathbb{R}^{M} \}_{m \in [R] }$ such that their union 
$\mathcal{U} \triangleq \{ \mathbf{v}_{l} \}_{l \in [L]} \cup \{ \mathbf{u}_{m}  \}_{m \in [R]}$ forms 
an ONB for $\mathbb{R}^{M}$. This implies in particular that we can represent the covariance matrix of the observation $\widetilde{\mathbf{C}}(\mathbf{x}) = \mathbf{C}(\mathbf{x}) + \sigma^{2} \mathbf{I}$ (cf.\ \eqref{equ_def_cov_matrix_observation_SPCM}), as 
\begin{equation}
\label{equ_expr_obs_cov_matrix_SDPCM}
\widetilde{\mathbf{C}}(\mathbf{x}) = \mathbf{U} \big[ \mathbf{D}'(\mathbf{x}) + \sigma^{2} \mathbf{I} \big] \mathbf{U}^{T}, 
\end{equation} 
where the columns of the orthonormal matrix $\mathbf{U} \in \mathbb{R}^{M \times M}$ are the vectors belonging to $\mathcal{U}$.
The diagonal matrix $\mathbf{D}'(\mathbf{x}) \in \mathbb{R}^{MÊ\times M}$ has a main diagonal consisting of two parts: the first part of length $R$ contains the entries $x_{k}$ where each entry is repeated $r_{k}$ times, 
while the second part consists of $M - R$ zeros, i.e.,
\begin{equation}
\label{equ_matrix_D_prime_SDPCM}
\mathbf{D}'(\mathbf{x}) = \begin{pmatrix} x_{1} \mathbf{I}_{r_{1}} & \mathbf{0} & \ldots & \ldots & \ldots & \ldots & \mathbf{0} \\ 
\mathbf{0} & x_{2} \mathbf{I}_{r_{2}} & \ddots &  \ddots & \ddots &\ddots &  \vdots \\ 
\vdots & \ddots& \ddots &  \ddots& \ddots &\ddots &  \vdots \\ 
\vdots & \ddots& \ddots &  x_{N} \mathbf{I}_{r_{N}} & \ddots &\ddots & \vdots \\ 
\vdots & \ddots & \ddots &   \ddots & \mathbf{0} &\ddots& \vdots 
\\ 
\vdots & \ddots & \ddots &   \ddots & \ddots &\ddots& \vdots \\ 
\mathbf{0} & \ldots & \ldots &   \ldots & \ldots &\ldots& \mathbf{0} 
\end{pmatrix}.
\end{equation} 


\subsection{The $\mathcal{D}$-restricted SDPCM}

As for the SPCM, we cannot associate a RKHS to the minimum variance problem $\mathcal{M}_{\text{SDPCM}}$ directly. Instead, 
we will work with the RKHS associated to the $\mathcal{D}$-restricted minimum variance problem $\mathcal{M}_{\mathcal{D}, \text{SDPCM}} = \mathcal{M}_{\text{SDPCM}}\big|_{\mathcal{D}}$, 
with a set $\mathcal{D} \subseteq \mathcal{X}_{S,+}$ that satisfies \eqref{equ_nec_suff_cond_D_SPCM_kernel_exists}. 

For the SDPCM, the condition \eqref{equ_nec_suff_cond_D_SPCM_kernel_exists} can be rewritten due to \eqref{equ_expr_obs_cov_matrix_SDPCM} as 
\begin{equation}
\label{equ_nec_suff_cond_D_SDPCM_kernel_exists}
 2 ( x_{i}+\sigma^{2})^{-1}  - (x_{0,i}+\sigma^{2} )^{-1} > 0 \quad \forall \mathbf{x}\in \mathcal{D}, \forall i \in [N]. 
\end{equation}
This is because a matrix $\mathbf{A} \in \mathbb{R}^{M \times M}$ is positive definite, i.e., $\mathbf{A} > \mathbf{0}$, if and only if $\mathbf{V} \mathbf{A} \mathbf{V}^{T} > \mathbf{0}$ 
for any orthonormal matrix $\mathbf{V} \in \mathbb{R}^{M \times M}$, i.e., satisfying $\mathbf{V}^{T} \mathbf{V} = \mathbf{I}$. 
The largest set $\mathcal{D} \subseteq \mathcal{X}_{S,+}$ which satisfies \eqref{equ_nec_suff_cond_D_SDPCM_kernel_exists} is obviously given by 
\begin{equation} 
\mathcal{D}_{0} \triangleq \{ \mathbf{x} \in \mathcal{X}_{S,+} \big|  x_{i} < 2 x_{0,i} + \sigma^{2} \} . 
\end{equation} 
Note that the set $\mathcal{D}_{0}$ depends on $\mathbf{x}_{0}$. 

Since the set $\mathcal{D}_{0}$ is the largest set satisfying \eqref{equ_nec_suff_cond_D_SDPCM_kernel_exists}, the corresponding minimum achievable variance $L_{\mathcal{M}_{\mathcal{D}_{0}, \text{SDPCM}}}$ 
will yield via \eqref{equ_lower_bound_spcm_restr_spcm} the tightest lower bound on $L_{\mathcal{M}_{\text{SDPCM}}}$ among all bounds 
that are obtained from $L_{\mathcal{M}_{\mathcal{D}, \text{SDPCM}}}$ with a set $\mathcal{D}$ that satisfies \eqref{equ_nec_suff_cond_D_SDPCM_kernel_exists}. 









\subsection{Variance Bounds for the SDPCM}

We will now derive a  lower bound on the minimum achievable variance $L_{\mathcal{M}_{\text{SDPCM}}}$ by lower bounding $L_{\mathcal{M}_{\mathcal{D}_{0}, \text{SDPCM}}}$, i.e., 
the minimum achievable variance of $\mathcal{M}_{\mathcal{D}_{0}, \text{SDPCM}}$, using the RKHS $\mathcal{H}(\mathcal{M}_{\mathcal{D}_{0}, \text{SDPCM}})$ which is 
associated to the kernel $R_{\mathcal{M}_{\mathcal{D}_{0}, \text{SDPCM}}}(\cdot, \cdot): \mathcal{D}_{0} \times \mathcal{D}_{0} \rightarrow \mathbb{R}$: 
\begin{align} 
\label{equ_def_kernel_SDPCM_D_0_def_part}
R_{\mathcal{M}_{\mathcal{D}_{0},\text{SDPCM}}} (\mathbf{x}_{1}, \mathbf{x}_{2})
& \triangleq \mathsf{E}_{\mathbf{x}_{0}}  \left\{ \rho_{\mathcal{M}_{\mathcal{D}_{0},\text{SDPCM}}}(\mathbf{y} , \mathbf{x}_{1}) \rho_{\mathcal{M}_{\mathcal{D}_{0},\text{SDPCM}}}(\mathbf{y},\mathbf{x}_{2}) \right\}. 
\end{align} 
Using \eqref{equ_kernel_D_restr_SPCM} and \eqref{equ_expr_obs_cov_matrix_SDPCM} and the identity $\detm{\mathbf{A} \mathbf{B}} = \detm{\mathbf{A}} \detm{\mathbf{B}}$ 
valid for any two matrices $\mathbf{A}, \mathbf{B} \in \mathbb{R}^{M \times M}$ \cite[p.\ 51]{golub96}, we obtain
\begin{align}
\label{equ_def_kernel_SDPCM_D_0}
R_{\mathcal{M}_{\mathcal{D}_{0},\text{SDPCM}}} (\mathbf{x}_{1}, \mathbf{x}_{2}) & = \big[ \detmb{ \widetilde{\mathbf{C}}(\mathbf{x}_{0}) } \big]^{1/2} \ist 
  \big[ \detmb{ \widetilde{\mathbf{C}}(\mathbf{x}_{1}) + \widetilde{\mathbf{C}}(\mathbf{x}_{2}) - \widetilde{\mathbf{C}}(\mathbf{x}_{1}) \widetilde{\mathbf{C}}^{-1}(\mathbf{x}_{0})\widetilde{\mathbf{C}}(\mathbf{x}_{2}) } \big]^{-1/2} 
 \nonumber \\[5mm]
 & =  \frac{ \prod\limits_{k \in [N]} (x_{0,k} \rmv\rmv+\rmv \sigma^{2})^{r_{k}} }{ \prod\limits_{k \in [N]} 
  \!\big[ (x_{0,k} \rmv\rmv+\rmv \sigma^{2})^2 - (x_{1,k} \!-\rmv\rmv x_{0,k})(x_{2,k} \!-\rmv\rmv x_{0,k})\big]^{r_k/2} }. 
\end{align} 
 Note that the kernel $R_{\mathcal{M}_{\mathcal{D}_{0}, \text{SDPCM}}}(\cdot, \cdot)$ given by \eqref{equ_def_kernel_SDPCM_D_0} is obviously differentiable up to any order $m$. 

The following result partly characterizes the RKHS $\mathcal{H}(\mathcal{M}_{\mathcal{D}_{0}, \text{SDPCM}})$. 
\begin{theorem}
\label{thm_charac_RKHS_SDPCM_D_0}
The RKHS $\mathcal{H}(\mathcal{M}_{\mathcal{D}_{0}, \emph{SDPCM}})$ is differentiable up to any order. Let $\mathcal{K} \subseteq [N]$ be an arbitrary set of $S$ different indices, i.e., $|\mathcal{K}| = S$, 
and let $\mathbf{p} \in \mathbb{Z}_{+}^{N}$ be an arbitrary multi-index with $\supp(\mathbf{p}) \subseteq \mathcal{K}$.
Then the function $g^{(\mathbf{p}, \mathcal{K})}(\cdot): \mathcal{D}_{0} \rightarrow \mathbb{R}$ defined as  
\vspace*{2mm}
\begin{align} 
\label{equ_def_g_p_thm_charac_RKHS_SDPCM_D_0}
g^{(\mathbf{p}, \mathcal{K})}(\mathbf{x}) & 
\triangleq \frac{ \partial^{\mathbf{p}} R_{\mathcal{M}_{\mathcal{D}_{0}, \emph{SDPCM}}}(\mathbf{x}, \mathbf{x}_{2})}{\partial \mathbf{x}_{2}^{\mathbf{p}}}\bigg|_{\mathbf{x}_{2} = \mathbf{x}_{0}^{\mathcal{K}}} \nonumber \\[4mm]
& = \prod_{k \in \supp(\mathbf{p})}  \frac{c_{k} (x_k - x_{0,k})^{p_{k}} }{ (x_{0,k} \rmv\rmv+\rmv \sigma^{2})^{2p_{k}} } 
 \prod_{k \in \supp(\mathbf{x}_{0}) \setminus \mathcal{K}} \frac{(x_{0,k} \rmv\rmv+\rmv \sigma^{2})^{r_{k}} }{\!\big[ (x_{0,k} \rmv\rmv+\rmv \sigma^{2})^2 + (x_{k} \!-\rmv\rmv x_{0,k})\rmv\rmv x_{0,k}\big]^{r_k/2} } ,
\end{align}
where $c_{k} \triangleq \prod_{l \in [p_{k}]} (r_{k}/2 + (l-1))$, is an element of $\mathcal{H}(\mathcal{M}_{\mathcal{D}_{0}, \emph{SDPCM}})$, i.e., 
\begin{equation} 
\label{equ_charac_RKHS_SDPCM_D_0_contains_part_der}
g^{(\mathbf{p},\mathcal{K})}(\cdot) \in \mathcal{H}(\mathcal{M}_{\mathcal{D}_{0}, \emph{SDPCM}}).
\end{equation}
The inner product of an arbitrary function $f(\cdot) \in \mathcal{H}(\mathcal{M}_{\mathcal{D}_{0}, \emph{SDPCM}})$ with $g^{(\mathbf{p}, \mathcal{K})}(\cdot)$ is given by 
\begin{equation} 
\label{equ_inner_prod_part_der_RKHS_Re_SDPCM}
\big\langle f (\cdot), g^{(\mathbf{p}, \mathcal{K})}(\cdot) \big\rangle_{\mathcal{H}(\mathcal{M}_{\mathcal{D}_{0}, \emph{SDPCM}})} =  \frac{ \partial^{\mathbf{p}} f(\mathbf{x})} {\partial \mathbf{x}^{\mathbf{p}}} \bigg|_{\mathbf{x} = \mathbf{x}_{0}^{\mathcal{K}}}. 
\end{equation} 
\end{theorem}
\begin{proof}
The statement follows straightforwardly from Theorem \ref{thm_der_repr_prop}, since the kernel $R_{\mathcal{M}_{\mathcal{D}_{0}, \text{SDPCM}}}(\cdot, \cdot)$ 
is differentiable up to any order.
\end{proof}

Based on Theorem \ref{thm_charac_RKHS_SDPCM_D_0}, we will now derive a lower bound on the minimum achievable variance $L_{\mathcal{M}_{\text{SDPCM}}}$: 
\begin{theorem}
\label{thm_sparse_lower_bound_SDPCM_ICASSP}
Consider the minimum variance problem $\mathcal{M}_{\emph{SDPCM}}=\left(\mathcal{E}_{\emph{SDPCM}},c(\cdot) \equiv 0 , \mathbf{x}_{0}Ê\right)$ 
with sparsity degree $S$. 
Let $\{ \mathbf{p}_{l} \in \mathbb{Z}_{+}^{N} \}_{l \in [L]}$ be a set of multi-indices such that 
$\supp(\mathbf{p}_{l}) \subseteq \mathcal{K}$, where $\mathcal{K}\subseteq [N]$ 
is an arbitrary set of $S$ different indices, i.e., $| \mathcal{K} | = S$. 
If the parameter function $g(\cdot):\mathcal{X}_{S,+} \rightarrow \mathbb{R}$ associated with $\mathcal{E}_{\emph{SDPCM}}$ is such that 
the partial derivatives $\frac{\partial^{\mathbf{p}_{l}} g(\mathbf{x})}{\partial \mathbf{x}^{\mathbf{p}_{l}}} \big|_{\mathbf{x} = \mathbf{x}_{0}^{\mathcal{K}}}$ exist for all $l \in [L]$, 
we have 
\begin{align} 
\label{equ_lower_bound_SPCM_ICASSP}
L_{\mathcal{M}_{\emph{SDPCM}}} \geq  \sum_{l \in [L]} \frac{1}{q_{l}(\mathbf{x}_{0})} \bigg[  \frac{\partial^{\mathbf{p}_{l}} g(\mathbf{x})}{\partial \mathbf{x}^{\mathbf{p}_{l}}}\bigg|_{\mathbf{x} = \mathbf{x}_{0}^{\mathcal{K}}}
\bigg]^{2} - \big[ g(\mathbf{x}_{0}) \big]^{2},
\end{align} 
where $q_{l}(\mathbf{x}_{0}) \triangleq 
  \frac{ \partial^{\mathbf{p}_{l}}  \partial^{\mathbf{p}_{l}} R_{\mathcal{M}_{\mathcal{D}_{0},\emph{SDPCM}}}(\mathbf{x}_{1}, \mathbf{x}_{2})}
  {\partial \mathbf{x}_{1}^{\mathbf{p}_{l}}\partial \mathbf{x}_{2}^{\mathbf{p}_{l}} }\bigg|_{\mathbf{x}_{1} = \mathbf{x}_{2} = \mathbf{x}_{0}^{\mathcal{K}}}$. 
Furthermore, there exists an unbiased estimator $\hat{g}(\cdot)$ whose variance at $\mathbf{x}_{0}$ achieves this bound, i.e., $v(\hat{g}(\cdot);\mathbf{x}_{0})=  \sum_{l \in [L]} \frac{1}{q_{l}(\mathbf{x}_{0})} \Big[  \frac{\partial^{\mathbf{p}_{l}} g(\mathbf{x})}{\partial \mathbf{x}^{\mathbf{p}_{l}}}\big|_{\mathbf{x} = \mathbf{x}_{0}^{\mathcal{K}}}\Big]^{2} - \big[ g(\mathbf{x}_{0}) \big]^{2}$, if and only if it can be written as 
\begin{equation}
\label{equ_nec_suff_cond_est_attains_bound_SDPCM_ICASSP} 
\hat{g}(\cdot) = \sum_{l \in [L]Ê} a_{l}   \frac{ \partial^{\mathbf{p}_{l}}  \rho_{\mathcal{M}_{\mathcal{D}_{0},\emph{SDPCM}}}(\cdot, \mathbf{x})}{\partial \mathbf{x}^{\mathbf{p}_{l}} }\bigg|_{\mathbf{x}  = \mathbf{x}_{0}^{\mathcal{K}}}, 
\end{equation} 
with suitable (non-random) coefficients $a_{l} \in \mathbb{R}$.
\end{theorem}

\begin{proof}
Since the bound in \eqref{equ_lower_bound_SPCM_ICASSP} is always finite, we assume without loss of generality that $g(\cdot)$ is estimable for $\mathcal{M}_{\text{SDPCM}}$, 
which implies via Theorem \ref{thm_par_set_reduction_classic_est_mve} that 
the restriction $g(\cdot)\big|_{\mathcal{D}_{0}}$ is estimable for $\mathcal{M}_{\mathcal{D}_{0}, \text{SDPCM}}$. 
Thus we can assume by Theorem \ref{thm_main_facts_RKHS_MVE} that the prescribed mean function 
$\gamma(\cdot): \mathcal{D}_{0} \rightarrow \mathbb{R}: \gamma(\mathbf{x}) = c(\mathbf{x}) + g(\mathbf{x})= g(\mathbf{x})$ belongs to the RKHS $\mathcal{H}(\mathcal{M}_{\mathcal{D}_{0},\text{SDPCM}})$. 

Let us consider the subspace $\mathcal{U} \triangleq \linspan \big\{ g^{(\mathbf{p}_{l},\mathcal{K})} (\cdot) \big\}_{l \in [L]}  \subseteq \mathcal{H}(\mathcal{M}_{\mathcal{D}_{0},\text{SDPCM}})$ spanned 
by the functions 
$g^{(\mathbf{p}_{l},\mathcal{K})} (\cdot) \in \mathcal{H}(\mathcal{M}_{\mathcal{D}_{0},\text{SDPCM}})$ as defined in Theorem \ref{thm_charac_RKHS_SDPCM_D_0}.
For the inner products between $g^{(\mathbf{p}_{l},\mathcal{K})} (\cdot)$ and $g^{(\mathbf{p}_{l'},\mathcal{K})} (\cdot)$, we obtain 
\vspace*{3mm}
\begin{align}
\label{equ_lower_bound_SPCM_ICASSP_inner_prod_g_p_general}
\big\langle g^{(\mathbf{p}_{l},\mathcal{K})} (\cdot) , g^{(\mathbf{p}_{l'},\mathcal{K})} (\cdot) \big\rangle_{ \mathcal{H}(\mathcal{M}_{\mathcal{D}_{0},\text{SDPCM}})}  
& \stackrel{\eqref{equ_inner_prod_part_der_RKHS_Re_SDPCM}}{=} 
 \frac{ \partial^{\mathbf{p}_{l'}} g^{(\mathbf{p}_{l},\mathcal{K})} (\mathbf{x}) }
  {\partial \mathbf{x}^{\mathbf{p}_{l'}} }\bigg|_{\mathbf{x} = \mathbf{x}_{0}^{\mathcal{K}}} \nonumber \\[4mm]
  & \hspace*{-50mm} \stackrel{\eqref{equ_def_g_p_thm_charac_RKHS_SDPCM_D_0}}{=} 
 \frac{ \partial^{\mathbf{p}_{l'}} } {\partial \mathbf{x}^{\mathbf{p}_{l'}}}
 \prod_{k \in \supp(\mathbf{p}_{l})}   \frac{c_{k} (x_{k} - x_{0,k})^{p_{k}}}{ (x_{0,k} \rmv\rmv+\rmv \sigma^{2})^{2p_{k}} } 
 \prod_{k \in \supp(\mathbf{x}_{0}) \setminus \mathcal{K}} \frac{(x_{0,k} \rmv\rmv+\rmv \sigma^{2})^{r_{k}} }{\!\big[ (x_{0,k} \rmv\rmv+\rmv \sigma^{2})^2 + (x_{k} \!-\rmv\rmv x_{0,k})\rmv\rmv x_{0,k}\big]^{r_k/2} }\bigg|_{\mathbf{x} = \mathbf{x}_{0}^{\mathcal{K}}}  \nonumber \\[4mm]
  & \hspace*{-50mm} \stackrel{(a)}{=}Ê\Big[ \prod_{k \in \supp(\mathbf{x}_{0}) \setminus \mathcal{K}} \frac{(x_{0,k} \rmv\rmv+\rmv \sigma^{2})^{r_{k}} }{\!\big[ (x_{0,k} \rmv\rmv+\rmv \sigma^{2})^2  \!-\rmv\rmv x_{0,k}^{2}\big]^{r_k/2} }  \Big]
\frac{ \partial^{\mathbf{p}_{l'}} } {\partial \mathbf{x}^{\mathbf{p}_{l'}}}
 \prod_{k \in \supp(\mathbf{p}_{l})}   \frac{c_{k} (x_{k} - x_{0,k})^{p_{k}}}{ (x_{0,k} \rmv\rmv+\rmv \sigma^{2})^{2p_{k}} } \bigg|_{\mathbf{x} = \mathbf{x}_{0}^{\mathcal{K}}},
\end{align}
where step $(a)$ is due to the fact that $\supp(\mathbf{p}_{l}),\supp(\mathbf{p}_{l'}) \subseteq \mathcal{K}$. 
Evaluating \eqref{equ_lower_bound_SPCM_ICASSP_inner_prod_g_p_general} for $l \neq l'$ reveals that the functions $\big\{ g^{(\mathbf{p}_{l},\mathcal{K})} (\cdot) \big\}_{l \in [L]}$ are orthogonal, i.e.
\begin{align} 
\label{equ_proof_ICASSP_SDPCM_orthogonality_q_l}
\big\langle g^{(\mathbf{p}_{l},\mathcal{K})} (\cdot) , g^{(\mathbf{p}_{l'},\mathcal{K})} (\cdot) \big\rangle_{ \mathcal{H}(\mathcal{M}_{\mathcal{D}_{0},\text{SDPCM}})} & = 
\delta_{l,l'}\big\langle g^{(\mathbf{p}_{l},\mathcal{K})} (\cdot) , g^{(\mathbf{p}_{l},\mathcal{K})} (\cdot) \big\rangle_{ \mathcal{H}(\mathcal{M}_{\mathcal{D}_{0},\text{SDPCM}})}  \nonumber \\[4mm] 
& \stackrel{\eqref{equ_inner_prod_part_der_RKHS_Re_SDPCM}}{=} 
 \delta_{l,l'}\frac{ \partial^{\mathbf{p}_{l'}} g^{(\mathbf{p}_{l},\mathcal{K})} (\mathbf{x}) }
  {\partial \mathbf{x}^{\mathbf{p}_{l'}} }\bigg|_{\mathbf{x} = \mathbf{x}_{0}^{\mathcal{K}}} \nonumber \\[4mm]
  & = \delta_{l,l'}  \frac{ \partial^{\mathbf{p}_{l}}  \partial^{\mathbf{p}_{l}} R_{\mathcal{M}_{\mathcal{D}_{0},\emph{SDPCM}}}(\mathbf{x}_{1}, \mathbf{x}_{2})}
  {\partial \mathbf{x}_{1}^{\mathbf{p}_{l}}\partial \mathbf{x}_{2}^{\mathbf{p}_{l}} }\bigg|_{\mathbf{x}_{1} = \mathbf{x}_{2} = \mathbf{x}_{0}^{\mathcal{K}}}  \nonumber \\[4mm]
 & =  \delta_{l,l'} q_{l}(\mathbf{x}_{0}).
\end{align}  
Therefore, an ONB for the subspace $\mathcal{U}\subseteq  \mathcal{H}(\mathcal{M}_{\mathcal{D}_{0},\text{SDPCM}})$ is given by the set $\Big\{ \frac{1}{\sqrt{q_{l}(\mathbf{x}_{0})}} g^{(\mathbf{p}_{l},\mathcal{K})} (\cdot) \Big\}_{l \in [L]}$.

Using the inner products (recall that $\gamma(\cdot) = g(\cdot)$) 
\begin{equation} 
\label{equ_inner_prod_part_der_gamm_SPCM_ICASSP_lower_bound}
\big\langle \gamma(\cdot),g^{(\mathbf{p}_{l},\mathcal{K})} (\cdot) \big\rangle_{ \mathcal{H}(\mathcal{M}_{\mathcal{D}_{0},\text{SDPCM}})} 
\stackrel{\eqref{equ_inner_prod_part_der_RKHS_Re_SDPCM}}{=} 
\frac{\partial^{\mathbf{p}_{l}} g(\mathbf{x})}{\partial \mathbf{x}^{\mathbf{p}_{l}}} \bigg|_{\mathbf{x} = \mathbf{x}_{0}^{\mathcal{K}}},
\end{equation}
the bound in \eqref{equ_lower_bound_SPCM_ICASSP} follows by projecting $\gamma(\cdot)$ onto the subspace $\mathcal{U}$, since 
\begin{align}
\label{equ_proof_lower_bound_SDPCM_ICASSP_lower_bound}
 L_{\mathcal{M}_{\text{SDPCM}}} & \stackrel{\eqref{equ_lower_bound_spcm_restr_spcm}}{\geq} 
 L_{\mathcal{M}_{\mathcal{D}_{0},\text{SDPCM}}} \stackrel{\eqref{equ_squared_norm_min_achiev_var}}{=} \| \gamma(\cdot) \|^{2}_{\mathcal{H}(\mathcal{M}_{\mathcal{D}_{0},\text{SDPCM}})}
 - \big[\underbrace{\gamma(\mathbf{x}_{0})}_{= g(\mathbf{x}_{0})}\big]^{2}  \nonumber \\[4mm] 
& \stackrel{\eqref{equ_squared_norm_orthog_proj_pythag_thm}}{\geq}  \| \mathbf{P}_{\mathcal{U}}\gamma(\cdot)  \|^{2}_{\mathcal{H}(\mathcal{M}_{\mathcal{D}_{0},\text{SDPCM}})} - \big[g(\mathbf{x}_{0})\big]^{2}  \nonumber \\[4mm]
&  \stackrel{\eqref{equ_norm_projection_finite_dim_subspace}}{=}  
\sum_{l \in [L]}   \frac{ \big\langle \gamma(\cdot), g^{(\mathbf{p}_{l},\mathcal{K})}(\cdot) \big\rangle_{\mathcal{H}(\mathcal{M}_{\mathcal{D}_{0},\text{SPCM}})}^{2}}{ \big\langle g^{(\mathbf{p}_{l},\mathcal{K})}(\cdot) , g^{(\mathbf{p}_{l},\mathcal{K})}(\cdot) \big\rangle_{\mathcal{H}(\mathcal{M}_{\mathcal{D}_{0},\text{SPCM}})}}
 - \big[g(\mathbf{x}_{0})\big]^{2}  \nonumber \\[4mm] 
&  \stackrel{\eqref{equ_proof_ICASSP_SDPCM_orthogonality_q_l},\eqref{equ_inner_prod_part_der_gamm_SPCM_ICASSP_lower_bound}}{=} \sum_{l \in [L]} \frac{1}{q_{l}(\mathbf{x}_{0})} \bigg[ \frac{\partial^{\mathbf{p}_{l}} g(\mathbf{x})}{\partial \mathbf{x}^{\mathbf{p}_{l}}}\bigg|_{\mathbf{x} = \mathbf{x}_{0}^{\mathcal{K}}}
\bigg]^{2} - \big[ g(\mathbf{x}_{0}) \big]^{2}.
\end{align} 

Let us now show that if an unbiased estimator for $\mathcal{M}_{\text{SDPCM}}$ achieves the bound \eqref{equ_lower_bound_SPCM_ICASSP}, i.e., its variance at $\mathbf{x}_{0}$ equals 
$\sum_{l \in [L]} \frac{1}{q_{l}(\mathbf{x}_{0})} \Big[  \frac{\partial^{\mathbf{p}_{l}} g(\mathbf{x})}{\partial \mathbf{x}^{\mathbf{p}_{l}}}\big|_{\mathbf{x} = \mathbf{x}_{0}^{\mathcal{K}}}\Big]^{2} - \big[ g(\mathbf{x}_{0}) \big]^{2}$ (which also implies that $g(\cdot)$ is estimable for  $\mathcal{M}_{\text{SDPCM}}$), it must be necessarily of the form \eqref{equ_nec_suff_cond_est_attains_bound_SDPCM_ICASSP}. 
To that end, we note that according to the above derivation of \eqref{equ_proof_lower_bound_SDPCM_ICASSP_lower_bound}, the bound \eqref{equ_lower_bound_SPCM_ICASSP} holds also for $L_{\mathcal{M}_{\mathcal{D}_{0},\text{SDPCM}}}$, i.e., 
\begin{equation} 
\label{equ_proof_lower_bound_ICASSP_SDPCM_existence_LMV_bound_LM_D0}
L_{\mathcal{M}_{\mathcal{D}_{0},\text{SDPCM}}} \geq  \sum_{l \in [L]} \frac{1}{q_{l}(\mathbf{x}_{0})} \bigg[  \frac{\partial^{\mathbf{p}_{l}} g(\mathbf{x})}{\partial \mathbf{x}^{\mathbf{p}_{l}}}\bigg|_{\mathbf{x} = \mathbf{x}_{0}^{\mathcal{K}}}
\bigg]^{2} -  \big[ g(\mathbf{x}_{0}) \big]^{2}.
\end{equation} 
Therefore, if an unbiased estimator exists whose variance at $\mathbf{x}_{0}$ achieves the bound \eqref{equ_lower_bound_SPCM_ICASSP}, it also achieves the bound \eqref{equ_proof_lower_bound_ICASSP_SDPCM_existence_LMV_bound_LM_D0}, and thus is the unique LMV estimator (which is guaranteed to exist in this case) for $\mathcal{M}_{\mathcal{D}_{0},\text{SDPCM}}$. 
The variance at $\mathbf{x}_{0}$ of this LMV estimator, which is of course also the minimum achievable variance $L_{\mathcal{M}_{\mathcal{D}_{0},\text{SDPCM}}}$, 
attains the bound \eqref{equ_lower_bound_SPCM_ICASSP} if and only if the projection $\mathbf{P}_{\mathcal{U}} \gamma(\cdot)$ 
coincides with $\gamma(\cdot)$, since this is necessary and sufficient for the second inequality in \eqref{equ_proof_lower_bound_SDPCM_ICASSP_lower_bound} to become an equality. 
This in turn, is the case if and only if $\gamma(\cdot) \in \mathcal{U}$, i.e., the prescribed mean function $\gamma(\cdot)$ 
can be written as the linear basis expansion $\sum_{l \in [L]} a_{l} g^{(\mathbf{p}_{l},\mathcal{K})} (\cdot)$ with coefficients $a_{l} \in \mathbb{R}$. 
However, in this case we can express the unique LMV estimator $\hat{g}^{(\mathbf{x}_{0})}(\cdot)$ for $\mathcal{M}_{\mathcal{D}_{0},\text{SDPCM}}$ via Theorem \ref{thm_isometry_RKHS_rhos_derivative_kernel} and Theorem \ref{thm_main_facts_RKHS_MVE} as 
\begin{equation}
\hat{g}^{(\mathbf{x}_{0})}(\cdot) \stackrel{\eqref{equ_lmv_estimator_general_congruence_L_M_RKHS}}{=} \mathsf{J} \big[ \gamma(\cdot) \big] = \sum_{l \in [L]Ê} a_{l}  \mathsf{J} \Big[   g^{(\mathbf{p}_{l},\mathcal{K})} (\cdot) \Big]
\stackrel{\eqref{equ_def_g_p_thm_charac_RKHS_SDPCM_D_0},\eqref{equ_isometry_RKHS_rhos_derivative_kernel}}{=}  \sum_{l \in [L]Ê} a_{l}   \frac{ \partial^{\mathbf{p}_{l}}  \rho_{\mathcal{M}_{\mathcal{D}_{0},\text{SDPCM}}}(\cdot, \mathbf{x})}{\partial \mathbf{x}^{\mathbf{p}_{l}} }\bigg|_{\mathbf{x}  = \mathbf{x}_{0}^{\mathcal{K}}}. 
\end{equation} 
Thus, if an estimator achieves the bound \eqref{equ_lower_bound_SPCM_ICASSP}, it is necessarily of the form \eqref{equ_nec_suff_cond_est_attains_bound_SDPCM_ICASSP}.  

Conversely, if there exists an estimator $\hat{g}(\cdot)$ of the form \eqref{equ_nec_suff_cond_est_attains_bound_SDPCM_ICASSP} we have by Theorem \ref{thm_isometry_RKHS_rhos} and Theorem \ref{thm_isometry_RKHS_rhos_derivative_kernel} that 
\begin{align}
\label{equ_mean_func_bound_ICASSP_SDPCM_sufficient_cond_LMV}
\gamma(\mathbf{x}) &  \stackrel{(a)}{=} \mathsf{E}_{\mathbf{x}} \{ \hat{g}(\mathbf{y}) \} =  \mathsf{E}_{\mathbf{x}} \Bigg\{ \sum_{l \in [L]Ê} a_{l}  \frac{ \partial^{\mathbf{p}_{l}}  \rho_{\mathcal{M}_{\mathcal{D}_{0},\text{SDPCM}}}(\mathbf{x}, \mathbf{x}_{2})}{\partial \mathbf{x}_{2}^{\mathbf{p}_{l}} }\bigg|_{\mathbf{x}_{2}  = \mathbf{x}_{0}^{\mathcal{K}}} \Bigg\}Ê\nonumber \\[4mm]
& =  \mathsf{E}_{\mathbf{x}} \Bigg\{ \sum_{l \in [L]Ê} a_{l}  \frac{ \partial^{\mathbf{p}_{l}}  \rho_{\mathcal{M}_{\mathcal{D}_{0},\text{SDPCM}}}(\mathbf{x}, \mathbf{x}_{2})}{\partial \mathbf{x}_{2}^{\mathbf{p}_{l}} }\bigg|_{\mathbf{x}_{2}  = \mathbf{x}_{0}^{\mathcal{K}}} \Bigg\}
\stackrel{\eqref{equ_def_g_p_thm_charac_RKHS_SDPCM_D_0},\eqref{equ_isometry_RKHS_rhos_derivative_kernel}}{=}  \mathsf{E}_{\mathbf{x}} \Bigg\{ \sum_{l \in [L]Ê} a_{l}  \mathsf{J}^{-1} \Big[   g^{(\mathbf{p}_{l},\mathcal{K})} (\cdot) \Big]  \Bigg\} \nonumber \\[4mm]
& \stackrel{\eqref{equ_expect_inversion_J_isometry_RKHS_rhos}}{=}  \sum_{l \in [L]Ê} a_{l}     g^{(\mathbf{p}_{l},\mathcal{K})} (\mathbf{x})  \quad \in \mathcal{U},
\end{align}  
where $(a)$ follows from the assumption that $\hat{g}(\cdot)$ is unbiased (note that we consider unbiased estimation since $c(\cdot) \equiv 0$).
Thus, the existence of an estimator $\hat{g}(\cdot)$ of the form \eqref{equ_nec_suff_cond_est_attains_bound_SDPCM_ICASSP} implies that $\gamma(\cdot) \in \mathcal{U}$, 
which, as already shown above, is necessary and sufficient for the variance of the estimator 
in \eqref{equ_nec_suff_cond_est_attains_bound_SDPCM_ICASSP},  i.e., $v(\hat{g}(\cdot),\mathbf{x}_{0})$ to achieve the bound \eqref{equ_lower_bound_SPCM_ICASSP}.  
\end{proof}

Note that by Theorem \ref{thm_equ_bias_param_function}, unbiased estimation of a parameter function $g(\mathbf{x})$ is equivalent to biased estimation of $x_{k}$ itself with bias $g(\mathbf{x}) - x_{k}$. 
Therefore, we have due to the definition of the minimum achievable variance 
the following implication of Theorem \ref{thm_sparse_lower_bound_SDPCM_ICASSP} for the variance of a biased estimator $\hat{x}_{k}(\mathbf{y})$ of $x_{k}$: 
\begin{corollary}
\label{cor_thm_sparse_lower_bound_SDPCM_ICASSP}
Consider an estimator $\hat{x}_{k}(\mathbf{y})$ which uses the observation $\mathbf{y}$ of a SDPCM with sparsity degree $S$ 
and denote its mean function by $m_{k}(\cdot): \mathcal{X}_{S,+} \rightarrow \mathbb{R}: m_{k} (\mathbf{x}) \triangleq \mathsf{E}_{\mathbf{x}} \big\{ \hat{x}_{k}(\mathbf{y}) \big\}$. 
Let $\{ \mathbf{p}_{l} \in \mathbb{Z}_{+}^{N} \}_{l \in [L]}$ be a set of multi-indices such that $\supp(\mathbf{p}_{l}) \subseteq \mathcal{K}$, where $\mathcal{K}\subseteq [N]$ 
is an arbitrary set of $S$ different indices, i.e., $| \mathcal{K} | = S$. If for a fixed parameter vector $\mathbf{x}_{0} \in \mathcal{X}_{S,+}$, the partial derivatives $\frac{\partial^{\mathbf{p}_{l}}m_{k}(\mathbf{x})}{\partial \mathbf{x}^{\mathbf{p}_{l}}} \big|_{\mathbf{x} = \mathbf{x}_{0}^{\mathcal{K}}}$ of the mean function exist, 
the variance $v(\hat{x}_{k}(\cdot); \mathbf{x}_{0})$ at $\mathbf{x}_{0}$ is lower bounded by 
\begin{equation} 
\label{equ_cor_thm_sparse_lower_bound_SDPCM_ICASSP}
v(\hat{x}_{k}(\cdot); \mathbf{x}_{0}) \geq \sum_{l \in [L]} \frac{1}{q_{l}(\mathbf{x}_{0})} \bigg[  \frac{\partial^{\mathbf{p}_{l}} m_{k}(\mathbf{x})}{\partial \mathbf{x}^{\mathbf{p}_{l}}}\bigg|_{\mathbf{x} = \mathbf{x}_{0}^{\mathcal{K}}}
\bigg]^{2} - \big[ m_{k}(\mathbf{x}_{0}) \big]^{2}
\end{equation} 
\end{corollary}

Let us now present another lower bound on $L_{\mathcal{M}_{\text{SDPCM}}}$ which depends only on the first-order partial derivatives of the parameter function $g(\cdot)$ associated 
with the minimum variance problem $\mathcal{M}_{\text{SDPCM}}=\left(\mathcal{E}_{\text{SDPCM}},c(\cdot) \equiv 0 , \mathbf{x}_{0}Ê\right)$. This bound  is derived by a slight modification of the proof of Theorem \ref{thm_sparse_lower_bound_SDPCM_ICASSP}.
\begin{theorem}
\label{thm_sparse_lower_bound_SDPCM_ICASSP_modifciaton}
Consider the minimum variance problem $\mathcal{M}_{\emph{SDPCM}}=\left(\mathcal{E}_{\emph{SDPCM}},c(\cdot) \equiv 0 , \mathbf{x}_{0}Ê\right)$ associated with 
the SDPCM $\mathcal{E}_{\emph{SDPCM}} = \left( \mathcal{X}_{S,+},f_{\emph{SDPCM}}(\mathbf{y}; \mathbf{x}), g(\mathbf{x})  \right)$
with sparsity degree $S$ and denote the value and index of the $S$-largest entry of $\mathbf{x}_{0} \in \mathcal{X}_{S,+}$ by $\xi_{0}$ and $j_{0}$, respectively. 
We furthermore define the sets $\mathcal{K}_{l} \triangleq \big\{ \{ l \} \cup \big\{ \supp(\mathbf{x}_{0}) \setminus \{ j_{0}  \} \big\} \big\}$ for $l \in [N] \setminus \supp(\mathbf{x}_{0})$ and 
$\mathcal{K}_{l}  \triangleq \supp(\mathbf{x}_{0})$ for $l \in \supp(\mathbf{x}_{0})$. 
If the parameter function $g(\cdot):\mathcal{X}_{S,+} \rightarrow \mathbb{R}$ associated with $\mathcal{E}_{\emph{SDPCM}}$ is such that 
the partial derivatives $b_{l} \triangleq \frac{\partial^{\mathbf{e}_{l}} g(\mathbf{x})}{\partial \mathbf{x}^{\mathbf{e}_{l}}} \big|_{\mathbf{x} = \mathbf{x}_{0}^{\mathcal{K}_{l}}}$ exist, we have
\begin{equation} 
\label{equ_lower_bound_SPCM_ICASSP_modification}
L_{\mathcal{M}_{\emph{SDPCM}}}  \geq 2 \!\! \sum_{l \in \supp(\mathbf{x}_{0})} \hspace*{-5mm} \frac{(x_{0,l} + \sigma^{2})^{2}}{r_{l}}b^{2}_{l} \,\, + \,\, 2
  \sigma^{4} \ist  \Bigg[ 1-  \frac{\xi_{0}^{2}}{(\xi_{0} + \sigma^{2})^{2}} \Bigg]^{\frac{r_{j_{0}}}{2}} \sum_{l \in [N] \setminus \supp(\mathbf{x}_{0})}  \frac{b^{2}_{l}}{r_{l}}.
\end{equation} 
\end{theorem}

\begin{proof} 
As for the derivation of Theorem \ref{thm_sparse_lower_bound_SDPCM_ICASSP}, we can assume that the prescribed mean function $\gamma(\cdot): \mathcal{D}_{0} \rightarrow \mathbb{R}: \gamma(\mathbf{x}) = g(\mathbf{x})$ 
belongs to the RKHS $\mathcal{H}(\mathcal{M}_{\mathcal{D}_{0},\text{SDPCM}})$. 
We then consider the subspace $\mathcal{U} \triangleq \linspan \big\{ \{ w_{0}(\cdot) \} \cup \{ g^{(\mathbf{e}_{l},\mathcal{K}_{l})} (\cdot) \}_{l \in [N]}Ê\big\}  \subseteq \mathcal{H}(\mathcal{M}_{\mathcal{D}_{0},\text{SDPCM}})$ 
which is spanned by the functions $w_{0}(\cdot) \triangleq R_{\mathcal{M}_{\mathcal{D}_{0},\text{SDPCM}}}(\cdot, \mathbf{x}_{0})$ and 
$g^{(\mathbf{e}_{l},\mathcal{K}_{l})} (\cdot) \in \mathcal{H}(\mathcal{M}_{\mathcal{D}_{0},\text{SDPCM}})$ as defined in Theorem \ref{thm_charac_RKHS_SDPCM_D_0}, i.e., 
$g^{(\mathbf{e}_{l},\mathcal{K}_{l})} (\cdot) = \frac{ \partial^{\mathbf{p}} R_{\mathcal{M}_{\mathcal{D}_{0}, \text{SDPCM}}}(\cdot, \mathbf{x}_{2})}{\partial \mathbf{x}_{2}^{\mathbf{e}_{l}}}\bigg|_{\mathbf{x}_{2} = \mathbf{x}_{0}^{\mathcal{K}_{l}}}$. 

The function $w_{0}(\cdot)$ is orthogonal to any function $g^{(\mathbf{e}_{l},\mathcal{K}_{l})} (\cdot)$ for $l \in [N]$, since 
\begin{align} 
\big\langle g^{(\mathbf{e}_{l},\mathcal{K}_{l})} (\cdot),w_{0}(\cdot) \big\rangle_{\mathcal{H}(\mathcal{M}_{\mathcal{D}_{0},\text{SDPCM}})} & =
\big\langle g^{(\mathbf{e}_{l},\mathcal{K}_{l})} (\cdot),R_{\mathcal{M}_{\mathcal{D}_{0},\text{SDPCM}}}(\cdot, \mathbf{x}_{0}) \big\rangle_{\mathcal{H}(\mathcal{M}_{\mathcal{D}_{0},\text{SDPCM}})} \nonumber \\[4mm]
&  \stackrel{\eqref{equ_reproduction_property}}{=}  
 \frac{ \partial^{\mathbf{e}_{l}} R_{\mathcal{M}_{\mathcal{D}_{0}, \text{SDPCM}}}(\mathbf{x}_{0}, \mathbf{x}_{2})}{\partial \mathbf{x}_{2}^{\mathbf{e}_{l}}}\bigg|_{\mathbf{x}_{2} = \mathbf{x}_{0}^{\mathcal{K}_{l}}} Ê\nonumber \\[4mm]
& \stackrel{\eqref{equ_kernel_one_arg_x_0_equal_1}}{=}  \frac{ \partial^{\mathbf{e}_{l}}1}  {\partial \mathbf{x}_{2}^{\mathbf{e}_{l}}} \bigg|_{\mathbf{x}_{2} = \mathbf{0}}  = 0. 
\end{align}
The inner products between the functions $\big\{ g^{(\mathbf{e}_{l},\mathcal{K}_{l})} (\cdot) \big\}_{l \in [N]}$ are calculated according to 
\eqref{equ_inner_prod_part_der_RKHS_Re_SDPCM} as 
\vspace*{2mm}
\begin{align} 
\label{equ_inner_prod_l_l_prime_x_0_SPCM_ICASSP_mod_lower_bound_starting}
\big\langle g^{(\mathbf{e}_{l},\mathcal{K}_{l})} (\cdot), g^{(\mathbf{e}_{l'},\mathcal{K}_{l'})} (\cdot) \big\rangle_{\mathcal{H}(\mathcal{M}_{\mathcal{D}_{0},\text{SDPCM}})} 
& \stackrel{\eqref{equ_inner_prod_part_der_RKHS_Re_SDPCM}}{=}  \frac{ \partial^{\mathbf{e}_{l}}} {\partial \mathbf{x}_{1}^{\mathbf{e}_{l}}} g^{(\mathbf{e}_{l'},\mathcal{K}_{l'})} (\mathbf{x}_{1})\bigg|_{\mathbf{x}_{1} = \mathbf{x}_{0}^{\mathcal{K}_{l}}}\nonumber \\[4mm]
& \hspace*{-40mm} \stackrel{\eqref{equ_def_g_p_thm_charac_RKHS_SDPCM_D_0}}{=} \frac{ \partial^{\mathbf{e}_{l}} } {\partial \mathbf{x}_{1}^{\mathbf{e}_{l}}}
 \frac{r_{l'} (x_{1,l'} - x_{0,l'})}{2 ( x_{0,l'}+\sigma^{2})^{2}} \prod_{k \in \supp(\mathbf{x}_{0}) 
 \setminus \mathcal{K}_{l'}} \frac{(x_{0,k} \rmv\rmv+\rmv \sigma^{2})^{r_{k}} }{\!\big[ (x_{0,k} \rmv\rmv+\rmv \sigma^{2})^2 + (x_{1,k} \!-\rmv\rmv x_{0,k})\rmv\rmv x_{0,k}\big]^{r_k/2} } 
 \bigg|_{\mathbf{x}_{1} = \mathbf{x}_{0}^{\mathcal{K}_{l}}},
\end{align} 
which is most conveniently evaluated by considering separately the case where $l,l' \in \supp(\mathbf{x}_{0})$ and the complementary case 
where either $l$ or $l'$ (or both) belong to $[N] \setminus \supp(\mathbf{x}_{0})$. 
Let us begin with the latter case, where due to the symmetry of the inner product, we can assume without loss of generality that $l \in [N] \setminus \supp(\mathbf{x}_{0})$ 
and $l' \in [N]$ is arbitrary. In this case, \eqref{equ_inner_prod_l_l_prime_x_0_SPCM_ICASSP_mod_lower_bound_starting} yields 
\begin{align} 
\label{equ_inner_prod_l_l_prime_x_0_SPCM_ICASSP_mod_lower_bound_latter_case}
&\big\langle g^{(\mathbf{e}_{l},\mathcal{K}_{l})} (\cdot), g^{(\mathbf{e}_{l'},\mathcal{K}_{l'})} (\cdot)  \big\rangle_{\mathcal{H}(\mathcal{M}_{\mathcal{D}_{0},\text{SDPCM}})}  \nonumber \\[4mm]
&= \frac{ \partial^{\mathbf{e}_{l}} } {\partial \mathbf{x}_{1}^{\mathbf{e}_{l}}}
 \frac{r_{l'} (x_{1,l'} - x_{0,l'})}{2 ( x_{0,l'}+\sigma^{2})^{2}} \prod_{k \in \supp(\mathbf{x}_{0}) 
 \setminus \mathcal{K}_{l'}} \frac{(x_{0,k} \rmv\rmv+\rmv \sigma^{2})^{r_{k}} }{\!\big[ (x_{0,k} \rmv\rmv+\rmv \sigma^{2})^2 + (x_{1,k} \!-\rmv\rmv x_{0,k})\rmv\rmv x_{0,k}\big]^{r_k/2} } 
 \bigg|_{\mathbf{x}_{1} = \mathbf{x}_{0}^{\mathcal{K}_{l}}}\nonumber \\[4mm]
& \stackrel{l \notin \supp(\mathbf{x}_{0})}{=} \Bigg[\prod_{k \in \supp(\mathbf{x}_{0}) 
 \setminus \mathcal{K}_{l'}} \frac{(x_{0,k} \rmv\rmv+\rmv \sigma^{2})^{r_{k}} }{\!\big[ (x_{0,k} \rmv\rmv+\rmv \sigma^{2})^2 + (x_{1,k} \!-\rmv\rmv x_{0,k})\rmv\rmv x_{0,k}\big]^{r_k/2} } 
 \bigg|_{\mathbf{x}_{1} = \mathbf{x}_{0}^{\mathcal{K}_{l}}} \Bigg] \frac{ \partial^{\mathbf{e}_{l}} } {\partial \mathbf{x}_{1}^{\mathbf{e}_{l}}}
 \frac{r_{l'} (x_{1,l'} - x_{0,l'})}{2 ( x_{0,l'}+\sigma^{2})^{2}} \bigg|_{\mathbf{x}_{1} = \mathbf{x}_{0}^{\mathcal{K}_{l}}} \nonumber \\[4mm]
&= \delta_{l,l'} \big\langle g^{(\mathbf{e}_{l},\mathcal{K}_{l})} (\cdot), g^{(\mathbf{e}_{l},\mathcal{K}_{l})} (\cdot)  \big\rangle_{\mathcal{H}(\mathcal{M}_{\mathcal{D}_{0},\text{SDPCM}})}.
\end{align}
In the former case, where $l,l' \in \supp(\mathbf{x}_{0})$ and therefore $\mathcal{K}_{l} = \mathcal{K}_{l'} = \supp(\mathbf{x}_{0})$, 
we obtain by \eqref{equ_inner_prod_l_l_prime_x_0_SPCM_ICASSP_mod_lower_bound_starting} that 
\begin{align} 
\label{equ_inner_prod_l_l_prime_x_0_SPCM_ICASSP_mod_lower_bound_former_case}
&\big\langle g^{(\mathbf{e}_{l},\mathcal{K}_{l})} (\cdot), g^{(\mathbf{e}_{l'},\mathcal{K}_{l'})} (\cdot)  \big\rangle_{\mathcal{H}(\mathcal{M}_{\mathcal{D}_{0},\text{SDPCM}})}  \nonumber \\[4mm]
&= \frac{ \partial^{\mathbf{e}_{l}} } {\partial \mathbf{x}_{1}^{\mathbf{e}_{l}}}
 \frac{r_{l'} (x_{1,l'} - x_{0,l'})}{2 ( x_{0,l'}+\sigma^{2})^{2}} \prod_{k \in \supp(\mathbf{x}_{0}) 
 \setminus \mathcal{K}_{l'}} \frac{(x_{0,k} \rmv\rmv+\rmv \sigma^{2})^{r_{k}} }{\!\big[ (x_{0,k} \rmv\rmv+\rmv \sigma^{2})^2 + (x_{1,k} \!-\rmv\rmv x_{0,k})\rmv\rmv x_{0,k}\big]^{r_k/2} } 
 \bigg|_{\mathbf{x}_{1} = \mathbf{x}_{0}^{\mathcal{K}_{l}}}\nonumber \\[4mm]
& \stackrel{\supp(\mathbf{x}_{0}) 
 \setminus \mathcal{K}_{l'}=\emptyset}{=}  \frac{ \partial^{\mathbf{e}_{l}} } {\partial \mathbf{x}_{1}^{\mathbf{e}_{l}}}
 \frac{r_{l'} (x_{1,l'} - x_{0,l'})}{2 ( x_{0,l'}+\sigma^{2})^{2}} \bigg|_{\mathbf{x}_{1} = \mathbf{x}_{0}^{\mathcal{K}_{l}}} \nonumber \\[4mm]
&= \delta_{l,l'} \big\langle g^{(\mathbf{e}_{l},\mathcal{K}_{l})} (\cdot), g^{(\mathbf{e}_{l},\mathcal{K}_{l})} (\cdot)  \big\rangle_{\mathcal{H}(\mathcal{M}_{\mathcal{D}_{0},\text{SDPCM}})}.
\end{align}
By combining \eqref{equ_inner_prod_l_l_prime_x_0_SPCM_ICASSP_mod_lower_bound_latter_case} and \eqref{equ_inner_prod_l_l_prime_x_0_SPCM_ICASSP_mod_lower_bound_former_case}, 
we obtain 
\begin{equation} 
\big\langle g^{(\mathbf{e}_{l},\mathcal{K}_{l})} (\cdot), g^{(\mathbf{e}_{l'},\mathcal{K}_{l'})} (\cdot) \big\rangle_{\mathcal{H}(\mathcal{M}_{\mathcal{D}_{0},\text{SDPCM}})} 
= \delta_{l,l'} q_{l}(\mathbf{x}_{0}),
\end{equation} 
with 
\begin{align} 
\label{equ_sqared_norm_q_x_0_SPCM_ICASSP_mod_lower_bound}
q_{l}(\mathbf{x}_{0}) & \triangleq   
\big\langle g^{(\mathbf{e}_{l},\mathcal{K}_{l})} (\cdot), g^{(\mathbf{e}_{l},\mathcal{K}_{l})} (\cdot)  \big\rangle_{\mathcal{H}(\mathcal{M}_{\mathcal{D}_{0},\text{SDPCM}})} \nonumber \\[4mm]
& =\frac{ \partial^{\mathbf{e}_{l}} } {\partial \mathbf{x}_{1}^{\mathbf{e}_{l}}}
 \frac{r_{l} (x_{1,l} - x_{0,l})}{2 ( x_{0,l}+\sigma^{2})^{2}} \prod_{k \in \supp(\mathbf{x}_{0}) 
 \setminus \mathcal{K}_{l}} \frac{(x_{0,k} \rmv\rmv+\rmv \sigma^{2})^{r_{k}} }{\!\big[ (x_{0,k} \rmv\rmv+\rmv \sigma^{2})^2 + (x_{1,k} \!-\rmv\rmv x_{0,k})\rmv\rmv x_{0,k}\big]^{r_k/2} } 
 \bigg|_{\mathbf{x}_{1} = \mathbf{x}_{0}^{\mathcal{K}_{l}}}  \nonumber \\[4mm] 
 & =\Bigg[  \prod_{k \in \supp(\mathbf{x}_{0}) 
 \setminus \mathcal{K}_{l}} \frac{(x_{0,k} \rmv\rmv+\rmv \sigma^{2})^{r_{k}} }{\!\big[ (x_{0,k} \rmv\rmv+\rmv \sigma^{2})^2 -\rmv\rmv x_{0,k}^{2}\big]^{r_k/2} } 
 \Bigg]\frac{ \partial^{\mathbf{e}_{l}} } {\partial \mathbf{x}_{1}^{\mathbf{e}_{l}}}
 \frac{r_{l} (x_{1,l} - x_{0,l})}{2 ( x_{0,l}+\sigma^{2})^{2}} \bigg|_{\mathbf{x}_{1} = \mathbf{x}_{0}^{\mathcal{K}_{l}}} \nonumber \\[4mm] 
& =  \begin{cases}  \frac{r_{l}}{2 (x_{0,l}+\sigma^{2})^{2}}, & \mbox{if } l \in \supp(\mathbf{x}_{0}) \\Ê
                            \frac{r_{l}}{2 \sigma^{4}} \bigg[  \frac{(\xi_{0} \rmv\rmv+\rmv \sigma^{2})^{2} }{ (\xi_{0} \rmv\rmv+\rmv \sigma^{2})^2 -\rmv\rmv \xi_{0}^{2} } 
 \bigg]^{r_{j_{0}}/2} & \mbox{if } l \in [N] \setminus \supp(\mathbf{x}_{0}). 
 \end{cases}
\end{align}  

Using the inner products (recall that $\gamma(\cdot) = g(\cdot)$)
\begin{equation} 
\label{equ_inner_prod_part_der_gamm_SPCM_ICASSP_mod_lower_bound}
\big\langle \gamma(\cdot),g^{(\mathbf{e}_{l},\mathcal{K}_{l})} (\cdot) \big\rangle_{ \mathcal{H}(\mathcal{M}_{\mathcal{D}_{0},\text{SDPCM}})} 
\stackrel{\eqref{equ_inner_prod_part_der_RKHS_Re_SDPCM}}{=} 
\frac{\partial^{\mathbf{e}_{l}} g(\mathbf{x})}{\partial \mathbf{x}^{\mathbf{e}_{l}}} \bigg|_{\mathbf{x} = \mathbf{x}_{0}^{\mathcal{K}_{l}}}
= b_{l},
\end{equation}  
the bound \eqref{equ_lower_bound_SPCM_ICASSP_modification} follows then 
by projecting $\gamma(\cdot)$ onto the subspace $\mathcal{U}$: 
\begin{align}
 L_{\mathcal{M}_{\text{SDPCM}}} & \stackrel{\eqref{equ_lower_bound_spcm_restr_spcm}}{\geq} 
 L_{\mathcal{M}_{\mathcal{D}_{0},\text{SDPCM}}} \stackrel{\eqref{equ_squared_norm_min_achiev_var}}{=} \| \gamma(\cdot) \|^{2}_{\mathcal{H}(\mathcal{M}_{\mathcal{D}_{0},\text{SDPCM}})}
 - \big[\underbrace{\gamma(\mathbf{x}_{0})}_{= g(\mathbf{x}_{0})}\big]^{2}  \nonumber \\[4mm] 
&  \stackrel{\eqref{equ_squared_norm_orthog_proj_pythag_thm}}{\geq}  \| \mathbf{P}_{\mathcal{U}}\gamma(\cdot)  \|^{2}_{\mathcal{H}(\mathcal{M}_{\mathcal{D}_{0},\text{SDPCM}})} - \big[g(\mathbf{x}_{0})\big]^{2}  \nonumber \\[4mm]
&  \stackrel{\eqref{equ_idendity_projection_subspace_spanned_union_orthogonal_sets}}{=}  
\frac{ \big\langle \gamma(\cdot), w_{0}(\cdot) \big\rangle_{\mathcal{H}(\mathcal{M}_{\mathcal{D}_{0},\text{SPCM}})}^{2}}{ \big\langle w_{0}(\cdot), w_{0}(\cdot) 
\big\rangle_{\mathcal{H}(\mathcal{M}_{\mathcal{D}_{0},\text{SPCM}})}} + \sum_{l \in [N]}  
 \frac{ \big\langle \gamma(\cdot), g^{(\mathbf{e}_{l},\mathcal{K}_{l})}(\cdot) \big\rangle_{\mathcal{H}(\mathcal{M}_{\mathcal{D}_{0},\text{SPCM}})}^{2}} 
 {\big\langle  g^{(\mathbf{e}_{l},\mathcal{K}_{l})}(\cdot), g^{(\mathbf{e}_{l},\mathcal{K}_{l})}(\cdot) \big\rangle_{\mathcal{H}(\mathcal{M}_{\mathcal{D}_{0},\text{SPCM}})}}  - \big[g(\mathbf{x}_{0})\big]^{2}  \nonumber \\[4mm]
& = \frac{\big\langle \gamma(\cdot), R_{\mathcal{M}_{\mathcal{D}_{0},\text{SDPCM}}}(\cdot, \mathbf{x}_{0})(\cdot) \big\rangle_{\mathcal{H}(\mathcal{M}_{\mathcal{D}_{0},\text{SPCM}})}^{2}}{ \big\langle R_{\mathcal{M}_{\mathcal{D}_{0},\text{SDPCM}}}(\cdot, \mathbf{x}_{0}), R_{\mathcal{M}_{\mathcal{D}_{0},\text{SDPCM}}}(\cdot, \mathbf{x}_{0}) 
\big\rangle_{\mathcal{H}(\mathcal{M}_{\mathcal{D}_{0},\text{SPCM}})} }  \nonumber \\[4mm] 
&  \hspace*{10mm}+ \sum_{l \in [N]}   \frac{\big\langle \gamma(\cdot), g^{(\mathbf{e}_{l},\mathcal{K}_{l})}(\cdot) \big\rangle_{\mathcal{H}(\mathcal{M}_{\mathcal{D}_{0},\text{SPCM}})}^{2}}{\big\langle  g^{(\mathbf{e}_{l},\mathcal{K}_{l})}(\cdot), g^{(\mathbf{e}_{l},\mathcal{K}_{l})}(\cdot) \big\rangle_{\mathcal{H}(\mathcal{M}_{\mathcal{D}_{0},\text{SPCM}})}}  - \big[g(\mathbf{x}_{0})\big]^{2}  \nonumber \\[4mm]
& \stackrel{\eqref{equ_reproduction_property},\eqref{equ_sqared_norm_q_x_0_SPCM_ICASSP_mod_lower_bound}}{=} \frac{\big[ \gamma(\mathbf{x}_{0}) \big]^{2}}{R_{\mathcal{M}_{\mathcal{D}_{0},\text{SDPCM}}}(\mathbf{x}_{0}, \mathbf{x}_{0})} + \sum_{l \in [N]}   \frac{\big\langle \gamma(\cdot), g^{(\mathbf{e}_{l},\mathcal{K}_{l})}(\cdot) \big\rangle_{\mathcal{H}(\mathcal{M}_{\mathcal{D}_{0},\text{SPCM}})}^{2}}{q_{l}(\mathbf{x}_{0})}  - \big[g(\mathbf{x}_{0})\big]^{2}  \nonumber \\[4mm]
& \stackrel{\eqref{equ_kernel_one_arg_x_0_equal_1}}{=} \frac{\big[ \gamma(\mathbf{x}_{0}) \big]^{2}}{1} + \sum_{l \in [N]}   \frac{\big\langle \gamma(\cdot), g^{(\mathbf{e}_{l},\mathcal{K}_{l})}(\cdot) \big\rangle_{\mathcal{H}(\mathcal{M}_{\mathcal{D}_{0},\text{SPCM}})}^{2}}{q_{l}(\mathbf{x}_{0})}  - \big[g(\mathbf{x}_{0})\big]^{2}  \nonumber \\[4mm]
& \stackrel{\gamma(\cdot)=g(\cdot)}{=}  \sum_{l \in [N]}   \frac{\big\langle \gamma(\cdot), g^{(\mathbf{e}_{l},\mathcal{K}_{l})}(\cdot) \big\rangle_{\mathcal{H}(\mathcal{M}_{\mathcal{D}_{0},\text{SPCM}})}^{2}}{q_{l}(\mathbf{x}_{0})}   \nonumber \\[4mm]
& \stackrel{\eqref{equ_sqared_norm_q_x_0_SPCM_ICASSP_mod_lower_bound},\eqref{equ_inner_prod_part_der_gamm_SPCM_ICASSP_mod_lower_bound}}{=} 2 \!\! \sum_{l \in \supp(\mathbf{x}_{0})} \hspace*{-4mm} \frac{(x_{0,l} + \sigma^{2})^{2}}{r_{l}}b^{2}_{l} + 2
\sigma^{4} \ist 
 \bigg[ 1 - \frac{  \xi_{0}^{2} } {(\xi_{0} \rmv\rmv+\rmv \sigma^{2})^{2} }
 \bigg]^{r_{j_{0}}/2}  \sum_{l \in [N] \setminus \supp(\mathbf{x}_{0})}  \frac{b^{2}_{l}}{r_{l}}.
\end{align} 
\end{proof} 
In addition to being easily evaluated, the bound in \eqref{equ_lower_bound_SPCM_ICASSP_modification} has two appealing properties. First, it is continuous with respect to $\mathbf{x}_{0}$, and 
second, it is always at least as tight as the bound given by Theorem \ref{thm_CRB_SPCM} when specialized to the SDPCM. Indeed, the bounds \eqref{equ_CRB_SPCM_not_full_sparsity} and 
\eqref{equ_CRB_SPCM_full_sparsity} of Theorem \ref{thm_CRB_SPCM} applied to the SDPCM yield $L_{\mathcal{M}_{\text{SDPCM}}} \geq \sum_{l \in [N]} 2 (x_{0,l} + \sigma^{2})^{2} \frac{b_{l}^{2}}{r_{l}}$ 
for $\| \mathbf{x}_{0} \|_{0} < S$ and $L_{\mathcal{M}_{\text{SDPCM}}} \geq \sum_{l \in \supp(\mathbf{x}_{0}) } 2 (x_{0,l} + \sigma^{2})^{2} \frac{b_{l}^{2}}{r_{l}}$ when $\| \mathbf{x}_{0} \|_{0} = S$. 
Comparing these expressions with \eqref{equ_lower_bound_SPCM_ICASSP_modification} (note that $\xi_{0} = 0$ when $\|\mathbf{x}_{0} \|_{0} < S$), it can be shown 
that Theorem \ref{thm_sparse_lower_bound_SDPCM_ICASSP_modifciaton} yields always a lower bound that is at least as tight, i.e., at least as high as 
the bound given by Theorem \ref{thm_CRB_SPCM} applied to the SDPCM.  

Similar to Corollary \ref{cor_thm_sparse_lower_bound_SDPCM_ICASSP}, we have the following corollary 
of Theorem \ref{thm_sparse_lower_bound_SDPCM_ICASSP_modifciaton}, which considers estimators for the SDPCM with an arbitrary bias function. 
\begin{corollary}
\label{cor_thm_sparse_lower_bound_SDPCM_ICASSP_modifciaton}
Consider an estimator $\hat{x}_{k}(\mathbf{y})$ which uses the observation $\mathbf{y}$ of a 
SDPCM with sparsity degree $S$ and denote its mean function by $m_{k}(\cdot): \mathcal{X}_{S,+} \rightarrow \mathbb{R}: m_{k} (\mathbf{x}) \triangleq \mathsf{E}_{\mathbf{x}} \big\{ \hat{x}_{k}(\mathbf{y}) \big\}$. 
Let us fix a parameter vector $\mathbf{x}_{0} \in \mathcal{X}_{S,+}$ and denote the value and index of the $S$ largest entry of $\mathbf{x}_{0}$ by $\xi_{0}$ and $j_{0}$, respectively. 
Furthermore, we define the sets $\mathcal{K}_{l} \triangleq \big\{ \{ l \} \cup \big\{ \supp(\mathbf{x}_{0}) \setminus \{ j_{0}  \} \big\} \big\}$ for $l \in [N] \setminus \supp(\mathbf{x}_{0})$ and 
$\mathcal{K}_{l}  \triangleq \supp(\mathbf{x}_{0})$ for $l \in \supp(\mathbf{x}_{0})$. Then, if $m_{k}(\cdot)$ is such that 
the partial derivatives $b_{l} \triangleq \frac{\partial^{\mathbf{e}_{l}} m_{k} (\mathbf{x})}{\partial \mathbf{x}^{\mathbf{e}_{l}}} \big|_{\mathbf{x} = \mathbf{x}_{0}^{\mathcal{K}_{l}}}$ exist, we have
\begin{align} 
\label{equ_cor_thm_sparse_lower_bound_SDPCM_ICASSP_modifciaton}
v(\hat{x}_{k}(\cdot); \mathbf{x}_{0}) \geq   2 \!\! \sum_{l \in \supp(\mathbf{x}_{0})} \hspace*{-5mm} \frac{(x_{0,l} + \sigma^{2})^{2}}{r_{l}}b^{2}_{l} \,\, + \,\, 2
  \sigma^{4} \ist  \Bigg[ 1-  \frac{\xi_{0}^{2}}{(\xi_{0} + \sigma^{2})^{2}} \Bigg]^{\frac{r_{j_{0}}}{2}} \sum_{l \in [N] \setminus \supp(\mathbf{x}_{0})}  \frac{b^{2}_{l}}{r_{l}}. Ê\\[-4mm]
\nonumber
\end{align} 
\end{corollary}

If we want to use Corollary \ref{cor_thm_sparse_lower_bound_SDPCM_ICASSP} or Corollary \ref{cor_thm_sparse_lower_bound_SDPCM_ICASSP_modifciaton} for the comparison of the actual variance behavior of a given estimation scheme $\hat{x}_{k}(\cdot)$
with the lower bounds \eqref{equ_cor_thm_sparse_lower_bound_SDPCM_ICASSP} and \eqref{equ_cor_thm_sparse_lower_bound_SDPCM_ICASSP_modifciaton}, respectively, we 
have to ensure that the partial derivatives of the mean function $m_{k}(\mathbf{x}) = \mathsf{E}_{\mathbf{x}} \{ \hat{x}_{k}(\mathbf{y}) \}$ of the given estimator $\hat{x}_{k}(\cdot)$ exist. 
That this is indeed the case for a very broad class of estimators $\hat{x}_{k}(\cdot)$ is stated in 

\begin{lemma}
\label{lem_cond_exist_partial_der_exist_spec_est_SDPCM}
Consider an estimator $\hat{x}_{k}(\mathbf{y})$ which uses the observation $\mathbf{y}$ of a SDPCM, and whose mean is $m_{k}(\mathbf{x})$.
If $\hat{x}_{k}(\cdot): \mathbb{R}^{M} \rightarrow \mathbb{R}$ is a Lebesgue-measurable function, and moreover  
\begin{equation}
\label{equ_bound_suff_cond_exist_par_der_exist_SDPCM_spec_est}
|\hat{x}_{k}(\mathbf{y})| \leq C \| \mathbf{y} \|_{2}^{L} \quad \quad \forall \mathbf{y} \in \mathbb{R}^{M},
\end{equation}
with arbitrary but fixed constants $C,L \in \mathbb{R}$, then 
the partial derivatives $\frac{ \partial^{\mathbf{p}} m_{k}(\mathbf{x})}{\partial \mathbf{x}^{\mathbf{p}}}$ exist for every multi-index $\mathbf{p} \in \mathbb{Z}_{+}^{N} \in [N]$. 
For the case $\mathbf{p} = \mathbf{e}_{l}$ with $l \in [N]$, the partial derivatives are given explicitly by
\begin{align}
\label{equ_expr_par_der_existing_SDPCM_spec_est}
\frac{ \partial^{\mathbf{e}_{l}} m_{k}(\mathbf{x})}{\partial \mathbf{x}^{\mathbf{e}_{l}}}= 
-\frac{r_{l}}{2(x_{l}+\sigma^{2})} m_{k}(\mathbf{x}) + \frac{1}{2(x_{l}+\sigma^{2})^{2}}  \mathsf{E}_{\mathbf{x}} \big\{ \mathbf{y}^{T} \mathbf{C}_{l} \mathbf{y}\hat{x}_{k}(\mathbf{y}) \big\}. \\[-5mm] 
\nonumber
\end{align} 
\end{lemma} 
\begin{proof}
Appendix \ref{chap_appendix_E}.
\end{proof}

\subsection{Special Case: Unbiased Spectrum Estimation}
\label{sec_unbiasecd_SDPCM}

Let us now consider the special case of the SDPCM obtained for the parameter function $g(\mathbf{x}) = x_{k}$, where $k \in [N]$ is an arbitrary but fixed index. 
Since the entries $x_{k}$ of the parameter vector $\mathbf{x}$ can be interpreted as mean powers of the signal vector $\mathbf{s}$ within the subspace $\linspan(\mathbf{C}_{k}) \subseteq \mathbb{R}^{M}$ 
(cf.\ \eqref{equ_cov_param_SPCM}), the resulting estimation problem is that of unbiased nonstationary spectrum estimation \cite{TimeFrequencyAnalysisBoashash,stoi97} 
if the basis matrices $\mathbf{C}_{k}$ correspond to well-localized regions in a time-frequency domain. 
The specialization of Theorem \ref{thm_sparse_lower_bound_SDPCM_ICASSP_modifciaton} to the parameter function $g(\mathbf{x}) = x_{k}$ yields  
\begin{corollary} 
\label{corr_SDCM_unbiased}
Consider the minimum variance problem  $\mathcal{M}_{\emph{SDPCM}}=\left(\mathcal{E}_{\emph{SDPCM}},c(\cdot) \equiv 0 , \mathbf{x}_{0}Ê\right)$ associated with the 
SDPCM $\mathcal{E}_{\emph{SDPCM}} = \left( \mathcal{X}_{S,+}, f_{\emph{SDPCM}}(\mathbf{y}; \mathbf{x}), g(\mathbf{x}) = x_{k} \right)$. 
For any estimator $\hat{x}_{k}(\cdot): \mathbb{R}^{M} \rightarrow \mathbb{R}$ that is unbiased, i.e., $\mathsf{E}_{\mathbf{x}} \{ \hat{x}_{k}(\cdot) \} =x_k$ for every $\mathbf{x} \in \mathcal{X}_{S,+}$, 
we have the following lower bound on its variance at $\mathbf{x}_{0}$: 
\begin{align} 
&v(\hat{x}_{k}(\cdot); \mathbf{x}_{0}) \geq \begin{cases}   
 \displaystyle \frac{2} {r_{k}}(x_{0,k} \rmv+\rmv \sigma^{2})^2 , & k \rmv\in\rmv \supp(\mathbf{x}_{0}) \\[2.5mm]Ê
 \displaystyle \frac{2} {r_{k}} \sigma^{4} \Bigg[ 1-  \frac{\xi_{0}^{2}}{(\xi_{0} + \sigma^{2})^{2}} \Bigg]^{\frac{r_{j_{0}}}{2}} \,, 
   &k \rmv\not\in\rmv \supp(\mathbf{x}_{0}) \ist, \end{cases}  \label{equ_corr_lower_bound}\\[-2mm]
& \nonumber \\[-7mm]
& \nonumber
\end{align}
where $\xi_{0}$ and $j_{0}$ denote the value and index, respectively, of the $S$-largest entry of $\mathbf{x}_{0}$. 
\end{corollary} 
\begin{proof}
This Corollary follows straightforwardly from Theorem \ref{thm_sparse_lower_bound_SDPCM_ICASSP_modifciaton} and the relation $b_{l} = \delta_{l,k}$ where $b_{l}$ as defined in Theorem  \ref{thm_sparse_lower_bound_SDPCM_ICASSP_modifciaton}.
\end{proof}

The lower bound \eqref{equ_corr_lower_bound} 
can be achieved at least for the case $| \{ k \} \cup \supp(\mathbf{x}_{0})| < S+1$. 
In this case, 
the estimator given by (cf.\ \eqref{equ_basis_matrix_sum_rank_1_projection_SDPCM})
\begin{equation} 
\hat{x}_{k}(\mathbf{y}) = \beta_k(\mathbf{y}) \rmv-\rmv \sigma^{2} , \quad \text{with} \;\; 
  \beta_k(\mathbf{y}) \triangleq \frac{1}{r_{k}} \rmv\sum_{iÊ\in [r_k]} \rmv\big( \mathbf{u}^{T}_{m_{k,i}} \ist \mathbf{y} \big)^{2} , 
\label{equ_tight_est_SDPCM}
\end{equation} 
is unbiased and its variance achieves the bound 
\eqref{equ_corr_lower_bound}. Indeed, 
we have for the mean  
\begin{align}
\label{equ_mean_tight_est_SDPCM}
\mathsf{E}_{\mathbf{x}} \big\{  \hat{x}_{k}(\mathbf{y})  \big \} & = \mathsf{E}_{\mathbf{x}} \big\{  \beta_k(\mathbf{y})- \sigma^{2} \big\}  
 = \mathsf{E}_{\mathbf{x}} \bigg\{  \frac{1}{r_{k}} \rmv\sum_{iÊ\in [r_k]} \rmv\big( \mathbf{u}^{T}_{m_{k,i}} \ist \mathbf{y} \big)^{2}  \bigg\} - \sigma^{2}  \nonumber \\[4mm]
 & =  \frac{1}{r_{k}} \rmv\sum_{iÊ\in [r_k]} \mathsf{E}_{\mathbf{x}} \big\{  \rmv\big( \mathbf{u}^{T}_{m_{k,i}} \ist \mathbf{y} \big)^{2}  \big\} - \sigma^{2}
  =  \frac{1}{r_{k}} \rmv\sum_{iÊ\in [r_k]} (x_{k} + \sigma^{2}) - \sigma^{2} \nonumber \\[4mm]  
 & =  \frac{1}{r_{k}} r_{k} (x_{k} + \sigma^{2}) - \sigma^{2} = x_{k},
\end{align} 
where we used the fact 
\begin{equation}
 \mathbf{u}^{T}_{m_{k,i}} \ist \mathbf{y}  \sim \mathcal{N}(0,\mathbf{u}^{T}_{m_{k,i}} \widetilde{\mathbf{C}}(\mathbf{x}) \mathbf{u}_{m_{k,i}})= \mathcal{N}(0, x_{k} + \sigma^{2}), 
\end{equation} 
which follows from $\mathbf{y} \sim \mathcal{N}(\mathbf{0}, \widetilde{\mathbf{C}}(\mathbf{x}))$, \eqref{equ_def_cov_matrix_observation_SPCM} and \eqref{equ_basis_matrix_sum_rank_1_projection_SDPCM}. 
Using the shorthand $z_{i} \triangleq  \mathbf{u}^{T}_{m_{k,i}} \ist \mathbf{y}$, the variance at an arbitrary parameter vector $\mathbf{x}_{0} \in \mathcal{X}_{S,+}$ can be calculated as  
\vspace*{1mm}
\begin{align} 
\label{equ_spdcm_variance_tight_est}
v( \hat{x}_{k}(\cdot); \mathbf{x}_{0}) & = v( \beta_{k}(\cdot); \mathbf{x}_{0})  \nonumber \\[4mm]
& = P (\beta_{k}(\cdot); \mathbf{x}_{0}) - \big[ \mathsf{E}_{\mathbf{x}_{0}} \{  \beta_{k}(\mathbf{y}) \} \big]^{2}  \nonumber \\[4mm]
& = \mathsf{E}_{\mathbf{x}_{0}} \bigg\{ \bigg[  \frac{1}{r_{k}} \sum_{iÊ\in [r_k]} z_{i}^{2}  \bigg]^{2} \bigg\} - \bigg[ \frac{1}{r_{k}}\sum_{iÊ\in [r_k]}  \mathsf{E}_{\mathbf{x}_{0}} \{ z^{2}_{i} \} \bigg]^{2} \nonumber \\[4mm] 
 & = \frac{1}{r^{2}_{k}} \sum_{i,i'Ê\in [r_k]}  \mathsf{E}_{\mathbf{x}_{0}} \big\{  z_{i}^{2} z_{i'}^{2} \big\}     - \frac{1}{r_{k}^{2}}\bigg[  \sum_{iÊ\in [r_k]}  \mathsf{E}_{\mathbf{x}_{0}} \{ z^{2}_{i} \} \bigg]^{2}.
 \end{align} 
 Using the fact that the random variables $\{ z_{i} \}_{i \in [r_k]}$ are i.i.d.\ with $z_{i} \sim \mathcal{N}(0,x_{k} + \sigma^{2})$, we obtain further
 \begin{align}
 \label{equ_var_tight_est_SDPCM}
 v( \hat{x}_{k}(\cdot); \mathbf{x}_{0}) 
& = \frac{1}{r^{2}_{k}} \Bigg[ \sum_{iÊ\in [r_k]}  \mathsf{E}_{\mathbf{x}_{0}} \big\{  z_{i}^{4} \big\}  +
 \sum_{\substack{i,i'Ê\in [r_k]\\ i \neq i'}}  \mathsf{E}_{\mathbf{x}_{0}} \big\{  z_{i}^{2} \big\} \mathsf{E}_{\mathbf{x}_{0}} \big\{  z_{i'}^{2} \big\}  
  - \bigg[ \sum_{iÊ\in [r_k]}  \mathsf{E}_{\mathbf{x}_{0}} \{ z^{2}_{i} \} \bigg]^{2} \Bigg]
  \nonumber \\[4mm] 
  &  = \frac{1}{r^{2}_{k}} \Bigg[ \sum_{iÊ\in [r_k]}  \mathsf{E}_{\mathbf{x}_{0}} \big\{  z_{i}^{4} \big\}  +
 \bigg[ \sum_{iÊ\in [r_k]}  \mathsf{E}_{\mathbf{x}_{0}} \{ z^{2}_{i} \} \bigg]^{2}
 -  \sum_{iÊ\in [r_k]} \big[ \mathsf{E}_{\mathbf{x}_{0}} \big\{  z_{i}^{2} \big\} \big]^{2}  
  - \bigg[\sum_{iÊ\in [r_k]}  \mathsf{E}_{\mathbf{x}_{0}} \{ z^{2}_{i} \} \bigg]^{2}  \Bigg] \nonumber \\[4mm] 
 & = \frac{1}{r^{2}_{k}}  \Bigg[  \sum_{iÊ\in [r_k]}  \mathsf{E}_{\mathbf{x}_{0}} \big\{  z_{i}^{4} \big\} 
 - \sum_{iÊ\in [r_k]} \big[ \mathsf{E}_{\mathbf{x}_{0}} \big\{  z_{i}^{2} \big\} \big]^{2} \Bigg]  \nonumber \\[4mm] 
  & = \frac{1}{r^{2}_{k}} \Bigg[ \sum_{iÊ\in [r_k]}  3 (x_{k} + \sigma^{2})^{2} -   \sum_{iÊ\in [r_k]} (x_{k} + \sigma^{2})^{2}\Bigg] \nonumber \\[4mm]
 & = \frac{2}{r_{k}} (x_{k} + \sigma^{2})^{2}, 
\end{align}
where we used a well-known identity for the fourth-order central moment of a zero-mean Gaussian random variable \cite{papoulis}. 
The estimator in \eqref{equ_tight_est_SDPCM} does not 
use the sparsity information and does not depend on 
$\mathbf{x}_{0}$. Moreover, since this estimator is the unique LMVU estimator at any $\| \mathbf{x}_{0} \|_{0} < S$, we have that 
if a UMVU estimator existed for the unbiased spectrum estimation problem $\mathcal{E}_{\text{SDPCM}} =  \left(\mathcal{X}_{S,+}, f_{\text{SDPCM}}(\mathbf{y}; \mathbf{x}), g(\mathbf{x}) = x_{k}\right)$, 
it would necessarily coincide with the nonsparse estimator in \eqref{equ_tight_est_SDPCM} and, thus, it would 
ignore the sparsity information $\mathbf{x} \in \mathcal{X}_{S,+}$. 
This is an instance of the general rule that unbiased estimation is usually not optimal for 
estimation problems with sparsity constraints. 

Let us define a ``signal-to-noise ratio'' (SNR) quantity as $\mbox{SNR}(\mathbf{x}_{0}) \triangleq \xi_{0}/\sigma^{2}\rmv\rmv$.
For $\mbox{SNR}(\mathbf{x}_{0}) \!\ll\! 1$, where $\big[ 1-  \frac{\xi_{0}^{2}}{(\xi_{0} + \sigma^{2})^{2}} \big]^{\frac{r_{j_{0}}}{2}}Ê\approx 1$, 
the lower bound \eqref{equ_corr_lower_bound} in Corollary \ref{corr_SDCM_unbiased} is approximately
$ \frac{2} {r_{k}}(x_{0,k} \rmv+\rmv \sigma^{2})^2$ for any $k$, which does not depend on
$S$ and moreover equals the variance of the 
unbiased estimator \eqref{equ_tight_est_SDPCM}. 
Since that estimator does not exploit any sparsity information, Corollary \ref{corr_SDCM_unbiased} suggests that, in the low-SNR regime, 
unbiased estimators cannot exploit the prior information that $\mathbf{x}$ is $S$-sparse, i.e., $\mathbf{x} \in \mathcal{X}_{S,+}$. 
However, in the high-SNR regime ($\mbox{SNR}(\mathbf{x}_{0}) \!\rightarrow\! \infty$), \eqref{equ_corr_lower_bound} becomes
$\frac{2} {r_{k}}(x_{0,k} \rmv+\rmv \sigma^{2})^2$ for $k \rmv\in\rmv \supp(\mathbf{x}_{0})$ and $0$ for $k \rmv\not\in\rmv \supp(\mathbf{x}_{0})$, which can be shown to equal
the variance of the oracle estimator that knows 
$\supp(\mathbf{x}_{0})$ (this oracle estimator yields $\hat{x}_{k} \!=\! x_{0,k} \!=\! 0$ 
for all $k \rmv\not\in\rmv \supp(\mathbf{x}_{0})$). 
The transition of the lower bound \eqref{equ_corr_lower_bound} from the low-SNR regime to the high-SNR regime 
has a polynomial characteristic; it is thus much slower than the exponential transition of the analogous 
lower bound for the SLM given by Theorem \ref{thm_bound_asilomar}.
This slow
transition suggests that the optimal unbiased estimator for low SNR---which ignores the sparsity information---
will also be a nearly optimal unbiased estimator over a relatively wide SNR range. 
This further suggests that, for unbiased covariance estimation based on the SDPCM, 
prior information of sparsity is not as helpful as for minimum variance estimation for the SLM.

Based on the specific estimator given by \eqref{equ_tight_est_SDPCM}, we can also show that for the SDPCM, 
similarily to the SLM (see Section \ref{sec_strict_sparstiy_SLM}), the strict sparsity requirement $\mathbf{x} \in \mathcal{X}_{S,+}$ is 
necessary in order to allow for the existence of unbiased spectrum estimators $\hat{x}_{k}(\cdot)$ whose variance at $\mathbf{x}_{0}$ is smaller than 
the variance of the nonsparse estimator in \eqref{equ_tight_est_SDPCM}. We will also show that the variance of the estimator in \eqref{equ_tight_est_SDPCM} is simultaneously 
the theoretical minimum variance that can be achieved without any sparsity constraints. 

Indeed, consider any minimum variance problem $\mathcal{M}' = \left(\mathcal{E}',c(\cdot)\equiv 0, \mathbf{x}_{0}Ê\right)$ associated with 
the estimation problem $\mathcal{E}' = \left(\mathcal{X}', f_{\text{SPCM}}(\mathbf{y}; \mathbf{x}), g(\mathbf{x}) = x_{k} \right)$ (with arbitrary $k \in [N]$), where the statistical model 
$f_{\text{SPCM}}(\mathbf{y}; \mathbf{x})$ (cf.\ \eqref{equ_def_SPCM_stat_model}) is defined using basis matrices $\mathbf{C}_{k}$ satisfying \eqref{equ_conds_basis_matrices_SDPCM}. 
Thus, the minimum variance problem $\mathcal{M}'$ is identical to the minimum variance problem $\mathcal{M}_{\text{SDPCM}}=\left(\mathcal{E}_{\text{SDPCM}},c(\cdot) \equiv 0 , \mathbf{x}_{0}Ê\right)$ 
with $\mathcal{E}_{\text{SDPCM}} =  \left(\mathcal{X}_{S,+}, f_{\text{SDPCM}}(\mathbf{y}; \mathbf{x}), g(\mathbf{x}) = x_{k}\right)$ 
except for the parameter set which is $\mathcal{X}'$ instead of $\mathcal{X}_{S,+}$. Let us assume that the parameter set $\mathcal{X}'$ 
contains an open ball $\mathcal{B}(\mathbf{x}_{0},r)$, i.e., $\mathcal{B}(\mathbf{x}_{0},r) \subseteq \mathcal{X}'$ with some radius $r > 0$. We can then apply Theorem \ref{thm_unconstr_CR}, i.e., 
the ordinary unconstrained CRB for a general nonsparse estimation problem with Gaussian observation $\mathbf{y} \sim \mathcal{N}(\mathbf{0}, \mathbf{C}(\mathbf{x}))$ \cite{kay}, 
to obtain the lower bound $v(\hat{x}_{k}(\cdot); \mathbf{x}_{0}) \geq \frac{2} {r_{k}}(x_{0,k} \rmv+\rmv \sigma^{2})^2$, which is achieved by the estimator in \eqref{equ_tight_est_SDPCM}. 
From this, we conclude that if the parameter set is chosen in particular as $\mathcal{X}' = \mathbb{R}^{N}_{+} \cap \mathcal{B}_{q}(S)$ with $q>0$ (cf.\ \eqref{equ_def_approximate_S_sparse_vecs_q_ball}), i.e., 
consisting of approximately $S$-sparse vectors, then the UMVU estimator is always given by the estimator in \eqref{equ_tight_est_SDPCM}, whose variance 
is in general strictly larger than the bound \eqref{equ_corr_lower_bound} for the case $k \notin \supp(\mathbf{x}_{0})$. In particular for the case $\| \mathbf{x}_{0} \|_{0} =S$, i.e., $\xi_{0} > 0$ we have 
\begin{equation} 
\frac{2} {r_{k}} \sigma^{4} \underbrace{\big[ 1-  \frac{\xi_{0}^{2}}{(\xi_{0} + \sigma^{2})^{2}} \big]^{\frac{r_{j_{0}}}{2}}}_{ < 1} < \frac{2} {r_{k}}(x_{0,k} \rmv+\rmv \sigma^{2})^2.
\end{equation} 

However, for the SDPCM with the strictly sparse parameter set $\mathcal{X}_{S,+}$, there exist in the general case unbiased 
estimators $\hat{x}_{k}(\cdot)$ for $\mathcal{M}_{\text{SDPCM}}=\left(\mathcal{E}_{\text{SDPCM}},c(\cdot) \equiv 0 , \mathbf{x}_{0}Ê\right)$ 
with $\mathcal{E}_{\text{SDPCM}} =  \left(\mathcal{X}_{S,+}, f_{\text{SDPCM}}(\mathbf{y}; \mathbf{x}), g(\mathbf{x}) = x_{k}\right)$, 
that have a variance that is strictly smaller than that of the estimator \eqref{equ_tight_est_SDPCM} at 
a certain parameter vector $\mathbf{x}_{0} \in \mathcal{X}_{S,+}$. 
In particular, consider the simplest configuration of the SDPCM where $\mathbf{C}_{k} = \sum_{i \in [r_{k}]} \mathbf{e}_{j_{k,i}} \mathbf{e}_{j_{k,i}}^{T}$, $N > 1$, and $S=1$. 
For this instance of the SDPCM, one can find the LMVU estimator $\hat{x}^{(\mathbf{x}_{0})}_{k}(\cdot)$ for the minimum variance problem 
$\mathcal{M}_{\text{SDPCM}}=\left(\mathcal{E}_{\text{SDPCM}},c(\cdot) \equiv 0 , \mathbf{x}_{0}Ê\right)$ (with parameter function $g(\mathbf{x}) = x_{k}$): 
\begin{theorem} 
\label{thm_LMVU_SDPCM_S_1}
Consider the minimum variance problem $\mathcal{M}_{\emph{SDPCM}}=\left(\mathcal{E}_{\emph{SDPCM}},c(\cdot) \equiv 0 , \mathbf{x}_{0}Ê\right)$, 
with $\mathbf{x}_{0} \in \mathcal{X}_{S,+}$, $N>1$, $S=1$, and $g(\mathbf{x}) = x_{k}$. 
Denoting by $\xi_{0}$ and $j_{0}$ the value and index, respectively, of the only nonzero entry of $\mathbf{x}_{0}$, the LMVU estimator for $\mathcal{M}_{\emph{SDPCM}}$ can be written as 
\begin{equation}
\label{equ_LMVU_SDPCM_S_1}
\hat{x}^{(\mathbf{x}_{0})}_{k}(\mathbf{y}) = \begin{cases} 
  \beta_k(\mathbf{y}) \rmv-\rmv \sigma^{2} , & k = j_{0} \\[.5mm]
  \alpha(\mathbf{y};\mathbf{x}_{0}) \ist \big( \beta_k(\mathbf{y}) \rmv-\rmv \sigma^{2} \big) \,, & k \neq j_{0} \,,  \end{cases}
\end{equation}
where 
\begin{equation} 
\alpha(\mathbf{y};\mathbf{x}_{0}) \triangleq a(\mathbf{x}_{0}) \exp\rmv\rmv\big(\!\rmv-\rmv\rmv r_{j_{0}} b(\mathbf{x}_{0}) \ist \beta_{j_{0}}(\mathbf{y})\big)
\end{equation} 
with 
$a(\mathbf{x}_{0}) \triangleq \Big[ \frac{ 2\xi_{0} + \sigma^{2} }{\xi_{0} + \sigma^{2}} \Big]^{r_{j_{0}}/2}$
, $b(\mathbf{x}_{0}) \triangleq \frac{1}{2}\ist\big(\frac{1}{\sigma^{2}} - \frac{1}{\xi_{0} + \sigma^{2}} \big)$, and $\beta_{k}(\mathbf{y})$ as defined in \eqref{equ_tight_est_SDPCM}. 
\end{theorem}
\begin{proof}
Appendix \ref{chap_appendix_B}. 
\end{proof}


\section{Comparison of the Bounds with Existing Estimators} 
\label{sec_comp_bounds_actual_var_SPCM}

We will now compare the lower bound \eqref{equ_cor_thm_sparse_lower_bound_SDPCM_ICASSP} (in Corollary \ref{cor_thm_sparse_lower_bound_SDPCM_ICASSP}) 
on the variance of estimators for the SDPCM, with the actual variance behavior of two practical estimation schemes.
The first is an ad-hoc adaptation of the hard-thresholding (HT) estimator \cite{Donoho94idealspatial} to SDPCM-based covariance estimation.
It is defined componentwise as (cf.\ \eqref{equ_tight_est_SDPCM})
\[
\hat{x}_{\text{HT},k}(\mathbf{y}) \ist\triangleq\ist \frac{1}{r_k} \!\sum_{i \in [r_k]} \!\varphi_T^2 \big( \mathbf{u}^{T}_{m_{k,i}} \ist \mathbf{y} \big)  - \sigma^{2} ,
\vspace{-.7mm}
\]
where $ \varphi_T \!:\rmv \mathbb{R} \!\to\! \mathbb{R}$ denotes the hard-thresholding function with threshold $T \rmv\ge\rmv 0$, 
i.e., 
\begin{equation}
\varphi_T(y)= \begin{cases} y \,,  & |y| \geq T\\[.5mm]
  0 \,, & \text{else} \ist.
\end{cases}
\end{equation} 
The second standard method 
is the maximum likelihood (ML) 
\vspace*{-.5mm}
estimator  
\[ 
\hat{\mathbf{x}}_{\text{ML}}(\mathbf{y}) \,\triangleq\, \argmax_{\mathbf{x}' \in \mathcal{X}_{S,+}} \ist  f_{\text{SDPCM}}( \mathbf{y};  \mathbf{x}' ) \,.
\vspace{-1mm}
\] 

For the ML estimator, we obtain  
\begin{align}
\hat{\mathbf{x}}_{\text{ML}}(\mathbf{y}) & \triangleq \argmax_{\mathbf{x}' \in \mathcal{X}_{S,+}}  f_{\text{SDPCM}}(\mathbf{y}; \mathbf{x}')  \nonumber \\[4mm]
& \stackrel{\eqref{equ_def_SDPCM_stat_model}}{=} \argmax_{\mathbf{x}' \in \mathcal{X}_{S,+}} \frac{1}{(2 \pi)^{M/2} [\detm{\widetilde{\mathbf{C}}(\mathbf{x}')}]^{1/2}}\exp \left(- \frac{1}{2} \mathbf{y}^{T}  \widetilde{\mathbf{C}}^{-1}(\mathbf{x}') \mathbf{y} \right) \nonumber \\[4mm]
& \stackrel{(a)}{=}  \argmax_{\mathbf{x}' \in \mathcal{X}_{S,+}} \log \Bigg[\frac{1}{[\detm{\widetilde{\mathbf{C}}(\mathbf{x}')}]^{1/2}}\exp \left(- \frac{1}{2} \mathbf{y}^{T}  \widetilde{\mathbf{C}}^{-1}(\mathbf{x}') \mathbf{y} \right)\Bigg] \nonumber \\[4mm] 
& = \argmax_{\mathbf{x}' \in \mathcal{X}_{S,+}}  \bigg\{ - \frac{1}{2} \mathbf{y}^{T}  \widetilde{\mathbf{C}}^{-1}(\mathbf{x}') \mathbf{y} - \frac{1}{2}\log\big[\detm{\widetilde{\mathbf{C}}(\mathbf{x}')}\big] \bigg\}
\nonumber \\[4mm] 
& = \argmax_{\mathbf{x}' \in \mathcal{X}_{S,+}}\big\{ -  \mathbf{y}^{T}  \widetilde{\mathbf{C}}^{-1}(\mathbf{x}') \mathbf{y} - \log\big[\detm{\widetilde{\mathbf{C}}(\mathbf{x}')}\big] \big\},
\end{align}
where the step $(a)$ follows from the fact that the logarithm is monotonically increasing on its domain $\mathbb{R}_{+} \setminus \{0\}$, i.e., it preserves the position of the maximum. 
Using \eqref{equ_expr_obs_cov_matrix_SDPCM} and \eqref{equ_matrix_D_prime_SDPCM}, we obtain further 
\begin{align} 
\label{equ_SDPCM_ML_est_expression_sum_h_k}
\hat{\mathbf{x}}_{\text{ML}}(\mathbf{y})  = \argmax_{\mathbf{x}' \in \mathcal{X}_{S,+}} \Bigg\{ - \sum_{k \in [N]}Ê r_{k} \bigg[ \frac{\beta_{k}(\mathbf{y})}{x'_{k}+\sigma^{2}} +  \log\big( x'_{k}+\sigma^{2} \big) \bigg]   \Bigg\}, 
\end{align}
with $\beta_{k}(\mathbf{y})$ as defined in \eqref{equ_tight_est_SDPCM}. 

We have 
\begin{theorem}
\label{thm_ml_est_SDPCM}
The ML estimator for the SDPCM is given componentwise as 
\begin{align} 
\label{equ_def_ML_est_SDPCM_expr}
\vspace*{-.5mm}
\hat{x}_{\emph{ML},k}(\mathbf{y}) = \begin{cases} 
  \beta_k(\mathbf{y}) \rmv-\rmv \sigma^{2} , & k \rmv\in\rmv \mathcal{L}_1 \!\cap\rmv \mathcal{L}_2 \\ 
  0 \,, &\mbox{else} \ist, \end{cases} \\[-5mm]
\nonumber
\end{align}
where $\beta_k(\mathbf{y})$ is defined in \eqref{equ_tight_est_SDPCM}, $\mathcal{L}_1$ consists of the $S$ indices $k \in [N]$ for which 
$r_{k} \big[ \beta_{k}(\mathbf{y})/\sigma^2 \! -Ê\! \log \ist (\beta_{k}(\mathbf{y})/\sigma^{2}) \rmv-\! 1 \big]$
is largest, and 
$\mathcal{L}_2$ consists of all indices $k$ for which $\beta_k(\mathbf{y}) \rmv\geq\rmv \sigma^{2}\rmv$.
\end{theorem}
\begin{proof}
Appendix \ref{chap_appendix_C}. 
\end{proof}

In what follows, we will denote by $L^{(m_{k}(\cdot))}(\mathbf{x}_{0})$ the bound \eqref{equ_cor_thm_sparse_lower_bound_SDPCM_ICASSP} in 
Corollary \ref{cor_thm_sparse_lower_bound_SDPCM_ICASSP} for an estimator $\hat{x}_{k}(\cdot)$ with mean function 
$\mathsf{E}_{\mathbf{x}} \big\{ \hat{x}_{k}(\mathbf{y}) \big\}=m_{k}(\cdot)$, using the index set $\mathcal{K} = \supp(\mathbf{x}_{0})$ when $k \in \supp(\mathbf{x}_{0})$ and 
$\mathcal{K} = \big\{ \supp(\mathbf{x}_{0})\setminus \{j_{0}\}Ê\big\} \cup \{k \}$ else, where $j_{0} \in [N]$ denotes the index of the $S$-largest entry of $\mathbf{x}_{0}$. 
We furthermore use the choices $L \!=\! 2$, $\mathbf{p}_{1} \!=\! \mathbf{0}$, and $\mathbf{p}_{2} \!=\rmv \mathbf{e}_k$ for the application of Corollary \ref{cor_thm_sparse_lower_bound_SDPCM_ICASSP}. 
In order to evaluate the bound $L^{(m_{k}(\cdot))}(\mathbf{x}_{0})$, we need to compute the first-order partial derivatives of the mean function $m_{k}(\cdot)$.
This is accomplished by using Lemma \ref{lem_cond_exist_partial_der_exist_spec_est_SDPCM} for the HT estimator and by using a finite-difference quotient approximation \cite{HeroUniformCRB} for the ML estimator, i.e., 
\begin{equation}
\frac{\partial \mathsf{E}_{\mathbf{x}} \big\{\hat{x}_{\text{ML},k}(\mathbf{y})\big\} }{\partial x_{l}} \approx \frac{ \mathsf{E}_{\mathbf{x}+ \Delta \mathbf{e}_{l}} \big\{\hat{x}_{\text{ML},k}(\mathbf{y})\big\} -\mathsf{E}_{\mathbf{x}} \big\{\hat{x}_{\text{ML},k}(\mathbf{y})\big\} }{\Delta}, 
\end{equation}  
where $\Delta \in \mathbb{R}_{+}$ is a small stepsize and the expectations are calculated using numerical integration.

Given an estimator $\hat{\mathbf{x}}(\cdot)$ whose mean is equal to $\mathbf{m}(\mathbf{x}) \triangleq \mathsf{E}_{\mathbf{x}}\big\{ \hat{\mathbf{x}}(\cdot) \big\}$, we obtain due to \eqref{equ_sum_min_var_scalar_min_var} a lower bound on the estimator variance $v(\hat{\mathbf{x}}(\cdot);\mathbf{x}_{0})$ by summing the quantities 
$L^{(m_{k}(\cdot))}(\mathbf{x}_{0})$, i.e., 
\begin{equation}
\label{equ_lower_bound_var_vec_est_numerical_SDPCM}
v(\hat{\mathbf{x}}(\cdot);\mathbf{x}_{0}) \geq L^{(\mathbf{m}(\cdot))}(\mathbf{x}_{0}) \triangleq \sum_{k \in [N]} \rmv\rmv L^{(m_{k}(\cdot))}(\mathbf{x}_{0}).
\end{equation}

\begin{figure}
\vspace{-1.5mm}
\centering
\psfrag{SNR}[c][c][.9]{\uput{3.4mm}[270]{0}{\hspace{-.1mm}$\mbox{SNR}$ [dB]}}
\psfrag{title}[c][c][.9]{\uput{2.5mm}[270]{0}{}}
\psfrag{x_0_01}[c][c][.9]{\uput{0.3mm}[270]{0}{$-20$}}
\psfrag{x_0_1}[c][c][.9]{\uput{0.3mm}[270]{0}{$-10$}}
\psfrag{x_1}[c][c][.9]{\uput{0.3mm}[270]{0}{$0$}}
\psfrag{x_10}[c][c][.9]{\uput{0.3mm}[270]{0}{$10$}}
\psfrag{x_100}[c][c][.9]{\uput{0.3mm}[270]{0}{$20$}}
\psfrag{x_1000}[c][c][.9]{\uput{0.3mm}[270]{0}{$30$}}
\psfrag{y_0}[c][c][.9]{\uput{0.1mm}[180]{0}{$0$}}
\psfrag{y_5}[c][c][.9]{\uput{0.1mm}[180]{0}{$5$}}
\psfrag{y_10}[c][c][.9]{\uput{0.1mm}[180]{0}{$10$}}
\psfrag{y_15}[c][c][.9]{\uput{0.1mm}[180]{0}{$15$}}
\psfrag{y_20}[c][c][.9]{\uput{0.1mm}[180]{0}{$20$}}
\psfrag{y_25}[c][c][.9]{\uput{0.1mm}[180]{0}{$25$}}
\psfrag{y_30}[c][c][.9]{\uput{0.1mm}[180]{0}{$30$}}
\psfrag{y_35}[c][c][.9]{\uput{0.1mm}[180]{0}{$35$}}
\psfrag{y_40}[c][c][.9]{\uput{0.1mm}[180]{0}{$40$}}
\psfrag{y_45}[c][c][.9]{\uput{0.1mm}[180]{0}{$45$}}
\psfrag{y_50}[c][c][.9]{\uput{0.1mm}[180]{0}{$50$}}
\psfrag{y_2}[c][c][.9]{\uput{0.1mm}[180]{0}{$2$}}
\psfrag{y_3}[c][c][.9]{\uput{0.1mm}[180]{0}{$3$}}
\psfrag{variance}[c][c][.9]{\uput{3mm}[90]{0}{\hspace*{-3mm}normalized variance/bound}}
\psfrag{data2}[l][l][0.8]{bound on $v(\hat{\mathbf{x}}_{\text{ML}}(\cdot);\mathbf{x}_0)$}
\psfrag{data1}[l][l][0.8]{$v(\hat{\mathbf{x}}_{\text{ML}}(\cdot);\mathbf{x}_0)$}
\psfrag{data3}[l][l][0.8]{$v(\hat{\mathbf{x}}_{\text{HT}}(\cdot); \mathbf{x}_0)$, $T \!=\! 3$}
\psfrag{data5}[l][l][0.8]{$v(\hat{\mathbf{x}}_{\text{HT}}(\cdot); \mathbf{x}_0)$, $T \!=\! 6$}
\psfrag{data7}[l][l][0.8]{$v(\hat{\mathbf{x}}_{\text{HT}}(\cdot); \mathbf{x}_0)$, $T \!=\! 9$}
\psfrag{data4}[l][l][0.8]{bound on $v(\hat{\mathbf{x}}_{\text{HT}}(\cdot); \mathbf{x}_0)$, $T \!=\! 3$}
\psfrag{data6}[l][l][0.8]{bound on $v(\hat{\mathbf{x}}_{\text{HT}}(\cdot); \mathbf{x}_0)$, $T \!=\! 6$}
\psfrag{data8}[l][l][0.8]{bound on $v(\hat{\mathbf{x}}_{\text{HT}}(\cdot); \mathbf{x}_0)$, $T \!=\! 9$}
\psfrag{ML}[l][l][0.8]{ML}
\psfrag{T=3}[l][l][0.8]{\uput{0mm}[180]{0}{HT\,\ist($T \!=\! 3$)}}
\psfrag{T=6}[l][l][0.8]{HT\,\ist($T \!=\! 6$)}
\psfrag{T=9}[l][l][0.8]{HT\,\ist($T \!=\! 9$)}
\centering
\hspace*{1.5mm}\includegraphics[height=8cm,width=15cm]{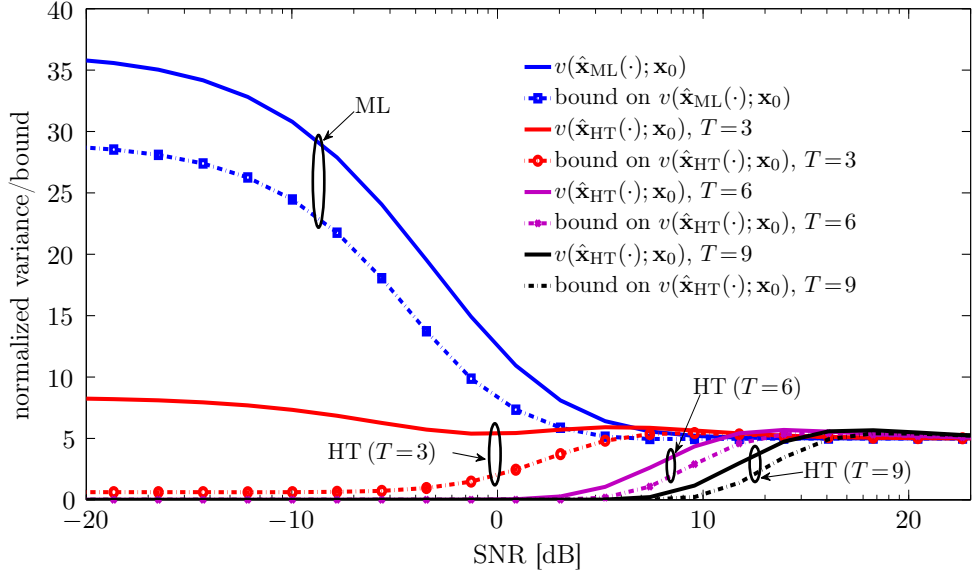}
\vspace{1mm}
  \caption{Normalized variance of the ML and HT estimators and corresponding lower bounds 
  versus the $\mbox{SNR}$, for the SDPCM 
  with $N\!=\!50$, $S\!=\!5$, $\sigma^{2} \!=\!1$, and $\mathbf{C}_{k} \!=\rmv \mathbf{e}_{k} \mathbf{e}_{k}^{T}$. } 
\label{fig_bounds_SDPCM}
\vspace*{3.5mm}
\end{figure}

For a numerical evaluation, we considered the SDPCM with $N\!\rmv=\!50$, $S\!=\!5$, and $\mathbf{C}_{k} \!=\rmv \mathbf{e}_{k}\mathbf{e}_{k}^{T}$.
We generated parameter vectors $\mathbf{x}_{0} = \text{SNR} \cdot \sigma^{2} \mathbf{x}_{1}$, where $\mathbf{x}_{1} \in \{0,1\}^{50}$, 
$\supp(\mathbf{x}_{1}) = [S]$, and $\text{SNR}$ varies between $0.01$ and $100$. (The fixed choice $\supp(\mathbf{x}_{0}) = [S]$ is justified by the fact that neither the variances of the ML and HT estimators nor the corresponding variance bounds 
depend on the location of $\supp(\mathbf{x}_{0})$.)
In Fig.\ \ref{fig_bounds_SDPCM}, we show the 
variance at $\mathbf{x}_{0}$, $v(\hat{\mathbf{x}}(\cdot); \mathbf{x}_0) \!=\! \sum_{k \in [N]} v(\hat{x}_{k}(\cdot);\mathbf{x}_0)$ 
(computed by means of numerical integration), for the HT estimator using various choices of 
$T$ and for the ML estimator. 
The variance is plotted versus $\mbox{SNR}$. 
Along with each variance curve, we display the corresponding lower 
bound $L^{(\mathbf{m}(\cdot))}(\mathbf{x}_{0})$, using for 
$\mathbf{m}(\cdot)$ the mean function of the respective estimator (HT or ML). 
(The mean functions of the HT and ML estimators were computed by means of numerical integration.)
In Fig.\ \ref{fig_bounds_SDPCM}, all variances and bounds are normalized by $2  \ist (\mbox{SNR}\rmv+\rmv1)^{2} \sigma^{4}$, which is 
the variance of an optimum unbiased oracle estimator for the SDPCM with $S=1$ and the nonzero entry of $\mathbf{x}_{0}$ being equal to $\mbox{SNR} \cdot \sigma^{2}$. 
The computation of the oracle estimator presupposes the knowledge of the position of the nonzero entry of $\mathbf{x}_{0}$, i.e., it knows $\supp(\mathbf{x}_{0})$. 

It can be seen from Fig.\ \ref{fig_bounds_SDPCM} that in the high-SNR regime, for both estimators, the gap between the variance and the corresponding lower bound 
is quite small. This indicates that the variances of both estimators are nearly minimal (for the respective bias functions). Indeed, for high SNR, all variance curves approach the variance of an optimum unbiased oracle estimator that knows the position of the nonzero entries of $\mathbf{x}_{0}$, i.e., it knows $\supp(\mathbf{x}_{0})$. The variance 
of this oracle estimator is given by $2  S \ist (\mbox{SNR}+1)^{2} \sigma^{4}$ (this expression is only valid for parameter vectors $\mathbf{x}_{0}$ with exactly $S$ nonzero entries whose values are all equal to $\mbox{SNR} \cdot \sigma^{2}$).
However, in the low-SNR regime, 
the variances of the 
estimators tend to be
significantly higher than the 
bounds. This means that there \emph{may} be estimators 
with the same bias function as that of the HT or ML estimator but a lower variance. However, the actual existence of such estimators 
is not shown by
our analysis.

\chapter{Conclusion and Outlook}

This thesis specialized the RKHS approach to minimum variance estimation \cite{Parzen59,Duttweiler73b} to estimation problems with sparsity constraints on 
the unknown parameter vector to be estimated. 
After a brief introduction to the main concepts of classical (non-Bayesian) estimation and Hilbert space theory - with emphasis on RKHS - we reviewed the RKHS approach to minimum variance estimation. This approach is based on two isometric Hilbert spaces, one of them being a RKHS, which can be associated in a natural manner to a minimum variance problem. 
We presented two fundamental facts about minimum variance estimation using the RHKS approach, which seem to be novel, to the best of the author's knowledge. 
The first fact is that the minimum achievable variance, viewed as a function of the parameter vector $\mathbf{x}_{0}$ at which the variance is minimized, 
is a lower semi-continuous function. The second fact is that the RKHS associated to a minimum variance problem remains unchanged if the observation of the 
minimum variance problem is replaced by a sufficient statistic. This fact is closely related in spirit to the Rao-Blackwell-Lehmann-Scheff{\'e} theorem of classical estimation. 

The core of this thesis are Chapter \ref{chap_SLM} and Chapter \ref{chap_SCM}, in which we considered classical estimation problems with sparsity constraints. 
In Chapter \ref{chap_SLM}, we considered the SLM, where the unknown sparse parameter vector is linearly distorted and corrupted by additive white Gaussian noise. 
The SLM can be used 
to model the recovery problem in a linear compressed sensing scheme.
We analyzed the minimum variance estimation problem using the RKHS approach and derived novel lower bounds on the estimator variance for a prescribed estimator bias. These lower bounds are obtained by an orthogonal projection of the prescribed mean function onto a subspace of the RKHS that is associated to the SLM. The bounds share the characteristic that they vary between two extreme cases. On the one hand, we have a low-SNR regime where the entries of the true parameter vector are small compared with the noise variance. Here, our bounds predict that the a-priori information of a sparse parameter vector does not help much for the SLM if the estimator bias is approximately zero. However, by allowing a nonzero bias it 
is possible to reduce the variance. On the other hand, we have a high-SNR regime where the nonzero entries of the sparse parameter vector are all 
considerable larger than the noise variance. Here, our bounds coincide with the CRB of an associated LGM that is obtained from the SLM if the 
support of the unknown parameter vector is known. Our bounds ehibit a very steep transition between these two regimes. In general, this transition has an exponential decay. 
For the special case of the SLM given by the SSNM, we strengthened the results and derived closed-form 
expressions for the minimum achievable variance and corresponding LMV estimator. This includes the (unbiased) Barankin bound and LMVU 
estimator for the SSNM. We also compared our lower bounds on the variance as well as the expressions for the minimum achievable variance to the actual variance behavior of popular estimation schemes for the SLM and SSNM. 

In Chapter \ref{chap_SCM}, we considered the SPCM, where the unknown sparse parameter vector determines the covariance matrix 
of a random signal vector of which a noisy version is observed. The SPCM is relevant, e.g., in the context of spectrum estimation for scene analysis in cognitive radio systems.
We discussed the fact that the RKHS approach is more involved for the SPCM than for the SLM.  
We presented one particular method of associating a RKHS to the SPCM. Based on this RKHS, we derived lower bounds on the estimator variance for estimators with 
a given bias function. As for the SLM, these bounds were obtained by projecting the corresponding mean function onto a subspace of the RKHS. 
For the special case of the SPCM which is given by the SDPCM, we strengthened our results and derived the minimum achievable variance and corresponding LMVU estimator 
for the case where the parameter vector is known to have at most one nonzero entry. As for the SLM and SSNM, also the lower bounds for the SPCM and SDPCM exhibit 
a transition from a low-SNR regime to a high-SNR regime, and the bounds for the high-SNR regime again coincide with the CRB of an associated estimation problem without sparsity constraints. 
However, a remarkable difference between the bounds for the 
SDPCM and for the SSNM is that the transition of the SDPCM bounds is only polynomial in the SNR whereas that of the SSNM bounds is exponential.
This further suggests that in terms of minimum variance estimation, for covariance estimation based on the SDPCM, prior information of sparsity is not as helpful as for estimating the parameter vector of the SSNM. 

Topics for future research that may be addressed using this thesis as a starting point include the following.
\begin{itemize}
\item {\bf Asymptotic Analysis.} Can we find (classes of) estimators that come close to the lower variance bounds asymptotically, 
i.e., when we consider not only one observation but an ensemble of i.i.d.\ observations whose number tends to infinity? 
Intuitively, the class of ML estimators should achieve the variance bounds asymptotically in the unbiased case. However, a rigorous proof of this claim seems to 
be nontrivial for the sparse parameter set $\mathcal{X}_{S}$. Indeed, most approaches to the analysis of the asymptotic behavior of ML estimators assume 
that the parameter set is an open subset of $\mathbb{R}^{N}$ \cite{LC,IbragimovBook,AsympVanderVaartBook}, which is not the case for the parameter set $\mathcal{X}_{S}$. 
Another popular class of estimators is given by the M-estimators or penalized maximum likelihood estimators, for which a characterization of their asymptotic behavior is available (cf.\ \cite{EldarUniformCRB,AsympVanderVaartBook,HuberRobustBook}). Moreover, under mild conditions, the class of M-estimators allows an efficient implementation via convex optimization techniques.

\item{\bf Finite-Rate-of-Innovation Signals.}
Another desirable generalization of our results would be to consider estimation of finite-rate-of-innovation signals, which are in some sense the continuous-time analogue 
of sparse vectors in the finite-dimensional setting \cite{ZvikaFRICRB}. 

\item{\bf Block Sparsity.} 
A popular generalization of the concept of strict sparsity that we used in our thesis is block or group sparsity \cite{BlockSparsityEldarTSP,EiwenSPAWC2010}, 
where a block structure is pre-specified for the parameter vector $\mathbf{x} \in \mathbb{R}^{N}$. 
Block sparsity then refers to the fact that only few blocks contain at least one nonzero entry. The traditional notion of strict sparsity is reobtained when the block size is equal to one. 
Block sparsity could be useful, e.g., in the context of spectrum estimation for sparse underspread processes \cite{jung-ssp09,jung-icassp09} if the process energy occurs in clusters in the time-frequency plane. 

\item{\bf Sparse Exponential Family.}
The SLM and SPCM are estimation problems with a statistical model that is a special case of the exponential family. 
It would be interesting to study estimation problems with sparsity constraints whose statistical model is the general exponential family. 
In particular, it would be interesting to study how the cumulant function of the exponential family relates to the RKHS and its geometric properties. 

\item{\bf Random System Matrix and Basis Matrices.} 

Another class of estimation problems that may be interesting to consider is obtained by a modification of the SLM and SPCM in which the system matrix $\mathbf{H}$
and basis matrices $\mathbf{C}_{k}$, respectively, are modeled as random matrices. In particular, 
a characterization of the RKHS for a SLM with a random system matrix $\mathbf{H}$, corresponding to the compressed sensing measurement model, would be desirable. 
Similarly, one could investigate the RKHS associated to a SPCM with random basis matrices $\mathbf{C}_{k}$. The resulting 
RKHS could be used to analyze the compressive spectrum estimators proposed in \cite{jung-ssp09,jung-icassp09}. 
\end{itemize}


\begin{appendices}

\chapter[]{Proof of Theorem \ref{thm_diag_bias_min_achiev_var_LMV}}
\label{chap_appendix_A}

\begin{proof} 
We give two different (though similar) proofs for the case $| \supp(\mathbf{x}_{0}) \cup \{k\}| < S+1$ and the complementary case $| \supp(\mathbf{x}_{0}) \cup \{k\}| = S+1$. 

\paragraph{Case I: $| \supp(\mathbf{x}_{0}) \cup \{k\}| < S+1$:}$~~$\\
Consider the functions 
\begin{equation} 
\label{equ_def_part_der_proof_SSNM_main_result_case_1}
h^{(l)}(\cdot): \mathcal{X}_{S} \rightarrow \mathbb{R}: h^{(l)}(\mathbf{x}) = \frac{\sigma^{l}}{\sqrt{l!}} \frac{\partial^{l \mathbf{e}_{k}} 
R_{\mathcal{M}_{\text{SSNM}}}(\mathbf{x}, \mathbf{x}_{2})}{\partial \mathbf{x}_{2}^{l \mathbf{e}_{k}}} \bigg|_{\mathbf{x}_{2} = \mathbf{x}_{0}} 
\stackrel{\eqref{equ_def_kernel_SSNM}}{=} \frac{1}{\sigma^{l}\sqrt{l!}} (x_{1,k} - x_{0,k})^{l}, 
\end{equation}
with $l \in \mathbb{Z}_{+}$. 
By Theorem \ref{thm_der_repr_prop} (note that $\mathcal{N}_{\mathbf{x}_{0}}^{\supp(l \mathbf{e}_{k})} \subseteq \mathcal{X}_{S}$) we have that the functions $\{ h^{(l)}(\cdot) \}_{l \in \mathbb{Z}_{+}}$ belong to $\mathcal{H}(\mathcal{M}_{\text{SSNM}})$ and moreover are orthonormal, since 
\begin{align} 
\big\langle h^{(l)}(\cdot), h^{(l')}(\cdot) \big\rangle_{\mathcal{H}(\mathcal{M}_{\text{SSNM}})}
& \stackrel{\eqref{equ_def_part_der_proof_SSNM_main_result_case_1}}{=} \bigg\langle \frac{\sigma^{l}}{\sqrt{l!}} \frac{\partial^{l \mathbf{e}_{k}} 
R_{\mathcal{M}_{\text{SSNM}}}(\mathbf{x}, \mathbf{x}_{2})}{\partial \mathbf{x}_{2}^{l \mathbf{e}_{k}}} \bigg|_{\mathbf{x}_{2} = \mathbf{x}_{0}},  \frac{\sigma^{l'}}{\sqrt{l'!}} \frac{\partial^{l' \mathbf{e}_{k}} 
R_{\mathcal{M}_{\text{SSNM}}}(\mathbf{x}, \mathbf{x}_{2})}{\partial \mathbf{x}_{2}^{l' \mathbf{e}_{k}}} \bigg|_{\mathbf{x}_{2} = \mathbf{x}_{0}} \bigg\rangle_{\mathcal{H}(\mathcal{M}_{\text{SSNM}})} \nonumber \\[4mm]
& \stackrel{\eqref{equ_der_reproduction_prop}}{=}    \frac{\sigma^{l}}{\sqrt{l!}} \frac{\sigma^{l'}}{\sqrt{l'!}} \frac{\partial^{l \mathbf{e}_{k}} \partial^{l' \mathbf{e}_{k}} 
R_{\mathcal{M}_{\text{SSNM}}}(\mathbf{x}_{1}, \mathbf{x}_{2})}{ \partial \mathbf{x}_{1}^{l \mathbf{e}_{k}}\partial \mathbf{x}_{2}^{l' \mathbf{e}_{k}}} \bigg|_{\mathbf{x}_{1} = \mathbf{x}_{2} = \mathbf{x}_{0}} \nonumber \\[4mm]
& \stackrel{\eqref{equ_def_kernel_SSNM}}{=}\frac{\sigma^{l}}{\sqrt{l!}} \frac{\sigma^{l'}}{\sqrt{l'!}}  \frac{\partial^{l \mathbf{e}_{k}} \partial^{l' \mathbf{e}_{k}} \exp \left( \frac{1}{\sigma^{2}} (\mathbf{x}_{1} - \mathbf{x}_{0})^{T} (\mathbf{x}_{2} - \mathbf{x}_{0}) \right)}{ \partial \mathbf{x}_{1}^{l \mathbf{e}_{k}}\partial \mathbf{x}_{2}^{l' \mathbf{e}_{k}}} \bigg|_{\mathbf{x}_{1} = \mathbf{x}_{2} = \mathbf{x}_{0}} \nonumber \\[4mm]
& = \frac{\sigma^{l}}{\sqrt{l!}} \frac{\sigma^{l'}}{\sqrt{l'!}}  \frac{1}{\sigma^{2l'}} \frac{\partial^{l \mathbf{e}_{k}} (x_{1,k} - x_{0,k})^{l'}}{ \partial \mathbf{x}_{1}^{l \mathbf{e}_{k}}} \bigg|_{\mathbf{x}_{1}  = \mathbf{x}_{0}}
= \delta_{l,l'}.
\end{align}
Obviously, then the function 
\begin{equation}
s_{N}(\cdot): \mathcal{X}_{S} \rightarrow \mathbb{R}: s_{N}(\mathbf{x}) \triangleq \sum_{l \in [N]} \frac{m_{l}\sigma^{l} }{\sqrt{l!}} h^{(l)}(\mathbf{x}) + m_{0}h^{(0)}(\mathbf{x})
\end{equation}
also belongs to $\mathcal{H}(\mathcal{M}_{\text{SSNM}})$ and has squared norm $\| s_{N}(\cdot)\|^{2}_{\mathcal{H}(\mathcal{M}_{\text{SSNM}})} = m_{0}^{2}+ \sum_{l \in [N]} \frac{m_{l}^{2}}{l!} \sigma^{2l}$. 

Now let us assume that \eqref{equ_cond_valid_bias_diagon_bias_SSNM} is satisfied. Since (for $N > N'$) 
\begin{equation} 
\| s_{N}(\cdot) - s_{N'}(\cdot)\|^{2}_{\mathcal{H}(\mathcal{M}_{\text{SSNM}})} =  \sum_{l \in [N]\setminus [N']} \frac{m_{l}^{2}}{l!} \sigma^{2l}, 
\end{equation}
we have from \eqref{equ_cond_valid_bias_diagon_bias_SSNM} that the sequence $\{ s_{N}(\cdot) \}_{N \rightarrow \infty}$ is a Cauchy sequence 
and therefore converges in the RKHS $\mathcal{H}(\mathcal{M}_{\text{SSNM}})$.
Since, by \eqref{equ_def_part_der_proof_SSNM_main_result_case_1} and \eqref{equ_prescr_mean_diag_bias}, this sequence moreover converges pointwise to the prescribed mean function $\gamma(\mathbf{x})$, we have by Theorem \ref{thm_point_convergent_impl_norm} that 
\begin{equation} 
\label{equ_proof_appendix_A_lim_N_s_N_gamma}
\lim_{N \rightarrow \infty} s_{N}(\cdot) = \gamma(\cdot), 
\end{equation}
which implies that $\gamma(\cdot) \in \mathcal{H}(\mathcal{M}_{\text{SSNM}})$ and in turn that the prescribed bias function $c(\cdot)$ is valid. We have then by Theorem \ref{thm_main_facts_RKHS_MVE} that 
\begin{equation}
L_{\mathcal{M}_{\text{SSNM}}} = \| \gamma(\cdot) \|^{2}_{\mathcal{H}(\mathcal{M}_{\text{SSNM}})} - \big[ \gamma(\mathbf{x}_{0}) \big]^{2} =
 \lim_{N \rightarrow \infty} \| s_{N}(\cdot)  \|^{2}_{\mathcal{H}(\mathcal{M}_{\text{SSNM}})} - m_{0}^{2}= \sum_{l \in \mathbb{N}}  \frac{m_{l}^{2}}{l!} \sigma^{2l}, 
\end{equation}
where we have used that $\gamma(\mathbf{x}_{0}) = m_{0}$ (according to \eqref{equ_prescr_mean_diag_bias}). 
Thus, we have verified the sufficiency of \eqref{equ_cond_valid_bias_diagon_bias_SSNM} for a bias function to be valid and for the expression \eqref{equ_expr_min_achiev_var_case_1_SSNM} for the minimum achievable variance. 

It remains to verify the necessity of \eqref{equ_cond_valid_bias_diagon_bias_SSNM} for a bias function to be valid. We will prove this by contradiction. 
Assume that \eqref{equ_cond_valid_bias_diagon_bias_SSNM} is not fulfilled, 
i.e., the sums $\sum_{l \in \mathcal{T}}  \frac{m_{l}^{2}}{l!} \sigma^{2l}$ with a finite index set $\mathcal{T} \subseteq \mathbb{Z}_{+}$ can be arbitrarily large, 
but the bias function $c(\cdot)$ is valid, i.e., $\gamma(\cdot) \in  \mathcal{H}(\mathcal{M}_{\text{SSNM}})$. By Theorem \ref{thm_orthog_proj_ineq}, we then have 
for any subspace $\mathcal{U} \triangleq \linspan\{ h^{(l)}(\cdot) \}_{l \in \mathcal{T}}$ with a finite set $\mathcal{T} \subseteq \mathbb{Z}_{+}$ that\footnote{Note that, trivially, $\{ h^{(l)}(\cdot) \}_{l \in \mathcal{T}}$ is an ONB for $\mathcal{U}$.} 
\begin{align}
\label{equ_proof_SSNM_sum_cond_valid_case_1}
 \| \gamma(\cdot) \|^{2}_{\mathcal{H}(\mathcal{M}_{\text{SSNM}})}
 & \stackrel{\eqref{equ_squared_norm_orthog_proj_pythag_thm}}{\geq}  \| \mathbf{P}_{\mathcal{U}} \gamma(\cdot) \|^{2}_{\mathcal{H}(\mathcal{M}_{\text{SSNM}})} \nonumber \\[4mm]
 & \stackrel{\eqref{equ_orth_proj_finite_dim_onb_squared_norm}}{=}\sum_{l \in \mathcal{T}}  \big\langle \gamma(\cdot), h^{(l)}(\cdot) \big\rangle^{2}_{\mathcal{H}(\mathcal{M}_{\text{SSNM}})} \nonumber \\[4mm]
 & \stackrel{\eqref{equ_der_reproduction_prop}}{=}\sum_{l \in \mathcal{T}}   \frac{\sigma^{2l}}{l!} \bigg[ \frac{\partial^{l \mathbf{e}_{k}} \gamma(\cdot)}{\partial \mathbf{x}^{l \mathbf{e}_{k}}} \bigg|_{\mathbf{x} = \mathbf{x}_{0}} \bigg]^{2} \nonumber \\[4mm]
 & \stackrel{\eqref{equ_prescr_mean_diag_bias}}{=}  \sum_{l \in \mathcal{T}}  \frac{m_{l}^{2}}{l!} \sigma^{2l},
\end{align}
where we also used Theorem \ref{thm_orth_proj_finite_dim_onb} and Theorem \ref{thm_der_repr_prop}. 
However, since \eqref{equ_cond_valid_bias_diagon_bias_SSNM} is not satisfied, i.e., we can make $ \sum_{l \in \mathcal{T}}  \frac{m_{l}^{2}}{l!} \sigma^{2l}$ arbitrarily large by using a suitable finite index set $\mathcal{T}$, 
it follows from \eqref{equ_proof_SSNM_sum_cond_valid_case_1} that the squared norm $ \| \gamma(\cdot) \|^{2}_{\mathcal{H}(\mathcal{M}_{\text{SSNM}})}$ is unbounded, which is impossible because $\gamma(\cdot) \in \mathcal{H}(\mathcal{M}_{\text{SSNM}})$. 

Finally, we will show \eqref{equ_expr_LMV_diag_bias_case_1_SSNM}. If the prescribed bias is valid, we can by \eqref{equ_proof_appendix_A_lim_N_s_N_gamma} represent the corresponding mean as 
\begin{equation} 
\label{equ_proof_appendix_A_lim_N_gamma_summe_h_l}
\gamma(\cdot) = \lim_{N \rightarrow \infty} s_{N}(\cdot)  = \sum_{l \in \mathbb{Z}_{+}} \frac{m_{l}}{\sqrt{l!}} \sigma^{l} h^{(l)}(\cdot).
\end{equation} 
From this, we obtain the LMV estimator in \eqref{equ_expr_LMV_diag_bias_case_1_SSNM} 
from \eqref{equ_def_part_der_proof_SSNM_main_result_case_1}, Theorem \ref{thm_der_repr_prop}, 
Theorem \ref{thm_isometry_RKHS_rhos_derivative_kernel}, Theorem \ref{thm_main_facts_RKHS_MVE} and the identity 
\begin{align}
\label{equ_proof_main_thm_diag_bias_SSNM_cases_identiy_case_1_hermite}
 & \frac{m_{l} \sigma^{2l}}{l!} \frac{ \partial^{l} \exp \left( - \frac{1}{2 \sigma^{2}} ( 2 y(x_{0} -x)+ x^{2} - x_{0}^{2}) \right)}{ \partial x^{l}} \bigg|_{x = x_{0}}  \nonumber \\[4mm]
 & = \frac{m_{l} \sigma^{2l}}{l!} \frac{ \partial^{l} \exp \left( - \frac{1}{2 \sigma^{2}} ( (x-y)^{2} - y^{2} + 2yx_0 - x_{0}^{2}) \right)}{ \partial x^{l}} \bigg|_{x = x_{0}}  \nonumber \\[4mm] 
& = \frac{m_{l} \sigma^{2l}}{l!}  \exp \left(  \frac{1}{2 \sigma^{2}}  (y-x_0)^{2} \right)\frac{ \partial^{l} \exp \left( - \frac{1}{2 \sigma^{2}}  (x-y)^{2} \right)}{ \partial x^{l}} \bigg|_{x = x_{0}}  \nonumber \\[4mm] 
& {\scriptstyle \bigg( x' \triangleq \frac{x - y}{\sigma}  \bigg)}  \nonumber \\[4mm]Ê
&  = \frac{m_{l} \sigma^{l}}{l!}  \exp \left( \frac{1}{2 \sigma^{2}}  (x_0-y)^{2} \right)\frac{ \partial^{l} \exp \left( - \frac{1}{2}  x'^{2} \right)}{ \partial x'^{l}} \bigg|_{x' = \frac{x_{0}-y}{\sigma}}  \nonumber \\[4mm] 
&  = \frac{m_{l} \sigma^{l}}{l!} (-1)^{l} \mathsf{H}_{l}\left( \frac{ x_{0} - y}{\sigma} \right) = \frac{m_{l} \sigma^{l}}{l!} \mathsf{H}_{l}\left( \frac{y- x_{0}}{\sigma} \right). 
\end{align}
Indeed, we have 
\begin{align} 
\hat{x}_{k}^{(\mathbf{x}_{0})} (\mathbf{y}) &\stackrel{\eqref{equ_lmv_estimator_general_congruence_L_M_RKHS}}{=} \mathsf{J} \big[ \gamma(\cdot) \big] \stackrel{\eqref{equ_proof_appendix_A_lim_N_gamma_summe_h_l}}{=}\mathsf{J} \bigg[ \sum_{l \in \mathbb{Z}_{+}} \frac{m_{l}}{\sqrt{l!}} \sigma^{l} h^{(l)}(\cdot)\bigg] =  \sum_{l \in \mathbb{Z}_{+}} \frac{m_{l}}{\sqrt{l!}} \sigma^{l} \mathsf{J} \big[ h^{(l)}(\cdot)\big]  \nonumber \\[4mm]
& \stackrel{\eqref{equ_def_part_der_proof_SSNM_main_result_case_1}}{=} \sum_{l \in \mathbb{Z}_{+}} \frac{m_{l}}{l!} \sigma^{2l}  \mathsf{J} 
\bigg[  \frac{\partial^{l \mathbf{e}_{k}} R_{\mathcal{M}_{\text{SSNM}}}(\cdot, \mathbf{x})}{\partial \mathbf{x}^{l \mathbf{e}_{k}}} \bigg|_{\mathbf{x} = \mathbf{x}_{0}}\bigg] \nonumber \\[4mm]
& \stackrel{\eqref{equ_isometry_RKHS_rhos_derivative_kernel}}{=} \sum_{l \in \mathbb{Z}_{+}} \frac{m_{l}}{l!} \sigma^{2l}    \frac{\partial^{l \mathbf{e}_{k}} \rho_{\mathcal{M}_{\text{SSNM}}}(\mathbf{y}, \mathbf{x})}{\partial \mathbf{x}^{l \mathbf{e}_{k}}} \bigg|_{\mathbf{x} = \mathbf{x}_{0}}
\nonumber \\[4mm]
& \stackrel{\eqref{equ_rho_M_SSNM}}{=} \sum_{l \in \mathbb{Z}_{+}} \frac{m_{l}\sigma^{2l}}{l!}  \frac{\partial^{l \mathbf{e}_{k}} \exp \left( \frac{1}{\sigma^{2}} \mathbf{y}^{T}(\mathbf{x}-\mathbf{x}_{0}) - \frac{1}{2 \sigma^{2}}(\| \mathbf{x} \|^{2}_{2} - \|Ê\mathbf{x}_{0} \|^{2}_{2}) \right)}{\partial \mathbf{x}^{l \mathbf{e}_{k}}} \bigg|_{\mathbf{x} = \mathbf{x}_{0}}\nonumber \\[4mm]
& \stackrel{\eqref{equ_proof_main_thm_diag_bias_SSNM_cases_identiy_case_1_hermite}}{=} \sum_{l \in \mathbb{Z}_{+}} \frac{m_{l} \sigma^{l}}{l!} \mathsf{H}_{l}\left( \frac{y_k - x_{0,k}}{\sigma} \right).
\end{align} 

\paragraph{Case II: $| \supp(\mathbf{x}_{0}) \cup \{k\}| = S+1$:}$~~$\\ 
We will assume without loss of generality\footnote{Indeed, 
consider the minimum variance problem $\mathcal{M}_{\text{SSNM}} = \left( \mathcal{E}_{\text{SSNM}},c(\cdot),\mathbf{x}_{0} \right)$ 
with observation $\mathbf{y} \in \mathbb{R}^{N}$ and an arbitrary support $\supp(\mathbf{x}_{0})$ and index $k \in [N] \setminus \supp(\mathbf{x}_{0})$. 
According to Section \ref{sec_invariance_mve}, we then have that the minimum variance problem $\mathcal{M}'$ that arises from $\mathcal{M}_{\text{SSNM}}$ by applying an invertible
permutation matrix $\mathbf{P} \in \{0,1\}^{N \times N}$ to the observation $\mathbf{y}$, is completely 
equivalent to $\mathcal{M}_{\text{SSNM}}$. 
By a suitable choice of the permutation matrix $\mathbf{P}$, the modified minimum variance problem $\mathcal{M}'$ coincides with 
$\mathcal{M}_{\text{SSNM}}=\left( \mathcal{E}_{\text{SSNM}},c(\cdot),\mathbf{x}'_{0} \right)$ where $\supp(\mathbf{x}'_{0}) = [S]$ and $k=S+1$.} that $k = S+1$ and $\supp(\mathbf{x}_{0}) = [S]$.
The proof for Case II is similar to that for Case I. However, a considerable difference is that 
the functions $h^{(l)}(\cdot)$ defined in \eqref{equ_def_part_der_proof_SSNM_main_result_case_1} do not belong to the RKHS $\mathcal{H}(\mathcal{M}_{\text{SSNM}})$ in this case. 

Consider the functions 
\begin{equation}
\label{equ_appendix_A_def_v_l}
v^{(l)}(\cdot): \mathcal{X}_{S} \rightarrow \mathbb{R}: v^{(l)}(\mathbf{x}) \triangleq \frac{1}{\sigma^{l} l!} x_{k}^{l} \nu_{\mathbf{x}_{0}}(\mathbf{x}), 
\end{equation}
where $l \in \mathbb{Z}_{+}$ and $\nu_{\mathbf{x}_{0}}(\mathbf{x}) = \exp\left( - \frac{1}{2 \sigma^{2}} \|\mathbf{x}_{0} \|^{2}_{2} + Ê\frac{1}{\sigma} \mathbf{x}^{T} \mathbf{x}_{0}\right)$ as defined in Theorem \ref{thm_general_min_achiev_var_LMV_SSNM}. 

Using the Taylor expansion $\exp(x) = \sum_{l \in \mathbb{Z}_{+}} \frac{x^{l}}{l!}$ of the exponential function \cite{RudinBook}, we obtain 
\begin{align}
\label{equ_proof_ssnm_main_case_2_def_part_der_expansion_onb}
v^{(l)}(\mathbf{x})  & = \frac{1}{\sigma^{l} l!}  x_{k}^{l} \nu_{\mathbf{x}_{0}}(\mathbf{x})  =\frac{1}{\sigma^{l} l!} \exp \left( - \frac{ \| \mathbf{x}_{0} \|^{2}_{2}}{2 \sigma^{2}} \right) x_{k}^{l} 
\sum_{\mathbf{p} \in \mathbb{Z}_{+}^{S}} 
 \frac{1}{\mathbf{p}!} \prod_{j \in [S]}  \left( \frac{ x_{0,j} x_{j}}{\sigma}\right)^{p_{j}} \nonumber \\[4mm]
& \stackrel{(a)}{=}\exp \left( - \frac{ \| \mathbf{x}_{0} \|^{2}_{2}}{2 \sigma^{2}} \right) \sum_{\mathbf{p} \in \mathcal{A}^{(l)}} \frac{1}{\mathbf{p}!} \widetilde{\mathbf{x}}_{0}^{\mathbf{p}} \mathbf{x}^{\mathbf{p}} \nonumber \\[4mm]
&=  \exp \left( - \frac{ \| \mathbf{x}_{0} \|^{2}_{2}}{2 \sigma^{2}} \right) \sum_{\mathbf{p} \in \mathcal{A}^{(l)}} \frac{\widetilde{\mathbf{x}}_{0}^{\mathbf{p}}}{\sqrt{\mathbf{p}!}}  g^{(\mathbf{p})} (\mathbf{x}),
\end{align}
where  $\widetilde{\mathbf{x}}_{0} \triangleq \frac{1}{\sigma}\big( x_{0,1},  \ldots, x_{0,S}, 1,Ê\ldots, 1 \big)^{T} \in \mathbb{R}^{N}$, $\mathcal{A}^{(l)} \triangleq \left\{ \mathbf{p} \in  \mathbb{Z}_{+}^{N} \cap \mathcal{X}_{S}  | p_{k} = l \mbox{, } \supp(\mathbf{p}) \subseteq [S+1]  \right\}$, and the functions $g^{(\mathbf{p})}(\cdot): \mathcal{X}_{S} \rightarrow \mathbb{R}: 
g^{(\mathbf{p})}(\mathbf{x}) = \frac{1}{\sqrt{\mathbf{p}!}} \frac{\partial^{\mathbf{p}} R_{e}(\mathbf{x}, \mathbf{x}_{2})}{\partial \mathbf{x}_{2}^{\mathbf{p}}} \big|_{\mathbf{x}_{2}=0}$ as 
defined in \eqref{equ_def_func_g_p_part_der_R_e} of Theorem \ref{thm_entire_charact_RKHS_SSNM}).
The step $(a)$ in \eqref{equ_proof_ssnm_main_case_2_def_part_der_expansion_onb} follows from the 
fact that the function $v^{(l)}(\mathbf{x})$ is defined only on the domain $\mathcal{X}_{S}$ on which every monomial $\mathbf{x}^{\mathbf{p}}$ ($\mathbf{p} \in \mathbb{Z}_{+}^{N}$) with $\| \mathbf{p}  \|_{0} > S$ vanishes.

By Theorem \ref{thm_entire_charact_RKHS_SSNM}, the functions $g^{(\mathbf{p})}(\cdot)$ belong to $\mathcal{H}(R_{e})$ and are orthonormal, i.e., 
\begin{equation}
\big\langle g^{(\mathbf{p})}(\cdot) , g^{(\mathbf{p}')}(\cdot) \big\rangle_{\mathcal{H}(R_{e})} = \delta_{\mathbf{p},\mathbf{p}'}.
\end{equation} 
Since the coefficient sequence of \eqref{equ_proof_ssnm_main_case_2_def_part_der_expansion_onb} satisfies
$\frac{\widetilde{\mathbf{x}}_{0}^{\mathbf{p}}}{\sqrt{\mathbf{p}!}} \in \ell^{2}(\mathbb{Z}_{+}^{N} \cap \mathcal{X}_{S})$, we then have also that $v^{(l)}(\cdot) \in \mathcal{H}(R_{e})$. 

Based on the orthonormality of the functions $\left\{ g^{(\mathbf{p})}(\cdot) \right\}_{\mathbf{p} \in \mathcal{A}^{(l)}}$, we can calculate the squared norm of $v^{(l)}(\cdot)$ as  
\begin{align} 
\label{equ_proof_main_SSNM_case_2_atom_squared_norm_part_1}
\| v^{(l)}(\cdot) \|_{\mathcal{H}(R_{e})}^{2} &  =  \exp \left( - \frac{ \| \mathbf{x}_{0} \|^{2}_{2}}{ \sigma^{2}} \right) \sum_{\mathbf{p} \in \mathcal{A}^{(l)}} 
\frac{\widetilde{\mathbf{x}}_{0}^{2\mathbf{p}}}{\mathbf{p}!}   \nonumber \\[4mm]
&  =   \exp \left( - \frac{ \| \mathbf{x}_{0} \|^{2}_{2}}{ \sigma^{2}} \right) \sum_{\mathbf{p} \in \mathcal{A}^{(l)}} \prod_{j \in [N]} \frac{\widetilde{x}_{0,j}^{2 p_j }}{p_j !}   \nonumber \\[4mm]Ê
  & \stackrel{(a)}{=}  \exp \left( - \frac{ \| \mathbf{x}_{0} \|^{2}_{2}}{ \sigma^{2}} \right) \sum_{m \in [S]} \sum_{\mathbf{p} \in \mathcal{A}^{(l)}_{m}} \prod_{j \in [N]} 
  \frac{\widetilde{x}_{0,j}^{2 p_j }}{p_j !},
\end{align} 
where 
\begin{align}
\mathcal{A}^{(l)}_{m} & \triangleq \big\{ \mathbf{p} \in  \mathbb{Z}_{+}^{N} \cap \mathcal{X}_{S}  \big| p_{n} > 0 \mbox{ for } n \in [m-1] \mbox{, } p_{m}= 0 \mbox{, } p_{k} = l \mbox{, } \supp(\mathbf{p}) \subseteq [S+1]  \big\}  \nonumber \\[4mm] 
& =  \big\{ \mathbf{p} \in \mathbb{Z}_{+}^{N} \big| p_{n} > 0 \mbox{ for } n \in [m-1] \mbox{, } p_{m}= 0 \mbox{, } p_{n} \geq 0 \mbox{ for } n \in [S] \setminus [m] \mbox{, }p_{k} = l \mbox{, } \nonumber \\[4mm]
& \hspace*{10mm}  \supp(\mathbf{p}) \subseteq [S+1]  \big\}. 
\end{align} 
The step $(a)$ in \eqref{equ_proof_main_SSNM_case_2_atom_squared_norm_part_1} follows from the fact that the sets $\big\{ \mathcal{A}^{(l)}_{m} \big\}_{m \in [S]}$ form a 
partition of the set $\mathcal{A}^{(l)}$, i.e., the sets are disjoint ($\mathcal{A}^{(l)}_{m} \cap \mathcal{A}^{(l)}_{m'} = \emptyset$ for $m \neq m'$) and their union 
equals $\mathcal{A}^{(l)}$, i.e., $\mathcal{A}^{(l)} = \bigcup_{m \in [S]} \mathcal{A}^{(l)}_{m}$.
Further developing \eqref{equ_proof_main_SSNM_case_2_atom_squared_norm_part_1} yields 
\begin{align}
\label{equ_proof_main_SSNM_case_2_atom_squared_norm}
\| v^{(l)}(\cdot) \|_{\mathcal{H}(R_{e})}^{2}   & = 
\exp \left( - \frac{ \| \mathbf{x}_{0} \|^{2}_{2}}{ \sigma^{2}} \right) \sum_{m \in [S]} \sum_{\mathbf{p} \in \mathcal{A}^{(l)}_{m}} \Bigg(\prod_{n \in [m-1]} 
  \frac{\widetilde{x}_{0,n}^{2 p_n }}{p_n !} \Bigg) \frac{\widetilde{x}_{0,m}^{2 p_m }}{p_m !}  \Bigg(\prod_{n' \in [S] \setminus [m]} 
  \frac{\widetilde{x}_{0,n'}^{2 p_n' }}{p_n' !}\Bigg) \frac{\widetilde{x}_{0,k}^{2 p_k }}{p_k !}  \nonumber \\[4mm]
  & \hspace*{-15mm} \stackrel{p_{m}=0\mbox{, } p_{k}=l}{=} \frac{1}{\sigma^{2l} l!} \exp \left( - \frac{ \| \mathbf{x}_{0} \|^{2}_{2}}{ \sigma^{2}} \right) \sum_{m \in [S]} \Bigg( \prod_{n \in [m-1]} \sum_{p_{n} \in \mathbb{N}} \frac{\widetilde{x}_{0,n}^{2 p_n }}{p_n !} \Bigg) 
  \Bigg(\prod_{n' \in [S] \setminus [m]} \sum_{p_{n'} \in \mathbb{Z}_{+}} \frac{\widetilde{x}_{0,n'}^{2 p_{n'} }}{p_{n'} !} \Bigg)\nonumber \\[4mm] 
  & \hspace*{-15mm} =  \frac{1}{\sigma^{2l} l!} \exp \left( - \frac{ \| \mathbf{x}_{0} \|^{2}_{2}}{ \sigma^{2}} \right) \sum_{m \in [S]} \prod_{n \in [m-1]} \Bigg[ \exp \bigg(  \frac{ x_{0,n}^{2}}{ \sigma^{2}} \bigg) - 1 \Bigg]
 \prod_{n' \in [S] \setminus [m]}  \exp \bigg(  \frac{ x_{0,n'}^{2}}{ \sigma^{2}} \bigg) \nonumber \\[4mm] 
 & \hspace*{-15mm} = \frac{1}{\sigma^{2l} l!}  \sum_{m \in [S]}\exp \bigg(  -\frac{ x_{0,m}^{2}}{ \sigma^{2}} \bigg)  \prod_{n \in [m-1]} \Bigg[ 1 -\exp \bigg(  -\frac{ x_{0,n}^{2}}{ \sigma^{2}} \bigg)  \Bigg].
\end{align} 
It can also be verified by 
the orthonormality (cf.\ Theorem \ref{thm_entire_charact_RKHS_SSNM}) of the functions $g^{(\mathbf{p})}(\cdot)$, that the functions $\{v^{(l)}(\cdot) \}_{l \in \mathbb{Z}_{+}}$ are orthogonal, i.e., $\big\langle v^{(l)}(\cdot) , v^{(l')}(\cdot) \big\rangle_{\mathcal{H}(R_{e})} = \delta_{l,l'} \| v^{(l)}(\cdot) \|^{2}_{\mathcal{H}(R_{e})}$.

Now consider the function 
\begin{equation} 
\label{equ_partial_sum_converges_image_mean_valid_bias_case_2_ssnm}
w_{N}(\cdot): \mathcal{X}_{S} \rightarrow \mathbb{R}: w_{N}(\mathbf{x}) \triangleq \sum_{l \in [N]} m_{l} \sigma^{2l} v^{(l)}(\mathbf{x})  + m_{0}v^{(0)}(\mathbf{x}), 
\end{equation} 
which obviously belongs to the RKHS $\mathcal{H}(R_{e})$ as the functions $v^{(l)}(\cdot)$ do. 
In the following we will use the relation 
\begin{equation}
\label{equ_appendix_a_relation_squared_norm_v_l_v_0}
\| v^{(l)} (\cdot) \|^{2}_{\mathcal{H}(R_{e})} = \frac{1}{\sigma^{2l} l!}\| v^{(0)} (\cdot) \|^{2}_{\mathcal{H}(R_{e})},
\end{equation}
which can be verified easily by \eqref{equ_proof_main_SSNM_case_2_atom_squared_norm}. 
Based on \eqref{equ_proof_main_SSNM_case_2_atom_squared_norm} 
and the orthogonality of the functions  $\{v^{(l)}(\cdot) \}_{l \in \mathbb{Z}_{+}}$ we obtain 
\begin{align}
\label{equ_expr_norm_sequ_gamma_case_2_ssnm_main_result}
\| w_{N}(\cdot) \|_{\mathcal{H}(R_{e})}^{2} & = m_{0}^{2}\| v^{(0)} (\cdot) \|^{2}_{\mathcal{H}(R_{e})} + \sum_{l \in [N]} m^{2}_{l} \sigma^{4l} \| v^{(l)} (\cdot) \|^{2}_{\mathcal{H}(R_{e})}    \nonumber \\[4mm]
& \stackrel{\eqref{equ_appendix_a_relation_squared_norm_v_l_v_0}}{=}   \bigg( m_{0}^{2} +  \sum_{l \in [N]} \frac{m_{l}^{2}}{l!} \sigma^{2l} \bigg) \| v^{(0)} (\cdot) \|^{2}_{\mathcal{H}(R_{e})}. 
\end{align} 

If \eqref{equ_cond_valid_bias_diagon_bias_SSNM} is satisfied, we have by \eqref{equ_expr_norm_sequ_gamma_case_2_ssnm_main_result} 
that the sequence $\{ w_{N}(\cdot) \}_{N \rightarrow \infty}$ is a Cauchy sequence and therefore converges in the RKHS $\mathcal{H}(R_{e})$. Indeed, by 
a similar calculation as used in \eqref{equ_expr_norm_sequ_gamma_case_2_ssnm_main_result}, we have for $N' > N$ that 
\begin{equation}
 \| w_{N'}(\cdot) - w_{N}(\cdot) \|_{\mathcal{H}(R_{e})}^{2} =  \| v^{(0)} (\cdot) \|^{2}_{\mathcal{H}(R_{e})} \sum_{l \in [N'] \setminus [N]} \frac{m_{l}^{2}}{l!} \sigma^{2l}, 
\end{equation} 
which must become arbitrarily small for $N',N$ sufficiently large since otherwise the sum in \eqref{equ_cond_valid_bias_diagon_bias_SSNM} cannot be finite.
Moreover we have that the sequence $\{ w_{N}(\cdot) \}_{N \rightarrow \infty}$ converges pointwise to the function $\tilde{\gamma}(\cdot): \mathcal{X}_{S} \rightarrow \mathbb{R}: \tilde{\gamma}(\mathbf{x}) \triangleq \gamma(\sigma \mathbf{x}) \nu_{\mathbf{x}_{0}}(\mathbf{x})=\mathsf{K}_{e}[\gamma(\cdot)]$. 
Here, we denote by $\mathsf{K}_{e}[\cdot]: \mathcal{H}(\mathcal{M}_{\text{SSNM}}) \rightarrow \mathcal{H}(R_{e})$ the congruence defined by \eqref{equ_def_congruence_SSNM} in Theorem \ref{thm_entire_charact_RKHS_SSNM}. 
This pointwise convergence can be seen from
\begin{align}
\lim_{N \rightarrow \infty}  w_{N}(\mathbf{x}) & \stackrel{\eqref{equ_partial_sum_converges_image_mean_valid_bias_case_2_ssnm}}{=}
\lim_{N \rightarrow \infty} \sum_{l \in [N]} m_{l} \sigma^{2l} v^{(l)}(\mathbf{x})  + m_{0}v^{(0)}(\mathbf{x}) \nonumber \\[4mm]
& \stackrel{\eqref{equ_appendix_A_def_v_l}}{=} \lim_{N \rightarrow \infty} \sum_{l \in [N]} m_{l} \sigma^{2l} \frac{1}{\sigma^{l} l!} x_{k}^{l} \nu_{\mathbf{x}_{0}}(\mathbf{x})  
+ m_{0} \nu_{\mathbf{x}_{0}}(\mathbf{x}) \nonumber \\[4mm]
& =  \nu_{\mathbf{x}_{0}}(\mathbf{x})\sum_{l \in \mathbb{Z}_{+}} \frac{m_{l} \sigma^{l} }{ l!} x_{k}^{l} \nonumber \\[4mm]
& \stackrel{\eqref{equ_prescr_mean_diag_bias}}{=}  \gamma(\sigma \mathbf{x}) \nu_{\mathbf{x}_{0}}(\mathbf{x}).
\end{align} 
Therefore, we have by Theorem \ref{thm_point_convergent_impl_norm} that 
\begin{equation} 
\label{equ_Cauchy_sequ_converges_image_mean_valid_bias_case_2_ssnm}
\lim_{N \rightarrow \infty} w_{N}(\cdot) = \tilde{\gamma}(\cdot) 
\end{equation}
and therefore $\tilde{\gamma}(\cdot) \in \mathcal{H}(R_{e})$. 
However, since $\tilde{\gamma}(\cdot)$ is the image of the prescribed mean function $\gamma(\cdot)$ under the congruence $\mathsf{K}_{e}[\cdot]: \mathcal{H}(\mathcal{M}_{\text{SSNM}}) \rightarrow \mathcal{H}(R_{e})$, we have that $\gamma(\cdot) \in \mathcal{H}(\mathcal{M}_{\text{SSNM}})$, i.e., 
that the prescribed bias function $c(\cdot)$ is valid. 

We then have also by Theorem \ref{thm_main_facts_RKHS_MVE} and Theorem \ref{thm_entire_charact_RKHS_SSNM} that 
\begin{align}
L_{\mathcal{M}_{\text{SSNM}}} 
& \stackrel{\eqref{equ_squared_norm_min_achiev_var}}{=} \|  \gamma(\cdot) \|_{\mathcal{H}(\mathcal{M}_{\text{SSNM}})}^{2} - \big[ \gamma(\mathbf{x}_{0}) \big]^{2}
\nonumber \\[4mm]
& \stackrel{\eqref{equ_def_congruence_SSNM}}{=}
 \|  \gamma(\sigma \mathbf{x}) \nu_{\mathbf{x}_{0}}(\mathbf{x}) \|_{\mathcal{H}(R_{e})}^{2}    - \big[ \gamma(\mathbf{x}_{0}) \big]^{2}Ê\nonumber \\[4mm]Ê
& =  \|  \tilde{\gamma}(\cdot)\|_{\mathcal{H}(R_{e})}^{2}    - \left( \gamma(\mathbf{x}_{0}) \right)^{2}Ê \nonumber \\[4mm]
& \stackrel{\eqref{equ_expr_norm_sequ_gamma_case_2_ssnm_main_result},\eqref{equ_Cauchy_sequ_converges_image_mean_valid_bias_case_2_ssnm}}{=} 
\lim_{N \rightarrow \infty}  \bigg( m_{0}^{2} +  \sum_{l \in [N]} \frac{m_{l}^{2}}{l!} \sigma^{2l} \bigg) \| v^{(0)} (\cdot) \|^{2}_{\mathcal{H}(R_{e})}- \big[ \gamma(\mathbf{x}_{0}) \big]^{2}Ê \nonumber \\[4mm]Ê
&  \stackrel{\eqref{equ_proof_main_SSNM_case_2_atom_squared_norm}}{=}    \bigg(  \sum_{l \in \mathbb{Z}_{+}} \frac{m_{l}^{2}}{l!} \sigma^{2l} \bigg) \sum_{m \in [S]}\exp \bigg(  -\frac{ x_{0,m}^{2}}{ \sigma^{2}} \bigg) 
 \prod_{n \in [m-1]} \Bigg[ 1 -\exp \bigg(  -\frac{ x_{0,n}^{2}}{ \sigma^{2}} \bigg)  \Bigg]    -  \big[ \gamma(\mathbf{x}_{0}) \big]^{2},
\end{align} 
which is equal to \eqref{equ_min_achiev_var_case_2_ssnm_main}. 

In order to verify the expression \eqref{equ_lmv_case_2_ssnm_main} for the LMV estimator, we first consider the function 
\begin{equation} 
\label{equ_appendix_a_def_r_l}
r^{(l)}(\cdot) \triangleq \mathsf{K}_{e}^{-1}[ v^{(l)}(\cdot)].
\end{equation} 
From \eqref{equ_proof_ssnm_main_case_2_def_part_der_expansion_onb}, \eqref{equ_def_congruence_SSNM}, and \eqref{equ_def_func_g_p_part_der_R_e} it follows that
\begin{align} 
\label{equ_proof_main_ssnm_case_2_expansion_atom_lmv}
r^{(l)}(\mathbf{x}) & =  \mathsf{K}_{e}^{-1}[ v^{(l)}(\cdot)] \nonumber \\[4mm]
& \stackrel{\eqref{equ_proof_ssnm_main_case_2_def_part_der_expansion_onb}}{=} \mathsf{K}_{e}^{-1} \bigg[   \exp \left( - \frac{ \| \mathbf{x}_{0} \|^{2}_{2}}{2 \sigma^{2}} \right) \sum_{\mathbf{p} \in \mathcal{A}^{(l)}} \frac{\widetilde{\mathbf{x}}_{0}^{\mathbf{p}}}{\sqrt{\mathbf{p}!}}  g^{(\mathbf{p})} (\mathbf{x}) \bigg]  \nonumber \\[4mm]
& \stackrel{\eqref{equ_def_congruence_SSNM}}{=} \frac{1}{\nu_{\mathbf{x}_{0}} \big(\frac{\mathbf{x}}{ \sigma} \big)}    \exp \left( - \frac{ \| \mathbf{x}_{0} \|^{2}_{2}}{2 \sigma^{2}} \right) 
\sum_{\mathbf{p} \in \mathcal{A}^{(l)}} \frac{\widetilde{\mathbf{x}}_{0}^{\mathbf{p}}}{\mathbf{p}!}   g^{(\mathbf{p})} \bigg(\frac{\mathbf{x}}{\sigma}\bigg) \nonumber \\[4mm]
& \stackrel{\eqref{equ_def_func_g_p_part_der_R_e}}{=}
\frac{1}{\nu_{\mathbf{x}_{0}} \big(\frac{\mathbf{x}}{ \sigma} \big)}    \exp \left( - \frac{ \| \mathbf{x}_{0} \|^{2}_{2}}{2 \sigma^{2}} \right) 
\sum_{\mathbf{p} \in \mathcal{A}^{(l)}} \frac{\widetilde{\mathbf{x}}_{0}^{\mathbf{p}}}{\mathbf{p}!}  \frac{\partial^{\mathbf{p}} R_{e}\left(\frac{\mathbf{x}}\sigma, \mathbf{x}_{2}Ê\right)}{ \partial \mathbf{x}_{2}^{\mathbf{p}}} \bigg|_{\mathbf{x}_{2}= \mathbf{0}} 
\nonumber \\[4mm] 
& \stackrel{\eqref{equ_def_R_e}}{=}   \exp \left( - \frac{1}{ \sigma^{2}} \mathbf{x}^{T} \mathbf{x}_{0} \right) 
\sum_{\mathbf{p} \in \mathcal{A}^{(l)}} \frac{\widetilde{\mathbf{x}}_{0}^{\mathbf{p}}}{\mathbf{p}!}  \frac{\partial^{\mathbf{p}} \exp \left( \frac{1}{\sigma} \mathbf{x}^{T} \mathbf{x}_{2} \right) }{ \partial \mathbf{x}_{2}^{\mathbf{p}}} \bigg|_{\mathbf{x}_{2}= \mathbf{0}}  \nonumber \\[4mm] 
& \stackrel{\eqref{equ_def_kernel_SSNM}}{=} \sum_{\mathbf{p} \in \mathcal{A}^{(l)}} \frac{\widetilde{\mathbf{x}}_{0}^{\mathbf{p}}}{\mathbf{p}!}  \frac{\partial^{\mathbf{p}}  \left[ R_{\mathcal{M}_{\text{SSNM}}}(\mathbf{x}, \sigma \mathbf{x}_{2})
 \exp \left( -\frac{1}{\sigma^{2}} \| \mathbf{x}_{0} \|^{2}_{2} + \frac{1}{\sigma} \mathbf{x}^{T}_{2}\mathbf{x}_{0} \right) \right] }{ \partial \mathbf{x}_{2}^{\mathbf{p}}} \bigg|_{\mathbf{x}_{2}= \mathbf{0}}.
\end{align} 
The image of $r^{(l)}(\cdot)$ under the congruence $\mathsf{J}[\cdot]: \mathcal{H}(\mathcal{M}_{\text{SSNM}}) \rightarrow \mathcal{L}(\mathcal{M}_{\text{SSNM}})$ as 
defined in Theorem \ref{thm_isometry_RKHS_rhos} can be calculated due to the expansion \eqref{equ_proof_main_ssnm_case_2_expansion_atom_lmv} and Theorem \ref{thm_isometry_RKHS_rhos_derivative_kernel} as 
{\allowdisplaybreaks \begin{align}
\mathsf{J}[ r^{(l)}(\cdot)](\mathbf{y}) & \stackrel{\eqref{equ_isometry_RKHS_rhos_derivative_kernel}}{=} \sum_{\mathbf{p} \in \mathcal{A}^{(l)}} \frac{\widetilde{\mathbf{x}}_{0}^{\mathbf{p}}}{\mathbf{p}!}  \frac{\partial^{\mathbf{p}}  \left[ \rho_{\mathcal{M}_{\text{SSNM}}}(\mathbf{y}, \sigma \mathbf{x}_{2})
 \exp \left( -\frac{1}{\sigma^{2}} \| \mathbf{x}_{0} \|^{2}_{2} + \frac{1}{\sigma} \mathbf{x}_{2}^{T}\mathbf{x}_{0} \right) \right] }{ \partial \mathbf{x}_{2}^{\mathbf{p}}} \bigg|_{\mathbf{x}_{2}= \mathbf{0}} \nonumber \\[4mm]Ê
 & \stackrel{\eqref{equ_rho_M_SSNM}}{=}   \exp \left( -\frac{1}{2 \sigma^{2}} \| \mathbf{x}_{0} \|^{2}_{2} \right) \sum_{\mathbf{p} \in \mathcal{A}^{(l)}} \frac{\widetilde{\mathbf{x}}_{0}^{\mathbf{p}}}{\mathbf{p}!}  \frac{\partial^{\mathbf{p}} 
 \exp \left(  \frac{1}{\sigma^{2}} \mathbf{y}^{T}(\sigma \mathbf{x}_{2} - \mathbf{x}_{0}) + \frac{1}{\sigma} \mathbf{x}_{2}^{T} \mathbf{x}_{0} - \frac{1}{2} \| \mathbf{x}_{2} \|^{2}_{2} \right) }
 { \partial \mathbf{x}_{2}^{\mathbf{p}}} \bigg|_{\mathbf{x}_{2}= \mathbf{0}}. 
 \end{align}}
 Using the partition $\mathcal{A}^{(l)} = \bigcup_{m \in [S]} \mathcal{A}^{(l)}_{m}$, this becomes 
{\allowdisplaybreaks \begin{align}
 \mathsf{J}[ r^{(l)}(\cdot)](\mathbf{y})    
&=  \exp \left( -\frac{1}{2 \sigma^{2}} \| \mathbf{x}_{0} \|^{2}_{2} \right) \sum_{m \in [S]}  \sum_{\mathbf{p} \in \mathcal{A}^{(l)}_{m}} \frac{\widetilde{\mathbf{x}}_{0}^{\mathbf{p}}}{\mathbf{p}!}  \frac{\partial^{\mathbf{p}} 
 \exp \left(  \frac{1}{\sigma^{2}} \mathbf{y}^{T}(\sigma \mathbf{x}_{2} - \mathbf{x}_{0}) + \frac{1}{\sigma} \mathbf{x}_{2}^{T} \mathbf{x}_{0} - \frac{1}{2} \| \mathbf{x}_{2} \|^{2}_{2} \right) }{ \partial \mathbf{x}_{2}^{\mathbf{p}}} \bigg|_{\mathbf{x}_{2}= \mathbf{0}}  \nonumber \\[4mm]
 & \hspace*{-10mm} =  \exp \left( -\frac{1}{2 \sigma^{2}} \| \mathbf{x}_{0} \|^{2}_{2} \right)  \nonumber \\[4mm]
 & \hspace*{-5mm} \times \sum_{m \in [S]}  \sum_{\mathbf{p} \in \mathcal{A}^{(l)}_{m}} 
 \prod_{n \in [N]} \frac{\widetilde{x}_{0,n}^{p_{n}}}{p_{n}!} \frac{\partial^{p_{n}} 
 \exp \left(  \frac{1}{\sigma^{2}} y_{n} (\sigma x_{2,n} - x_{0,n}) + \frac{1}{\sigma} x_{2,n} x_{0,n} - \frac{1}{2}  x_{2,n}^{2}\right) }{ \partial x_{2,n}^{p_n}} \bigg|_{x_{2,n}= 0} \nonumber \\[4mm]Ê
   & \hspace*{-10mm} =  \exp \left( -\frac{1}{2 \sigma^{2}} \| \mathbf{x}_{0} \|^{2}_{2} \right) \sum_{m \in [S]} \sum_{\mathbf{p} \in \mathcal{A}^{(l)}_{m}} 
 \frac{\widetilde{x}_{0,m}^{p_{m}}}{p_{m}!} \frac{\partial^{p_{m}} 
 \exp \left(  \frac{1}{\sigma^{2}} y_{m} (\sigma x_{2,m} - x_{0,m}) + \frac{1}{\sigma} x_{2,m} x_{0,m} - \frac{1}{2}  x_{2,m}^{2}\right) }{ \partial x_{2,m}^{p_m}} \bigg|_{x_{2,m}= 0}  \nonumber \\[4mm]
& \hspace*{-5mm} \times \Bigg [ \prod_{n \in [m-1]}  \sum_{p_{n} \in \mathbb{N}} \frac{\widetilde{x}_{0,n}^{p_{n}}}{p_{n}!} \frac{\partial^{p_{n}} 
 \exp \left(  \frac{1}{\sigma^{2}} y_{n} (\sigma x_{2,n} - x_{0,n}) + \frac{1}{\sigma} x_{2,n} x_{0,n} - \frac{1}{2}  x_{2,n}^{2}\right) }{ \partial x_{2,n}^{p_n}} \bigg|_{x_{2,n}= 0} \Bigg]   \nonumber \\[4mm]Ê
 & \hspace*{-5mm} \times \Bigg [ \prod_{n' \in [S] \setminus [m]}  \sum_{p_{n'} \in \mathbb{Z}_{+}} \frac{\widetilde{x}_{0,n'}^{p_{n'}}}{p_{n'}!} \frac{\partial^{p_{n'}} 
 \exp \left(  \frac{1}{\sigma^{2}} y_{n'} (\sigma x_{2,n'} - x_{0,n'}) + \frac{1}{\sigma} x_{2,n'} x_{0,n'} - \frac{1}{2}  x_{2,n'}^{2}\right) }{ \partial x_{2,n'}^{p_{n'}}} \bigg|_{x_{2,n'}= 0} \Bigg] \nonumber \\[4mm] 
 &  \hspace*{-5mm} \times \frac{\widetilde{x}_{0,k}^{p_{k}}}{p_{k}!} \frac{\partial^{p_{k}} 
 \exp \left(  \frac{1}{\sigma^{2}} y_{k} (\sigma x_{2,k} - x_{0,k}) + \frac{1}{\sigma} x_{2,k} x_{0,k} - \frac{1}{2}  x_{2,k}^{2}\right) }{ \partial x_{2,k}^{p_k}} \bigg|_{x_{2,k}= 0}  \nonumber \\[4mm]Ê
  & \hspace*{-10mm} =  \exp \left( -\frac{1}{2 \sigma^{2}} \| \mathbf{x}_{0} \|^{2}_{2} \right) \sum_{m \in [S]} \exp \left( -  \frac{1}{ \sigma^{2}} y_{m} x_{0,m}  \right)   \nonumber \\[4mm]
& \hspace*{-5mm} \times \Bigg [ \prod_{n \in [m-1]}  \sum_{p_{n} \in \mathbb{N}} \frac{\widetilde{x}_{0,n}^{p_{n}}}{p_{n}!} \frac{\partial^{p_{n}} 
 \exp \left(  \frac{1}{\sigma^{2}} y_{n} (\sigma x_{2,n} - x_{0,n}) + \frac{1}{\sigma} x_{2,n} x_{0,n} - \frac{1}{2}  x_{2,n}^{2}\right) }{ \partial x_{2,n}^{p_n}} \bigg|_{x_{2,n}= 0} \Bigg]   \nonumber \\[4mm]Ê
 & \hspace*{-5mm} \times \Bigg[ \prod_{n' \in [S] \setminus [m]}  \sum_{p_{n'} \in \mathbb{Z}_{+}} \frac{\widetilde{x}_{0,n'}^{p_{n'}}}{p_{n'}!} \frac{\partial^{p_{n'}} 
 \exp \left(  \frac{1}{\sigma^{2}} y_{n'} (\sigma x_{2,n'} - x_{0,n'}) + \frac{1}{\sigma} x_{2,n'} x_{0,n'} - \frac{1}{2}  x_{2,n'}^{2}\right) }{ \partial x_{2,n'}^{p_{n'}}} \bigg|_{x_{2,n'}= 0} \Bigg]  \nonumber \\[4mm] 
 &  \hspace*{-5mm} \times \frac{1}{l! \sigma^{l}} \frac{\partial^{l} 
 \exp \left(  \frac{1}{\sigma} y_{k} x_{2,k} - \frac{1}{2}  x_{2,k}^{2}\right) }{ \partial x_{2,k}^{l}} \bigg|_{x_{2,k}= 0}. 
 \end{align} 
 Using the Taylor expansion of the exponential function \cite{RudinBook}, we obtain further
 \begin{align} 
 \label{equ_expr_case_2_main_ssnm_lmv_atom_congruence_part_2}
 \mathsf{J}[ r^{(l)}(\cdot)](\mathbf{y})   
 &= \exp \left( -\frac{1}{2 \sigma^{2}} \| \mathbf{x}_{0} \|^{2}_{2} \right) \sum_{m \in [S]} \exp \left( -  \frac{1}{ \sigma^{2}} y_{m} x_{0,m}  \right) \nonumber \\[4mm]
& \hspace*{-5mm}  \times \prod_{n \in [m-1]} \Bigg[\exp \left( \frac{1}{2 \sigma^{2}} {x}_{0,n}^{2} \right) - \exp \left( -  \frac{1}{ \sigma^{2}} y_{n} x_{0,n}  \right) \Bigg]  
 \prod_{n' \in [S] \setminus [m]} \exp \left( \frac{1}{2 \sigma^{2}} {x}_{0,n'}^{2} \right)  \nonumber \\[4mm] 
 &  \hspace*{-5mm} \times  \frac{1}{l! \sigma^{l}} \frac{\partial^{l} 
 \exp \left(  \frac{1}{\sigma} y_{k} x_{2,k} - \frac{1}{2}  x_{2,k}^{2}\right) }{ \partial x_{2,k}^{l}} \bigg|_{x_{2,k}= 0}. 
 \end{align} 
Using the identity 
\begin{align}
 \frac{\partial^{l} \exp \left( \frac{1}{\sigma} y x - \frac{1}{2} x^{2} \right) }{\partial x^{l}} \bigg|_{x=0}  & =  
\frac{\partial^{l} \exp \left( - \frac{1}{2} \left( x - \frac{y}{\sigma} \right)^{2} + \frac{1}{2} \frac{y^{2}}{\sigma^{2}} \right) }{\partial x^{l}} \bigg|_{x=0} \nonumber \\[4mm]Ê
   \Big( x' & \triangleq \frac{y}{\sigma} -x  \Big) \nonumber \\[4mm] 
& = \exp \left( \frac{1}{2} \frac{y^2}{\sigma^{2}} \right)  \frac{\partial^{l} \exp \left(  - \frac{1}{2} x'^{2} \right) }{\partial x'^{l}} \bigg|_{x'=-\frac{y}{\sigma}} =  \mathsf{H}_{l} \left(  \frac{y}{\sigma} \right), 
\end{align}
we finally have that
\begin{align}
 \label{equ_expr_case_2_main_ssnm_lmv_atom_congruence}
 \mathsf{J}[ r^{(l)}(\cdot)](\mathbf{y})   
 &Ê\stackrel{\eqref{equ_expr_case_2_main_ssnm_lmv_atom_congruence_part_2}}{=}  \exp \left( -\frac{1}{2 \sigma^{2}} \| \mathbf{x}_{0} \|^{2}_{2} \right) \sum_{m \in [S]} \exp \left( -  \frac{1}{ \sigma^{2}} y_{m} x_{0,m}  \right)   \nonumber \\[4mm]
& \hspace*{-5mm} \times \prod_{n \in [m-1]} \Bigg[\exp \left( \frac{1}{2 \sigma^{2}} {x}_{0,n}^{2} \right) - \exp \left( -  \frac{1}{ \sigma^{2}} y_{n} x_{0,n}  \right) \Bigg]  
 \prod_{n' \in [S] \setminus [m]} \exp \left( \frac{1}{2 \sigma^{2}} {x}_{0,n'}^{2} \right)  \nonumber \\[4mm] 
&  \hspace*{-5mm} \times \frac{1}{l! \sigma^{l}}  \mathsf{H}_{l} \left( \frac{y_{k}}{\sigma} \right) \nonumber \\[4mm] 
 & \hspace*{-10mm} =  \frac{1}{l! \sigma^{l}}  \mathsf{H}_{l} \left( \frac{y_{k}}{\sigma} \right)\sum_{m \in [S]} \exp \left( -  \frac{1}{ 2 \sigma^{2}} (x_{0,m}^{2} + 2y_{m} x_{0,m})  \right) \prod_{n \in [m-1]} \Bigg[ 1- \exp \left( -  \frac{1}{ 2\sigma^{2}} (x_{0,n}^{2} + 2y_{n} x_{0,n} )  \right) \Bigg]  ,
\end{align}}

By Theorem \ref{thm_main_facts_RKHS_MVE}, we have that if the prescribed bias function $c(\cdot)$ is valid, i.e., $\gamma(\cdot) \in  \mathcal{H}(\mathcal{M}_{\text{SSNM}})$ and in turn $\tilde{\gamma}(\cdot) \in \mathcal{H}(R_{e})$, the corresponding 
LMV estimator is given by 
\begin{align} 
\hat{x}_{k}^{(\mathbf{x}_{0})}(\mathbf{y}) & \stackrel{\eqref{equ_lmv_estimator_general_congruence_L_M_RKHS}}{=} \mathsf{J}[ \gamma(\cdot)] = \mathsf{J}[ \mathsf{K}_{e}^{-1}[\tilde{\gamma}(\cdot)]]  \stackrel{\eqref{equ_Cauchy_sequ_converges_image_mean_valid_bias_case_2_ssnm}}{=} \mathsf{J}\bigg[ \mathsf{K}_{e}^{-1}\bigg[ \lim_{N \rightarrow \infty} w_{N}(\cdot)\bigg]\bigg]   
 \stackrel{\eqref{equ_partial_sum_converges_image_mean_valid_bias_case_2_ssnm}}{=}  \mathsf{J}\bigg[ \mathsf{K}_{e}^{-1}\bigg[ \lim_{N \rightarrow \infty} \sum_{l \in [N]\cup \{0\} } m_{l} \sigma^{2l} v^{(l)}(\cdot)  \bigg]\bigg]  \nonumber \\[4mm]
 & =  \lim_{N \rightarrow \infty} \sum_{l \in[N]\cup \{0\} } m_{l} \sigma^{2l}  \mathsf{J}\big[ \mathsf{K}_{e}^{-1}\big[v^{(l)}(\cdot) \big]\big]  
 \stackrel{\eqref{equ_appendix_a_def_r_l}}{=}  \lim_{N \rightarrow \infty} \sum_{l \in [N]\cup \{0\} } m_{l} \sigma^{2l} \mathsf{J}\big[ r^{(l)}(\cdot)  \big] \nonumber  \\[4mm]Ê
& \stackrel{\eqref{equ_expr_case_2_main_ssnm_lmv_atom_congruence}}{=} \sum_{l \in \mathbb{Z}_{+}} \frac{m_{l}}{l! }\sigma^{l}  \mathsf{H}_{l} \left(\frac{y_{k}}{\sigma} \right)\sum_{m \in [S]} \exp \left( -  \frac{1}{ 2 \sigma^{2}} (x_{0,m}^{2} + 2y_{m} x_{0,m})  \right)Ê \nonumber \\[4mm]
& \hspace*{5mm} \times \prod_{n \in [m-1]} \Big[ 1- \exp \big( -  \frac{1}{ 2\sigma^{2}} (x_{0,n}^{2} + 2y_{n} x_{0,n} )  \big) \Big]  ,
\end{align} 
which is equal to \eqref{equ_lmv_case_2_ssnm_main}. 

It remains to verify the necessity of the condition \eqref{equ_cond_valid_bias_diagon_bias_SSNM} for a bias function $c(\cdot)$ to be valid. To that end, assume that \eqref{equ_cond_valid_bias_diagon_bias_SSNM} is not fulfilled, 
i.e., the sums $\sum_{l \in \mathcal{T}}  \frac{m_{l}^{2}}{l!} \sigma^{2l}$ can be made arbitrary large by using 
a suitable finite index set $\mathcal{T} \subseteq \mathbb{Z}_{+}$, but the bias function $c(\cdot)$ is valid, i.e., $\gamma(\cdot) \in  \mathcal{H}(\mathcal{M}_{\text{SSNM}})$. However, 
$\gamma(\cdot) \in  \mathcal{H}(\mathcal{M}_{\text{SSNM}})$ implies, due to the congruence $\mathsf{K}_{e}[\cdot]$ in \eqref{equ_def_congruence_SSNM} of Theorem \ref{thm_entire_charact_RKHS_SSNM}, that we also have 
$\tilde{\gamma}(\cdot) =  \gamma(\sigma \mathbf{x}) \nu_{\mathbf{x}_{0}}(\mathbf{x})= \mathsf{K}_{e}[\gamma(\cdot)] \in \mathcal{H}(R_{e})$. 
By Theorem \ref{thm_orthog_proj_ineq}, we then have 
for any subspace $\mathcal{U} \triangleq \linspan\{ v^{(l)}(\cdot) \}_{l \in \mathcal{T}}$ that\footnote{An ONB for $\mathcal{U}$ is given by 
$\big\{  v^{(l)}(\cdot) / \| v^{(l)}(\cdot) \|_{\mathcal{H}(R_{e})} \big\}_{l \in \mathcal{T}}$ since 
the functions $v^{(l)}(\cdot)$ are orthogonal.} 
\begin{align}
\label{equ_proof_SSNM_sum_cond_valid_case_2}
 \| \gamma(\cdot) \|^{2}_{\mathcal{H}(\mathcal{M}_{\text{SSNM}})} & = \|\mathsf{K}_{e}[\gamma(\cdot)]  \|^{2}_{\mathcal{H}(R_{e})} = \| \tilde{\gamma}(\cdot) \|^{2}_{\mathcal{H}(R_{e})}  \stackrel{\eqref{equ_squared_norm_orthog_proj_pythag_thm}}{\geq}  \| \mathbf{P}_{\mathcal{U}} \tilde{\gamma}(\cdot) \|^{2}_{\mathcal{H}(R_{e})}  \nonumber \\[4mm]
 & \stackrel{\eqref{equ_orth_proj_finite_dim_onb_squared_norm}}{=}  \sum_{l \in \mathcal{T}} \bigg[ \big\langle \tilde{\gamma}(\cdot), v^{(l)}(\cdot) / \| v^{(l)}(\cdot) \|_{\mathcal{H}(R_{e})} \big\rangle_{\mathcal{H}(R_{e})} \bigg]^{2} \nonumber \\[4mm]
  & \stackrel{\eqref{equ_Cauchy_sequ_converges_image_mean_valid_bias_case_2_ssnm}}{=}  \sum_{l \in \mathcal{T}} \bigg[ \big\langle \lim_{N \rightarrow \infty} w_{N}(\cdot), v^{(l)}(\cdot) / \| v^{(l)}(\cdot) \|_{\mathcal{H}(R_{e})} \big\rangle_{\mathcal{H}(R_{e})} \bigg]^{2}\nonumber \\[4mm]
  & \stackrel{(a)}{=}  \sum_{l \in \mathcal{T}} \bigg[  \lim_{N \rightarrow \infty} \big\langle w_{N}(\cdot), v^{(l)}(\cdot) / \| v^{(l)}(\cdot) \|_{\mathcal{H}(R_{e})} \big\rangle_{\mathcal{H}(R_{e})} \bigg]^{2}Ê\nonumber \\[4mm]
  &  \stackrel{\eqref{equ_partial_sum_converges_image_mean_valid_bias_case_2_ssnm}}{=}  \sum_{l \in \mathcal{T}} \bigg[  \frac{m_{l}\sigma^{2l}} {  \| v^{(l)}(\cdot) \|_{\mathcal{H}(R_{e})}}  \big\langle v^{(l)}(\cdot), v^{(l)}(\cdot)  \big\rangle_{\mathcal{H}(R_{e})} \bigg]^{2} \nonumber \\[4mm]
  &  =  \sum_{l \in \mathcal{T}} \bigg[ m_{l}\sigma^{2l}  \| v^{(l)}(\cdot) \|_{\mathcal{H}(R_{e})} \bigg]^{2} \nonumber \\[4mm]
 &    \stackrel{\eqref{equ_appendix_a_relation_squared_norm_v_l_v_0}}{=}  \sum_{l \in \mathcal{T}} \bigg[  \frac{m_{l} \sigma^{l}}{\sqrt{l!} } \| v^{(0)}(\cdot) \|_{\mathcal{H}(R_{e})}  \bigg]^{2} 
  =  \sum_{l \in \mathcal{T}}  \frac{m_{l}^{2}}{l!} \sigma^{2l} \| v^{(0)}(\cdot) \|_{\mathcal{H}(R_{e})}^{2},
\end{align}
where in $(a)$ the change of the order of taking the limit $\lim_{N \rightarrow \infty}$ and computing the inner product with $v^{(l)}(\cdot) / \| v^{(l)}(\cdot) \|_{\mathcal{H}(R_{e})}$, respectively, is due to the 
fact that the inner product is a continuous mapping (see, e.g., \cite[Theorem 4.6]{RudinBook}). 
However, since \eqref{equ_cond_valid_bias_diagon_bias_SSNM} is not satisfied, it follows from \eqref{equ_proof_SSNM_sum_cond_valid_case_2} 
that the squared norm $ \| \gamma(\cdot) \|^{2}_{\mathcal{H}(\mathcal{M}_{\text{SSNM}})}$ is unbounded, which is impossible because $\gamma(\cdot) \in \mathcal{H}(\mathcal{M}_{\text{SSNM}})$. 
\end{proof}

\chapter[]{Proof of Theorem \ref{thm_sparse_lower_bound_SPCM}}
\label{chap_appendix_D}

\begin{proof}

The statement is derived by considering the minimum variance problem $\mathcal{M}_{\mathcal{D},\text{SPCM}}=\mathcal{M}_{\text{SPCM}}\big|_{\mathcal{D}}$, 
where the set $\mathcal{D}$ is chosen as 
\begin{equation} 
\label{equ_def_set_D_proof_RIP_lower_bound_SPCM}
\mathcal{D} = \big( \mathcal{B}(\mathbf{0},r) \cup \mathcal{B}(\mathbf{x}_{0},r)  \big) \cap \big \{ \mathbf{x} \in \mathcal{X}_{S,+} \big| \supp(\mathbf{x}) \subseteq \mathcal{K} \big\},
\end{equation}
where 
\begin{equation} 
\label{equ_proof_RIP_bound_SPCM_index_set_K}
\mathcal{K} \triangleq \big \{ \{ l \} \cup \supp(\mathbf{x}_{0}) \big\} =\{i_{1},\ldots,i_{K} \}
\end{equation} 
 with $K \leq S+1$ and an arbitrary but (for the rest of the proof) fixed index $l \in [N] \setminus \supp(\mathbf{x}_{0})$. 
 Furthermore, the radius $r$ is sufficiently small 
such that Theorem \ref{thm_cond_d_restr_SPCM_ball} applies for the two sets $\mathcal{B}(\mathbf{0},r) \cap \mathcal{X}_{S,+}$ and $\mathcal{B}(\mathbf{x}_{0},r) \cap \mathcal{X}_{S,+}$ 
separately. Therefore, the condition \eqref{equ_nec_suff_cond_D_SPCM_kernel_exists} is also satisfied for their union $\big( \mathcal{B}(\mathbf{0},r) \cup \mathcal{B}(\mathbf{x}_{0},r)  \big) \cap \mathcal{X}_{S,+}$ and in turn also for the set $\mathcal{D}$ given in \eqref{equ_def_set_D_proof_RIP_lower_bound_SPCM}, since 
this set is contained in the union. Thus, the set $\mathcal{D}$ given in \eqref{equ_def_set_D_proof_RIP_lower_bound_SPCM} satisfies \eqref{equ_nec_suff_cond_D_SPCM_kernel_exists} and therefore 
the RKHS $\mathcal{H}(\mathcal{M}_{\mathcal{D},\text{SPCM}})$ exists. 
Since we can assume that $g(\cdot)$ is estimable for $\mathcal{M}_{\text{SDPCM}}$, which implies via Theorem \ref{thm_par_set_reduction_classic_est_mve} that 
the restriction $g(\cdot)\big|_{\mathcal{D}}$ is estimable for $\mathcal{M}_{\mathcal{D}, \text{SDPCM}}$, we have that the prescribed mean function $\gamma(\cdot): \mathcal{D} \rightarrow \mathbb{R}: \gamma(\mathbf{x}) = c(\mathbf{x}) + g(\mathbf{x})= g(\mathbf{x})$ belongs to the RKHS $\mathcal{H}(\mathcal{M}_{\mathcal{D},\text{SPCM}})$. 

The outline of the rest of the proof is as follows: We will construct a subspace 
\begin{equation}
\mathcal{U} \triangleq \linspan \{ w_{0}(\cdot), w_{1}(\cdot) \} \subseteq \mathcal{H}(\mathcal{M}_{\mathcal{D},\text{SPCM}})
\end{equation}
that is spanned by two orthogonal functions $w_{0}(\cdot)$ and $w_{1}(\cdot)$, respectively. We can then invoke Theorem \ref{thm_main_facts_RKHS_MVE} and Theorem \ref{thm_orthog_proj_ineq} to lower bound the 
minimum achievable variance $L_{\mathcal{M}_{\mathcal{D},\text{SPCM}}}$ by projecting the mean function $\gamma(\cdot)$ on the subspace $\mathcal{U}$. However, we will not compute the squared norm of the projection $\mathbf{P}_{\mathcal{U}} \gamma(\cdot)$ explicitly, but lower bound it by deriving an upper bound on the squared norm of the function $w_{1}(\cdot)$. 
 The so obtained lower bound on  $L_{\mathcal{M}_{\mathcal{D},\text{SPCM}}}$ yields then via \eqref{equ_lower_bound_spcm_restr_spcm} a lower bound on  $L_{\mathcal{M}_{\text{SPCM}}}$.

Consider the matrix $\mathbf{H}^{\mathcal{K}} \triangleq \big( \mathbf{H}_{i_{1}},\cdots, \mathbf{H}_{i_{K}} \big) \in \mathbb{R}^{M \times  Q}$, 
where $Q \triangleq \sum_{k \in \mathcal{K}} r_{k}$, and its thin SVD 
\begin{equation} 
\label{equ_svd_H_K_lower_bound_not_supp_SPCM}
\mathbf{H}^{\mathcal{K}} = \mathbf{U}_{1} \mathbf{\Sigma}_{1} \mathbf{V}_{1}^{T}.
\end{equation} 
Due to the assumption that the basis matrices $\big\{ \mathbf{C}_{k} \big\}_{k \in [N]}$ satisfy the RIP of order $S+1$, the diagonal entries of the diagonal matrix $\mathbf{\Sigma}_{1} \in \mathbb{R}^{ Q \times Q}$, i.e., 
the singular values of $\mathbf{H}^{\mathcal{K}}$, belong to the interval $[1-\delta_{S+1},1+\delta_{S+1}]$. 
This also implies that the matrix $\mathbf{H}^{\mathcal{K}}$ has full column rank, i.e., 
$\rank(\mathbf{H}^{\mathcal{K}}) =  Q$, and that the matrix 
$\mathbf{\Sigma}_{1}\mathbf{V}_{1}^{T} \in \mathbb{R}^{ Q \times Q}$ appearing in the thin SVD \eqref{equ_svd_H_K_lower_bound_not_supp_SPCM} is invertible. 
For every $\mathbf{x} \in \mathcal{D}$, we have due to \eqref{equ_svd_H_K_lower_bound_not_supp_SPCM} that 
\begin{align} 
\label{equ_subspace_decomp_matrix_proof_D_K}
\widetilde{\mathbf{C}}(\mathbf{x}) & \stackrel{\eqref{equ_def_cov_matrix_observation_SPCM},\eqref{equ_cov_param_SPCM_1}}{=} \mathbf{H} \mathbf{D}^{\mathcal{K}}(\mathbf{x}) \mathbf{H}^{T} + \sigma^{2} \mathbf{I} = \sum_{k \in \mathcal{K}} \mathbf{H}_{k} x_{k} \mathbf{H}_{k}^{T} + \sigma^{2} \mathbf{I}  \nonumber \\[4mm]
&  = \mathbf{H}^{\mathcal{K}}\mathbf{D}^{\mathcal{K}}(\mathbf{x})\big(\mathbf{H}^{\mathcal{K}}\big)^{T} + \sigma^{2} \mathbf{I}  \nonumber \\[4mm]
&  =\mathbf{U}_{1} \mathbf{\Sigma}_{1} \mathbf{V}_{1}^{T} \mathbf{D}^{\mathcal{K}}(\mathbf{x})\mathbf{V}_{1} \mathbf{\Sigma}_{1} \mathbf{U}_{1}^{T} + \sigma^{2} \mathbf{P}^{\mathcal{K}}+ \sigma^{2}(\mathbf{I} - \mathbf{P}^{\mathcal{K}})   \nonumber \\[4mm]
&  =\mathbf{U}_{1} \mathbf{\Sigma}_{1} \mathbf{V}_{1}^{T} \mathbf{D}^{\mathcal{K}}(\mathbf{x})\mathbf{V}_{1} \mathbf{\Sigma}_{1} \mathbf{U}_{1}^{T} + \sigma^{2} \mathbf{U}_{1} \mathbf{U}_{1}^{T} + \sigma^{2}(\mathbf{I} - \mathbf{P}^{\mathcal{K}})   \nonumber \\[4mm]
& = \mathbf{U}_{1} \mathbf{\Sigma}_{1} \mathbf{V}_{1}^{T}  \mathbf{C}^{\mathcal{K}}(\mathbf{x}) \mathbf{V}_{1} \mathbf{\Sigma}_{1} \mathbf{U}_{1}^{T} + \big(\mathbf{I} - \mathbf{P}^{\mathcal{K}}\big) \sigma^{2}.
\end{align}
Here, we used
\begin{equation}
\label{equ_def_RIP_SPCM_bound_D_K}
\mathbf{D}^{\mathcal{K}}(\mathbf{x}) \triangleq \begin{pmatrix} x_{i_{1}} \mathbf{I}_{r_{i_{1}}} & \mathbf{0} & \ldots & \mathbf{0} \\ 
\mathbf{0} & x_{i_{2}} \mathbf{I}_{r_{i_{2}}} & \ddots  &  \vdots \\ 
\vdots & \ddots& \ddots &   \mathbf{0} \\ 
\mathbf{0} & \ldots&  \mathbf{0} &  x_{i_{K}} \mathbf{I}_{r_{i_{K}}} \\ 
\end{pmatrix},
\end{equation} 
\vspace*{3mm}
\begin{equation} 
\label{equ_def_transformed_cov_matrix_model_SPCM}
\mathbf{C}^{\mathcal{K}}(\mathbf{x}) \triangleq \mathbf{D}^{\mathcal{K}}(\mathbf{x}) + \sigma^{2} \mathbf{V}_{1} \mathbf{\Sigma}_{1}^{-2} \mathbf{V}_{1}^{T} =\mathbf{D}^{\mathcal{K}}(\mathbf{x})  + \sigma^{2} \mathbf{I} + \mathbf{E}_{1},
\end{equation} 
\vspace*{-1mm}
\begin{equation} 
\mathbf{E}_{1} \triangleq \sigma^{2} \big( \mathbf{V}_{1} \mathbf{\Sigma}_{1}^{-2} \mathbf{V}_{1}^{T} - \mathbf{I} \big),
\end{equation}
and 
\begin{equation} 
\mathbf{P}^{\mathcal{K}} \triangleq \mathbf{H}^{\mathcal{K}} \big(\mathbf{H}^{\mathcal{K}}\big)^{\dagger} = \mathbf{U}_{1} \mathbf{\Sigma}_{1} \mathbf{V}_{1}^{T}\mathbf{V}_{1}\mathbf{\Sigma}^{-1}_{1} \mathbf{U}^{T}_{1}
=  \mathbf{U}_{1} \mathbf{U}^{T}_{1} 
\end{equation} 
denotes the orthogonal projection matrix on the subspace 
$\linspan(\mathbf{H}^{\mathcal{K}}) \subseteq \mathbb{R}^{M}$. 
Note that 
\begin{equation} 
(\mathbf{I} - \mathbf{P}^{\mathcal{K}}) \mathbf{U}_{1} =  \mathbf{U}_{1} -  \mathbf{U}_{1} \mathbf{U}^{T}_{1}  \mathbf{U}_{1} =  \mathbf{U}_{1} - \mathbf{U}_{1}  = \mathbf{0}.
\end{equation}  
Since, by \eqref{equ_subspace_decomp_matrix_proof_D_K}, the matrix $\widetilde{\mathbf{C}}(\mathbf{x})$ is the sum of two matrices 
which act (viewed as a linear transformation) independently on the two complementary subsets 
$\linspan(\mathbf{U}_{1}), \linspan \big(\mathbf{I}-\mathbf{P}^{\mathcal{K}} \big) \subseteq \mathbb{R}^{M}$ 
we have by \cite{golub96,HalmosFiniteDimVecSpace} that
\begin{align} 
\label{equ_det_tilde_observation_space_support_space_SPCM}
\detmb{\widetilde{\mathbf{C}}(\mathbf{x})} & = \detmb{ \mathbf{\Sigma}_{1} \mathbf{V}_{1}^{T}  \mathbf{C}^{\mathcal{K}}(\mathbf{x}) \mathbf{V}_{1} \mathbf{\Sigma}_{1}}\sigma^{2 (M -Q)}.
\end{align} 

The $2$-norm of the matrix $\mathbf{E}_{1}$ can be bounded as 
\begin{equation} 
\label{equ_proof_lower_bound_RIP_SPCM_matrix_pertub}
\| \mathbf{E}_{1} \|_{2} \stackrel{(a)}{\leq} \sigma^{2} \frac{4 \delta_{S+1}}{1-2\delta_{S+1}} \stackrel{(b)}{\leq}  \sigma^{2} 5 \delta_{S+1}.
\end{equation} 
Here, step $(b)$ follows from $\delta_{S+1} < 1/32$. 
In order to verify step $(a)$ in \eqref{equ_proof_lower_bound_RIP_SPCM_matrix_pertub}, note that for an arbitrary vector $\mathbf{x} \in \mathbb{R}^{N}$ with $\| \mathbf{x} \|_{2}=1$, 
we obtain (note that the diagonal elements of $\mathbf{\Sigma}_{1}^{-2}$ belong to the interval 
$\big[(1+\delta_{S+1})^{-2}, (1-\delta_{S+1})^{-2} \big]$)
\begin{equation} 
\label{equ_proof_RIP_bound_SPCM_quadratic_form_E_1}
\mathbf{x}^{T} \mathbf{E}_{1} \mathbf{x} =  \sigma^{2} \big(Ê \underbrace{\mathbf{x}^{T} \mathbf{V}_{1} \mathbf{\Sigma}_{1}^{-2} \mathbf{V}_{1}^{T} \mathbf{x}}_{ \in \big[(1+\delta_{S+1})^{-2}, (1-\delta_{S+1})^{-2} \big]}- \underbrace{\mathbf{x}^{T}\mathbf{x}}_{=1} \big)  \in  \sigma^{2}  \big[(1+\delta_{S+1})^{-2}-1, (1-\delta_{S+1})^{-2}-1 \big].
\end{equation}
From \eqref{equ_proof_RIP_bound_SPCM_quadratic_form_E_1} it follows via the EVD of the symmetric matrix $\mathbf{E}_{1}$ \cite{golub96} that the eigenvalues $\lambda_{k} (\mathbf{E}_{1})$ belong to the interval $\sigma^{2} \big[(1+\delta_{S+1})^{-2}-1, (1-\delta_{S+1})^{-2}-1 \big]$, i.e., 
\begin{equation} 
\label{equ_proof_RIP_bound_SPCM_eigvals_E_1_belong_interval}
\lambda_{k} (\mathbf{E}_{1}) \in  \sigma^{2}  \big[(1+\delta_{S+1})^{-2}-1, (1-\delta_{S+1})^{-2}-1 \big].
\end{equation} 
This further implies that
\begin{align}
\| \mathbf{E}_{1} \|_{2} & = \max\big\{ \big|\lambda_{Q} (\mathbf{E}_{1})\big|, \big|\lambda_{1}(\mathbf{E}_{1})\big| \big\} \leq  \big|\lambda_{Q} (\mathbf{E}_{1})\big|+ \big|\lambda_{1}(\mathbf{E}_{1})\big| \nonumber \\[4mm]
& \stackrel{\eqref{equ_proof_RIP_bound_SPCM_eigvals_E_1_belong_interval}}\leq \sigma^{2} \big[(1-\delta_{S+1})^{-2} - 1\big]- \sigma^{2}\big[ (1+\delta_{S+1})^{-2}-1\big]  \nonumber \\[4mm]
& = \sigma^{2} \big[ (1-\delta_{S+1})^{-2} - (1+\delta_{S+1})^{-2}\big] = \sigma^{2} \frac{(1+\delta_{S+1})^{2}-(1-\delta_{S+1})^{2}}{(1-\delta_{S+1})^{2}(1+\delta_{S+1})^{2}} \nonumber \\[4mm]
& = \sigma^{2}  \frac{4\delta_{S+1}}{(1-\delta_{S+1})^{2}(1+\delta_{S+1})^{2}} 
 \leq \sigma^{2}  \frac{4\delta_{S+1}}{(1-\delta_{S+1})^{2}} =  \sigma^{2}  \frac{4\delta_{S+1}}{1-2\delta_{S+1}+\delta_{S+1}^{2}}  \nonumber \\[4mm]
 & \leq  \sigma^{2}  \frac{4\delta_{S+1}}{1-2\delta_{S+1}} ,
\end{align}
where step $(a)$ follows from the fact that $\mathbf{E}_{1}$ is symmetric, i.e., $\mathbf{E}^{T}_{1} = \mathbf{E}_{1}$ (cf.\ \cite{golub96}).
It is important to note that the matrix $\mathbf{E}_{1}$ appearing in \eqref{equ_def_transformed_cov_matrix_model_SPCM} does not depend on $\mathbf{x}$.  

From \eqref{equ_subspace_decomp_matrix_proof_D_K}, we have via elementary linear algebra \cite{golub96} that the kernel $R_{\mathcal{M}_{\mathcal{D},\text{SPCM}}}(\cdot,\cdot)$ in \eqref{equ_kernel_D_restr_SPCM} can be rewritten as 
\vspace*{3mm}
\begin{align} 
R_{\mathcal{M}_{\mathcal{D},\text{SPCM}}}(\mathbf{x}_{1}, \mathbf{x}_{2}) & = \big[ \detmb{ \widetilde{\mathbf{C}}(\mathbf{x}_{0}) } \big]^{1/2} \ist 
  \big[ \detmb{ \widetilde{\mathbf{C}}(\mathbf{x}_{1}) + \widetilde{\mathbf{C}}(\mathbf{x}_{2}) - \widetilde{\mathbf{C}}(\mathbf{x}_{1}) 
  \widetilde{\mathbf{C}}^{-1}(\mathbf{x}_{0})\widetilde{\mathbf{C}}(\mathbf{x}_{2}) } \big]^{-1/2}   \nonumber \\[4mm]
  & \hspace*{-20mm} \stackrel{\eqref{equ_subspace_decomp_matrix_proof_D_K},\eqref{equ_det_tilde_observation_space_support_space_SPCM}}{=} \big[ \detmb{ \mathbf{\Sigma}_{1} \mathbf{V}_{1}^{T} \mathbf{C}^{\mathcal{K}}(\mathbf{x}_{0}) \mathbf{V}_{1} \mathbf{\Sigma}_{1}  } \big]^{1/2} \ist \times \nonumber \\[1mm]
  & \big[ \detmb{ \mathbf{\Sigma}_{1} \mathbf{V}_{1}^{T} \big(  \mathbf{C}^{\mathcal{K}}(\mathbf{x}_{1}) + \mathbf{C}^{\mathcal{K}}(\mathbf{x}_{2}) - \mathbf{C}^{\mathcal{K}}(\mathbf{x}_{1}) 
 \big[ \mathbf{C}^{\mathcal{K}}(\mathbf{x}_{0})\big]^{-1} \mathbf{C}^{\mathcal{K}}(\mathbf{x}_{2})\big) \mathbf{V}_{1} \mathbf{\Sigma}_{1} } \big]^{-1/2}   \nonumber \\[4mm]
  & \hspace*{-20mm} \stackrel{(a)}{=} \big[ \detmb{ \mathbf{C}^{\mathcal{K}}(\mathbf{x}_{0}) } \big]^{1/2} \ist 
  \big[ \detmb{  \mathbf{C}^{\mathcal{K}}(\mathbf{x}_{1}) + \mathbf{C}^{\mathcal{K}}(\mathbf{x}_{2}) - \mathbf{C}^{\mathcal{K}}(\mathbf{x}_{1}) 
 \big[ \mathbf{C}^{\mathcal{K}}(\mathbf{x}_{0})\big]^{-1} \mathbf{C}^{\mathcal{K}}(\mathbf{x}_{2}) } \big]^{-1/2},
  \label{equ_rewrite_kernel_SPCM_RIP_bound}
\end{align}
where $(a)$ is due to the identity $\detmb{ \mathbf{M} \mathbf{N} } = \detmb{\mathbf{M}} \detmb{\mathbf{N}}$ for two square matrices $\mathbf{M}, \mathbf{N}$ \cite{HalmosFiniteDimVecSpace}.

Let us now define the subspace $\mathcal{U} \triangleq \linspan \{ w_{0}(\cdot), w_{1}(\cdot) \} \subseteq \mathcal{H}(\mathcal{M}_{\mathcal{D},\text{SPCM}})$ spanned by 
the two functions $w_{0}(\cdot) \triangleq R_{\mathcal{M}_{\mathcal{D},\text{SPCM}}}(\cdot,\mathbf{x}_{0}) \in  \mathcal{H}(\mathcal{M}_{\mathcal{D},\text{SPCM}})$ and 
\begin{equation} 
w_{1}(\cdot) \triangleq \frac{ \partial^{\mathbf{e}_{l}} R_{\mathcal{M}_{\mathcal{D},\text{SPCM}}} (\cdot, \mathbf{x}_{2})}  {\partial \mathbf{x}_{2}^{\mathbf{e}_{l}}} \bigg|_{\mathbf{x}_{2} = \mathbf{0}},
\end{equation}  
where $l \in [N] \setminus \supp(\mathbf{x}_{0})$ is the specific index which is used in the definition of the set $\mathcal{K}$ in \eqref{equ_proof_RIP_bound_SPCM_index_set_K}. 
By Theorem \ref{thm_der_repr_prop} (using the function $g^{(\mathbf{p})}_{\mathbf{x}_{c}}(\cdot)$ with $\mathbf{p}Ê\!=\! \mathbf{e}_{l}$ and $\mathbf{x}_{c} \!=\! \mathbf{0}$), 
we have that $w_{1} (\cdot)  \in \mathcal{H}(\mathcal{M}_{\mathcal{D},\text{SPCM}})$ 
and moreover 
\begin{align}
\label{equ_inner_prod_zero_eval_par_der_zero_RIP_bound}
\big\langle  w_{0}(\cdot), w_{1}(\cdot) \big\rangle_{\mathcal{H}(\mathcal{M}_{\mathcal{D},\text{SPCM}})} 
\stackrel{\eqref{equ_der_reproduction_prop}}{=}  \frac{ \partial^{\mathbf{e}_{l}} R_{\mathcal{M}_{\mathcal{D},\text{SPCM}}} (\mathbf{x}_{2},\mathbf{x}_{0})}  {\partial \mathbf{x}_{2}^{\mathbf{e}_{l}}} \bigg|_{\mathbf{x}_{2} = \mathbf{0}}
\stackrel{\eqref{equ_kernel_one_arg_x_0_equal_1}}{=}  \frac{ \partial^{\mathbf{e}_{l}}1}  {\partial \mathbf{x}_{2}^{\mathbf{e}_{l}}} \bigg|_{\mathbf{x}_{2} = \mathbf{0}}  = 0. 
\end{align}
Therefore, an ONB for $\mathcal{U}$ is given by 
$\big\{ \frac{w_{0}(\cdot)}{\|  w_{0}(\cdot) \|_{\mathcal{H}(\mathcal{M}_{\mathcal{D},\text{SPCM}})} }, \frac{w_{1}(\cdot)}{\|  w_{1}(\cdot) \|_{\mathcal{H}(\mathcal{M}_{\mathcal{D},\text{SPCM}})} } \big\}$. 

By projecting the prescribed mean function $\gamma(\mathbf{x})=g(\mathbf{x}) + c(\mathbf{x}) = g(\mathbf{x})$ of $\mathcal{M}_{\mathcal{D},\text{SPCM}}$ onto the subspace $\mathcal{U}$, we obtain 
\begin{align}
 L_{\mathcal{M}_{\text{SPCM}}}  & \stackrel{\eqref{equ_lower_bound_spcm_restr_spcm}}{\geq}  L_{\mathcal{M}_{\mathcal{D},\text{SPCM}}} \stackrel{\eqref{equ_squared_norm_min_achiev_var}}{=} 
 \| \gamma(\cdot) \|^{2}_{\mathcal{H}(\mathcal{M}_{\mathcal{D},\text{SPCM}})} - \big[ \gamma(\mathbf{x}_{0}) \big]^{2} \stackrel{\eqref{equ_lower_bound_RKHS_projection}}{\geq}  \| \mathbf{P}_{\mathcal{U}} \gamma(\cdot) \|^{2}_{\mathcal{H}(\mathcal{M}_{\mathcal{D},\text{SPCM}})} - \big[ \gamma(\mathbf{x}_{0}) \big]^{2}  \nonumber \\[4mm]
 & \stackrel{(a)}{=} 
\frac{ \langle \gamma(\cdot), w_{0}(\cdot) \rangle_{\mathcal{H}(\mathcal{M}_{\mathcal{D},\text{SPCM}})}^{2}}{ \langle w_{0}(\cdot), w_{0}(\cdot) \rangle_{\mathcal{H}(\mathcal{M}_{\mathcal{D},\text{SPCM}})}} +  
\frac{\langle \gamma(\cdot), w_{1}(\cdot) \rangle_{\mathcal{H}(\mathcal{M}_{\mathcal{D},\text{SPCM}})}^{2}}{\langle w_{1}(\cdot), w_{1}(\cdot) \rangle_{\mathcal{H}(\mathcal{M}_{\mathcal{D},\text{SPCM}})}}  - \big[ \gamma(\mathbf{x}_{0})  \big]^{2} \nonumber \\[4mm]
 & =
\frac{ \langle \gamma(\cdot),  R_{\mathcal{M}_{\mathcal{D},\text{SPCM}}}(\cdot,\mathbf{x}_{0}) \rangle_{\mathcal{H}(\mathcal{M}_{\mathcal{D},\text{SPCM}})}^{2}}{ \langle  R_{\mathcal{M}_{\mathcal{D},\text{SPCM}}}(\cdot,\mathbf{x}_{0}),  R_{\mathcal{M}_{\mathcal{D},\text{SPCM}}}(\cdot,\mathbf{x}_{0}) \rangle_{\mathcal{H}(\mathcal{M}_{\mathcal{D},\text{SPCM}})}} +  
\frac{\langle \gamma(\cdot), w_{1}(\cdot) \rangle_{\mathcal{H}(\mathcal{M}_{\mathcal{D},\text{SPCM}})}^{2}}{\langle w_{1}(\cdot), w_{1}(\cdot) \rangle_{\mathcal{H}(\mathcal{M}_{\mathcal{D},\text{SPCM}})}}  - \big[ \gamma(\mathbf{x}_{0})  \big]^{2} \nonumber \\[4mm]
 & \stackrel{\eqref{equ_reproduction_property},\eqref{equ_kernel_one_arg_x_0_equal_1}}{=}
\frac{ \big[\gamma(\mathbf{x}_{0}) \big]^{2}}{ 1} +  
\frac{\langle \gamma(\cdot), w_{1}(\cdot) \rangle_{\mathcal{H}(\mathcal{M}_{\mathcal{D},\text{SPCM}})}^{2}}{\langle w_{1}(\cdot), w_{1}(\cdot) \rangle_{\mathcal{H}(\mathcal{M}_{\mathcal{D},\text{SPCM}})}}  - \big[ \gamma(\mathbf{x}_{0}) \big]^{2} \nonumber \\[4mm]
 & \stackrel{\gamma(\cdot) =g(\cdot)}{=}
\frac{ \big[g(\mathbf{x}_{0}) \big]^{2}}{ 1} +  
\frac{\langle g(\cdot), w_{1}(\cdot) \rangle_{\mathcal{H}(\mathcal{M}_{\mathcal{D},\text{SPCM}})}^{2}}{\langle w_{1}(\cdot), w_{1}(\cdot) \rangle_{\mathcal{H}(\mathcal{M}_{\mathcal{D},\text{SPCM}})}}  - \big[ g(\mathbf{x}_{0})  \big]^{2} \nonumber \\[4mm]
& \stackrel{(b)}{=} b_{l}^{2}\|  w_{1}(\cdot) \|^{-2}_{\mathcal{H}(\mathcal{M}_{\mathcal{D},\text{SPCM}})},
\label{equ_proof_lower_bound_RIP_projection_SPCM}
\end{align}
where $b_{l} \triangleq \frac{\partial g(\mathbf{x})}{\partial x_{l}} \big|_{\mathbf{x} = \mathbf{0}}$; step $(a)$ follows from Theorem \ref{thm_norm_projection_finite_dim_subspace_union_orthog_subspaces}, and step $(b)$ is 
due to Theorem \ref{thm_der_repr_prop} (using the function $g^{(\mathbf{p})}_{\mathbf{x}_{c}}(\cdot)$ with $\mathbf{p}Ê\!=\! \mathbf{e}_{l}$ and $\mathbf{x}_{c} \!=\! \mathbf{0}$).

The last and central step of the proof derives an upper bound on the squared norm $\|  w_{1}(\cdot) \|^{2}_{\mathcal{H}(\mathcal{M}_{\mathcal{D},\text{SPCM}})}$, 
which combined with \eqref{equ_proof_lower_bound_RIP_projection_SPCM} yields 
the bound in \eqref{equ_lower_bound_SPCM_RIP}.  
By combining Theorem \ref{thm_der_repr_prop} (using the function $g^{(\mathbf{p})}_{\mathbf{x}_{c}}(\cdot)$ with $\mathbf{p}Ê\!=\! \mathbf{e}_{l}$ and $\mathbf{x}_{c} \!=\! \mathbf{0}$), Lemma \ref{lem_der_det_function}, and Lemma \ref{lem_derivation_inverse_matrix} we obtain
\begin{align} 
\label{equ_proof_lower_bond_RIP_SPCM_squared_norm_expr}
\|  w_{1}(\cdot) \|^{2}_{\mathcal{H}(\mathcal{M}_{\mathcal{D},\text{SPCM}})} & = \big\langle  w_{1}(\cdot), w_{1}(\cdot) \big\rangle_{\mathcal{H}(\mathcal{M}_{\mathcal{D},\text{SPCM}})}  
 \stackrel{\eqref{equ_der_reproduction_prop}}{=} \frac{\partial^{\mathbf{e}_{l}} \partial^{\mathbf{e}_{l}} R_{\mathcal{M}_{\mathcal{D},\text{SPCM}}} (\mathbf{x}_{1}, \mathbf{x}_{2})}  {\partial \mathbf{x}_{1}^{\mathbf{e}_{l}} \partial \mathbf{x}_{2}^{\mathbf{e}_{l}}} \bigg|_{\mathbf{x}_{1} = \mathbf{x}_{2} = \mathbf{0}} \nonumber \\[4mm]
& \hspace*{-30mm} \stackrel{\eqref{equ_rewrite_kernel_SPCM_RIP_bound}}{=} \big[ \detmb{ \mathbf{C}^{\mathcal{K}}(\mathbf{x}_{0}) } \big]^{1/2}  \frac{ \partial^{\mathbf{e}_{l}} \partial^{\mathbf{e}_{l}} 
  \big[ \detmb{  \mathbf{C}^{\mathcal{K}}(\mathbf{x}_{1}) + \mathbf{C}^{\mathcal{K}}(\mathbf{x}_{2}) - \mathbf{C}^{\mathcal{K}}(\mathbf{x}_{1}) 
  \big[ \mathbf{C}^{\mathcal{K}}(\mathbf{x}_{0})\big]^{-1} \mathbf{C}^{\mathcal{K}}(\mathbf{x}_{2}) } \big]^{-1/2} }  {\partial \mathbf{x}_{1}^{\mathbf{e}_{l}} \partial \mathbf{x}_{2}^{\mathbf{e}_{l}}} \bigg|_{\mathbf{x}_{1} = \mathbf{x}_{2} = \mathbf{0}} \nonumber \\[4mm]
  & \hspace*{-30mm} \stackrel{\eqref{equ_lem_der_det_function},\eqref{equ_def_projection_P_l}}{=} -\frac{1}{2} \big[ \detmb{ \mathbf{C}^{\mathcal{K}}(\mathbf{x}_{0}) } \big]^{1/2}    \nonumber \\[4mm]Ê
  & \hspace*{-25mm} \timesÊ \frac{ \partial^{\mathbf{e}_{l}}}{\partial \mathbf{x}_{1}^{\mathbf{e}_{l}}}
   \frac{  \trace \big\{ \big[  \mathbf{C}^{\mathcal{K}}(\mathbf{x}_{1}) + \mathbf{C}^{\mathcal{K}}(\mathbf{0}) - \mathbf{C}^{\mathcal{K}}(\mathbf{x}_{1}) 
  \big[ \mathbf{C}^{\mathcal{K}}(\mathbf{x}_{0})\big]^{-1} \mathbf{C}^{\mathcal{K}}(\mathbf{0})  \big]^{-1} \big[  \mathbf{I} - \mathbf{C}^{\mathcal{K}}(\mathbf{x}_{1}) 
 \big[ \mathbf{C}^{\mathcal{K}}(\mathbf{x}_{0})\big]^{-1}   \big] \mathbf{P}_{l} \big\} }    
  {\big[ \detmb{  \mathbf{C}^{\mathcal{K}}(\mathbf{x}_{1}) + \mathbf{C}^{\mathcal{K}}(\mathbf{0}) - \mathbf{C}^{\mathcal{K}}(\mathbf{x}_{1}) 
  \big[ \mathbf{C}^{\mathcal{K}}(\mathbf{x}_{0})\big]^{-1} \mathbf{C}^{\mathcal{K}}(\mathbf{0}) } \big]^{1/2}} \bigg|_{\mathbf{x}_{1} =  \mathbf{0}} \nonumber \\[4mm] 
    & \hspace*{-30mm} \stackrel{\eqref{equ_derivation_inverse_matrix},\eqref{equ_def_projection_P_l}}{=}  \Bigg[ \frac{ \detmb{ \mathbf{C}^{\mathcal{K}}(\mathbf{x}_{0}) }}{\detmb{ \mathbf{N} }} \Bigg]^{1/2} \nonumber \\[4mm]   
    & \hspace*{-5mm} \times  \Bigg[ \frac{1}{4} \trace \big\{ \mathbf{N}^{-1}  \big[\mathbf{I} -   \mathbf{C}^{\mathcal{K}}(\mathbf{0}) \big[ \mathbf{C}^{\mathcal{K}}(\mathbf{x}_{0})\big]^{-1}\big] \mathbf{P}_{l} \big\}
      \trace\big\{ \mathbf{N}^{-1}  \mathbf{P}_{l}  \big[\mathbf{I} -    \big[ \mathbf{C}^{\mathcal{K}}(\mathbf{x}_{0})\big]^{-1}\mathbf{C}^{\mathcal{K}}(\mathbf{0})\big]  \big\}
      \nonumber \\[4mm]
    & + \frac{1}{2}  \tracem{\mathbf{N}^{-1}   \mathbf{P}_{l} \big[\mathbf{I} -   \big[ \mathbf{C}^{\mathcal{K}}(\mathbf{x}_{0})\big]^{-1}  \mathbf{C}^{\mathcal{K}}(\mathbf{0}) \big]  \mathbf{N}^{-1} \big[\mathbf{I} -   \mathbf{C}^{\mathcal{K}}(\mathbf{0}) \big[ \mathbf{C}^{\mathcal{K}}(\mathbf{x}_{0})\big]^{-1} \big]  \mathbf{P}_{l} }  \nonumber \\[4mm]Ê
        &+ \frac{1}{2}  \trace \big\{\mathbf{N}^{-1}  \mathbf{P}_{l} \big[ \mathbf{C}^{\mathcal{K}}(\mathbf{x}_{0})\big]^{-1} \mathbf{P}_{l} \big\}  \Bigg],
\end{align}
where 
\begin{equation} 
\mathbf{N} \triangleq  2  \mathbf{C}^{\mathcal{K}}(\mathbf{0}) - \mathbf{C}^{\mathcal{K}}(\mathbf{0}) \big[ \mathbf{C}^{\mathcal{K}}(\mathbf{x}_{0})\big]^{-1} \mathbf{C}^{\mathcal{K}}(\mathbf{0})
\end{equation}
and 
\begin{equation} 
\label{equ_def_projection_P_l}
\mathbf{P}_{l} \triangleq \sum_{l' \in \mathcal{I}_{l}} \mathbf{e}_{l'} \mathbf{e}^{T}_{l'} = \frac{\partial^{\mathbf{e}_{l}}}{\partial \mathbf{x}^{\mathbf{e}_{l}}}\mathbf{C}^{\mathcal{K}}(\mathbf{x}),
\end{equation} 
with the index set 
\begin{equation} 
\label{equ_index_set_I_l_proof_non_supp_SPCM}
\mathcal{I}_{l} \triangleq \bigg\{ \sum_{k \in [l-1]} r_{k}+1, \sum_{k \in [l-1]} r_{k}+2,\ldots,\sum_{k \in [l-1]} r_{k}+r_{l}-1 \bigg\}.
\end{equation} 
The index set $\mathcal{I}_{l}$ consists of the indices of the main diagonal entries of $\mathbf{D}^{\mathcal{K}}(\mathbf{x})$ 
in \eqref{equ_diag_matrix_param_vec_SPCM} that correspond to $x_{l}$ as defined by \eqref{equ_diag_matrix_param_vec_SPCM}. 
Note that $\mathbf{P}_{l}$ is a psd orthogonal projection matrix, i.e., $\mathbf{P}_{l} \mathbf{P}_{l} = \mathbf{P}_{l}$.  

The various terms in \eqref{equ_proof_lower_bond_RIP_SPCM_squared_norm_expr} can be bounded by using Lemma \ref{lem_matrix_inversion_perturb}, the assumption $\delta_{S+1} < 1/32$, and the Definition in \equref{equ_def_transformed_cov_matrix_model_SPCM}, as follows. 
First, observe that by Lemma \ref{lem_eigvalues_sum_symmetric_matrices}, the eigenvalues
$\lambda_{l'_{k}} ( \mathbf{C}^{\mathcal{K}}(\mathbf{x}_{0}))$ of the matrix $\mathbf{C}^{\mathcal{K}}(\mathbf{x}_{0})$ (cf.\ \eqref{equ_def_transformed_cov_matrix_model_SPCM}), 
that correspond to the $k$th largest entry $x_{0,k}$ of $\mathbf{x}_{0}$, are upper bounded as 
\begin{align}
\label{equ_proof_RIP_Bound_SPCM_eignval_widetilde_C_K_x_0}
 \lambda_{l'_{k}} ( \mathbf{C}^{\mathcal{K}}(\mathbf{x}_{0})) & \stackrel{\eqref{equ_double_ineq_eigvals_sum}}{\leq} \lambda_{l'_{k}}(\mathbf{D}^{\mathcal{K}}(\mathbf{x}_{0}) + \sigma^{2} \mathbf{I}) + \lambda_{1}( \mathbf{E}_{1} ) \nonumber \\[4mm]
 & \stackrel{\eqref{equ_any_eigval_lower_mag_norm_2}}{\leq}  \lambda_{l'_{k}}(\mathbf{D}^{\mathcal{K}}(\mathbf{x}_{0}) + \sigma^{2} \mathbf{I}) + \| \mathbf{E}_{1} \|_{2}
 = x_{0,k} + \sigma^{2}+ \| \mathbf{E}_{1} \|_{2} \nonumber \\[4mm]
 & \stackrel{\eqref{equ_proof_lower_bound_RIP_SPCM_matrix_pertub}}{\leq}  x_{0,k}+ \sigma^{2} (1+5\delta_{S+1}).
\end{align} 
Note that there are $r_{k}$ eigenvalues $\lambda_{l'_{k}} ( \mathbf{C}^{\mathcal{K}}(\mathbf{x}_{0}))$ that correspond to $x_{0,k}$.
From this it follows that 
\begin{align} 
\label{equ_proof_RIP_lower_bound_SPCM_various_bound_1_part_1} 
 \big[ \detmb{ \mathbf{C}^{\mathcal{K}}(\mathbf{x}_{0}) } \big]^{1/2} & \stackrel{\eqref{equ_rel_eigvals_det}}{=}\bigg[ \prod_{k \inÊ\mathcal{K}} \prod_{l'_{k}} \lambda_{l'_{k}} ( \mathbf{C}^{\mathcal{K}}(\mathbf{x}_{0}))\bigg]^{1/2}
\stackrel{\eqref{equ_proof_RIP_Bound_SPCM_eignval_widetilde_C_K_x_0}}{\leq} \prod_{k \in \mathcal{K}} \big[x_{0,k}+ \sigma^{2} (1+5\delta_{S+1})\big]^{r_{k}/2} \nonumber \\[4mm]
& \stackrel{\eqref{equ_proof_RIP_bound_SPCM_index_set_K}}{=} \big[ \sigma^{2} (1+5\delta_{S+1})\big]^{r_{l}/2} \prod_{k \in \supp(\mathbf{x}_{0})} \big[x_{0,k}+ \sigma^{2} (1+5\delta_{S+1})\big]^{r_{k}/2}.
\end{align}

Similar to the bound \eqref{equ_proof_RIP_Bound_SPCM_eignval_widetilde_C_K_x_0}, we also have 
\begin{align}
\label{equ_proof_RIP_Bound_SPCM_min_eigval_C_K}
 \lambda_{l'} ( \mathbf{C}^{\mathcal{K}}(\mathbf{x}_{0})) & \stackrel{ \eqref{equ_def_transformed_cov_matrix_model_SPCM}}{=}  
  \lambda_{l'} ( \mathbf{D}^{\mathcal{K}}(\mathbf{x}_{0}) + \sigma^{2} \mathbf{V}_{1} \mathbf{\Sigma}_{1}^{-2} \mathbf{V}_{1}^{T})  \nonumber \\[4mm]
 & \stackrel{\eqref{equ_double_ineq_eigvals_sum}}{\geq}   \lambda_{l'} ( \sigma^{2} \mathbf{V}_{1} \mathbf{\Sigma}_{1}^{-2} \mathbf{V}_{1}^{T})
 +   \lambda_{Q} ( \mathbf{D}^{\mathcal{K}}(\mathbf{x}_{0})) \nonumber \\[4mm]
 & \geq  \lambda_{l'} ( \sigma^{2} \mathbf{V}_{1} \mathbf{\Sigma}_{1}^{-2} \mathbf{V}_{1}^{T}) = \sigma^{2} \lambda_{l'} (  \mathbf{\Sigma}_{1}^{-2} ) \geq \sigma^{2}(1+\delta_{S+1})^{-2},
\end{align}
and
\begin{align}
\label{equ_proof_RIP_Bound_SPCM_eignval_widetilde_C_K_0}
 \lambda_{l'} ( \mathbf{C}^{\mathcal{K}}(\mathbf{0})) & \stackrel{ \eqref{equ_def_transformed_cov_matrix_model_SPCM}}{=}  
  \lambda_{l'} ( \mathbf{D}^{\mathcal{K}}(\mathbf{0}) + \sigma^{2} \mathbf{V}_{1} \mathbf{\Sigma}_{1}^{-2} \mathbf{V}_{1}^{T})  \nonumber \\[4mm]
 & \stackrel{\eqref{equ_def_RIP_SPCM_bound_D_K}}{=} 
 \lambda_{l'} ( \sigma^{2} \mathbf{V}_{1} \mathbf{\Sigma}_{1}^{-2} \mathbf{V}_{1}^{T}) = \sigma^{2} \lambda_{l'} (  \mathbf{\Sigma}_{1}^{-2} ) \in \sigma^{2} [(1+\delta_{S+1})^{-2},(1-\delta_{S+1})^{-2}],
\end{align}
where the last step in \eqref{equ_proof_RIP_Bound_SPCM_min_eigval_C_K} and \eqref{equ_proof_RIP_Bound_SPCM_eignval_widetilde_C_K_0}, respectively, follow from the fact that the basis matrices satisfy the RIP of order $S+1$. From \eqref{equ_proof_RIP_Bound_SPCM_min_eigval_C_K} and 
we can conclude that 
\begin{equation} 
\label{equ_proof_RIP_SPCM_inv_C_K_norm}
\big\| \big[ \mathbf{C}^{\mathcal{K}}(\mathbf{x}_{0})\big]^{-1} \big\|_{2} \leq (1+\delta_{S+1})^{2} \sigma^{-2}
\end{equation} 
and in turn 
\begin{equation} 
\label{equ_proof_RiP_SPCM_3_factor_prod_bound_norm_2}
\big\| \mathbf{C}^{\mathcal{K}}(\mathbf{0}) \big[ \mathbf{C}^{\mathcal{K}}(\mathbf{x}_{0})\big]^{-1} \mathbf{C}^{\mathcal{K}}(\mathbf{0}) \big\|_{2} \leq 
\big\| \mathbf{C}^{\mathcal{K}}(\mathbf{0}) \big\|^{2}_{2} \,\, \big\|  \big[ \mathbf{C}^{\mathcal{K}}(\mathbf{x}_{0})\big]^{-1} \big\|_{2} \stackrel{\eqref{equ_proof_RIP_SPCM_inv_C_K_norm},\eqref{equ_proof_RIP_Bound_SPCM_eignval_widetilde_C_K_0}}\leq \sigma^{2} \frac{  (1+\delta_{S+1})^{2}}{(1-\delta_{S+1})^{4}}.
\end{equation} 
This yields further 
\begin{align} 
\label{equ_proof_RIP_SPCM_lower_bound_eig_vals_N}
\lambda_{l'}(\mathbf{N}) & \stackrel{\eqref{equ_double_ineq_eigvals_sum}}{\geq} \lambda_{l'}( 2  \mathbf{C}^{\mathcal{K}}(\mathbf{0}))+ \lambda_{Q}(- \mathbf{C}^{\mathcal{K}}(\mathbf{0}) \big[ \mathbf{C}^{\mathcal{K}}(\mathbf{x}_{0})\big]^{-1} \mathbf{C}^{\mathcal{K}}(\mathbf{0}) \nonumber \\[4mm]
 & \stackrel{\eqref{equ_any_eigval_lower_mag_norm_2}}{\geq}  \lambda_{l'}( 2  \mathbf{C}^{\mathcal{K}}(\mathbf{0}))- \big\| \mathbf{C}^{\mathcal{K}}(\mathbf{0}) \big[ \mathbf{C}^{\mathcal{K}}(\mathbf{x}_{0})\big]^{-1} \mathbf{C}^{\mathcal{K}}(\mathbf{0}) \big\|_{2} \nonumber \\[4mm]
  & \stackrel{\eqref{equ_proof_RiP_SPCM_3_factor_prod_bound_norm_2}}{\geq}  \lambda_{l'}( 2  \mathbf{C}^{\mathcal{K}}(\mathbf{0}))- \sigma^{2} \frac{  (1+\delta_{S+1})^{2}}{(1-\delta_{S+1})^{4}} \nonumber \\[4mm] 
    & \stackrel{\eqref{equ_proof_RIP_Bound_SPCM_eignval_widetilde_C_K_0}}{\geq}  2  \sigma^{2} (1+\delta_{S+1})^{-2}- \sigma^{2} \frac{  (1+\delta_{S+1})^{2}}{(1-\delta_{S+1})^{4}}  \nonumber \\[4mm] 
    & =  \sigma^{2} \frac{1}{(1+\delta_{S+1})^{2}} \bigg( 2 - \frac{(1+\delta_{S+1})^{4}}{(1-\delta_{S+1})^{4}} \bigg)= \sigma^{2} \beta,
\end{align}
and in turn
\begin{align} 
\label{equ_proof_RIP_lower_bound_SPCM_various_bound_1_part_2} 
 \big[ \detmb{ \mathbf{N} } \big]^{1/2} & \stackrel{\eqref{equ_rel_eigvals_det}}{=}\bigg[ \prod_{k \inÊ\mathcal{K}} \prod_{l'_{k}} \lambda_{l'_{k}} ( \mathbf{N}) \bigg]^{1/2}  
 \stackrel{\eqref{equ_proof_RIP_SPCM_lower_bound_eig_vals_N}}{\geq} \big( \sigma^{2} \beta \big)^{Q/2}.
\end{align} 
Note that \eqref{equ_proof_RIP_SPCM_lower_bound_eig_vals_N} also implies that 
\begin{equation} 
\label{equ_proof_RIP_SPCM_bound_norm_inv_N}
\| \mathbf{N}^{-1} \|_{2}  = \frac{1}{\min_{l' \in [Q]} |\lambda_{l'}(\mathbf{N})|} \leq \sigma^{2} \beta.
\end{equation}


Since we have (cf.\ \eqref{equ_def_RIP_SPCM_bound_D_K}) $\big\| \big( \mathbf{D}^{\mathcal{K}}(\mathbf{x}_{0}) + \sigma^{2} \mathbf{I} \big)^{-1} \big\|_{2} \leq  1/\sigma^{2}$, %
it holds that 
\begin{equation} 
\big\| \big( \mathbf{D}^{\mathcal{K}}(\mathbf{x}_{0}) + \sigma^{2} \mathbf{I} \big)^{-1} \cdot \mathbf{E}_{1} \big\|_{2} \leq \big\| \big( \mathbf{D}^{\mathcal{K}}(\mathbf{x}_{0}) + \sigma^{2} \mathbf{I} \big)^{-1} \big\|_{2}\cdot \|\mathbf{E}_{1} \|_{2} \stackrel{\eqref{equ_proof_lower_bound_RIP_SPCM_matrix_pertub}}{\leq}  5 \delta_{S+1} \stackrel{\delta_{S+1} < 1/32}{<} 1.
\end{equation}
This allows us to involve Lemma \ref{lem_matrix_inversion_perturb}, with the choices $\mathbf{A} = \mathbf{D}^{\mathcal{K}}(\mathbf{x}_{0}) + \sigma^{2} \mathbf{I}$ and 
$\mathbf{E} =  \mathbf{C}^{\mathcal{K}}(\mathbf{x}_{0})-  \big[ \mathbf{D}^{\mathcal{K}}(\mathbf{x}_{0}) + \sigma^{2} \mathbf{I} \big] = \mathbf{E}_{1}$ 
(cf.\ \eqref{equ_def_transformed_cov_matrix_model_SPCM}), to obtain 
\begin{equation}
\label{equ_C_K_inverse} 
\big[ \mathbf{C}^{\mathcal{K}}(\mathbf{x}_{0})\big]^{-1}= \big[ \mathbf{D}^{\mathcal{K}}(\mathbf{x}_{0}) + \sigma^{2} \mathbf{I} \big]^{-1} + \mathbf{E}' 
\end{equation} 
with an error term $\mathbf{E}'$ that satisfies 
\begin{align}
\label{equ_error_norm_bound_C_K_inverse} 
\| \mathbf{E}' \|_{2} & = \big\| \big[ \mathbf{C}^{\mathcal{K}}(\mathbf{x}_{0})\big]^{-1}-\big[ \mathbf{D}^{\mathcal{K}}(\mathbf{x}_{0}) + \sigma^{2} \mathbf{I} \big]^{-1} \big\|_{2} 
\nonumber \\[4mm]
& \stackrel{\eqref{equ_matrix_inversion_perturb}}{\leq}  \frac{\big\| \big[ \mathbf{D}^{\mathcal{K}}(\mathbf{x}_{0}) + \sigma^{2} \mathbf{I} \big]^{-1}Ê\big\|^{2}_{2} 
\| \mathbf{E}_{1} \|_{2}}{1 - \big\| \big[ \mathbf{D}^{\mathcal{K}}(\mathbf{x}_{0}) + \sigma^{2} \mathbf{I} \big]^{-1}Ê\big\|_{2} \| \mathbf{E}_{1} \|_{2}} 
\leq  \frac{1}{\sigma^{4}} \frac{ \| \mathbf{E}_{1} \|_{2}}{1 - \| \mathbf{E}_{1} \|_{2} / \sigma^{2} }  \nonumber \\[4mm]
& \stackrel{\eqref{equ_proof_lower_bound_RIP_SPCM_matrix_pertub}}{\leq} \frac{1}{\sigma^{4}} \frac{ \sigma^{2} 5 \delta_{S+1} }{1 - 5 \delta_{S+1}}
 \stackrel{\delta_{S+1} < 1/32}{\leq}  \frac{1}{\sigma^{2}} \frac{  5 \delta_{S+1} }{1 - 5/32}
\leq \frac{1}{\sigma^{2}} 6 \delta_{S+1}.
\end{align} 
In what follows, we will also need the identity \cite{golub96}
\begin{equation} 
\label{equ_identiy_trace_ONB}
\tracem{\mathbf{A}}  = \sum_{l' \in [Q]} \mathbf{u}_{l'}^{T} \mathbf{A} \mathbf{u}_{l'}
\end{equation} 
valid for any matrix $\mathbf{A} \in \mathbb{R}^{Q \times Q}$ and 
any set $\big\{ \mathbf{u}_{l'} \in \mathbb{R}^{T} \big\}_{l' \in [Q]}$ of vectors that form an ONB for $\mathbb{R}^{Q}$. 
We then obtain
\begin{align}
\label{equ_proof_RIP_lower_bound_SPCM_various_bound_2} 
 & \big| \trace \big\{ \mathbf{N}^{-1}  \big[\mathbf{I} -   \mathbf{C}^{\mathcal{K}}(\mathbf{0})\big[ \mathbf{C}^{\mathcal{K}}(\mathbf{x}_{0})\big]^{-1}\big] \mathbf{P}_{l} \big\} \big|   \nonumber \\[4mm]
 & \stackrel{\eqref{equ_identiy_trace_ONB}}{=} \bigg|  \sum_{l' \in [N]} \mathbf{e}_{l'}^{T} \mathbf{N}^{-1}  \big[\mathbf{I} -   \mathbf{C}^{\mathcal{K}}(\mathbf{0}) \big[ \mathbf{C}^{\mathcal{K}}(\mathbf{x}_{0})\big]^{-1} \big] \mathbf{P}_{l} \mathbf{e}_{l'} Ê \bigg|
   \stackrel{\eqref{equ_def_projection_P_l}}{=}  \bigg| \sum_{l' \in \mathcal{I}_{l}} \mathbf{e}_{l'}^{T} \mathbf{N}^{-1}  \big[\mathbf{I} -   \mathbf{C}^{\mathcal{K}}(\mathbf{0}) \big[ \mathbf{C}^{\mathcal{K}}(\mathbf{x}_{0})\big]^{-1} \big] \mathbf{e}_{l'} Ê \bigg|   \nonumber \\[4mm]
   & \stackrel{\eqref{equ_def_transformed_cov_matrix_model_SPCM},\eqref{equ_C_K_inverse}}{=}  \bigg| \sum_{l' \in \mathcal{I}_{l}} \mathbf{e}_{l'}^{T} \mathbf{N}^{-1}  \bigg[\mathbf{I} -   (\sigma^{2} \mathbf{I} + \mathbf{E}_{1})   \big[( \mathbf{D}^{\mathcal{K}}(\mathbf{x}_{0}) + \sigma^{2} \mathbf{I} )^{-1}
    + \mathbf{E}'\big]\bigg]  \mathbf{e}_{l'}\bigg|  \nonumber \\[4mm]
   & = \bigg| \sum_{l' \in \mathcal{I}_{l}} \mathbf{e}_{l'}^{T} \mathbf{N}^{-1}  (\mathbf{I} -   \sigma^{2} ( \mathbf{D}^{\mathcal{K}}(\mathbf{x}_{0}) + \sigma^{2} \mathbf{I} )^{-1} 
   - \mathbf{E}_{1}( \mathbf{D}^{\mathcal{K}}(\mathbf{x}_{0}) + \sigma^{2} \mathbf{I} )^{-1} - \sigma^{2}  \mathbf{E}' - \mathbf{E}_{1} \mathbf{E}')  \mathbf{e}_{l'}\bigg|.
   \end{align}
Using the assumption that $l \notin \supp(\mathbf{x}_{0})$, implying that $( \mathbf{D}^{\mathcal{K}}(\mathbf{x}_{0}) + \sigma^{2} \mathbf{I} )^{-1} \mathbf{e}_{l'} = \mathbf{e}_{l'} /\sigma^{2}$ for 
every $l' \in \mathcal{I}_{l}$, this further becomes
\begin{align}
\label{equ_proof_RIP_lower_bound_SPCM_various_bound_2} 
 & \big| \trace \big\{ \mathbf{N}^{-1}  \big[\mathbf{I} -   \mathbf{C}^{\mathcal{K}}(\mathbf{0})\big[ \mathbf{C}^{\mathcal{K}}(\mathbf{x}_{0})\big]^{-1}\big] \mathbf{P}_{l} \big\} \big|   \nonumber \\[4mm]
   &  = \bigg| \sum_{l' \in \mathcal{I}_{l}} \mathbf{e}_{l'}^{T} \mathbf{N}^{-1}  \big[\mathbf{e}_{l'} -  \mathbf{e}_{l'} - \mathbf{E}_{1}( \mathbf{D}^{\mathcal{K}}(\mathbf{x}_{0}) 
   + \sigma^{2} \mathbf{I} )^{-1}\mathbf{e}_{l'} - \sigma^{2}  \mathbf{E}'\mathbf{e}_{l'} - \mathbf{E}_{1} \mathbf{E}'  \mathbf{e}_{l'} \big]\bigg|
  \nonumber \\[4mm]
   & = \bigg| \sum_{l' \in \mathcal{I}_{l}} \mathbf{e}_{l'}^{T} \mathbf{N}^{-1}  \big[ -\mathbf{E}_{1} ( \mathbf{D}^{\mathcal{K}}(\mathbf{x}_{0}) + \sigma^{2} \mathbf{I} )^{-1}\mathbf{e}_{l'} 
   - \sigma^{2}  \mathbf{E}'\mathbf{e}_{l'} - \mathbf{E}_{1} \mathbf{E}'  \mathbf{e}_{l'}\big]\bigg|
  \nonumber \\[4mm]
   &  \leq \sum_{l' \in \mathcal{I}_{l}} \bigg|  \mathbf{e}_{l'}^{T} \mathbf{N}^{-1} 
    Ê\big[ -\mathbf{E}_{1}( \mathbf{D}^{\mathcal{K}}(\mathbf{x}_{0}) + \sigma^{2} \mathbf{I} )^{-1}\mathbf{e}_{l'} - \sigma^{2}  \mathbf{E}'\mathbf{e}_{l'} - \mathbf{E}_{1} \mathbf{E}'  \mathbf{e}_{l'} \big]\bigg|
      \nonumber \\[4mm]
   &\stackrel{(b)}{\leq}  r_{l}  \frac{1}{\sigma^{2} \beta} ( 5 \delta_{S+1} +6Ê\delta_{S+1}+ 30 \delta_{S+1}^{2}) \stackrel{\delta_{S+1} < 1/32}{\leq}  r_{l}  \frac{1}{\sigma^{2} \beta} 12 \delta_{S+1},
\end{align} 
where $(b)$ follows from the inequality $\| \mathbf{A} \mathbf{B} \|_{2} \leq \|  \mathbf{A}\|_{2} \| \mathbf{B} \|_{2}$, which is valid for any two matrices 
$\mathbf{A},Ê\mathbf{B} \in \mathbb{R}^{Q \times Q}$ \cite{golub96}. 
For the step $(b)$ we also used \eqref{equ_proof_lower_bound_RIP_SPCM_matrix_pertub}, \eqref{equ_proof_RIP_SPCM_bound_norm_inv_N}, and \eqref{equ_error_norm_bound_C_K_inverse}. 

Similar to \eqref{equ_proof_RIP_lower_bound_SPCM_various_bound_2}, one can also derive the inequalities 
\begin{align}
\label{equ_proof_RIP_lower_bound_SPCM_various_bound_3} 
& \big|Ê\trace \big\{ \mathbf{N}^{-1} \mathbf{P}_{l}  \big[ \mathbf{I} -    \big[ \mathbf{C}^{\mathcal{K}}(\mathbf{x}_{0})\big]^{-1}\mathbf{C}^{\mathcal{K}}(\mathbf{0}) \big]  \big\} \big| \leq r_{l} 
\frac{1}{\sigma^{2} \beta} 12 \delta_{S+1}Ê\nonumber \\[4mm]Ê
& \big|Ê  \trace \big\{\mathbf{N}^{-1}  \mathbf{P}_{l} \big[\mathbf{I} -   \big[ \mathbf{C}^{\mathcal{K}}(\mathbf{x}_{0})\big]^{-1}  \mathbf{C}^{\mathcal{K}}(\mathbf{0})\big]  \mathbf{N}^{-1} \big[ \mathbf{I} -  \mathbf{C}^{\mathcal{K}}(\mathbf{0})\big[ \mathbf{C}^{\mathcal{K}}(\mathbf{x}_{0})\big]^{-1} \big] \mathbf{P}_{l} \big\} \big| \leq r_{l}   (12  \delta_{S+1} )^{2} \frac{1}{ \sigma^{4} \beta^{2}} \nonumber \\[4mm]Ê
& |Ê  \tracem{\mathbf{N}^{-1} \mathbf{P}_{l} \big[ \mathbf{C}^{\mathcal{K}}(\mathbf{x}_{0})\big]^{-1}\mathbf{P}_{l}} |  \leq r_{l} \frac{1}{\beta \sigma^{4}} (6 \delta_{S+1}+ 1).
\end{align}
The bound in \eqref{equ_lower_bound_SPCM_RIP} follows then finally by inserting the bounds \eqref{equ_proof_RIP_lower_bound_SPCM_various_bound_1_part_1}, 
\eqref{equ_proof_RIP_lower_bound_SPCM_various_bound_1_part_2}, \eqref{equ_proof_RIP_lower_bound_SPCM_various_bound_2}, \eqref{equ_proof_RIP_lower_bound_SPCM_various_bound_3} into \eqref{equ_proof_lower_bond_RIP_SPCM_squared_norm_expr} and using \eqref{equ_proof_lower_bound_RIP_projection_SPCM}. 
\end{proof}

\chapter[]{Proof of Lemma \ref{lem_cond_exist_partial_der_exist_spec_est_SDPCM}}
\label{chap_appendix_E}
\begin{proof}
For the following, we restrict ourselves to the case where $\mathbf{p} = \mathbf{e}_{l}$ with an arbitrary index $l \in [N]$ and introduce the function $\mathbf{z}(\mathbf{x}): \mathbb{R}_{+}^{N} \rightarrow \mathbb{R}_{+}^{N}$ defined 
elementwise by $z_{k}(\mathbf{x}) = \frac{1}{x_{k} + \sigma^{2}}$. Moreover, we define the modified mean function $\widetilde{m}_{k}(\mathbf{z})$ by requiring that $\widetilde{m}_{k}(\mathbf{z}(\mathbf{x})) = m_{k}(\mathbf{x})$. 
We have 
\begin{align} 
\label{equ_proof_part_der_SDPCM_1}
\widetilde{m}_{k}(\mathbf{z}(\mathbf{x})) = m_{k}(\mathbf{x}) & = \mathsf{E}_{\mathbf{x}} \{ \hat{x}_{k}(\mathbf{y})  \} =   \int_{\mathbf{y} \in \mathbb{R}^{M}} \hat{x}_{k}(\mathbf{y}) f_{\text{SDPCM}}(\mathbf{y}; \mathbf{x}) d \mathbf{y}\nonumber \\[4mm] 
&\hspace*{-30mm}  \stackrel{\eqref{equ_def_SDPCM_stat_model}}{=}  \frac{1}{(2 \pi)^{M/2} [\detm{\widetilde{\mathbf{C}}(\mathbf{x})}]^{1/2}} \int_{\mathbf{y} \in \mathbb{R}^{M}} \hat{x}_{k}(\mathbf{y}) \exp \left( - \frac{1}{2} \mathbf{y}^{T} \widetilde{\mathbf{C}}^{-1}(\mathbf{x}) \mathbf{y} \right) d \mathbf{y}\nonumber \\[4mm] 
&\hspace*{-30mm} \stackrel{\eqref{equ_basis_matrix_sum_rank_1_projection_SDPCM},\eqref{equ_expr_obs_cov_matrix_SDPCM}}{=}  \frac{1}{(2 \pi)^{M/2} [\detm{\widetilde{\mathbf{C}}(\mathbf{x})}]^{1/2}} \nonumber \\[4mm]
& \hspace*{-25mm} \times \int_{\mathbf{y} \in \mathbb{R}^{M}} \hat{x}_{k}(\mathbf{y}) \exp \left( - \frac{1}{2} \bigg[\sum_{k \in [N]} z_{k}(\mathbf{x}) \sum_{i \in [r_{k}]} \big[\mathbf{u}^{T}_{m_{k,i}} \mathbf{y}\big]^{2} + \frac{1}{\sigma^{2}}\mathbf{y}^{T} \mathbf{P}_{n}\mathbf{y}\bigg] \right) d \mathbf{y} \nonumber \\[4mm]
&\hspace*{-30mm}  \stackrel{\eqref{equ_expr_obs_cov_matrix_SDPCM}}{=}  \frac{1}{(2 \pi)^{M/2}}  \frac{1}{\prod_{k \in [N]} (x_{k} + \sigma^{2})^{r_{k}/2}}  \frac{1}{\sigma^{M-R}}  \nonumber \\[4mm]
& \hspace*{-25mm}  \times \int_{\mathbf{y} \in \mathbb{R}^{M}} \hat{x}_{k}(\mathbf{y}) \exp \left( - \frac{1}{2} \bigg[\sum_{k \in [N]} z_{k}(\mathbf{x}) \sum_{i \in [r_{k}]} \big[\mathbf{u}^{T}_{m_{k,i}} \mathbf{y}\big]^{2} + \frac{1}{\sigma^{2}}\mathbf{y}^{T} \mathbf{P}_{n}\mathbf{y}\bigg] \right) d \mathbf{y}\nonumber \\[4mm]
&\hspace*{-30mm}  =  \frac{1}{(2 \pi)^{M/2}\sigma^{M-R}}  \prod_{k \in [N]} (z_{k}(\mathbf{x}))^{r_{k}/2}    \nonumber \\[4mm]
& \hspace*{-25mm}  \times \int_{\mathbf{y} \in \mathbb{R}^{M}} \hat{x}_{k}(\mathbf{y}) \exp \left( - \frac{1}{2} \bigg[\sum_{k \in [N]} z_{k}(\mathbf{x}) \sum_{i \in [r_{k}]} \big[\mathbf{u}^{T}_{m_{k,i}} \mathbf{y}\big]^{2} + \frac{1}{\sigma^{2}}\mathbf{y}^{T} \mathbf{P}_{n}\mathbf{y}\bigg] \right) d \mathbf{y}\nonumber \\[4mm]
&\hspace*{-30mm}  \stackrel{\tilde{\mathbf{y}}=\mathbf{U} \mathbf{y}}{=}  \frac{1}{(2 \pi)^{M/2}\sigma^{M-R}}  \prod_{k \in [N]} (z_{k}(\mathbf{x}))^{r_{k}/2}    \nonumber \\[4mm]
& \hspace*{-25mm}  \times \int_{\tilde{\mathbf{y}} \in \mathbb{R}^{M}} \hat{x}_{k}(\mathbf{U}^{T} \tilde{\mathbf{y}}) \exp \left( - \frac{1}{2} \bigg[\sum_{k \in [N]} z_{k}(\mathbf{x}) \sum_{i \in [r_{k}]} \tilde{y}_{m_{i,k}}^{2} + \frac{1}{\sigma^{2}} \sum_{k \in [M]\setminus[R]} \tilde{y}_{k}^{2}\bigg]  \right) d \tilde{\mathbf{y}} , 
\end{align}
where $\mathbf{P}_{n}$ denotes the projection matrix on the subspace of $\mathbb{R}^{M}$ which is orthogonal to the subspace $\bigcup_{k \in [N]} \linspan(\mathbf{C}_{k})$ with $\mathbf{C}_{k}$ being the basis matrices of the SDPCM. 
The matrix $\mathbf{U}$, used for the substitution $\mathbf{y} \rightarrow \tilde{\mathbf{y}}$, is the orthonormal matrix that appears in \eqref{equ_expr_obs_cov_matrix_SDPCM}.
The existence of the integrals in \eqref{equ_proof_part_der_SDPCM_1} 
follows from the dominated convergence theorem \cite{HalmosMeasure,RudinBookPrinciplesMatheAnalysis,RudinBook} since the function 
$\hat{x}_{k}(\mathbf{y})  \exp \left( - \frac{1}{2} \mathbf{y}^{T} \widetilde{\mathbf{C}}^{-1}(\mathbf{x}) \mathbf{y} \right)$
is measurable, since $\hat{x}_{k}(\mathbf{y})$ is assumed measurable \cite{HalmosMeasure}. 
Moreover, because of \eqref{equ_bound_suff_cond_exist_par_der_exist_SDPCM_spec_est}, the above function is upper bounded in magnitude (dominated) by the function 
$C \| \mathbf{y} \|_{2}^{L} \exp \left( - \frac{1}{2} \mathbf{y}^{T} \widetilde{\mathbf{C}}^{-1}(\mathbf{x}) \mathbf{y} \right)$, 
which is obviously integrable. 
Let us now calculate the partial derivative $\frac{ \partial^{\mathbf{e}_{l}} \widetilde{m}_{k}(\mathbf{z})}{\partial \mathbf{z}^{\mathbf{e}_{l}}}$. 
\begin{align}
\label{equ_proof_part_der_tilde_m_SDPCM}
\frac{ \partial^{\mathbf{e}_{l}} \widetilde{m}_{k}(\mathbf{z})}{\partial \mathbf{z}^{\mathbf{e}_{l}}}& \stackrel{\eqref{equ_proof_part_der_SDPCM_1}}{=} \frac{1}{(2 \pi)^{M/2}\sigma^{M-R}}\frac{ \partial^{\mathbf{e}_{l}}}{\partial \mathbf{z}^{\mathbf{e}_{l}}} \prod_{k \in [N]} z_{k}^{r_{k}/2}    \nonumber \\[4mm]
& \hspace*{-5mm}  \times \int_{\tilde{\mathbf{y}} \in \mathbb{R}^{M}} \hat{x}_{k}(\mathbf{U}^{T} \tilde{\mathbf{y}}) 
\exp \left( - \frac{1}{2} \bigg[\sum_{k \in [N]} z_{k} \sum_{i \in [r_{k}]} \tilde{y}_{m_{i,k}}^{2} + \frac{1}{\sigma^{2}} \sum_{k \in [M]\setminus[R]} \tilde{y}_{k}^{2}\bigg]  \right) d \tilde{\mathbf{y}} 
 \nonumber \\[4mm]
& \hspace*{-15mm} \stackrel{(a)}{=} \frac{1}{(2 \pi)^{M/2}\sigma^{M-R}} \nonumber \\[4mm]
 & \hspace*{-10mm}  \Bigg\{ \Bigg[ \frac{ \partial^{\mathbf{e}_{l}}}{\partial \mathbf{z}^{\mathbf{e}_{l}}} \prod_{k \in [N]} z_{k}^{r_{k}/2} \Bigg]  \int_{\tilde{\mathbf{y}} \in \mathbb{R}^{M}} \hat{x}_{k}(\mathbf{U}^{T} \tilde{\mathbf{y}})
  \exp \left( - \frac{1}{2} \bigg[\sum_{k \in [N]} z_{k} \sum_{i \in [r_{k}]} \tilde{y}_{m_{i,k}}^{2} + \frac{1}{\sigma^{2}} \sum_{k \in [M]\setminus[R]} \tilde{y}_{k}^{2} \bigg] \right) d \tilde{\mathbf{y}}  \nonumber \\[4mm]
 &\hspace*{-10mm}  + \prod_{k \in [N]} (z_{k})^{r_{k}/2} \Bigg[ \frac{ \partial^{\mathbf{e}_{l}}}{\partial \mathbf{z}^{\mathbf{e}_{l}}}  \int_{\tilde{\mathbf{y}} \in \mathbb{R}^{M}} \hat{x}_{k}(\mathbf{U}^{T} \tilde{\mathbf{y}})
  \exp \left( - \frac{1}{2} \bigg[\sum_{k \in [N]} z_{k}\sum_{i \in [r_{k}]} \tilde{y}_{m_{i,k}}^{2} + \frac{1}{\sigma^{2}} \sum_{k \in [M]\setminus[R]} \tilde{y}_{k}^{2}\bigg]  \right) d \tilde{\mathbf{y}} \bigg] \Bigg\}
   \nonumber \\[4mm]
   & \hspace*{-15mm} \stackrel{(b)}{=} \frac{1}{(2 \pi)^{M/2}\sigma^{M-R}} \nonumber \\[4mm]
 & \hspace*{-10mm}  \Bigg\{ \Bigg[ \frac{ \partial^{\mathbf{e}_{l}}}{\partial \mathbf{z}^{\mathbf{e}_{l}}} \prod_{k \in [N]} z_{k}^{r_{k}/2} \Bigg]  \int_{\tilde{\mathbf{y}} \in \mathbb{R}^{M}} \hat{x}_{k}(\mathbf{U}^{T} \tilde{\mathbf{y}})
  \exp \left( - \frac{1}{2} \bigg[\sum_{k \in [N]} z_{k} \sum_{i \in [r_{k}]} \tilde{y}_{m_{i,k}}^{2} + \frac{1}{\sigma^{2}} \sum_{k \in [M]\setminus[R]} \tilde{y}_{k}^{2} \bigg] \right) d \tilde{\mathbf{y}}  \nonumber \\[4mm]
 &\hspace*{-10mm}  + \prod_{k \in [N]} (z_{k})^{r_{k}/2} \Bigg[  \int_{\tilde{\mathbf{y}} \in \mathbb{R}^{M}} \hat{x}_{k}(\mathbf{U}^{T} \tilde{\mathbf{y}}) \frac{ \partial^{\mathbf{e}_{l}}}{\partial \mathbf{z}^{\mathbf{e}_{l}}}
  \exp \left( - \frac{1}{2} \bigg[\sum_{k \in [N]} z_{k}\sum_{i \in [r_{k}]} \tilde{y}_{m_{i,k}}^{2} + \frac{1}{\sigma^{2}} \sum_{k \in [M]\setminus[R]} \tilde{y}_{k}^{2}\bigg]  \right) d \tilde{\mathbf{y}} \bigg] \Bigg\}
   \nonumber \\[4mm]
& \hspace*{-15mm}= \frac{1}{(2 \pi)^{M/2}\sigma^{M-R}} \nonumber \\[4mm]
 & \hspace*{-10mm}Ê\times  \Bigg\{ \Bigg[\int_{\tilde{\mathbf{y}} \in \mathbb{R}^{M}} \hat{x}_{k}(\mathbf{U}^{T} \tilde{\mathbf{y}})
  \exp \left( - \frac{1}{2} \bigg[\sum_{k \in [N]} z_{k} \sum_{i \in [r_{k}]} \tilde{y}_{m_{i,k}}^{2} + \frac{1}{\sigma^{2}} \sum_{k \in [M]\setminus[R]} \tilde{y}_{k}^{2} \bigg]  \right) d \tilde{\mathbf{y}}\Bigg]\frac{r_{l}}{2z_{l}} \prod_{k \in [N]} z_{k}^{r_{k}/2}  \nonumber \\[4mm]
 &\hspace*{-20mm}  +\Bigg[   \int_{\tilde{\mathbf{y}} \in \mathbb{R}^{M}}  \frac{-\sum_{i \in [r_{l}]} \tilde{y}_{m_{i,l}}^{2}}{2}\hat{x}_{k}(\mathbf{U}^{T} \tilde{\mathbf{y}})
  \exp \left( - \frac{1}{2} \bigg[\sum_{k \in [N]} z_{k} \sum_{i \in [r_{k}]} \tilde{y}_{m_{i,k}}^{2} + \frac{1}{\sigma^{2}} \sum_{k \in [M]\setminus[R]} \tilde{y}_{k}^{2} \bigg]  \right) d \tilde{\mathbf{y}} \Bigg]  \prod_{k \in [N]} z_{k}^{r_{k}/2} \Bigg\},
\end{align} 
where $(a)$ follows from the product rule for differentiation \cite{RudinBookPrinciplesMatheAnalysis} and $(b)$ results from a change of the order of differentiation and integration 
\cite[Theorem 1.5.8]{LC}, which can be verified, e.g., by a standard argument using the dominated convergence theorem \cite{FundmentExpFamBrown,RudinBookPrinciplesMatheAnalysis,HalmosMeasure}.  
From \eqref{equ_proof_part_der_tilde_m_SDPCM}, we obtain by the chain rule for differentiation \cite[Theorem 9.15]{RudinBookPrinciplesMatheAnalysis} that   
\begin{align}
 \frac{ \partial^{\mathbf{e}_{l}} m_{k}(\mathbf{x})}{\partial \mathbf{x}^{\mathbf{e}_{l}}} 
 & =  \frac{ \partial^{\mathbf{e}_{l}} \widetilde{m}_{k}(\mathbf{z})}{\partial \mathbf{z}^{\mathbf{e}_{l}}} \frac{ \partial z_{l}(\mathbf{x})}{\partial x_{l}}
 =  \frac{ \partial^{\mathbf{e}_{l}} \widetilde{m}_{k}(\mathbf{z})}{\partial \mathbf{z}^{\mathbf{e}_{l}}} \frac{-1}{(x_{l}+\sigma^{2})^{2}} \nonumber \\[4mm]
 & \hspace*{-15mm} \stackrel{\eqref{equ_proof_part_der_tilde_m_SDPCM}}{=} \frac{1}{(2 \pi)^{M/2}\sigma^{M-R}}\cdot \frac{-1}{(x_{l}+\sigma^{2})^{2}} \nonumber \\[4mm]
 & \hspace*{-10mm}Ê\times  \Bigg\{\Bigg[\int_{\tilde{\mathbf{y}} \in \mathbb{R}^{M}} \hat{x}_{k}(\mathbf{U}^{T} \tilde{\mathbf{y}})
  \exp \left( - \frac{1}{2} \bigg[\sum_{k \in [N]} z_{k} \sum_{i \in [r_{k}]} \tilde{y}_{m_{i,k}}^{2} + \frac{1}{\sigma^{2}} \sum_{k \in [M]\setminus[R]} \tilde{y}_{k}^{2} \bigg]  \right) d \tilde{\mathbf{y}}\Bigg]  \frac{r_{l}}{2z_{l}} \prod_{k \in [N]} z_{k}^{r_{k}/2} \nonumber \\[4mm]
 &\hspace*{-20mm}  + \Bigg[   \int_{\tilde{\mathbf{y}} \in \mathbb{R}^{M}}  \frac{-\sum_{i \in [r_{l}]} \tilde{y}_{m_{i,l}}^{2}}{2}\hat{x}_{k}(\mathbf{U}^{T} \tilde{\mathbf{y}})
  \exp \left( - \frac{1}{2} \bigg[\sum_{k \in [N]} z_{k} \sum_{i \in [r_{k}]} \tilde{y}_{m_{i,k}}^{2} + \frac{1}{\sigma^{2}} \sum_{k \in [M]\setminus[R]} \tilde{y}_{k}^{2} \bigg]  \right) d \tilde{\mathbf{y}} \Bigg] \prod_{k \in [N]} z_{k}^{r_{k}/2} \Bigg\} \nonumber \\[4mm]
   & \hspace*{-15mm} \stackrel{\eqref{equ_basis_matrix_sum_rank_1_projection_SDPCM},\eqref{equ_expr_obs_cov_matrix_SDPCM}}{=} \frac{1}{(2 \pi)^{M/2}\sigma^{M-R}}\cdot\frac{-1}{(x_{l}+\sigma^{2})^{2}} \nonumber \\[4mm]
 & \hspace*{-10mm}  \times \Bigg\{\bigg[ \int_{\mathbf{y} \in \mathbb{R}^{M}} \hat{x}_{k}(\mathbf{y})
  \exp \left( - \frac{1}{2} \mathbf{y}^{T} \widetilde{\mathbf{C}}^{-1}(\mathbf{x}) \mathbf{y} \right) d \mathbf{y}\bigg] \frac{r_{l}}{2z_{l}} \prod_{k \in [N]} z_{k}^{r_{k}/2} \nonumber \\[4mm]
 & + \Bigg[   \int_{\mathbf{y} \in \mathbb{R}^{M}}  \frac{-\mathbf{y}^{T} \mathbf{C}_{l} \mathbf{y}}{2} \hat{x}_{k}(\mathbf{y})
  \exp \left( - \frac{1}{2} \mathbf{y}^{T} \widetilde{\mathbf{C}}^{-1}(\mathbf{x}) \mathbf{y} \right) d \mathbf{y} \Bigg] \prod_{k \in [N]} z_{k}^{r_{k}/2} \Bigg\} \nonumber \\[4mm]
   & \hspace*{-15mm} \stackrel{\eqref{equ_expr_obs_cov_matrix_SDPCM}}{=}   \frac{1}{(2 \pi)^{M/2} [\detm{\widetilde{\mathbf{C}}(\mathbf{x})}]^{1/2}}\cdot \frac{-1}{(x_{l}+\sigma^{2})^{2}} \nonumber \\[4mm]
 & \hspace*{-10mm}  \times \Bigg\{\Bigg[ \int_{\mathbf{y} \in \mathbb{R}^{M}} \hat{x}_{k}(\mathbf{y})
  \exp \left( - \frac{1}{2} \mathbf{y}^{T} \widetilde{\mathbf{C}}^{-1}(\mathbf{x}) \mathbf{y} \right) d \mathbf{y}\Bigg] \frac{r_{l}}{2z_{l}}  \nonumber \\[4mm]
 & + \Bigg[   \int_{\mathbf{y} \in \mathbb{R}^{M}}  \frac{-\mathbf{y}^{T} \mathbf{C}_{l} \mathbf{y}}{2}  \hat{x}_{k}(\mathbf{y})
  \exp \left( - \frac{1}{2} \mathbf{y}^{T} \widetilde{\mathbf{C}}^{-1}(\mathbf{x}) \mathbf{y} \right) d \mathbf{y} \Bigg]  \Bigg\} \nonumber \\[4mm]
   & \hspace*{-15mm} =   \frac{1}{(2 \pi)^{M/2} [\detm{\widetilde{\mathbf{C}}(\mathbf{x})}]^{1/2}}\nonumber \\[4mm]
 & \hspace*{-10mm}  \times \Bigg\{\bigg[ \int_{\mathbf{y} \in \mathbb{R}^{M}} \hat{x}_{k}(\mathbf{y})
  \exp \left( - \frac{1}{2} \mathbf{y}^{T} \widetilde{\mathbf{C}}^{-1}(\mathbf{x}) \mathbf{y} \right) d \mathbf{y}\bigg] \frac{-r_{l}}{2(x_{l}+\sigma^{2})}  \nonumber \\[4mm]
 & +  \frac{-1}{(x_{l}+\sigma^{2})^{2}}  \int_{\mathbf{y} \in \mathbb{R}^{M}}  \frac{-\mathbf{y}^{T} \mathbf{C}_{l} \mathbf{y}}{2}  \hat{x}_{k}(\mathbf{y})
  \exp \left( - \frac{1}{2} \mathbf{y}^{T} \widetilde{\mathbf{C}}^{-1}(\mathbf{x}) \mathbf{y} \right) d \mathbf{y}  \Bigg\} 
  \nonumber \\[4mm]
   & \hspace*{-15mm} =  -\frac{r_{l}}{2(x_{l}+\sigma^{2})} \mathsf{E}_{\mathbf{x}} \big\{ \hat{x}_{k}(\mathbf{y}) \big\} + \frac{1}{2(x_{l}+\sigma^{2})^{2}}  \mathsf{E}_{\mathbf{x}} \big\{ \mathbf{y}^{T} \mathbf{C}_{l} \mathbf{y}\hat{x}_{k}(\mathbf{y}) \big\}  \nonumber \\[4mm]
      & \hspace*{-15mm} =  -\frac{r_{l}}{2(x_{l}+\sigma^{2})} m_{k}(\mathbf{x}) + \frac{1}{2(x_{l}+\sigma^{2})^{2}}  \mathsf{E}_{\mathbf{x}} \big\{ \mathbf{y}^{T} \mathbf{C}_{l} \mathbf{y}\hat{x}_{k}(\mathbf{y}) \big\}.  \nonumber \\[4mm]
\end{align}
The existence of the higher order partial derivatives $\frac{ \partial^{\mathbf{p}} m_{k}(\mathbf{x})}{\partial \mathbf{x}^{\mathbf{p}}}$ can then be verified recursively. 
\end{proof}

\chapter[]{Proof of Theorem \ref{thm_LMVU_SDPCM_S_1}}
\label{chap_appendix_B}

\begin{proof} 
We prove the statement separately for the two complementary cases where $k = j_{0}$ and the case where $k \neq j_{0}$.
First, consider the case $k = j_{0}$, where we have obviously that $| \{k\} \cup \supp(\mathbf{x}_{0}) | = | \{k \}| = 1 < S+1$. This implies that the estimator \eqref{equ_LMVU_SDPCM_S_1} equals the estimator \eqref{equ_tight_est_SDPCM}, which is unbiased (cf.\ \eqref{equ_mean_tight_est_SDPCM}) and whose variance achieves the bound \eqref{equ_corr_lower_bound} (for the case $k \in \supp(\mathbf{x}_{0})$) of Corollary \ref{corr_SDCM_unbiased} (cf.\ \eqref{equ_spdcm_variance_tight_est}). Therefore, the estimator \eqref{equ_LMVU_SDPCM_S_1} must be 
the LMVU for the case $k = j_{0}$.

Next, we consider the complementary case $k \neq j_{0}$. 
We calculate the mean $\mathsf{E}_{x} \{ \hat{x}_{k}^{(\mathbf{x}_{0})}(\mathbf{y}) \}$ separately for parameter vectors $\mathbf{x}$ such that $\supp(\mathbf{x}) = j_{0}$ 
and parameter vectors such that $\supp(\mathbf{x}) \neq j_{0}$. 
\begin{itemize}
\item $\supp(\mathbf{x}) \neq  j_{0}$: 
We obtain for the mean of $\hat{x}_{k}^{(\mathbf{x}_{0})}(\mathbf{y})$
\begin{align} 
\label{equ_appendix_B_k_neq_j_0_supp_x_neq_j_0}
\mathsf{E}_{\mathbf{x}} \big\{\hat{x}^{(\mathbf{x}_{0})}_{k}(\mathbf{y}) \big\} & = \mathsf{E}_{\mathbf{x}} \big\{   \alpha(\mathbf{y};\mathbf{x}_{0}) \ist \big( \beta_k(\mathbf{y}) \rmv-\rmv \sigma^{2} \big) \big\}
\stackrel{(a)}{=}  \mathsf{E}_{\mathbf{x}} \big\{   \alpha(\mathbf{y};\mathbf{x}_{0}) \big\}  \mathsf{E}_{\mathbf{x}} \big\{   \beta_k(\mathbf{y}) \rmv-\rmv \sigma^{2}  \big\} \nonumber \\[4mm] 
 & \stackrel{\eqref{equ_mean_tight_est_SDPCM}}{=} x_{k}  \mathsf{E}_{\mathbf{x}} \big\{   \alpha(\mathbf{y};\mathbf{x}_{0}) \big\}  = x_{k}  \mathsf{E}_{\mathbf{x}} \big\{  a(\mathbf{x}_{0}) \exp\rmv\rmv\big(\!\rmv-\rmv\rmv r_{j_{0}} b(\mathbf{x}_{0}) \ist \beta_{j_{0}}(\mathbf{y})\big)\big\} 
 \nonumber \\[4mm] 
 & =x_{k} a(\mathbf{x}_{0}) \mathsf{E}_{\mathbf{x}} \big\{   \exp\rmv\rmv\big(\!\rmv-\rmv\rmv r_{j_{0}} b(\mathbf{x}_{0}) \ist \beta_{j_{0}}(\mathbf{y})\big) \big\}  \nonumber \\[4mm]
& \stackrel{\eqref{equ_tight_est_SDPCM}}{=} x_{k} a(\mathbf{x}_{0}) \mathsf{E}_{\mathbf{x}} \big\{   \exp\rmv\rmv\big(\!\rmv-\rmv\rmv r_{j_{0}} b(\mathbf{x}_{0}) \ist\frac{1}{r_{j_{0}}} \rmv\sum_{iÊ\in [r_{j_{0}}]} \rmv\big( \mathbf{u}^{T}_{m_{j_{0},i}} \ist \mathbf{y} \big)^{2}\big) \big\} 
 \nonumber \\[4mm]
 & \stackrel{(b)}{=} x_{k} a(\mathbf{x}_{0}) \prod_{i \in [r_{j_{0}}]} \mathsf{E}_{\mathbf{x}} \big\{   \exp\rmv\rmv\big(\!\rmv-\rmv\rmv b(\mathbf{x}_{0}) z_{i}^{2}\big) \big\}, 
\end{align} 
where $z_{i} \triangleq  \mathbf{u}^{T}_{m_{j_{0},i}} \ist \mathbf{y}$ (cf.\ \eqref{equ_basis_matrix_sum_rank_1_projection_SDPCM}). Note that the random variables $\{ z_{i} \}_{i \in [r_{j_{0}}]}$ 
are i.i.d.\ with $z_{i} \sim \mathcal{N}(0 , \sigma^{2})$, which can be verified by \eqref{equ_expr_obs_cov_matrix_SDPCM} and the fact that we consider the case $\supp(\mathbf{x}) \neq j_{0}$.
The steps $(a)$ and $(b)$ are due to the statistical independence of the random variables $\{ \beta_{l}(\mathbf{y}) \}_{l \in [N]}$,  which can be verified by \eqref{equ_expr_obs_cov_matrix_SDPCM}. 
We now use the identity 
\begin{equation} 
\label{equ_proof_LMVU_SDPCM_exp_expect}
\mathsf{E} \{ \exp(- B z^{2}) \} =  \big[ (2B + c^{-2}) c^{2}Ê\big ]^{-1/2}, 
\end{equation} 
which is valid for a Gaussian random variable $z \sim \mathcal{N}(0,c^{2})$ with arbitrary variance $c^{2}$ \cite{papoulis}, to obtain from \eqref{equ_appendix_B_k_neq_j_0_supp_x_neq_j_0} that 
\begin{align}
\mathsf{E}_{\mathbf{x}} \big\{\hat{x}^{(\mathbf{x}_{0})}_{k}(\mathbf{y}) \big\}  &  = x_{k} a(\mathbf{x}_{0}) \prod_{i \in [r_{j_{0}}]} \big[ (2b(\mathbf{x}_{0}) + \sigma^{-2}) \sigma^{2}Ê\big ]^{-1/2} \nonumber \\[4mm]
& = x_{k} a(\mathbf{x}_{0}) \big[ (2b(\mathbf{x}_{0}) + \sigma^{-2}) \sigma^{2}Ê\big ]^{-r_{j_{0}}/2} 
   \nonumber \\[4mm] 
 & =x_{k} a(\mathbf{x}_{0}) \big[ (2 (\sigma^{-2} - (\xi_{0} + \sigma^{2})^{-1})/2  + \sigma^{-2}) \sigma^{2}Ê\big ]^{-r_{j_{0}}/2}  \nonumber \\[4mm]
  & =x_{k} a(\mathbf{x}_{0}) \big[ 1 - \sigma^{2}(\xi_{0} + \sigma^{2})^{-1} + 1 Ê\big ]^{-r_{j_{0}}/2}  \nonumber \\[4mm]
    & =x_{k} a(\mathbf{x}_{0}) \big[2 - \sigma^{2}(\xi_{0} + \sigma^{2})^{-1}  Ê\big ]^{-r_{j_{0}}/2}  \nonumber \\[4mm]
        & =x_{k} a(\mathbf{x}_{0}) \big[2(\xi_{0} + \sigma^{2}) - \sigma^{2} \big ]^{-r_{j_{0}}/2} (\xi_{0}+\sigma^{2})^{r_{j_{0}}/2}  \nonumber \\[4mm]
& =  x_{k} a(\mathbf{x}_{0}) (2\xi_{0}+\sigma^{2})^{-r_{j_{0}}/2}(\xi_{0}+\sigma^{2})^{r_{j_{0}}/2} 
  = x_{k}.
\end{align}
\item $\supp(\mathbf{x}) =  j_{0}$: 
In this case, we have $k \neq \supp(\mathbf{x})$ (since $k \neq j_{0}$ is assumed) and thus $x_{k}=0$. Hence, we obtain for the mean of $\hat{x}_{k}^{(\mathbf{x}_{0})}(\mathbf{y})$ 
\begin{align} 
\mathsf{E}_{\mathbf{x}} \big\{ \hat{x}^{(\mathbf{x}_{0})}_{k}(\mathbf{y}) \big\} & = \mathsf{E}_{\mathbf{x}} \big\{   \alpha(\mathbf{y};\mathbf{x}_{0}) \ist \big( \beta_k(\mathbf{y}) \rmv-\rmv \sigma^{2} \big) \big\} \nonumber \\[4mm]
& \stackrel{(a)}=  \mathsf{E}_{\mathbf{x}} \big\{   \alpha(\mathbf{y};\mathbf{x}_{0}) \big\}  \mathsf{E}_{\mathbf{x}} \big\{  \beta_k(\mathbf{y}) \rmv-\rmv \sigma^{2}  \big\} \stackrel{(b)}{=} 0 = x_{k}, 
\end{align}
where $(a)$ follows from the statistical independence of the random variables $\{ \beta_{l}(\mathbf{y}) \}_{l \in [N]}$ and $(b)$ is due to 
$\mathsf{E}_{\mathbf{x}} \big\{  \beta_k(\mathbf{y}) \rmv-\rmv \sigma^{2}  \big\} = x_{k}=0$ (again see \eqref{equ_mean_tight_est_SDPCM}). 
\end{itemize}
To summarize, we have shown that for the case $k \neq j_{0}$, the estimator \eqref{equ_LMVU_SDPCM_S_1} is unbiased for every $\mathbf{x} \in \mathcal{X}_{S=1,+}$. 

Finally, we show that the variance at $\mathbf{x}_{0}$ of the estimator \eqref{equ_LMVU_SDPCM_S_1} achieves the lower bound \eqref{equ_corr_lower_bound} (for the case $k \notin \supp(\mathbf{x}_{0})$) of Corollary \ref{corr_SDCM_unbiased}. 
Indeed, we obtain under the assumption $k \neq j_{0} = \supp(\mathbf{x}_{0})$ and using the shorthand $z_{i} \triangleq  \mathbf{u}^{T}_{m_{j_{0},i}} \ist \mathbf{y}$  
\begin{align} 
\label{equ_appendix_b_k_neq_j_0_variance_hat_x_part_1}
v(\hat{x}^{(\mathbf{x}_{0})}_{k}(\cdot); \mathbf{x}_{0}) & \stackrel{(a)}{=} P(\hat{x}^{(\mathbf{x}_{0})}_{k}(\cdot); \mathbf{x}_{0}) \nonumber \\[4mm]
& = \mathsf{E}_{\mathbf{x}_{0}} \big\{ \alpha^{2}(\mathbf{y};\mathbf{x}_{0}) \ist \big( \beta_k(\mathbf{y}) \rmv-\rmv \sigma^{2} \big)^{2} \big\} \nonumber \\[4mm]
& \stackrel{(b)}{=}  \mathsf{E}_{\mathbf{x}_{0}} \big\{ \alpha^{2}(\mathbf{y};\mathbf{x}_{0}) \big\}  \mathsf{E}_{\mathbf{x}_{0}} \big\{  \big( \beta_k(\mathbf{y}) \rmv-\rmv \sigma^{2} \big)^{2} \big\}.
\end{align} 
Here, for step $(a)$ we used the fact that the estimator $\hat{x}_{k}^{(\mathbf{x}_{0})}(\cdot)$ is unbiased as has been proven just before and step $(b)$ is due the statistical 
independence of the random variables $\{ \beta_{l}(\mathbf{y}) \}_{l \in [N]}$.
We have $\mathsf{E}_{\mathbf{x}_{0}} \big\{  \big( \beta_k(\mathbf{y}) \rmv-\rmv \sigma^{2} \big)^{2} \big\} = \frac{2 \sigma^{4}}{r_{k}}$ because of \eqref{equ_var_tight_est_SDPCM}. Thus, we obtain from \eqref{equ_appendix_b_k_neq_j_0_variance_hat_x_part_1} that  
\begin{align}
v(\hat{x}^{(\mathbf{x}_{0})}_{k}(\cdot); \mathbf{x}_{0})& = \frac{2 \sigma^{4}}{r_{k}}  \mathsf{E}_{\mathbf{x}_{0}} \big\{ \alpha^{2}(\mathbf{y};\mathbf{x}_{0}) \big\} 
\nonumber \\[4mm] 
& \stackrel{(a)}{=}  \frac{2 \sigma^{4}}{r_{k}}  a^{2}(\mathbf{x}_{0}) \prod_{i \in [r_{j_{0}}]} \mathsf{E}_{\mathbf{x}_{0}} \big\{  \exp\rmv\rmv\big(\!\rmv-\rmv\rmv 2b(\mathbf{x}_{0}) z_{i}^{2}\big) \big \} \nonumber \\[3mm] 
& \stackrel{\eqref{equ_proof_LMVU_SDPCM_exp_expect}}{=} \frac{2 \sigma^{4}}{r_{k}} a^{2}(\mathbf{x}_{0})  \prod_{i \in [r_{j_{0}}]} \big[ 2 \big[ \sigma^{-2} - (\xi_{0} + \sigma^{2})^{-1}\big] + \big(\xi_{0} + \sigma^{2}\big)^{-1} \big]^{-1/2} \big(\xi_{0} + \sigma^{2}\big)^{-1/2} \nonumber \\[4mm]
& = \frac{2 \sigma^{4}}{r_{k}} a^{2}(\mathbf{x}_{0})  \prod_{i \in [r_{j_{0}}]} \big[ 2  \sigma^{-2} \big(\xi_{0} + \sigma^{2}\big) - 1 \big]^{-1/2} \nonumber \\[4mm]  
& = \frac{2 \sigma^{4}}{r_{k}} a^{2}(\mathbf{x}_{0}) \big[ 2  \sigma^{-2} \big(\xi_{0} + \sigma^{2}\big) - 1 \big]^{-r_{j_{0}}/2} \nonumber \\[4mm]  
& =  \frac{2 \sigma^{4}}{r_{k}}a^{2}(\mathbf{x}_{0})  \sigma^{r_{j_{0}}} \big(2\xi_{0} + \sigma^{2}\big)^{-r_{j_{0}}/2}   \nonumber \\[4mm]
& = \frac{2 \sigma^{4}}{r_{k}}  \bigg[ \frac{2 \xi_{0} + \sigma^{2}}{\xi_{0} + \sigma^{2}} \bigg]^{r_{j_{0}}}  \sigma^{r_{j_{0}}}  \big(2\xi_{0} + \sigma^{2}\big)^{-r_{j_{0}}/2}  
\nonumber \\[4mm] 
& =  \frac{2 \sigma^{4}}{r_{k}}   \frac{\big( 2 \xi_{0}\sigma^{2} + \sigma^{4}\big)^{r_{j_{0}}/2}}{\big( \xi_{0} + \sigma^{2} \big)^{r_{j_{0}}}} \nonumber \\[4mm]
& =  \frac{2 \sigma^{4}}{r_{k}}   \frac{\big[  (\xi_{0}+\sigma^{2})^{2} - \xi_{0}^{2} \big]^{r_{j_{0}}/2}}{\big( \xi_{0} + \sigma^{2} \big)^{r_{j_{0}}}}  \nonumber \\[4mm]
& =  \frac{2} {r_{k}} \sigma^{4} \Bigg[ 1-  \frac{\xi_{0}^{2}}{(\xi_{0} + \sigma^{2})^{2}} \Bigg]^{\frac{r_{j_{0}}}{2}}, 
\end{align} 
which is \eqref{equ_corr_lower_bound} (for the case $k \notin \supp(\mathbf{x}_{0})$).
Here, the step $(a)$ is due to the fact that the random variables 
$\{ z_{i} \}_{i \in [r_{j_{0}}]}$ are i.i.d.\ with $z_{i} \sim \mathcal{N}(0 , \sigma^{2}+ \xi_{0})$. 
\end{proof}

\chapter[]{Proof of Theorem \ref{thm_ml_est_SDPCM}}
\label{chap_appendix_C}

\begin{proof} 
Consider the function 
\begin{equation} 
\label{equ_appendix_C_def_h_k_x_function}
h_{k}(x): \mathbb{R}_{+} \rightarrow \mathbb{R}: h_{k}(x) \triangleq  -r_{k} \bigg[ \frac{\beta_{k}(\mathbf{y})}{x+\sigma^{2}} + \log\big( x+\sigma^{2} \big) \bigg].
\end{equation} 
This function is continuous on its domain and differentiable at every point $x > 0$, with derivative 
\begin{equation}
\label{proof_spdcm_ML_derivative_h_k}
h'_{k}(x) = r_{k}   \bigg[  \frac{\beta_{k}(\mathbf{y})}{(x+\sigma^{2})^2} -  \frac{1}{x+\sigma^{2}} \bigg].
\end{equation}
As can be verified easily, if $\beta_{k}(\mathbf{y}) < \sigma^{2}$, then the derivative in \eqref{proof_spdcm_ML_derivative_h_k} is always negative, which implies 
via the mean value theorem \cite{RudinBookPrinciplesMatheAnalysis} that the function $h_{k}(x)$ is monotonically decreasing and attains its maximum at $x=0$. Thus, we have 
\begin{equation}
\label{equ_max_h_k_beta_k_lower_sigma_2}
\beta_{k}(\mathbf{y}) < \sigma^{2} \quad \Rightarrow \quad \max_{x \in \mathbb{R}_{+}} h_{k}(x) = h_{k}(0) =  -r_{k} \bigg[ \frac{\beta_{k}(\mathbf{y})}{\sigma^{2}} +  \log\big( \sigma^{2} \big) \bigg] .
\end{equation}
On the other hand, if  $\beta_{k}(\mathbf{y}) \geq \sigma^{2}$, then we can distinguish two cases: (i) for $xÊ\in [0, \beta_{k}(\mathbf{y}) - \sigma^{2}]$, the derivative $h'_{k}(x)$ is nonnegative, implying that 
$h_{k}(x)$ is monotonically increasing; (ii) for $x \geq \beta_{k}(\mathbf{y}) - \sigma^{2}$, the derivative is negative, implying that $h_{k}(x)$ is monotonically decreasing. These 
two facts imply that the maximum of $h_{k}(x)$ must occur at $x=\beta_{k}(\mathbf{y})-\sigma^{2}$, i.e., 
\begin{equation}
\label{equ_max_h_k_beta_k_geq_sigma_2}
\beta_{k}(\mathbf{y}) \geq \sigma^{2} \quad \Rightarrow \quad \max_{x \in \mathbb{R}_{+}} h_{k}(x) = h_{k}(\beta_{k}(\mathbf{y})- \sigma^{2}) =  -r_{k} \big[ 1+  \log\big( \beta_{k}(\mathbf{y}) \big) \big] .
\end{equation}
In what follows we will use the shorthand 
\begin{equation}
\label{equ_appendix_C_def_H_sum_h_k} 
H(\mathbf{x}) \triangleq \sum_{k \in [N]} h_{k}(x_{k}), 
\end{equation}
which allows us to rewrite the ML estimator (cf.\ \eqref{equ_SDPCM_ML_est_expression_sum_h_k}) as 
\begin{equation} 
\label{equ_appendix_C_expr_ML_est_H}
\hat{\mathbf{x}}_{\text{ML}}(\mathbf{y}) = \argmax_{\mathbf{x}' \in \mathcal{X}_{S,+}} \Bigg\{ \sum_{k \in [N]}Êh_{k}(x'_{k})  \Bigg\} = \argmax_{\mathbf{x}' \in \mathcal{X}_{S,+}} H(\mathbf{x}').
\end{equation} 

Let us now show by contradiction that $\supp(\hat{\mathbf{x}}_{\text{ML}}(\mathbf{y})) \subseteq \mathcal{L}_{2}$, i.e., $\hat{x}_{k,\text{ML}}(\mathbf{y})  = 0$ for every index $k \notin \mathcal{L}_{2}$. 
Let $\tilde{\mathbf{x}} \triangleq\hat{\mathbf{x}}_{\text{ML}}(\mathbf{y})$, and assume that the ML estimate $\tilde{\mathbf{x}}$ has a nonzero entry $\tilde{x}_{k'}$ 
with index $k'  \notin \mathcal{L}_{2}$. 
Consider then a new parameter vector $\mathbf{x}_{2} \in \mathcal{X}_{S,+}$ which is obtained from $\tilde{\mathbf{x}}$ by zeroing the entry with index $k'$ and retaining the rest. We can 
then write the difference of the objective values $H(\mathbf{x})$ for the parameter vectors $\tilde{\mathbf{x}}$ and $\mathbf{x}_{2}$ as 
\begin{align}
\label{equ_appendix_C_diff_H_tilde_x_H_x_2_part_1}
H(\tilde{\mathbf{x}}) - H(\mathbf{x}_{2}) & \stackrel{\eqref{equ_appendix_C_def_H_sum_h_k}}{=} \Bigg\{ \sum_{k \in [N]} h_{k}(\tilde{x}_{k}) \Bigg\}  -\Bigg\{ \sum_{k \in [N] } h_{k}(x_{2,k}) \Bigg\} \nonumber \\[4mm]
& = \Bigg\{ \sum_{k \in [N] \setminus \{k'\}} h_{k}(\tilde{x}_{k}) \Bigg\}  + h_{k'}(\tilde{x}_{k'}) -\Bigg\{ \sum_{k \in [N] \setminus \{k'\}} h_{k}(x_{2,k}) \Bigg\}- h_{k'}(x_{2,k'}) \nonumber \\[4mm]
& = h_{k'}(\tilde{x}_{k'}) - h_{k'}(x_{2,k'}) \nonumber \\[4mm]
& =  h_{k'}(\tilde{x}_{k'}) - h_{k'}(0). 
\end{align} 
Because $k' \notin \mathcal{L}_{2}$ we have $\beta_{k'}(\mathbf{y}) < \sigma^{2}$ and in turn that the derivative \eqref{proof_spdcm_ML_derivative_h_k} (for $k=k'$) is strictly negative for all $x \in \mathbb{R}_{+}$. 
It then follows from \eqref{equ_appendix_C_diff_H_tilde_x_H_x_2_part_1} and \eqref{equ_max_h_k_beta_k_lower_sigma_2} that
\begin{align}
H(\tilde{\mathbf{x}}) - H(\mathbf{x}_{2})< 0.
\end{align} 
Thus, we have that the new parameter vector $\mathbf{x}_{2}$ yields a larger value of the objective $H(\mathbf{x})$ than $\tilde{\mathbf{x}}$. However, this contradicts the 
assumption that $\tilde{\mathbf{x}}$ is the ML estimator, i.e., the parameter vector in $\mathcal{X}_{S,+}$ that maximizes $H(\mathbf{x})$. Thus, we have shown that $\supp(\hat{\mathbf{x}}_{\text{ML}}(\mathbf{y})) \subseteq \mathcal{L}_{2}$, i.e., $\hat{x}_{k,\text{ML}}(\mathbf{y})  = 0$ for every index $k \notin \mathcal{L}_{2}$.

Let $\mathcal{I}$ be a given set of indices such that $\mathcal{I} \subseteq \mathcal{L}_{2}$ and $|\mathcal{I}| \leq S$, and let us denote by $H_{\mathcal{I}}$ the maximum value of $H(\mathbf{x}')$ under the additional constraint that $\supp(\mathbf{x}') \subseteq \mathcal{I}$, i.e., 
\begin{align}
\label{equ_ml_sdpcm_maximizer_restricted}
H_{\mathcal{I}} &  \triangleq \max_{\substack{\mathbf{x}' \in \mathbb{R}_{+}^{N} \\ \supp(\mathbf{x}') \subseteq \mathcal{I}}} H(\mathbf{x}')  \nonumber \\[4mm]
& = \max_{\substack{\mathbf{x}' \in \mathbb{R}_{+}^{N} \\ \supp(\mathbf{x}') \subseteq \mathcal{I}}}   \sum_{k \in [N]} h_{k}(x'_{k}) \nonumber \\[4mm]
& =  \max_{\substack{\mathbf{x}' \in \mathbb{R}_{+}^{N} \\ \supp(\mathbf{x}') \subseteq \mathcal{I}}} \sum_{k \notin \mathcal{I}} h_{k}(0) +  \sum_{k \in \mathcal{I}} h_{k}(x'_{k})  \nonumber \\[4mm]
&\stackrel{\eqref{equ_max_h_k_beta_k_geq_sigma_2}}{=}\sum_{k \notin \mathcal{I}} h_{k}(0)  + \sum_{k \in \mathcal{I}} \bigg\{  -r_{k} \big[ 1+  \log\big( \beta_{k}(\mathbf{y}) \big) \big] \bigg\}.  \\[-4mm]
\end{align} 
 
Up to now, we have verified that the ML estimate $\tilde{\mathbf{x}}$ must satisfy $\supp(\tilde{\mathbf{x}}) \subseteq \mathcal{L}_{2}$. 
We will now show that the vector $\mathbf{x}_{3}$ given by \eqref{equ_def_ML_est_SDPCM_expr} yields the ML estimate, i.e., $H(\tilde{\mathbf{x}}) = H(\mathbf{x}_{3})$. Note 
that there may be multiple ML estimates $\tilde{\mathbf{x}}$ which yield the global maximum in \eqref{equ_SDPCM_ML_est_expression_sum_h_k}; we will show that $\mathbf{x}_{3}$ yields one of these solutions. 
By the definition of $\mathbf{x}_{3}$, 
\begin{align}
\label{equ_appendix_c_def_H_x_3}
H(\mathbf{x}_{3}) & \stackrel{\eqref{equ_appendix_C_def_H_sum_h_k}}{=}  \sum_{k \in [N]} h_{k}(x_{3,k}) \nonumber \\[4mm]
 & \stackrel{\eqref{equ_def_ML_est_SDPCM_expr}}{=} \sum_{k \notin \mathcal{L}_{1} \cap \mathcal{L}_{2}} h_{k}(0)  + \sum_{k \in \mathcal{L}_{1} \cap \mathcal{L}_{2}} h_{k}(\beta_{k}(\mathbf{y}) - \sigma^{2}) \nonumber \\[4mm]
& \stackrel{\eqref{equ_appendix_C_def_h_k_x_function}}{=} \sum_{k \notin \mathcal{L}_{1} \cap \mathcal{L}_{2}} h_{k}(0)  + \sum_{k \in \mathcal{L}_{1} \cap \mathcal{L}_{2}} \bigg\{  -r_{k} \big[ 1+  \log\big( \beta_{k}(\mathbf{y}) \big) \big] \bigg\} \nonumber \\[4mm]
& \stackrel{\eqref{equ_ml_sdpcm_maximizer_restricted}}{=} H_{\mathcal{L}_{1} \cap \mathcal{L}_{2}}.
\end{align} 

Consider then an arbitrary parameter vector $\mathbf{x}' \in \mathcal{X}_{S,+}$ which satisfies $\supp(\mathbf{x}') \subseteq \mathcal{L}_{2}$ and denote 
its support by $\mathcal{I}' Ê\triangleq \supp(\mathbf{x}')$. Obviously, we have that $\mathcal{I}' \subseteq \mathcal{L}_{2}$ and $|\mathcal{I}'| \leq S$ which implies 
that 
\begin{equation}
\label{equ_appendix_C_cardinality_I_prime_bound}
| \mathcal{I}' | \leq \min \{ S , |\mathcal{L}_{2}| \}. 
\end{equation} 
The difference of the objective values $H(\mathbf{x})$ for the parameter 
vectors $\mathbf{x}_{3}$ and $\mathbf{x}'$ then satisfies  
\begin{align}
\label{equ_appendix_C_diff_H_x_3_H_x_prime}
 H(\mathbf{x}_{3}) - H(\mathbf{x}') & \stackrel{\eqref{equ_appendix_c_def_H_x_3}}{=}  H_{\mathcal{L}_{1} \cap \mathcal{L}_{2}} - H(\mathbf{x}')
\stackrel{\eqref{equ_ml_sdpcm_maximizer_restricted}}{\geq} 
H_{\mathcal{L}_{1} \cap \mathcal{L}_{2}} - H_{\mathcal{I}'} \nonumber \\[4mm]
 &  \stackrel{\eqref{equ_ml_sdpcm_maximizer_restricted}}{=} 
 \sum_{k \in \big(\mathcal{L}_{1}\cap \mathcal{L}_{2}\big)} \bigg\{- r_{k} \big[ 1+  \log\big( \beta_{k}(\mathbf{y}) \big) \big]\bigg\}  + \sum_{k \notin  \big(\mathcal{L}_{1}\cap \mathcal{L}_{2}\big)}
 h_{k}(0) \nonumber \\[4mm]
& \hspace*{10mm}-  \sum_{k \in  \mathcal{I}' }\bigg\{ - r_{k} \big[ 1+  \log\big( \beta_{k}(\mathbf{y}) \big) \big]\bigg\}  - \sum_{k \notin \mathcal{I}'}
 h_{k}(0)\nonumber \\[4mm]
 &  = 
 \sum_{k \in \big(\mathcal{L}_{1}\cap \mathcal{L}_{2} \cap \mathcal{I}' \big)} \bigg\{- r_{k} \big[ 1+  \log\big( \beta_{k}(\mathbf{y}) \big) \big]\bigg\} +
  \sum_{k \in \big(\mathcal{L}_{1}\cap \mathcal{L}_{2}\big) \setminus \mathcal{I}'} \bigg\{- r_{k} \big[ 1+  \log\big( \beta_{k}(\mathbf{y}) \big) \big]\bigg\}  \nonumber \\[4mm]
 & \hspace*{10mm}  + \sum_{k \in  [N] \setminus \big((\mathcal{L}_{1}\cap \mathcal{L}_{2}) \cup \mathcal{I}'\big)}h_{k}(0)  + \sum_{k \in    \mathcal{I}' \setminus\big(\mathcal{L}_{1}\cap \mathcal{L}_{2}\big)}
 h_{k}(0) \nonumber \\[4mm]
& \hspace*{10mm}-  \sum_{k \in \big(\mathcal{L}_{1}\cap \mathcal{L}_{2} \cap \mathcal{I}' \big)} \bigg\{- r_{k} \big[ 1+  \log\big( \beta_{k}(\mathbf{y}) \big) \big]\bigg\} -
  \sum_{k \in   \mathcal{I}' \setminus \big(\mathcal{L}_{1}\cap \mathcal{L}_{2}\big)} \bigg\{- r_{k} \big[ 1+  \log\big( \beta_{k}(\mathbf{y}) \big) \big]\bigg\}  \nonumber \\[4mm]
& \hspace*{10mm} -  \sum_{k \in  [N] \setminus \big((\mathcal{L}_{1}\cap \mathcal{L}_{2}) \cup \mathcal{I}'\big)}h_{k}(0)  - \sum_{k \in  \big(\mathcal{L}_{1}\cap \mathcal{L}_{2}\big) \setminus  \mathcal{I}' }
 h_{k}(0)\nonumber \\[4mm]
 &  \stackrel{\eqref{equ_ml_sdpcm_maximizer_restricted}}{=} 
 \sum_{k \in \big(\mathcal{L}_{1}\cap \mathcal{L}_{2}\big) \setminus \mathcal{I}'} \bigg\{ - r_{k} \big[ 1+  \log\big( \beta_{k}(\mathbf{y}) \big) \big] -h_{k}(0) \bigg\} \nonumber \\[4mm]
& \hspace*{10mm}-  \sum_{k \in  \mathcal{I}' \setminus\big(\mathcal{L}_{1}\cap \mathcal{L}_{2}\big) } \bigg\{ - r_{k} \big[ 1+  \log\big( \beta_{k}(\mathbf{y}) \big) \big] -h_{k}(0) \bigg \} \nonumber \\[4mm]
&  = \sum_{k \in \big(\mathcal{L}_{1}\cap \mathcal{L}_{2}\big) \setminus \mathcal{I}'}  \bigg \{ - r_{k} \big[ 1+  \log\big( \beta_{k}(\mathbf{y}) \big) \big] +
 r_{k} \bigg[ \frac{\beta_{k}(\mathbf{y})}{\sigma^{2}} +  \log\big( \sigma^{2} \big) \bigg] \bigg\} \nonumber \\[4mm]
& \hspace*{10mm} -  \sum_{k \in  \mathcal{I}' \setminus\big(\mathcal{L}_{1}\cap \mathcal{L}_{2}\big) } \bigg\{- r_{k} \big[ 1+  \log\big( \beta_{k}(\mathbf{y}) \big) \big] +
r_{k} \bigg[ \frac{\beta_{k}(\mathbf{y})}{\sigma^{2}} +  \log\big( \sigma^{2} \big) \bigg] \bigg \} \nonumber \\[4mm]
&  = \sum_{k \in \big(\mathcal{L}_{1}\cap \mathcal{L}_{2}\big) \setminus \mathcal{I}'} r_{k} \bigg[ \frac{\beta_{k}(\mathbf{y})}{\sigma^2} \! -Ê\! \log \ist \bigg( \frac{\beta_{k}(\mathbf{y})}{\sigma^{2}} \bigg) \rmv-\! 1 \bigg]\nonumber \\[4mm]
&Ê\hspace*{10mm} - \sum_{k' \in  \mathcal{I}' \setminus\big(\mathcal{L}_{1}\cap \mathcal{L}_{2}\big) }r_{k'} \bigg[ \frac{\beta_{k'}(\mathbf{y})}{\sigma^2} \! -Ê\! \log \ist \bigg( \frac{\beta_{k'}(\mathbf{y})}{\sigma^{2}} \bigg) \rmv-\! 1 \bigg] 
 \end{align}
Note that by definition $| \mathcal{L}_{1} \cap \mathcal{L}_{2}| = \min \{ S, |\mathcal{L}_{2}| \}$, implying due to \eqref{equ_appendix_C_cardinality_I_prime_bound} that
\begin{equation}
 | \mathcal{L}_{1} \cap \mathcal{L}_{2}| \geq |\mathcal{I}'|, 
\end{equation}   
and in turn 
\begin{equation} 
\label{equ_appendix_c_card_L1_L2_I_prime}
 |  \big(\mathcal{L}_{1}\cap \mathcal{L}_{2}\big) \setminus \mathcal{I}' | \geq | \mathcal{I}' \setminus\big(\mathcal{L}_{1}\cap \mathcal{L}_{2}\big)|. 
\end{equation}
From \eqref{equ_appendix_c_card_L1_L2_I_prime}, we have that there exists an injective index map $\Omega(k'): \mathcal{I}' \setminus\big(\mathcal{L}_{1}\cap \mathcal{L}_{2}\big) \rightarrow  \big(\mathcal{L}_{1}\cap \mathcal{L}_{2}\big) \setminus \mathcal{I}' $, 
such that for each summand (with index $k'  \in \mathcal{I}' \setminus\big(\mathcal{L}_{1}\cap \mathcal{L}_{2}\big)$) of the second summand in 
\eqref{equ_appendix_C_diff_H_x_3_H_x_prime}, we can associate a distinct summand (with index $\Omega(k') \in \big(\mathcal{L}_{1}\cap \mathcal{L}_{2}\big) \setminus \mathcal{I}'$) 
of the first sum in \eqref{equ_appendix_C_diff_H_x_3_H_x_prime}.
However, for any $\Omega(k') \in \big(\mathcal{L}_{1}\cap \mathcal{L}_{2}\big) \setminus \mathcal{I}'$ (implying that $\Omega(k') \in \mathcal{L}_{1} \cap \mathcal{L}_{2}$) 
and $k'  \in \mathcal{I}' \setminus\big(\mathcal{L}_{1}\cap \mathcal{L}_{2}\big)$ (implying that $k' \in \mathcal{L}_{2} \setminus \mathcal{L}_{1}$, since we assume that $\mathcal{I}' \subseteq \mathcal{L}_{2}$) we have 
\begin{equation}
\label{equ_appendix_C_ineq_Omage_k_rime_h_k_modified}
 r_{\Omega(k')} \bigg[ \frac{\beta_{\Omega(k')}(\mathbf{y})}{\sigma^2} \! -Ê\! \log \ist \bigg( \frac{\beta_{\Omega(k')}(\mathbf{y})}{\sigma^{2}} \bigg) \rmv-\! 1 \bigg] \geq r_{k'} \bigg[ \frac{\beta_{k'}(\mathbf{y})}{\sigma^2} \! -Ê\! \log \ist \bigg( \frac{\beta_{k'}(\mathbf{y})}{\sigma^{2}} \bigg) \rmv-\! 1 \bigg].
\end{equation} 
Consider the inequality 
\begin{equation}
x - ( \log(x)  + 1) \geq 0, 
\end{equation} 
valid for every $x \geq 1$, which can be verified by inspecting the derivative of the left hand side.
This inequality implies that for any index $k \in \mathcal{L}_{2}$ (i.e., $\beta_{k}(\mathbf{y}) \geq \sigma^{2}$) we have 
\begin{equation}
\label{equ_appendix_C_expr_geq_0_b_k_log}
  r_{k} \bigg[ \frac{\beta_{k}(\mathbf{y})}{\sigma^2} \! -Ê\! \log \ist \bigg( \frac{\beta_{k}(\mathbf{y})}{\sigma^{2}} \bigg) \rmv-\! 1 \bigg] \geq 0.
\end{equation} 
Let us introduce the shorthand $\mathcal{I}'' \triangleq  \big[ \big(\mathcal{L}_{1}\cap \mathcal{L}_{2}\big) \setminus \mathcal{I}'\big] \setminus \Omega\big( \mathcal{I}' \setminus\big(\mathcal{L}_{1}\cap \mathcal{L}_{2}\big) \big)$, and observe that 
\begin{equation} 
\label{equ_appendix_C_set_relation_Omega_1}
\mathcal{I}'' \cup  \Omega\big( \mathcal{I}' \setminus\big(\mathcal{L}_{1}\cap \mathcal{L}_{2}\big) \big) =   \big(\mathcal{L}_{1}\cap \mathcal{L}_{2}\big) \setminus \mathcal{I}',
\end{equation} and
\begin{equation}
\label{equ_appendix_C_set_relation_Omega_2}
\mathcal{I}'' \cap  \Omega\big( \mathcal{I}' \setminus\big(\mathcal{L}_{1}\cap \mathcal{L}_{2}\big) \big) =\emptyset,
\end{equation}
which is due to the obvious fact that $\Omega\big( \mathcal{I}' \setminus\big(\mathcal{L}_{1}\cap \mathcal{L}_{2}\big) \big) \subseteq \big(\mathcal{L}_{1}\cap \mathcal{L}_{2}\big) \setminus \mathcal{I}'$. 
We then have from \eqref{equ_appendix_C_diff_H_x_3_H_x_prime} that 
\begin{align}
  H(\mathbf{x}_{3}) - H(\mathbf{x}') & \stackrel{\eqref{equ_appendix_C_diff_H_x_3_H_x_prime}}{=} 
    \sum_{k \in \big(\mathcal{L}_{1}\cap \mathcal{L}_{2}\big) \setminus \mathcal{I}'} r_{k} \bigg[ \frac{\beta_{k}(\mathbf{y})}{\sigma^2} \! -Ê\! \log \ist \bigg( \frac{\beta_{k}(\mathbf{y})}{\sigma^{2}} \bigg) \rmv-\! 1 \bigg]\nonumber \\[4mm]
&Ê\hspace*{10mm} - \sum_{k' \in  \mathcal{I}' \setminus\big(\mathcal{L}_{1}\cap \mathcal{L}_{2}\big) }r_{k'} \bigg[ \frac{\beta_{k'}(\mathbf{y})}{\sigma^2} \! -Ê\! \log \ist \bigg( \frac{\beta_{k'}(\mathbf{y})}{\sigma^{2}} \bigg) \rmv-\! 1 \bigg]   \nonumber \\[4mm]
  &  \stackrel{\eqref{equ_appendix_C_set_relation_Omega_1},\eqref{equ_appendix_C_set_relation_Omega_2}}{=} \sum_{k \in \mathcal{I}''} r_{k} \bigg[ \frac{\beta_{k}(\mathbf{y})}{\sigma^2} \! -Ê\! \log \ist \bigg( \frac{\beta_{k}(\mathbf{y})}{\sigma^{2}} \bigg) \rmv-\! 1 \bigg]\nonumber \\[4mm]
& \hspace*{10mm} +  \sum_{k \in \Omega\big( \mathcal{I}' \setminus\big(\mathcal{L}_{1}\cap \mathcal{L}_{2}\big)\big)} r_{k} \bigg[ \frac{\beta_{k}(\mathbf{y})}{\sigma^2} \! -Ê\! \log \ist \bigg( \frac{\beta_{k}(\mathbf{y})}{\sigma^{2}} \bigg) \rmv-\! 1 \bigg]\nonumber \\[4mm]
&Ê\hspace*{10mm} - \sum_{k' \in  \mathcal{I}' \setminus\big(\mathcal{L}_{1}\cap \mathcal{L}_{2}\big) }r_{k'} \bigg[ \frac{\beta_{k'}(\mathbf{y})}{\sigma^2} \! -Ê\! \log \ist \bigg( \frac{\beta_{k'}(\mathbf{y})}{\sigma^{2}} \bigg) \rmv-\! 1 \bigg]  
 \nonumber \\[4mm]
  & \hspace*{-20mm}  = \sum_{k \in \mathcal{I}''} r_{k} \bigg[ \frac{\beta_{k}(\mathbf{y})}{\sigma^2} \! -Ê\! \log \ist \bigg( \frac{\beta_{k}(\mathbf{y})}{\sigma^{2}} \bigg) \rmv-\! 1 \bigg]\nonumber \\[4mm]
&Ê\hspace*{-15mm} +  \sum_{k' \in  \mathcal{I}' \setminus\big(\mathcal{L}_{1}\cap \mathcal{L}_{2}\big) } \Bigg\{ r_{\Omega(k')} \bigg[ \frac{\beta_{\Omega(k')}(\mathbf{y})}{\sigma^2} \! -Ê\! \log \ist \bigg( \frac{\beta_{\Omega(k')}(\mathbf{y})}{\sigma^{2}} \bigg) \rmv-\! 1 \bigg] -  r_{k'} \bigg[ \frac{\beta_{k'}(\mathbf{y})}{\sigma^2} \! -Ê\! \log \ist \bigg( \frac{\beta_{k'}(\mathbf{y})}{\sigma^{2}} \bigg) \rmv-\! 1 \bigg] \Bigg\} \nonumber \\[4mm]
& \stackrel{\eqref{equ_appendix_C_ineq_Omage_k_rime_h_k_modified}}{\geq}  \sum_{k \in \mathcal{I}''} r_{k} \bigg[ \frac{\beta_{k}(\mathbf{y})}{\sigma^2} \! -Ê\! \log \ist \bigg( \frac{\beta_{k}(\mathbf{y})}{\sigma^{2}} \bigg) \rmv-\! 1 \bigg]\nonumber \\[4mm]
& \stackrel{\eqref{equ_appendix_C_expr_geq_0_b_k_log}}{\geq} 0.
\end{align}
Thus, we have verified that $H(\mathbf{x}_{3})   \geq H(\mathbf{x}')$ for all $\mathbf{x}' \in \mathcal{X}_{S,+}$ such that $\supp(\mathbf{x}') \subseteq \mathcal{L}_{2}$. 
This means that the vector $\mathbf{x}_{3}$ given by \eqref{equ_def_ML_est_SDPCM_expr} yields the global maximum in \eqref{equ_SDPCM_ML_est_expression_sum_h_k}, i.e., 
it yields an ML estimate.
\end{proof}


\end{appendices}

\bibliographystyle{IEEEtran}
\bibliography{LitAJPHD}

\end{document}